%% file: These.tex
\documentclass[10pt,A4paper]{book}

 \usepackage[latin1]{inputenc}
  \usepackage[french]{varioref}
  \usepackage{psfig}
  \usepackage{epsfig}
  \usepackage{color}
  \usepackage{amsfonts,amsmath,amssymb}
 \usepackage{longtable}

 \usepackage[francais,english]{babel}

\newcommand\english{\selectlanguage{english}}
\newcommand\francais{\selectlanguage{francais}}

  \usepackage{xypic}
  \xyoption{all}

  \usepackage{makeidx}
  \makeindex

\usepackage{fancyhdr}
\fancypagestyle{plain}{%
        \fancyhead{} 
        
}

 \input{Preambule.tex}

\title{ 
    {\bf SUR LES $\ai$-CAT…GORIES}
      }
    \author{ {\bf Kenji LefËvre-Hasegawa}\\
       UFR de MathÈmatiques,
UniversitÈ Paris 7 \\
2, place Jussieu, 
75251 Paris Cedex 05, 
France\\
lefevre@math.jussieu.fr \\
\hbox{www.math.jussieu.fr/$\tilde{\space}$ lefevre}}
\date{mars 2002}
\begin{document}

 \begin{titlepage}
 \thispagestyle{empty}

\thispagestyle{empty}
\begin{center}

{\Large{\bf{UNIVERSITE PARIS 7 - DENIS DIDEROT}}}

\vspace{0.2cm}

{\Large{\bf{U.F.R. de MathÈmatiques}}}

\vspace{1.5cm}

{\large{AnnÈe : 2003}} \hfill {\large{N∞ \ $\begin{array}{|c|c|c|c|c|c|c|c|c|c|}
      \ & \ & \ & \ & \ & \ & \ & \ & \ & \ \\
      \hline
      \end{array}$}}
\vspace{1.5cm}

{\Large{\sc{ThËse de Doctorat}}}
\vspace{0.2cm}

{\large{\bf{SpÈcialitÈ : MathÈmatiques}}}

\vspace{0.8cm}

{\large{prÈsentÈe par}}

\vspace{0.3cm}

{\Large{\it{Kenji LefËvre-Hasegawa}}}

\vspace{0.8cm}

{\large{pour obtenir le grade de}}

\vspace{0.3cm}

{\Large{\sc{Docteur de l'UniversitÈ Paris 7}}}

\vspace{2.2cm}

\begin{tabular}{c}
\hline \\[-5pt]
{\huge{\bf{Sur les $\ai$-catÈgories}}} \\[7pt]
\hline
\end{tabular}

\vspace{2.2cm}

{\large{Soutenue le 6 novembre 2003 devant le jury composÈ de :\\
M.~Johannes {\sc{Huebschmann}} , {\it{rapporteur}}\\
M.~Bernhard {\sc{Keller}}, {\it{directeur}} \\
M.~Pierre {\sc{Cartier}}\\
M.~Alain {\sc{ProutÈ}}\\
M.~RaphaÎl {\sc{Rouquier}} \\
M.~Alexander {\sc{Zimmermann}} \\
}}

\end{center}

  \end{titlepage}

 \thispagestyle{empty}

 \francais
    \tableofcontents

\francais
\chapter*{Remerciements}
\addcontentsline{toc}{chapter}{\hspace{.65cm}Remerciements}
\fancyhead[LO]{}
\fancyhead[RE]{}
\francais

\input{Remerciements.tex}

\francais
\chapter*{Abstract/RÈsumÈ}
\addcontentsline{toc}{chapter}{\hspace{.65cm}Abstract/RÈsumÈ}
\fancyhead[LO]{}
\fancyhead[RE]{}
\francais


   \input{Resume.tex}

\francais

\chapter*{Introduction}
\addcontentsline{toc}{chapter}{\hspace{.65cm}Introduction}
\fancyhead[LO]{Sur les $\ai$-catÈgories}
\fancyhead[RE]{Introduction}
\francais

\input{Introduction.tex}



\francais
 \chapter{ThÈorie de l'homotopie des $\ai$-algËbres}
\fancyhead[LO]{\rightmark}
\fancyhead[RE]{\leftmark}
\label{chapitre_Homot_aialg}

     \input{Homot_aialg.tex}

\francais
\chapter{ThÈorie de l'homotopie des polydules}
\label{chapitre_Homot_aimod}

   \input{Homot_aimod.tex}

\francais
 \chapter{UnitÈs ‡ homotopie prËs et unitÈs strictes}
\label{chapitre_Unite}

    \input{Unite.tex}

\francais
 \chapter{CatÈgorie dÈrivÈe}
\label{chapitre_Cat_der}

    \input{Cat_der.tex}


\francais
\chapter{$\ai$-catÈgories et $\ai$-foncteurs}
\label{chapitre_ai_cat}

    \input{Ai_cat.tex}

\francais
\chapter{Torsion d'$\ai$-structures}
\label{chapitre_torsion}

\input{Torsion.tex}

\francais
\chapter{L'$\ai$-foncteur de Yoneda et les objets tordus}
\label{chapitre_objets_tordus}

\input{Objets_tordus.tex}

\francais
\chapter{L'$\ai$-catÈgorie des $\ai$-foncteurs}
\label{chapitre_ai_func}

      \input{Ai_func.tex}

\francais
\chapter{Les $\ai$-Èquivalences}
\label{chapitre_ai_equiv}

   \input{Ai_equiv.tex}

\francais
\renewcommand{\thechapter}{A}
 \renewcommand{\theproposition}{\thechapter.\arabic{proposition}}
\chapter{CatÈgories de modËles} 
\label{section_rappels_cmf}

\input{Appendice_A.tex}

\francais
\renewcommand{\thechapter}{B}
\renewcommand{\theproposition}{\thesection.\arabic{proposition}}
\chapter{ThÈorie de l'obstruction}
\label{section_theorie_obstruction}

\input{Appendice_B.tex}

\francais
  \addcontentsline{toc}{chapter}{\hspace{.65cm}Bibliographie}

\def\cprime{$'$}


\francais
 \chapter*{Notations}
\fancyhead[LO]{Notations}
\fancyhead[RE]{Notations}
\addcontentsline{toc}{chapter}{\hspace{.65cm}Notations}

\input{Notation_index.tex}


\fancyhead[LO]{}
\fancyhead[RE]{}


\begin{theindex}

  \item {$\ai$-alg\`ebre}, 25
    \subitem topologique, 153
  \item {$\ai$-cat\'egorie}, 135
    \subitem tordue, 148, 158
  \item {$\ai$-cog\`ebre}, 27
  \item {$\ai$-\'equivalence}, 196
  \item {$\ai$-foncteur}, 136
    \subitem de Yoneda, 162
    \subitem de Yoneda g\'en\'eralis\'e, 186
    \subitem pleinement fid\`ele, 162
    \subitem tordu, 149, 159
  \item {$\ai$-isomorphie}, 195
  \item {$\ai$-module}, 79
  \item {$\ai$-morphisme}, 25, 80
    \subitem strict, 25, 80
  \item {$\ai$-pr\'e-triangul\'ee (cat\'egorie)}, 171
  \item {$\ai$-quasi-isomorphisme}, 26, 80
  \item {$\rm{A}_n$-alg\`ebre}, 25
  \item {$\rm{A}_n$-module}, 79
  \item {$\rm{A}_n$-morphisme}, 25, 79
  \item acyclique (cocha\IeC {\^\i }ne), 68
  \item adjonction de Quillen, 204
  \item admissible
    \subitem cocha\IeC {\^\i }ne, 67
    \subitem cog\`ebre, 34
    \subitem comodule, 71
    \subitem filtration, 34
  \item {alg\`ebre}, 22
    \subitem diff\'erentielle gradu\'ee, 23
    \subitem enveloppante, 84, 127
    \subitem gradu\'ee, 23
    \subitem libre, 22
    \subitem presque libre, 23
    \subitem {r\'eduite}, 62
    \subitem tensorielle augment\'ee, 83
    \subitem tensorielle r\'eduite, 22
  \item augmentation, 62, 81

  \indexspace

  \item {${\mathbb A}$-${\mathbb B}$-bimodule}, 133
  \item bipolydule, 92, 137
    \subitem tordu, 151, 159

  \indexspace

  \item {cat\'egorie}
    \subitem {$\ai$-pr\'e-triangul\'ee}, 171
    \subitem de Frobenius, 76
    \subitem de mod\`eles, 202
    \subitem {d\'eriv\'ee}, 76, 88, 94, 118, 126, 129
    \subitem diff\'erentielle gradu\'ee, 135
    \subitem homotopique, 203
    \subitem stable, 76
    \subitem triangul\'ee alg\'ebrique, 172
  \item co-augment\'ee (cog\`ebre), 64
  \item co-augmentation, 64
  \item co-induction, 65
  \item co-unitaire, 64, 65
  \item cocha\IeC {\^\i }ne tordante, 31, 67
    \subitem acyclique, 68
    \subitem admissible, 67
    \subitem g\'en\'eralis\'ee, 124
    \subitem universelle, 69
  \item cocompl\`ete (cog\`ebre), 23
  \item cocomplet (comodule), 66
  \item {cod\'erivation}, 23, 65, 109, 183
  \item cofibrant, 203
  \item cofibration, 202
  \item cofibration standard (de $\alg$), 37
  \item cofibration standard (de $\Modu A$), 72
  \item {cog\`ebre}, 23
    \subitem admissible, 34
    \subitem {co-augment\'ee}, 64
    \subitem cocompl\`ete, 23
    \subitem {r\'eduite}, 64
    \subitem tensorielle r\'eduite, 24
  \item colibre (comodule), 66
  \item comodule, 65
    \subitem admissible, 71
    \subitem cocomplet, 66
    \subitem colibre, 66
    \subitem diff\'erentiel gradu\'e, 65
    \subitem gradu\'e, 65
    \subitem presque colibre, 66
  \item compact (objet), 172
  \item compl\'etion, 153
  \item complexe, 21
    \subitem filtr\'e, 34
  \item {c\^one d'un $\ai$-morphisme}, 90
  \item construction bar, 28, 68, 83, 93
  \item construction cobar, 30, 68
  \item contractant, 153
  \item contraction, 54
  \item corestriction, 65
  \item cylindre, 202

  \indexspace

  \item {d\'erivation}, 22, 63
  \item diff\'erentielle tordue, 67

  \indexspace

  \item {\'el\'ement tordant}, 145, 154
  \item {\'equivalence de Quillen}, 204
  \item {\'equivalence faible}, 202
    \subitem {d'$\ai$-foncteurs}, 192
  \item exhaustive (filtration), 34
  \item extension, 164
  \item extension scind\'ee, 164

  \indexspace

  \item fibrant, 203
  \item fibration, 202
  \item filtration, 34
    \subitem admissible, 34
    \subitem {$C$-primitive}, 35
    \subitem exhaustive, 34
    \subitem primitive, 23, 35, 66
  \item foncteur
    \subitem de Quillen, 204
    \subitem {d\'eriv\'e}, 204
    \subitem standard, 113, 114

  \indexspace

  \item {g\'en\'erateur}, 172

  \indexspace

  \item {$H$-unitaire}, 119, 123
  \item hauteur, 164
  \item Hochschild
    \subitem complexe, 214
    \subitem complexe r\'eduit, 215
  \item homotopie, 22--24, 50, 63, 65
    \subitem \`a gauche, 203
    \subitem {$\ai$-morphismes}, 26, 80, 81

  \indexspace

  \item identit\'e
    \subitem d'un objet, 135
    \subitem le foncteur, 136
  \item induction, 63
  \item irr\'eductible, 22
  \item isomorphie, 195

  \indexspace

  \item lemme clef, 140
  \item lemme de perturbation, 54
  \item libre
    \subitem alg\`ebre, 22
    \subitem module, 63

  \indexspace

  \item Maurer-Cartan (\'equation), 100, 146
  \item {mod\`ele diff\'erentiel gradu\'e}, 171
  \item {mod\`ele minimal}, 54, 105
  \item module, 63
    \subitem diff\'erentiel gradu\'e, 63
    \subitem gradu\'e, 63
    \subitem libre, 63
    \subitem presque libre, 63
  \item {${\mathbb A}$-module}, 133

  \indexspace

  \item objet de chemins, 203
  \item objet filtr\'e, 34
  \item objet gradu\'e, 20
  \item objet tordu, 163, 165
  \item obstruction, 207

  \indexspace

  \item perturbation, 54
  \item polydule, 79, 137
    \subitem topologique, 153
  \item presque colibre (comodule), 66
  \item presque libre
    \subitem {alg\`ebre}, 23
    \subitem module, 63
  \item {$n$-primitifs}, 23, 66
  \item produit tensoriel tordu, 67, 125

  \indexspace

  \item quasi-isomorphisme filtr\'e, 34

  \indexspace

  \item {$R$-$\ai$-cat\'egorie}, 154
  \item r\'eduite (alg\`ebre), 62
  \item r\'eduite (cog\`ebre), 64
  \item r\'esolution cofibrante, 203
  \item r\'esolution fibrante, 203
  \item {r\'eduction}, 81
  \item rel\`evement, 201
  \item restriction, 63

  \indexspace

  \item satur\'ee (classe), 202
  \item {scind\'ee (extension)}, 164
  \item {strict ($\ai$-morphisme)}, 25, 80
  \item strictement unitaire
    \subitem {$\ai$-alg\`ebre}, 81
    \subitem bipolydule, 92
    \subitem {\'el\'ement de $\Hom_{\Nunc}(f_1,f_2)$}, 185
    \subitem homotopie, 98
    \subitem polydule, 81
  \item suspension, 21, 79

  \indexspace

  \item {tensorielle augment\'ee (alg\`ebre)}, 83
  \item tensorielle r\'eduite (alg\`ebre), 22
  \item tensorielle r\'eduite (cog\`ebre), 24
  \item tensoriellement nilpotent, 146
  \item topologie, 153
  \item torsion, 148, 149, 151, 158, 159
  \item triviale (cofibration), 202
  \item triviale (fibration), 202

  \indexspace

  \item {unit\'e}
    \subitem homologique, 98
    \subitem stricte, 81
  \item universelle (cocha\IeC {\^\i }ne), 69

  \indexspace

  \item {Yoneda ($\ai$-foncteur de)}, 162
  \item {Yoneda g\'en\'eralis\'e ($\ai$-foncteur de)}, 186

\end{theindex}


\newpage

\input{Page_dos.tex}

\end{document}

%% file: Preambule.tex
\newtheorem{proposition}{Proposition}[section]
\newtheorem{exemple}[proposition]{Exemple}
\newtheorem{theoreme}[proposition]{ThÈorËme}
\newtheorem{definition}[proposition]{DÈfinition}
\newtheorem{lemme}[proposition]{Lemme}
\newtheorem{corollaire}[proposition]{Corollaire}
\newtheorem{remarque}[proposition]{Remarque}

\newtheorem{notation}[proposition]{Notation}

\newcommand{\indexnotation}[1]{\label{#1}}

  \makeatletter
 \renewcommand{\theproposition}{\thesubsection.\arabic{proposition}}
\@addtoreset{proposition}{subsection} \makeatother
  \makeatother

\newcommand{\dem}{ \emph{\mbox{DÈmonstration :\ \ \ }}}
\newcommand{\findem}{\hspace*{.01cm}\hfill $\Box$}

\newcommand{\diagnum}[2]{({\rm #1}) \quad \begin{array}{c} #2 \end{array}}

\newcommand{\ra}{\rightarrow} \newcommand{\la}{\leftarrow}
\newcommand{\Ra}{\Rightarrow}

 \newcommand{\lra}{\leftrightarrow}
\newcommand{\arr}[1]{\stackrel{#1}{\longrightarrow}}

\newcommand{\iso}{\simeq}

\newcommand{\larr}[1]{\stackrel{#1}{\longleftarrow }}

\newcommand{\ca}{{\mathcal A}} \newcommand{\cb}{{\mathcal B}}
\newcommand{\cc}{{\mathcal C}} \newcommand{\cd}{{\mathcal D}}
\newcommand{\ce}{{\mathcal E}} \newcommand{\cf}{{\mathcal F}}
\newcommand{\cg}{{\mathcal G}} \newcommand{\ch}{{\mathcal H}}

\newcommand{\cn}{{\mathcal N}} \newcommand{\co}{{\mathcal O}}

\newcommand{\calt}{{\mathcal T}}

 \renewcommand{\sf}[1]{\mathsf{#1}}

\newcommand{\Ind}{\opname{Ind}\nolimits}
\newcommand{\Res}{\opname{Res}\nolimits}
\newcommand{\colim}{\opname{colim}\nolimits}

\newcommand{\Z}{{\mathbf{Z}}} 
\newcommand{\I}{{\mathbf{I}}} 
 \newcommand{\N}{{\mathbf{N}}}
\newcommand{\corps}{\mathbb{K}}

\newcommand{\Obj}{{\opname{Obj\,}\nolimits}}
 
\newcommand{\id}{{\mathbf{I}}} 
\newcommand{\Id}{\mathbf{1}} 
\newcommand{\im}{\opname{Im}\nolimits}

\newcommand{\coker}{\opname{cok}\nolimits}
\renewcommand{\ker}{\opname{ker}\nolimits}
\newcommand{\vect}{\mathsf{Vect}}
\newcommand{\Func}{{\mathsf{Func}_\infty}}
\newcommand{\Nunc}{{\mathsf{Nunc}_\infty}}

\newcommand{\Ext}{{\mathsf{Ext}}}
\newcommand{\Hom}{{\mathsf{Hom}}}
\newcommand{\Homi}{\raisebox{0ex}[2ex][0ex]{$\overset{\,\infty}{
                              \raisebox{0ex}[1ex][0ex]{$\mathsf{Hom}$}
                                                                 }$}}
\newcommand{\coder}{{\mbox{$\mathsf{coder}$}}}

\newcommand{\opname}[1]{\operatorname{\mathsf{#1}}}

\newcommand{\alg}{\mbox{${\mathsf{Alg}}$}}     
\newcommand{\alga}{\mbox{${\mathsf{Alga}}$}}   
\newcommand{\aia}{\mbox{$\mathsf{Alg}_{\infty}$}} 
\newcommand{\aiaa}{\mbox{$\mathsf{Alga}_{\infty}$}} 
\newcommand{\aiModu}{\opname{Mod_\infty}\nolimits } 
\newcommand{\aiMod}{\opname{Nod_\infty}\nolimits } 
\newcommand{\aiModst}{\opname{Nod_\infty^{strict}}\nolimits } 
\newcommand{\aiModust}{\opname{Mod_\infty^{strict}}\nolimits } 
\newcommand{\Modu}{\opname{Mod}\nolimits }           
\newcommand{\Modush}{\opname{Modsh}\nolimits }           
\newcommand{\unite}{{\eta}} 
\newcommand{\augmentation}{{\mbox{$\cal \epsilon$}}} 
\newcommand{\free}{\opname{free}\nolimits }          

\newcommand{\aic}{\mbox{$\mathsf{Cog}_{\infty}$}} 
\newcommand{\ct}{{T^c}} 
\newcommand{\ctr}{{{\b{T^c}}}} 
\newcommand{\pctr}[1]{\mbox{$\b{T^c_{[#1]}}$}} 
 \newcommand{\coModu}{\opname{Com}\nolimits  }        
\newcommand{\coModcu}{\opname{Comc}\nolimits  }        
\newcommand{\colibre}{\opname{prcol}\nolimits  }        
\newcommand{\cocog}{\mbox{${\mathsf{Cogc}}$}}  
\newcommand{\cocoga}{\mbox{${\mathsf{Cogca}}$}}  
\newcommand{\cog}{{\mbox{${\mathsf{Cog}}$}}}   
\newcommand{\coga}{{\mbox{${\mathsf{Coga}}$}}}   
\newcommand{\cogtr}{{\mbox{${\mathsf{Cogtr}}$}}}   
\newcommand{\counite}{{\mbox{$\eta$}}} 
\newcommand{\coaugmentation}{{\mbox{$\cal \epsilon$}}} 

\newcommand{\aicat}{\mathsf{cat_\infty}}        
\newcommand{\ainat}{\mathsf{nat_\infty}}        

 \newcommand{\TW}{{\mathsf{tw}}}            

\newcommand{\Ba}{{B^+\!}}
\newcommand{\Oma}{{\Omega^+}\!}

 \newcommand{\Tria}{\opname{Tria}\nolimits}
 \newcommand{\tria}{\opname{tria}\nolimits}
 \newcommand{\rp}{{\mathbf{p}}} 

 \newcommand{\weq}{\mbox{\mbox{\rm $\ce  q$}}}
 \newcommand{\fib}{\mbox{\mbox{\rm $\cf ib$}}}
 \newcommand{\cof}{{\mbox{\mbox{\rm $\cc of$}}}}
 \newcommand{\Ho}{\opname{Ho}\nolimits}

 \newcommand{\Ass}{\mbox{${\mathsf{Ass}}$}}
 \newcommand{\gr}{{\cg {r}}} 

 \newcommand{\mor}[1]{\mathsf{mor}(#1)}
 \newcommand{\cod}[1]{\mathsf{cod}(#1)}
 \newcommand{\del}[1]{\delta (#1)}
 \newcommand{\der}[1]{\mathsf{der}(#1)}
 \newcommand{\ogr}{\sf{Gr}}

\newcommand{\tso}{\odot}

\newcommand{\tsi}{\stackrel{\infty}{\ts}}
\newcommand{\tpo}[1]{^{\odot #1}}
\newcommand{\ts}{\otimes}
\newcommand{\tp}[1]{{}^{\ts #1}}
\newcommand{\tbox}{\,\Box\,}
\newcommand{\prim}[1]{_{[ #1]}} 
\newcommand{\tw}{\ts_\tau}   

\newcommand{\ai}{\mbox{${\rm A}_{\infty}$}}
\renewcommand{\b}[1]{\overline{#1}}
\renewcommand{\hat}[1]{\widehat{#1}}
\newcommand{\ad}{{\mathsf{ad}}}
\newcommand{\mf}[1]{{\mathfrak{#1}}}
\newcommand{\?}{{\rule[.4mm]{2.4mm}{0.1mm}}\,}
\newcommand{\E}{{U}\!}
\newcommand{\lrus}[3]{{}_{#1}{#2}_{#3}} 
\newcommand{\+}{{}^{\sz +}\!}

\newcommand{\pw}{{}^{\sz \wedge \!}}
\newcommand{\sz}{\scriptsize}
\newcommand{\si}{{\omega}}
\renewcommand{\epsilon}{\varepsilon}
\newcommand{\op}{{}^{op}}
\newdir{ >}{{}*!/-8pt/\dir{>}} 

%% file: Remerciements.tex
Les annÈes de thËse s'accumulent et les soutenances se succËdent avec leur
suite infinie de remerciements. Finalement, tous les mots finissent usÈs
par le protocole, vecteur de distance retenue ou de pudeurs obligÈes.
Alors comment interprÈter l'ordre des remerciements? Il y a toujours le directeur
en premier, puis les rapporteurs, suivis du jury et des autres. La
hiÈrarchie se dÈcline, impeccable, sans pulsions, trop sage. Ici aussi
je vais dÈrouler cette carte classique des remerciements : je n'ai pas su inventer une
forme qui m'aurait semblait plus adÈquate. Mais je tenais quand mÍme ‡
Ècrire, ‡ insister et rÈpÈter que s'il ne fallait qu'une ligne, Áa serait celle-ci :\\

{\em
Merci Monsieur Keller. De tout mon coeur, merci. Vous seul devinez peut-Ítre ce que je vous dois.}\\

Je voudrais que cette ligne soit lue distinctement, rÈpÈtÈe avec conviction,
dans une respiration de gratitude. Je voudrais qu'elle s'imprime et
qu'elle se souvienne de la distance mesurÈe qui fonde le respect
que nous avons l'un pour l'autre. AllÈluia ! Le hasard qui nous
a menÈs l'un ‡ l'autre fut une belle route.\\

Bien. Si ces remerciements ne tiennent pas en une ligne, c'est parce que
chacun, et j'insiste sur cette Èvidence, mÈrite la ligne qui est la sienne.\\

L'intervention des rapporteurs, J.~D.~Stasheff et J.~Huebschmann, fut tout sauf mineure.
Par leurs remarques, ils ont grandement aidÈ ‡ la lisibilitÈ du texte.
Je me permets de leur exprimer mes remerciements les plus objectivement sincËres.\\

Je me plais ‡ recevoir dans mon jury P.~Cartier. R.~Rouquier (qui est aussi
dans le jury) fait partie des personnes que je croisais avec plaisir dans le couloir des mathÈmatiques.
J'ai lu avec intÈrÍt la thËse de doctorat d'A.~ProutÈ (qui est dans le jury).
Elle traitait d'$\ai$-structures. J'espËre que la suite que je prÈsente ici lui fera plaisir. \\

Le beau sourire de M.~Wasse, le professionnalisme ‡ visage humain de l'Èquipe informatique,
la fÈminitÈ spontanÈe de M.~Douchez ont su adoucir les tracas et obligations diverses qui savent
immanquablement ponctuer les existences mathÈmatiques.\\

Et puis, une fin ‡ ces lignes, une impossible fin : il y a dans mon coeur 
un point prÈcis qu'occupent ceux que j'aime et ceux qui m'aiment (des personnes parfois disjointes).
Ce lieu ne saurait Ítre dÈcrit ici... Trois petits points de suspension, pirouettes, cacahouËtes
et puis s'en vont !

%% file: Resume.tex
{\noindent \bf Abstract\\}

\noindent We study (not necessarily connected) $\Z$-graded $\ai$-algebras and their $\ai$-modules.
Using the cobar and the bar construction and Quillen's homotopical algebra, we describe the localisation
of the category of $\ai$-algebras with respect to $\ai$-quasi-isomorphisms. We then adapt these methods
to describe the derived category $\cd_\infty A$ of an augmented $\ai$-algebra $A$.
The case where $A$ is not endowed with an augmentation is treated differently.
Nevertheless, when $A$ is strictly unital, its derived category can be described in the same way as
in the augmented case. Next, we compare two different notions of $\ai$-unitarity :
strict unitarity and homological unitarity.
We show that, up to homotopy, there is no difference between these two notions.
We then establish a formalism which allows us to view $\ai$-categories as $\ai$-algebras in suitable
monoidal categories. We generalize the fundamental constructions of category theory to this setting : 
Yoneda embeddings, categories of
functors, equivalences of categories... We show that any algebraic
triangulated category $\calt$ which admits a set of
generators is $\ai$-pretriangulated,
that is to say, $\calt$ is equivalent to $H^0 \TW \ca$, where $\TW \ca$ is the $\ai$-category of twisted objets of a certain
$\ai$-category $\ca$. 

Thus we give detailed proofs of a part of the results on homological algebra
which M.~Kontsevich stated in his course
``Triangulated categories and geometry'' \cite{Kontsevich98}.\\

{\noindent \bf RÈsumÈ\\}

\noindent 
Nous Ètudions les $\ai$-algËbres $\Z$-graduÈes (non nÈcessairement connexes)
et leurs $\ai$-modules. En utilisant les constructions bar et cobar ainsi que les outils de l'algËbre
homotopique de Quillen, nous dÈcrivons la localisation de la catÈgorie des
$\ai$-algËbres par rapport aux $\ai$-quasi-isomorphismes.
Nous adaptons ensuite ces mÈthodes pour dÈcrire la catÈgorie dÈrivÈe
$\cd_\infty A$ d'une $\ai$-algËbre augmentÈe $A$.
Le cas o˘ $A$ n'est pas muni d'une augmentation est traitÈ diffÈremment.
NÈanmoins, lorsque $A$ est strictement unitaire, sa catÈgorie dÈrivÈe peut Ítre dÈcrite
de la mÍme maniËre que dans le cas augmentÈ.
Nous Ètudions ensuite deux variantes de la notion
d'unitaritÈ pour les $\ai$-algËbres : l'unitaritÈ stricte et l'unitaritÈ homologique. 
Nous montrons que d'un point de vue homotopique, il n'y a pas de diffÈrence
entre ces deux notions. Nous donnons ensuite un formalisme qui permet de
dÈfinir les $\ai$-catÈgories comme des $\ai$-algËbres dans certaines catÈgories monoÔdales.
Nous gÈnÈralisons ‡ ce cadre les constructions fondamentales de la thÈorie des catÈgories :
le foncteur de Yoneda,
les catÈgories de foncteurs, les Èquivalences de catÈgories...
Nous montrons que toute catÈgorie triangulÈe algÈbrique engendrÈe par un ensemble d'objets
est $\ai$-prÈtriangulÈe, c'est-‡-dire qu'elle est Èquivalente ‡ $H^0\TW \ca$, o˘ $\TW \ca$ est
l'$\ai$-catÈgorie des objets tordus d'une certaine $\ai$-catÈgorie $\ca$.

Nous dÈmontrons ainsi une partie des ÈnoncÈs d'algËbre homologique presentÈs par
M.~Kontsevich pendant son cours ``CatÈgories triangulÈes et gÈomÈtrie''  \cite{Kontsevich98}.\\[5cm]

\begin{center}
\input{k_5.pstex_t}\\[.5cm]
l'associaËdre $K_5$
\end{center}

%% file: k_5.pstex_t
\begin{picture}(0,0)%
\includegraphics{k_5.pstex}%
\end{picture}%
\setlength{\unitlength}{3947sp}%
\begingroup\makeatletter\ifx\SetFigFont\undefined%
\gdef\SetFigFont#1#2#3#4#5{%
  \reset@font\fontsize{#1}{#2pt}%
  \fontfamily{#3}\fontseries{#4}\fontshape{#5}%
  \selectfont}%
\fi\endgroup%
\begin{picture}(1824,1524)(289,-973)
\end{picture}

%% file: Introduction.tex
\noindent 
Nous renvoyons ‡ \cite{Keller01a} et \cite{Keller01b} pour une introduction
aux $\ai$-algËbres et leurs modules. Cette thËse contient, entre autres, les
dÈmonstrations dÈtaillÈes des ÈnoncÈs de \cite{Keller01a}. Hormis \cite{Kontsevich98}
et \cite{Keller01a}, nous nous sommes appuyÈs principalement sur l'article de V.~Hinich
\cite{Hinich01} et sur les travaux suivants : \cite{Stasheff63a}, \cite{Proute85},
\cite{Getzler90}, \cite{Huebschmann91a}, \cite{Gugenheim91},  \cite{Markl96}, \cite{Hinich97c}. 
Certains des rÈsultats de cette thËse ont ÈtÈ obtenus rÈcemment et de maniËre
indÈpendante par K.~Fukaya \cite{Fukaya01b}, P.~Seidel \cite{Seidel02}, A.~Lazarev \cite{Lazarev02},
V.~Lyubashenko \cite{Lyubashenko02} et M.~Kontsevich et Y.~Soibelmann \cite{Kontsevich02a},
\cite{Kontsevich02}. \\

{\bf \noindent Structures strictes et structures ‡ homotopie prËs}\\

\noindent Les structures de l'algËbre classique que l'on appelle strictes,
par exemple les algËbres associatives, commutatives ou les algËbres de Lie, se sont
rÈvÈlÈes insuffisantes en topologie car elles ne sont pas compatibles avec l'homotopie.
Ainsi, si $X$ est un espace de lacets et $Y$ est un espace topologique homotope ‡ $X$,
il n'est pas toujours possible de transfÈrer la structure de H-espace (qui est stricte)
de $X$ vers $Y$. C'est pour pallier ce dÈfaut que J.~Stasheff \cite{Stasheff63a}
a introduit la notion de structure $\ai$, qui est un assouplissement
de celle de semigroupe topologique. Les structures $\ai$ font partie des
{\em structures ‡ homotopie prËs}, c'est ‡ dire des structures dont
``le dÈfaut de strictitude'' est contrÙlÈ de maniËre
cohÈrente par des {\em homotopies d'ordre supÈrieur}. Pour certaines structures ‡ homotopie prËs,
des homotopies d'ordre supÈrieur
sont connues de longue date comme les opÈrations de Steenrod \cite{Steenrod47}, \cite{Steenrod52}
ou les produits de
Massey. Les structures ‡ homotopie prËs se comportent bien relativement aux Èquivalences d'homotopie :
si un objet (topologique ou
diffÈrentiel graduÈ) est muni d'une structure ‡ homotopie prËs, on peut
sous certaines conditions la translater sur un autre objet quand 
ce dernier est homotope ‡ l'objet de dÈpart. La premiËre partie de cette thËse
traitera des $\ai$-structures algÈbriques, c'est-‡-dire des $\ai$-structures
dans le cadre de l'algËbre diffÈrentielle graduÈe. Dans la deuxiËme partie nous
Ètudierons leurs gÈnÈralisations au cadre catÈgorique. \\

{\bf \noindent $\ai$-structures algÈbriques}\\

\noindent Soit $\corps$ un corps.
Une {\em $\ai$-algËbre} \cite{Stasheff63b} est un $\corps$-espace vectoriel $\Z$-graduÈ
$A$ muni de morphismes graduÈs  
\[
m_i : A \tp i \ra A, \quad i \geq 1,
\]
de degrÈ $2-i$, vÈrifiant des Èquations dont la premiËre dit que $m_1$ est une
diffÈrentielle, la deuxiËme que $m_1$ est une dÈrivation pour la {\em multiplication
} $m_2$ et la troisiËme
\[
m_2(m_2 \ts \Id) - m_2 (\Id \ts m_2) = \del{m_3} 
\]
que le dÈfaut d'associativitÈ de $m_2$ est mesurÈ par
le bord de $m_3$ dans l'espace diffÈrentiel graduÈ $\Hom (A\tp 3,A)$. Intuitivement, 
une $\ai$-algËbre est donc une ``algËbre diffÈrentielle graduÈe dont
le dÈfaut d'associativitÈ est contrÙlÈ (au sens fort)
par des homotopies d'ordre supÈrieur''. Si $A$ et $A'$ sont deux $\ai$-algËbres, un {\em $\ai$-morphisme}
$f : A \ra A'$
est une suite de morphismes graduÈs
\[
f_i : A \tp i \ra A' , \quad i \geq 1,
\]
de degrÈ $1-i$, vÈrifiant des Èquations dont les premiËres affirment que $f_1$ est un morphisme de complexes
qui est compatible aux multiplications $m_2$ et $m_2'$ ‡ une homotopie $f_2$ prËs.
De la mÍme maniËre, si $f$ et $g$ sont des $\ai$-morphismes $A \ra A'$,
une {\em homotopie} $h$ entre $f$ et $g$ est une
suite de morphismes
\[
h_i : A \tp i \ra A', \quad i \geq 1,
\]
de degrÈ $-i$, qui vÈrifient des Èquations dont les deux
premiËres affirment que $h_1$ est une homotopie entre les ``morphismes d'algËbres diffÈrentielles
graduÈes''
\[
f_1 \quad \mbox{et} \quad g_1 : (A,m_1,m_2) \ra  (A',m_1',m_2').
\]
Soit $A$ une $\ai$-algËbre.
Un {\em $A$-polydule} (appelÈ {\em $\ai$-module sur $A$} dans la littÈrature)
est un $\corps$-espace vectoriel
$\Z$-graduÈ $M$ muni de morphismes graduÈs
\[
m_i^M : M \ts A \tp{i-1} \ra M , \quad i \geq 1,
\]
de degrÈ $2-i$, vÈrifiant des Èquations dont les premiËres affirment que $m^M_1$ est une diffÈrentielle
et que $m^M_2$ dÈfinit une action de l'algËbre (fortement homotopiquement) associative $A$
dont la compatibilitÈ ‡ la multiplication
de $A$ est contrÙlÈe par l'homotopie d'ordre supÈrieur $m_3^M$.
Comme pour les $\ai$-algËbres, on a des $\ai$-morphismes entres $A$-polydules et des homotopies
entre les $\ai$-morphismes.\\

{\bf \noindent Lien avec la thÈorie des opÈrades}\\

\noindent Certains arguments de la thËse sont liÈs ‡ la thÈorie des opÈrades
(par exemple, la thÈorie de l'obstruction des $\ai$-algËbres (\ref{section_obstruction})).
Nous n'utiliserons pas explicitement le formalisme des opÈrades dans nos ÈnoncÈs
(et leurs dÈmonstrations), lui prÈfÈrant une approche naÔve. Rappelons nÈanmoins quelques
faits et rÈfÈrences sur ce sujet.

Les complexes cellulaires de Stasheff  $\{K_i \times \Sigma_i\}_{i\geq 2}$ (voir \cite{Stasheff63a})
forment une opÈrade topologique \cite{May72}.
Les complexes de chaÓnes qui leur sont associÈs forment donc une opÈrade
diffÈrentielle graduÈe. C'est l'opÈrade des $\ai$-algËbres.
Les opÈrades diffÈrentielles graduÈes ont ÈtÈ ÈtudiÈe abondamment au dÈbut des annÈes 90
\cite{Hinich93}, \cite{Getzler94}, \cite{Ginzburg94} pour expliciter 
le lien entre les structures strictes et les structures
‡ homotopie prËs \cite{Ginzburg94}, \cite{Markl96}, \cite{Markl99}, \cite{Markl00}.
En ce qui concerne les structures $\ai$, on retiendra des opÈrades deux rÈsultats :
{\em l'opÈrade des $\ai$-algËbres est
le modËle minimal cofibrant au sens de M.~Markl \cite{Markl96}
de l'opÈrade des algËbres associatives $\Ass$} ; {\em le dual de Koszul $\Ass^!$ de $\Ass$ est la co-opÈrade
des cogËbres co-associatives.}\\

{\bf \noindent Chapitre \ref{chapitre_Homot_aialg} : une thÈorie de l'homotopie des $\ai$-algËbres.}\\

\noindent Rappelons pour commencer un rÈsultat de H.~J.~Munkholm.
Soit $\mathsf{DA}$ la catÈgorie 
des algËbres diffÈrentielles graduÈes
(vÈrifiant certaines conditions sur la graduation et sur la connexitÈ)
et $\Ho \mathsf{DA}$ la localisation de $\mathsf{DA}$ par rapport aux quasi-isomorphismes.
Soit $\mathsf{DASH}$ la catÈgorie
des algËbres diffÈrentielles graduÈes dont les morphismes sont les $\ai$-morphismes.
Utilisant les idÈes de J.~Stasheff et S.~Halperin \cite{Stasheff70},
H.~J.~Munkholm \cite{Munkholm78} (voir aussi \cite{Munkholm76})
a montrÈ, premiËrement, que {\em la relation d'homo\-topie sur $\Hom_{\mathsf{DA}}(A,A')$, $A,A' \in DA$,
(qui n'est pas une relation d'Èquivalence en gÈnÈral) s'Ètend en une relation
sur les espaces de morphismes $\Hom_{\mathsf{DASH}}(A,A')$ qui est une relation d'Èquivalence
quelles que soient $A$ et $A'$,} et deuxiËmement, que
{\em la catÈgorie $\Ho \mathsf{DA}$ est Èquivalente au quotient de
$\mathsf{DASH}$ par cette relation d'Èquivalence}.
En d'autres termes, quitte ‡ augmenter le nombre de morphismes entre algËbres
diffÈrentielles graduÈes, la localisation par rapport
aux quasi-isomorphismes est Èquivalente au passage au quotient par rapport ‡ l'homotopie.
Dans la premiËre partie de ce chapitre, nous gÈnÈraliserons les rÈsultats de \cite{Munkholm78} 
aux $\ai$-algËbres. Un {\em $\ai$-quasi-isomorphisme}
$f$ est un $\ai$-morphisme tel que $f_1$ est un quasi-isomorphisme.
Nous montrons les rÈsultats suivants :\\

 \noindent (ThÈorËme de l'homotopie) {\it La relation d'homotopie sur 
les $\ai$-morphismes est une relation d'Èquiva\-lence
(\ref{corollaire_cmf_aia}~{\it a})}. \\

\noindent (ThÈorËme des $\ai$-quasi-isomorphismes) {\it 
Tout $\ai$-quasi-isomorphisme d'$\ai$-algËbres est inversible ‡ homotopie prËs
(\ref{corollaire_cmf_aia}~{\it b})}.\\

\noindent 
L'analogue topologique du thÈorËme des $\ai$-quasi-isomorphismes
est d˚ ‡ M.~Fuchs \cite{Fuchs76} (voir aussi \cite{Fuchs65}).
Dans sa thËse \cite{Proute85}, A.~ProutÈ a montrÈ les deux thÈorËmes
sous des conditions sur la graduation ou la connexitÈ (voir aussi \cite{Kadeishvili87}).
La nÈcÈssitÈ de gÈnÈraliser ces rÈsultats est due fait que
dans les constructions de K.~Fukaya et al.~de
$\ai$-algËbres ($\ai$-catÈgories), des composantes non nulles peuvent apparaÓtre en
tout degrÈ entier.
Dans le cas gÈnÈral, nous dÈduirons les thÈorËmes ci-dessus des rÈsultats suivants~: 
{\em la construction bar $B$
est une Èquivalence de catÈgories entre $\aia$, la catÈgorie des $\ai$-algËbres, et
la sous-catÈgorie des objets cofibrants et fibrants d'une catÈgorie de modËles $\cocog$
de cogËbres (\ref{theoreme_cmf_cocog}). La construction bar fait correspondre l'homotopie 
des $\ai$-morphismes ‡
l'homotopie ‡ gauche de $\cocog$ entre morphismes entre objets cofibrants et fibrants
(\ref{proposition1_homotopie}) et les $\ai$-quasi-isomorphismes aux Èquivalences faibles
(\ref{proposition1_cmf_aia}).}

La catÈgorie $\cocog$ en question est la catÈgorie des cogËbres diffÈrentielles
graduÈes cocomplËtes. Soit  $\alg$ la catÈgorie des algËbres diffÈrentielles graduÈes et
$\Omega : \cocog \ra \alg$ la construction cobar.
La structure de catÈgorie de modËles de $\cocog$ (\ref{theoreme_cmf_cocog}.~{\it a})
est telle que le couple de foncteurs adjoints
\[
(\Omega, B) : \cocog \ra \alg,
\]
est une Èquivalence de Quillen (\ref{theoreme_cmf_cocog}.~{\it b}).
L'utilisation de ce couple de foncteurs adjoints pour Ètudier la catÈgorie
$\alg [Qis^{-1}]$ remonte aux annÈes 70 avec les travaux de 
D.~Husemoller, J.~C.~Moore et J.~Stasheff \cite{Husemoller74} (voir aussi \cite{Eilenberg66}).
Ils considËrent les algËbres augmentÈes diffÈrentielles, graduÈes positivement d'une part et
les cogËbres co-augmentÈes diffÈrentielles, graduÈes positivement et connexes d'autre part, et montrent que
la localisation de la catÈgorie des algËbres par rapport aux quasi-isomorphismes est
Èquivalente ‡ la localisation de la catÈgorie des cogËbres par rapport aux quasi-isomorphismes.
Sans les hypothËses sur la graduation ou la connexitÈ, leur ÈnoncÈ n'est plus vrai.
Dans le cas gÈnÈral (\ref{theoreme_cmf_cocog}),
nous devons remplacer la classe des quasi-isomorphismes de $\cocog$ par une
classe de morphismes (appelÈs Èquivalences faibles) qui est
strictement contenue dans celle des quasi-isomorphismes
(\ref{proposition_equiv_qis}.~{\it c}). Nous montrons qu'entre deux cogËbres
graduÈes positivement, les Èquivalences faibles sont exactement
les quasi-isomorphismes (voir \ref{proposition_equiv_qis}.~{\it e}). Nos rÈsultats gÈnÈralisent
donc \cite[Chap.~II, Thm.~4.4 et Thm.~4.5]{Husemoller74}.

Notre dÈmonstration du fait que $\cocog$ admet une structure de catÈgorie de modËles (\ref{theoreme_cmf_cocog})
suit les idÈes de V.~Hinich \cite{Hinich01} inspirÈes de celles de Quillen
\cite{Quillen67}, \cite{Quillen69}. 
Nous relevons la structure de catÈgorie de modËles
de $\alg$ le long de l'adjonction $(\Omega , B)$. 
Cette adjonction est du mÍme type que l'adjonction entre la catÈgorie des algËbres de Lie
diffÈrentielles graduÈes et la catÈgorie des cogËbres cocommutatives diffÈrentielles
graduÈes en homotopie rationnelle.
Elle provient de la dualitÈ de Koszul entre l'opÈrade $\Ass$ et la co-opÈrade des cogËbres
co-associatives.

La caractÈrisation des objets fibrants de $\cocog$
peut Ítre interprÈtÈe comme une consÈquence du fait que l'opÈrade des $\ai$-algËbres est
le modËle minimal cofibrant au sens de M.~Markl \cite{Markl96}
de l'opÈrade des algËbres associatives. Ce fait implique
que l'obstruction ‡ la construction par rÈcurrence des
morphismes graduÈs $m_i $, $i\geq 1,$ dÈfinissant une $\ai$-structure sur un objet graduÈ $A$
est de la forme ``$m_{n+1}$ doit tuer un certain cocycle (construit ‡ partir des $m_i,$ $1 \leq i \leq n$)''
(voir  \ref{extension_nstructure_algebres}).
La condition qui mesure l'obstruction ‡ la construction
par rÈcurrence des $\ai$-morphismes est du mÍme type (\ref{extension_nstructure_morphismes}).
Nous appelons l'Ètude de ces obstructions la {\em thÈorie de l'obstruction} des $\ai$-algËbres.
Cette thÈorie est l'objet de l'appendice (\ref{section_obstruction}).

A la fin du chapitre \ref{chapitre_Homot_aialg}, nous redÈmontrerons (\ref{theoreme_transfert_structures})
la ``compatibilitÈ des structures $\ai$ ‡ l'homotopie'' : {\em soit $A$ une $\ai$-algËbre et
\[
g : (V,d) \ra (A,m^A_1)
\]
une Èquivalence d'homotopie de complexes.
Il existe une structure d'$\ai$-algËbre sur $V$ telle que $m_1^V$ est Ègale ‡ $d$ et telle que
$V$ et $A$ sont homotopes en tant qu'$\ai$-algËbres. }
Ce rÈsultat est bien connu. T.~Kadeishvili \cite{Kadeishvili80} et A.~ProutÈ \cite{Proute85}
l'ont montrÈ dans le cas o˘ $d = 0$ et sous des hypothËses sur la graduation et la connexitÈ
en utilisant la mÈthode des obstructions. Le cas gÈnÈral est d˚
‡ V.~K.~A.~M.~Gugenheim, L.~A.~Lambe et J.~Stasheff \cite{Gugenheim91} qui
utilisent ``l'astuce du tenseur'' inventÈe par J.~Huebschmann \cite{Huebschmann86}.
Le point essentiel de leur dÈmonstration est
que le lemme de perturbation \cite{Gugenheim72}
est compatible ‡ une structure additionnelle (de cogËbre dans notre cas).
Sur ce sujet, voir aussi  \cite{Huebschmann91a}, \cite{Gugenheim89} et les rappels historiques de la section
\ref{section_transfert_structures}.
Notre dÈmonstration de la ``compatibilitÈ ‡ l'homotopie''
(section \ref{theoreme_transfert_structures}) sera basÈe sur la
thÈorie des obstructions (\ref{section_obstruction}).
La ``compatibilitÈ ‡ l'homotopie'' implique que toute $\ai$-algËbre $A$ admet un {\em modËle minimal},
i.~e.~une structure $\ai$ sur l'homologie $H^*A$ telle que $H^*A$ et $A$
sont homotopes en tant qu'$\ai$-algËbres (\ref{corollaire_modele_minimal}).
Le lien entre un certain modËle minimal
obtenu par notre mÈthode et celui obtenu par le lemme de perturbation \cite{Gugenheim91}
est dÈcrit en (\ref{lemme_modele_minimal_I}).

La ``minimalitÈ'' du modËle $H^*A$ ci-dessus se rÈfËre au fait que la cogËbre tensorielle
$B(H^*A)$ est un modËle minimal (au sens de H.~J.~Baues et J.-M.~Lemaire \cite{Baues77}) de
la cogËbre $BA$.\\

{\bf \noindent Chapitre \ref{chapitre_Homot_aimod} :
une thÈorie de l'homotopie des polydules.}\\

\noindent 
Soit $A$ une $\ai$-algËbre {\em augmentÈe}. Rappelons que dans cette thËse les structures
communÈment appelÈes $\ai$-modules sur $A$ sont appelÈes {\em $A$-polydules}
(``poly'' car la structure est donnÈe par plusieurs
multiplications).

Le but de ce chapitre est de dÈcrire la catÈgorie dÈrivÈe
$\cd_\infty A$ dont les objets sont les $A$-polydules strictement
unitaires. On adapte pour cela les mÈthodes d'algËbre homotopique du chapitre 1 aux
$A$-polydules.
La catÈgorie dÈrivÈe d'une $\ai$-algËbre qui n'est pas munie d'une augmentation
sera ÈtudiÈe au chapitre \ref{chapitre_Cat_der}.

Soit $C$ une cogËbre diffÈrentielle graduÈe co-augmentÈe cocomplËte et
$\coModcu C$ la catÈgorie
des $C$-comodules diffÈrentiels graduÈs co-unitaires cocomplets.
Nous construisons (\ref{theoreme_cmf_coModcu}) une structure de catÈgorie de modËles sur $\coModcu C$
qui est telle que, si $A$ est une algËbre diffÈrentielle graduÈe augmentÈe et $ \tau : C \ra A$ une cochaÓne
tordante admissible acyclique, le couple de foncteurs adjoints
``produits tensoriels tordus'' (\ref{section_cochaines_tordantes})
\[
(\? \tw A , ? \tw C) : \coModcu C \ra \Modu A
\]
est une Èquivalence de Quillen. 
Les foncteurs ``produits tensoriels tordus'' remplacent ici les constructions bar et cobar du chapitre
prÈcÈdent. La catÈgorie homotopique $\Ho \coModcu C$ (voir appendice \ref{section_rappels_cmf})
est donc Èquivalente ‡ la catÈgorie dÈrivÈe
\[
\cd A = \Ho \Modu A.
\]
Dans \cite{Husemoller74}, D.~Husemoller, J.~C.~Moore et J.~Stasheff ont dÈmontrÈ un rÈsultat (le
thÈorËme 5.15) un peu plus gÈnÈral mais sous des hypothËses sur la graduation et la connexitÈ.
Nous ne considÈrerons pas ici les algËbres et les cogËbres Ètendues (voir \cite{Husemoller74}), nous
restreignant ‡ l'Ètude sÈparÈe des (co)algËbres et de leurs (co)modules.
Remarquons juste que notre rÈsultat (\ref{theoreme_cmf_coModcu})
gÈnÈralise la spÈcialisation du thÈorËme 5.15 de
\cite{Husemoller74} ‡ la sous-catÈgorie formÈe des algËbres Ètendues $(M,A,0)$, o˘
$A$ est une algËbre fixÈe et $M$ un $A$-module, et ‡ son image dans la catÈgorie des cogËbres Ètendues.

Nous Ètudions ensuite les objets fibrants de $\coModcu C$ pour une certaine classe de cogËbres $C.$
Soit $A$ une $\ai$-algËbre augmentÈe et 
$\aiModu A$ la catÈgorie des $A$-polydules strictement unitaires dont les morphismes
sont les $\ai$-morphismes strictement unitaires.
Notons $\Ba A$ la construction bar co-augmentÈe de la rÈduction $\b A$ de $A.$
Lorsque $C$ est une cogËbre isomorphe ‡ $\Ba A$, nous montrons (\ref{proposition2_categorie_derivee})
‡ l'aide de la thÈorie de l'obstruction (\ref{section_obstruction_polydules})
qu'{\em un objet de $\coModcu C$ est fibrant si et seulement si il est facteur direct
d'un objet presque colibre}. Comme tous les objets de $\coModcu C$
sont cofibrants, la sous-catÈgorie des objets cofibrants et fibrants est
l'image essentielle de la construction bar des $A$-polydules strictement unitaires.
Nous en dÈduisons (\ref{corollaire3_categorie_derivee}) 
que la catÈgorie dÈrivÈe
\[
\cd_\infty A = \aiModu A [\{Qis\}^{-1}]
\]
est  Èquivalente au quotient de la catÈgorie $\aiModu A$ par la relation d'homotopie (ceci montre
le thÈorËme des $\ai$-quasi-isomorphismes pour les $A$-polydules).
La structure triangulaire de $\cd_\infty A$ sera ÈtudiÈe dans la section
(\ref{section_structure_triangulee}).

Dans la section
\ref{section_categorie_derivee_augmentee_bipolydules}, nous Ètudions, par les
mÍmes mÈthodes, la catÈgorie dÈrivÈe des bipolydules (appelÈs {\em $\ai$-bimodules} dans
la littÈrature) strictement unitaires sur deux $\ai$-algËbres augmentÈes. 
Nous utiliserons les rÈsultats de cette section dans la seconde partie
de la thËse qui concerne les $\ai$-catÈgories.\\

{\bf \noindent Chapitre \ref{chapitre_Unite} : les unitÈs.}\\

\noindent 
Une $\corps$-algËbre associative $(A,\mu)$ est {\em unitaire} si
elle est munie d'un morphisme $\unite : \corps \ra A$ vÈrifiant les relations
\[
\mu (\unite \ts \Id ) = \Id \quad \mbox{et} \quad  \mu (\Id \ts \unite) = \Id.
\]
Il existe plusieurs relËvements de la notion d'unitaritÈ aux
$\ai$-algËbres. Nous en Ètudions deux : l'{\em unitaritÈ stricte} (dÈj‡ prÈsente dans
la version topologique de J.~Stasheff \cite{Stasheff63a}) et l'{\em unitaritÈ homologique}.
L'unitaritÈ stricte est la notion qui nous permettra de gÈnÈraliser
certaines propriÈtÈs classiques des algËbres unitaires aux $\ai$-algËbres.
L'unitaritÈ homologique, plus gÈnÈrale, apparaÓt dans les exemples gÈomÈtriques
\cite{Fukaya93}. Nous montrons que d'un point
de vue homotopique il n'y a pas de diffÈrence entre ces deux relËvements possibles de la notion d'unitaritÈ.
Plus prÈcisÈment, nous montrerons le rÈsultat suivant :
soit $\big(\aia\big)_{hu}$ 
la catÈgorie des $\ai$-algËbres homologiquement unitaires dont les morphismes
sont les $\ai$-morphismes homologiquement unitaires et $\big(\aia\big)_{su}$
la catÈgorie des $\ai$-algËbres strictement unitaires dont les morphismes
sont les $\ai$-morphismes strictement unitaires.
{\em Les catÈgories $\big(\aia\big)_{hu}$ et $\big(\aia\big)_{su}$
deviennent Èquivalentes aprËs passage ‡ l'homotopie} (\ref{corollaire_ai-algebres_strictement_unitaires}).
La dÈmonstration de ce rÈsultat sera basÈe sur une thÈorie de l'obstruction
des $\ai$-structures minimales (\ref{section_obstruction_unite}) et sur l'existence d'un
modËle minimal strictement unitaire pour les $\ai$-algËbres strictement unitaires
(\ref{proposition_modele_minimal_strictement_unitaire}). 

RÈcemment, K.~Fukaya \cite{Fukaya01}, \cite{Fukaya01a},
P.~Seidel \cite{Seidel02}, A.~Lazarev \cite{Lazarev02} et V.~Lyubashenko \cite{Lyubashenko02}
ont ÈtudiÈ le problËme des unitÈs de maniËre indÈpendante. Le relËvement de la notion d'unitaritÈ
de V.~Lyubashenko se spÈcialise ‡ notre notion d'unitaritÈ homologique si on travaille sur un corps
(V.~Lyubashenko travaille sur un anneau commutatif quelconque).\\

{\bf \noindent Chapitre \ref{chapitre_Cat_der} : la catÈgorie dÈrivÈe.}\\

\noindent Ici, nous dÈfinissons la catÈgorie dÈrivÈe d'une $\ai$-algËbre quelconque $A$ (non nÈcessairement
strictement unitaire). Nous montrerons que, lorsque $A$ est strictement unitaire,
sa catÈgorie dÈrivÈe admet les quatre descriptions suivantes
(\ref{theoreme_categorie_derivee_strictement_unitaire}) :
\english \begin{itemize}
\item[{\rm D1.}] 
la sous-catÈgorie triangulÈe $\Tria A$ de la catÈgorie dÈrivÈe $\cd_\infty (A\+)$
(o˘ $A\+$ est l'augmentation de $A$ et $\cd_\infty (A\+)$ est dÈfinie
dans le chapitre~\ref{chapitre_Homot_aimod}),
\item[{\rm D2.}] la catÈgorie
\[
\ch_\infty A = \aiModu A/\!\sim\,
\]
o˘ $\aiModu A$ est la catÈgorie des $A$-polydules strictement unitaires et $\sim$ est la relation d'homotopie,
\item[{\rm D3.}] la catÈgorie localisÈe
\[
\big(\aiModu A \big)[Qis^{-1}]
\]
o˘ $Qis$ est la classe des $\ai$-quasi-isomorphismes de $\aiModu A$,
\item[{\rm D4.}] la catÈgorie localisÈe
\[
\big(\aiModust A \big)[Qis^{-1}]
\]
o˘ $\aiModust A$ est la sous-catÈgorie non pleine de $\aiModu A$ dont les morphismes sont
les $\ai$-morphismes stricts.
\end{itemize} \francais
Nous montrerons (\ref{lemme_categorie_derivee_algebre_differentielle_graduee_unitaire})
que si $A$ est une algËbre diffÈrentielle graduÈe unitaire, la catÈgorie dÈrivÈe $\cd A$ 
(voir par exemple \cite{Keller94})
est Èquivalente aux catÈgories dÈfinies ci-dessus.\\

{\bf \noindent Chapitre \ref{chapitre_ai_cat} : prÈliminaires sur les $\ai$-catÈgories.}\\

\noindent La notion d'$\ai$-catÈgorie est une gÈnÈralisation naturelle de celle
d'$\ai$-algËbre. Au dÈbut des annÈes 90, les travaux de K.~Fukaya \cite{Fukaya93}
(voir aussi \cite{Fukaya01a}) ont montrÈ qu'elle apparaÓt naturellement dans l'Ètude
de l'homologie de Floer. InspirÈ par ces travaux, M.~Kontsevich, dans son exposÈ
\cite{Kontsevich94} au congrËs international, a donnÈ une interprÈtation conjecturale
de la symÈtrie miroir comme l'``ombre'' d'une Èquivalence
entre les catÈgories dÈrivÈes de deux $\ai$-catÈgories d'origine
gÈomÈtrique  (voir aussi \cite{Polishchuk98} o˘ cette conjecture
a ÈtÈ dÈmontrÈe pour les courbes elliptiques). Dans la suite de cette thËse,
nous gÈnÈralisons au cadre $\ai$-catÈgorique les constructions fondamentales de la
thÈorie des catÈgories : le foncteur de Yoneda,
les catÈgories de foncteurs, les Èquivalences de catÈgories, etc.~, et dÈmontrons certains
des rÈsultats ÈnoncÈs ou implicites dans \cite{Kontsevich98}. Nous utiliserons ou adapterons pour cela
certaines mÈthodes de la premiËre partie de la thËse.\\

Une $\ai$-catÈgorie est une $\ai$-algËbre avec plusieurs objets,
et rÈciproquement, une $\ai$-algËbre est une $\ai$-catÈgorie avec un objet.
Les problËmes soulevÈs
par l'augmentation du nombre d'objets sont nombreux et la gÈnÈralisation des
rÈsultats des chapitres prÈcÈdents est parfois trËs technique (par exemple pour
l'homotopie entre $\ai$-morphismes). Nous introduisons
une bicatÈgorie $\sf C$ dont les objets sont les ensembles. Comme $\sf C$ est
une bicatÈgorie, pour tout ensemble $\mathbb O$,
la catÈgorie des morphismes $\sf C(\mathbb O,\mathbb O)$ est une catÈgorie monoÔdale
(voir \cite[Chap.~XII, ß6]{MacLane98}).
Nous dÈfinissons (\ref{definition_ai-categorie})
une petite $\ai$-catÈgorie dont l'ensemble des objets est en bijection avec
un ensemble $\mathbb O$ comme une $\ai$-algËbre dans $\sf C(\mathbb O,\mathbb O).$
Nous dÈfinissons ensuite les $\ai$-foncteurs et les catÈgories diffÈrentielles
graduÈes $\cc_\infty \ca$ et $\cc_\infty (\ca,\cb)$ de $\ca$-polydules et $\ca$-$\cb$-bipolydules
strictement unitaires ($\ca$ et $\cb$ sont des $\ai$-catÈgories strictement unitaires).
Un lemme clef qui sera utile pour la
construction de l'$\ai$-foncteur de Yoneda (chapitre \ref{chapitre_objets_tordus})
est dÈmontrÈ en (\ref{lemme_clef}).\\

{\bf \noindent Chapitre \ref{chapitre_torsion} : la torsion d'$\ai$-catÈgories.}\\

\noindent Dans ce chapitre, nous gÈnÈralisons aux $\ai$-algËbres une technique de torsion bien connue en
thÈorie des dÈformations (pour un panorama, voir par exemple \cite{Huebschmann99}).
Soit $\ca$ une $\ai$-catÈgorie.
ConsidÈrons l'Èquation de Maurer-Cartan gÈnÈralisÈe
\[
\sum_{i=1}^\infty m_i\big(x \ts \hdots \ts x\big) = 0.
\]
Nous montrons (\ref{section_torsion_ai-categories_I}
et \ref{section_torsion_ai-categories_II}) qu'une solution
$x$ de cette Èquation (lorsqu'elle a un sens) donne 
une nouvelle $\ai$-catÈgorie $\ca_x$ appelÈe la torsion de $\ca$
par $x.$ La torsion des $\ai$-algËbres est due ‡ K.~Fukaya qui l'a introduite
(ainsi que celle des $\mathrm{L}_\infty$-algËbres)
dans \cite{Fukaya01a} et \cite{Fukaya01b} pour l'Ètude des $\ai$-dÈformations.
Nos formules pour les compositions tordues $m_i^x$, $i\geq 1$, de $\ca_x$
sont les mÍmes (‡ des signes Èquivalents prËs)
que dans \cite{Fukaya01a} mais la dÈmonstration du fait qu'elles dÈfinissent
bien une structure d'$\ai$-catÈgorie est diffÈrente.
Nous dÈcrivons ensuite la torsion des $\ai$-foncteurs
(\ref{section_torsion_ai-foncteurs_I} et \ref{section_torsion_ai-foncteurs_II})
et des (bi)polydules (\ref{section_torsion_bipolydules_I} et \ref{section_torsion_bipolydules_II})
par des solutions de l'Èquation de Maurer-Cartan.
Nous montrons aussi que si un $\ai$-foncteur $f$ induit un quasi-isomorphisme dans
les espaces de morphismes, sa torsion $f_x$ induit elle aussi un quasi-isomorphisme
dans les espaces de morphismes (\ref{proposition_torsion_ai-foncteurs_equiv_faible_I}). 

La torsion sera utile dans les chapitres \ref{chapitre_objets_tordus} et \ref{chapitre_ai_func}.\\

{\bf \noindent Chapitre \ref{chapitre_objets_tordus} : l'$\ai$-foncteur de Yoneda et les objets
tordus.}\\

\noindent Soit $\ca$ une catÈgorie. Rappelons que le foncteur de Yoneda est le foncteur
\[
\ca \ra \Modu \ca, \quad A \mapsto \Hom_\ca(\?,A).
\]
Dans ce chapitre, nous relevons ce foncteur en un $\ai$-foncteur (\ref{definition_ai-foncteur_de_Yoneda})
\[
y :\ca \ra \cc_\infty \ca, \quad A \mapsto \Hom_\ca(\?,A),
\]
o˘ $\ca$ est une $\ai$-catÈgorie. Si $\ca$ est strictement unitaire, nous montrons que
l'$\ai$-foncteur $y$ est strictement unitaire et qu'il se factorise par l'$\ai$-catÈgorie des objets tordus
$\TW \ca$ (\ref{theoreme_factorisation_Yoneda}). Les compositions de l'$\ai$-catÈgorie
$\TW \ca$ sont obtenues par torsion (chapitre \ref{chapitre_torsion}).

Si $\cg$ est une catÈgorie diffÈrentielle graduÈe unitaire, la catÈgorie (diffÈrentielle graduÈe)
des objets tordus est due ‡ A.~I.~Bondal et M.~M.~Kapranov \cite{Bondal91}
(ils la notent $\mbox{Pr-Tr}\+\,\cg$). Le but de \cite{Bondal91} est de pallier
un dÈfaut des axiomes des catÈgories triangulÈes pour dÈcrire les catÈgories dÈrivÈes \cite{Verdier77} :
le cÙne n'est pas fonctoriel.
PlutÙt que les catÈgories triangulÈes, ils considËrent les catÈgories prÈ-triangulÈes dÈcrites ‡ l'aide
de la catÈgorie des objets tordus et montrent l'Èquivalence de catÈgories suivante :
soit $\ce$ une catÈgorie prÈ-triangulÈe ($H^0\ce$ est alors triangulÈe).
Soit $\cg$ une sous-catÈgorie pleine de $\ce$. La sous-catÈgorie triangulÈe $\tria \cg \subset H^0\ce$
engendrÈe par $\cg$ est Èquivalente ‡ la catÈgorie triangulÈe $H^0(\mbox{Pr-Tr}^+\cg)$.
Dans le cas $\ai$, on a les mÍmes rÈsultats : nous montrons  (\ref{section_equivalence_tria_ca_H^0TWca})
que si $\ca$ est une $\ai$-catÈgorie strictement unitaire, les catÈgories
\[
H^0\TW \ca\quad \mbox{et} \quad \tria \ca \subset \cd_{\infty}\ca
\]
sont Èquivalentes (comme annoncÈ dans \cite{Kontsevich94}). De plus, nous montrons
(section \ref{section_categories_stables}) que toute catÈgorie
triangulÈe algÈbrique qui est engendrÈe par un ensemble d'objets
est $\ai$-prÈ-triangulÈe, i.~e.~elle est Èquivalente
‡ $H^0\TW \ca$, pour une certaine $\ai$-catÈgorie $\ca$.

Soit $\ca$ une $\ai$-catÈgorie strictement unitaire.
La catÈgorie $\cc_\infty \ca$ est diffÈrentielle graduÈe et l'$\ai$-foncteur de Yoneda
$y : \ca \ra \cc_\infty \ca$  induit
(\ref{lemme_Yoneda_quasi-isomorphisme}) un quasi-isomorphisme dans les espaces de morphismes.
Nous en dÈduisons que l'image $y(\ca) \subset \cc_\infty \ca$ est une catÈgorie diffÈrentielle
graduÈe unitaire qui est quasi-isomorphe ‡ $\ca$. Ceci montre que d'un point de vue homologique,
l'Ètude des $\ai$-catÈgories strictement unitaires (et mÍme homologiquement unitaires, par le
chapitre \ref{chapitre_Unite})  revient ‡ l'Ètude des catÈgories diffÈrentielles graduÈes unitaires. 
Concernant les  catÈgories diffÈrentielles
graduÈes et leurs catÈgories dÈrivÈes, on renvoie ‡ \cite{Keller94},~\cite{Keller99}.\\

{\bf \noindent Chapitre \ref{chapitre_ai_func} : l'$\ai$-catÈgorie des $\ai$-foncteurs.}\\

\noindent Soit $\ca$ et $\cb$ deux $\ai$-catÈgories strictement unitaires. Nous dÈfinissons
(\ref{section_Nunc} et \ref{section_Func}) une
$\ai$-catÈgorie $\Func (\ca,\cb)$ dont les objets sont les $\ai$-foncteurs
strictement unitaires $\ca \ra \cb$. La difficultÈ consiste ‡ dÈfinir les compositions
supÈrieures des morphismes entre $\ai$-foncteurs. Nous utiliserons pour cela
la mÈthode de la torsion du chapitre \ref{chapitre_torsion}.
Cette $\ai$-catÈgorie est fonctorielle
en $\ca$ et $\cb$ (\ref{section_fonctorialite_Nunc}).
Nous en dÈduisons une $2$-catÈgorie $\aicat$ dont les objets sont les
petites $\ai$-catÈgories strictement unitaires
et les espaces de morphismes $\ca \ra \cb$ sont les catÈgories 
\[
\aicat (\ca,\cb) = H^0\Func(\ca,\cb), \quad \ca, \cb \in \Obj \aicat.
\]
Nous caractÈrisons (\ref{proposition_equivalence_faible_d'ai-foncteurs})
ensuite les ÈlÈments
\[
H \in \Hom_{\Func(\ca,\cb)} (f,g), \quad f, g : \ca \ra \cb
\]
qui deviennent des isomorphismes
$f \ra g$ dans la catÈgorie $\aicat(\ca,\cb)$.
La dÈmonstration de cette caractÈrisation utilisera l'existence d'un {\em $\ai$-foncteur
de Yoneda gÈnÈralisÈ} (\ref{section_ai-foncteur_de_Yoneda_generalise})
\[
z : \Func (\ca,\cb) \ra \cc_\infty (\ca,\cb)
\]
qui induit un quasi-isomorphisme dans les espaces de morphismes.

L'$\ai$-catÈgorie $\Func(\ca,\cb)$ a ÈtÈ construite indÈpendamment par
K.~Fukaya \cite{Fukaya01a}, V.~Lyubashenko \cite{Lyubashenko02}  et
M.~Kontsevich et Y.~Soibelman \cite{Kontsevich02}, \cite{Kontsevich02a}.
Bien qu'obtenues par des mÈthodes diffÈrentes, les compositions de $\Func(\ca,\cb)$ de
\cite{Lyubashenko02} sont les mÍmes que les nÙtres.
\\

{\bf \noindent Chapitre \ref{chapitre_ai_equiv} : les $\ai$-Èquivalences.}\\

\noindent Soit $\ca$ une $\ai$-catÈgorie strictement unitaire. Dans (\ref{section_ai-isomorphie}),
nous relevons la notion d'isomorphisme de $H^0 \ca$ ‡ $\ca$. Nous montrons ensuite
qu'un $\ai$-foncteur $f : \ca \ra \cb$ est une
$\ai$-Èquivalence si et seulement si $f_1$ est un quasi-isomorphisme et 
s'il induit une Èquivalence de catÈgories (au sens classique) entre
$H^0\ca$ et $H^0 \cb$ (\ref{section_ai-equivalences}).
D'autres dÈmonstrations de cette caractÈrisation (annoncÈe dans \cite{Kontsevich98})
se trouvent dans \cite{Fukaya01a} et
\cite{Lyubashenko02}.

%% file: Homot_aialg.tex
{\noindent \bf Introduction}\\

\noindent Rappelons trois rÈsultats classiques
sur les $\ai$-algËbres :
\english \begin{itemize}
\item[{\it 1.}]{\it (Relation d'homotopie)}
La relation d'homotopie sur les $\ai$-morphis\-mes est une relation d'Èquivalence
(\ref{corollaire_cmf_aia}~{\it a}).
\item[{\it 2.}] {\it ($\ai$-quasi-isomorphisme)} 
Tout $\ai$-quasi-isomorphisme d'$\ai$-algËbres est inversible ‡ homotopie
prËs (\ref{corollaire_cmf_aia}~{\it b}).
\item[{\it 3.}] {\it (ModËle minimal)} Toute $\ai$-algËbre admet un modËle minimal
(\ref{corollaire_modele_minimal}).
\end{itemize}\francais
Dans la littÈrature, les rÈsultats {\it 1} et {\it 2} sont dÈmontrÈs pour les $\ai$-algËbres vÈrifiant
certaines conditions sur leur graduation ou leur connexitÈ
(voir les rÈfÈrences figurant dans le corps du chapitre).
Le but de ce chapitre est de les gÈnÈraliser aux $\ai$-algËbres quelconques.\\

{\noindent \bf Plan du chapitre}\\

\noindent Le chapitre est divisÈ en quatre sections.
Dans la section \ref{section_notations}, on fixe les notations
et on dÈfinit les algËbres libres et les cogËbres tensorielles.

Dans la section \ref{section_definition_ai-algebres}, on dÈfinit les $\ai$-algËbres, les
$\ai$-morphismes et les homotopies entre $\ai$-morphismes.
On rappelle les constructions bar et cobar (\ref{construction_bar_cobar}). 

Dans la section \ref{section_cmf_cocog}, nous montrons le rÈsultat principal (\ref{theoreme_cmf_cocog})
de ce chapitre :\\

{\em \noindent La catÈgorie $\cocog$ des cogËbres diffÈrentielles graduÈes cocomplËtes admet une structure
de catÈgorie de modËles qui la rend Quillen-Èquivalente ‡ la catÈgorie de modËles $\alg$ des
algËbres diffÈrentielles graduÈes. Tous les objets de $\cocog$ sont cofibrants et les objets
fibrants de $\cocog$ sont ceux qui, en tant que cogËbres graduÈes, sont isomorphes ‡ des cogËbres
tensorielles rÈduites}.\\

\noindent 
La dÈmonstration du fait que la catÈgorie $\cocog$ admet une telle structure
nous a ÈtÈ inspirÈe du travail de V.~Hinich \cite{Hinich01}.
Nous considÈrons des objets filtrÈs et Ètudions dans ce cadre les propriÈtÈs des constructions
bar et cobar.
La caractÈrisation des objets cofibrants sera immÈdiate car les cofibrations sont
les injections. La caractÈrisation des objets fibrants sera un rÈsultat plus
profond, consÈquence du thÈorËme (\ref{theoreme_cmf_aia}) :
{\em la catÈgorie des $\ai$-algËbres $\aia$ admet une structure de ``catÈgorie
de modËles sans limites'' dont la classe des Èquivalences faibles est formÈe
des $\ai$-quasi-isomorphismes.}

Notre dÈmonstration de ce rÈsultat sera entiËrement basÈe sur la thÈorie de l'obstruction
(voir appendice \ref{section_obstruction}). Elle peut donc Ítre interprÈtÈe comme une
consÈquence du fait que l'opÈrade des $\ai$-algËbres est un modËle
cofibrant minimal au sens de M.~Markl \cite{Markl96}
pour l'opÈrade des algËbres associatives.

Les $\ai$-algËbres s'identifient par la construction bar aux objets
fibrants et cofibrants de $\cocog$. Les rÈsultats {\it 1} et {\it 2} citÈs plus
haut apparaÓtront alors comme des cas particuliers de rÈsultats fondamentaux de l'algËbre
homotopique de Quillen (voir appendice \ref{section_rappels_cmf}).

Dans la section \ref{section_transfert_structures}, nous remontrons (\ref{corollaire_modele_minimal})
le rÈsultat {\it 3} (modËle minimal).
Notre dÈmonstration utilisera la thÈorie de l'obstruction.
Ensuite, nous comparons (\ref{lemme_modele_minimal_I}) un modËle minimal obtenu ainsi avec celui
obtenu gr‚ce au lemme de perturbation (voir par ex.~\cite{Huebschmann91a}).

\section{Rappels et notations} \label{section_notations}

\subsection{Objets diffÈrentiels graduÈs}

\noindent Nous fixons des notations que nous utiliserons tout au long de ce chapitre.\\

{\noindent \bf La catÈgorie de base}\\ \label{categorie_de_base}

\noindent Soit $\corps$ \indexnotation{corps} un corps. Soit $\sf C$ \indexnotation{C0}
une catÈgorie $\corps$-linÈaire abÈlienne, semi-simple,
cocomplËte, aux colimites filtrantes exactes (i.e.~une $\corps$-catÈgorie de Grothen\-dieck
semi-simple). Nous supposons en outre que $\sf C$ est munie d'une structure de
catÈgorie monoÔdale $\corps$-bilinÈaire donnÈe par un foncteur \indexnotation{ts}
\[
\ts : \sf C \times \sf C  \ra \sf C,
\]
un objet neutre \indexnotation{e}  $e$,  et des contraintes d'associativitÈ et d'unitaritÈ
\[
X\ts (Y \ts Z) \iso (X \ts Y) \ts Z,  \hspace{1cm} X \ts e \iso X  \iso e \ts X, \hspace{1cm} X,Y,Z \in \sf C.
\]
Nous supposons que pour tout objet $X$ de $\sf C$, les foncteurs $X \ts \?$ et $\? \ts X$
sont exacts et commutent aux colimites filtrantes.

La catÈgorie des $\corps$-espaces vectoriels vÈrifie bien s˚r ces hypothËses. La raison pour
laquelle nous travaillons dans un cadre plus gÈnÈral est l'apparition naturelle d'autres
exemples dans l'Ètude des $\ai$-catÈgories (voir le chapitre \ref{chapitre_ai_cat}).\\

{\noindent \bf Objets graduÈs}\\

\noindent Un {\em objet graduÈ} \index{objet graduÈ}
(sur $\sf C$) est une suite $M = (M^p)_{p\in \Z}$ d'objets de $\sf C$.
Soit deux objets graduÈs $M$ et $L$. 
La {\em catÈgorie}  $\gr \sf C$ \indexnotation{grC} des objets graduÈs a pour 
espace des morphismes de $M$ dans $L$
l'espace vectoriel $\Z$-graduÈ de composantes
\[
\Hom_{\gr \sf C} (M,L)^r = \prod_p \Hom_{\scriptsize \sf C}(M^p, L^{p+r}), \hspace{1cm} r \in \Z.
\]
On appelle {\em morphismes graduÈs de degrÈ $r$} les ÈlÈments de la $r$-iËme composante.
Le {\em produit tensoriel} de deux objets graduÈs $M$ et $L$ a pour
composantes
\[
(M\ts L)^n = \bigoplus_{p+q=n} M^p \ts L^q, \hspace{1cm} n \in \Z.
\]
Soit $f : M \ra M'$ et $g : L \ra L'$ deux morphismes graduÈs de degrÈ $r$ et $s$.
Le {\em produit tensoriel}
\[
f \ts g : M \ts L \ra M' \ts L'
\] 
est le morphisme de degrÈ $r+s$ dont la $n$-iËme composante est induite par les morphismes
\[
(-1)^{ps} f^p \ts g^q : M^p \ts L^q \ra M'^{p+r} \ts L'^{q+s}, \hspace{1cm} p+q = n.
\]
L'ÈlÈment neutre pour le produit tensoriel graduÈ est l'objet graduÈ dont toutes
les composantes sont nulles sauf la $0$-iËme, qui vaut $e$. Nous le notons aussi $e.$
La catÈgorie $\gr \sf C$ est ainsi munie d'une structure de catÈgorie monoÔdale.
On dÈfinit le {\em foncteur suspension} \index{suspension} \indexnotation{S}
$S : \gr \sf C \ra \gr \sf C$ par
\[
(SM)^i = M^{i+1}, \hspace{1cm} i \in \Z.
\]
Nous notons \indexnotation{s}
\[
s_M : M \ra SM
\]
le {\em morphisme graduÈ fonctoriel} de degrÈ $-1$ de composantes
\[
s_M^i = \Id_{M^i} : M^i \ra (SM)^{i-1}, \hspace{1cm} i\in \Z.
\]
Le morphisme $s^{-1}$ est notÈ $\si.$ \indexnotation{si} Remarquons l'ÈgalitÈ
\[
\si \tp i \circ s \tp i = (-1)^{\frac{i(i-1)}{2}} \Id_{M \tp i}.
\]

{\noindent \bf Objets diffÈrentiels graduÈs}\\ \label{complexes}

%
%
%

\noindent Un {\em objet diffÈrentiel graduÈ} (ou {\em complexe})\index{complexe}
 est un couple $(M,d)$, o˘ $M$ est un  objet
graduÈ et  $d$ est une {\em diffÈrentielle}, c'est-‡-dire un
endomorphisme de $M$ de degrÈ $+1$, tel que $d^2 = 0$.
Le {\em sous-objet}  $Z^iM = \ker d^i$ de $M^i$ 
est l'objet des {\em cycles} \label{cycle} de degrÈ $i$ du complexe $M$.
Le {\em sous-objet}  $B^iM = \im  d^{i-1}$ de $Z^iM$ 
est l'objet des {\em bords} \label{bord} de degrÈ $i$ du complexe $M$.
Si $(M,d_M)$ et
$(L,d_L)$ sont deux complexes, nous munissons l'espace des morphismes graduÈs
$\Hom_{\gr \sf C} (M,L)$ de la diffÈrentielle  $\delta$ \label{delta} 
\indexnotation{del} dont les
composantes sont
\[
\begin{array}{rccl}
\delta^r  : & \Hom_{\gr \sf C} (M,L)^r & \ra & \Hom_{\gr \sf C} (M,L)^{r+1}, \hspace{1cm} r \in \Z.\\
& f & \mapsto &  d_L \circ f - (-1)^r f \circ d_M
\end{array}
\]
La {\em catÈgorie} $\cc \sf C$ \indexnotation{ccC}
a pour objets les complexes  et pour espaces de morphismes
\[
\Hom_{\sz \cc \sf C} (M,L) = Z^0 (\Hom_{\gr \sf C} (M,L),\delta).
\]
Si $M$ et $L$ sont deux complexes, on munit le produit tensoriel graduÈ
$M \ts L$ de la diffÈrentielle 
\[
d_{M\ts L} = d_M \ts \Id_L + \Id_M \ts d_L.
\]
Nous avons ainsi muni $\cc \sf C$ d'une structure de catÈgorie monoÔdale
d'objet neutre l'objet graduÈ $e$ muni de la diffÈrentielle nulle.
Si $M$ est un complexe, nous munissons sa suspension $SM$ de la diffÈrentielle
\[
d_{SM}=-s_M \circ  d_M \circ s_M^{-1}.
\]

%
%
%
Le foncteur {\em homologie} $H : \cc \sf C \ra \gr \sf C$ envoie un complexe $M$
sur l'objet graduÈ $HM$ de composantes
\[
H^i M = Z^iM / B^iM, \hspace{1cm} i\in \Z.
\]
Un {\em quasi-isomorphisme} de $\cc \sf C$
est un morphisme qui induit un isomorphisme en homologie. 
Un complexe est {\em acyclique} s'il est quasi-isomorphe ‡ l'objet
nul.
%
%
%
Deux morphismes de complexes $f,g : M\ra L$ sont
{\em homotopes} \index{homotopie} s'il existe un
morphisme $r : M \ra L$ de degrÈ $-1$ tel que $\del r = f - g$.
L'homotopie est une relation d'Èquivalence.
La {\em catÈgorie}  $\ch \sf C$ a pour objets les complexes et pour
espaces de morphismes de $M$ dans $L$ les classes d'homotopie de morphismes de
la catÈgorie $\cc \sf C$ :
\[
\Hom_{\scriptsize \ch \sf C} (M,L) = H^0(\Hom_{\gr \sf C} (M,L),\delta).
\]
Nous notons encore $H : \ch \sf C \ra \gr \sf C$ le foncteur induit par
le foncteur homologie.
%
%
%
%
%
%
%
%

\subsection{AlgËbres et cogËbres}

{\noindent \bf  AlgËbres}\\

Soit $\sf M$ l'une des catÈgories  $\sf C$, $\gr \sf C$ ou $\cc \sf C$.
Une {\em algËbre}\index{algebre@{algËbre}} $(A,\mu)$ dans $\sf M$
est un objet $A$ muni d'une
multiplication $\mu : A \ts A \ra A$ associative (et de degrÈ $0$ si $\sf M = \gr \sf C$). 
DÈfinissons $\mu^{(2)} = \mu$, et pour tout $n\geq 3$,
$ \mu^{(n)} :  A\tp n \ra A$ par
\[
\mu^{(n)} = \mu (\Id \ts \mu^{(n-1)}).
\]
Pour $n\geq 1$, on appelle $\coker \mu^{(n+1)}$ l'{\em algËbre des $n$-irrÈductibles}
\index{irrÈductible} de $A$.

Soit $f,g: A\ra B$ deux morphismes d'algËbres.
Une {\em $(f,g)$-dÈrivation} \index{derivation@{dÈrivation}}
est un morphisme $D : A\ra B$  vÈrifiant la rËgle de Leibnitz
\[
D \circ \mu = \mu \circ (f\ts D + D \ts g).
\]
Une {\em dÈrivation} de l'algËbre $A$  est une $(\Id_A ,\Id_A)$-dÈrivation.

%
%
%
%
%
%
%
%
%
Soit $V$ un objet graduÈ de $\sf M.$ L'{\em algËbre tensorielle rÈduite sur $V$}
\index{tensorielle rÈduite (algËbre)} \index{algebre@{algËbre}!tensorielle rÈduite} \indexnotation{bTV}
est l'objet 
\[
\b TV = \bigoplus_{i\geq 1} V\tp i
\]
muni de la multiplication $\mu$ dont les composantes
\[
V \tp i \ts V\tp j \ra V\tp{i+j} \ra \b TV
\]
sont donnÈes par la contrainte d'associativitÈ de la catÈgorie monoÔdale
$\sf M.$ 
Une algËbre $A$ de $\sf M$ est {\em
libre}\index{algebre@{algËbre}!libre} \index{libre!algËbre}
si elle est isomorphe ‡ $\b TV$ pour un objet $V$ de $\sf M.$
Nous avons alors $V \iso \coker \mu_A.$

\begin{lemme}[propriÈtÈ universelle de l'algËbre tensorielle]
\label{algebres_libres} \hspace*{1cm}\\
Soit $(A,\mu)$ une algËbre. Pour $n\geq 1$, nous notons $j_n : V\tp n \ra \b T(V)$
l'injection canonique.
\begin{enumerate}
\item[a.] L'application $f \mapsto f \circ j_1$ est une bijection de l'ensemble
des morphismes d'algËbres $\b T(V) \ra A$ sur l'ensemble des morphismes $V \ra A$ de $\sf M$
(de degrÈ $0$ si $\sf M = \gr \sf C$). 
L'application inverse associe ‡ $g : V \ra A$ le morphisme
d'algËbres $\mor g : \b TV \ra A$ dont la $n$-iËme  composante est 
\[
 V\tp n \arr{g\tp n} A \tp n \arr{\mu ^{(n)}} A, \hspace{1cm} n \geq 1.
\]
\item[b.] Soit $f,g : A \ra B$ deux morphismes d'algËbres.
L'application $D \mapsto D \circ j_1$ est une bijection de l'ensemble des
$(f,g)$-dÈrivations sur l'ensemble des morphismes $V \ra A$ de $\sf M$. 
 L'application inverse associe ‡ $h : V \ra A$ la $(f,g)$-dÈrivation
$\der h : \b TV \ra A$ dont dont la $n$-iËme  composante est 
\[
\mu^{(n)}\circ  \left(  \sum_{l+1+j=n}(f^{\ts l}\ts h\ts g^{\ts
j})\right), \hspace{1cm} n \geq 1.
\]
\hspace*{1cm} \findem
\end{enumerate}
\end{lemme}
Une {\em algËbre graduÈe} \index{algebre@{algËbre}!graduÈe} (resp.~{\em diffÈrentielle graduÈe})
\index{algebre@{algËbre}!diffÈrentielle graduÈe} \label{definition_alg} est une algËbre de la catÈgorie $\gr \sf C$
(resp.~de la catÈgorie $\cc \sf C$). 
Nous notons $\alg$ \indexnotation{alg} la catÈgorie des
algËbres diffÈrentielles graduÈes. Un morphisme de $\alg$
est un {\em quasi-isomorphisme} s'il induit un isomorphisme
en homologie.
Une algËbre diffÈrentielle graduÈe est
{\em presque libre}\index{algebre@{algËbre}!presque libre} \index{presque libre!algebre@{algËbre}}
si elle est libre en tant qu'algËbre graduÈe.
Deux morphismes $f,g :A \ra B$ de $\alg$ sont
{\em homotopes}{\index{homotopie}} s'il
existe  une $(f,g)$-dÈrivation  $H :A \ra B$ graduÈe de degrÈ
$-1$ telle que
\[
f - g = dH + Hd.
\]
Il rÈsultera de la proposition \ref{proposition_cmf}
appliquÈe ‡ l'exemple \ref{cmf_algebres} que, si l'algËbre $A$ est presque libre,
la relation d'homotopie est
une relation d'Èquiva\-lence sur l'ensemble des morphismes d'algËbres de $A$ dans $B$. \\
%
%
%
%
%
%
%
%

{\noindent \bf CogËbres}\\

Une {\em cogËbre} \index{cogebre@{cogËbre}} dans  $\sf M$
est la  donnÈe d'un objet $C$  muni  d'une comultiplication
$\Delta:C\ra C\ts C$ co-associative, i.~e.~$(\Delta \ts \Id)\Delta = (\Id \ts
\Delta) \Delta $. DÈfinissons $\Delta^{(2)} = \Delta$ et, pour tout $n\geq 3$,
$\Delta ^{(n)} : C \ra C  \tp{n}$ par
\[
\Delta ^{(n)} = (\Id \tp{n-2} \ts \Delta)\circ \Delta^{(n-1)}.
\]
Soit $n \geq 1.$ Le noyau $C \prim{n} = \ker \Delta^{(n+1)}$ est une sous-cogËbre de $C$ ; nous
l'appelons la {\em sous-cogËbre des $n$-primitifs} \label{primitifs} \indexnotation{prim}
\index{primitifs@{$n$-primitifs}}
de $C.$ La suite croissante de sous-cogËbres
\[
C \prim{1} \subset C \prim{2} \subset C \prim{3} \subset \cdots
\]
est la {\em filtration primitive}
\index{filtration!primitive} de la cogËbre
$C.$ La cogËbre $C$ est  {\em cocomplËte}\index{cocomplËte (cogËbre)}
\index{cogebre@{cogËbre}!cocomplËte}
si l'on~a
\[
\colim C \prim{i} = C.
\]
Soit $f$ et $g:C\ra B$ deux morphismes  de cogËbres. Une {\em
$(f,g)$-codÈrivation}\index{coderivation@{codÈrivation}}
est un morphisme $ D:C\ra B$  vÈrifiant l'identitÈ duale de
la rËgle de Leibniz
\[
\Delta \circ D = (f\ts D+D \ts g)\circ \Delta .
\] 
Une {\em codÈrivation} de
$C$ est une $(\Id_C ,\Id_C)$-codÈrivation. \\

%
%
%
%
%
%
%
%
Soit $V$ un objet de $\sf M$. La {\em cogËbre tensorielle rÈduite sur $V$}
\index{tensorielle rÈduite (cogËbre)} \index{cogebre@{cogËbre}!tensorielle rÈduite} \indexnotation{ctr}
est l'objet
\[
\ctr V = \bigoplus_{i\geq
1} V \tp i
\]
muni de la comultiplication dont la $n$-iËme composante
\[
V\tp{n} \arr{} \oplus_{i+j=n}V \tp i \ts V\tp j \arr{} \ctr V \ts \ctr V,
\]
est la somme des morphismes $V\tp{n} \ra V \tp i \ts V\tp j$ donnÈs par
la contrainte d'associativitÈ de la structure monoÔdale de $\sf M.$
Remarquons que si $C$ est isomorphe ‡ une cogËbre
tensorielle rÈduite, elle est isomorphe ‡ $\ctr (C \prim{1})$. Les
cogËbres tensorielles  rÈduites sont cocomplËtes.  \\
\begin{lemme}[propriÈtÈ universelle de la cogËbre tensorielle]
\label{cogebres_tensorielles} \hspace*{1cm} \\
Soit $C$ une cogËbre cocomplËte. Pour $n \geq 1$, nous notons
$p_n:\ctr (V)\ra V^{\ts n}$ la
projection canonique.
\english
\begin{itemize}
\item[a.] L'application $f \mapsto p_1\circ f$ est une bijection de l'ensemble
des morphismes de cogËbres
sur l'ensemble des morphismes  $C \ra V$ de $\sf M$ (de degrÈ $0$ si $\sf M = \gr \sf C$).
L'application inverse associe ‡ $g:C \ra V$ le morphisme de cogËbres
$\mor g : C \ra \ctr V$  dont la $n$-iËme composante est
\[
C \arr{\Delta^{(n)}} C \tp n \arr{g \tp n} V\tp n, \hspace{1cm} n\geq 1.
\]

\item[b.] Soit $f,g : C \ra \ctr V$ deux morphismes de cogËbres.
L'application $D\mapsto p_1\circ D$ est une bijection de l'ensemble des
$(f,g)$-codÈrivations $C \ra \ctr V$ sur l'ensemble des morphismes $C \ra V$.
L'application inverse associe  ‡ $h:C\ra V$ la
$(f,g)$-codÈrivation $\cod{h} : C \ra \ctr V$ 
dont la $n$-iËme composante est
\[
\left( \sum_{l+1+j=n}(f^{\ts l}\ts h\ts
g^{\ts j})\right) \circ \Delta^{(n)}, \hspace{1cm} n\geq 1.
\]
\hspace*{1cm} \findem
\end{itemize}
\francais
\end{lemme}

\begin{remarque}{\em
L'isomorphisme canonique 
\[
e \arr{\sim} e \ts e
\]
fait de $C = e$ une cogËbre. Elle n'est pas cocomplËte. Aucun morphisme
non nul $C \ra V$ ne se relËve en un morphisme de cogËbres $C\ra \ctr V.$
}\end{remarque}

Nous notons $\cog $ \label{cog} la catÈgorie des cogËbres diffÈrentielles
graduÈes et  $\cocog$ \label{cogc} la sous-catÈgorie de $\cog$ formÈe des cogËbres
cocomplËtes.
Deux morphismes $f,g :C \ra B$ de cogËbres diffÈrentielles graduÈes
sont {\em homotopes}\index{homotopie} s'il
existe  une $(f,g)$-codÈrivation graduÈe  $H :C \ra B$ de degrÈ $-1$ 
telle que 
\[
f - g = dH + Hd.
\]
Il rÈsultera du thÈorËme \ref{theoreme_cmf_cocog} et du lemme \ref{lemme1_cmf}
que, si la cogËbre graduÈe sous-jacente ‡ $B$ est isomorphe ‡ une cogËbre graduÈe
tensorielle rÈduite, l'homotopie est une relation d'Èquivalence sur l'ensemble des morphismes
de cogËbres de $C$ dans $B.$
%
%
%
%
%
%
%
%
%
%
%
%

\section{$\ai $-algËbres et $\ai$-cogËbres} \label{section_definition_ai-algebres}


%
%
%
%
\subsection{DÈfinitions}

\begin{definition} \label{definition_ai-algebre} {\em 
Soit $n$ un entier $\geq 1.$ Une {\em $\rm{A}_n$-algËbre}
\index{A-n algËbre@{$\rm{A}_n$-algËbre}}
 est un objet $A$ de $\gr \sf C$ muni d'une famille de morphismes
graduÈs
\[
m_i:A^{\ts i}\ra A, \hspace{1cm} 1 \leq i \leq n,
\] 
de degrÈ $2-i$ vÈrifiant, pour tout $1 \leq m\leq n$, l'Èquation
\[
 (*_m) \hspace{1cm}\sum (-1)^{jk+l}m_i(\Id^{\ts j}\ts m_k\ts \Id^{\ts l})= 0 
\]
dans $\Hom_{\gr \sf C} (A\tp{m}, A)$, o˘ les entiers $i,j,k,l$ sont  tels que
$j+k+l=m$ et $i=j+1+l$.
Une {\em $\ai$-algËbre} \index{A-infini algebre@{$\ai$-algËbre}}
 (ou {\em algËbre fortement homotopiquement
associative}) 
est un objet $A$ de $\gr \sf C$ muni de morphismes graduÈs $m_i:A^{\ts i}\ra A$, $i\geq 1$,
de degrÈ $2-i$ vÈrifiant l'Èquation $(*_m)$ pour tout $m\geq 1$. 
}
\end{definition}
\begin{definition} \label{definition_ai-morphisme} {\em 
Un {\em $\rm{A}_n$-morphisme} d'$\rm{A}_n$-algËbres
\index{A-n morphisme@{$\rm{A}_n$-morphisme}}
$f:A \ra B$  est une famille de morphismes graduÈs 
\[
f_i:A^{\ts i}\ra B, \hspace{1cm} 1 \leq i \leq n ,
\]
de degrÈ  $1-i$ vÈrifiant, pour tout $1 \leq m \leq n$,
l'Èquation
\[
 (**_m) \hspace{1cm}\sum (-1)^{jk+l}f_i(\Id^{\ts j}\ts m_k\ts \Id^{\ts l}) = \sum (-1)^s
m_r (f_{i_1}\ts \hdots \ts f_{i_r}) 
\]
dans $\Hom_{\gr \sf C} (A\tp{m},B)$, o˘ les entiers $i,j,k,l$ dans la somme de
gauche sont tels que $j+k+l=m$ et $i=j+1+l$ et
\[
s = \sum_{2\leq u\leq r}\Big((1-i_u)\sum_{1\leq v \leq u}i_v\Big).
\]
Un $\rm{A}_n$-morphisme $f$
 est {\em strict}
\index{strict (A-infini morphisme)@{strict ($\ai$-morphisme)}}
\index{A-infini morphisme@{$\ai$-morphisme}!strict}
si $f_i = 0$ pour tout $i \geq 2$.
 La {\em composition} d'un $\rm{A}_n$-morphisme $f:A \ra B$ avec un $\rm{A}_n$-morphisme
$g:B \ra C$ est dÈfinie par
\[
(gf)_m = \sum_r \sum_{i_1+..+i_r=m} (-1)^s g_r(f_{i_1}\ts \hdots \ts
f_{i_r})
\]
en tant que morphisme de $A\tp{m}$ sur $C$, o˘ $s$ est le mÍme signe que
prÈcÈdemment.
L'{\em identitÈ} de l'$\rm{A}_n$-algËbre $A$ est le $\rm{A}_n$-morphisme
tel que $f_1 = \Id_A$ et $f_i = 0$ si $2 \leq i\leq m$.
Un {\em $\ai$-morphisme} 
\index{A-infini morphisme@{$\ai$-morphisme}}
est une famille de morphismes graduÈs $f_i:A^{\ts i}\ra B$, $i \geq 1$, de degrÈ $1-i$
vÈrifiant l'Èquation $(**_m)$ pour tout $m \geq 1.$
Pour les $\ai$-algËbres, les composantes de la composition et
de l'identitÈ sont dÈfinies
par les mÍmes formules que pour les $\rm{A}_n$-algËbres..
}
\end{definition}

Il rÈsultera de la section \ref{construction_bar_cobar} qu'on obtient
ainsi une catÈgorie.
Notons $\aia$ \indexnotation{aia} la {\em catÈgorie des $\ai$-algËbres}. De mÍme, pour
tout $n \geq 1$, on obtient la {\em catÈgorie $\alg_n$ des $\rm{A}_n$-algËbres}.
\\[.2cm]
\begin{remarque}\label{remarque1_definition_aia} {\em
La dÈfinition des $\ai$-algËbres implique les
formules suivantes  qui expliquent l'autre appellation d'une
${\ai}$-algËbre: algËbre fortement homotopiquement associative.
L'ÈgalitÈ
\[
(*_1)  \hspace{1cm} m_1 \,m_1=0
\]
montre que $(A,m_1)$ est un complexe. L'ÈgalitÈ
\[
(*_2)  \hspace{1cm} m_1 \, m_2 = m_2 \,(m_1 \ts \Id + \Id\ts m_1)
\]
de morphismes  $A\tp{2} \ra A$ signifie que la diffÈrentielle $m_1$ est une
dÈrivation  pour la  multiplication $m_2$.
L'ÈgalitÈ
\[
(*_3)  \hspace{.5cm}  m_2 (m_2 \ts \Id - \Id\ts m_2) = 
 m_1 m_3 + m_3\, (m_1 \ts \Id \ts \Id + \Id\ts m_1 \ts \Id + \Id \ts \Id\ts m_1)
\] 
de morphismes $A\tp{3} \ra A$ exprime que le
dÈfaut d'associativitÈ de $m_2$ est Ègal au bord de
$m_3$  dans le complexe 
\[
\Hom_{\gr \sf C}((A,m_1)\tp{3},(A,m_1)).
\]
Ceci signifie que l'objet graduÈ $A$ muni de la multiplication $m_2$
est une {algËbre dont la multiplication est associative ‡ homotopie prËs}.

 De mÍme, la dÈfinition d'un $\ai$-morphisme $f:A\ra B$ implique les
formules suivantes. L'ÈgalitÈ 
\[
(**_1)  \hspace{1cm} f_1 m_1= m_1 f_1
\]
signifie que le morphisme graduÈ
$f_1$ est un morphisme de complexes. L'ÈgalitÈ 
\[
(**_2)   \hspace{1cm} f_1 m_2 = m_2 \, (f_1 \ts f_1) + m_1 f_2 + f_2 \, (m_1 \ts \Id
+ \Id\ts m_1)
\]
signifie que le
dÈfaut de compatibilitÈ de $f_1$ aux multiplications de $A$ et $B$
est mesurÈ par le bord de $f_2$ dans 
\[
\Hom_{\gr \sf C}((A,m_1)\tp{2},(B,m_1)).
\]
}
\end{remarque}

\begin{remarque}\label{remarque_ai-structure_complexe}
{\em 
Si $(V,d)$ est un complexe,
les morphismes 
\[
m_1 = d, \hspace{1cm} m_i = 0\quad  \rm{ pour }\quad i \geq 2
\]
dÈfinissent une structure d'$\ai$-algËbre sur $V$.
La catÈgorie $\cc \sf C$ des complexes est
une sous-catÈgorie non pleine de $\aia$.}
\end{remarque}

\begin{remarque}\label{remarque_ai-structure_alg}{\em 
Si $((A,d),m)$ est une algËbre diffÈrentielle graduÈe,
les morphismes 
\[
m_1 = d, \hspace{1cm} m_2 = m, \hspace{1cm} m_i = 0\quad  \rm{ pour } \quad i \geq 3
\]
dÈfinissent une structure d'$\ai$-algËbre sur $A$.
RÈciproquement, si dans une $\ai$-algËbre
$A$, les multiplications $m_i$ sont nulles pour $i\geq 3$, le complexe
$(A,m_1)$ muni de la multiplication $m_2 : A \ts A \ra A$ est une algËbre diffÈrentielle
graduÈe. La catÈgorie $\alg$ des algËbres diffÈrentielles graduÈes est
une sous-catÈgorie non pleine de $\aia$.}
\end{remarque}

\begin{definition}
{\em Un {\em $\ai$-quasi-isomorphisme}
\index{A-infini quasi-isomorphisme@{$\ai$-quasi-isomorphisme}}
 $f$ est un $\ai$-morphisme tel
que $f_{1}$ est un quasi-isomorphisme de complexes.}
\end{definition}

\begin{definition}\label{definition_homotopie_ai-morphismes}
{\em
Soit $A$ et $A'$ deux $\ai$-algËbres.
Soit $f$ et $g$ deux $\ai$-morphismes $A \ra A'.$
Une {\em homotopie entre $f$ et $g$}
\index{homotopie!A-infini-morphismes@{$\ai$-morphismes}}
est une famille de morphismes
\[
h_i:A^{\ts i}\ra B, \hspace{1cm} 1 \leq i ,
\]
de degrÈ  $-i$ vÈrifiant, pour tout $1 \leq n$,
l'Èquation $(***_n)$
\[
\begin{array}{rcl}
 f_n - g_n & = & \sum (-1)^s
m_{r+1+t} (f_{i_1}\ts \hdots \ts f_{i_r} \ts h_{k} \ts g_{j_1}\ts \hdots \ts g_{i_t})\\ 
& & + \sum (-1)^{jk+l}h_i(\Id^{\ts j}\ts m_k\ts \Id^{\ts l})
\end{array}
\]
dans $\Hom_{\gr \sf C} (A\tp{n},B)$, o˘ la somme des entiers $i_1, 
\hdots i_r, k, j_1, \hdots j_t$ vaut $n$, o˘ $j+k+l = n$ et
o˘
\[
s = t + \sum_{1\leq \alpha \leq t}(1-j_{\alpha})(n-\sum_{u\geq \alpha}j_u) + k \sum_{1 \leq u \leq r}i_u
+ \sum_{2\leq \alpha \leq r}(1-i_{\alpha})\sum_{u < \alpha}i_u.
\]
Deux $\ai$-morphismes d'$\ai$-algËbres $f$ et $g$ sont {\em homotopes} si il existe une homotopie
entre $f$ et $g$. 
}
\end{definition}

%
%
%
%
%
%
%
%
%
%
%
%
%
%
%
\begin{definition}\label{definition_ai-cogebre} {\em 
 Une {\em $\ai$-cogËbre}\index{A-infini cogËbre
@{$\ai$-cogËbre}} (ou {\em
cogËbre fortement homotopiquement co-associative}) est un objet  $C$ de $\gr \sf C$
muni d'une famille de morphismes graduÈs
\[
\Delta_i : C \ra C \tp i, \hspace{1cm} i \geq 1,
\]
de degrÈ $2-i$ telle que le morphisme
\[
S^{-1}C \arr{} \prod_{i\geq 1}(S^{-1}C)\tp i
\]
dont les composantes sont les
\[
- \si \tp i \circ \Delta_i \circ s \hspace{1cm} (\mbox{o˘}  \hspace{.3cm} \si = s^{-1})
\]
se factorise par le monomorphisme
\[
\bigoplus_{i\geq 1}(S^{-1}C)\tp i \arr{} \prod_{i\geq 1}(S^{-1}C)\tp i
\]
et que, pour tout $m\geq 1$, on a
\[
\sum (-1)^{i+jk}(\Id \tp i \ts \Delta_j \ts \Id \tp k)\Delta_l = 0,
\]
o˘ les les entiers $i,j,k,l$ dans la somme de
gauche sont tels que $i+j+k=m$ et $l=i+1+k$. 

La construction cobar ci-dessous permettra de mieux comprendre cette dÈfinition.
}
\end{definition}

\subsection{Constructions bar et cobar} \label{construction_bar_cobar}
La construction bar est due ‡ S.~Eilenberg et S.~Mac Lane \cite{Eilenberg53} pour les algËbres diffÈrentielles
graduÈes (voir aussi \cite{Cartan55}) et ‡ J. Stasheff \cite{Stasheff63b} pour les $\ai$-algËbres.
Elle permet, entre autres, de reformuler la dÈfinition des $\ai$-structures.
Elle donne aussi une explication (\ref{remarque_signes}) des signes apparaissant dans les Èquations
$(*_m)$ de la dÈfinition des $\ai$-algËbres.
La construction cobar est l'analogue de la construction bar
dans le cas des $\ai$-cogËbres \cite{Adams56}.\\

{\noindent \bf Construction bar}\\

Soit $A$ un objet graduÈ. Soit une famille de morphismes graduÈs
\[
m_i : A \tp i \ra A, \hspace{1cm} i\geq 1,
\]
de degrÈ $2-i.$ Pour tout $i\geq 1$,
nous dÈfinissons une bijection
\[
\begin{array}{rcl}
\Hom_{\gr \sf C} (A\tp i ,A) &\ra &\Hom_{\gr \sf C}((SA)\tp i, SA) \\
m_i & \mapsto & b_i
\end{array}
\]
par la relation
\[
b_i = - s \circ m_i \circ \si \tp i \hspace{1cm} \mbox{o˘}  \hspace{1cm} \si = s^{-1}.
\]
Remarquons que le morphisme
$b_i$ est de degrÈ $+1.$\\

Soit $\ctr (SA)$ la cogËbre tensorielle rÈduite graduÈe sur $SA$.
En vertu du lemme \ref{cogebres_tensorielles}, le morphisme
\[
\bigoplus_{i\geq 1} (SA) \tp i \ra SA
\]
de composantes les $b_i$ se relËve en une unique codÈrivation
\[
b : \ctr (SA) \arr{} \ctr (SA).
\]

\begin{lemme}[J.~Stasheff \cite{Stasheff63b}] \label{bar}
Les propositions suivantes sont Èquivalentes :
\begin{enumerate}
\item[a.] Les $m_i$ dÈfinissent une structure de $\ai$-algËbre sur $A$.
\item[b.] Pour chaque $m\geq 1$, l'Èquation suivante est vÈrifiÈe
\[
\sum_{\substack{j+k+l = m\\j+1+l = i}} b_i(\Id \tp j \ts b_k \ts \Id \tp l) = 0.
\]
\item[c.] La codÈrivation $b$ est une diffÈrentielle, i.~e.~$b^2 = 0$.
\end{enumerate} 
\end{lemme}
\dem L'Èquivalence entre les deux premiers points rÈsulte des
ÈgalitÈs suivantes dans $\Hom_{\gr \sf C}(A\tp i ,SA)$ 
\[
\begin{array}{ccl}
b_i \circ (\Id\tp{j} \ts b_k \ts \Id\tp{l}) & = & s m_i \si \tp i \circ (\Id\tp{j} \ts
s m_k \si \tp k \ts \Id\tp{l})\\ & = &
(-1)^l\ s m_i \circ (\si \tp{j} \ts (m_k\circ \si \tp{k}) \ts \si \tp{l})\\ & = &
(-1)^{l+jk}\
s m_i \circ (\Id\tp{j} \ts m_k \ts \Id\tp{l}) \circ \si \tp n. 
\end{array}
\]
Comme la codÈrivation $b$ est de degrÈ impair, son carrÈ est encore une codÈrivation. Par
le lemme \ref{cogebres_tensorielles}, nous avons donc $b^2 = 0$ si et seulement si $p_1b^2 = 0.$
Ceci montre l'Èquivalence des deux derniers points.
\findem \\

\begin{remarque}[signes]{\em \label{remarque_signes}
Le choix de la bijection $m_i \lra b_i$  n'est
pas canonique. Un autre choix donnerait d'autres signes dans les Èquations
$(*_m)$ de la dÈfinition \ref{definition_ai-algebre}.
}
\end{remarque}

\begin{definition}\label{definition_construction_bar}{\em
La cogËbre diffÈrentielle graduÈe $(\ctr(SA),b)$ associÈe ‡ une
$\ai$-algËbre $A$ est notÈe $BA$ et s'appelle la {\em construction bar}
\index{construction bar} \indexnotation{B} \indexnotation{b} de $A.$}
\end{definition}

Soit $A$ et $A'$ deux $\ai$-algËbres. 
Pour tout $i\geq 1$, nous dÈfinissons une bijection
\[
\begin{array}{rcl}
\Hom_{\gr \sf C} (A\tp i ,A') &\ra &\Hom_{\gr \sf C}((SA)\tp i, SA') \\
\end{array}
\]
par la relation
\[
\si \circ F_i = (-1)^{|F_i|} f_i \circ \si \tp i
\]
o˘ $F_i$ est un morphisme graduÈ de degrÈ $|F_i|.$
Soit une famille de morphismes graduÈs
\[
f_i : A \tp i \ra A', \hspace{1cm} i\geq 1,
\]
de degrÈ $1-i.$
Soit
\[
F : BA \arr{} BA'
\]
le morphisme de cogËbres graduÈes qui relËve le morphisme
\[
\bigoplus_{i\geq 1}(SA) \tp i \arr{} SA'
\]
de composantes les $F_i.$ Une dÈmonstration similaire ‡ celle du lemme
\ref{bar} montre que les $f_i$ dÈfinissent un morphisme d'$\ai$-algËbres si
et seulement si $F$ est compatibles aux diffÈrentielles.
Ainsi, les Èquations $(**_m)$ sont la traduction du fait que la $(F,F)$-codÈrivation
$F \circ b_{BA} - b_{BA'} \circ F$ s'annule.\\


Soit $A$ et $A'$ deux $\ai$-algËbres.
Soit $f$ et $g$ deux $\ai$-morphismes d'$\ai$-algËbres. Notons $F$ et $G$ les
morphismes de cogËbres $BA \ra BA'$ correspondant ‡ $f$ et $g.$ Soit
$H : BA \ra BA'$ une $(F,G)$-codÈrivation de degrÈ $-1$. Elle est dÈterminÈe
(\ref{cogebres_tensorielles}) par sa composÈe avec la projection sur $SA'$
\[
p_1 \circ H : BA \ra SA'.
\]
dont les composantes  sont notÈes
\[
H_i : (SA)\tp i \ra SA', \quad i \geq 1.
\]
Les bijections $F_i \lra f_i$ envoient les morphismes $H_i$ sur des morphismes
$h_i : A \tp i \ra A'$, $i\geq 1$. Ceci dÈfinit une bijection de
l'ensemble des $(F,G)$-codÈrivations de degrÈ $-1$ vers le produit des espaces de morphismes graduÈs
$A \tp i \ra A'$, $i\geq 1$, de degrÈ $-i.$
Cette bijection envoie une homotopie
$H : BA \ra BA'$ entre les morphismes de cogËbres $F$ et $G$ sur l'homotopie entre
les $\ai$-morphismes $f$ et $g$ dÈfinie par la famille
\[
h_i : A \tp i \ra A', \quad i \geq 1.
\]
Les Èquations $(***_m)$ de la dÈfinition \ref{definition_homotopie_ai-morphismes}
sont la traduction de l'Èquation $F - G = \del H.$\\

On obtient ainsi un foncteur $B : \aia \ra \cocog$ appelÈ
le foncteur {\em construction bar}. Il envoie des $\ai$-morphismes
homotopes sur des morphismes homotopes de cogËbres.
La construction bar induit une Èquivalence entre la catÈgorie des $\ai$-algËbres et la
sous-catÈgorie pleine de $\cocog$ formÈe
des cogËbres diffÈrentielles graduÈes dont la cogËbre graduÈe sous-jacente est
isomorphe ‡ une cogËbre graduÈe tensorielle rÈduite.\\

Soit $V$ un objet graduÈ et $n\geq 1$. La sous-cogËbre des $n$-primitifs de $\ctr V$ a
pour espace graduÈ sous-jacent
\[
\bigoplus_{1 \leq i \leq n} V\tp i.
\]
Nous notons $\pctr{n}V$ \indexnotation{pctr}
cette cogËbre. Un raisonnement analogue ‡ celui que nous venons de faire
pour les $\ai$-algËbres permet de construire un foncteur pleinement fidËle
\[
B_n :\alg_n \ra \cocog
\]
qui envoie une $\rm{A}_n$-algËbre $A$ sur la cogËbre
diffÈrentielle graduÈe $(\pctr n (SA),b)$, o˘ $b$ est la diffÈrentielle
construite ‡ l'aide de la bijection  $b_i \lra m_i.$ \\

%
%
%
%
%
%
%
%
%
%
%
%
%
%
%
%
%
%
%
%
%
%
%
%
%
%
%
%
{\noindent \bf Construction cobar}\\

Soit $C$ un objet graduÈ. 
Pour $i\geq 1$, dÈfinissons la bijection
\[
\begin{array}{rcl}
\Hom_{\gr \sf C} (C , C \tp i) & \arr{\sim} & \Hom_{\gr \sf C} (S^{-1}C,(S^{-1}C)\tp i)\\
\Delta_i & \mapsto & {D_i}
\end{array} 
\]
par la relation
\[
D_i = - \si \tp i \circ \Delta_i \circ s.
\]
Soit une famille de morphismes graduÈs
\[
\Delta_i : C \ra C \tp i, \hspace{1cm} i \geq 1,
\]
de degrÈ $2-i$ telle que le morphisme
\[
S^{-1}C \arr{} \prod_{i\geq 1}(S^{-1}C)\tp i
\]
dont les composantes sont les $D_i$, $i\geq 1$,
se factorise par le monomorphisme
\[
\bigoplus_{i\geq 1}(S^{-1}C)\tp i \arr{} \prod_{i\geq 1}(S^{-1}C)\tp i.
\] 
Gr‚ce au lemme \ref{algebres_libres}, le morphisme graduÈ $S^{-1}C \ra \b TS^{-1}C$
obtenue ainsi
s'Ètend en une unique dÈrivation d'algËbres de $\b TS^{-1}C.$
A l'aide du lemme \ref{algebres_libres}, nous montrons qu'on a $D^2 = 0$ si et
seulement si les $\Delta_i$ dÈfinissent une structure d'$\ai$-cogËbre sur $C.$
Ainsi, les diffÈrentielles de l'algËbre $\b TS^{-1}C$ sont en bijection avec
les structures d'$\ai$-cogËbre
sur l'objet graduÈ $C$.

\begin{definition}\label{definition_construction_cobar}{\em 
On note $\Omega C$ l'algËbre diffÈrentielle graduÈe $(\b TS^{-1}C,D)$ associÈe ‡
une $\ai$-cogËbre $C$. Elle s'appelle la {\em construction cobar}
\index{construction cobar} \indexnotation{Omega} de $C$.}
\end{definition}
La catÈgorie $\aic$ \indexnotation{aic} des $\ai$-cogËbres a pour objets les $\ai$-cogËbres. On dÈfinit
ses morphismes de telle maniËre que la construction cobar
\[
\Omega : \aic \arr{} \alg
\]
devienne un foncteur pleinement fidËle. La catÈgorie des cogËbres diffÈrentielles
graduÈes s'identifie alors ‡ une sous-catÈgorie (non pleine) de la catÈgorie des
$\ai$-cogËbres.

On note encore $B$ (resp.~$\Omega$) la restriction de la construction
bar (resp.~cobar) aux algËbres (resp.~cogËbres cocomplËtes) diffÈrentielles graduÈes.
\begin{lemme}\label{adjonction_bar_cobar}
Le foncteur $\Omega : \cocog \ra \alg$ est adjoint ‡ gauche
au foncteur $B : \alg \ra \cocog.$
\end{lemme}

\dem Ce lemme est bien connu.
Soit $A$ une algËbre et $C$ une cogËbre cocomplËte. Il s'agit de montrer
que nous avons un isomorphisme fonctoriel
\[
\Hom_{\sz \cocog}(C,BA) \arr{\sim} \Hom_{\sz \alg}(\Omega C,A).
\]
Soit $F : C \ra BA$ un morphisme de cogËbres. Comme $BA$ est une cogËbre tensorielle
rÈduite en tant que cogËbre graduÈe, la donnÈe de $F$ Èquivaut (\ref{cogebres_tensorielles})
‡ celle de
\[
f = p_1F : C \ra SA.
\]
Posons $\tau = \si \circ f.$
La condition $d_{BA} \circ F - F \circ d_C = 0$ se traduit par le fait
que $\tau$ est une
{\em 
cochaÓne tordante}\index{cochaÓne tordante}, c'est-‡-dire que l'on a
\[
 d_A \circ \tau + \tau \circ d_C +  m   \circ \tau \tp 2 \circ \Delta = 0.
\]
Le morphisme graduÈ  $f' = \tau \circ s$ s'Ètend de maniËre unique (\ref{algebres_libres})
en un morphisme
d'algËbres $F' : \Omega C \ra A$ car $\Omega C$ est
libre sur $S^{-1}C$ en tant qu'algËbre graduÈe.
La compatibilitÈ de $F'$ ‡ la diffÈrentielle est Èquivalente au
fait que $\tau$ est une cochaÓne tordante. \findem

%
%
%
%
%
%
%
%
%
%
%
%
%
%
%
%
\section{$\cocog$ comme catÈgorie de modËles} \label{section_cmf_cocog}

%
%

{\noindent \bf Plan de la section}\\

\noindent Cette section est divisÈe en cinq sous-sections.

Dans la premiËre sous-section (\ref{section_cmf_alg}), nous rappelons \cite{Hinich97c}
la structure de catÈgorie de modËles sur
la catÈgorie $\alg$ des algËbres diffÈrentielles graduÈes. Nous ÈnonÁons
le thÈorËme principal (\ref{theoreme_cmf_cocog}) et en dÈduisons {\em
thÈorËme des $\ai$-quasi-isomorphismes} (\ref{corollaire_cmf_aia}.~{\it a}) et le {\em thÈorËme de l'homotopie}
(\ref{corollaire_cmf_aia}.~{\it b}).

Dans la deuxiËme sous-section \ref{section_theoreme_cmf_cocog}, nous montrons le
thÈorËme principal (\ref{theoreme_cmf_cocog}). Pour la caractÈrisation des objets fibrants
de $\cocog$, nous aurons besoin de
certains rÈsultats de la sous-section suivante.

Dans la sous-section \ref{section_cmf_aia}, nous montrons que la catÈgorie $\aia$ admet une structure
de ``catÈgorie de modËles sans limites'' (\ref{theoreme_cmf_aia}). Nous montrons ensuite que
la construction bar $B : \aia \ra \cocog$
est compatible aux structures de catÈgorie de modËles (``sans limites'')
de $\aia$ et $\cocog$ (\ref{proposition1_cmf_aia}).

Dans la sous-section \ref{section_homotopie}, nous comparons 
l'homotopie ‡ gauche (au sens des catÈgories de modËles)
avec l'homotopie ``au sens classique'' sur les morphismes de cogËbres
diffÈrentielles graduÈes cocomplËtes.

Dans la sous-section \ref{section_equiv_qis}, nous comparons
les Èquivalences faibles de $\cocog$ avec les quasi-isomorphismes de $\cocog$.

\subsection{Le thÈorËme principal}
\label{section_theoreme_cmf_cocog}

Le lecteur qui n'est pas familier avec les catÈgories de modËles au sens de
Quillen trouvera dans l'appendice \ref{section_rappels_cmf} quelques rappels de certains ÈnoncÈs clefs
et les rÈfÈrences classiques.\\

\noindent{\bf La catÈgorie de modËles $\alg$}\\
\label{section_cmf_alg}

Dans la catÈgorie $\alg$ des algËbres diffÈrentielles $\Z$-graduÈes (\ref{alg}),
considÈrons les trois classes de morphismes suivantes :
\english \begin{itemize}
\item[-] la classe $Qis$ des quasi-isomorphismes,
\item[-] la classe $\fib$ des morphismes $f : A \ra B$ tels que
$f^n$ est un Èpimorphisme pour tout $n \in \Z,$
\item[-] la classe $\cof$ des morphismes qui ont la propriÈtÈ de relËvement
‡ gauche par rapport aux morphismes appartenant ‡ $Qis \cap \fib$.
\end{itemize} \francais
Soit $\sf E$  l'une des sous-catÈgories pleines de $\alg$ dont les objets sont
respectivement
\begin{description}
\item[(I)] les algËbres $A$ telles que $A^p = 0$ pour tout $p > 0$,
\item[(II)] les algËbres $A$ telles que $A^p = 0$ pour tout $p \leq 0$.
\end{description}
H. Munkholm a dÈmontrÈ dans \cite{Munkholm78} que $\sf E$  devient
une catÈgorie de modËles si on la munit de $\sf E \cap Qis$, $\sf E \cap \fib$ 
 et de la
classe des morphismes de $\sf E$ qui ont la propriÈtÈ de relËvement 
‡ gauche par rapport aux morphismes de $\sf E \cap Qis \cap \fib$.
Le rÈsultat de H.~Munkholm a ÈtÈ renforcÈ par V.~Hinich :

\begin{theoreme}[Hinich \cite{Hinich97c}] \label{cmf_algebres}
La catÈgorie $\alg$ munie des classes de morphismes dÈfinies ci-dessus
est une catÈgorie de modËles. Les algËbres cofibrantes sont
les algËbres qui sont isomorphes ‡ une algËbre presque libre.
Toutes les algËbres sont fibrantes. \findem
\end{theoreme}
 Le cas plus gÈnÈral o˘ l'anneau de base n'est pas
un corps est d˚ ‡ J.~F.~Jardine \cite{Jardine97}. S. Schwede et B. Shipley \cite{SchwedeShipley02}
ont gÈnÈralisÈ ces rÈsultats ‡ des catÈgories d'algËbres dans des catÈgories de
modËles monoÔdales\footnote{``monoÔdales'' se rapporte ‡ ``catÈgories'' ("modËle" est masculin).}.\\

%
%
%
%
%
%
%
%
%
%
%
%
%
%
%
%
%
%
%


\noindent{\bf Le thÈorËme principal et ses consÈquences}\\

Dans la catÈgorie $\cocog$ des cogËbres cocomplËtes
diffÈrentielles graduÈes, nous considÈrons les trois classes de morphismes
suivantes :

\english \begin{itemize}
\item[-] la classe $\weq $ des {\em Èquivalences
faibles} est formÈe des morphismes $f : C\ra D$ tels que
 $\Omega F : \Omega C \ra \Omega D$ est un quasi-isomorphisme 
d'algËbres,

\item[-] la classe $\cof$ des {\em cofibrations} est formÈe des morphismes
$f : C \ra D$ qui, en tant que morphismes de complexes, sont des
monomorphismes,

\item[-] la classe $\fib$ des {\em fibrations} est formÈe des morphismes
qui ont la propriÈtÈ de relËvement ‡ droite par rapport aux cofibrations triviales.
\end{itemize} \francais

Il s'avËre que la classe des Èquivalences faibles est strictement incluse dans la
classe des quasi-isomorphismes de cogËbres (voir \ref{section_equiv_qis}). 
De l'autre
cÙtÈ, il est bien connu (et nous le redÈmontrerons, voir la proposition \ref{proposition_equiv_qis})
qu'un quasi-isomorphisme entre cogËbres cocomplËtes
est une Èquivalence faible si les deux cogËbres sont concentrÈes en degrÈs $<-1$ ou en degrÈs
$\geq 0$.

\begin{theoreme} \label{theoreme_cmf_cocog}
\english \begin{itemize}
\item[a.]
La catÈgorie $\cocog$ munie des trois classes de morphismes ci-dessus est une catÈgorie de modËles.
Tous ses objets sont cofibrants. Un objet de $\cocog$ est fibrant si et seulement
si sa cogËbre graduÈe sous-jacente est isomorphe ‡ une cogËbre
tensorielle rÈduite.
\item[b.]
Munissons la catÈgorie $\alg$ de la structure de catÈgorie de
modËles du thÈorËme \ref{cmf_algebres}.
La paire de foncteurs adjoints $(\Omega, B)$ de $\cocog$ dans $\alg$ est une
Èquivalence de Quillen.
\end{itemize} \francais
\end{theoreme}

\dem Voir la section suivante \ref{section_demonstration_thm_principal}. \findem \\

Le point {\it b} du thÈorËme renforce 
des thÈorËmes classiques (voir \cite{Moore70}, \cite[th.~4.4 et 4.5]{Husemoller74}).
Il semble Ítre nouveau sous la forme que nous donnons. Notre dÈmonstration est
une adaptation de celle de Hinich \cite{Hinich01}, basÈe ‡ son tour sur celle
de Quillen \cite{Quillen69}. Le fait que les foncteurs
bar et cobar induisent des Èquivalences inverses
l'une de l'autre dans les catÈgories homotopiques est non trivial mais sa
dÈmonstration n'est pas trËs difficile.  DÈduisons maintenant, gr‚ce aux techniques d'algËbre
homotopique de Quillen (voir appendice \ref{section_rappels_cmf})
le {\em thÈorËme des $\ai$-quasi-isomorphismes}, le {\em thÈorËme de l'homotopie}
et la gÈnÈralisation du thÈorËme \cite[Thm.~6.2]{Munkholm78} de H.~J.~Munkholm.

\begin{corollaire} \label{corollaire_cmf_aia}
\english \begin{itemize}
\item[a.] La relation d'homotopie (voir \ref{definition_homotopie_ai-morphismes})
dans $\aia$ est une relation d'Èqui\-valence.
\item[b.] Un quasi-isomorphisme d'$\ai$-algËbres est une Èquivalence d'homotopie (i.~e.~un
isomorphisme dans la catÈgorie quotient de $\aia$ par la relation d'homotopie).
\item[c.] Soit $\mathsf{dash}$ la sous-catÈgorie pleine  de $\aia$ formÈe des algËbres diffÈren\-tielles
graduÈes. Notons $\sim$ la relation d'homotopie sur $\mathsf{dash}$. L'inclu\-sion 
$\alg \hookrightarrow \mathsf{dash}$ induit une Èquivalence
\[
\alg[Qis^{-1}] \arr{\sim} \mathsf{dash}/\!\sim.
\] 
\end{itemize} \francais
\end{corollaire}

L'idÈe du point {\it c} remonte ‡ J.~Stasheff et S.~Halperin \cite{Stasheff70}. Il a ÈtÈ
dÈmontrÈ sous les conditions {\bf (I)} ou {\bf (II)} (voir en haut de cette section) par H.~J.~Munkholm
\cite{Munkholm78}. Les points {\it a} et {\it b} sont connus (surtout parmi les
spÈcialistes de l'homotopie rationnelle) depuis le dÈbut des annÈes 80, au moins pour des
$\ai$-algËbres connexes (i.~e.~concentrÈes en degrÈs homologiques $\geq 1$),
voir par exemple A.~ProutÈ \cite[chap.~4]{Proute85} ou T.~V.~Kadeishvili
\cite{Kadeishvili87}.\\

\dem 
Par la section \ref{construction_bar_cobar}, nous savons que deux morphismes
d'$\ai$-algËbres 
\[
f,g :A\ra A'
\]
sont homotopes si et seulement si $Bf$ et $Bg$ sont des
morphismes de cogËbres homotopes. Par le thÈorËme principal (\ref{theoreme_cmf_cocog}),
la cogËbre $BA'$ est fibrante dans $\cocog$ et tout objet de $\cocog$ est cofibrant.
Acceptons provisoirement (voir \ref{proposition1_homotopie} plus bas)
le rÈsultat suivant : la relation d'homotopie au sens classique
sur $\Hom_{\sz \cocog}(BA,BA')$ est Ègale ‡ la
relation d'homotopie ‡ gauche pour la catÈgorie de modËles $\cocog$.

{\it a.} C'est le lemme \ref{lemme1_cmf} appliquÈ ‡ la catÈgorie de
modËles fermÈe $\cocog$.

{\it b.} Un $\ai$-quasi-isomorphisme $f : A \ra A'$ induit (voir \ref{proposition1_cmf_aia}
plus bas) un morphisme 
\[
Bf : BA \ra BA',
\]
qui est une Èquivalence faible de $\cocog$ entre objets
fibrants et cofibrants. Il est donc inversible ‡ homotopie prËs dans $\cocog$
(voir la proposition \ref{proposition_cmf}).

{\it c.} Par le thÈorËme principal \ref{theoreme_cmf_cocog}, le foncteur
 $B$ induit une Èquivalence
\[
\alg[Qis^{-1}] = \Ho \alg \arr{\sim} \Ho \cocog.
\]
Nous avons l'Èquivalence (voir \ref{proposition_cmf})
\[
\cocog_{\sf{cf}}/\!\sim \, \arr{\sim} \Ho \cocog.
\]
Le foncteur $B$ prend ses valeurs dans $\cocog_{\sf{cf}}$.
Il induit donc une Èquivalence
\[
\alg[Qis^{-1}] \arr{\sim} \cocog_{\sf{cf}}/\!\sim.
\]
Son image est isomorphe ‡ $\sf{dash}/\!\sim$. \findem

\subsection{DÈmonstration du thÈorËme principal}
\label{section_demonstration_thm_principal}

Notre dÈmonstration du thÈorËme principal \ref{theoreme_cmf_cocog} nÈcessite l'Ètude prÈalable des algËbres
et des cogËbres filtrÈes.\\

%
%
%
%
%
%
%
%
%
%
%
%

 {\noindent \bf Objets filtrÈs}\\

Soit $\sf M$ l'une des catÈgories $\gr \sf C$ ou $\cc \sf C$.
Une {\em filtration  d'un objet $X$} \index{filtration} de $\sf M$  est une suite croissante
\[
X_0 \subset X_1 \subset \cdots \subset X_i \subset X_{i+1}\subset \cdots, \quad i\in \N
\]
de sous-objets de $X$.
Elle est {\em exhaustive} \index{exhaustive (filtration)} \index{filtration!exhaustive}si l'on a
\[
\colim X_i = X.
\]
Elle est {\em admissible} \index{admissible!filtration} \index{filtration!admissible}
si elle est exhaustive et si $X_0 = 0$.
Un {\em objet filtrÈ} \index{objet filtrÈ}de $\sf M$ est un objet de $\sf M$ muni d'une
filtration. Soit $X$ et $Y$ deux objets filtrÈs.
L'{\em objet graduÈ $\ogr X$ associÈ ‡ $X$} est dÈfini par la suite d'objets de $\sf M$
\[
\ogr_0 = X_0,\quad {\ogr}_iX = X_i/X_{i-1},\quad i\geq 1.
\]
Un morphisme $f : X \ra Y$ de $\sf M$ est un {\em morphisme d'objets filtrÈs}
si on a
\[
f(X_i) \subset Y_i
\]
pour tout $i \in \N$.
Le produit tensoriel $X \ts Y$ 
est muni de la filtration dÈfinie par la suite
\[
(X \ts Y)_i = \sum_{p+q=i} X_p \ts Y_q, \quad i \in \N. 
\]
Ceci munit la catÈgorie des objets filtrÈs de $\sf M$ d'une structure
de catÈgorie monoÔdale dont l'ÈlÈment neutre est l'objet $e$ muni de la filtration
$e_i = e$, $i \in \N.$ 
La suspension $SX$ de l'objet de $\sf M$ sous-jacent ‡ $X$ est munie de la
filtration donnÈe par $(SM)_i = SM_i$, $i\in \N$.\\

Un {\em complexe filtrÈ} \index{complexe!filtrÈ} est un objet filtrÈ de $\cc \sf C.$

\begin{definition}{\em 
Soit $X$ et $Y$ deux complexes filtrÈs.
Un morphisme  $f : X\ra Y$ est un {\em quasi-isomorphisme filtrÈ} \index{quasi-isomorphisme filtrÈ}
si  les morphismes 
\[
\ogr_i C \ra \ogr_i D,  \quad i\in \N,
\]
induits par $f$
sont des quasi-isomorphismes de complexes.
}
\end{definition}

Une {\em algËbre filtrÈe} (resp.~{\em cogËbre filtrÈe})
est une algËbre (resp.~cogËbre) dans la catÈgorie des complexes filtrÈs.
Une {\em cogËbre filtrÈe admissible} \index{admissible!cogËbre} \index{cogebre@{cogËbre}!admissible}
est une cogËbre $C$ munie d'une filtration
admissible. Notons qu'on a alors
\[
\Delta C_{i+1} \subset C_i \ts C_i, \quad i\in \N.
\]
Nous montrerons   (\ref{lemme4_cmf_cocog}) que tout  quasi-isomorphisme filtrÈ entre cogËbres filtrÈes
admissibles est une Èquivalence faible de $\cocog$.\\

Soit $C$ une cogËbre filtrÈe, cocomplËte en tant que cogËbre.
La filtration de $C$ induit une filtration sur chaque puissance tensorielle
de $S^{-1}C$. Nous obtenons ainsi une filtration d'algËbre sur la construction
cobar $\Omega C$.
Soit $C$ et $D$ deux cogËbres filtrÈes cocomplËtes. Munissons les constructions cobar $\Omega C$
 et $\Omega D$ des filtrations induites par celles de $C$ et $D$.
La construction cobar envoie un morphisme de cogËbres filtrÈes $f : C\ra D$ sur un
morphisme d'algËbres filtrÈes $\Omega f : \Omega C \ra \Omega D.$

Soit $A$ une algËbre filtrÈe. La filtration de $A$ induit une filtration de
cogËbre sur la construction bar $BA$ de $A.$ Soit $A$ et $A'$ deux algËbres
filtrÈes. Munissons les constructions bar $BA$ et $BA'$ des filtrations induites par
celles de $A$ et de $A'.$ La construction bar envoie un morphisme d'algËbres
filtrÈes $f :A \ra A'$ sur un morphisme de cogËbres filtrÈes $Bf :BA \ra BA'.$

Soit $C$ une cogËbre cocomplËte. La {\em filtration primitive} \index{filtration!primitive}
de la cogËbre $C$ est dÈfinie
par la suite des sous-cogËbres des $i$-primitifs  $C \prim{i}$, $i\geq 1$, complÈtÈe par $C \prim{0} = 0$.
Comme la catÈgorie de base $\sf C$ est semi-simple, la filtration primitive de $C$ est une
filtration de cogËbre. 
Elle est admissible
et induit une filtration sur $\Omega C$, qui induit
une filtration sur la construction bar $B \Omega C.$ Nous appelons cette derniËre
filtration la  {\em filtration $C$-primitive } \index{filtration!C-primitive@{$C$-primitive}} de $B\Omega C$.

%
%
%
%
%
%
%
%
%
%
%
%

\begin{lemme}\label{lemme4_cmf_cocog}
Un quasi-isomorphisme filtrÈ de cogËbres filtrÈes admissibles est une Èquivalence faible.
\end{lemme}

\dem Soit $C$ et $D$ deux cogËbres admissibles et $f : C \ra D$ un quasi-isomorphisme
filtrÈ. Nous allons montrer que le morphisme d'algËbre
\[
\Omega f : \Omega C \ra \Omega D
\]
est un quasi-isomorphisme filtrÈ pour les filtrations de $\Omega C$ et $\Omega D$
induites par celles de $C$ et $D$. On rappelle que la diffÈrentielle de $\Omega C$ est l'unique
codÈrivation $d$ qui relËve le morphisme 
\[
S^{-1}C \ra \Omega C
\]
de composantes non nulles $\si d_C s$ et
$\si  \Delta s\tp 2.$
Munissons  $\Omega C$ de la filtration
induite par celle de $C$. Soit $i \geq 1.$ Comme la filtration de $C$ est admissible,
$\ogr_i (C\tp j) = 0$ si $j > i.$
Munissons
\[
\ogr_i \Omega  C = \ogr_i \Big(\bigoplus_{1\leq j\leq i}C\tp j \Big)
\]
de la filtration
\[
F_l = \ogr_i\Big(\bigoplus_{i-l\leq j\leq i}C\tp j\Big), \quad l\geq 0.
\]
La contribution de $\si \Delta s\tp 2$ dans la diffÈrentielle $d$ de $\ogr_i \Omega  C$
fait dÈcroÓtre la filtration. Ainsi, seul le morphisme $\si d_C s$ contribue
‡ la diffÈrentielle de l'objet graduÈ associÈ au $F_l$, $l\geq 1$.
Le morphisme
\[
\ogr_i \Omega  C  \arr{}
\ogr_i \Omega  D 
\]
est  filtrÈ pour cette filtration et il induit clairement un quasi-isomorphisme
dans les objets graduÈs. \findem


\begin{lemme}\label{lemme1_cmf_cocog}
\english \begin{itemize}
\item[a.] Soit $A$ et $A'$ deux algËbres diffÈrentielles graduÈes.
La construction bar envoie un quasi-isomorphisme d'algËbres $f :A \ra A'$
sur un quasi-isomorphisme filtrÈ $f : BA \ra BA'$ pour la filtration primitive.
 
\item[b.] Soit $A$ une algËbre diffÈrentielle graduÈe. Le morphisme d'adjonction
\[
\phi : \Omega B A \arr{} A
\]
est un quasi-isomorphisme d'algËbres.
\item[c.] Soit $C$ une cogËbre cocomplËte. Munissons $C$ de la filtration 
primitive et $B\Omega C$ de la filtration $C$-primitive. Le morphisme d'adjonction
\[
\psi : C \arr{} B \Omega C
\]
est un quasi-isomorphisme filtrÈ.
\end{itemize} \francais
\end{lemme}

\dem {a.} La filtration primitive de $BA$ a pour objet graduÈ associÈ
\[
\ogr_i(BA) = (SA)\tp i,\quad i \in \N.
\]
Par le thÈorËme de K¸nneth, un quasi-isomorphisme $f : A \ra A'$ induit un quasi-isomorphisme
dans ces sous-quotients.

{b.} Nous allons donner des filtrations exhaustives sur $A$ et $\Omega B A$ de faÁon ‡
ce que le morphisme d'adjonction devienne un quasi-isomorphisme filtrÈ.
Soit la filtration de $A$ dÈfinie par $A_i = A,$ $i \geq 1$ et
$A_0 = 0.$ Munissons $\Omega BA$ de la filtration induite  par
la filtration primitive de $BA.$ Le morphisme d'adjonction
\[
\phi : \Omega B A \arr{} A
\]
est clairement un morphisme filtrÈ. Il induit un morphisme 
\[
\ogr_i (\Omega B A) \arr{} \ogr_i A ,\quad i \in \N,
\]
dans les objets graduÈs qui est l'identitÈ de $A$ si $i =1$, et qui
est nul si $i \geq 2.$ Pour montrer que le morphisme d'adjonction
est un quasi-isomorphisme, il suffit donc de montrer que, pour $i \geq 2$,
le complexe $\ogr_i (\Omega B A)$
est contractile. Soit le complexe $V = SA$.
Remarquons que nous avons un isomorphisme de complexes
\[
\bigoplus_{i\geq 1} \ogr_i (\Omega B A) \arr{\sim } \Omega \ctr{V}
\]
 qui identifie ‡ la composante
$\ogr_i (\Omega B A)$, $i \geq 1$, ‡ la somme des
\[
S^{-1}V\tp{i_1} \ts \hdots \ts S^{-1}V\tp{i_k} \subset (S^{-1}\ctr{V})\tp k,
\]
o˘ $k \geq 1$ et o˘ $i_1 + \hdots + i_k = i.$
Soit $i \geq 2.$
Soit  le morphisme graduÈ $r :\ogr_i  (\Omega B A)\ra \ogr_i  (\Omega B A)$
de degrÈ $-1$ donnÈ par les morphismes
\[
S^{-1}V\tp{ i_1} \ts S^{-1}V\tp{i_2}\ts \hdots \ts S^{-1}V\tp{ i_{k} }
\ra S^{-1}V\tp{ i_1+i_2} \ts \hdots \ts S^{-1}V\tp {i_{k}}
\]
que nous dÈfinissons comme nuls si $i_1 \neq 1$ et valant $\eta \circ (s \ts \Id \tp k)$ sinon ;
ici $\eta$ est l'isomorphisme naturel
\[
V \ts S^{-1}V\tp{i_2} \ts \hdots \ts S^{-1}V\tp{i_{k}} \arr{\sim}S^{-1}V\tp{1+i_2} \ts \hdots 
\ts S^{-1}V\tp{i_{k}}.
\]
Nous vÈrifions que le morphisme graduÈ $r$ est une homotopie contractante du complexe $\ogr_i  (\Omega B A)$.

{\it c.} Nous devons montrer que le morphisme de complexes
\[
\psi : \ogr C \ra \ogr (B\Omega C)
\]
est un quasi-isomorphisme. Posons $W = \ogr (S^{-1}C)$. Comme $C$ est
admissible, la comultiplication de $\ogr C$ est nulle et
\[
\ogr (B\Omega C) \arr{\sim} B \Omega (\ogr C)
\]
est la somme des complexes
\[
V_i = \bigoplus SW \tp{i_1} \ts  \hdots \ts SW\tp{i_k}, \quad i \geq 1,
\]
o˘ $k \geq 1$ et $i_1 + \hdots + i_k = i.$ La composÈe du morphisme
\[
\ogr C \ra \ogr B \Omega C
\]
avec la projection sur les $V_i$ est nulle si $i \geq 2$ et c'est l'identitÈ
de $\ogr C$ si $i = 1$. Il reste ‡ montrer que $V_i$ est contractile pour
$i \geq 2$.
Soit $i \geq 2.$
Soit  $r : V_i \ra V_i$ le morphisme graduÈ de degrÈ $-1$
dÈfini par les
 morphismes 
\[
\xymatrix{
SW \tp{i_1 + i_2} \ts  \hdots \ts SW\tp{i_k} \ar[r] & 
SW \tp{i_1} \ts SW \tp{i_2} \ts  \hdots \ts SW\tp{i_k}}
\]
que nous dÈfinissons comme nuls si $i_1 \neq 1$ et, comme $\eta  \circ (s\ts \Id \tp {i-1})$ sinon ; ici
$\eta$ est l'isomorphisme naturel
\[
 S^{-1}W\tp{1+i_2} \ts \hdots 
\ts S^{-1}W\tp{i_{k}}\arr{\sim} W \ts S^{-1}W\tp{i_2} \ts \hdots \ts S^{-1}W\tp{i_{k}}.
\]
Nous vÈrifions que le morphisme $r$ est une homotopie contractante de $V_i.$
\findem \\

{\noindent \bf DÈmonstration du thÈorËme principal \ref{theoreme_cmf_cocog}}\\

CommenÁons par quelques lemmes prÈliminaires.

\begin{lemme} \label{lemme2_cmf_cocog}
Soit $C$ une cogËbre et $C'$ une sous-cogËbre de $C$ telle que
$\Delta C \subset C' \ts C'.$
La construction
cobar envoie l'inclusion $C' \hookrightarrow C$ sur une cofibration
standard (\ref{cofibration_standard}).
\end{lemme}

Pour dÈmontrer ce lemme et le suivant, nous aurons besoin de la description suivante \cite{Hinich97c}
des cofibrations de $\alg$ : notons $A^\sharp$ le complexe sous-jacent
‡ une algËbre diffÈrentielle graduÈe $A$ et $FV = \b TV$ l'algËbre diffÈrentielle
graduÈe libre sur le complexe $V$.
Soit $A$ une algËbre diffÈrentielle graduÈe et $M$ un complexe.
Soit $\alpha : M \ra A^{\sharp}$ un morphisme de complexes.
On note $C(\alpha)$ le cÙne de $\alpha$ dans la
catÈgorie $\cc\sf C$. Notons
$A\langle M,\alpha\rangle $ \label{A<M,a>} la colimite dans $\alg$ du diagramme
\[
A \la F(A^{\sharp}) \ra FC(\alpha).
\]

\begin{definition}{\em  \label{cofibration_standard}
Un morphisme $f : A \ra B$ est une {\em cofibration standard}
\index{cofibration standard (de $\alg$)}
 s'il est la colimite
d'une suite de morphismes composÈs
\[
A = A_0 \ra A_1 \ra \hdots \ra A_{n-1} \ra A_n, \quad n\geq 1,
\]
o˘ toutes les flËches $A_i \ra A_{i+1}$ sont donnÈes par les morphismes canoniques 
\[
A_i \ra A_i\langle M_i,\alpha_i\rangle = A_{i+1}
\]
pour des morphismes de complexes $\alpha_i : M_i \ra A_i^{\sharp}$. 
Une {\em cofibration standard triviale}
est une cofibration standard telle que
tous les complexes $M_i$ sont contractiles (i.~e.~isomorphes ‡ $0$ dans 
$\ch \sf C.$)
}
\end{definition}

Les faits suivants sont dÈmontrÈs dans $\cite{Hinich97c}$ :
Toute cofibration est rÈtract d'une cofibration
standard. De mÍme, toute cofibration triviale est rÈtract d'une
cofibration standard triviale.\\

{\em DÈmonstration du lemme \ref{lemme2_cmf_cocog}} :
Soit $E$ le conoyau dans la catÈgorie des complexes
de l'inclusion $C' \hookrightarrow C.$ Choisissons une section de $C \ra E$
dans la catÈgorie graduÈe pour obtenir un isomorphisme
\[
C' \oplus E \arr{\sim} C
\]
d'objets graduÈs.
En tant qu'algËbre graduÈe,
la construction cobar $\Omega C = \Omega (C' \oplus E )$ est
isomorphe au coproduit d'algËbres graduÈes
\[
FS^{-1} C' \amalg FS^{-1} E,
\]
o˘ $F = \b T$ comme dans (\ref{section_cmf_alg}).
La diffÈrentielle de $\Omega C$ est induite par 
la comultiplication de $C$ et la diffÈrentielle du complexe  $C$.
Selon la dÈcomposition $C = C' \oplus E$, la comultiplication de $C$ est donnÈe par deux
composantes 
\[
\Delta_{C'} : C' \ra C' \ts C' \quad \mbox{et} \quad \Delta_{E} : E \ra C' \ts C',
\]
et la diffÈrentielle de $C$ est donnÈe par la diffÈrentielle de $C'$, celle de $E$ et un morphisme
graduÈ $d : E \ra C'$ de degrÈ $+1$.
Soit le morphisme de complexes
\[
[D_1,D_2] : S^{-2}E \arr{}  S^{-1}C'\oplus \left(S^{-1}C' \ts S^{-1}C' \right)
\]
dont les composantes sont dÈfinies par  $s \tp 2 \circ D_2  = \Delta_E \circ s^2$ et
par $s \circ D_1  = d \circ s^2.$
Nous notons
\[
D : S^{-2}E \arr{} FS^{-1} C' \amalg FS^{-1} E
\]
sa composition avec l'injection de $ S^{-1}C'\oplus \left(S^{-1}C' \ts S^{-1}C' \right)$
dans $ FS^{-1} C' \amalg FS^{-1} E.$
Par construction, l'algËbre diffÈrentielle graduÈe
\[
\Omega C'\langle S^{-2}E, D \rangle
\]
est l'algËbre graduÈe $ FS^{-1} C' \amalg FS^{-1} E$ dont la diffÈrentielle est
induite par la comultiplication de $C'$, les
diffÈrentielles des complexes $C'$ et $E$,  le morphisme $\Delta_{E}$ et le morphisme
$d$. Elle est donc isomorphe ‡ $\Omega C$ en tant qu'algËbre diffÈrentielle graduÈe. 
\findem

\begin{lemme} \label{lemme3_cmf_cocog}
\english \begin{itemize}
\item[a.] La construction cobar prÈserve les cofibrations et les
Èquivalences faibles.
\item[b.] La construction bar prÈserve les fibrations et les
Èquivalences faibles.
\end{itemize} \francais
\end{lemme}

\dem {\it a.} Soit $i : C \rightarrowtail D$ une cofibration de cogËbres.
Soit la filtration de $D$ dÈfinie par la suite des $D_i = i(C) + D\prim i$, $i\in \N$.
Remarquons que $D_0$ est isomorphe ‡ $C$ et que, pour tout $i \geq 1$, on a
\[
\Delta (D_{i+1}) \subset D_i \ts D_i.
\]
Nous pouvons donc appliquer le lemme \ref{lemme2_cmf_cocog}. Il certifie
que $\Omega D_i \ra \Omega D_{i+1}$ est une cofibration standard.
Le morphisme $\Omega C \ra \Omega D$ est la composition dÈnombrable
des cofibrations standard $\Omega D_i \rightarrowtail \Omega D_{i+1}$. Il est donc aussi une
cofibration standard. La construction cobar prÈserve les Èquivalences faibles
par dÈfinition des Èquivalences faibles de $\cocog$.

{\it b.} Soit $p : A \twoheadrightarrow A'$ une fibration d'algËbres. Le morphisme
$Bf$ est une fibration s'il vÈrifie la propriÈtÈ de relËvement ‡ droite par rapport
aux cofibrations triviales $i : C \ra D$ de cogËbres. Gr‚ce ‡ l'adjonction entre les constructions
bar et cobar, cette propriÈtÈ est Èquivalente au fait que $\Omega i$ ait la propriÈtÈ de
relËvement ‡ gauche par rapport ‡ $p.$ Mais ceci est toujours vrai par le point {\it a.}
Le morphisme $Bf$ est donc une fibration de $\cocog$.

Soit $f : A \arr{\sim} A'$ un quasi-isomorphisme d'algËbres. Nous voulons montrer que
$Bf$ est une Èquivalence faible, c'est-‡-dire que $\Omega Bf$ est un quasi-isomorphisme.
Gr‚ce au point {\it b} du lemme \ref{lemme1_cmf_cocog}, les flËches verticales
du diagramme commutatif
\[
\xymatrix{ A \ar[r]^f & A' \\
\Omega B A \ar[u] \ar[r] & \Omega BA' \ar[u]
}
\]
sont des quasi-isomorphismes. Par la propriÈtÈ de saturation des quasi-isomor\-phismes,
le morphisme $\Omega Bf$ est aussi un quasi-isomorphisme. \findem \\

%
%
%
%
%
%
%
%
%
%
%
%

{\em DÈmonstration du point {\em a} du thÈorËme \ref{theoreme_cmf_cocog}} :\\

(CM1) Les colimites de diagrammes finis de cogËbres sont donnÈes par les colimites des dia\-grammes
de complexes sous-jacents. Les constructions de produits et d'Ègalisateurs dans la catÈgorie
des cogËbres cocomplËtes sont duales de celles de coproduits et de co-Ègalisateur
dans la catÈgorie des algËbres, qui sont dÈcrites dans $\cite[3.3]{Munkholm78}$.\\

(CM2) Ceci est une consÈquence de la dÈfinition des Èquivalences faibles et
 de l'axiome (CM2) pour la structure
de catÈgorie de modËles sur $\alg$.\\

(CM3) Les cofibrations sont stables par rÈtract car elles sont les mono\-morphismes.
Les Èquiva\-lences faibles aussi car le foncteur $\Omega$ envoie un rÈtract sur un
rÈtract. Pour les fibrations, on rappelle qu'un morphisme $p$ est une fibration
s'il a la propriÈtÈ de relËvement ‡ droite par rapport aux cofibrations triviales.
On vÈrifie qu'un rÈtract d'un tel morphisme $p$ a la mÍme propriÈtÈ
de relËvement.\\

(CM4) AprËs (CM5).\\

(CM5) {\em factorisation : }\\
Soit $f : C \ra D$ un morphisme de $\cocog.$ Par l'axiome (CM5)  pour la structure
de catÈgorie de modËles sur $\alg$,  nous avons une factorisation de $\Omega f$
en
\[
\xymatrix{
  \Omega C \ar@{ >->}[rd]_i \ar[rr]^{f} & & \Omega D\\
 & A \ar@{->>}[ur]_p,}
\]
o˘ la cofibration $i$ (resp.~la fibration $p$) de $\alg$ est un quasi-isomorphisme.
Ainsi, le morphisme $B \Omega f:B \Omega C \ra B \Omega D$ se factorise
en $Bp \circ Bi .$
Soit le diagramme suivant dans $\cocog$
\[
\xymatrix{     & BA\prod_{B\Omega D}D \ar[dr]^q
\ar@{}[dddr]|{cart.}  \ar[dd]|\hole   &   \\ C 
\ar[rr]^(.7){f} \ar[dd]^{\in \ce q}  & & D
\ar[dd]^{\in \ce q}\\ &  BA \ar[dr]^{Bp} &  \\ 
B\Omega C \ar[ur]^{Bi} \ar[rr]_{B\Omega f}&  & B\Omega D. }
\]
Comme il est commutatif, le morphisme
$f : C \ra D$ et la composition 
\[
C \ra B \Omega C \arr{Bi} BA
\]
 dÈterminent
un morphisme  $ \tilde \imath : C \ra BA\prod_{B\Omega D}D $. Nous allons montrer
que
\[
\xymatrix{
 & BA\prod_{B\Omega D}D \ar[dr]^q
   &   \\ C  \ar[ur]^{\tilde \imath} \ar[rr]_{f} & & D
}
\]
fournit une factorisation du morphisme $f$ dans $\cocog$, o˘  $\tilde{\imath}$  est une
cofibration et $q$ est une fibration. Nous montrerons ensuite que
la cofibration $\tilde \imath$ (resp.~la fibration $q$) est triviale.

D'aprËs le point {\it b} du lemme \ref{lemme3_cmf_cocog}, le morphisme
$Bp$ est une fibration dans $\cocog$.
La projection $q : BA\prod_{B\Omega D}D \ra D $ est aussi une fibration car 
les fibrations sont stables par changement de base.
Admettons pour l'instant que nous savons que
$BA\prod_{B\Omega D}D \ra BA$ est cofibration (voir le
lemme \ref{lemmecle_cmf_cocog} ci-dessous). Le morphisme
$\tilde{\imath}$  est un monomorphisme (c'est-‡-dire une Èquivalence faible dans $\cocog$)
puisque la composition
\[
C \ra B \Omega C \arr{Bi} BA
\]
en est une.
Il reste ‡  montrer que la cofibration $\tilde{\imath}$ (resp.~la fibration~$q$)
 est une Èquivalence faible dans $\cocog$.  Admettons pour l'instant que nous savons que
$BA\prod_{B\Omega D}D \ra BA$ est Èquivalence faible (voir le
lemme \ref{lemmecle_cmf_cocog} ci-dessous).
Nous savons par le point {\it b} du lemme \ref{lemme3_cmf_cocog}
que le morphisme $Bi$ (resp.~$Bp$)  est une Èquivalence faible.
Comme le morphisme $C \ra B\Omega C$ (resp.~$D \ra B\Omega D$) est
une Èquivalence faible, $\tilde{\imath}$ (resp.~$q$) en est aussi une par
la propriÈtÈ de saturation de la classe des Èquivalences faibles de $\cocog$.\\

(CM4) {\em relËvement : }\\
{\it a.}
Soit le diagramme
commutatif dans $\cocog$
\[
\xymatrix{ E \ar@{ >->}[d]_u \ar@{->}[r] & C \ar@{->>}[d]^t \\
F\ar@{->}[r]  &  D }
\]
o˘ $t$ est une fibration triviale et $u$ une
cofibration.
Nous cherchons un morphisme $\alpha $ tel que les deux triangles du diagramme
\[
\xymatrix{ E \ar@{ >->}[d]_u \ar@{->}[r] & C \ar@{->>}[d]^t \\ F
\ar@{->}[r] \ar@{.>}[ur]^\alpha &  D}
\]
soient commutatifs.
En utilisant la construction de la dÈmonstration de (CM5), nous
 factorisons $t$ 
en $q\circ \tilde{\imath}$, o˘ le morphisme $q:BA \prod_{B\Omega D}D \ra D$
est une fibration  et o˘ le morphisme $\tilde{\imath}:C \ra BA \prod_{B\Omega D}D$ est
une cofibration. Par la propriÈtÈ de saturation de la classe $\weq$, les morphismes
$\tilde \imath$ et $q$  sont tous les deux des Èquivalences faibles.
Les fibrations Ètant les morphismes ayant la
propriÈtÈ de relËvement  ‡ droite par rapport aux cofibrations
triviales, il existe un relËvement $r : BA  
\prod_{B\Omega D}D \ra C$ dans
le diagramme de $\cocog$
\[
\xymatrix{ C \ar@{ >->}[d]_{\tilde \imath} \ar@{->}[r]^{\Id} &
C\ar@{->>}[d]^t \\
BA \prod_{B\Omega D}D \ar@{->}[r]_(.7)q &  D. }
\]
 Il nous suffit
donc de trouver un relËvement dans le diagramme
\[
\xymatrix{ E \ar@{ >->}[d]_u \ar@{->}[r] & BA \prod_{B\Omega D}D
\ar@{->>}[d]^q \\ F \ar@{->}[r]
\ar@{.>}[ur] &  D,}
\]
ou de maniËre Èquivalente dans le diagramme.
\[
\xymatrix{ E \ar@{ >->}[d]_u \ar@{->}[r] & BA \prod_{B\Omega D}D
\ar@{->>}[d] \ar@{->}[r] & BA \ar@{->>}[d]^{Bp}\\ F \ar@{->}[r]
\ar@{.>}[urr] &  D \ar@{}[ur]|{cart.}\ar@{->}[r] & B\Omega D.}
\]
Un tel relËvement existe gr‚ce ‡ l'adjonction entre $\Omega$ et $B$
et gr‚ce ‡ l'axiome de relËvement (CM4) de la structure de catÈgorie de modËles
fermÈe sur $\alg$. \\
%

%
%
%
%
%
%
%
%
%
%
%
%

{\noindent \bf Objets cofibrants et fibrants}\\

Tous les objets de $\cocog$ sont cofibrants puisque les
cofibrations sont les monomorphismes.

Montrons qu'un objet de $\cocog$
est fibrant si et seulement si il
est isomorphe, en tant que cogËbre graduÈe, ‡ une cogËbre tensorielle
rÈduite.

 Soit $C$ un objet fibrant de $\cocog$. Par l'axiome de relËvement (CM4), la
cofibration triviale $\psi : C \ra B \Omega C$ admet une rÈtraction $r$
dans $\cocog$. Notons $p_1 : B\Omega C \ra (B\Omega C)\prim 1$ la projection
canonique et posons $p_1^C = r \prim 1 \circ p_1 \circ \psi$. nous vÈrifions
facilement que le morphisme $p_1^C : C \ra C \prim 1$ est universel parmi
les morphismes d'objets graduÈs $C' \ra C\prim 1$, o˘ $C'$ est une cogËbre
graduÈe cocomplËte. Ainsi $p_1^C$ induit un isomorphisme de cogËbres
graduÈes
\[
C \arr{\sim} \ctr({C\prim 1}).
\]

La rÈciproque utilise les rÈsultats de la section \ref{section_cmf_aia}.
EnonÁons les deux rÈsultats de cette section qui nous seront utiles ici.\\
(\ref{theoreme_cmf_aia}) {\em La catÈgorie $\aia$ peut Ítre munie d'une structure de catÈgorie de modËles
dont la classe des Èquivalences faibles est exactement celle des $\ai$-quasi-isomorphismes et dont
la classe des cofibrations (resp.~des fibrations) est formÈe des morphismes $f : A \ra A'$, o˘ $A, A'$ sont des $\ai$-algËbres,
tel que $f_1$ est un monomorphisme (resp.~un Èpimorphisme)}.\\
(\ref{proposition1_cmf_aia}.~{\it a}) {\em Un morphisme $f$ est une Èquivalence faible de $\aia$ si et
seulement si sa construction bar $B f$ est une Èquivalence faible de $\cocog$}.\\
Notre dÈmonstration de (\ref{theoreme_cmf_aia}) est basÈe sur la thÈorie de l'obstruction
(voir \ref{section_obstruction}).
On peut donc interprÈter la rÈciproque que nous allons montrer comme une
consÈquence du fait que l'opÈrade des $\ai$-algËbres est le modËle minimal cofibrant
au sens de M.~Markl \cite{Markl96} de celle des algËbres associatives (voir l'introduction
‡ l'appendice \ref{section_obstruction}).\\

Supposons que $C$ est une cogËbre isomorphe,
en tant que cogËbre graduÈe,
‡ une cogËbre tensorielle rÈduite. Nous voulons montrer qu'elle est fibrante. On rappelle que
la sous-catÈgorie de $\cocog$ formÈe de telles cogËbres
est Èqui\-valente ‡ la catÈgorie $\aia$ des $\ai$-algËbres.
La cogËbre $B\Omega C$ appartient elle aussi ‡ cette sous-catÈgorie.
Le morphisme $C \ra B \Omega C$ est une Èquivalence faible de $\cocog$.
Par la proposition (\ref{proposition1_cmf_aia}.~{\it a}), il induit un
quasi-isomorphisme dans les primitifs.
L'axiome (CM4) du thÈorËme (\ref{theoreme_cmf_aia}) nous donne un relËvement dans le diagramme
\[
\xymatrix{
C \ar@{ >->}[d] \ar[r]^\Id &  C \ar[d]\\
B \Omega C \ar@{..>}[ur] \ar[r] & 0.
}
\]
La cogËbre $C$ est donc un rÈtract de $B \Omega C$. Comme la construction
bar conserve les fibrations et comme $\Omega C$ est
une algËbre fibrante, la cogËbre $B \Omega C$ est une cogËbre fibrante.
Le rÈtract d'une cogËbre fibrante
Ètant aussi fibrant, la cogËbre $C$ est fibrante.\\

{\em DÈmonstration du point {\em b} du thÈorËme \ref{theoreme_cmf_cocog}} : 

C'est un corollaire du lemme \ref{lemme1_cmf_cocog} qui nous dit
que les morphismes d'adjonction $C \ra B\Omega C$, o˘ $C$ est une cogËbre,
et $\Omega B A\ra A$,
o˘ $A$ est une algËbre, sont des Èquivalences faibles dans $\cocog$ et dans $\alg.$ \findem \\

Le lemme suivant complËte la dÈmonstration ci-dessus.

%
%
%
%
%
%
%
%
%
%
%
%

\begin{lemme} \label{lemmecle_cmf_cocog}
Soit $A$ une algËbre et $D$ une cogËbre. Soit une fibration $p : A \twoheadrightarrow \Omega D$ de $\alg$.
Le morphisme  $j : BA \prod_{B\Omega D}D \ra BA $ de cogËbres du diagramme cartÈsien 
\[
\xymatrix{ BA \prod_{B\Omega D}D
\ar[d]_{j} \ar@{->}[r] & D \ar[d]\\ 
  BA  \ar@{}[ur]|{cart.}\ar@{->}[r]_{Bp} & B\Omega D.}
\]
est une cofibration triviale de $\cocog$. 
\end{lemme}

\dem Nous allons donner des filtrations sur les cogËbres 
\[
BA \prod_{B\Omega D}D \quad \quad \mbox{ et} \quad \quad BA
\]
 telles
qu'elles soient des cogËbres filtrÈes admissibles et telles que $j$ soit un quasi-isomorphisme
filtrÈ. 

Soit la suite exacte de complexes
\[
0 \ra K \ra A \stackrel{p}{\twoheadrightarrow} \Omega D \ra 0.
\]
Comme l'algËbre $\Omega D$ est libre, nous avons un scindage
de $p$ dans la catÈgorie des algËbres graduÈes. La diffÈrentielle
de 
\[
A \arr{\sim } K \oplus \Omega D
\]
est alors donnÈe par une matrice
\[\left[
\begin{array}{cc}
d_K & d' \\ 0 & d_{\Omega D}
\end{array}
\right].
\]
Le scindage nous  donne des isomorphismes
de cogËbres graduÈes
\[
BA  \arr{\sim} BK \prod B\Omega D,
\]
\[
BA \prod_{B\Omega D}D \arr{\sim} BK \prod D.
\]
Munissons la cogËbre $B\Omega D$ de la filtration $D$-primitive.
Nous dÈfinissons des filtrations sur $BA$ et $BA \prod_{B\Omega D}D$
par les suites
\[
(BA)_j =  \sum_{p+q = j} (BK)\prim p \prod (B\Omega D)_q, \quad j \in \N,
\]
\[
(BA \prod_{B\Omega D}D)_j = \sum_{p+q = j} (BK)\prim p \prod D\prim q ,\quad j \in \N.
\]
Elles sont admissibles et respectent les diffÈrentielles des cogËbres $BA$ et $BA \prod_{B\Omega D}D$.
Pour ces filtrations, le morphisme $j$ est un morphisme filtrÈ. 
Soit $j \geq 1$.
En tant qu'objet graduÈ, le complexe
$\ogr(BA)$ est la somme des
\[
\diagnum{I}{
\ogr (B\Omega D)\ts K\tp{p_1}\ts \hdots  \ts
\ogr (B\Omega D)\ts K\tp{p_{k}}}, \quad k \geq 1.
\]
La diffÈrentielle
de $\ogr(BA)$ est construite ‡ partir des diffÈrentielles de $K$, $\ogr D$ et
du morphisme $d' : \Omega D \ra K.$
En tant qu'objet graduÈ, le complexe 
\[
\ogr (BA \prod_{B\Omega D} D) \arr{\sim} \ogr(BA) \prod \ogr D
\]
est la somme  des
\[\diagnum{II}{
\ogr D\ts K\tp{p_1}\ts \hdots  \ts
\ogr D\ts K\tp{p_{k}}}, \quad k \geq 1.
\]
La diffÈrentielle
de $\ogr (BA \prod_{B\Omega D} D) $
est construite ‡ partir des diffÈrentielles de $K$, $\ogr (B\Omega D)$ et
du morphisme $d' : \Omega D \ra K.$ Ainsi, la diffÈrentielle ``naÔve'' sur la
somme des termes (I), respectivement (II), est perturbÈe par la contribution
de $d' : \Omega D \ra K.$ Pour montrer que $j$ induit nÈanmoins un
quasi-isomorphisme entre les sommes, nous introduisons une filtration
supplÈmentaire telle que dans les objets graduÈs associÈs, la contribution de
$d' : \Omega D \ra K$ s'annule.
Soit la filtration $F_l\ogr (BA)$, $l \in \N$, de $\ogr (BA)$  induite par la suite
\[
(BA)_{l} = BK \prod (B\Omega D)\prim l, \quad l \in \N. 
\]
Soit la filtration $F_l \ogr (BA \prod_{B\Omega D}D)$, $l \in \N$, de $ \ogr (BA \prod_{B\Omega D}D)$
dont le \mbox{$l$-iËme} sous-objet, $l \in \N$, est la somme des objets de type (II)
comprenant un nombre de termes $\ogr D$ infÈrieur ou Ègal ‡ $l$.
Le morphisme 
\[
\ogr j : \ogr (BA \prod_{B\Omega D}D) \ra \ogr (BA)
\]
induit des morphismes
\[
F_l\ogr (BA \prod_{B\Omega D}D) \ra F_l \ogr (BA), \quad l \in \N.
\]
Il induit donc un morphisme
entre les objets graduÈs associÈs aux filtrations selon l'indice $l$. Ce dernier a
pour composantes les morphismes de complexes (avec les diffÈrentielles ``naÔves'')
\[
\xymatrix{
\ogr_{q_1}D\ts K\tp{p_1}\ts \hdots  \ts \ogr_{q_{k-1}}D \ts K\tp{p_{k}}
 \ar[d] \\
 \ogr_{q_1} (B\Omega D)\ts K\tp{p_1}\ts \hdots  \ts
\ogr_{q_{k-1}}(B\Omega D) \ts K\tp{p_{k}}}
\]
qui sont des quasi-isomorphismes (voir \ref{lemme1_cmf_cocog}). 
Le morphisme 
\[
\ogr j : \ogr(BA \prod_{B\Omega D}D) \arr{\sim} \ogr(BA)
\]
est donc un quasi-isomorphisme.
Nous venons ainsi de montrer
que $j$ est un quasi-isomorphisme filtrÈ de cogËbres admissibles.
Par le lemme \ref{lemme4_cmf_cocog}, le morphisme $j$ est une Èquivalence faible.
Il est une cofibration car il est clairement un monomorphisme. \findem

%
%
%
%
%
%
%
%
%
%
%
%
%
%
%
%

\subsection{$\aia$ comme ``catÈgorie de modËles sans limites''}
\label{section_cmf_aia}

Dans la catÈgorie $\aia$ des $\ai$-algËbres, nous considÈrons les trois
classes de morphismes suivantes~:

\english \begin{itemize}
\item[-] la classe $\weq$ est formÈe des {\em Èquivalences faibles}, c'est-‡-dire des
morphismes $f : A \ra A'$ tels que $f_1$ est un quasi-isomorphisme,
\item[-] la classe $\cof$ est formÈe des {\em cofibrations}, c'est-‡-dire des morphismes
 $f : A \ra A'$ tels que $f_1$  un monomorphisme,
\item[-] la classe $\fib$ est formÈe des {\em fibrations}, c'est-‡-dire des morphismes
 $f : A \ra A'$ tels que $f_1$  un Èpimorphisme.
\end{itemize} \francais

\begin{theoreme} \label{theoreme_cmf_aia}
La catÈgorie $\aia$, munie des trois classes dÈfinies ci-dessus, vÈrifie
l'axiome  (A) ci-dessous et les axiomes (CM2) -- (CM5) de la dÈfinition
\ref{definition_cmf}. Tous les objets sont fibrants et cofibrants. 
\end{theoreme}

\english \begin{itemize}
\item[(A)] Soit $ q : A \twoheadrightarrow A'$ une fibration et $f:A'' \ra A'$ un
morphisme. Il existe  un produit fibrÈ au-dessus de
\[
\xymatrix{A \ar@{->>}[r]^q & A' & A'' \ar[l]_f}.
\]
\end{itemize} \francais
L'axiome (A) est un affaiblissement de l'axiome (CM1) de la dÈfinition
\ref{definition_cmf}. 
Notre dÈmonstration de ce thÈorËme est entiËrement basÈe sur la thÈorie de l'obstruction
(\ref{section_obstruction}).

\begin{lemme} \label{lemme1_cmf_aia}
Soit $A$ une $\ai$-algËbre et $K$ un complexe considÈrÈ comme
$\ai$-algËbre (\ref{remarque_ai-structure_complexe}).
Supposons que le complexe $K$ est contractile.
Soit $g : (A,m^A_1)\ra (K,m^K_1)$ un morphisme de complexes.
Il existe un morphisme d'$\ai$-algËbres
\[
f : A \arr{} K
\]
tel que $f_1 = g$.
\end{lemme}

\dem 
Nous construisons  par rÈcurrence les morphismes
\[
f_i : A \tp i \ra K, \quad i \geq 1.
\]
 Soit $f_1 = g$.
Supposons que nous avons dÈj‡ construit des morphismes $f_i$, $1 \leq i \leq n$, 
qui dÈfinissent un ${\rm A}_n$-morphisme $A \ra K$. Nous cherchons un morphisme
$f_{n+1}$ dont le bord est le cycle $-r(f_1,\hdots ,f_n)$, i.~e.
\[
\del{f_{n+1}} + r(f_1,\hdots ,f_n) = 0 \quad \mbox{(voir 
\ref{extension_nstructure_morphismes})}.
\]
Comme $(K,m^K_1)$ est contractile, il existe bien un tel morphisme
$f_{n+1}$.
  \findem

\begin{lemme} \label{lemme2_cmf_aia}
\begin{enumerate}
\item[a.] Soit $j :A \ra D$ une cofibration de $\aia$. Il existe une $\ai$-algËbre $D'$
 et un isomorphisme d'$\ai$-algËbres
 $k:D \ra D'$ tels que
  la composition $k\circ j:A \ra D'$ est un morphisme strict.
\item[b.] Soit $q : C \ra E$ une fibration de $\aia$. Il existe  une $\ai$-algËbre $C'$
 et un isomorphisme $l : C' \ra C$ tels que la composition $q\circ l:C' \ra E$
 est un morphisme strict.
\end{enumerate}
\end{lemme}

\dem 

{\it a.} 
Nous construisons par rÈcurrence 
des morphismes
\[
k_i : D \tp i \ra D,\quad i\geq 1,
\]
homogËnes de degrÈ $1-i$
tels que $k \circ j$ est un morphisme strict.
Posons $k_1 = \Id_{D}$. Supposons que nous avons dÈj‡ construit des 
morphismes $k_i$, $1 \leq i \leq n$, tels que l'Èquation
\[
(eq_m) \quad \sum_{1 \leq l \leq m} \ \sum_{\sum i_r = m} (-1)^s k_l \circ (j_{i_1} \ts \hdots \ts j_{i_l}) = 0,
\quad 2 \leq m  \leq n,
 \]
o˘ $s$ est le signe apparaissant dans \ref{definition_ai-morphisme}, est vÈrifiÈe pour tout
$ 2 \leq m \leq n$.
Soit $r$ une rÈtraction dans $\gr \sf C$ de $j_1 : A \ra D$. 
Soit $k_{n+1} $ le morphisme  dÈfini par la somme
\[
- \left[\sum_{1 \leq l \leq n} \ \sum_{\sum i_r = n+1}
(-1)^s  k_l \circ (j_{i_1} \ts \hdots \ts j_{i_l})\right]\circ r\tp{n+1}.
\]
La suite $(k_1, \hdots , k_{n+1})$ vÈrifie l'Èquation $(eq_m)$ pour $2 \leq m \leq n+1$.
Comme $k_1$ est un isomorphisme d'objets graduÈs, les morphismes $k_i$, $i \geq 1$, induisent
un isomorphisme
\[
K : \ctr (SD) \arr{\sim} \ctr (SD).
\]
Nous dÈfinissons $D'$ comme l'$\ai$-algËbre dont l'objet graduÈ sous-jacent est $D$ et dont
les multiplications $m'_i$, $i\geq 1$, sont dÈfinies gr‚ce aux bijections
$m'_i \lra b'_i$ (voir \ref{construction_bar_cobar}), par les ÈgalitÈs
\[
b'_i = (K\circ b \circ K^{-1})_i,\quad i \geq 1.
\]
Alors le morphisme $ k: D \ra D'$ est clairement un isomorphisme de $\aia$ et la composÈe $k \circ j$
est stricte par construction de $k$.

 {\it b.} La dÈmonstration est similaire. Il faut utiliser une section de $q_1$ au lieu
 d'une rÈtraction de $j_1.$ \findem\\

{\em DÈmonstration du thÈorËme \ref{theoreme_cmf_aia} : } \\

(A) : Soit $ q : A \twoheadrightarrow A'$ une fibration et
$f :  A'' \ra A'$ un morphisme de $\aia$. La construction bar envoie
les morphismes $q$ et $f$ sur les morphismes $Q : BA \ra BA'$ et $F : BA'' \ra BA'$.
Nous allons montrer que le produit fibrÈ de $\cocog$ au-dessus
de 
\[
\xymatrix{BA \ar@{->>}[r]^Q & BA' & BA'' \ar[l]_F}
\]
est encore une cogËbre tensorielle rÈduite dans $\gr \sf C.$
Une section de $Q_1$ dans  $\gr \sf C$
induit un isomorphisme 
\[
SA \arr{\sim} SA'  \oplus K,
\]
o˘ $K$ est le noyau de $Q_1.$ 
Le produit fibrÈ $BA \prod_{BA'} BA''$ est
isomorphe, en tant que cogËbre graduÈe, ‡ 
\[
\ctr K \prod \ctr (SA'') \arr{\sim}\ctr (K \oplus SA'').
\]

(CM2) et (CM3) :  immÈdiat.\\

(CM4) {\em relËvement} :  Soit un diagramme d'$\ai$-algËbres
\[  \diagnum{I}{
\xymatrix{
A \ar[r]^f \ar@{ >->}[d]_j & C \ar@{->>}[d]^q \\
D \ar[r]_g & E,}}
\]
o˘ $q$ est une fibration et $j$ est une cofibration.
 Par le lemme
\ref{lemme2_cmf_aia}, quitte ‡ remplacer ce diagramme par
un diagramme isomorphe, 
nous pouvons supposer que les morphismes $j$ et $q$ sont stricts.
Supposons que la fibration $q$  (resp.~la cofibration $j$)
est triviale. Nous cherchons un relËvement
$\alpha$ qui rend commutatifs les deux triangles du diagramme
\[
\diagnum{I^+}{
\xymatrix{
A \ar[r]^f \ar@{ >->}[d]_j & C \ar@{->>}[d]^q \\
D \ar[r]_g \ar[ur]^\alpha & E.
}}
\]
Nous allons construire par rÈcurrence les morphismes correspondants
\[
\alpha_i : D\tp i \ra C, \quad i \geq 1.
\]
Gr‚ce au point {\it a} de l'axiome (CM4) pour la catÈgorie de modËles
$\cc \sf C$, il existe un relËvement $\alpha_1$ qui rend commutatifs les deux triangles
\[\diagnum{II}{
\xymatrix{
(A,m^A_1) \ar[r]^{f_1} \ar@{ >->}[d]_{j_1} & (C,m^C_1) \ar@{->>}[d]^{q_1} \\
(D,m^D_1) \ar@{..>}[ur]^{\alpha_1}\ar[r]_{g_1} & (E,m^E_1).}}
\] 
Supposons que nous avons dÈj‡ construit des morphismes $\alpha_i$, $1 \leq i \leq n$,
tels que le diagramme $(I^+)$ commute dans la catÈgorie des ${\rm A}_n$-algËbres.
Nous devons trouver un $\alpha_{n+1}$ tel que 
\[ 
\begin{array}{rll}
(1) & &\del{\alpha_{n+1}} + r({\alpha_{1}},\hdots ,{\alpha_{n}}) = 0, \quad
\mbox{(voir \ref{extension_nstructure_morphismes})}\\
(2) & & \alpha_{n+1} \cdot j_{1}\tp{n+1} = f_{n+1},  \\
(3) & &  q_{1}\cdot 
\alpha_{n+1} = g_{n+1}.
\end{array}
\]
Choisissons une solution $\beta$ de $(2)$ et $(3)$. Par exemple, 
si $\rho$ est une rÈtraction de $j_1$ et $\sigma$ une section de $q_1$
dans $\gr \sf C$, nous pouvons choisir
\[
\beta = f_{n+1}  \rho \tp{n+1} + \sigma  g_{n+1} - \sigma q_1 f_{n+1} \rho \tp{n+1}.
\]
Le morphisme $j$ est strict. Par le lemme \ref{fonctorialite_obstruction},
nous avons donc
\[
(\delta(\beta) + r(\alpha_1,\hdots ,\alpha_n))\circ j_1  =
\delta(\beta\circ j_1) + r(\alpha_1 \circ j_1,\hdots ,\alpha_n\circ j_1\tp n),
\]
et le terme de droite est Ègal ‡
\[
\delta(f_{n+1}) + r(f_1 \circ j_1,\hdots ,f_n) = 0.
\]
De mÍme, nous avons $q_1 \circ (\delta(\beta) + r(\alpha_1,\hdots ,\alpha_n)) = 0.$
Le cycle  
$\del{\beta} + r({\alpha_{1}},\hdots ,{\alpha_{n}})$ se
factorise donc en
\[
D\tp{n+1} \arr{p} \coker j_1\tp{n+1} \arr{c'}\ker q_1 \arr{i} C,
\]
o˘ $p$ est la projection canonique et $i$ l'injection canonique.
Comme $\ker q_1$ (resp.~$\coker (j_1\tp{n+1})$) est contractile, le cycle $c'$ est le bord
d'un morphisme $h'$. Le morphisme $\alpha_{n+1} = \beta - i \circ h' \circ p$
vÈrifie alors les Èquations $(1)$, $(2)$ et $(3).$

\begin{remarque}{\em \label{remarque_theoreme_cmf_aia}
La dÈmonstration de l'axiome de relËvement (CM4) montre que pour tout relËvement
$\alpha_1$ dans la catÈgorie $\cc \sf C$ du diagramme $(\mathrm{II})$,
on a un relËvement  $\alpha : D \ra C$ dans le diagramme $(\mathrm I)$.
}\end{remarque}
 
(CM5) {\em factorisation : } Soit $f : A \ra B$ un morphisme d'$\ai$-algËbres.

{\it a.} Soit $C = B \oplus S^{-1}B $ le cÙne de l'identitÈ de $S^{-1}B.$
ConsidÈrons le complexe $C$ comme une $\ai$-algËbre (voir \ref{remarque_ai-structure_complexe}).
Soit $j : A \ra A \prod C$ le morphisme d'$\ai$-algËbres de composantes $\Id_A$ et $0$.
Le morphisme  $q_1 : A \oplus C \ra B$ 
de  composantes le morphisme $f$ et la projection canonique $C \ra B$ 
est un relËvement du diagramme de $\cc \sf C$
\[
\xymatrix{
A \ar[r]^{f_1} \ar@{ >->}[d]^{j_1} & B \ar@{->>}[d] \\
A \oplus C \ar@{..>}[ur] \ar[r] & 0.
}
\]
La remarque \ref{remarque_theoreme_cmf_aia} appliquÈe au point {\it a} de l'axiome (CM4) nous
donne un relËvement dans le diagramme de $\aia$
\[
\xymatrix{
A \ar[r]^f \ar@{ >->}[d]^j & B \ar@{->>}[d] \\
A \prod C \ar@{..>}[ur]^q \ar[r] & 0.
}
\]
Dans la factorisation $f = q \circ j$,
$j$ est une cofibration triviale et $q$ est une fibration.

{\it b.} Soit $C= SA \oplus A$ le cÙne de l'identitÈ du complexe $(A,m_1)$.
ConsidÈrons $C$ comme une $\ai$-algËbre.
Par le lemme \ref{lemme1_cmf_aia},
il existe un morphisme d'$\ai$-algËbres $i : A \ra C$ tel que $i_1$ est l'injection
canonique $A \ra C.$
Soit $j : A \ra B \prod C$ le morphisme d'$\ai$-algËbres de composantes $f$ et $i.$ C'est une
cofibration triviale. Notons $q$ la projection canonique $B \prod C \ra B.$
C'est une fibration et le
morphisme $f$ se factorise en $q\circ j.$

\findem\\

{\noindent \bf Liens entre la ``catÈgories de modËles sans limites''
$\aia$ et la catÈgorie de modËles $\cocog$}\\

Soit $\cogtr$ la sous-catÈgorie de
$\cocog$ formÈe des cogËbres qui sont tensorielles rÈduites en tant que cogËbres graduÈes.
La construction bar induit un isomorphisme de catÈgories $\aia \ra \cogtr.$
Munissons $\cogtr$ de la structure de ``catÈgorie de modËles sans limites'' donnÈe par
cet isomorphisme. Les Èquivalences faibles (resp.~les cofibrations, resp.~les
fibrations) sont donc les morphismes $F : (\ctr V,b)\ra (\ctr V',b')$ qui induisent dans les primitifs
 un quasi-isomorphisme $F_1 : (V,b_1)\ra (V',b'_1)$ (resp.~un monomorphisme, resp.~un Èpimorphisme).

\begin{proposition} \label{proposition1_cmf_aia}
Soit $A$ et $A'$ deux $\ai$-algËbres.
\english \begin{itemize}
\item[a.]
Un morphisme  $f : BA \ra BA'$ est une Èquivalence faible
de $\cogtr$ si et seulement si c'est une Èquivalence faible de
$\cocog$.
\item[b.] Un morphisme $j : BA \ra BA'$ est une cofibration de $\cogtr$
 si et seulement si c'est une cofibration de
$\cocog$. 
\item[c.]  Un morphisme $q : BA \ra BA'$ est une fibration de $\cogtr$
 si et seulement si c'est une fibration de
$\cocog$.  
\end{itemize} \francais
\end{proposition}

CommenÁons par un lemme.

\begin{lemme} \label{lemme3_cmf_aia} 
Soit $A$ une $\ai$-algËbre.
Le morphisme $\phi : BA \ra  B \Omega BA$ est une Èquivalence faible de $\cogtr$.
\end{lemme}

\dem  Il s'agit de montrer
que le morphisme $\phi \prim 1$ est un quasi-isomorphisme ou, de faÁon Èquivalente,
que le morphisme
\[
S^{-1}\phi\prim 1 : (A,m_1) \ra \Omega B A 
\]
est un quasi-isomorphisme. Le morphisme $S^{-1}\phi\prim 1$ est l'injection canonique
de $A$ dans $\Omega BA.$ Munissons $\Omega BA$ de la filtration 
induite par la filtration
primitive de $BA.$ Tout comme ‡ la fin de la dÈmonstration du point {\it b}
du lemme \ref{lemme1_cmf_cocog}, nous montrons que 
\[
\ogr_0 (\Omega BA) = A \quad \mbox{et} \quad \ogr_i (\Omega BA) = 0 \quad \mbox{pour} \quad i\geq 1.
\]
\hspace*{1cm} \findem

{\em DÈmonstration de la proposition \ref{proposition1_cmf_aia}}

{\it a.}
Soit $f : BA \ra BA'$ une Èquivalence faible de $\cogtr$. Le morphisme
$f$ est clairement un quasi-isomorphisme filtrÈ pour la filtration primitive.
Il est donc une Èquivalence faible de $\cocog$.
Supposons que $f$ est une Èquivalence faible
de $\cocog$. Par dÈfinition des Èquivalences faibles de $\cocog$, le morphisme
$\Omega f$ est un quasi-isomorphisme et, par suite, le morphisme
$B\Omega f$ est une Èquivalence faible de $\cogtr.$
Par le lemme \ref{lemme3_cmf_aia},
les deux flËches horizontales du diagramme commutatif 
\[
\xymatrix{BA \ar[r] \ar[d]_f & B \Omega BA \ar[d]^{B\Omega f}\\
BA' \ar[r] & B\Omega B'A',
}
\]
sont des Èquivalences faibles de $\cogtr$, et $f$ est donc aussi une Èquivalence faible de $\cogtr.$

{\it b.} Comme les cofibrations de $\cocog$ sont les monomorphismes,
une cofibration de $\cogtr$ est  une
cofibration de $\cocog$. RÈciproquement,
si $j : BA \ra BA'$ est une cofibration de $\cocog$, sa restriction aux primitifs $(BA) \prim 1
= SA$ est un monomorphisme. Comme nous avons $f ((BA) \prim 1) \subset (BA') \prim 1$, le morphisme
$j \prim 1 : SA \ra SA'$ est un monomorphisme et $j$ est donc une cofibration de
$\cogtr$.

{\it c.}
On rappelle que  les fibrations d'une catÈgorie de modËles
sont les morphismes ayant la propriÈtÈ de relËvement
par rapport aux cofibrations triviales. Ce fait rÈsulte
des axiomes (CM5) et (CM3) et vaut donc aussi pour
$\cogtr.$  Par les points {\it a} et {\it b},
une fibration de $\cocog$ est une fibration
de $\cogtr$. Supposons que $q$ est fibration de $\cogtr.$ 
Soit le diagramme de $\cocog$
\[
\xymatrix{C \ar@{ >->}[d]_j \ar[r]^f & BA \ar[d]^q\\
C' \ar[r]_g & BA',
}
\] 
o˘ $j$ est une cofibration triviale de $\cocog.$ Nous cherchons un relËvement de
$g$ relatif ‡ $f$. Dans le diagramme de $\cocog$ ci-dessous
\[
\xymatrix{C \ar@{ >->}[dd]_j \ar[rr]^f \ar@{ >->}[rd]_{\phi} & & BA \ar[dd]^q\\
& B\Omega C  \ar@{ >->}[dd]& \\
C' \ar[rr]^(.7)g|(.475){\ }  \ar@{ >->}[rd] & & BA',\\
& B\Omega C'  & 
}
\]
$\phi$ est une cofibration triviale de $\cocog$ et $BA$ est fibrant dans $\cocog$.
Nous avons donc une factorisation de $f$ en $f' \circ \phi$ pour un morphisme
$f' : B \Omega C \ra BA$. Comme $\Omega j$ est un monomorphisme et un quasi-isomorphisme,
le morphisme $B\Omega j$ est une cofibration  triviale de $\cogtr$. L'objet $BA'$
Ètant cofibrant dans $\cogtr$, le morphisme $q \circ f'$ se factorise en $g' \circ B \Omega f$
pour un morphisme $g' : B\Omega C' \ra BA'$.
Il suffit donc de trouver un relËvement de $g'$ relatif ‡ $f'.$
Par l'axiome (CM4) pour la catÈgorie $\cogtr$, il en existe un. \findem

%
%
%
%
%
%
%
%
%
%
%

\subsection{Homotopie au sens classique} \label{section_homotopie}
\index{homotopie}

Soit $C$ et $C'$ deux cogËbres cocomplËtes.
Soit $f$ et $g$ deux morphismes de cogËbres $C \ra C'$.
Ils sont {\em homotopes au sens
classique} 
s'il existe une $(f,g)$-codÈrivation $h : C \ra C'$ de degrÈ
$-1$ tel que $\del h = f - g$. 
 Nous comparons cette notion ‡ la notion
d'homotopie au sens des catÈgories de modËles (voir l'appendice \ref{section_rappels_cmf}). \\

\begin{proposition} \label{proposition1_homotopie}
Soit $C$ et $C'$ deux cogËbres cocomplËtes et $f,g$ deux morphismes $C \ra C'.$
\english \begin{itemize}
\item[a.] Si $f$ et $g$ sont homotopes au sens classique, ils sont homotopes
‡ gauche (voir la dÈfinition \ref{definition_homotopie_cmf}).
\item[b.] Si la cogËbre $C'$ est fibrante, alors $f$ et $g$ sont homotopes au sens
classique si et seulement si ils sont homotopes ‡ gauche.
\end{itemize} \francais
\end{proposition}

\dem \\
{\it a.}
Nous allons construire un cylindre $C \wedge I$ pour la cogËbre
$C$, puis nous allons montrer que la notion
d'homotopie classique est Èquivalente ‡ la notion de $C \wedge I$-homotopie ‡ gauche.

Nous notons $I$ le complexe dont la composante de degrÈ $0$ est
$e \oplus e$, la composante de degrÈ $-1$ est $e$, toutes les autres
composantes sont nulles. Nous notons $e_0$ et $e_1$ les composantes
de $I_0$. La diffÈrentielle $d : I \ra I$ est donnÈe par
\[
d_1 = \left[\begin{matrix}\Id \\ -\Id
\end{matrix}\right] : e\ra e_0\oplus e_1  .
\] 
Soit $\Delta :  I \ra I \ts I$  le morphisme dont les composantes 
non nulles sont les morphismes
\[
e_0 \arr{\sim} e_0 \ts e_0,\quad e_1 \arr{\sim} e_1 \ts e_1,
\quad e \arr{\sim} e_0\ts e, \quad e \arr{\sim} e\ts e_1
\]
donnÈs par la
contrainte d'unitaritÈ de la catÈgorie monoÔdale de base (\ref{categorie_de_base}).
Ceci dÈfinit sur $I$ une structure de
cogËbre co-associative diffÈrentielle graduÈe.

Soit $C$ une cogËbre cocomplËte.
Le produit tensoriel $C\ts I$ de complexe hÈrite naturellement
d'une structure de cogËbre diffÈrentielle graduÈe par la
comultiplication $C\ts I \ra (C\ts I)\ts C\ts I \iso C\ts C\ts I\ts
I$. Elle est cocomplËte. Nous notons $C_0$ et $C_1$ les composantes de $C \coprod C.$
Nous dÈfinissons le cylindre $C \wedge I = C \ts I$ pour $C$ par les deux morphismes
de cogËbres diffÈrentielles graduÈes $i$ et $p$
\[
C_0\coprod C_1 \arr{i} C\ts I \arr{p} C,
\]
o˘ le morphisme $i$ a pour composantes non nulles
\[
C_0 \arr{\sim} C \ts e_0, \quad C_1 \arr{\sim} C \ts e_1,
\]
et o˘
le morphisme $p$ a pour composantes non nulles
\[
C \ts e_0 \arr{\sim} C, \quad C \ts e_1 \arr{\sim} C,
\]
donnÈes par les contraintes d'unitaritÈ de la catÈgorie de base.
Le morphisme $i$ est une
cofibration et le morphisme $p$ une Èquivalence faible.\\
Soit $C'$ une cogËbre cocomplËte. Soit $f , g$ et $h$ trois morphismes 
graduÈs $C \ra C'$ respectivement de degrÈ zÈro, zÈro et -1.\\
Soit le morphisme graduÈ de degrÈ nul $H :C\ts I \ra B$ dont les composantes sont les trois
morphismes graduÈs
\[
 C \ts e_0 \iso C \arr{f} C', \quad C \ts e_1 \iso C \arr{g} C' 
\]
\[
\mbox{et} \quad C \ts e \iso C \arr{h} C'.
\]
Le morphisme $H : C\ts I \ra C'$ est un morphisme de cogËbres si et
seulement~si
\english \begin{itemize}
\item[-] les morphismes $f$ et $g$ sont des morphismes de cogËbres $C \ra C'$,
\item[-] le morphisme $h : C \ra C'$  est une $(f,g)$-codÈrivation.
\end{itemize} \francais
Il est compatible aux diffÈrentielles si et seulement si
\english \begin{itemize}
\item[-] les morphismes $f$ et $g$ sont des morphismes de complexes $C \ra C'$
\item[-] le morphisme $h : C \ra C'$
rÈalise une homotopie entre les morphismes de complexes $f$ et $g$.
\end{itemize} \francais
Pour finir, nous vÈrifions que le morphisme $H$ est bien une $C \wedge I$-homotopie entre
$f$ et $g.$\\

{\it b.}
Soit $C'$ une cogËbre fibrante.
Soit $f$ et $g : C \ra C'$ deux morphismes  homotopes au sens
des catÈgories de modËles. Notons toujours $C \wedge I$ le cylindre
construit ci-dessus. Par le lemme \ref{lemme1_cmf}, il existe
une $C \wedge I$-homotopie ‡ gauche $H : C \wedge I \ra C'$ entre
$f$ et $g$. Par la dÈmonstration du point {\it a},
il existe une homotopie $h : C \ra C'$
au sens classique entre $f$ et $g.$ \findem

%
%
%
%
%
%
%
%
%
%
%
%
%
%
%
%

%
%
%
%
%
%
%
%
%
%
%
%
%
%
%
%

\subsection{Equivalences faibles et quasi-isomorphismes} \label{section_equiv_qis}

Nous notons $Qis$ la classe des quasi-isomorphismes de $\cocog$ \indexnotation{Qis} et
$\mbox{\it Qisf\,}$ \indexnotation{Qisf} la classe des morphismes $f : C \ra D$ de $\cocog$ tels que
$C$ et $D$ admettent des filtrations admissibles pour lesquelles
$f$ est un quasi-isomorphisme filtrÈ.

Cette section est consacrÈe ‡ la comparaison des trois classes $\weq$, $Qis$ et $\mbox{\it Qisf\,}$.
Nous allons montrer en particulier les inclusions  suivantes 
\[
\mbox{\it Qisf\,} \subseteq \weq \subset Qis.
\]

\begin{proposition} \label{proposition_equiv_qis}
\english \begin{itemize}
\item[a.] Nous avons l'inclusion $\mbox{\it Qisf\,} \subseteq \weq$.
De l'autre cÙtÈ, le foncteur canonique
\[
\cocog[\mbox{\it Qisf\,}^{-1}] \arr{}\cocog[\weq^{-1}] = \Ho \cocog
\]
est une Èquivalence.
\item[b.] Les Èquivalences faibles de $\cocog$ sont des quasi-isomorphismes.
\item[c.] La classe $\weq$ est strictement incluse dans la classe $Qis.$
\item[d.] Soit $C$ et $D$ deux objets de $\cocog$ concentrÈs en degrÈs $< -1.$
Tout quasi-isomorphisme de cogËbres $C \ra D$ est une Èquivalence faible.
\item[e.] Soit $C$ et $D$ deux objets de $\cocog$ concentrÈs en degrÈs $\geq 0.$
Tout quasi-isomorphisme de cogËbres $C \ra D$ est une Èquivalence faible.
\end{itemize} \francais
\end{proposition}

\dem 

{\it a.}
On rappelle (\ref{lemme4_cmf_cocog}) qu'un quasi-isomorphisme filtrÈ de cogËbres
est une Èquivalence faible de $\cocog.$ Il nous faut donc montrer que  les Èquivalences faibles
deviennent des isomorphismes dans la catÈgorie localisÈe $\cocog[\mbox{\it Qisf\,}^{-1}]$.
Soit $f : C \ra C'$ une Èquivalence faible de $\cocog.$ Le morphisme
\[
\Omega f : \Omega C \ra \Omega C'
\]
est donc un quasi-isomorphisme d'algËbre. Par le lemme
\ref{lemme1_cmf_cocog}, le morphisme $B \Omega f : B\Omega C \ra B \Omega C'$ est
un quasi-isomorphisme filtrÈ. On rappelle (\ref{lemme1_cmf_cocog}) que les morphismes d'adjonction
$C\ra B \Omega C$ et $D \ra B \Omega D$ sont des quasi-isomorphismes filtrÈs.
Nous dÈduisons du diagramme commutatif de $\cocog$
\[
\xymatrix{
C \ar[r] \ar[d]_f & B \Omega C \ar[d]^{B\Omega f} \\
C'\ar[r]  & B \Omega C' }
\]
que le morphisme $f$
devient un isomorphisme dans la catÈgorie 
$\cocog [\mbox{\it Qisf\,}^{-1}].$

{\it b.} Les quasi-isomorphismes filtrÈs sont des quasi-isomorphismes.
La propriÈtÈ de saturation de la classe $Qis$, appliquÈ au diagramme ci-dessus
montre qu'une Èquivalence faible
est un quasi-isomorphisme.

{\it c.}
Nous allons construire un exemple de cogËbre qui est acyclique mais qui
n'est pas faiblement Èquivalente ‡ la cogËbre nulle. 

Soit $A$ une algËbre unitaire non nulle de la catÈgorie de base $\sf C$.
ConsidÈrons $A$ comme une algËbre associative (dont on oublie l'unitÈ), c'est-‡-dire
comme un objet de la catÈgorie $\alg$ du thÈorËme (\ref{theoreme_cmf_cocog})
Comme $A$ n'est pas quasi-isomorphe ‡ l'algËbre nulle,
la cogËbre $BA = (\ctr SA,b)$ n'est pas faiblement Èquivalente ‡ la cogËbre nulle (\ref{theoreme_cmf_cocog}, {\it b}).
Or, elle est bien quasi-isomorphe ‡ la cogËbre nulle :
en effet, le complexe sous-jacent ‡ $S^{-1}BA$ est le complexe
\[
\cdots \ra A \ts A \ts A \ra A \ts A \ra A \ra 0,
\]
qui est isomorphe ‡ la rÈsolution bar de l'algËbre $A$.
Ce complexe est acyclique car $A$ est unitaire (voir \cite[IX.6]{Cartan56} o˘ ce complexe
se nomme ``rÈsolution standard'').

{\it d.} Soit $C$ et $D$ deux cogËbres cocomplËtes concentrÈes en degrÈs $<-1$.
Nous allons montrer que le morphisme $\Omega f : \Omega C \ra \Omega D$ est un
quasi-isomorphisme de $\alg.$ Munissons $\Omega C$ (resp.~$\Omega D$) de la filtration
dÈcroissante donnÈe par
\[
F_l \Omega C = \bigoplus_{p \geq l} (S^{-1}C)\tp p  \quad \left(\mbox{resp.~}
F_l \Omega D = \bigoplus_{p \geq l} (S^{-1}D)\tp p \right), \quad l \in \N.
\]
Par notre hypothËse, le morphisme $\Omega f$ induit des quasi-isomorphismes
dans les sous-quotients de ces filtrations. Il en rÈsulte que, pour tout $n \in \N$,
il induit un isomorphisme en $H^{-n}$, car nous avons
\[
(F_l  \Omega C)^n = (F_l  \Omega D)^n = 0 \quad \mbox{pour} \quad l > n,
\]
d'aprËs l'hypothËse sur $C$ et $D$.

{\it e.} La dÈmonstration est la mÍme que pour le point {\it d}. Il suffit de remarquer
que le complexe $S^{-1}C$ est concentrÈ en
degrÈs $>0$ (au lieu de $<0$). \hspace*{1cm}\findem

%
%
%
%
%
%
%
%
%
%
%
%
%
%
%
%

\section{Transfert de structures le long d'Èquivalences d'homo\-topie} \label{section_transfert_structures}

Le but de cette section est de (re)montrer le {\em thÈorËme du modËle minimal} 
(\ref{corollaire_modele_minimal}).

\subsection{ModËle minimal}

\begin{theoreme} \label{theoreme_transfert_structures}
Soit $A$  une $\ai$-algËbre. 
Soit une Èquivalence d'homotopie dans $\cc \sf C$
\[
g : (V,d) \ra (A,m^A_1),
\]
o˘ $(V,d)$ est un complexe.
Il existe une structure d'$\ai$-algËbre sur $V$ telle que $m^V_1 = d$
et un morphisme d'$\ai$-algËbres 
\[
f: V \ra A
\]
telle que $f_1 = g$. 
\end{theoreme}

Ce rÈsultat est connu depuis les annÈes 70 dans le cas d'une $\ai$-algËbre connexe
(i.~e.~concentrÈes en degrÈs homologiques $\geq 1$) et d'un complexe $(V,d)$
o˘ la diffÈrentielle $d$ est nulle
($V$ est alors isomorphe ‡ $H^*A$).
Il y a deux mÈthodes pour montrer ce thÈorËme, celle utilisant la ``mÈthode des obstructions''
 \cite{Chen77a}, \cite{Chen77b}, \cite{Kadeishvili80}, \cite{Smirnov80}, \cite{Gugenheim82}
et celle utilisant ``l'astuce du tenseur'' \cite{Huebschmann86}, \cite{Gugenheim86},  \cite{Gugenheim89},
\cite{Gugenheim91}, \cite{Huebschmann91a}, \cite{Merkulov99}, \cite{Kontsevich01}.
L'article \cite{Johansson01} prÈsente l'unification de ces diffÈrentes mÈthodes. Nous donnons ici
une dÈmonstration utilisant les obstructions. \\

\dem Par l'axiome (CM5) pour la  catÈgorie de modËles  $\cc \sf C$,
le morphisme $g$ se factorise en $q\circ j$, o˘ $q$ est une fibration triviale et o˘ $j$
est une cofibration triviale. Il suffit donc de montrer
le thÈorËme dans le cas o˘ l'Èquivalence d'homotopie
est un Èpimorphisme et dans le cas o˘ elle est un monomorphisme.

Supposons que $g$ est une fibration triviale de $\cc \sf C.$ Soit $K$ le noyau de $g.$ 
Comme $K$ est contractile, nous
pouvons scinder $g$ dans la catÈgorie des complexes. Ce scindage induit un isomorphisme de complexes
\[
V \arr{\sim} K \oplus A
\]
par lequel  le morphisme $g$ s'identifie ‡ la projection $K \oplus A \ra A.$
ConsidÈrons $K$ comme une $\ai$-algËbre (voir \ref{remarque_ai-structure_complexe}).
Munissons l'objet graduÈ sous-jacent ‡ $V$ de la structure d'$\ai$-algËbre de $K \prod A.$
Le morphisme $f$ est le morphisme canonique $K \prod A \ra A$ de $\aia.$

Supposons que $g$ est une cofibration triviale de $\cc \sf C.$ Soit $K$ le conoyau de $g.$
Comme il est contractile, on
peut scinder $g$ dans la catÈgorie des complexes. Ce scindage induit un isomorphisme de $\cc \sf C$
\[
A \arr{\sim} K \oplus V
\]
par lequel le morphisme $g$ s'identifie ‡ l'injection $V \ra K \oplus V.$
ConsidÈrons $K$ comme une $\ai$-algËbre.
Par le lemme \ref{lemme1_cmf_aia}, il existe un morphisme d'$\ai$-algËbres $h : A \ra K$ tel que
$h_1$ est la projection $K \oplus V \ra K$ dans $\cc \sf C.$ Gr‚ce ‡ l'axiome (A) du thÈorËme
\ref{theoreme_cmf_aia}, morphisme $h$ admet un noyau dans la catÈgorie $\aia$. L'objet
graduÈ sous-jacent ‡ $\ker h$ est $V$. Nous avons ainsi muni $V$ d'une structure d'$\ai$-algËbre
telle que $m^V_1$ est la diffÈrentielle de $V.$ Le morphisme canonique
$V \ra A$ est tel que $f_1 = g$. \findem \\

\newpage

{\noindent \bf ModËle minimal}

\begin{definition}\label{definition_ai-algebre_minimale}{\em
Une $\ai$-algËbre est {\em minimale}
si $m_1 = 0$. Soit $A$ une $\ai$-algËbre.
Un {\em modËle minimal pour $A$} \index{modele minimal@{modËle minimal}}
est un $\ai$-quasi-isomorphisme d'$\ai$-algËbres
$A' \ra A$
o˘ $A'$ est minimale.}
\end{definition}

\begin{remarque}{\em
Cette utilisation du terme ``modËle minimal'', due ‡ M.~Kontsevich, est diffÈrente
de l'usage conventionnel en homotopie rationnelle (modËle minimal de Sullivan).
On peut la justifier par le fait que la construction bar
$BA'$ est un modËle minimal au sens de H.~J.~Baues et J.-M.~Lemaire \cite{Baues77} de
la cogËbre $BA$. Remarquons qu'un modËle minimal de $BA$ ne donne pas en
gÈnÈral un modËle minimal de $A$ : soit $(\ctr SV,b)$ une cogËbre 
tensorielle rÈduite sur $SV$ dont la diffÈrentielle $b$ induit
zÈro dans les $1$-primitifs ; si $(\ctr SV,b)$ est un modËle minimal 
de $BA$, c'est-‡-dire si on a un {\em quasi-isomorphisme} de cogËbres
\[
F : (\ctr SV,b) \ra BA,
\]
l'$\ai$-algËbre $V$ telle que $BV = (\ctr SV,b)$
n'est pas en gÈnÈral un modËle minimal pour l'$\ai$-algËbre $A.$
Cependant, si $F$ est une {\em Èquivalence faible} de $\cocog$,
$V$ est un modËle minimal de l'$\ai$-algËbre $A$.
}\end{remarque}

\begin{corollaire} \label{corollaire_modele_minimal}
Soit $A$ une $\ai$-algËbre. Il existe une structure d'$\ai$-algËbre sur son homologie
$H^*A$ telle que
\english \begin{itemize}
\item[a.] $m_1 = 0$ et $m_2$ est induite par $m^A_2$,
\item[b.] il existe un morphisme d'$\ai$-algËbres $H^*A \ra A$ relevant l'identitÈ de
$H^*A$.
\end{itemize} \francais
Cette structure est unique ‡ un isomorphisme (non unique) prËs.
\end{corollaire}

\dem Comme la catÈgorie de base $\sf C$ est semi-simple,
nous avons un isomorphisme dans la catÈgorie des complexes
\[
(A,m^A_1) \arr{\sim} H^*A \oplus K
\]
pour un complexe contractile $K.$ Le rÈsultat est dÈduit du thÈorËme
\ref{theoreme_transfert_structures} appliquÈ ‡ l'injection canonique
\[
g : H^* A \arr{} A.
\]
L'unicitÈ de la structure provient du fait qu'un morphisme $f$ entre $\ai$-algËbres
minimales est un quasi-isomorphisme si et seulement si $f_1$ est un isomorphisme si et
seulement si $f$ est un isomorphisme.
\findem\\

\subsection{Lien avec le lemme de perturbation}
\index{lemme de perturbation}

Une {\em perturbation} \index{perturbation} 
$\delta$ de la diffÈrentielle $d$ d'un complexe filtrÈ $W$
est un morphisme graduÈ $\delta : W \ra W$
de degrÈ $+1$ qui diminue la filtration et
tel que $d + \delta$ est encore une diffÈrentielle, c'est-‡-dire, tel que
\[
d \circ \delta + \delta \circ d + \delta^2 = 0.
\]
Une {\em contraction} \index{contraction}
\cite{Eilenberg53} (voir aussi \cite{Huebschmann91a} et les
rÈfÈrences donnÈes dans \cite{Huebschmann91a})
\[
\big(\xymatrix{
V \ar@<-1ex>[r]_{i} & W \ar@<-1ex>[l]_{\rho}
}, H \big)
\]
est donnÈe par deux complexes $V$ et $W$, deux morphismes de complexes
$i : V \ra W$ et $\rho : V \ra W$ et un morphisme graduÈ $H : W \ra W$
de degrÈ $-1$
tels que
\[
\rho \circ i = \Id_V, \quad i \circ \rho = \Id_W + \delta (H), \quad H \circ i = 0, \quad
\rho \circ H = 0 \quad   \mbox{et} \quad H^2 = 0.
\]
On dit aussi que $W$ {\em se contracte sur} $V.$
Si  les complexes sont filtrÈs,
la contraction est {\em filtrÈe} si les morphismes sont filtrÈs relativement
‡ ces filtrations.

Soit $V$ et $W$ des complexes munis de filtrations exhaustives et soit
\[
\big(\xymatrix{
(V,d_V) \ar@<-1ex>[r]_{i} & (W,d_W) \ar@<-1ex>[l]_{\rho}
}, H \big)
\]
une contraction filtrÈe et $\delta$ une perturbation de la diffÈrentielle $d_W.$
Le lemme de perturbation (\cite{Gugenheim72}, \cite{Huebschmann91a})
donne une nouvelle diffÈrentielle $d_V^\delta$ de $V$ et des morphismes $i^\delta$, $\rho^\delta$ et
$H^\delta$ tels que
\[
\big(\xymatrix{
(V,d_V^\delta) \ar@<-1ex>[r]_(.4){i^\delta} & (W,d_W + \delta) \ar@<-1ex>[l]_(.6){\rho^\delta}
}, H^\delta \big)
\]
est une contraction filtrÈe. Supposons que la contraction filtrÈe ci-dessus est une
{\em contraction filtrÈe de cogËbres} : les objets
$V$ et $W$ sont des cogËbres diffÈrentielles graduÈes filtrÈes,
les morphismes $i$ et $\rho$ sont des morphismes de cogËbres filtrÈes,
$H$ est une $\Id$-$(i \rho)$-codÈrivation filtrÈe de $W$. Supposons aussi
que la perturbation $\delta$ est une perturbation d'une diffÈrentielle de cogËbres,
i.~e.~$\delta$ est une $\Id$-$\Id$-codÈrivation de $W$.
Le lemme de perturbation produit alors une contraction de cogËbres
(\cite{Huebschmann91a},
\cite{Gugenheim86},  \cite{Gugenheim89}, \cite{Gugenheim91}, \cite{Merkulov99}).\\

Soit $A$ une $\ai$-algËbre et soit
\[
\xymatrix{
0 \ar[r] & (V,d_V) \ar@<-1ex>[r]_{i} & (A,m_1) \ar@<-1ex>[r]_{p} \ar@<-1ex>[l]_{\rho}
& (K,d_K) \ar@<-1ex>[l]_{\sigma} \ar[r] & 0
}
\]
 une suite exacte scindÈe de complexes telle que
\[
\rho \circ \sigma = 0 \quad \mbox{et} \quad i \circ \rho + \sigma \circ p = \Id_A.
\]
Soit $h$ une homotopie contractante de $K$ telle que $h^2 = 0$.
A partir de ces donnÈes, nous avons deux maniËres naturelles
de dÈfinir une structure d'$\ai$-algËbre sur $V$ et un $\ai$-morphisme
\[
V \ra A
\]
dont le premiËre composante est $i.$\\

{\em PremiËre mÈthode : le lemme de perturbation}\\
Nous appliquons le lemme de perturbation ‡ la contraction filtrÈe et ‡ la perturbation
de cogËbres
\[
\xymatrix{
\Big(
\ctr S(V,d_V) \ar@<-1ex>[r]_{F} & \ctr S(A,m_1) \ar@<-1ex>[l]_{R}
}, H\Big) \quad \mbox{et} \quad \delta : \ctr SA \ra \ctr SA,
\]
o˘ $F = \ctr Si$, $R =\ctr S\rho$,
$H$ est l'unique $\Id$-$(FR)$-codÈrivation relevant $\sigma \circ h \circ p$  et 
$\delta = b - b_1$ (ici $b$ est la diffÈrentielle de $BA$).
Nous obtenons une nouvelle diffÈrentielle $b'$ sur $\ctr SV$ et un morphisme de cogËbres
\[
F^\delta : (\ctr SV,b') \ra (\ctr SA,b).
\]
Nous obtenons une structure d'$\ai$-algËbre sur $V$
(notons cette $\ai$-algËbre $V^\delta$) et un $\ai$-morphisme
\[
f^\delta : V^\delta \ra A.
\]

{\em DeuxiËme mÈthode : le noyau d'un $\ai$-morphisme $g$}\\
DÈfinissons par rÈcurrence des morphismes
\[
g_i :  A \tp i \ra K, \quad i \geq 1,
\]
en posant
\[
g_1 = p \quad \mbox{et} \quad g_i = - h \circ r(g_1,\hdots,g_{i-1}) , \quad i \geq 2,
\]
o˘ $r(g_1,\hdots,g_{i-1})$ est le cycle du lemme (\ref{extension_nstructure_morphismes}).
Le lemme (\ref{extension_nstructure_morphismes}) montre qu'ils dÈfinissent
un $\ai$-morphisme $g : A \ra K$ (o˘ $K$ est le complexe $K$
considÈrÈ comme $\ai$-algËbre). L'axiome (A) du thÈorËme (\ref{theoreme_cmf_aia}) montre
qu'il existe un noyau de $g$ dans la catÈgorie $\aia$
\[
V^g = \ker g \ra A.
\]
Comme l'objet graduÈ sous-jacent de l'$\ai$-algËbre $V^g$ est $V$,
cela dÈfinit une $\ai$-structure sur $V$ et un $\ai$-morphisme
\[
f^g : V^g \ra A.
\]

\begin{lemme}{\em  \label{lemme_modele_minimal_I}
Nous avons un isomorphisme $\theta : V^\delta \ra V^g$ tel que
$\theta_1 = \Id$ et  $f^\delta  = f^g \circ \theta$.
}\end{lemme}

\dem 
Rappelons les descriptions de l'$\ai$-structure de $V^\delta$
et de $f^\delta$ en terme d'arbres due ‡ M.~Kontsevich et Y.~Soibelman
\cite[6.4]{Kontsevich01}. 

L'$\ai$-structure de $V^\delta$ est dÈfinie par les formules suivantes :
\[
m^\delta_1 = 0, \quad m^\delta_2 = \rho \circ m_2 \circ (i \ts i), \quad m^\delta_i = \sum_{T \in \mathcal T}
(-1)^s m_{i,T}, \quad i \geq 3,
\]
o˘ $s$ et $T$, $\mathcal T$ et $m_{i,T}$ sont dÈfinis ainsi :
ConsidÈrons l'ensemble $\mathcal T$ des arbres planaires orientÈs $T$ avec $i+1$ sommets terminaux
(la racine et les feuilles),
tels que l'aritÈ $|v|$ de tout sommet interne $v \in T$ (i.~e.~le nombre de flËches arrivant
‡ $v$) est $\geq 2.$
Pour dÈcrire le morphisme
\[
m_{i,T} : (V^\delta) \tp i \ra V^\delta, \quad i \geq 3, \quad T \in \mathcal T,
\]
on a besoin de considÈrer l'arbre $\b T$ construit ‡ partir de $T$ en
rajoutant un sommet interne au milieu de chaque arÍte interne. L'arbre $\b T$
est ainsi constituÈ de deux types de sommets internes : 
les {\em anciens} qui correspondent aux sommets internes de $T$
et les {\em nouveaux} que l'on vient de rajouter.
On colorie les sommets de $\b T$ par les morphismes suivants :
\english\begin{itemize}
\item[-] $\rho$ sur la racine,
\item[-] $i$ sur les feuilles,
\item[-] $m_{|v|}$ sur les sommets internes anciens $v$ (dont l'aritÈ est $|v|$),
\item[-] $H$ sur les sommets internes nouveaux.
\end{itemize}\francais
A chaque arbre $\b T$ ainsi coloriÈ, on associe le morphisme $m_{i,T}$
qui consiste ‡ composer les coloriages en descendant le long de l'arbre des feuilles vers
la racine. Voici un exemple :
\begin{center}
\input{konts1.pstex_t}\hspace{1.5cm}
 \input{konts2.pstex_t}.
\end{center}
Le morphisme $m_{6,T}$ vaut
\[
\rho \circ m_2 \circ (H \ts \Id) \circ (m_3 \ts \Id ) \circ (\Id \ts H \ts \Id \tp 2)
\circ (\Id \ts m_3 \ts \Id \tp 2)\circ (i \tp 6).
\]
Le signe $(-1)^s$ associÈ ‡ $T$ est donnÈ par l'ÈgalitÈ
\[
\begin{array}{c}
\rho \circ m_2 \circ (H \ts \Id) \circ (m_3 \ts \Id ) \circ (\Id \ts H \ts \Id \tp 2)
\circ (\Id \ts m_3 \ts \Id \tp 2)\circ (i \tp 6) \circ (\si \tp 6) = \\
(-1)^{s}\si \circ \rho' \circ b_2 \circ (H' \ts \Id) \circ (b_3 \ts \Id ) \circ (\Id \ts H' \ts \Id \tp 2)
\circ (\Id \ts b_3 \ts \Id \tp 2)\circ (i \tp 6),
\end{array}
\]
o˘
\[
\rho' = s \circ \rho \circ \si ,\quad  H' = - s \circ H \circ \si
\quad \mbox{et} \quad i' = s\circ i \circ \si.
\]
Le signe dans le cas gÈnÈral s'obtient de la mÍme maniËre.

Le morphisme $f^\delta : V^\delta \ra A$ est donnÈ par les formules
\[
f^\delta_1 = i , \quad f^\delta_i = \sum_{T \in \mathcal T} (-1)^s f_{i,T}, \quad i \geq 2,
\]
o˘ les morphismes $f_{i,T}$ et le signe $s$ sont construits de la mÍme maniËre en coloriant
la racine de l'arbre $\b T$ par $H$ au lieu de $\rho$.
La remarque (\ref{remarque1_lemme_perturbation}) ci-dessous montrera que
les morphismes $m^\delta_i$, $i\geq 1,$ et $f^\delta_i$, $i\geq 1$,
dÈfinissent bien des $\ai$-structures.

Remarquons que les signes ci-dessus sont tels que
\[
b^\delta_i = \sum_{T \in \mathcal T} b_{i,T} \quad \mbox{et} \quad F^\delta_i = \sum_{T \in \mathcal T}
F_{i,T},
\quad i\geq 1,
\]
o˘ $b_{i,T}$ et $F_{i,T}$ sont obtenus en coloriant les sommets des arbres $\b T$ par des $b_i$
(resp.~$i'$, $\rho'$, $H'$) sur les sommets qui
Ètaient prÈcÈdemment de couleur $m_i$ (resp.~$i$, $\rho$, $H$).

Nous allons  maintenant expliciter l'$\ai$-morphisme
\[
g : A \ra K
\]
en terme d'arbres.
Un calcul facile (nous utilisons le fait que $h^2 = 0$)
montre que le morphisme $g_i$, $i\geq 1$, est donnÈ par les formules
\[
g_1 = p \quad \mbox{et} \quad g_i = - p \circ h \circ  m_i, \quad i \geq 2.
\]
Comme $h\circ p = p\circ H$, les morphismes $g_i$ correspondent aux arbres coloriÈs (ils n'appartiennent
pas nÈcessairement ‡ $\mathcal T$)
\begin{center}
\input{konts2.4.pstex_t}\hspace{2cm} \input{konts2.5.pstex_t}.
\end{center}
Le signe intervenant dans la formule pour $g$ implique
les ÈgalitÈs
\[
G_1 = p' \quad \mbox{et} \quad G_i = - p' \circ H' \circ  b_i, \quad i \geq 2,
\]
o˘ $p' = s \circ p \circ \si$.

Montrons que la composition $g \circ f^\delta$ est nulle.
Il suffit de montrer les ÈgalitÈs
\[
\sum_{\sum \alpha_k =n} G_{i}(F^\delta_{\alpha_1} \ts \hdots \ts F^\delta_{\alpha_i}) = 0, \quad n \geq 1.
\]
Soit $n \geq 1$. Comme les $G_i$ et les $F^\delta_{\alpha_k}$ sont des sommes de compositions
associÈes ‡ des arbres coloriÈs,
la somme ci-dessus est la somme de compositions associÈes aux arbres coloriÈs concatÈnÈs.
Nous vÈrifions que les arbres coloriÈs concatÈnÈs intervenant dans les sommes
\[
\sum_{\sum \alpha_k =n,\ i\geq 2} G_{i}(F^\delta_{\alpha_1} \ts \hdots \ts F^\delta_{\alpha_i})
\quad \mbox{et} \quad G_1 \circ F^\delta_n 
\]
sont les mÍmes. Dans la premiËre somme, le signe intervenant
devant chaque composition associÈe ‡ un arbre coloriÈ concatÈnÈ est nÈgatif
car, pour $i \geq 2$, nous avons $G_i = - p' \circ H' \circ  b_i$. Dans la seconde somme,
il est positif car $G_1 = p'$. Nous avons donc $G \circ F^\delta = 0.$
Le morphisme $f^\delta$ se factorise en $f^g \circ \theta$.
Comme $f^\delta_1 = f^g_1$, nous avons $\theta_1 = \Id_V$. Il s'ensuit
que $\theta : V^\delta \ra V^g$ est un isomorphisme.
 \findem

\begin{remarque}{\em \label{remarque1_lemme_perturbation}
La dÈmonstration montre  que les morphismes $m^\delta_i$, $i\geq 1,$ et $f^\delta_i$, $i\geq 1$,
dÈfinis en terme d'arbres dÈfinissent bien des $\ai$-structures
(voir \cite[6.4]{Kontsevich01} pour
une autre preuve).
}\end{remarque}

\begin{remarque}{\em
Si $A$ est une algËbre diffÈrentielle graduÈe, l'$\ai$-morphis\-me
\[
g : A \ra K
\]
n'a que deux composantes non nulles $g_1$ et $g_2.$ La complexitÈ
des formules pour
les $m^\delta_i$, $i\geq 1$, provient donc de la complexitÈ
des formules des $f^\delta_i$, $i\geq 1$, dÈfinissant le noyau de $g$ dans $\aia$
\[
f^\delta : V^\delta \hookrightarrow A.
\]
}\end{remarque}

\begin{remarque}{\em \label{remarque3_lemme_perturbation}
Le lemme de perturbation nous donne, outre $V^\delta$ et $f^\delta$, une
contraction d'$\ai$-algËbres
\[
\big(\xymatrix{
V^\delta \ar@<-1ex>[r]_(.6){f^\delta} & A \ar@<-1ex>[l]_(.4){q^\delta}
}, H^\delta \big).
\]
Remarquons que l'$\ai$-morphisme $q^\delta$ est le conoyau de l'$\ai$-morphisme
\[
j : K \ra A
\]
donnÈ par les formules
\[
j_1 = \sigma, \quad j_i = - m_i \circ (\sigma \tp i) \circ (h \ts \Id \tp{i-1}), \quad i\geq 2.
\]
\findem
}\end{remarque}

\begin{remarque}{\em
Soit $V$ et $W$ des complexes munis de filtrations exhaustives et soit
\[
\big(\xymatrix{
(V,d_V) \ar@<-1ex>[r]_{i} & (W,d_W) \ar@<-1ex>[l]_{\rho}
}, H \big)
\]
une contraction filtrÈe de complexes.
Alors il existe une suite exacte scindÈe de complexes
\[
\xymatrix{
0 \ar[r] & (V,d_V) \ar@<-1ex>[r]_{i} & (W,d_W) \ar@<-1ex>[r]_{p} \ar@<-1ex>[l]_{\rho}
& (K,d_K) \ar@<-1ex>[l]_{\sigma} \ar[r] & 0
}
\]
telle que
\[
\rho \circ \sigma = 0 \quad \mbox{et} \quad i \circ \rho + \sigma \circ p = \Id_A
\]
et une homotopie contractante $h$ de $K$ telle que
\[
h^2 = 0 \quad \mbox{et} \quad H = \sigma \circ h \circ p.
\]
Le complexe contractile $K$ est donc un facteur direct de $W.$
Soit $\delta$ une perturbation de la diffÈrentielle $d_W.$
Le lemme de perturbation produit une contraction filtrÈe de complexes
\[
\big(\xymatrix{
(V,d_V^\delta) \ar@<-1ex>[r]_(.4){i^\delta} & (W,d_W + \delta) \ar@<-1ex>[l]_(.6){\rho^\delta}
}, H^\delta \big).
\]
Un calcul montre que les morphismes
\[
(p - pH\delta) : (W,d_W + \delta) \ra (K,d_K) \quad \mbox{et} \quad
(\sigma - \delta H \sigma)  :  (K,d_K) \ra (W,d_W + \delta)
\]
sont des morphismes de complexes et qu'ils sont
le conoyau et le noyau de $i^\delta$ et $\rho^\delta$.
La composition
\[
(p - p H \delta) \circ (\sigma - \delta H \sigma )  : (K,d_K) \ra (K,d_K)
\]
induit un isomorphisme dans les objets graduÈs associÈs ‡ la filtration.
C'est donc un isomorphisme. Le complexe contractile$(K,d_K)$ est donc aussi
un facteur direct du complexe perturbÈ  $(W,d_W + \delta)$ et l'inclusion
\[
\sigma : K \ra W
\]
est ``perturbÈe'' en $\sigma - \delta H \sigma$ pour devenir compatible ‡
$d_W + \delta.$
}
\end{remarque}

%% file: konts1.pstex_t
\begin{picture}(0,0)%
\includegraphics{konts1.pstex}%
\end{picture}%
\setlength{\unitlength}{3947sp}%
\begingroup\makeatletter\ifx\SetFigFont\undefined%
\gdef\SetFigFont#1#2#3#4#5{%
  \reset@font\fontsize{#1}{#2pt}%
  \fontfamily{#3}\fontseries{#4}\fontshape{#5}%
  \selectfont}%
\fi\endgroup%
\begin{picture}(1593,2489)(558,-1757)
\put(953,-1460){\makebox(0,0)[lb]{\smash{\SetFigFont{12}{14.4}{\rmdefault}{\mddefault}{\updefault}{\color[rgb]{0,0,0}$T$}%
}}}
\end{picture}

%% file: konts2.pstex_t
\begin{picture}(0,0)%
\includegraphics{konts2.pstex}%
\end{picture}%
\setlength{\unitlength}{3947sp}%
\begingroup\makeatletter\ifx\SetFigFont\undefined%
\gdef\SetFigFont#1#2#3#4#5{%
  \reset@font\fontsize{#1}{#2pt}%
  \fontfamily{#3}\fontseries{#4}\fontshape{#5}%
  \selectfont}%
\fi\endgroup%
\begin{picture}(1594,2644)(586,-1790)
\put(2101,-586){\makebox(0,0)[lb]{\smash{\SetFigFont{12}{14.4}{\rmdefault}{\mddefault}{\updefault}{\color[rgb]{0,0,0}$i$}%
}}}
\put(1651,689){\makebox(0,0)[lb]{\smash{\SetFigFont{12}{14.4}{\rmdefault}{\mddefault}{\updefault}{\color[rgb]{0,0,0}$i$}%
}}}
\put(1126,689){\makebox(0,0)[lb]{\smash{\SetFigFont{12}{14.4}{\rmdefault}{\mddefault}{\updefault}{\color[rgb]{0,0,0}$i$}%
}}}
\put(1077,-541){\makebox(0,0)[lb]{\smash{\SetFigFont{12}{14.4}{\rmdefault}{\mddefault}{\updefault}{\color[rgb]{0,0,0}$m_3$}%
}}}
\put(1636,-1205){\makebox(0,0)[lb]{\smash{\SetFigFont{12}{14.4}{\rmdefault}{\mddefault}{\updefault}{\color[rgb]{0,0,0}$m_2$}%
}}}
\put(1084,175){\makebox(0,0)[lb]{\smash{\SetFigFont{12}{14.4}{\rmdefault}{\mddefault}{\updefault}{\color[rgb]{0,0,0}$m_3$}%
}}}
\put(1420,-875){\makebox(0,0)[lb]{\smash{\SetFigFont{12}{14.4}{\rmdefault}{\mddefault}{\updefault}{\color[rgb]{0,0,0}$H$}%
}}}
\put(1101,-180){\makebox(0,0)[lb]{\smash{\SetFigFont{12}{14.4}{\rmdefault}{\mddefault}{\updefault}{\color[rgb]{0,0,0}$H$}%
}}}
\put(725,689){\makebox(0,0)[lb]{\smash{\SetFigFont{12}{14.4}{\rmdefault}{\mddefault}{\updefault}{\color[rgb]{0,0,0}$i$}%
}}}
\put(586, 85){\makebox(0,0)[lb]{\smash{\SetFigFont{12}{14.4}{\rmdefault}{\mddefault}{\updefault}{\color[rgb]{0,0,0}$i$}%
}}}
\put(1734, 26){\makebox(0,0)[lb]{\smash{\SetFigFont{12}{14.4}{\rmdefault}{\mddefault}{\updefault}{\color[rgb]{0,0,0}$i$}%
}}}
\put(1695,-1790){\makebox(0,0)[lb]{\smash{\SetFigFont{12}{14.4}{\rmdefault}{\mddefault}{\updefault}{\color[rgb]{0,0,0}$\rho$}%
}}}
\put(593,-1520){\makebox(0,0)[lb]{\smash{\SetFigFont{12}{14.4}{\rmdefault}{\mddefault}{\updefault}{\color[rgb]{0,0,0}$\b T$ colori\'e}%
}}}
\end{picture}

%% file: konts2.4.pstex_t
\begin{picture}(0,0)%
\includegraphics{konts2.4.pstex}%
\end{picture}%
\setlength{\unitlength}{3947sp}%
\begingroup\makeatletter\ifx\SetFigFont\undefined%
\gdef\SetFigFont#1#2#3#4#5{%
  \reset@font\fontsize{#1}{#2pt}%
  \fontfamily{#3}\fontseries{#4}\fontshape{#5}%
  \selectfont}%
\fi\endgroup%
\begin{picture}(109,1254)(728,-678)
\put(728,-620){\makebox(0,0)[lb]{\smash{\SetFigFont{12}{14.4}{\rmdefault}{\mddefault}{\updefault}{\color[rgb]{0,0,0}$\hspace*{.1cm}$}%
}}}
\put(767,-167){\makebox(0,0)[lb]{\smash{\SetFigFont{12}{14.4}{\rmdefault}{\mddefault}{\updefault}{\color[rgb]{0,0,0}$p$}%
}}}
\put(761,411){\makebox(0,0)[lb]{\smash{\SetFigFont{12}{14.4}{\rmdefault}{\mddefault}{\updefault}{\color[rgb]{0,0,0}$\Id$}%
}}}
\end{picture}

%% file: konts2.5.pstex_t
\begin{picture}(0,0)%
\includegraphics{konts2.5.pstex}%
\end{picture}%
\setlength{\unitlength}{3947sp}%
\begingroup\makeatletter\ifx\SetFigFont\undefined%
\gdef\SetFigFont#1#2#3#4#5{%
  \reset@font\fontsize{#1}{#2pt}%
  \fontfamily{#3}\fontseries{#4}\fontshape{#5}%
  \selectfont}%
\fi\endgroup%
\begin{picture}(1144,1521)(560,-708)
\put(560,648){\makebox(0,0)[lb]{\smash{\SetFigFont{12}{14.4}{\rmdefault}{\mddefault}{\updefault}{\color[rgb]{0,0,0}$\Id$}%
}}}
\put(943,644){\makebox(0,0)[lb]{\smash{\SetFigFont{12}{14.4}{\rmdefault}{\mddefault}{\updefault}{\color[rgb]{0,0,0}$\Id$}%
}}}
\put(1288,644){\makebox(0,0)[lb]{\smash{\SetFigFont{12}{14.4}{\rmdefault}{\mddefault}{\updefault}{\color[rgb]{0,0,0}$\Id$}%
}}}
\put(1622,640){\makebox(0,0)[lb]{\smash{\SetFigFont{12}{14.4}{\rmdefault}{\mddefault}{\updefault}{\color[rgb]{0,0,0}$\Id$}%
}}}
\put(1081,-263){\makebox(0,0)[lb]{\smash{\SetFigFont{12}{14.4}{\rmdefault}{\mddefault}{\updefault}{\color[rgb]{0,0,0}$H$}%
}}}
\put(1073,156){\makebox(0,0)[lb]{\smash{\SetFigFont{12}{14.4}{\rmdefault}{\mddefault}{\updefault}{\color[rgb]{0,0,0}$m_i$}%
}}}
\put(1092,-650){\makebox(0,0)[lb]{\smash{\SetFigFont{12}{14.4}{\rmdefault}{\mddefault}{\updefault}{\color[rgb]{0,0,0}$p$}%
}}}
\end{picture}

%% file: Homot_aimod.tex
{\bf \noindent Introduction}\\

\noindent Soit $A$ une $\ai$-algËbre augmentÈe. Rappelons que dans cette thËse les structures
communÈment appelÈes $\ai$-modules sur $A$ sont appelÈes {\em $A$-polydules}
(``poly'' car la structure est donnÈe par plusieurs
multiplications).
Le but de ce chapitre est de dÈcrire la catÈgorie dÈrivÈe
$\cd_\infty A$ dont les objets sont les $A$-polydules strictement
unitaires. Nous utiliserons pour cela les outils de l'algËbre homotopique
de Quillen (voir l'appendice \ref{section_rappels_cmf})
en adaptant les mÈthodes du chapitre 1 aux
polydules.
La catÈgorie dÈrivÈe d'une $\ai$-algËbre quelconque sera ÈtudiÈe au chapitre \ref{chapitre_Cat_der}.\\

{\bf \noindent Plan du chapitre}\\

\noindent Ce chapitre est divisÈ en deux parties. 

La premiËre partie qui composÈe des sections 
(\ref{section_reduction_augmentation}) et (\ref{section_cmf_coModcu}) ne traitera pas
des $\ai$-structures proprement dites. Dans la premiËre section (\ref{section_reduction_augmentation}),
on dÈfinit les (co)modules diffÈrentiels
graduÈs (co)unitaires. Dans la section \ref{section_cmf_coModcu},
nous dÈmontrons le thÈorËme (\ref{theoreme_cmf_coModcu}) :

{\em \noindent 
 Soit $C$ une cogËbre diffÈrentielle graduÈe cocomplËte co-augmentÈe.
La catÈgorie $\coModcu C$ des $C$-comodules diffÈrentiels graduÈs co-unitaires cocomplets admet
une unique structure de catÈgorie de modËles telle que, pour toute 
algËbre diffÈrentielle graduÈe augmentÈe $A$ et toute cochaÓne
tordante admissible acyclique $\tau : C \ra A$, le couple de foncteurs
adjoints
\[
(? \tw A, \? \tw C) : \coModcu C \ra \Modu A
\]
est une Èquivalence de Quillen. Tous les objets de $\coModcu C$ sont cofibrants.
}

\noindent Nous caractÈrisons ensuite l'acyclicitÈ des cochaÓnes tordantes
(\ref{caracterisation_cochaines_acycliques}).\\

La deuxiËme partie est consacrÈe aux $\ai$-structures concernÈes par ce chapitre :
les (bi)polydules strictement unitaires sur des $\ai$-algËbres augmentÈes.
Dans la section \ref{section_definition_ai-modules}, on dÈfinit les polydules,
leurs suspensions, les $\ai$-morphismes et les homotopies entre $\ai$-morphis\-mes.
Nous dÈfinissons ensuite la notion d'unitaritÈ stricte pour les $\ai$-structures.
Cette notion sera ÈtudiÈe plus prÈcisÈment dans le chapitre \ref{chapitre_Unite}.
On rappelle ensuite les constructions bar et cobar et l'algËbre enveloppante.
Dans la section \ref{section_categorie_derivee_augmentee}, nous affinons le thÈorËme
(\ref{theoreme_cmf_coModcu}) prÈcitÈ. Nous montrons que, si la
cogËbre $C$ est isomorphe, en tant que cogËbre graduÈe, ‡ une cogËbre tensorielle
co-augmentÈe, les objets fibrants de $\coModcu C$ sont exactement les facteurs directs des
$C$-comodules presque colibres.
En particulier, dans le cas o˘ $C$ est Ègale ‡ la construction bar d'une $\ai$-algËbre
augmentÈe $A$, la catÈgorie des objets fibrants et cofibrants de $\coModcu C$ est l'image
essentielle par la
construction bar des $A$-polydules strictement unitaires.
Nous dÈduirons de ce rÈsultat plusieurs descriptions de la catÈgorie dÈrivÈe
\[
\cd_\infty A = \aiModu A [Qis^{-1}],
\]
o˘ $\aiModu A$ dÈsigne la catÈgorie des $A$-polydules strictement unitaires.

Dans la section \ref{section_categorie_derivee_augmentee_bipolydules},
nous Ètudions la catÈgorie dÈrivÈe $\cd_\infty (A,A')$ des
bipolydules strictement unitaires sur $A$ et $A'$, deux $\ai$-algËbres augmentÈes.
Les mÈthodes Ètant similaires, les dÈtails seront omis.
Les bipolydules seront utiles dans l'Ètude des
$\ai$-catÈgories.

\section{Rappels et notations} \label{section_reduction_augmentation}

Soit $(\sf C,\ts,e)$ une $\corps$-catÈgorie de Grothendieck semi-simple monoÔdale
et $\sf C'$ \indexnotation{C1}
une $\corps$-catÈgorie de Grothendieck semi-simple (non nÈcessairement 
monoÔdale).
Nous supposons que la catÈgorie monoÔdale $\sf C$ agit ‡
droite sur $\sf C'$, i.~e.~$\sf C'$ est munie d'un foncteur
\[
\sf C' \times \sf C \ra \sf C', \quad 
(M,A) \mapsto M \ts A
\]
tel que
\[
\Hom_{\sf C'}(M,M') \times \Hom{\sf C}(A,A') \ra \Hom_{\sf C'}(M\ts A,M' \ts A'),
\]
o˘ $A,A'$ sont dans $\sf C$ et $M,M'$ sont dans $\sf C',$ est $\corps$-bilinÈaire. Nous demandons en
outre que cette action soit associative et unitaire ‡ des isomorphismes donnÈs prËs
(voir \cite[chap.~XI]{MacLane98}).

\subsection{Modules sur une algËbre augmentÈe} \label{section_modules_augmentes}

Soit $\sf M$ (resp.~$\sf M'$) l'une des catÈgories  $\gr \sf C$ ou $\cc \sf C$
(resp.~$\gr \sf C'$ ou $\cc \sf C'$) dÈfinies ‡ la section \ref{categorie_de_base}.
La catÈgorie $\sf M$ est monoÔdale et agit clairement sur $\sf M'$. \\

{\noindent \bf AlgËbres augmentÈes, rÈduites}\\

Une algËbre  $(A,\mu)$ dans $\sf M$ est {\em unitaire}  si elle est munie
d'un morphisme $\eta : e \ra A$ tel que $\mu (\Id \ts \eta) = \mu
(\eta \ts \Id) = \Id$. On appelle le morphisme $\unite$ l'{\em unitÈ}  \indexnotation{unite}
de $A.$
Si $A$ et $A'$ sont des algËbres unitaires,
un {\em morphisme} d'algËbres unitaires $f :A \ra A'$ est un morphisme d'algËbres
$f$ tel que $f  \unite_A = \unite_{A'}.$
Le morphisme $e \ts e \ra e$ donnÈ par la contrainte d'unitaritÈ de
la catÈgorie de base (\ref{categorie_de_base}) dÈfinit une structure d'algËbre unitaire
sur l'objet neutre
$e$.
Une algËbre $A$ est {\em augmentÈe} 
si elle est unitaire et munie d'un morphisme
d'algËbres unitaires \indexnotation{augmentation_mor}
\[
\epsilon : A \ra e.
\]
Le morphisme $\epsilon$ s'appelle l'{\em augmentation} \index{augmentation} de $A$.
Si $A$ et $A'$ sont des algËbres augmentÈes,
un {\em morphisme} d'algËbres augmentÈes $f :A \ra A'$ est un morphisme d'algËbres
unitaires $f$ tel que $\epsilon_{A'}  f  = \epsilon_A.$

Si $A$ est une algËbre augmentÈe de $\sf M$,
l'{\em algËbre rÈduite} \index{rÈduite (algËbre)} \index{algebre@{algËbre}!reduite@{rÈduite}}
\indexnotation{reduction}
$\b A $ associÈe ‡ $A$ est  le noyau de l'augmentation.
Si $A$ est une algËbre de $\sf M$,
l'{\em algËbre augmentÈe} \indexnotation{augmentation}
associÈe ‡ $A$ est l'algËbre $A\+$ dont l'objet sous-jacent
est $e \oplus A$, et dont la multiplication est dÈfinie par les morphismes
\[
e \ts e \ra e, \quad e \ts A \ra A,\quad A \ts e \ra A\quad \mbox{et} \quad A \ts A \arr{\mu} A,
\]
o˘ les trois premiers morphismes sont donnÈs par la contrainte d'unitaritÈ de la
catÈgorie de base. L'augmentation de $A\+$ est la projection canonique
$A\+ \ra e$.
On note $\alga$ \indexnotation{alga} la {\em catÈgorie} des algËbres augmentÈes de $\cc \sf C.$
Le foncteur 
\[
 \alg \arr{} \alga, \quad A \mapsto A\+,
\]
est une Èquivalence dont le quasi-inverse est le foncteur $A \mapsto \b A.$\\

{\noindent \bf Modules}\\

Soit $A$ une algËbre dans $\sf M$. Un {\em $A$-module (‡ droite) \index{module}
dans $\sf M'$} est un objet $M$ de $\sf M'$
muni d'un morphisme $\mu^M : M \ts A \ra M$ (de degrÈ $0$ si $\sf M' = \gr \sf C'$)
tel que
\[
\mu^M (\mu^M \ts \Id) = \mu^M (\Id \ts \mu^A).
\]
On appelle $\mu^M$ la {\em multiplication de $M$.}
Si $M$ et $N$ sont deux modules, un {\em morphisme} de modules $f : M \ra N$
est un morphisme $f$ tel que 
\[
f  \mu^M = \mu^N (f \ts \Id).
\]
Si l'algËbre $A$ est unitaire, un $A$-module $M$ est {\em unitaire} 
si on a 
\[
\mu^M (\Id \ts \unite^A) = \Id_M.
\]
Soit $A$ une algËbre graduÈe (resp.~diffÈrentielle graduÈe).
Un $A$-{\em module graduÈ} \index{module!graduÈ} (resp.~{\em diffÈrentiel graduÈ})
\index{module!diffÈrentiel graduÈ} est un $A$-module dans
la catÈgorie $\gr \sf C'$ (resp.~$\cc \sf C'$).
Si $A$ est une algËbre diffÈrentielle graduÈe,
un $A$-module diffÈrentiel graduÈ est donc un objet $M$ de $\gr \sf C'$, muni d'une multiplication
$\mu^M : M\ts A\ra M$ et d'une diffÈrentielle $d^M : M \ra M$ telle que
\[
d^M (\mu^M) = \mu^M (d^M \ts \Id_A + \Id_M \ts d^A).
\]
Si $(M,\mu^M)$ est un $A$-module graduÈ, une {\em dÈrivation de modules}
\index{derivation@{dÈrivation}} est un morphisme $d^M : M \ra M$ vÈrifiant l'Èquation ci-dessus.
Une {\em diffÈrentielle de module} est une dÈrivation de degrÈ $+1$ et de carrÈ nul.
Si $A$ est une algËbre diffÈrentielle graduÈe unitaire,
on note $\Modu A$ \label{Modu}la {\em catÈgorie} des $A$-modules dif\-fÈrentiels graduÈs
unitaires.

Soit $f :  A \ra A'$ un morphisme de $\alg$. La {\em restriction le long de $f$}
\index{restriction} d'un $A'$-module $M$
est le $A$-module dont l'objet sous-jacent est $M$ et dont la
multiplication est $\mu^M (f \ts \Id).$ Le {\em $A'$-module induit par $f$}
\index{induction} d'un $A$-module $M$ a pour objet sous-jacent
$M \ts_{A} A'$ et pour multiplication  $\Id \ts \mu^{A'}$.
Soit $A$ une algËbre augmentÈe et soit  $i : \b A \ra A$ l'injection canonique. 
Le {\em foncteur restriction}
est une Èquivalence de $\Modu A$ sur la catÈgorie des modules diffÈrentiels graduÈs sur $\b A$,
son quasi-inverse
est le {\em foncteur induction}.\\

Soit $A$ une algËbre diffÈrentielle graduÈe et $M$ et $N$ deux modules diffÈrentiels graduÈs.
Si $f$ et $g$ sont deux morphismes $M \ra N$, une {\em homotopie} \index{homotopie}
entre $f$ et $g$ est un morphisme graduÈ de $A$-modules $h : M \ra N$ de degrÈ $-1$ tel que
$h \circ d + d \circ h = f - g.$ Deux morphismes $f$ et $g$ sont {\em homotopes} s'il existe une homotopie
entre $f$ et $g.$\\

{\noindent \bf Modules libres}\\

Soit $A$ une algËbre de $\sf M$. Soit $V$ un objet de $\sf M'.$ Le morphisme $\Id_V \ts \mu^A$
dÈfinit une structure de $A$-module sur $V \ts A$.
Un $A$-module $M$ est {\em libre sur $V$} \index{libre!module} \index{module!libre}s'il existe un isomorphisme
de $A$-modules $M \arr{\sim} V \ts A$.
Un module diffÈrentiel graduÈ est {\em presque libre} \index{module!presque libre} \index{presque libre!module}
s'il est libre en tant que module graduÈ.

\begin{lemme} \label{lemme_modules_libres}
Soit $A$ un objet de $\alga$.
Soit $M$ un objet de $\Modu A$ et $V$ un objet de $\gr \sf C'$.
\english \begin{itemize}
\item[a.] L'application $f \mapsto f (\Id \ts \unite)$ est une bijection de l'ensemble
des morphismes de modules graduÈs $V\ts A \ra M$ vers l'ensemble des morphismes
graduÈs $V \ra M$. L'application inverse associe ‡ $g : V \ra M$ le morphisme
de modules 
\[
V \ts A \arr{g \ts \Id} M \ts A \arr{\mu^M} M.
\]
\item[b.] L'application $d \mapsto d (\Id \ts \unite)$ est une bijection de l'ensemble $\ce$
des dÈrivations du modules graduÈs $V\ts A$ vers l'ensemble des morphismes
graduÈs $g : V \ra V \ts A$.
L'application inverse associe ‡ $g : M \ra N$ la diffÈrentielle
\[
\Id \ts d^A + (\Id \ts \mu^A)(g \ts \Id).
\]
Cette bijection fait correspondre le sous-ensemble de $\ce$
formÈ des diffÈrentielles de modules au morphismes de degrÈ $+1$ tel que 
\[
(\Id_V \ts \mu^A)(g \ts \Id)g + (\Id \ts d^A)g = 0.
\]
\end{itemize} \francais \findem
\end{lemme}

\subsection{Comodules co-augmentÈs}

{\noindent \bf CogËbres co-augmentÈes, rÈduites}\\

Une cogËbre $(C,\Delta)$ de $\sf M$ est {\em co-unitaire} \index{co-unitaire} si elle est munie d'un morphisme
$\counite : C \ra e$ tel que $(\Id \ts \counite)\Delta = (\counite \ts \Id )\Delta = \Id$. Le morphisme
$\counite$ s'appelle la {\em co-unitÈ} \indexnotation{counite} de $C$. Si $C$ et $C'$ sont
deux cogËbres co-unitaires, un {\em morphisme} de cogËbres co-unitaires $f : C \ra C'$ est un morphisme
de cogËbres $f$ tel que $\counite_{C'}  f = \counite_C.$ 
Le morphisme $e \ra e \ts e$ donnÈ par la contrainte d'unitaritÈ de la catÈgorie de base
dÈfinit une structure de cogËbre co-unitaire sur l'objet neutre $e.$
Une cogËbre $C$ est {\em co-augmentÈe}  si elle est munie d'un
morphisme de cogËbres co-unitaires \indexnotation{coaugmentation_mor}
\[
\coaugmentation : e \ra C.
\]
Le morphisme $\coaugmentation$ s'appelle la {\em co-augmentation} \index{co-augmentation} de la cogËbre $C.$
Si $C$ et $C'$ sont deux cogËbres co-augmentÈes, un {\em morphisme} de cogËbres co-augmentÈes
$f : C\ra C'$ est un morphisme de cogËbres unitaires $f$ tel que $f  \coaugmentation_C = 
\coaugmentation_{C'}$.

Si $C$ est une cogËbre co-augmentÈe de $\sf M$, la {\em cogËbre rÈduite} \index{rÈduite (cogËbre)}
\index{cogebre@{cogËbre}!reduite@{rÈduite}}
$\b C$ \indexnotation{reduction_cog} est le conoyau
de la co-augmentation.
Si $C$ est une cogËbre de $\sf M$, la {\em cogËbre co-augmentÈe} \index{co-augmentÈe (cogËbre)}
\index{cogebre@{cogËbre}!co-augmentee@{co-augmentÈe}}
$C\+$ \indexnotation{coaugmentation} est la cogËbre 
dont l'objet sous-jacent est $C \oplus e$ et dont la comultiplication est le morphisme dÈfini par
les composantes
\[
e \ra e \ts e,\quad C \ra e \ts C, \quad C \ra C \ts e \quad \mbox{et} \quad C \arr{\Delta}C \ts C,
\]
o˘ les trois premiers morphismes sont dÈfinis par la contrainte d'unitaritÈ de la catÈgorie de base.
La co-augmentation de $C\+$ est l'injection canonique $e \ra C\+.$ 
Si $V$ est un objet graduÈ de $\sf C$, on note $\ct V$ la cogËbre \indexnotation{ct} $(\ctr V)\+.$
Soit $\cocoga$ \label{cocoga} la {\em catÈgorie}
des cogËbres co-augmentÈes de $\cc \sf C$ dont les cogËbres rÈduites sont
cocomplËtes.
Le foncteur
\[
\cocog \ra \cocoga, \quad C \mapsto C\+,
\]
est une Èquivalence dont le quasi-inverse est le foncteur $C \ra \b C.$\\

{\bf \noindent Comodules}\\

Soit $C$ une cogËbre de $\sf M.$ Un $C$-{\em comodule (‡ droite) dans $\sf M'$ } \index{comodule}
est un objet graduÈ $N$ de $\sf M'$ muni
d'un morphisme $\Delta^N : N \ra N \ts C$ (de degrÈ $0$ si $\sf M' = \gr \sf C'$) tel que
\[
(\Id \ts \Delta^C) \Delta^N = (\Delta^N \ts \Id) \Delta^N.
\]
Si $N$ et $N'$ sont deux $C$-comodules, un {\em morphisme} de $C$-comodules $f : N \ra N'$
est un morphisme de $\sf M'$ tel que $\Delta^{N'}  f = (f \ts \Id)  \Delta^N$.
Si la cogËbre $C$ est co-unitaire, un $C$-comodule $N$ est
{\em co-unitaire} \index{co-unitaire} si $\Delta^N (\Id \ts \counite) = \Id_N.$

Soit $C$ une cogËbre graduÈe (resp.~diffÈrentielle graduÈe).
Un $C$-{\em comodule graduÈ} \index{comodule!graduÈ} (resp.~{\em diffÈrentiel graduÈ})
\index{comodule!diffÈrentiel graduÈ} est un $C$-comodule dans la
catÈgorie $\gr \sf C'$ (resp.~$\cc \sf C'$).
Si $C$ est une cogËbre diffÈrentielle graduÈe,
un comodule diffÈren\-tiel graduÈ est donc un objet $N$ de $\gr \sf C'$, muni d'une comultiplication
$\Delta^N : N \ra N \ts C$ et d'une diffÈrentielle $d^N : N \ra N$ telle que
\[
\Delta^N d^N  = (d^N \ts \Id_A + \Id_N \ts d^A) \Delta^N.
\]
Si $(N,\Delta^N)$ est un $C$-comodule graduÈ, une {\em codÈrivation de comodules}
\index{coderivation@{codÈrivation}} est un morphisme $d^N : N \ra N$ vÈrifiant l'Èquation ci-dessus.
Une {\em diffÈrentielle de comodule} est une codÈrivation de degrÈ $+1$ et de carrÈ nul.
Si la cogËbre $C$ est co-unitaire, on note $\coModu C$ \indexnotation{comC}
la catÈgorie des comodules diffÈrentiels graduÈs co-unitaires.

Soit $f : C \ra C'$ un morphisme de $\cog.$ La {\em corestriction le long de $f$}
\index{corestriction} d'un
$C$-comodule $N$ est le $C'$-comodule dont l'objet sous-jacent est $N$ et dont
la comultiplication est $(\Id \ts f) \Delta^N.$ 
Le {\em $C$-comodule co-induit par $f$} \index{co-induction}
associÈ ‡ un $C'$-comodule $N$ a pour objet sous-jacent le noyau
\[
\ker ( N\ts C \arr{u} N \ts C' \ts C),
\]
o˘ $u = \Delta^N  \ts \Id_C - (\Id_N \ts f \ts \Id_C)(\Id_N \ts \Delta^C)$,
et pour comultiplication le morphisme induit par $\Id^N \ts \Delta^C : N \ts C \ra N\ts C \ts C$.

Soit $C$ une cogËbre co-augmentÈe et soit  $p : C \ra \b C$ la projection canonique. 
Le {\em foncteur corestriction}
est une Èquivalence de la catÈgorie $\coModu C$ sur la catÈgorie des $\b C$-comodules diffÈrentiels
graduÈs.
Son quasi-inverse
est le {\em foncteur co-induction}.\\

Soit $C$ une cogËbre diffÈrentielle graduÈe et soit $N$ et $N'$ deux comodules diffÈrentiels graduÈs.
Si $f$ et $g$ sont deux morphismes $N \ra N'$, une {\em homotopie} \index{homotopie}
entre $f$ et $g$ est un morphisme graduÈ de $C$-comodules $h : N \ra N'$ de degrÈ $-1$ tel que
$h \circ d + d \circ h = f - g.$ Deux morphismes $f$ et $g$ sont {\em homotopes} s'il existe une homotopie
entre $f$ et $g.$\\

{\noindent \bf Comodules cocomplets}\\

Soit $C$ une cogËbre co-augmentÈe de $\sf M$ et $N$ un $C$-comodule co-unitaire dans $\sf M'$.
On dÈfinit $\Delta^{(2)} = \Delta^N$
et, pour tout $n\geq 3$, on dÈfinit
$\Delta ^{(n)} : N \ra N  \ts C\tp{n-1}$ par
\[
\Delta ^{(n)} = (\Id \tp{n-2} \ts \Delta^C) \Delta^{(n-1)}.
\]
Soit $n \geq 1.$ Le noyau $N\prim n$ \indexnotation{primN} du morphisme
\[
N \arr{\Delta^{(n+1)}} N \ts C\tp{n} \arr{\Id \ts p\tp{n}} N \ts \b C \tp{n}
\]
(o˘ $p : C \ra \b C$ est la projection canonique) est un sous-comodule de $N.$
Il s'appelle le {\em sous-comodule des $n$-primitifs}
\index{primitifs@{$n$-primitifs}} de $N.$
Pour $n = 1$, on obtient le sous-comodule des {\em primitifs} de $N$.
La suite croissante de sous-comodules
\[
N_{[1]} \subset N_{[2]} \subset N_{[3]} \subset \cdots
\]
est la {\em filtration primitive}
\index{filtration!primitive} du comodule
$N.$
Si $C$ est un objet de $\cocoga$,
un $C$-comodule diffÈrentiel graduÈ co-unitaire $N$ est {\em cocomplet} \index{cocomplet (comodule)}
\index{comodule!cocomplet}
si sa filtration primitive est exhaustive. On note $\coModcu C$ \indexnotation{comcC} la catÈgorie des
comodules cocomplets.\\

{\noindent \bf Comodules colibres}\\

Soit $C$ une cogËbre co-augmentÈe dans $\sf M$. Soit $V$ un objet de $\sf M'.$ Le
morphisme 
\[
\Id \ts \Delta^C : V\ts C \ra V \ts C \ts C
\]
munit $V \ts C$
d'une structure de $C$-comodule. Son sous-comodule des primitifs est le comodule $V \ts e.$
Pour $n \geq 2$, son sous-comodule des $n$-primitifs est le $C$-comodule
$V \ts C\prim {n-1}.$ Le $C$-comodule $V\ts C$ est
donc cocomplet si $C$ est un objet de $\cocoga$. Un $C$-module
$N$ est {\em colibre sur $V$} \index{colibre (comodule)} \index{comodule!colibre} s'il existe
un isomorphisme de $C$-comodules $N \arr{\sim} V \ts C.$
Si $C$ est un objet de $\cocoga$, un comodule diffÈrentiel graduÈ est
{\em presque colibre} \index{comodule!presque colibre} \index{presque colibre (comodule)}
s'il est libre en tant que comodule graduÈ.
La {\em sous-catÈgorie} de $\coModcu C$ formÈe des objets presque colibres est notÈe $\colibre C.$

\begin{lemme} \label{lemme_comodules_colibres}
Soit $C$ une cogËbre diffÈrentielle graduÈe co-unitaire,
$N$ un objet dans $\coModu C$ et $V$ un objet graduÈ.
\english \begin{itemize}
\item[a.] L'application $f \mapsto (\Id \ts \counite^C) f$ est une bijection de l'ensemble
des morphismes de comodules graduÈs $N \ra V\ts C$ sur l'ensemble des morphismes graduÈs
$N \ra V.$ L'application inverse envoie $g : N \ra V$ sur le morphisme de $C$-comodules
\[
N \arr{\Delta} N \ts C \arr{g \ts \Id} V \ts C.
\]
\item[b.]
L'application $d \mapsto (\Id \ts \counite^C) d$ est une bijection de l'ensemble 
\[
\coder(V \ts C)
\]
des codÈrivations du comodules $V \ts C$ sur l'ensemble des morphismes graduÈs
$g : V \ts C \ra V$.
L'application inverse envoie $g$ sur la codÈrivation
\[
(g \ts \Id )(\Id_V \ts \Delta^C) + \Id_V \ts d^C.
\]
Cette bijection fait correspondre les diffÈrentielles de comodules
aux morphismes graduÈs
de degrÈ $+1$ tels que
\[
g (\Id_V \ts d^C) + g (g \ts \Id_C )(\Id_V \ts \Delta^C) = 0.
\]
\end{itemize} \francais \findem
\end{lemme}

%
%
%
%
%
%
%
%
%
%
%
%
%
%

\section{$\coModcu C$ comme catÈgorie de modËles}  \label{section_cmf_coModcu}
%
%
%
%
%
%
%
%
%
%
%
%
%
%

\subsection{CochaÓne tordante et produits tensoriels tordus} 
\label{section_cochaines_tordantes}

\begin{definition}{\em 
Soit $C$ une cogËbre diffÈrentielle graduÈe et $A$ une al\-gËbre diffÈrentielle
graduÈe. Une {\em cochaÓne tordante}\index{cochaÓne tordante} \indexnotation{tau}
est un morphisme graduÈ $\tau : C \ra A$
de degrÈ $+1$ tel que 
\[
d_A \tau + \tau d_C + m (\tau \ts \tau) \Delta = 0.
\]
Si $f: A \ra A'$ est un morphisme dans $\alg$ (resp.~si $g : C' \ra C$ est un morphisme
dans $\cog$) la composition $f \circ \tau$ (resp.~$\tau \circ g$) est encore une cochaÓne
tordante. Ainsi, une cochaÓne tordante $\tau : C \ra A$ induit une cochaÓne tordante
$\tau \+ = i \circ \tau \circ p : C\+ \ra A\+$, o˘ $i$ est l'injection canonique $A \ra A\+$
et $p$ la projection canonique $C\+ \ra C$.
Soit $A$ un objet de $\alga$ et $C$ un objet de $\coga$. Une cochaÓne tordante
$ C \ra A$ est {\em admissible} \index{admissible!cochaÓne} \index{cochaÓne tordante!admissible}
si elle est induite par une cochaÓne tordante $\b C \ra \b A$.
}\end{definition}

Soit $A$ une algËbre diffÈrentielle graduÈe augmentÈe et $C$ une cogËbre diffÈrentielle
graduÈe co-augmentÈe.
Soit $\tau : C \ra A$ une cochaÓne tordante admissible.
Soit $M$ un objet de $\Modu A$. Soit le morphisme $t_\tau :  M \ts C \ra M \ts C$ dÈfini comme la composition
\[
M \ts C \arr{\Id \ts \Delta} M \ts C \ts C \arr{\Id \ts \tau \ts \Id} M \ts A \ts C
\arr{\mu^M \ts \Id} M \ts C.
\]
Comme $\tau$ est une cochaÓne tordante, la somme
\[
b_\tau = b + t_\tau : M\ts C \arr{} M \ts C
\]
o˘ $b$ est la diffÈrentielle du produit tensoriel $M\ts C$,
donne une diffÈrentielle sur le $C$-comodule graduÈ co-unitaire $M \ts C.$
Le produit tensoriel $M \ts C$ muni de la {\em diffÈrentielle tordue (par $\tau$) $b_\tau$}
\index{diffÈrentielle tordue} \index{produit tensoriel tordu} est notÈ $M \tw C.$
Si $M$ et $M'$ sont deux objets de $\Modu A$, un morphisme $f : M \ra M'$
induit un morphisme de $C$-comodules graduÈs co-unitaires
$f \ts \Id_C: M \tw C\ra M' \tw C$ compatible aux diffÈrentielles. 
Nous obtenons ainsi un {\em foncteur}\indexnotation{RM}
\[
R_\tau : \Modu A \ra \coModu C, \hspace{1cm} M \mapsto M \tw C.
\]
Lorsqu'il n'y aura pas d'ambiguitÈ nous noterons ce foncteur $R.$

De maniËre duale, si $N$ est $C$-comodule diffÈrentiel graduÈ co-unitaire,
le {\em morphisme} $T_\tau$ est dÈfini
comme la composition
\[
N \ts A \arr{\Delta^N \ts \Id} N \ts C \ts A \arr{\Id \ts \tau \ts \Id} N \ts A \ts A
\arr{\Id \ts \mu^A} N \ts C.
\]
La somme de la diffÈrentielle $D$ du produit tensoriel $N \ts A$ et du morphisme $T_\tau$ dÈfinit
une nouvelle diffÈrentielle sur le $A$-module graduÈ unitaire $N \ts A$. 
Le produit tensoriel $N \ts A$ muni de
la {\em diffÈrentielle tordue (par $\tau$)  $D_\tau = D + T_\tau$} 
est notÈ $N \tw A.$ \indexnotation{produit_tens_tordu}
Si $N$ et $N'$ sont deux
objets de $\coModcu C$, un morphisme $f : N \ra N'$ induit un morphisme de $A$-modules graduÈs unitaires
$f \ts \Id_A : N \tw A \ra N' \tw A$ compatible aux diffÈrentielles.
Nous obtenons ainsi un {\em foncteur}  \indexnotation{LN}
\[
L_\tau : \coModu C \ra \Modu A, \hspace{1cm} N \ra N \tw A
\]
que nous noterons $L$ lorsqu'il n'y aura pas d'ambiguitÈ.

\begin{lemme}
Le foncteur $L : \coModu C \ra \Modu A$ est adjoint ‡ gauche au foncteur $R : \Modu A \ra \coModu C.$
\end{lemme}

\dem  Soit $N$ un objet de $\coModu C$ et $M$ un objet de $\Modu A$. Nous allons donner la bijection
fonctorielle
\[
\phi : \Hom_{\sz \Modu A} (LN,M) \arr{} \Hom_{\sz \coModu C}(N,RM) .
\]
Soit $f : LN \ra M$ un morphisme de $\Modu A$. Par le lemme \ref{lemme_modules_libres}, il est
dÈterminÈ par sa composition $\alpha = f \circ (\Id_N \ts \unite^A) : N \ra M.$ Par le lemme
\ref{lemme_comodules_colibres}, le morphisme $\alpha$ dÈtermine ‡ son tour un morphisme graduÈ
de $C$-comodules co-unitaires
$\phi (f) : N \ra RM$ tel que $(\Id \ts \counite^C) \phi(f) = \alpha.$ On vÈrifie que la condition 
$b_\tau \phi (f) - \phi  (f)  d_N = 0$ Èquivaut ‡ la condition $d_M f - f D_\tau = 0.$ \findem

\begin{definition}{\em 
Une cochaÓne tordante admissible $\tau : C \ra A$ est {\em acyclique} \index{acyclique (cochaÓne)}
\index{cochaÓne tordante!acyclique}si, 
pour tout objet $M$ de $\Modu A$, le morphisme d'adjonction
\[
\phi : LR M \ra M
\]
est un quasi-isomorphisme (voir la proposition \ref{caracterisation_cochaines_acycliques}
ci-dessous pour des conditions Èquivalentes).
}\end{definition}

\begin{notation}[Construction bar et cobar]{\em
Soit $A$ un objet de $\alga$. Nous notons $\Ba A$ la cogËbre co-augmentÈe
$(B\b A)\+$, o˘ $\b A$ est l'algËbre rÈduite associÈe ‡ $A$. Attention ‡ ne pas
confondre les cogËbres co-augmentÈes $\Ba A$ et $(BA)\+.$ \index{construction bar}
\indexnotation{Ba1}
Soit $C$ un objet de $\cocoga$. Nous notons $\Oma C$ l'algËbre augmentÈe
$(\Omega \b C)\+$, o˘ $\b C$ est le cogËbre rÈduite associÈe ‡ $C$. Elle n'est pas
isomorphe ‡ $(\Omega C)\+$. \index{construction cobar} \indexnotation{Oma}
\indexnotation{Oma1}
}\end{notation}

\begin{lemme}
\label{cochaines_universelles} 
\english \begin{itemize}
\item[a.] Soit $A$ un objet de $\alga$. Soit $p : B\b A \ra S\b A$ la projection
canonique. La composition
\[ 
\tau_A : \Ba A \ra B\b A \arr{\si \circ p} \b A \ra A, 
\]
o˘ la premiËre flËche est la projection canonique et la derniËre l'injection canonique,
est une cochaÓne tordante admissible. La cochaÓne
$\tau_A$ est universelle parmi les cochaÓnes tordantes admissibles
de but $A$, i.~e.~si $C$ est un objet de $\coga$ et $\tau : C \ra A$ est une cochaÓne tordante admissible, il
existe un unique morphisme $g_\tau$ tel que $\tau_A \circ g_\tau = \tau.$

\item[b.] De faÁon duale, nous associons ‡ un objet $C$ de $\cocoga$ une cochaÓne tordante admissible
\[ 
\tau_C : C \ra \b C \arr{i \circ \si} \Omega \b C \ra \Oma C
\]
o˘ $i : S^{-1}\b C \ra \Omega \b C$ est l'injection canonique.
La cochaÓne
$\tau_C$ est universelle parmi les cochaÓnes tordantes admissibles
de source $C$, i.~e.~si $\tau : C \ra A$ est une cochaÓne tordante admissible, il existe un unique morphisme
$f_\tau$ tel que $f_\tau \circ \tau_C = \tau.$

\end{itemize} \francais
\end{lemme}

\dem 
Soit $C$ un objet de $\cocoga$ et $A$ un objet de $\alga$. 
Soit $\tau : C \ra A$ un morphisme graduÈ de degrÈ $+1$ dont la composition avec la co-augmentation
de $C$ et l'augmentation de $A$ est nulle.
Soit 
\[
f_\tau : \Oma C \ra A
\]
le morphisme graduÈ d'algËbres augmentÈes relevant (\ref{algebres_libres})
la composition $\tau \circ s$ et 
\[
g_\tau : C \ra \Ba A
\]
le morphisme graduÈ de cogËbres co-augmentÈes relevant
(\ref{cogebres_tensorielles}) la composition $s \circ \tau$. 
Par la dÈmonstration du lemme \ref{adjonction_bar_cobar}, le morphisme graduÈ $\tau$ est
une cochaÓne tordante si et seulement si $f_\tau$ est compatible aux diffÈrentielles
si et seulement si $g_\tau$ est compatible aux diffÈrentielles.

{\it a.} La composition $\si \circ p : B\b A \ra \b A$ est une cochaÓne tordante car le relËvement
(\ref{cogebres_tensorielles}) de $p : B \b A \ra S\b A$
est l'identitÈ de la cogËbre $B\b A$ (et cette derniËre
commute Èvidemment ‡ la diffÈrentielle de $B \b A$). L'universalitÈ est immÈdiate. 

{\it b.} Idem.
\findem

\begin{definition}{\em 
On appelle $\tau_A$ la {\em cochaÓne tordante universelle de $A$} \index{universelle (cochaÓne)}
\index{cochaÓne tordante!universelle}
\indexnotation{tauA} et $\tau_C$
la {\em cochaÓne tordante universelle de $C$.}
}\end{definition}

\begin{remarque}{\em
Dans \cite{Husemoller74}, le foncteur
\[
R_{\tau_A} : \Modu A \ra \coModcu \Ba A, \hspace{1cm} M \mapsto M \ts_{\tau_A} \Ba A,
\]
est notÈ $B_A M.$ \findem }\end{remarque}

\noindent Notons
\[
\Res : \Modu A \ra \Modu \Oma C
\] 
le foncteur restriction le long de $f_\tau$ et 
\[
\Ind :  \Modu \Oma C \ra \Modu A
\]
le foncteur induction.
On sait que $(\Ind , \Res)$ est une paire de foncteurs adjoints de la catÈgorie
$\Modu \Oma C$ vers la catÈgorie $\Modu A$. 
Notons
\[
\Res^{op} : \coModcu C \ra \coModcu \Ba A
\] 
le foncteur corestriction le long de $g_\tau$ et 
\[
\Ind^{op} :  \coModcu \Ba A \ra \coModcu C
\]
le foncteur co-induction.
On sait que $(\Res^{op} , \Ind^{op})$ est une paire de foncteurs adjoints de la catÈgorie
$\coModcu C$ vers la catÈgorie $\coModcu \Ba A$.

\begin{lemme} \label{lemme_composition_adjonction}
\english \begin{itemize}
\item[a.] La paire de foncteurs adjoints $(L_\tau, R_\tau)$ de la catÈgorie $\Modu A$ vers la
catÈgorie  $\coModcu C$ est
la composition de la paire $(\Ind , \Res)$ avec la paire
$(L_{\tau_C}, R_{\tau_C}).$
\item[b.] La paire de foncteurs adjoints $(L_\tau, R_\tau)$ de la catÈgorie $\Modu A$ vers la
catÈgorie  $\coModcu C$ est
la composition de la paire $(L_{\tau_A}, R_{\tau_A})$ avec la paire $(\Res^{op} , \Ind^{op})$.
\findem
\end{itemize} \francais
\end{lemme}

\begin{lemme} \label{lemme_acyclicite_cochaines_universelles} 
\english \begin{itemize}
\item[a.] Soit $A$ un objet de $\alga$. La cochaÓne tordante universelle $\tau_A$ est
acyclique.
\item[b.] Soit $C$ un objet de $\cocoga$. La cochaÓne tordante universelle $\tau_C$
est acyclique.
\end{itemize} \francais
\end{lemme}

\dem 
{\it a.} Soit $M$ un objet de $\Modu A$.
Montrons que $LRM = (M \ts \Ba A \ts A, d)$ est une rÈsolution (dite rÈsolution bar normalisÈe)
de $M$ 
\[
bar_A(M) = \cdots \ra M \ts \b A \tp i \ts A \ra \cdots \ra M \ts \b A \ts A \ra M\ts A,
\]
et que le morphisme $\Phi$ correspondant au morphisme $bar_A(M) \ra M$ est un quasi-isomorphisme.
Comme dans le cas o˘ $M$ est concentrÈ en degrÈ $0,$ (voir \cite[IX.6]{Cartan56} o˘ ce
complexe se nomme le complexe standard normalisÈ) les morphismes
\[
h_{i-1} = \Id \tp i \ts p \ts \epsilon : M \ts \b A \tp {i-1} \ts A \ra M \ts \b A \tp i \ts A,
\]
o˘ $p$ est la projection canonique, dÈfinissent une homotopie contractante du complexe
\[
\cdots \ra M \ts \b A \tp i \ts A \ra \cdots \ra M \ts \b A \ts A \ra M\ts A \ra M \ra 0.
\]

{\it b.} Soit $M$ un objet de $\Modu{\Oma C}$.
Montrons que $\Phi : LRM \ra M$ est un quasi-isomorphisme filtrÈ.  Munissons $\Omega \b C$ de la filtration
induite par la filtration primitive de $\b C$ considÈrÈ comme cogËbre. Nous avons alors une
filtration de $\Oma C$ dÈfinie par la suite
\[
\left(\Oma C\right)_i = (\Omega \b C)_i \oplus e, \quad i\geq 0. 
\]
Munissons $C$, considÈrÈ comme objet de $\coModu C$, de sa filtration primitive
de $C$-module (nous la complÈtons par $C\prim 0 = e$). 
Munissons $M$ de la filtration dÈfinie par la suite $M_i = M$, $i \geq 0.$
Ces filtrations induisent sur $LRM = (M \ts C \ts \Oma C)$ une filtration
de complexes.
Le morphisme $\Phi : LRM \ra M$ devient un morphisme filtrÈ pour ces filtrations.
Il induit un morphisme
\[
\ogr_0 (LRM) \ra \ogr_0 M
\]
qui est l'identitÈ de $M.$
Comme $\ogr_i M = 0$ pour tout $i \geq 1$, il nous suffit de montrer que
\[
 \ogr_i (LRM), \quad i\geq 1,
\]
est contractile. Soit $i \geq 1.$ Par construction, nous avons un isomorphisme d'objets graduÈs
\[
 \ogr_i (LRM) =  M \ts e \ts \ogr_i\Oma C \oplus 
\bigg(\bigoplus_{i_1 + i_2 =i\above 0pt i_1 \neq 0}
M \ts \ogr_{i_1}C \ts \ogr_{i_2}\Oma C\bigg).
\]
La diffÈrentielle a pour matrice
\[
\left[\begin{array}{cc}
0 &  \rho \\
0 & 0
\end{array}\right]
\] 
o˘ $\rho$ est le morphisme induit par $T_{\tau_C}$
\[
\bigoplus_{i_1 + i_2 =i\above 0pt
i_1 \neq 0} M \ts \ogr_{i_1}C \ts \ogr_{i_2}\Oma C 
\arr{} M \ts e \ts \ogr_i\Oma C.
\]
Ce dernier est un isomorphisme car il est induit par l'isomorphisme
\[
\bigoplus_{i_1 + i_2 =i\above 0pt
i_1 \neq 0} \ogr_{i_1}C \ts \ogr_{i_2}\Oma C 
\arr{} \ogr_i\Oma C
\]
\findem

%
%
%
%
%
%
%
%
%
%
%
%
%

\subsection{$\coModcu C$ comme catÈgorie de modËles}

Soit $C$ un objet de $\cocoga$. Dans cette section, nous allons munir $\coModcu C$
d'une structure de catÈgorie de modËles. Nous commenÁons par rappeler la structure
de catÈgorie de modËles sur $\Modu A$, o˘ $A$ est un objet de $\alga$ et nous
ÈnonÁons ensuite le thÈorËme principal (\ref{theoreme_cmf_coModcu}). Nous ne dÈtaillerons
pas toute sa dÈmonstration car elle est similaire ‡ celle de (\ref{theoreme_cmf_cocog}).
Seul les points qui diffËrent seront dÈveloppÈs.\\

%
%
%
%
%
%
%
%
%
%
%
%
%

{\noindent \bf Rappels sur la catÈgorie $\Modu A$}\\

Soit $A$ une algËbre diffÈrentielle graduÈe unitaire.
Dans la catÈgorie $\Modu A$,
considÈrons les trois classes de morphismes suivantes
\english \begin{itemize}
\item[-] la classe $Qis$ des quasi-isomorphismes,
\item[-] la classe $\fib$ des morphismes $f : M \ra M'$ tels que
$f^n$ est un Èpimorphis\-me pour tout $n \in \Z,$
\item[-] la classe $\cof$ des morphismes qui ont la propriÈtÈ de relËvement
‡ gauche par rapport aux morphismes appartenant ‡ $Qis \cap \fib$.
\end{itemize} \francais

\begin{theoreme}[Hinich \cite{Hinich97c}]\label{theoreme_cmf_Modu}
La catÈgorie $\Modu A$ munie des classes de morphismes dÈfinies ci-dessus
est une catÈgorie de modËles. 
Tous les objets sont fibrants. Les objet cofibrants sont dÈcrits dans
la remarque \ref{remarque_objets_cofibrants_Modu} ci-dessous.
\end{theoreme}


{\noindent \bf Le thÈorËme principal}\\

Soit $A$ un objet de $\alga$ et $C$ un objet de $\cocoga$.
Soit $\tau : C \ra A$ une cochaÓne tordante admissible acyclique.
Dans la catÈgorie $\coModcu C$ des comodules
diffÈrentiels graduÈs co-unitaires cocomplets,
nous considÈrons les trois classes de morphismes suivantes :

\english \begin{itemize}
\item[-] la classe $\weq $ des {\em Èquivalences
faibles} est formÈe des morphismes $f : N\ra N'$ tels que
 $L f : LN \ra LN'$ est un quasi-isomorphisme 
de modules,

\item[-] la classe $\cof$ des {\em cofibrations} est formÈe des morphismes
$f : N \ra N'$ qui, en tant que morphismes de complexes, sont des
monomorphismes,

\item[-] la classe $\fib$ des {\em fibrations} est formÈe des morphismes
qui ont la propriÈtÈ de relËvement ‡ droite par rapport aux cofibrations triviales.
\end{itemize} \francais

\begin{theoreme} \label{theoreme_cmf_coModcu}
\english \begin{itemize}
\item[a.]
La catÈgorie $\coModcu C$ munie des trois classes de morphismes ci-dessus est une catÈgorie de modËles.
Tous ses objets sont cofibrants. Un objet de $\coModcu C$ est fibrant si et seulement si il est un facteur
direct d'un objet $RM$, o˘ $M$ est un objet de $\Modu A$.
\item[b.]
Munissons la catÈgorie $\Modu A$ de la structure de catÈgorie de
modËles du thÈorËme \ref{theoreme_cmf_Modu}.
La paire de foncteurs adjoints $(L, R)$ de $\coModcu C$ dans $\Modu A$ est une
Èquivalence de Quillen.
\item[c.] La structure de catÈgorie de modËles sur $\coModcu C$ ne dÈpend pas de la cochaÓne
tordante admissible acyclique $\tau$.
\end{itemize} \francais
\end{theoreme}

En particulier, la catÈgorie $\Ho \coModcu C$ est Èquivalente ‡ la catÈgorie dÈrivÈe
$\cd A$ (voir la dÈfinition dans \ref{structure_triangulee_Ho_coModcu}).
Le thÈorËme \ref{theoreme_cmf_coModcu} et le lemme \ref{lemme_acyclicite_cochaines_universelles} implique
le corollaire suivant :

\begin{corollaire} \label{corollaire_theoreme_cmf_coModcu}
La catÈgorie $\coModcu C$ admet une unique structure de catÈgorie de modËles
telle que pour toute cochaÓne tordante admissible acyclique $\tau : C \ra A$, o˘ $A$
est un objet de $\alga$, le couple de foncteurs adjoints $(L,R)$ est une Èquivalence
de Quillen.
\findem
\end{corollaire}

\begin{definition}{\em \label{definition_cmf_canonique_coModcu}
Nous appelons la structure de catÈgorie de modËles sur $\coModcu C$ du
corollaire la {\em structure  canonique}.
}\end{definition}

Pour montrer le thÈorËme \ref{theoreme_cmf_coModcu}, nous avons besoin 
(comme pour la dÈmonstration du thÈorËme \ref{theoreme_cmf_cocog}) d'introduire des
filtrations.\\

Si l'algËbre $A$ (resp.~la cogËbre $C$) est filtrÈe,
un {\em $A$-module diffÈrentiel graduÈ filtrÈ} (resp.~{\em $C$-comodule diffÈrentiel graduÈ filtrÈ})
est un $A$-module (resp.~$C$-comodule) dans la catÈgorie des complexes filtrÈs.
Un $C$-comodule filtrÈ $M$ est {\em admissible} \index{admissible!comodule}
\index{comodule!admissible}
si sa filtration est exhaustive et si $M_0 = 0.$ Par dÈfinition, tous les objets de $\coModcu C$,
munis de leur filtration primitive sont admissibles.

\begin{lemme} \label{lemme1_cmf_coModcu}
Si $C$ est munie d'une filtration exhaustive de cogËbres telle que $C_0 = e$,
un quasi-isomorphisme filtrÈ de $C$-comodules admissibles est une Èquivalence faible.
\end{lemme}

\dem  Soit $f : N \ra N'$ un quasi-isomorphisme filtrÈ de $C$-comodules admissibles.
La filtration de $N$ induit une filtration de $A$-module dÈfinie par la suite
\[
(LN)_i = N_i \ts A, \quad i \geq 0.
\]
La diffÈrentielle de $(LN)_i$, $i\geq 0$, est la somme de la diffÈrentielle
du produit tensoriel $N_i \ts A$ et
de la contribution de $D_\tau$. Comme la filtration de $N$ est admissible
et que la cochaÓne $\tau : C \ra A$ est admissible, la contribution de $D_\tau$
fait dÈcroÓtre la filtration de $LN.$ Ainsi, la diffÈrentielle de
\[
\ogr LN \arr{\sim} \ogr N \ts A
\]
est celle du produit tensoriel $\ogr N \ts A$ et le morphisme $Lf$ est bien
un quasi-isomorphisme de $A$-modules.  \findem

\begin{lemme} \label{lemme0_cmf_coModcu}
\english \begin{itemize}
\item[a.] Soit $M$ et $M'$ deux objets de $\Modu A$. Le foncteur $R$ envoie un quasi-isomorphisme
$f : M \ra M'$ sur une Èquivalence faible $Rf : RM \ra RM'$ dans $\coModcu C$.
\item[b.] Soit $M$ un objet de $\Modu A$. Le morphisme d'adjonction
\[
\Phi : LR M \arr{} M
\]
est un quasi-isomorphisme de $A$-modules.
\item[c.] Soit $N$ un objet de $\coModcu C$. Le morphisme d'adjonction
\[
\Psi : N \arr{} RL N
\]
est une Èquivalence faible de $\coModcu C$.
\end{itemize} \francais
\end{lemme}

\dem 

{\it b.} La cochaÓne $\tau $ est acyclique.

{\it a.} Le morphisme $Rf$ est une Èquivalence faible si et seulement si $LRf$ est un quasi-isomorphisme.
Par le point {\it b}, $\Phi$ est un quasi-isomorphisme. Par ailleurs, on~a 
\[
\Phi_M \circ f = LRf \circ \Phi_{M'}.
\]
La saturation des quasi-isomorphismes dans $\Modu A$ nous donne le rÈsultat.

{\it c.} Nous voulons montrer que $\Psi $ est une Èquivalence faible, c'est-‡-dire que $L\Psi : LN \ra LRLN$
est un quasi-isomorphisme. Nous savons que 
\[
\Phi_{LN} \circ L\Psi_N = \Id_{LN}
\]
et que $\Phi$ est un
quasi-isomorphisme. Le morphisme $L \Psi$ est donc aussi un quasi-isomorphisme. 
\findem\\

Rappelons la description de \cite{Hinich97c} des cofibrations de $\Modu A$.
Les {\em cofibrations standard} (resp.~{\em triviales)} de $\Modu A$
\index{cofibration standard (de $\Modu A$)} sont dÈfinies
comme dans la dÈfinition \ref{cofibration_standard},
‡ la diffÈrence prËs que $M^\sharp$ dÈsigne le complexe
sous-jacent ‡  un objet $M$ de $\Modu A$ et que $FV$ dÈsigne le module diffÈrentiel graduÈ libre sur 
un complexe $V.$ Nous avons alors la mÍme description (voir juste dessous \ref{cofibration_standard}) des
 cofibrations (resp.~triviales) de $\Modu A$ ‡ partir des cofibrations standard
(resp.~triviales).

\begin{lemme} \label{lemme2_cmf_coModcu}
Soit $N$ un objet de $\coModcu C$ et $N'$ un sous-objet de $N$ tel que
$\Delta N \subset N \ts e \oplus N' \ts  C.$
Le foncteur $L$ envoie l'inclusion $N' \hookrightarrow N$ sur une cofibration
standard.
\end{lemme}

\dem 
Soit $E$ le conoyau de l'inclusion $N' \hookrightarrow N$. Choisissons un scindage
dans la catÈgorie des objets graduÈs
\[
N \arr{\sim} N' \oplus E.
\]
Selon cette dÈcomposition, la comultiplication $\Delta^N$ est donnÈe par deux composantes
\[
\Delta^{N'} : N' \ra N' \ts C \quad \mbox{et} \quad \Delta^E = \left[\begin{array}{c}\Delta^E_1
\\ \Delta^E_2 \end{array}\right]
: E \arr{}  N \ts e \oplus N' \ts C,
\]
et la diffÈrentielle est donnÈe par la diffÈrentielle de $N'$, celle de $E$ et un morphisme
\[
d' : E \arr{} N'.
\]
Nous avons un isomorphisme d'objets graduÈs
\[
LN \arr{\sim} LN' \oplus LE.
\]
La diffÈrentielle est la somme de celle de $LN' \oplus LE$, du morphisme 
\[
d' \ts \Id : E \ts A \ra N' \ts A
\]
et du morphisme $d_\tau'$ qui est la composition
\[
E \ts A \arr{\Delta^E_2 \ts \Id} N' \ts C \ts A \arr{\Id \ts \tau \ts \Id }
N' \ts A \ts A \arr{\Id \ts \mu^A} N' \ts A.
\]
Remarquons qu'il n'y a pas de contribution de $\Delta^E_1$ car la cochaÓne $\tau$
est admissible. Posons 
\[D' = (d' \ts \Id + d'_\tau)s : S^{-1}E \ra N' \ts A.
\]
Nous vÈrifions que 
$LN$ est isomorphe ‡
\[
LN' \langle S^{-1}E , D' \rangle.
\]
\findem

\begin{lemme} \label{lemme3_cmf_coModcu}
\english \begin{itemize}
\item[a.] Le foncteur $L$ prÈserve les cofibrations et les
Èqui\-valences faibles.
\item[b.] Le foncteur $R$ prÈserve les fibrations et les
Èquivalences faibles.
\end{itemize} \francais
\end{lemme}

\dem 

{\it a.} Soit $j : N' \rightarrowtail N$ une cofibration de $\coModcu C$. Soit
la filtration de $N$ donnÈe par la suite
\[
N_i = j(N') + N\prim i, \quad i\geq 0,
\]
o˘ $N\prim i$, $i\geq 1$, est la filtration primitive de $N$ (complÈtÈe par $N_0 = 0$).
Remarquons que, pour tout $i \geq 1$, nous avons
\[
\Delta N_i \subset N_i \ts e \oplus N_{i-1} \ts C.
\]
Nous pouvons donc appliquer le lemme \ref{lemme2_cmf_coModcu}. Il certifie que $LN_i \ra LN_{i+1}$
est une cofibration standard. Le morphisme $Lj : LN' \ra LN$ est ainsi la composition dÈnombrable
des cofibrations standard $LN_i \ra LN_{i+1}$, il est donc une cofibration.
Par dÈfinition des Èquivalences faibles dans $\coModcu C$, le foncteur $L$ prÈserve les
Èquivalences faibles.

{\it b.} Par le point {\it a} et l'adjonction $(L,R,\phi)$ de $\coModcu C$ dans $\Modu A$, le
foncteur $R$ conserve les fibrations. Le fait qu'il conserve les Èquivalences faibles est le
point {\it a} du lemme \ref{lemme0_cmf_coModcu}. \findem

\begin{lemme} \label{lemmecle_cmf_coModcu}
Soit $M$ un objet de $\Modu A$ et $N$ un objet de $\coModcu C$.
Soit une fibration $p : M \twoheadrightarrow L N$ de $\Modu A$.
Le morphisme  $j : RM \prod_{RL N}N \ra RM $ de comodules du diagramme cartÈsien 
\[
\xymatrix{ RM \prod_{RLN} N
\ar[d]_{j} \ar@{->}[r] & N \ar[d]^{\Psi}\\ 
  RM  \ar@{}[ur]|{cart.}\ar@{->}[r]_{Rp} & RLN.}
\]
est une cofibration triviale de $\coModcu C$. 
\end{lemme}

\dem Soit $K$ le noyau de $p$. Nous avons des isomorphismes d'objets graduÈs
\[
RM \arr{\sim} RK \oplus RLN, \quad RM \prod_{RLN} N \arr{\sim} RK \oplus N.
\]
Le morphisme $j$ s'Ècrit alors
\[
\left[\begin{array}{cc} \Id & * \\ 0 & \Psi \end{array}\right].
\]
Nous avons donc un diagramme de $\Modu A$
\[
\xymatrix{
0 \ar[r] & LRK \ar[r] \ar[d]_{L\Id} &L\Big( RM \prod_{RLN} N \Big)  \ar[d]^{Lj} \ar[r] & LN \ar[d]^{L\Psi}
\ar[r] & 0\\
0 \ar[r] & LRK \ar[r] & RM  \ar[r] & LRLN \ar[r] & 0,
}
\]
o˘ les lignes sont exactes et o˘ la flËche verticale de droite et celle de gauche sont des
quasi-isomorphismes. Le morphisme $Lj$ est donc un quasi-isomorphisme, et $j$ est une Èquivalence
faible de $\coModcu C$. Il est clairement un monomorphisme, donc une cofibration de $\coModcu C$.
\findem\\

{\noindent \bf DÈmonstration du thÈorËme \ref{theoreme_cmf_coModcu}}\\

Par les lemmes ci-dessus, 
la dÈmonstration du fait que les classes $\weq$, $\cof$ et $\fib$ dÈfinissent une
structure de catÈgorie de modËles  est la mÍme que celle du thÈorËme \ref{theoreme_cmf_cocog}.\\

{\noindent \bf Objets cofibrants et objets fibrants de $\coModcu C$}\\

Tous les objets de $\coModcu C$ sont cofibrants puisque les cofibrations sont les
monomorphismes.

Montrons qu'un objet de $\coModcu C$ est fibrant si et seulement si il est un facteur
direct d'un objet $RM$, o˘ $M$ est un objet de $\Modu A$. Nous rappelons (\ref{theoreme_cmf_Modu})
que tous les objets de $\Modu A$ sont fibrants.
Par le lemme \ref{lemme3_cmf_coModcu}, l'image du foncteur $R$ est donc formÈ d'objets
fibrants de $\coModcu C.$ Ainsi, tous les objets de la forme $RM$ et leurs facteurs
directs sont fibrants. 
RÈciproquement si $N$ est fibrant, par l'axiome (CM4), le
morphisme $\Psi : N \ra RLN$ (qui est une cofibration triviale) est scindÈ. 
L'objet $N$ est donc un facteur direct de $RLN$
\begin{remarque} \label{remarque_objets_cofibrants_Modu}{\em
La dualisation de cette dÈmonstration montre que les objets cofibrants de
$\Modu A$ sont les facteurs directs des $LN$, $N \in \coModcu C$.}
\end{remarque}
 
Le point {\it b} du thÈorËme \ref{theoreme_cmf_coModcu}
est un corollaire du lemme \ref{lemme1_cmf_coModcu}. Il nous reste ‡ montrer le point
{\it c.}\\

{\noindent \bf UnicitÈ de la structure de catÈgorie de modËles sur $\coModcu C$}\\

Soit $A'$ un objet de $\alga.$ Soit 
$\tau' : A' \ra C$ une cochaÓne tordante admissible acyclique. Nous voulons montrer
que la structure de
catÈgorie de modËles sur $\coModcu C$ (dÈfinie au point {\it a} de \ref{theoreme_cmf_coModcu})
relative ‡ $\tau$ est la mÍme que celle relative ‡ $\tau'.$

Il suffit de le montrer dans le cas o˘ $\tau'$ est la cochaÓne universelle $\tau_C.$
Nous allons montrer que les classes des cofibrations et les classes des Èquivalences
faibles relatives aux deux structures coÔncident.
C'est vrai pour les cofibrations puisqu'elles sont les monomorphismes.
Nous rappelons (\ref{lemme_composition_adjonction}) que 
la paire de foncteurs adjoints $(L_\tau, R_\tau)$ de $\Modu A$ vers $\coModcu C$ est
la composition de la paire $(\Ind , \Res)$ avec la paire
$(L_{\tau_C}, R_{\tau_C}).$ 
Comme le foncteur $\Res$  induit une
Èquivalence entre les localisations
de $\Modu A$ et $\Modu \Oma C$ par rapport aux quasi-isomorphismes
(voir \cite[exple 6.1]{Keller94}), 
les Èquivalences faibles des deux
structures sur $\coModcu C$ coÔncident par le point {\it b} du thÈorËme \ref{theoreme_cmf_coModcu}.
\findem \\

{\noindent \bf Quasi-isomorphismes filtrÈs et Èquivalences faibles}\\

Nous notons $\mbox{\it Qisf\,}$ la classe des morphismes $f : N \ra N'$ tels que
$C$ admet une filtration exhaustive de cogËbre telle que $C_0 = e$ et tels que $N$ et $N'$
admettent des filtrations admissibles de $C$-comodules pour lesquelles $f$ est un quasi-isomorphisme
filtrÈ. Le lemme \ref{lemme1_cmf_coModcu} montre que nous avons une inclusion
\[
\mbox{\it Qisf\,} \subset \weq.
\]
On rappelle (voir appendice \ref{section_rappels_cmf}) que la {\em catÈgorie homotopique} $\Ho \coModcu C$
est la localisation
\[
\Big(\coModcu C\Big)[\weq^{-1}].
\]

\begin{lemme} \label{lemme_Qisf_weq}
Le foncteur canonique
\[
\Big(\coModcu C\Big) [\mbox{\it Qisf\,}^{-1}] \arr{\sim} \Ho \coModcu C
\]
est une Èquivalence.
\end{lemme}

\dem La dÈmonstration est similaire ‡ celle du point {\it a} de la proposition \ref{proposition_equiv_qis}.
Nous vÈrifions que le morphisme d'adjonction 
\[
\Psi : N \ra R_{\tau_C} L_{\tau_C}N
\]
est
un morphisme quasi-isomorphisme filtrÈ pour la filtration primitive sur $N$ et la filtration
sur $R_{\tau_C} L_{\tau_C}N$ induite par les filtrations primitives de $N$ et $C$.
Le morphisme $R_{\tau_C} L_{\tau_C}f$
est clairement un quasi-isomorphisme filtrÈ. La propriÈtÈ
de saturation des quasi-isomorphismes filtrÈs appliquÈe ‡
l'ÈgalitÈ $RLf \circ \Psi_N = \Psi_{N'} \circ f$ nous donne le rÈsultat.  \findem

%
%
%
%
%
%
%
%
%
%
%
%
%

\newpage

\subsection{Structure triangulÈe sur $\Ho\coModcu C$} \label{structure_triangulee_Ho_coModcu}

{\bf \noindent Rappel sur la structure triangulÈe sur $\Ho \Modu A$}\\

Rappelons qu'une catÈgorie de Frobenius \index{categorie@{catÈgorie}!de Frobenius}
est une catÈgorie exacte au
sens de Quillen \cite{Quillen73} qui possËde assez d'injectifs et assez
de projectifs et dont la classe des projectifs coÔncide avec celle des
injectifs. Il est connu \cite{Heller60}, \cite{Happel87}, \cite{Keller87}
que le quotient d'une catÈgorie de Frobenius $\ca$ par l'idÈal des morphismes
se factorisant par un projectif est une catÈgorie triangulÈe \cite{Verdier77}.
On l'appelle la {\em catÈgorie stable} associÈe ‡ $\ca.$ 
\index{categorie@{catÈgorie}!stable}

Soit $A$ une algËbre diffÈrentielle graduÈe unitaire. La catÈgorie $\Modu A$,
munie de la {\em classe} $\ce$ formÈe des suites exactes
\[
0 \ra M' \arr{f} M \arr{g} M'' \ra 0
\]
qui sont scindÈes dans la catÈgorie de modules graduÈs, est une catÈgorie
exacte.
La classe des objets injectifs est formÈe des complexes de la forme
\[
IM = \Big(M\oplus SM,\left[\begin{array}{cc}0 & \si \\ 0 & 0\end{array}\right]\Big),
\quad M \in \Modu A.
\]
Elle coÔncide avec la classe des objets projectifs. La catÈgorie $\Modu A$ est donc
une catÈgorie de Frobenius.
Nous notons $\ch A$ la catÈgorie stable associÈe ‡ $\Modu A.$
Elle est une catÈgorie triangulÈe. Son foncteur suspension est le foncteur $M \mapsto SM$.
Ses triangles standard proviennent des suites exactes de $\ce$.
Les quasi-isomorphismes de $\Modu A$ sont
exactement les morphismes $f$ dont l'image $\b f$ par le foncteur canonique $\Modu A \ra \ch A$
s'insËre dans un triangle
\[
N \ra M \arr{\b f} M' \ra SN,
\]
o˘ $N$ est acyclique.
La {\em catÈgorie dÈrivÈe} $\cd A$ \index{categorie@{catÈgorie}!derivee@{dÈrivÈe}}
\indexnotation{cdA}
est la localisation de la catÈgorie $\ch A$
par rapport aux quasi-isomorphismes.
Les {\em triangles standard} de $\cd A$ sont
l'image par le foncteur
\[
Q : \ch A\arr{} \cd A
\]
des triangles standard de $\ch A.$
La catÈgorie dÈrivÈe $\cd A$, munie de l'endofoncteur suspension est triangulÈe pour la
classe des {\em triangles distinguÈs}, i.~e.~les triangles isomorphes ‡ des triangles standard.
Si $f$ est un morphisme de
$\Modu A$, on note $C(f)$ \indexnotation{Cf} son cÙne. Si
\[
0 \ra M' \arr{i} M \arr{p} M'' \ra 0
\]
est une suite exacte (non nÈcessairement scindÈe) de $\Modu A$,
le morphisme $[p,0] : C(i) \ra M''$ est un quasi-isomorphisme et la
suite
\[
M' \arr{Q \b i} M \arr{Q \b p} M'' \arr{\delta} SM',
\] 
o˘ le morphisme $\delta$ est le morphisme de $\cd A$ dÈfini par
\[
M'' \larr{[p,0]} C(i) \arr{[-1,0]} SM',
\]
est un triangle distinguÈ de $\cd A.$
\\

{\noindent \bf Structure triangulÈe sur $\Ho \coModcu C$}\\

Soit $C$ un objet de $\cocoga.$ La catÈgorie $\coModcu C$, munie de la classe $\cf$
des suites exactes courtes qui sont scindÈes dans la catÈgorie des comodules
graduÈs, est une catÈgorie de Frobenius dont la classe des objets injectifs
est formÈe des objets
\[
IN = \Big(N\oplus SN,\left[\begin{array}{cc}0 & \si \\ 0 & 0\end{array}\right]\Big),
\quad N \in \coModcu C.
\]
Nous notons $\ch C$ la catÈgorie stable associÈe.
Elle est triangulÈe. Son foncteur suspension est $N \mapsto SN$.
Les suites exactes de $\cf$ donnent lieu aux {\em triangles standard}. Les triangles
distinguÈs sont les triangles isomorphes aux triangles standard.

Soit $\tau : C \ra A$ une
cochaÓne tordante admissible acyclique o˘ $A$ est un objet de $\alga$.
Les foncteurs $L$ et $R$ forment un couple de foncteurs exacts entre
les catÈgories $\coModcu C$ et $\Modu A$ et prÈservent l'injectivitÈ.
Ils induisent donc un couple de foncteurs adjoints triangulÈs entre
les catÈgories stables $\ch C$ et $\ch A.$ La {\em catÈgorie dÈrivÈe}
$\cd C$ \indexnotation{cdC}
est la catÈgorie localisÈe $\big(\ch C\big)[\weq^{-1}].$ Elle est clairement
isomorphe ‡ la catÈgorie $\Ho \coModcu C$.
Rappelons (thm \ref{theoreme_cmf_coModcu}) que les foncteurs $R$ et $L$
(dÈfinis en \ref{section_cochaines_tordantes})
induisent des Èquivalences
inverses l'une de l'autre entre les catÈgories localisÈes
\[
\cd A = \big(\ch A\big)[Qis^{-1}] \quad \mbox{et} \quad \big(\ch C\big)[\weq^{-1}] = \cd C.
\]
En particulier, le systËme multiplicatif $\weq$ est compatible aux triangles
de $\ch C$ car il est l'image rÈciproque du systËme multiplicatif des
isomorphismes de $\cd A$ par le foncteur triangulÈ composÈ
\[
\ch C \arr{L} \ch A \arr{} \cd A.
\]
Il en rÈsulte que $\cd C$ porte une structure triangulÈe canonique et que les
Èquivalences induites entre $\cd A$ et $\cd C$ sont des foncteurs triangulÈs.

\subsection{CaractÈrisation de l'acyclicitÈ des cochaÓnes tordantes}

Nous rappelons que le foncteur $\?\+ : \cocog \ra \cocoga$ est une Èquivalence de catÈgories
(\ref{cocoga}).
Munissons $\cocoga$ de la structure de catÈgorie de modËles induite par celle de
$\cocog$ (voir \ref{theoreme_cmf_cocog}).

\begin{proposition} \label{caracterisation_cochaines_acycliques}
Soit $A$ un objet de $\alga$ et $C$ un objet de $\cocoga$. Soit $\tau : C \ra A$
une cochaÓne tordante admissible. Les conditions suivantes sont Èquivalentes.
\english \begin{itemize}
\item[a.] La cochaÓne tordante $\tau$ est acyclique, i.~e.~si $M$ est un objet de $\Modu A$, le morphisme d'adjonction 
\[
\Phi : LRM \ra M
\] est un quasi-isomorphisme de $\Modu A$.
\item[b.] Si $N$ est un objet de $\coModcu C$, le morphisme d'adjonction 
\[
\Psi : N \ra RLN
\]
est une Èquivalence faible de $\coModcu C$.
\item[c.] Le morphisme d'adjonction
\[
LRA = A \tw C \tw A \arr{\Phi_A} A 
\]
est un quasi-isomorphisme de $\Modu A$.
\item[d.] Le morphisme
\[
\unite_A \ts \coaugmentation_C : e \ra A \tw C
\]
est une Èquivalence faible de $\coModcu C.$
\item[e.] Le morphisme d'algËbres $f_\tau$ (\ref{cochaines_universelles}) est un quasi-isomorphisme.
\item[f.] Le morphisme de cogËbres $g_\tau$ (\ref{cochaines_universelles}) est une Èquivalence de $\cocoga.$
\end{itemize} \francais
\end{proposition}

\dem 

$a \Rightarrow b.$  C'est une consÈquence du point {\it b} du thÈorËme \ref{theoreme_cmf_coModcu}.

$a \Rightarrow c.$ C'est clair.

$b \Rightarrow d.$ Nous avons l'ÈgalitÈ $\Psi_e = \unite_A \ts \coaugmentation_C.$

$c \Rightarrow a.$ La sous-catÈgorie de $\cd A$ formÈe des objets $M$ tel
que
\[
\Phi : LRM \ra M
\]
est un quasi-isomorphisme est une sous-catÈgorie triangulÈe aux sommes infinies contenant
$A$ par hypothËse. Elle coÔncide donc (voir \cite[4.2]{Keller94}) avec $\cd A.$

$d \Rightarrow e.$
Rappelons que $\tau_C : C \ra \Oma C$ est acyclique (\ref{lemme_acyclicite_cochaines_universelles}).
Cela implique que le  morphisme
\[
L_{\tau_C}e = \Oma C \arr{} L_{\tau_C}(A \tw C) = L_{\tau_C}R_{\tau_C}\Res A
\]
et le morphisme d'adjonction
\[
L_{\tau_C}R_{\tau_C}\Res A \ra \Res A
\]
sont des quasi-isomorphismes.
Le morphisme $f_\tau $ est un quasi-isomorphisme car il est Ègal ‡ la composÈe
\[
\Oma C \arr{} L_{\tau_C}R_{\tau_C}\Res A \arr{} \Res A.
\]

$e \Leftrightarrow f.$ C'est le point {\it b} du thÈorËme \ref{theoreme_cmf_cocog}.

$e \Rightarrow a.$ Comme la cochaÓne $\tau_C$ est acyclique, le morphisme d'adjonction
\[
L_{\tau_C}R_{\tau_C} M = M \tw C \ts_{\tau_C} \Oma C \ra M
\]
est un quasi-isomorphisme. Par ailleurs, il est Ègal ‡ la composÈe
\[
M \tw C \ts_{\tau_C} \Oma C \arr{\phi_M} M \tw C \tw A \arr{\Phi_M} M.
\]
Il nous suffit donc de montrer que le morphisme $\phi_M$ induit par
le morphisme $f_\tau$ est un quasi-isomorphisme.
Munissons le comodule $M \tw C$ de sa filtration primitive.
Nous avons alors
\[
\ogr (M \tw C) = M \ts \ogr C 
\]
et des filtrations induites sur $M \tw C  \tw A$ et $M \tw C \tw \Oma C$ qui vÈrifient
\[
\ogr (M \tw C  \tw A) = M \ts \ogr C \ts A \quad \mbox{et} \quad
\ogr (M \tw C \ts_{\tau_C} \Oma C) = M \ts \ogr C \ts \Oma C.
\]
Pour ces filtrations, le morphisme $\phi_M$ est un morphisme filtrÈ et
il induit des quasi-isomorphismes dans les objets graduÈs
car $f_\tau$ est un quasi-isomorphisme. Il est donc un
quasi-isomorphisme.\findem

%
%
%
%
%
%
%
%
%
%
%
%
%
%
%
%
%
%

\section{Polydules} \label{section_definition_ai-modules}

\subsection{DÈfinitions}
%
%
%
%
%
%
%
%
%
%
\begin{definition} {\em 
Soit $A$ une $\rm{A}_n$-algËbre. Un {\em $\rm{A}_n$-module sur $A$ dans la catÈgorie $\gr \sf C'$}
\index{A-n module@{$\rm{A}_n$-module}} est un objet
graduÈ $M$ dans $\gr \sf C'$ muni d'une famille de morphismes graduÈs
\[
m^M_i :  M \ts A\tp{i-1} \ra M, \quad  1 \leq i\leq n,
\]
de degrÈ $2-i$, telle qu'une Èquation $(*'_m)$ de la mÍme forme que
l'Èquation $(*_m)$ de la dÈfinition
\ref{definition_ai-algebre} est vÈrifiÈe pour tout $1\leq m \leq n$.
Dans l'Èquation 
$(*'_m)$, pour $j > 0$, les termes
\[
m_i(\Id^{\ts j}\ts m_k\ts \Id^{\ts l})
\]
de l'Èquation $(*_m)$ doivent Ítre interprÈtÈs comme
\[
m^M_i(\Id^{\ts j}\ts m_k\ts \Id^{\ts l}): M \ts A\tp{m-1} \ra M,
\]
et, pour $j=0$, comme
\[
m_i^M( m_k^M \ts \Id^{\ts l}) : M \ts A\tp{m-1} \ra M.
\]}
\end{definition}

\begin{definition}{\em \label{definition_ai-module}
Soit $A$ une $\ai$-algËbre. Un {\em $A$-polydule dans $\gr \sf C'$} (dans la littÈrature, cette
structure est communÈment appelÈe un {\em $\ai$-module sur $A$})
\index{A-infini module@{$\ai$-module}} \index{polydule}
est un objet
graduÈ $M$ muni d'une famille de morphismes graduÈs
\[
m^M_i :  M \ts A\tp{i-1} \ra M, \quad  1 \leq i,
\]
de degrÈ $2-i$, telle que l'Èquation $(*'_m)$ est vÈrifiÈe pour tout $1 \leq m.$
}
\end{definition}

\begin{definition}\label{definition_suspension_ai-module}{\em 
La {\em suspension} \index{suspension}
$SM$ d'un $A$-polydule est le $A$-polydule dont l'objet graduÈ sous-jacent est la suspension
$SM$ et dont les multiplications
sont dÈfinies par 
\[
m_i^{SM} = (-1)^{i} s\circ m_i^M \circ (\si \ts \Id\tp{i-1}) , \quad i\geq 1.
\]
La section \ref{construction_bar_cobar_mod} nous certifiera que ceci dÈfinit bien un $A$-polydule.
}\end{definition}

\begin{definition}{\em
Soit $A$ une $\rm{A}_n$-algËbre, et $M$ et $N$ deux $\rm{A}_n$-modules sur $A$.
Un {\em $\rm{A}_n$-morphisme} de $\rm{A}_n$-modules
\index{A-n morphisme@{$\rm{A}_n$-morphisme}}
$f:M \ra N$  est une famille de
morphismes graduÈs de $\sf C'$
\[
f_i :  M \ts A\tp{i-1} \ra N, \quad 1 \leq i \leq n,
\]
de degrÈ $1-i$, vÈrifiant, pour tout $1\leq m \leq n$, l'ÈgalitÈ
\[
(**'_m) \quad \sum (-1)^{jk+l}f_i(\Id^{\ts j}\ts m_k\ts \Id^{\ts l}) = \sum
 m_{s+1} (f_{r}\ts \Id\tp s)
\]
dans $ \Hom_{\sz \gr \sf C'}(M \ts A\tp{m-1},N)$, o˘  $j+k+l=m$, $i=j+k+1$ et $r+s=m$. 
Un $\rm{A}_n$-morphisme $f$ est {\em strict} si $f_i = 0$ pour tout $i \geq 2.$
Soit $M$, $N$ et $T$ trois $\rm{A}_n$-modules sur $A$.
Soit $g : M\ra N $ et $f : N \ra T$ deux $\rm{A}_n$-morphismes de $\rm{A}_n$-modules.
La {\em composition} $f \circ g : M \ra T$ est dÈfinie par la suite
\[
(f\circ g)_i = \sum_{k+l = i} f_{1+l} (g_k \ts \Id \tp l), \quad 1 \leq i \leq n.
\]
}
\end{definition}

\begin{definition}\label{definition_ai-morphisme_mod} {\em
Soit $A$ une $\ai$-algËbre et $M$ et $N$ deux $A$-polydules
Un {\em $\ai$-morphisme} \index{A-infini morphisme@{$\ai$-morphisme}}
$f : M \ra N$ est une famille de morphismes graduÈs
\[
f_i :  M \ts A\tp{i-1} \ra N, \quad  1 \leq i,
\]
de degrÈ $1-i,$ telle que l'Èquation $(**'_m)$ est vÈrifiÈe pour tout $1 \leq m.$
La {\em composition} des $\ai$-morphismes est dÈfinie par les mÍmes formules que celle
de la composition d'$\rm{A}_n$-morphismes. Un $\ai$-morphisme $f$ est {\em strict} si $f_i = 0$ pour tout $i \geq 2.$
\index{strict (A-infini morphisme)@{strict ($\ai$-morphisme)}}
\index{A-infini morphisme@{$\ai$-morphisme}!strict}
}
\end{definition}

Il rÈsultera de la section \ref{construction_bar_cobar_mod} que nous obtenons bien ainsi une
catÈgorie. Nous la notons $\aiMod A$. \indexnotation{aiMod}
La lettre $\mathsf N$ remplace la lettre $\mathsf M$
de $\Modu$ et se rapporte au $\mathsf N$on dans ``$\ai$-module $\mathsf N$on unitaires''.
Notons $\aiModst A$ \indexnotation{aiModst}
la sous-catÈgorie de $\aiMod A$ dont les objets sont les $A$-polydules
et dont les morphismes sont les $\ai$-morphismes stricts.

\begin{remarque}{\em
Soit $A$ une $\ai$-algËbre.
De maniËre analogue ‡ la remarque \ref{remarque1_definition_aia}, si $M$ est
un $A$-polydule, 
\english \begin{itemize}
\item[-] $(M,m_1)$ est un complexe;
\item[-] le morphisme $m^M_2 : M \ts A \ra M$ dÈfinit une action ‡ homotopie
prËs de l'algËbre
fortement homotopiquement  associative (\ref{remarque1_definition_aia}) $A$
sur $M$. Le dÈfaut de compatibilitÈ de la multiplication $m^A_2$ et de l'action
$m^M_2$ est Ègal au bord de $m^M_3$ dans
\[
(\Hom_{\sz \gr \sf C'}(M \ts A\tp 2,M), \delta),
\]
o˘ $\delta$ est dÈfini ‡ l'aide de $m^M_1$ et $m^A_1.$
\item[-] Si $f : M \ra N$ est un $\ai$-morphisme de $A$-polydules, le morphisme $f_1$ est
un morphisme de complexes $(M,m^M_1) \ra (N,m^N_1).$
\end{itemize} \francais
}\end{remarque}

\begin{remarque}{\em
Soit $A$ une $\ai$-algËbre.
Les morphismes $m^A_i$, $i \geq 1$,
dÈfinissent une structure de $A$-polydule sur l'objet sous-jacent ‡ $A$.
}\end{remarque}

\begin{remarque}{\em 
Soit $A$ un objet de $\alg$ et $(M,d^M,\Delta^M)$ un $A$-module diffÈrentiel graduÈ.
Les morphismes
\[
m^M_1 = d^M, \quad m^M_2 = \Delta^M, \quad m^M_i =  0\quad  \rm{ pour }\quad i \geq 3
\]
dÈfinissent sur l'objet sous-jacent ‡ $M$ une structure de $A$-polydule.
La catÈgorie des $A$-modules diffÈrentiels graduÈs est une sous-catÈgorie non pleine
de la catÈgorie des $A$-polydules.
}\end{remarque}

\begin{definition}{\em
Soit $A$ une $\ai$-algËbre et $M$ et $N$ deux $A$-polydules.
Un $\ai$-morphisme de $A$-polydules $f : M \ra N$ est un {\em $\ai$-quasi-isomorphisme}
\index{A-infini quasi-isomorphisme@{$\ai$-quasi-isomorphisme}}
si $f_1$ est un quasi-isomorphisme de complexes.
}\end{definition}

\begin{definition}\label{definition_homotopie_ai-morphismes_mod}
{\em
Soit $A$ une $\ai$-algËbre et $M$ et $N$ deux $A$-polydules.
Soit $f$ et $g$ deux $\ai$-morphismes $M \ra N.$
Une {\em homotopie entre $f$ et $g$}
\index{homotopie!A-infini-morphismes@{$\ai$-morphismes}}
est une famille de morphismes
\[
h_i:M\ts A^{\ts i-1}\ra N, \hspace{1cm} 1 \leq i ,
\]
de degrÈ  $-i$ vÈrifiant, pour tout $1 \leq m$,
l'Èquation 
\[
\begin{array}{c}(***'_m) \\ \hspace*{1cm}\end{array} \quad
\begin{array}{rcl}
 f_m - g_m & = & \sum (-1)^{s}
m_{1+s} (h_{r} \ts \Id \tp s)\\ 
& & + \sum (-1)^{jk+l}h_i(\Id^{\ts j}\ts m_k\ts \Id^{\ts l})
\end{array}
\]
dans $\Hom_{\gr \sf C'} (M \ts A \tp{m-1},N)$, o˘ $r+s = m$ et $j+k+l = m$.
Deux $\ai$-morphismes d'$\ai$-algËbres $f$ et $g$ sont {\em homotopes} s'il existe une homotopie
entre $f$ et $g$. 
}
\end{definition}


\subsection{UnitÈs strictes, augmentations et rÈductions} \label{section_augmentation_reduction}

Dans ce chapitre, nous Ètudierons les polydules {\em strictement unitaires} sur des $\ai$-algËbres
{\em augmentÈes}. Nous dÈfinirons donc ici un type d'unitaritÈ pour les $\ai$-structures : l'{\em
unitaritÈ stricte}. Cette structure nous permettra de gÈnÈraliser certaines
propriÈtÈs des modules unitaires aux polydules.
La pertinence de cette notion d'unitaritÈ
relativement ‡ l'homotopie des $\ai$-structures fera l'objet du chapitre
\ref{chapitre_Unite}.

\begin{definition}{\em \label{definition_unite_stricte}
Une $\ai$-algËbre $A$ est {\em strictement unitaire} \index{strictement unitaire!A-infini algebre@{$\ai$-algËbre}}
si elle est munie d'un morphisme graduÈ
$\unite : e \ra A$ de degrÈ $0$ tel que $m_i(\Id \hdots \Id \ts \unite \ts \Id \hdots \Id) = 0$ 
pour tout $i \neq 2$ et 
\[
m_2 (\Id_A \ts \unite) = m_2 (\unite \ts \Id_A ) = \Id_A.
\]
Le morphisme $\unite$ s'appelle l'{\em unitÈ} (stricte) \index{unite@{unitÈ}!stricte} de $A.$
Si $A$ et $A'$ sont deux $\ai$-algËbres strictement unitaires, un  $\ai$-morphisme
$f : A \ra A'$ est {\em strictement unitaire} si $f_1 \unite^A  = \unite^{A'}$
et $f_i (\Id \hdots \Id \ts \unite \ts \Id \hdots \Id) = 0$ 
pour tout $i \geq 2.$}
\end{definition}

Par la remarque \ref{remarque_ai-structure_alg},
une algËbre diffÈrentielle graduÈe unitaire est une $\ai$-algËbre strictement
unitaire. En particulier, l'algËbre $e$ est une $\ai$-algËbre strictement
unitaire.

\begin{definition}{\em
Une $\ai$-algËbre $A$ est {\em augmentÈe} si elle est strictement unitaire et 
munie d'un $\ai$-morphisme
strict d'$\ai$-algËbres strictement unitaires $\epsilon : A \ra e.$ Le morphisme
$\epsilon$ s'appelle l'{\em augmentation} \index{augmentation} de $A.$}
\end{definition}

L'{\em $\ai$-algËbre rÈduite} \index{reduction@{rÈduction}} $\b A$ est le noyau de $\epsilon.$
Soit $A$ une $\ai$-algËbre. L'{\em $\ai$-algËbre augmentÈe} $A^+$
a pour objet sous-jacent $A \oplus e$, ses multiplications $m_i$, $i \geq 1$, sont telles que
l'injection canonique $e \ra A \oplus e$ est l'unitÈ stricte et telles qu'elles coÔncident
avec $m^A_i,$ $i\geq 1$, sur $A$. Son augmentation
est la projection canonique $A \oplus e \ra e$. Nous notons $\aiaa$ 
\indexnotation{aiaa}la {\em catÈgorie}
des $\ai$-algËbres augmentÈes. Le {\em foncteur augmentation} \index{augmentation} $\aia \ra \aiaa$ est
une Èquivalence dont le quasi-inverse est le {\em foncteur rÈduction.}

\begin{definition}{\em \label{definition_polydules_strictement_unitaires}
Soit $A$ une $\ai$-algËbre strictement unitaire.
Un $A$-polydule $M$ est {\em strictement unitaire}
\index{strictement unitaire!polydule} si
$m^M_i(\Id_M \ts \Id \hdots \Id \ts \unite \ts \Id \hdots \Id) = 0$ 
pour tout $i \geq 3$ et 
\[
m^M_2 (\Id_M \ts \unite) = \Id_M.
\]
Un morphisme {\em strictement unitaire} de $A$-polydules strictement unitaires est un $\ai$-morphisme
$f$ de $A$-polydules
tel que
\[
f_i(\Id_M \ts \Id \hdots \Id \ts \unite \ts \Id \hdots \Id) = 0 , \quad i \geq 2.
\]
Si $f$ et $g$ sont deux morphismes strictement unitaires, une homotopie $h$
entre $f$ et $g$ est {\em strictement unitaire} si
\[
h_i(\Id_M \ts \Id \hdots \Id \ts \unite \ts \Id \hdots \Id) = 0 , \quad i \geq 2.
\]
Si $h$ est une homotopie strictement unitaire entre deux morphismes strictement
unitaires $f$ et $g$, on dit que $f$ et $g$ sont {\em homotopes} (relativement ‡ $h$)
\index{homotopie!A-infini-morphismes@{$\ai$-morphismes}} et on note $f \sim g$.
Nous notons $\aiModu A$ \indexnotation{aiModu} la {\em catÈgorie} des $A$-polydules strictement unitaires dont les
morphismes sont les morphismes strictement unitaires et 
$\aiModust A$ \indexnotation{aiModust} la {\em catÈgorie} des $A$-polydules strictement unitaires
dont les morphismes sont les $\ai$-morphismes stricts et strictement unitaires.
}\end{definition}

Si $A$ est une $\ai$-algËbre et $M$ un $A$-polydule, $M^+$ est l'{\em $A^+$-polydule}
(strictement unitaire) qui a pour objet sous-jacent $M$ et dont
la multiplication $m^{M^+}_i$, $i \geq 1$,
est telle que, restreinte ‡ $A$, elle coÔncide avec $m^M_i$, $i \geq 1$ (en particulier le
$m_1$ ne change pas).
Ceci dÈfinit un isomorphisme
\[
{}\+ : \aiMod A \arr{\sim} \aiModu A\+
\]
compatible ‡ l'homotopie.
Le quasi-inverse est donnÈ par le foncteur qui envoie $M$ sur
le $\b A$-polydule $\b M$ dont l'objet sous-jacent est $M$ et dont la multiplication
$m^{\b M}_i$, $i\geq 2$, est la restriction de $m^M_i$, $i\geq 2$, ‡
$M \ts A\tp{i-1}$.\\

\subsection{Construction bar} \label{construction_bar_cobar_mod}

Les dÈmonstrations de cette section Ètant presque identiques ‡
celles de la section \ref{construction_bar_cobar}, nous nous contentons
d'Ènoncer les rÈsultats.\\

{\noindent \bf Construction bar des polydules}\\

Soit $A$ et $M$ deux objets graduÈs. Pour chaque $i \geq 1$, nous dÈfinissons
une bijection
\[
\begin{array}{rcl}
\Hom_{\sz \gr \sf C'}(M \ts A \tp{i-1} ,M) & \ra &  \Hom_{\sz \gr \sf C'}(SM \ts (SA) \tp{i-1},SM)\\
m^M_i & \mapsto & b^M_i
\end{array}
\]
par la relation 
\[
\si \circ b^M_i = -  m^M_i \circ \si \tp i \hspace{1cm} (\mbox{o˘}  \hspace{.3cm} \si = s^{-1}).
\]
Soit $A$ une $\ai$-algËbre.
Nous rappelons (\ref{lemme_comodules_colibres})
qu'une diffÈrentielle $b^M$ sur le $(BA)\+$-comodule (co-unitaire) graduÈ $SM \ts (BA)\+$ est dÈterminÈe
par la composition
\[
(\Id \ts \counite^{(BA)\+}) \circ b^M : SM \ts (BA)\+ \ra SM
\]
dont nous notons les composantes $b^M_i$, $i\geq 1$.
Les bijections $m^M_i \lra b^M_i$ induisent une bijection de l'ensemble des structures
de $A$-polydule sur $M$  sur l'ensemble
des diffÈrentielles $b^M$ sur le $(BA)\+$-comodule graduÈ $SM \ts (BA)\+$.

Soit $A$, $M$ et $N$ trois objets graduÈs. Pour chaque $i \geq 1$, nous dÈfinissons
une bijection
\[
\begin{array}{rcl}
\Hom_{\sz \gr \sf C'}(M \ts A \tp{i-1} ,M) & \ra &  \Hom_{\sz \gr \sf C'}(SM \ts (SA) \tp{i-1},SM)\\
f_i & \ra & F_i
\end{array}
\]
par les relations  
\[
\si \circ F_i = (-1)^{|F_i|} f_i  \circ \si \tp i, \quad i \geq 1,
\]
o˘ $F_i$ est un morphisme graduÈ de degrÈ $|F_i|.$
Soit $A$ une $\ai$-algËbre.
On rappelle (\ref{lemme_comodules_colibres})
qu'un morphisme graduÈ de $(BA)\+$-comodules (co-unitaires)
\[
F : SM \ts (BA)\+ \ra SN \ts (BA)\+
\]
est dÈterminÈ par la composition
\[
(\Id \ts \counite^{(BA)\+}) \circ F : SM \ts (BA)\+ \ra SM
\]
dont nous notons les composantes $F_i$, $i\geq 1$.
Les bijections $f_i \lra F_i$ induisent une bijection du produit des ensembles de morphismes
graduÈs 
\[
f_i : M \ts A \tp {i-1} \ra N, \quad i\geq 1,
\]
de degrÈ $1-i+n$, sur l'ensemble
des morphismes graduÈs de $(BA)\+$-comodules $F : SM \ts (BA)\+ \ra SN \ts (BA)\+$
de degrÈ $n.$ Si $M$ et $N$ sont des $A$-polydules,
cette bijection envoie bijectivement l'ensemble des familles dÈfinissant un $\ai$-morphisme $f : M \ra N$ sur
l'ensemble des morphismes diffÈrentiels graduÈs de $(BA)\+$-comodules
\[
F : SM \ts (BA)\+  \ra SN \ts (BA)\+ .
\]
Si $f$ et $g$ sont deux $\ai$-morphismes de $A$-polydules, la mÍme bijection envoie
bijectivement l'ensemble des homotopies entre $f$ et $g$ sur l'ensemble des homotopies
entre les morphismes de $(BA)\+$-comodules $F$ et $G$ correspondant ‡ $f$ et $g.$

Ceci nous donne un foncteur
\[
\aiMod A \ra \coModcu (BA)\+, \quad M \mapsto (SM \ts (BA)\+,b^M).
\]

{\noindent \bf Construction bar des polydules strictement unitaires sur une $\ai$-algËbre augmentÈe}\\

Soit $A$ une $\ai$-algËbre augmentÈe. Nous notons $\Ba A$ la cogËbre co-augmentÈe
$(B\b A)\+$, o˘ $\b A$ est l'$\ai$-algËbre rÈduite associÈe ‡ $A$. Attention ‡ ne pas
confondre les cogËbres co-augmentÈes $\Ba A$ et $(BA)\+.$ \index{construction bar}
\indexnotation{Ba}

Par la section \ref{section_augmentation_reduction},
le foncteur $N \mapsto \b N$ est un isomorphisme de catÈgories
\[
\aiModu A \arr{\sim} \aiMod \b A.
\]
Le foncteur composÈ
\[
B_A : \aiModu A \arr{\sim} \aiMod \b A \ra \coModcu \Ba A
\]
est appelÈ le foncteur {\em construction bar}\index{construction bar}\indexnotation{BM}.
Nous le noterons souvent $B$.
La suspension $SM$ d'un polydule est envoyÈe
par le construction bar sur $B SN = (S^2 N \ts \Ba A,b^{SN})$. Nous vÈrifions que ce dernier est isomorphe
‡ $SBN$. Le foncteur construction bar
envoie des $\ai$-morphismes homotopes sur des morphismes homotopes de comodules et il induit une Èquivalence
entre la catÈgorie $\aiModu A$ et la sous-catÈgorie $\colibre \Ba A$ de $\coModcu \Ba A$ formÈe des objets
presque colibres.\\

\subsection{AlgËbre Enveloppante} \label{section_algebre_enveloppante}

Dans cette section,
nous dÈfinissons l'algËbre enveloppante $\E A$
d'une $\ai$-algËbre augmentÈe $A$ puis montrons que la catÈgorie $\Modu \E A$
est isomorphe ‡ la catÈgorie $\aiModust A$. \\

Soit $V$ un espace graduÈ (resp.~diffÈrentiel graduÈ). 
L'{\em algËbre tensorielle (augmentÈe)}  \indexnotation{TV}
\index{tensorielle augmentee (algebre)@{tensorielle augmentÈe (algËbre)}}
\index{algebre@{algËbre}!tensorielle augmentÈe}
$TV$ est l'augmentation $(\b TV)\+$
de l'algËbre tensorielle rÈduite. Soit $i : V \ra TV$ l'injection canonique.

\begin{lemme} \label{lemme1_algebre_enveloppante}
Soit $M$ un objet graduÈ.
L'application $\mu^M \mapsto \mu^M (\Id  \ts i)$ est une bijection de l'ensemble
des structures de $TV$-module unitaire sur $M$ sur l'ensemble des morphismes graduÈs
\[
M \ts V \ra M
\] de degrÈ $0$. L'application inverse associe ‡ $g$ la multiplication
\[
\mu : M \ts TV \ra M
\]
dont la composante $M \ts e \ra M$ est l'identitÈ et la composante
$M \ts V \tp i \ra M$ est le morphisme $g \circ (g \ts \Id) \circ \cdots \circ (g \ts \Id \tp{i-1})$.
\end{lemme}

\begin{definition}{\em \label{definition_algebre_enveloppante}
Soit $A$ une $\ai$-algËbre augmentÈe.
L'{\em algËbre enveloppante} \index{algebre@{algËbre}!enveloppante}
de $A$ est l'algËbre diffÈrentielle
graduÈe $\E A = \Oma \Ba  A$, c'est-‡-dire l'algËbre $(\Omega B \b A)\+$.} \indexnotation{EA}
\end{definition}
\begin{lemme}\label{lemme_algebre_enveloppante}
L'$\ai$-morphisme $A \ra \E A$ donnÈ par
le morphisme d'ad\-jonction 
\[
\Ba A \ra \Ba \E A = \Ba \Oma \Ba  A
\]
est un $\ai$-quasi-isomorphisme. Il est universel parmi les $\ai$-morphismes de $A$
vers une algËbre diffÈrentielle graduÈe.
\findem
\end{lemme}

\dem 
C'est un $\ai$-quasi-isomorphisme par le lemme \ref{lemme3_cmf_aia}.
L'universalitÈ est immÈdiate gr‚ce ‡ l'adjonction $(\Omega, B)$.
\findem \\


\begin{lemme} \label{A-polydules_EA-modules}
Nous avons un isomorphisme de catÈgories
\[
i : \Modu \E A \ra \aiModust A, \quad M \ra S^{-1}M.
\]
\end{lemme}

\dem  Soit $M$ un objet graduÈ. Nous allons montrer que les structures de $\E A$-module unitaire
sur $SM$ sont les structures de $A$-polydule strictement unitaire sur $M$.
Soit $m^M_1$ une diffÈrentielle sur $M$ et soit
\[
m^M_i : M \ts A \tp{i-1} \ra M, \quad i \geq 2,
\]
des morphismes graduÈs de degrÈ $2-i$. Nous dÈfinissons ‡ l'aide 
des bijections $m^M_i \lra b^M_i$ de la section \ref{construction_bar_cobar},
un morphisme
\[
g : SM \ts (B\b A) \ra SM.
\]
Par le lemme \ref{lemme1_algebre_enveloppante}, le morphisme
\[
SM \ts S^{-1}(B\b A) \arr{\Id \ts s} SM \ts (B\b A) \arr{g} SM
\]
se relËve en une structure $\mu^{\E}$ de $\Oma \Ba A$-module graduÈ unitaire sur $SM.$
Nous vÈrifions que $(SM,\mu^{\E},Sm_1)$
dÈfinit un module diffÈrentiel graduÈ unitaire si et seulement si les $m^M_i$, $i\geq 1$,
dÈfinissent une structure de $A$-polydule strictement unitaire sur $M$.
Si $SM$ et $SN$ sont deux $\E A$-modules,
les morphismes de $\E A$-modules $SM \ra SN$ s'identifient clairement aux $\ai$-morphismes
stricts de $A$-polydules
$M \ra N$. \findem

%
%

%
%
%
%
%
%
%
%
%
%
%
%
%

\newpage

\section{CatÈgorie dÈrivÈe d'une $\ai$-algËbre augmentÈe} \label{section_categorie_derivee_augmentee}

{\bf Introduction}\\

\noindent Soit $A$ une $\ai$-algËbre augmentÈe.
Le but de cette section est de montrer que la catÈgorie dÈrivÈe
\[
\cd_\infty A = \aiModu A [Qis^{-1}]
\]
est Èquivalente aux catÈgories
\[
\ch_\infty A = \aiModu A/\!\sim\, \quad \mbox{et} \quad \big(\aiModust A\big) [Qis^{-1}]
\]
o˘ $\sim$ est la relation d'homotopie.
La catÈgorie dÈrivÈe d'une $\ai$-algËbre quelconque
est ÈtudiÈe au chapitre \ref{chapitre_Cat_der}.\\

{\noindent \bf Plan de la section}\\

\noindent Cette section est divisÈe en trois sous-sections.
Dans la sous-section \ref{section_objets_fibrants_coModcu}, nous dÈmontrons le thÈorËme de
l'homotopie et celui des $\ai$-quasi-isomorphismes pour les polydules. Pour cela, nous
caractÈriserons les objets fibrants de la catÈgorie de modËles
$\coModcu \Ba A$ :
{\em ils sont exactement les facteurs directs des objets presque colibres} et nous montrons que
les thÈorËmes ci-dessus  apparaissent alors comme des
cas particuliers de rÈsultats fondamentaux de l'algËbre homotopique de Quillen
(voir appendice \ref{section_rappels_cmf}). 
Dans la sous-section \ref{section_definition_categorie_derivee_augmentee_polydules},
nous montrons les Èquivalences annoncÈes dans l'introduction ci-dessus (toujours gr‚ce ‡
l'algËbre homotopique de Quillen).
Dans la section \ref{section_structure_triangulee}, nous
Ètudions la structure triangulÈe de $\cd_\infty A$.

%
%
%
%
%
%
%
%
%
%
%
%
%
%
%

\subsection{Objets fibrants de $\coModcu \Ba A$}
\label{section_objets_fibrants_coModcu}

Soit $A$ un objet de $\aiaa.$ Le but de cette section est de montrer la proposition suivante :

\begin{proposition} \label{corollaire2_categorie_derivee}
\english \begin{itemize}
\item[a.] La relation d'homotopie (\ref{definition_polydules_strictement_unitaires}) dans $\aiModu A$
est une relation d'Èqui\-valence compatible ‡ la composition.
\item[b.] Un $\ai$-quasi-isomorphisme de $A$-polydules est une Èquivalence d'homo\-topie.
\item[c.] Soit $A'$ un objet de $\alga$. Soit $\Modush A'$ la sous-catÈgorie pleine de $\aiModu A'$ formÈe
des $A'$-modules diffÈrentiels graduÈs unitaires. Notons $\sim$ la relation d'homotopie sur $\Modush A'$.
L'inclusion $\Modu A' \hookrightarrow \Modush A'$ induit une Èquivalence
\[
\cd A' \arr{\sim} \Modush A'/\! \sim .
\]\end{itemize} \francais
\end{proposition}

\begin{remarque}{\em
Le point {\it c} reste vrai mÍme dans le cas o˘ l'algËbre diffÈrentielle graduÈe unitaire $A'$ n'est
pas augmentÈe (voir \ref{lemme_categorie_derivee_algebre_differentielle_graduee_unitaire}).
}\end{remarque}

\dem La dÈmonstration est identique ‡ celle du corollaire \ref{corollaire_cmf_aia}. Elle procËde de
la mÍme maniËre en utilisant (‡ la place du thÈorËme principal \ref{theoreme_cmf_cocog})
le thÈorËme \ref{theoreme_cmf_coModcu} et la proposition \ref{proposition2_categorie_derivee} ci-dessous.
\findem \\

{\noindent \bf Un raffinement de la caractÈrisation des objets fibrants du thÈorËme
\ref{theoreme_cmf_coModcu}}\\

Soit $C$ un objet de $\cocoga.$
Munissons la catÈgorie $\coModcu C$ de sa structure canonique de catÈgorie de modËles
(\ref{definition_cmf_canonique_coModcu}).
Soit $\tau : C \ra A'$ une cochaÓne tordante admissible
acyclique, o˘ $A'$ est un objet de $\alga$ (il existe toujours une telle cochaÓne gr‚ce au lemme
\ref{lemme_acyclicite_cochaines_universelles}).
Le thÈorËme \ref{theoreme_cmf_coModcu} dit que les objets fibrants de
$\coModcu C$ sont les facteurs directs d'objets de la forme $R_\tau M$, o˘
$M$ est un objet de $\Modu A'$. En particulier,
les objets fibrants sont facteurs directs d'objets presque colibres de $\coModcu C.$
Montrons que la rÈciproque est vrai pour certaines cogËbres :

\begin{proposition}\label{proposition2_categorie_derivee}
Soit $C$ un objet de $\cocoga$ qui est isomorphe, en tant que cogËbre graduÈe, ‡ une cogËbre
tensorielle.
Les objets fibrants de $\coModcu C$ sont exactement les facteurs directs d'objets presque colibres.
\end{proposition}

En particulier, puisque la cogËbre $C$ est isomorphe ‡ la construction bar $\Ba A$ d'un objet $A$ de $\aiaa$
les objets fibrants de $\coModcu C$ sont exactement les facteurs directs de comodules qui sont l'image
par la construction bar d'un $A$-polydule. La dÈmonstration de ce rÈsultat est reportÈe
‡ la fin de cette section. Nous dÈmontrons au prÈalable quelques propositions.\\

{\noindent \bf $\aiModu A$ comme ``catÈgorie de modËles sans limites''}\\

Soit $A$ une $\ai$-algËbre augmentÈe.
Dans la catÈgorie $\aiModu A$, nous considÈrons les trois
classes de morphismes suivantes :

\english \begin{itemize}
\item[-] la classe $\weq$ est formÈe des {\em Èquivalences faibles}, c'est-‡-dire des
$\ai$-quasi-isomorphismes,
\item[-] la classe $\cof$ est formÈe des {\em cofibrations}, c'est-‡-dire des $\ai$-morphismes
 $j : M \ra M'$ tels que $j_1$ est un monomorphisme,
\item[-] la classe $\fib$ est formÈe des {\em fibrations}, c'est-‡-dire des $\ai$-morphismes
 $q : M \ra M'$ tels que $q_1$ est un Èpimorphisme.
\end{itemize} \francais

\begin{theoreme} \label{theoreme_cmf_aiModu}
La catÈgorie $\aiModu A$, munie des trois classes dÈfinies ci-dessus, vÈrifie
l'axiome  (A) du thÈorËme \ref{theoreme_cmf_aia}
 et les axiomes (CM2) -- (CM5) de la dÈfinition
\ref{definition_cmf}. Tous les objets sont fibrants et cofibrants. 
\end{theoreme}

\dem Elle est identique ‡ celle de \ref{theoreme_cmf_aia} car basÈe sur les lemmes d'obstruction
(voir appendice \ref{section_obstruction_polydules}). \findem \\

{\noindent \bf Liens entre la ``catÈgorie de modËles sans limites'' $\aiModu A$ et la catÈgorie de modËles
$\coModcu  \Ba A$ }

\begin{proposition} \label{proposition1_categorie_derivee}
Soit $M$ et $M'$ deux objets de $\aiModu A$.
\english \begin{itemize}
\item[a.]
Un $\ai$-morphisme  $f : M \ra M'$ est un $\ai$-quasi-isomorphisme
de $\aiModu A$ si et seulement si le morphisme $B f : B M \ra B M'$ est une Èquivalence faible de
$\coModcu \Ba A$.
\item[b.] Un $\ai$-morphisme  $j : M \ra M'$ est une cofibration
de $\aiModu A$ si et seulement si $B j : B M \ra B M'$ est une cofibration de
$\coModcu \Ba A$.
\item[c.]  Un $\ai$-morphisme  $q : M \ra M'$ est une fibration
de $\aiModu A$ si et seulement si $B q : B M \ra B M'$ est une fibration de
$\coModcu \Ba A$.
\end{itemize} \francais
\end{proposition}

\dem 
Soit $\E A$ l'algËbre enveloppante de $A$.
Rappelons (\ref{cochaines_universelles}) que la cochaÓne tordante
universelle
\[
\tau : \Ba A \ra \Oma \Ba A  = \E A
\]
est acyclique.
Par le corollaire \ref{corollaire_theoreme_cmf_coModcu}, nous avons une Èquivalence de Quillen
\[
(L,R) : \coModcu \Ba A \ra \Modu \E A.
\]
{\it a.} Si $f$ est un $\ai$-quasi-isomorphisme, le morphisme
$B f$ est un quasi-isomorphisme filtrÈ pour les filtrations primitives.
Par le lemme \ref{lemme1_cmf_coModcu},
il est une Èquivalence faible de $\coModcu \Ba A$. Supposons que $B f$ est
une Èquivalence faible de $\coModcu \Ba A$. Soit le diagramme de $\coModcu \Ba A$
\[
\xymatrix{B M \ar[d]_{B f} \ar[r] & RL B M \ar[d]^{RL B f} \\
B M' \ar[r] & RL B M'.}
\]
Comme $R = B i$, ce diagramme est l'image par $B$ d'un diagramme
\[
\xymatrix{M \ar[d]_{f} \ar[r] & iL B M \ar[d]^{iL B f} \\
M' \ar[r] & iL B M'.}
\]
Comme $Bf$ est une Èquivalence faible de $\coModcu \Ba A$, le morphisme $LB f$ est
un quasi-isomorphisme de $\Modu \E A$. Le morphisme (strict) $iLB f$ est
donc un $\ai$-quasi-isomorphisme dans $\aiModu A.$
Le lemme \ref{lemme1_categorie_derivee} ci-dessous montre que les flËches horizontales
du diagramme ci-dessus reprÈsentent des $\ai$-quasi-isomorphismes. Par la propriÈtÈ de
saturation des $\ai$-quasi-isomorphismes dans $\aiModu A$, $f$ est donc un $\ai$-quasi-isomorphisme.

{\it b} et {\it c.} MÍme dÈmonstration que pour la proposition \ref{proposition1_cmf_aia}.
\findem

\begin{lemme}\label{lemme1_categorie_derivee}
Soit $M$ un objet de $\aiModu A$.
Le morphisme d'adjonction $B M \ra RLB M$ induit un quasi-isomorphisme
dans les primitifs.
\end{lemme}

\dem Il s'agit de montrer que le morphisme
\[
SM \ra SM \ts \Ba A \ts \E A
\]
est un quasi-isomorphisme. Notons $C$ la cogËbre $\Ba A$. Nous rappelons que
par dÈfinition $\Oma C = \Omega \b C.$
Il faut montrer que
\[
SM \ra SM \ts C \ts \Oma C
\]
est un quasi-isomorphisme. Munissons $\Oma C$  de la filtration
induite par la filtration primitive de $\b C$ considÈrÈ comme cogËbre. Nous avons alors une
filtration de $\Oma C$ dÈfinie par la suite
\[
\left(\Oma C\right)_i = (\Omega \b C)_i \oplus e, \quad i\geq 0. 
\]
Munissons $C$, considÈrÈ comme objet de $\coModu C$, de sa filtration primitive
de $C$-module (on la complËte par $C\prim 0 = e$). 
Munissons $M$ de la filtration dÈfinie par la suite $M_i = M$, $i \geq 0.$
Ces filtrations induisent sur $SM \ts C \ts \Oma C$ une filtration
de complexes.
Tout comme ‡ la fin de la dÈmonstration du point {\it b} du lemme
\ref{lemme_acyclicite_cochaines_universelles}, nous montrons que
\[
\ogr_0 (SM \ts C \ts \Oma C) = SM , \quad \ogr_i (SM \ts C \ts \Oma C) = 0 \quad
\mbox{pour} \quad i \geq 1.
\]
\findem \\

{\em DÈmonstration de la proposition \ref{proposition2_categorie_derivee}} : Nous pouvons supposer que $C$ est Ègal ‡ $\Ba A$, pour $A$ une $\ai$-algËbre augmentÈe.
Soit $\tau$ la cochaÓne tordante universelle de $\Ba A$. Nous savons que les objets fibrants de
$\coModcu \Ba A$ sont les facteurs directs d'objets de la forme $R M = M \tw \Ba A$, o˘
$M$ est un objet de $\Modu \Oma \Ba A$. Ils sont donc des facteurs directs des objets presque colibres.
RÈciproquement, si $N$ est un objet presque colibre,
il est isomorphe ‡ l'image par la construction bar d'un objet $M$ de $\aiModu A$. Ce dernier Ètant
fibrant dans $\aiModu A$, l'objet $N$ est fibrant dans $\coModcu \Ba A$ par le point {\it c} de la
proposition \ref{proposition1_categorie_derivee}.\findem

\subsection{CatÈgorie dÈrivÈe $\cd_\infty A$}
\label{section_definition_categorie_derivee_augmentee_polydules}

Dans cette section, nous dÈfinissons la catÈgorie dÈrivÈe $\cd_\infty A$ et en donnons plusieurs
descriptions.\\

Le point {\it a.} de la proposition \ref{corollaire2_categorie_derivee} montre que la dÈfinition suivante
a un sens.

\begin{definition}\label{definition_categorie_derivee_augmentee}{\em
Soit $A$ une $\ai$-algËbre augmentÈe. Nous notons $\ch_\infty A$
la catÈgorie $\aiModu A/\! \sim \,,$ o˘ $\sim$ est la relation d'homotopie
(voir \ref{definition_polydules_strictement_unitaires}).
La {\em catÈgorie dÈrivÈe $\cd_\infty A$ de $\aiModu A$}
\index{categorie@{catÈgorie}!derivee@{dÈrivÈe}} est la localisation par rapport aux $\ai$-quasi-isomorphismes
de la catÈgorie $\aiModu A.$
}\end{definition}

La proposition (\ref{corollaire2_categorie_derivee}) entraÓne le rÈsultat suivant :

\begin{corollaire} \label{corollaire3_categorie_derivee}
La projection canonique
\[
\ch_\infty A \ra \cd_\infty A
\]
est un isomorphisme.
\end{corollaire}

\dem Les $\ai$-quasi-isomorphismes Ètant des Èquivalences d'homotopie,
la projection canonique
\[
\ch_\infty A \ra \big(\ch_\infty A\big)[\weq^{-1}] \iso \cd_\infty A
\]
est une Èquivalence. \findem

\begin{lemme}\label{lemme2_categorie_derivee}
La composition des foncteurs (voir \ref{A-polydules_EA-modules})
\[
J : \Modu \E A \arr{i} \aiModust A \hookrightarrow \aiModu A
\]
induit un isomorphisme $\cd \E A \ra \cd_\infty A.$
\end{lemme}

\dem Nous avons un diagramme commutatif
\[
\xymatrix{
\Modu \E A \ar[r]^{J} \ar[d]_{R} & \aiModust A \ar[d]^{\mathsf{incl.}} \\
\coModcu \Ba A & \aiModu A \ar[l]^{B}
}
\]
et les foncteurs $J$, $R$ et $B$ induisent des Èquivalences entre les catÈgories
\[
\cd \E A, \quad \cd C, \quad (\aiModust A)[\weq^{-1}]\quad \mbox{et} \quad \cd_\infty A.
\]
 \findem

%
%
%
%
%
%
%
%
%
%
%
%
%
%
%
%
%

\subsection{Structure triangulÈe sur $\cd_\infty A$}
\label{section_structure_triangulee}

{\noindent \bf Suites exactes de $\aiModu A$}\\  \label{suites_exactes_aiModu}

Le foncteur
\[
i : \Modu \E A \ra \aiModu A, \quad SM \mapsto M,
\]
identifie (voir \ref{A-polydules_EA-modules}) la catÈgorie $\Modu \E A$
‡ la sous-catÈgorie $\aiModust A$ de $\aiModu A$. Il envoie
la suspension d'un $\E A$-module sur la suspension d'un $A$-polydule
(voir \ref{definition_suspension_ai-module}).
Il identifie les suites exactes courtes de $\Modu \E A$ qui
sont scindÈes dans la catÈgorie des modules graduÈs aux suites
de $\aiModu A$ formÈes d'$\ai$-morphismes stricts
\[
(*) \hspace*{2cm} M' \arr{j} M \arr{q}  M'',
\]
telles que
\[
0 \ra M' \arr{j_1} M \arr{q_1}  M'' \ra 0
\]
est une suite exacte de $\cc \sf C'$ et telles qu'il existe une rÈtraction $\rho$ de $j_1$ 
dans $\gr \sf C'$ telle que, pour tout $i \geq 2$,
\[
\rho m^{M}_i = m^{M'}_i (\rho \ts \Id \tp{i-1}).
\]

{\noindent \bf Structure triangulÈe sur $\cd_\infty A$}\\

Nous munissons la catÈgorie dÈrivÈe $\cd_\infty A$ de
l'unique structure triangulÈe (unique ‡ Èquivalence triangulÈe prËs) pour laquelle l'Èquivalence
\[
J : \cd \E A \ra \cd_\infty A
\]
du lemme \ref{lemme2_categorie_derivee} est triangulÈe.
Comme les foncteurs
\[
R : \cd \E A \ra \cd \Ba A \quad \mbox{et} \quad B : \cd_\infty A \ra \cd \Ba A
\]
sont des foncteurs triangulÈs nous en dÈduisons le thÈorËme suivant.

\begin{theoreme}
La structure triangulÈe sur $\cd_\infty A$ a pour endo\-foncteur suspension celui dÈfini en
 \ref{definition_suspension_ai-module}. Les triangles distinguÈs
sont exactement ceux qui sont isomorphes aux triangles provenant de suites exactes de la
forme $(*)$ de $\aiModu A.$
\findem
\end{theoreme}

{\noindent \bf CÙne d'un $\ai$-morphisme.}\\

Si $f : M \ra M'$ est un $\ai$-morphisme de $A$-polydules, son {\em cÙne} 
\index{cone d'un $\ai$-morphisme@{cÙne d'un $\ai$-morphisme}}
$C(f)$ est
le $A$-polydule $M' \oplus SM$ dont les multiplications
\[
m^{C(f)}_i : (M' \oplus SM) \ts A \tp{i-1} \ra  M' \oplus SM, \quad i \geq 1,
\]
sont donnÈes par les morphismes 
\[
m^{M'}_i, \quad m^{SM}_i\  \mbox{(voir 
\ref{definition_suspension_ai-module})} \quad \mbox{et} \quad
f_i \circ (\si \ts \Id \tp{i-1}).
\]
La construction bar envoie $C(f)$ sur le cÙne de $Bf.$

\begin{lemme}
Soit $A$ un objet de $\alga$. L'inclusion 
\[
\Modu A \hookrightarrow \aiModu A
\]
induit une Èquivalence triangulÈe
\[
\cd A \ra \cd_\infty A.
\]
\end{lemme}

\dem Comme $A \ra \E A$ (voir \ref{lemme3_cmf_aia}) est un quasi-isomor\-phisme,
nous avons une Èquivalence triangulÈe entre la catÈgorie
$\cd A$ et la catÈgorie $\cd \E A$. L'inclusion (\ref{A-polydules_EA-modules})
\[
i : \Modu \E A \ra \aiModu A
\]
induisant une Èquivalence triangulÈe
de $\cd \E A$ sur $\cd_\infty A,$ nous en dÈduisons le rÈsultat.
 \findem

\section{CatÈgorie dÈrivÈe des bipolydules (le cas augmentÈ)}

\label{section_categorie_derivee_augmentee_bipolydules}

{\bf \noindent Introduction }\\

Soit $A$ et $A''$ deux $\ai$-algËbres augmentÈes.
Dans cette section, nous dÈfinissons la catÈgorie dÈrivÈe $\cd_\infty (A,A'')$ des $A$-$A''$-bipolydules
strictement unitaires et nous en donnons plusieurs descriptions.
Le cas o˘ $A$ et $A''$ sont quelconques sera traitÈ au chapitre 
\ref{chapitre_Cat_der}.\\

{\bf \noindent Notations}\\ 

Soit $(\sf C,\ts,e)$ et $(\sf C'',\ts,e)$ deux $\corps$-catÈgories de Grothendieck semi-simples monoÔdales
et $\sf C'$ une $\corps$-catÈgorie de Grothendieck semi-simple (non nÈcessairement 
monoÔdale). Nous supposons que $\sf C$ est tressÈe (voir \cite[Chap.~XI]{MacLane98}).
Nous notons $\ts\op$ le produit tensoriel de $\sf C$ dÈfini par
\[
A \ts\op B = B \ts A.
\]
Supposons que la catÈgorie monoÔdale $\sf C$ agit ‡ gauche sur $\sf C'$ et la catÈgorie 
 monoÔdale $\sf C''$ agit ‡ droite sur $\sf C'$ de maniËre compa\-tible, i.~e.~$\sf C'$ est munie de deux
foncteurs ($\corps$-bilinÈaires sur les espaces de morphismes)
\[
\begin{array}{rcl}
\sf C' \times \sf C''  & \ra & \sf C',\\
(M',A'') & \mapsto & M' \ts A''
\end{array} \quad \mbox{et} \quad
\begin{array}{rcl}
\sf C \times \sf C'  & \ra & \sf C',\\
(A,M') & \mapsto & A \ts M'
\end{array}
\]
associatifs et unitaires ‡ des isomorphismes donnÈs prËs
(voir \cite[Chap.~XI]{MacLane98}) et tels que
\[
(A \ts M') \ts A'' = A \ts (M' \ts A'').
\]
Nous supposons en outre qu'on a une $\corps$-catÈgorie de Grothendieck
semi-simple monoÔdale
$\sf C \ts \sf C''$, munie d'un foncteur monoÔdal
\[
(\sf C,\ts\op) \times (\sf C'',\ts) \ra \sf C \ts \sf C'', \quad (A,A'') \mapsto A \ts A'',
\]
\[
\Hom_{\sf C}(A,B) \times \Hom_{\sf C''}(A'',B'') \ra \Hom_{\sf C \ts \sf C''}((A \ts A''),(B \ts B'')),
\]
bilinÈaire sur les espaces de morphismes, d'une action sur $C'$ et d'un isomorphisme
\[
M \ts (A \ts A'')  = A \ts M \ts A''.
\]

L'exemple suivant apparaÓt naturellement dans l'Ètude des $\ai$-catÈgories
(\ref{section_categories_de_base_aicat}).

\begin{exemple}{\em 
Soit $\mathbb A$ et $\mathbb B$ deux ensembles considÈrÈs comme des catÈgories
discrËtes. 
Nous notons $\sf C( \mathbb A,\mathbb B)$ la catÈgorie des foncteurs
\[
\mathbb B\op \times \mathbb A \ra \vect \corps.
\]
Posons
\[
\sf C  =  \sf C(\mathbb A, \mathbb A), \quad
\sf  C'  =  \sf C(\mathbb A, \mathbb B) \quad \mbox{et} \quad
\sf  C''  =  \sf C(\mathbb B, \mathbb B).
\]
Les produits tensoriels au-dessus de $\mathbb A$ et de $\mathbb B$ dÈfinissent
les structures de catÈgories monoÔdales (tressÈes) sur $\sf C$ et $\sf C''$
et les actions de $\sf C$ et $\sf C''$ sur $\sf C'$.
La catÈgorie $\sf C \ts \sf C''$ est la catÈgorie
$\sf C( \mathbb A \times \mathbb B, \mathbb A \times \mathbb B)$
des foncteurs
\[
\big(\mathbb A \times \mathbb B\big)\op \times \big(\mathbb A \times \mathbb B\big) \ra \vect \corps.
\]
Le foncteur
\[
\sf C( \mathbb A,\mathbb A) \times \sf C( \mathbb B,\mathbb B) \ra
\sf C( \mathbb A \times \mathbb B, \mathbb A \times \mathbb B)
\]
envoie $(L,M)$ sur le foncteur
\[
(A,B,A',B') \mapsto L(A,A') \ts_{\corps} M(B,B').
\]
}\end{exemple}

\subsection{DÈfinitions des bipolydules} \label{section_definition_bipolydules}

Soit $A$ et $A''$ deux $\ai$-algËbres de $\sf C$ et $\sf C''.$

\begin{definition}{\em \label{definition_A_n-A_m-bimodules}
Un $\mathrm A_n$-$\mathrm A_{n'}$-{\em bimodule sur $A$ et $A''$}
est un objet de $\gr \sf C'$
muni d'une famille de morphismes graduÈs dans $\gr \sf C'$
\[
m_{i,j} : A \tp i  \ts M \ts A'' \tp j \ra M,\quad
0 \leq i \leq n, \quad 0 \leq j \leq n',
\]
de degrÈ $1-i-j$, telles qu'une Èquation $(*''_{r,t})$  de la mÍme
forme que l'Èquation $(*_{r+1+t})$, $r+1+t \geq 1$, de la dÈfinition
\ref{definition_ai-algebre} est vÈrifiÈe pour tous $0 \leq r \leq n$ et $0 \leq t \leq n'.$ 
Si $M$ et $M'$ sont deux $\mathrm A_n$-$\mathrm A_{n'}$-bimodule sur $A$ et $A''$, un {\em morphisme}
\[
f : M \ra M'
\]
est une famille de morphismes graduÈs dans $\gr \sf C'$
\[
f_{i,j} : A \tp i  \ts M \ts A'' \tp j \ra M',\quad
0 \leq i \leq n, \quad 0 \leq j \leq n',
\]
de degrÈ $-i-j$,
vÈrifiant les ÈgalitÈs $(**'')_{r,t}$, $0 \leq r \leq n$ et $0 \leq t \leq n'$, 
des morphismes
\[
A^{\ts r}  \ts M \ts A''^{\tso t} \ra M', \quad 0 \leq r \leq n, \quad 0 \leq t \leq n',
\]
\[
\begin{array}{c}
\sum (-1)^{\alpha(-i-j)} m_{\alpha,\beta}(\Id \tp \alpha \ts f_{k,l} \ts \Id \tp \beta) = \hspace*{4cm} \\
\hspace*{1cm} \sum (-1)^{j + i(|m_\bullet|)} f_{\bullet,\bullet}(\Id \tp i \ts m_{\bullet} \ts \Id \tp j)
\end{array}
\]
o˘ $|m_\bullet|$ est le degrÈ de $m_\bullet$; il faut interprÈter convenablement les $m_\bullet$
par des $m^A_{\bullet}$, $m^{A''}_{\bullet}$ ou $m_{\bullet,\bullet}$ selon leur place.
La composition $g\circ f$ de deux morphismes $f$ et $g$ est dÈfinie par la suite
\[
(g \circ f)_n = \sum (-1)^{\alpha(-i-j)}
g_{i,j}(\Id \tp \alpha \ts f_{k,l} \ts \Id \tp \beta), \quad n \geq 1.
\]
}\end{definition}

\begin{definition}{\em  
Un $A$-$A''$-{\em bipolydule dans $\sf C'$} \index{bipolydule} (appelÈ communÈment {\em $\ai$-bimodule sur $A$ et $A''$} dans la littÈrature)
est un objet de $\gr \sf C'$
muni d'une famille de morphismes graduÈs dans $\gr \sf C'$
\[
m_{i,j} : A \tp i  \ts M \ts A'' \tp j \ra M,\quad
i,j\geq 0,
\]
de degrÈ $1-i-j$, telles que l'Èquation $(*''_{n,n'})$, $n,n' \geq 0$ est vÈrifiÈe.
Si $M$ et $M'$ sont deux $A$-$A''$-polydules, un {\em morphisme}
\[
f : M \ra M'
\]
est une famille de morphismes graduÈs dans $\gr \sf C'$
tels que l'ÈgalitÈ $(**'')_{n,n'}$, $n +1 +n'\geq 1$, est vÈrifiÈe.
La composition $g\circ f$ de deux $\ai$-morphismes $f$ et $g$ est dÈfinie
par les mÍmes formules que dans le cas des morphismes de $\mathrm A_n$-$\mathrm A_{n'}$-bimodules
sur $A$ et $A''$. Nous obtenons ainsi une {\em catÈgorie} $\aiMod (A,A'')$.
\indexnotation{aiModAA'} La lettre $\mathsf N$ de
$\aiMod$ remplace la lettre $\mathsf M$ dans $\aiModu$ et se rapporte au $\mathsf N$ dans
``$\ai$-bimodules $\mathsf N$on (nÈcessairement) unitaires''.
}\end{definition}

Nous supposons dÈsormais que $A$ et ${A''}$ sont
augmentÈes.

\begin{definition}{\em  \label{definition_bipolydules}
Un $A$-${A''}$-{bipolydule} est {\em strictement unitaire}
\index{strictement unitaire!bipolydule}
si pour tous $i,j \geq 0$, on a
\[
m_{i,j} (\Id \tp \alpha \ts \unite \ts \Id \tp \beta) = 0, \quad \alpha \neq i, \quad 
(i,j) \notin \{(0,1),(1,0)\}
\]
et
\[
m_{1,0} \circ (\unite \ts \Id)  = m_{0,1} \circ (\Id \ts \unite) = \Id.
\]
Nous notons $\aiModu (A,{A''})$ \indexnotation{aiModuAA'}
la catÈgorie des $A$-${A''}$-bipolydules strictement unitaires.
Elle est isomorphe ‡ la catÈgorie des $\b A$-$\b {A''}$-bipolydules, o˘ $\b A$ et $\b {A''}$ sont
les rÈductions de $A$ et ${A''}.$ 
}\end{definition}

{\noindent \bf Construction bar}\\ 
\label{section_construction_bar_bipolydules}

Nous dÈfinissons des bijections 
\[
\begin{array}{rcl}
\Hom ((S\b A)\tp{i} \ts SM \ts (S\b {A''})\tp{j}, SM) &
\arr{\sim} & \Hom (\b A \tp i \ts M \ts \b {A''} \tp{j}, M),\\
m_{i,j} & \mapsto & b_{i,j}
\end{array}
\]
\[
\begin{array}{rcl}
\Hom ((S\b A)\tp{i} \ts SM \ts (S\b {A''})\tp{j}, SM) &
\arr{\sim} & \Hom (\b A \tp i \ts M \ts \b {A''} \tp{j}, M),\\
f_{i,j} & \mapsto & F_{i,j}
\end{array}
\]
par les
relations
\[
\si \circ b_{i,j} = -  m_{i,j} \circ \si \tp{i+1+j} \quad \mbox{et} \quad
\si \circ F_{i,j} = (-1)^{|F_{i,j}|} f_{i,j} \circ \si \tp{i+1+j}.
\]
Ces bijections dÈfinissent le foncteur {\em construction bar}, pleinement fidËle,
\index{construction bar}
\[
B : \aiModu (A,{A''}) \arr{} \coModcu (\Ba A,\Ba {A''}),
\]
o˘ $\coModcu (\Ba A,\Ba {A''})$ est la {\em catÈgorie} des objets de $\gr \sf C'$
munis de structures de $\Ba A$-$\Ba A''$-bicomodule
co-unitaire diffÈrentiel graduÈ cocomplet.
Son image est formÈe des objets qui sont presque colibres.

\subsection{CatÈgorie dÈrivÈe des $\ai$-bimodules}

Soit $A$ et $A''$ deux $\ai$-algËbres augmentÈes dans $\sf C$ et $\sf C''.$
Dans cette section, nous dÈfinissons la catÈgorie dÈrivÈe des
$A$-${A''}$-bipolydules strictement unitaires, puis nous en donnons plusieurs descriptions.\\

{\noindent \bf Structure de catÈgorie de modËles sur $\coModcu (\Ba A,\Ba {A''})$}\\

Notons $(\Ba A)\op$ la cogËbre opposÈe de $\Ba A$ dÈfinie ‡ l'aide du tressage de~$\sf C$.
L'objet  $(\Ba A)\op \ts\Ba {A''}$ de
$\sf C \ts \sf C''$ est une cogËbre diffÈrentielle graduÈe cocomplËte.
Notons qu'elle n'est pas cotensorielle en gÈnÈral. La catÈgorie
$\coModcu ((\Ba A)\op \ts \Ba {A''})$ est munie de sa structure canonique
de catÈgorie de modËles (\ref{definition_cmf_canonique_coModcu}).
La catÈgorie $\coModcu (\Ba A,\Ba {A''})$
devient une catÈgorie de modËles gr‚ce ‡ l'isomorphisme de catÈgories
\[
\coModcu (\Ba A,\Ba {A''}) \ra \coModcu ((\Ba A)\op \ts \Ba {A''}).
\]
Nous allons maintenant montrer que les objets fibrants de
$\coModcu (\Ba A,\Ba {A''})$ sont exactement les facteurs directs d'objets presque
colibres.\\

{\noindent \bf Une cochaÓne tordante acyclique}\\

Notons $(\E A)\op$ l'algËbre opposÈe de $\E A$ dÈfinie ‡ l'aide du tressage de~$\sf C$.
L'objet  $(\E A)\op \ts\E {A''}$ de
$\sf C \ts \sf C''$ est une algËbre diffÈrentielle graduÈe.
Munissons la catÈgorie $\Modu ((\E A)\op \ts \E {A''})$ de la structure
de catÈgorie de modËles du thÈorËme \ref{theoreme_cmf_Modu}.
Soit $\Modu (\E A, \E {A''})$ la {\em catÈgorie} des
bimodules diffÈrentiels graduÈs unitaires. \indexnotation{ModuAA'}
La catÈgorie $\Modu (\E A,\E {A''})$
devient une catÈgorie de modËles gr‚ce ‡ l'isomorphisme de catÈgories
\[
\Modu (\E A,\E {A''}) \ra \Modu ((\E A)\op \ts \E {A''}).
\]
Nous allons construire une cochaÓne tordante admissible acyclique
\[
\tau :  (\Ba A)\op \ts \Ba {A''} \ra (\E A)\op \ts \E {A''}.
\]
Il s'ensuivra (\ref{corollaire_theoreme_cmf_coModcu}) que le couple de foncteurs adjoints associÈ ‡ $\tau$ 
(voir \ref{section_cochaines_tordantes})
\[
(L,R) : \coModcu \big((\Ba A)\op \ts \Ba {A''}\big) \ra \Modu \big((\E A)\op \ts \E {A''}\big)
\]
est une Èquivalence de Quillen .\\

La cochaÓne tordante universelle (\ref{cochaines_universelles})
\[
\tau_{\Ba A} : \Ba A \ra \Oma \Ba A = \E A
\]
induit une cochaÓne tordante
\[
\tau_{\Ba A}' : (\Ba A)\op \ra (\E A)\op.
\]
Nous vÈrifions que
\[
\tau = \tau_{\Ba A} \ts \counite \circ \unite + \counite \circ \unite \ts \tau_{\Ba {A''}}
:  (\Ba A)\op \ts \Ba {A''} \ra (\E A)\op \ts \E {A''},
\]
o˘ les symboles $\unite$ dÈsignent les (co)unitÈs de $\Ba A$, $\Ba {A''}$,
$\E A$ et $\E {A''}$, est une cochaÓne
tordante admissible. Par le critËre d'acyclicitÈ des cochaÓnes tordantes
(\ref{caracterisation_cochaines_acycliques}),
l'objet de $\sf C \ts \sf C''$
\[
\begin{array}{l}
\Big((\Ba A)\op \ts \Ba {A''}\Big) \tw \Big((\E A)\op \ts \E {A''}\Big) = \\
\hspace*{2cm}\Big( (\Ba A)\op \ts_{\tau'_{\Ba A}} (\E A)\op \ts
(\Ba {A''})\op \ts_{\tau_{\Ba {A''}}} \E {A''}\Big)
\end{array}
\]
est quasi-isomorphe ‡ $e_{\sf C}\ts e_{\sf C''} = e_{\sf C \ts \sf C''}$.
La cochaÓne tordante $\tau$ est donc acyclique.\\

{\noindent \bf Objets fibrants de $\coModcu \big((\Ba A)\op \ts \Ba {A''}\big)$}\\

Comme dans le cas des polydules sur une $\ai$-algËbre augmentÈe (voir \ref{theoreme_cmf_aiModu}),
nous montrons gr‚ce ‡ la thÈorie de l'obstruction (\ref{section_obstruction_bipolydules})
que la catÈgorie des $A$-$A''$-bipolydules est munie d'une structure de
``catÈgorie de modËles sans limites'' : les Èquivalences faibles, les cofibrations et les fibrations sont
dÈfinies de la mÍme maniËre que dans le cas des $A$-polydules (\ref{theoreme_cmf_aiModu}).
Par le mÍme raisonnement que celui de la preuve de la proposition \ref{proposition2_categorie_derivee},
nous montrons que les objets fibrants de la catÈgorie de modËles $\coModcu (\Ba A,\Ba A'')$
sont exactement les facteurs directs des comodules presque colibres.\\

{\noindent \bf La catÈgorie dÈrivÈe}\\

La construction bar
\[
B : \aiModu (A,{A''}) \ra \coModcu (\Ba A, \Ba {A''})
\]
est un foncteur pleinement fidËle. La clÙture par rÈtracts de son image est la
sous-catÈgorie des objets fibrants et cofibrants.
La proposition \ref{proposition_cmf} et la compatibilitÈ de la construction bar ‡ l'homotopie et
aux Èquivalences faibles montre que la dÈfinition suivante a un sens.

\begin{definition}
\label{definition_categorie_derivee_augmentee_bipolydules}{\em
La {\em catÈgorie} $\ch_\infty (A,A'')$ est la catÈgorie $\aiModu (A,A'')/\! \sim$, o˘ $\sim$
est la relation d'homotopie.
La {\em catÈgorie dÈrivÈe $\cd_\infty (A,A'')$}
\index{categorie@{catÈgorie}!derivee@{dÈrivÈe}} est la loca\-lisation par rapport aux $\ai$-quasi-isomorphismes
de la catÈgorie $\aiModu (A,A'').$
}\end{definition}
Par la proposition \ref{proposition_cmf}, nous avons un isomorphisme
\[
\ch_\infty (A,A'') \ra \cd_\infty (A,A'').
\]
Nous avons un foncteur pleinement fidËle
\[
I : \Modu (\E A,\E {A''}) \ra \aiModust (A,{A''}), \quad M \ra S^{-1}M,
\]
o˘ $\aiModust (A,{A''})$ \indexnotation{aiModustAA'}
est la catÈgorie des $A$-${A''}$-polydules strictement
unitaires dont les morphismes sont les $\ai$-morphismes stricts. L'image de ce foncteur est formÈe
des $A$-$A''$-bipolydules $M$ dont les morphismes
\[
m_{i,j} : A \tp i \ts M \ts A'' \tp j \ra M, \quad i,j \geq 0,
\]
sont nuls si les deux entiers $i$ et $j$ sont diffÈrents de $0.$
Rappelons que le foncteur analogue dans le cas des polydules est un isomorphisme
(\ref{A-polydules_EA-modules}).

\begin{lemme}
La composition des foncteurs 
\[
J : \Modu (\E A,\E A'') \arr{I} \aiModust (A,A'') \hookrightarrow \aiModu (A,A'')
\]
induit une Èquivalence $\cd (\E A,\E A'') \ra \cd_\infty (A,A'').$
\end{lemme}

\dem 
Nous avons un diagramme commutatif
\[
\xymatrix{
\Modu (\E A,\E A'') \ar[r]^{I} \ar[d]_{R} & \aiModust (A,A'') \ar@{^(->}[d] \\
\coModcu (\Ba A,\Ba A'') & \aiModu (A,A'') \ar[l]^{B}
}
\]
o˘ $R$ et $B$ induisent des Èquivalences dans les catÈgories dÈrivÈes.
Cela montre que le foncteur induit par $J$ est pleinement fidËle. Montrons
qu'il est essentiellement surjectif.
Soit $M$ un $A$-$A''$-bipolydule. Le morphisme d'adjonction
\[
BM \ra RL BM = \Ba A \ts_{\tau_{\Ba A}} \E A \ts_{\tau_{\Ba A}} BM \ts_{\tau_{\Ba A''}} \E A'' 
\ts_{\tau_{\Ba A''}} \Ba A''
\]
est une Èquivalence faible. Le bicomodule $RL BM$ est la construction bar du $A$-$A''$-polydule
\[
M' = S^{-1}\big(\E A \ts_{\tau_{\Ba A}} BM \ts_{\tau_{\Ba A''}} \E A''\big).
\]
Nous avons alors un $\ai$-quasi-isomorphisme de $A$-$A''$-bipolydules
\[
M \ra  M'
\]
et, comme $M'$ est dans l'image de $J$, nous avons le rÈsultat.\findem

%% file: Unite.tex
{\bf \noindent Introduction}\\

\noindent Les $\ai$-espaces de \cite{Stasheff63a} sont munis d'unitÈs strictes.
Dans le cadre algÈbrique, la notion correspondante a ÈtÈ dÈfinie en
(\ref{definition_unite_stricte}).
Lorsque $A$ est une $\ai$-algËbre strictement unitaire, certaines propriÈtÈs des algËbres associatives
unitaires pourront Ítre gÈnÈralisÈe ‡ $A$.
Par exemple, nous montrerons l'analogue de l'isomorphisme
\[
M \ts_{B} {B} \ra M,
\]
lorsque $B$ est une algËbre associative unitaire et $M$ un $B$-module unitaire
(voir la gÈnÈralisation en \ref{lemme2_adjonction_foncteurs_standard} dans le
chapitre \ref{chapitre_Cat_der}).
Cependant, les $\ai$-algËbres (en fait $\ai$-catÈgories) apparaissant en
gÈomÈtrie \cite{Fukaya93} ne sont pas strictement unitaires mais {\em homologiquement
unitaires}, i.~e.~$H^*A$ munie de la multiplication induite par $m_2$ est une algËbre
graduÈe unitaire.
Le but de ce chapitre est de montrer que d'un point de vue homotopique, il n'y a pas de
diffÈrence entre les unitÈs strictes et les unitÈs homologiques. Plus prÈcisÈment,
nous montrerons que
{\em la sous-catÈgorie des $\ai$-algËbres homologiquement unitaires dont les morphismes
sont les $\ai$-morphismes homologiquement unitaires et la sous-catÈgorie des $\ai$-algËbres
strictement unitaires dont les morphismes
sont les $\ai$-morphismes strictement unitaires deviennent Èquivalentes aprËs passage ‡
l'homotopie} (\ref{corollaire_ai-algebres_strictement_unitaires}).\\

{\bf \noindent Plan du chapitre}\\

\noindent Ce chapitre est divisÈ en trois sections. Dans la section \ref{section_definition_unites},
nous dÈfinissons 
les unitÈs homologiques relatives aux $\ai$-structures.
Dans la section \ref{section_strictification_unitaire_aia}, nous montrons le rÈsultat ÈnoncÈ
ci-dessus.
Dans la section \ref{section_strictification_unitaire_(bi)polydules},
nous comparons les diffÈrents types de compatibilitÈs
aux unitÈs des (bi)polydules.

\section{DÈfinitions}\label{section_definition_unites}

Soit $\sf C$ une catÈgorie de base telle que dans le chapitre \ref{chapitre_Homot_aialg}.
Soit $A$ une $\ai$-algËbre sur $\sf C$ et soit 
\[
\mu : H^*A \ts H^*A \ra H^*A
\]
le morphisme induit par $m_2$.

\begin{definition}
{\em 
Un morphisme $\unite^A : e \ra A$ dans $\gr \sf C$ est une {\em unitÈ homologique}\index{unite@{unitÈ}!homologique}
si $m_1 \circ \unite = 0$ et s'il induit une unitÈ
pour l'algËbre graduÈe associative $(H^*A,\mu)$. 
Si $A$ est munie d'une unitÈ homologique, nous dirons qu'elle est {\em homologiquement unitaire}.
Si $A$ et $A'$ sont deux $\ai$-algËbres homologiquement unitaires, un $\ai$-morphisme
$f : A \ra A'$ est {\em homologiquement unitaire}
si $f_1$ induit un morphisme unitaire
\[
 H^*A \arr{\sim} H^*A'.
\]
}\end{definition}

\begin{remarque}{\em
L'unitÈ $e \ra A$ d'une $\ai$-algËbre strictement unitaire (\ref{definition_unite_stricte}) est
clairement une unitÈ homologique.
Un morphisme strictement unitaire d'$\ai$-algËbre strictement unitaire est homologiquement unitaire.
}\end{remarque}

On trouve dans les travaux de K.~Fukaya \cite{Fukaya01}
et V.~Lyubashenko \cite{Lyubashenko02}
d'autres relËvements de la notion d'unitaritÈ. Une $\ai$-algËbre munie d'une
``unitÈ homotopique'' (dÈfinie dans  \cite{Fukaya01} ‡
l'aide d'homotopies supÈrieures, voir aussi \cite{Fukaya01a}) donne une
``$\ai$-algËbre unitale'' au sens de \cite{Lyubashenko02}.
Le relËvement de la notion d'unitaritÈ de V.~Lyubashenko \cite{Lyubashenko02} se spÈcialise
‡ notre notion d'unitaritÈ homologique si on travaille sur un corps
(V.~Lyubashenko travaille sur un anneau commutatif quelconque). Remarquons que
l'unitaritÈ homologique n'est pas du type {\em ``‡ homotopie prËs''} : elle n'est pas dÈfinie
‡ l'aide d'homotopies supÈrieures vÈrifiant des conditions de cohÈrence. Elle est cependant une
notion valide puisque (comme nous le verrons dans ce chapitre)
la localisation de la catÈgorie des $\ai$-algËbres homologiquement unitaires par rapport
aux $\ai$-quasi-isomorphismes est Èquivalente ‡ la localisation de la catÈgorie des algËbres
unitaires par rapport aux quasi-isomorphismes.

\begin{definition}{\em
Si $f$ et $f'$ sont deux morphismes homotopiquement unitaires $A \ra A',$ une homotopie
$h$ entre $f$ et $f'$ est {\em strictement unitaire} si  \index{strictement unitaire!homotopie}
\[
h_i (\Id \tp j \ts \unite \ts \Id \tp l) = 0, \quad i \geq 1 \mbox{ et } j+1+l = i.
\]
}\end{definition}

\begin{remarque}{\em
 Si $A$ est une $\ai$-algËbre homologiquement unitaire et $H^*A$ est
un modËle minimal pour $A$ (\ref{corollaire_modele_minimal}), l'unitÈ homologique $\unite^A$ induit une
unitÈ homologique $\unite^{H^*A} : e \ra H^*A$ qui vÈrifie en outre
\[
m^{H^*A}_2(\unite^{H^*A} \ts \Id) = m^{H^*A}_2 (\Id \ts \unite^{H^*A}) = \Id.
\]

 Soit $f : A \ra A'$ un morphisme homologiquement
unitaire et $H^*A$ et $H^*A'$ des modËles minimaux de $A$ et $A'.$
Nous rappelons (\ref{corollaire_modele_minimal})
qu'il existe
des $\ai$-quasi-isomorphismes
\[
i : H^*A \ra A \quad \mbox{et} \quad i' : H^*A' \ra A'.
\]
Par le point {\it b} du corollaire \ref{corollaire_cmf_aia}, il
existe un inverse ‡ homotopie prËs $p'$ de $i'.$
Le morphisme  $g = p' \circ f_1 \circ i$ vÈrifie en outre $g_1 \unite_{H^*A} = \unite_{H^*A'}.$
}\end{remarque}

\section{$\ai$-algËbres homologiquement unitaires} \label{section_strictification_unitaire_aia}

\noindent Cette section est divisÈe en quatre sous-sections.

Dans la sous-section \ref{section_stritification_unitaire_ai-algebres},
nous donnons deux dÈmonstrations du fait que toute $\ai$-algËbre minimale homologiquement
unitaire est isomorphe ‡ une $\ai$-algËbre strictement unitaire.
 La premiËre de ces
dÈmonstrations est inspirÈe de la thÈorie des dÈformations des algËbres graduÈes et
n'est valable qu'en caractÈristique nulle.
La seconde est basÈe sur la thÈorie de l'obstruction des $\ai$-algËbres
minimales (voir l'appendice \ref{section_obstruction_unite}).

Dans les sous-sections \ref{section_stritification_unitaire_ai-morphismes} et 
\ref{section_stritification_unitaire_homotopies}, nous dÈmontrons, ‡ l'aide de la thÈorie
de l'obstruction, qu'on peut rendre strictement unitaire tout $\ai$-morphisme homologiquement
unitaire entre $\ai$-algËbres strictement unitaires et toute homotopie entre $\ai$-morphismes.

Dans la sous-section \ref{section_modele_minimal_strictement_unitaire},
nous montrons que toute $\ai$-algËbre strictement unitaire $A$
admet un modËle minimal strictement unitaire $A'$ et des $\ai$-quasi-isomorphismes
strictement unitaires
\[
A' \ra A \quad \mbox{et} \quad A \ra A'.
\]
Nous dÈduirons de ce rÈsultat et des sous-sections prÈcÈdentes le rÈsultat principal de
ce chapitre (\ref{corollaire_ai-algebres_strictement_unitaires}) : la catÈgorie $\big(\aia\big)_{hu}$
des $\ai$-algËbres homologiquement unitaires dont les morphismes
sont les $\ai$-morphismes homologiquement unitaires et sa sous-catÈgorie non pleine
$\big(\aia\big)_{su}$
des $\ai$-algËbres
strictement unitaires dont les morphismes
sont les $\ai$-morphismes strictement unitaires deviennent Èquivalentes aprËs passage
‡ l'homotopie.

\subsection{Strictification unitaire des $\ai$-algËbres}
\label{section_stritification_unitaire_ai-algebres}

%
%
%
%
%
\begin{theoreme}[A.~Lazarev \cite{Lazarev02}, P.~Seidel \cite{Seidel02}] \label{theoreme_unite_ai-algebre}
Toute $\ai$-algËbre minimale
homologiquement unitaire est isomorphe
‡ une $\ai$-algËbre minimale strictement unitaire.
\end{theoreme}
Le thÈorËme a ÈtÈ dÈmontrÈ de faÁon indÈpendante
par P.~Seidel \cite{Seidel02}, qui utilise la mÍme mÈthode que nous, ainsi que
par A. Lazarev \cite{Lazarev02}. Notre premiËre dÈmonstration utilisera 
les dÈformations et n'est valable qu'en
caractÈristique zÈro. Elle nous donne l'existence de l'$\ai$-algËbre minimale strictement unitaire.
La seconde dÈmonstration est basÈe sur les lemmes d'obstruction de l'appendice 
\ref{section_obstruction_unite}. Elle prÈcise les choix possibles de l'$\ai$-algËbre minimale
strictement unitaire.

Les deux dÈmonstrations sont liÈes :
pour un $m_2$ donnÈ, le complexe de Hochschild $C^*(A,A)$ (voir l'appendice
\ref{section_obstruction_unite}) contrÙle l'obstruction ‡ la construction par
rÈcurrence des $m_i$, $i \geq 3$, d'une structure d'$\ai$-algËbre minimale sur $A$ et il
est aussi l'algËbre de Lie diffÈrentielle graduÈe qui dÈcrit le problËme des dÈformations de l'algËbre
$(A,m_2).$ Nous renvoyons aux articles
\cite{Schlessinger85} et \cite{Kontsevich00} concernant ce point.

\begin{corollaire} \label{corollaire_unite_ai-algebre}
Toute $\ai$-algËbre homologiquement unitaire est homotopiquement Èquivalente
‡ une $\ai$-algËbre strictement unitaire.
\end{corollaire}

\dem Soit $A$ une $\ai$-algËbre homologiquement unitaire et soit
$A'$ un modËle minimal de $A$. Nous savons que $A$ et $A'$ sont homotopiquement Èquivalents.
Le rÈsultat se dÈduit alors du thÈorËme \ref{theoreme_unite_ai-algebre}
appliquÈ ‡ $A'.$ \findem \\

\begin{remarque}{\em
Nous montrerons ‡ la fin de ce chapitre (\ref{proposition_modele_minimal_strictement_unitaire}) que toute
$\ai$-algËbre strictement unitaire $A$ admet un modËle minimal
strictement unitaire $A'$ tel que l'$\ai$-quasi-morphisme
\[
A' \ra A
\]
est strictement unitaire.
}\end{remarque}

\newpage

{\noindent \em PremiËre dÈmonstration du thÈorËme \ref{theoreme_unite_ai-algebre} : }\\

{\noindent \bf Rappel sur les dÈformations }\\

Supposons que la caractÈristique de $\corps$ est nulle.
Soit $(\mathfrak g,\delta,[\?,\?])$ une $\corps$-algËbre de Lie diffÈrentielle graduÈe nilpotente,
i.~e.~il existe un entier $N \geq 1$ tel que
\[
\ad X_1 \ad X_2 \hdots \ad X_N = 0, \quad X_1,\hdots ,X_N \in \mathfrak{g}.
\]
On note $\mathsf{MC}(\mathfrak g)$ les ÈlÈments  $X \in \mathfrak g$ de degrÈ $+1$
qui sont solutions de l'Èquation de Maurer-Cartan \index{Maurer-Cartan (Èquation)}
\[
\del{X} + \frac{1}{2}[X,X] = 0.
\]
Soit $\Gamma$ le groupe nilpotent associÈ ‡ $\mathfrak{g}^0$.
Il agit sur $\mathfrak{g}^1$ par transformations affines, c'est-‡-dire, par l'exponentiation
de l'action de son algËbre de Lie
\[
g.x = \delta(g) + [g,x], \quad g \in \mathfrak{g}^0, x \in \mathfrak{g}^1.
\]
Cette action conserve $\mathsf{MC}(\mathfrak g)$ et on a l'ensemble
\[
\mathsf{MC}(\mathfrak g)/\!\sim\, = \mathsf{MC}(\mathfrak g)/\Gamma.
\]
On rappelle  \cite{Goldman90} le rÈsultat suivant.
\begin{theoreme} \label{theoreme_Goldman_Milson}
Si $\mathfrak h$ est une algËbre de Lie diffÈrentielle graduÈe nilpotente,
une Èquivalence d'homotopie $f : \mathfrak h \ra \mathfrak g$
induit une bijection
\[
\mathsf{MC}(\mathfrak h)/\!\sim\, \arr{\sim} \mathsf{MC}(\mathfrak g)/\!\sim\,.
\]
\findem
\end{theoreme}
Si $\mathfrak g'$ est une algËbre de Lie pronilpotente
(i.~e.~qui est la limite d'algËbres nilpotentes $\mathfrak{g_i}$,
$i\geq 0$) on dÈfinit
\[
\mathsf{MC}(\mathfrak g')= \lim \mathsf{MC}(\mathfrak g_i)
\quad \mbox{et} \quad
\mathsf{MC}(\mathfrak g')/\!\sim\, = \lim \Big(\mathsf{MC}(\mathfrak g_i)/\Gamma_i\Big).
\]

{\noindent \bf Lien avec les $\ai$-algËbres}\\

Soit $(A,\mu)$ une $\corps$-algËbre graduÈe associative unitaire. L'application
\[
(D,D') \mapsto [D,D'] = D\circ D' - (-1)^{pq}D' \circ D,
\]
o˘ $D$ et $D'$ sont homogËnes de degrÈ $p$ et $q$,
munit le complexe $(\coder (BA)\+,\delta)$ d'une structure d'algËbre de Lie
diffÈrentielle graduÈe. Notons $LA$ cette algËbre de Lie.
Nous avons un isomorphisme de complexes
\[
LA \ra SC(A,A),
\]
o˘ $C(A,A)$ est le complexe de Hochschild (voir l'appendice \ref{section_obstruction_unite}).
Il envoie le crochet de Lie de $LA$ sur le crochet de Gerstenhaber \cite{Gerstenhaber63}.
Soit $L^{\geq n}A  \subset LA$, $n \geq 3$, la sous-algËbre de Lie
\[
S \Big( \prod_{i\geq n}\Hom_{\sz \gr \sf C}(A\tp i,A) \Big).
\]
Les sous-algËbres $L^{\geq n}A$, $n \geq 4$, sont des idÈaux de $L^{\geq 3}A$ et
nous avons
\[
L^{\geq 3}A = \lim_{n \geq 4} \mathfrak{g}_n,
\]
o˘ $\mathfrak{g}_n$ est l'algËbre $L^{\geq 3}A /L^{\geq n}A.$
Comme nous avons 
\[
[L^{\geq n}A,L^{\geq n'}A] \subset L^{\geq n+n'-1}A, \quad n,n' \geq 1,
\] 
les algËbres de Lie $\mathfrak{g}_n$ sont nilpotentes et $L^{\geq 3}A$ est pronilpotente.
Le sous-complexe rÈduit $S\b C(A,A)$ est une sous-algËbre de Lie de $LA$ pour le crochet de Gerstenhaber.
Nous la notons $\b LA$.
Rappelons que l'inclusion $\b LA \hookrightarrow L A$
est une Èquivalence d'homotopie (voir \cite[Chap.~IX]{Cartan56}).
Par le thÈorËme \ref{theoreme_Goldman_Milson}, nous avons une bijection
\[
\Theta : \mathsf{MC}(\b L^{\geq 3}A)/\!\sim\, \arr{\sim} \mathsf{MC}(L^{\geq 3}A)/\!\sim\,,
\]
o˘ $\b L^{\geq 3}A = \b LA \cap L^{\geq 3}A.$
Un ÈlÈment $b' \in L^{\geq 3}A$ est dans $\mathsf{MC}(L^{\geq 3}A)$
si et seulement si $b = b'+b_2$ (o˘ $b_2$
correspond ‡ $m_2 = \mu$) est une diffÈrentielle
de $(BA)\+.$  En d'autres termes, nous avons une bijection entre 
$\mathsf{MC}(L^{\geq 3}A)$ et l'ensemble
des structures d'$\ai$-algËbre minimale sur $A$ dont la multiplication $m_2$ vaut $\mu$.
Sous cette bijection, les classes d'Èquivalence de $\mathsf{MC}(L^{\geq 3}A)$
correspondent aux classes d'isomorphie de structures $\ai$ minimales tel que $m_2$ vaut $\mu$.
Remarquons qu'un ÈlÈment $b'' \in \mathsf{MC}(L^{\geq 3}A)$ appartient ‡ la sous-algËbre $\b L^{\geq 3}A$
si et seulement si l'$\ai$-structure correspondant ‡ $b''$ est
strictement unitaire sur $A$. Nous dÈduisons alors de la bijection $\Theta$ que
toute $\ai$-structure (dont le $m_2$ vaut $\mu$) homologiquement unitaire sur $A$
est isomorphe ‡ une $\ai$-structure strictement unitaire.\\

{\noindent \em DeuxiËme dÈmonstration du thÈorËme \ref{theoreme_unite_ai-algebre} : }\\

La caractÈristique de $\corps$ est quelconque.

\begin{lemme}
Soit $A$ une $\ai$-algËbre minimale. Soit $n$ un entier $\geq 2$ et
\[
f_n : A \tp n \ra A
\]
un morphisme graduÈ de degrÈ $1-n$.
Il existe une $\ai$-algËbre minimale $A'$, $\ai$-isomorphe ‡ $A$,
dont l'objet graduÈ sous-jacent est $A$ et dont
les multiplications $m'_i$, $i\geq 2$, sont telles
que
\[
m'_i = m_i \quad \mbox{si}\quad  i\leq n \quad \mbox{et} \quad m'_{n+1} = m_{n+1} + \delta_{Hoch}(f_n).
\]
\end{lemme}

\dem Soit le morphisme de cogËbres graduÈes
\[
F : BA \ra BA
\]
dÈterminÈ par la suite
\[
(\Id_{SA}, 0 , \hdots ,0 , F_n, 0 \hdots),
\]
o˘ $F_n$ est donnÈ par la bijection $F_n \lra f_n$ de la section \ref{construction_bar_cobar}.
Le morphisme $F$ est un isomorphisme. Posons 
\[
b' = F \circ b^A \circ F^{-1}.
\]
C'est une diffÈrentielle sur $\ctr SA$. La cogËbre $(\ctr SA,b')$ est donc la construction bar
d'une $\ai$-algËbre $A'$, $\ai$-isomorphe ‡ $A$,
dont l'objet graduÈ sous-jacent est $A$. Il reste ‡ vÈrifier les conditions sur les
multiplications. La matrice 
du morphisme de cogËbres graduÈes
\[
F : \ctr (SA) = \bigoplus_{p\geq 1} (SA)\tp p \arr{\sim} \ctr(SA) = \bigoplus_{q\geq 1} (SA)\tp q
\]
est triangulaire supÈrieure et sa
diagonale est formÈe d'identitÈs. La matrice de $F^{-1}$
est donc de la mÍme forme. De plus la restriction de $F$
‡ 
\[
\pctr{n-1} SA  = \bigoplus_{1\leq p \leq n-1} (SA)^p
\]
est l'identitÈ. Il en est donc de mÍme
pour son inverse. La matrice de la diffÈrentielle $b^A$ est strictement triangulaire supÈrieure
puisque $b^A_1$ est nul.
Le calcul montre alors que
\[
b'_i = F_1 b^A_i (F^{-1})_1 ,\quad \mbox{pour}\quad i \leq n,
\]
\[
b'_{n+1}  =  F_1 b^A_{n+1} (F^{-1})_1 +
F_1 b^A_2 (F^{-1})_n +
F_n b^A_2 (F^{-1})_1.
\]
Nous dÈduisons le rÈsultat des ÈgalitÈs 
\[
(F^{-1})_n = -F_n \quad \mbox{et} \quad F_1 = F^{-1}_1 =\Id_{SA}.
\] \findem\\

DÈmontrons maintenant le thÈorËme \ref{theoreme_unite_ai-algebre}.
Nous raisonnons par rÈcurrence sur $n$. Soit $n\geq 2$. Supposons que $A$ est une
$\ai$-algËbre telle que, pour tout $3 \leq i \leq n,$ on a
\[
m_i(\Id \tp j \ts \unite \ts \Id \tp k) = 0, \quad j + k = n.
\]
Ceci est Èquivalent ‡ demander que les $m_i$, $3 \leq i \leq n$, soient des
ÈlÈments du sous-complexe de Hochschild rÈduit $\b C(A,A)$ (voir \ref{section_obstruction_unite}).
Montrons que nous pouvons construire
une $\ai$-algËbre $A'$, $\ai$-isomorphe ‡ $A$, dont l'objet graduÈ sous-jacent est $A$ et dont
les multiplications $m'_i$, $3 \leq i \leq n+1,$ sont des ÈlÈments de $\b C(A,A)$.
Par hypothËse sur les $m_i$, $3 \leq i \leq n$, le cycle de Hochschild
$r(m_3, \cdots , m_{n-1})$ du lemme \ref{lemme_obstruction_structure_minimale} appartient ‡ $\b C(A,A)$.
Comme $A$ est une $\ai$-algËbre, nous savons par le lemme \ref{lemme_obstruction_structure_minimale}
que
\[
\delta_{Hoch}(m_{n+1}) + r(m_3, \cdots , m_{n}) = 0
\]
et que l'ÈlÈment $r(m_3, \cdots , m_{n})$ est un cycle de Hochschild.
Ainsi, l'ÈlÈment
\[
(m_{n+1} , s r(m_3, \cdots , m_{n}))
\]
du cÙne $C$
sur l'inclusion $\b C(A,A) \hookrightarrow C(A,A)$ est un cycle.
Comme $C$ est acyclique, cet ÈlÈment
est le bord d'un ÈlÈment
$(f_{n},sm'_{n+1})$. En d'autres termes, il existe des ÈlÈments
\[
m'_{n+1} \in \Hom_{\sz \gr \sf C} (\b A \tp{n+1} ,A)\quad \mbox{et} \quad
f_{n} \in \Hom_{\sz \gr \sf C} (A \tp{n} ,A)
\]
tels que
\[
\delta_{Hoch} (f_n) + m'_{n+1} = m_{n} \quad \mbox{et} \quad
\delta_{Hoch}(m'_{n+1}) + r(m_3, \cdots , m_{n}) = 0.
\]
Par le lemme prÈcÈdent appliquÈ ‡ l'$\ai$-algËbre $A$ et
au morphisme $-f_n$, il existe une $\ai$-algËbre $A'$, $\ai$-isomorphe ‡ $A$,
telle que nous avons, pour tout $3 \leq i \leq n+1,$
\[
m'_i(\Id \tp j \ts \unite \ts \Id \tp k) = 0, \quad j + k = n.
\]
\findem

\subsection{Strictification unitaire des $\ai$-morphismes}
\label{section_stritification_unitaire_ai-morphismes}

\begin{theoreme} \label{theoreme_unite_ai-morphisme}
Un morphisme d'$\ai$-algËbres minimales strictement unitaires qui est homologiquement unitaire
 est homotope ‡ un morphisme strictement unitaire.
\end{theoreme}

\begin{lemme} \label{lemme2_unite}
Soit $A$ et $A'$ deux $\ai$-algËbres minimales et $f : A \ra A'$ un
$\ai$-morphisme. Soit $n$ un entier $\geq 2$ et 
\[
h_n : A \tp n \ra A
\]
un morphisme graduÈ de degrÈ~$-n$.
Il existe un $\ai$-morphisme $f' : A \ra A'$ homotope ‡ $f$ tel que
\[
f'_i = f_i \quad \mbox{si} \quad i\leq n \quad \mbox{et} \quad f'_{n+1} =  f_{n+1} - \delta_{Hoch}(h_n).
\]
\end{lemme}

\dem 
Nous allons construire un morphisme $f'$ tel que la 
suite
\[
(0,\hdots, 0,h_n,0,\hdots)
\]
dÈfinisse une homotopie
$h$ entre $f$ et $f'$.
Nous construisons les $f'_i$ par rÈcurrence sur $i$.
Soit $i\geq 1.$ Supposons qu'il existe un ${\rm{A}}_i$-morphisme
$f' : A \ra A'$ tel que $h$ dÈfinit une homotopie entre $f$ et $f'$ en
tant que ${\rm A}_i$-morphisme. Posons
\[
\begin{array}{rcl}
f'_{i+1} &= & f_{i+1} -  \sum (-1)^s
m_{r+1+t} (f_{i_1}\ts \hdots \ts f_{i_r} \ts h_{k} \ts f'_{j_1}\ts \hdots \ts f'_{i_t})\\ 
& & - \sum (-1)^{jk+l}h_z(\Id^{\ts j}\ts m_k\ts \Id^{\ts l}),
\end{array}
\]
o˘ $s$ est le signe apparaissant dans \ref{definition_homotopie_ai-morphismes}.
Par construction, la suite des
\[
(f'_1,\hdots,f'_i,f'_{i+1})
\]
dÈfinit un ${\rm A}_{i+1}$-morphisme homotope ‡ $f$. 
Le morphisme $f'$ ainsi construit vÈrifie clairement les conditions souhaitÈes sur les $f'_i$, 
$1\leq i\leq n+1$.
\findem\\

{\noindent \em DÈmonstration du thÈorËme \ref{theoreme_unite_ai-morphisme} :}\\

Soit $A$ et $A'$ deux $\ai$-algËbres minimales strictement unitaires et
\[
f : A \ra A'
\]
un $\ai$-morphisme homologiquemement unitaire. Nous cherchons un morphisme $f'$ homotope ‡ $f$ tel que
les morphismes $f'_i$, $i\geq 1$, vÈrifient
\[
f'_i (\Id \tp j \ts \unite \ts \Id \tp l) = 0, \quad i \geq 2 \mbox{ et } j+1+l = i.
\]
Construisons les $f'_i$, $1 \leq i \leq n$,  par rÈcurrence
sur $n.$ Soit $n\geq 1$. Supposons qu'on a un morphisme
$f$, tel que les morphismes $f_i$, $2 \leq i\leq n$, vÈrifient la condition prÈcitÈe.
En utilisant les mÍmes arguments que pour le thÈorËme
\ref{theoreme_unite_ai-algebre} dans lesquels nous remplaÁons le complexe
$C(A,A)$ par le complexe $C(A,A')$ et le lemme d'obstruction
\ref{lemme_obstruction_structure_minimale} par le lemme
\ref{lemme_obstruction_structure_minimale_morphisme}, nous trouvons
qu'il existe deux ÈlÈments
\[
f'_{n+1} \in \Hom_{\sz \gr \sf C} (\b A \tp{n+1} ,A')\quad \mbox{et} \quad
h_{n} \in \Hom_{\sz \gr \sf C} (A \tp{n} ,A')
\]
tels que
\[
\delta_{Hoch} (h_n) + f'_{n+1} = f_{n} \quad \mbox{et} \quad
\delta_{Hoch}(f'_{n+1}) + r(f_2, \cdots , f_{n}) = 0.
\]
Par le lemme \ref{lemme2_unite} appliquÈ ‡ $f$ et $h_n$, il existe un morphisme
$f'$ homotope ‡ $f$ dont les morphismes $f_i$, $2 \leq i \leq n+1$, vÈrifient
les Èquations
\[
f_i (\Id \tp j \ts \unite \ts \Id \tp l) = 0, \quad i \geq 2 \mbox{ et } j+1+l = i.
\]
\findem

\subsection{Strictification unitaire des homotopies}
\label{section_stritification_unitaire_homotopies}

\begin{theoreme} \label{theoreme_unite_homotopies}
Soit $A$ et $A'$ deux $\ai$-algËbres minimales strictement unitaires.
Si $f$ et $g$ sont deux $\ai$-morphismes strictement unitaires homotopes $A \ra A'$
il existe une homotopie strictement unitaire entre $f$ et $g.$
\end{theoreme}

\begin{lemme}\label{lemme3_unite}
Soit $A$ et $A'$ deux $\ai$-algËbres minimales.
Soit $f$ et $g$ sont deux $\ai$-morphismes homotopes $A \ra A'$ et $h$ une homotopie
de $f$ vers $g$. Soit $n \geq 2$ et
\[
\rho_n : A \tp n \ra A'
\]
un morphisme graduÈ de degrÈ $-n-1.$
Il existe une homotopie $h'$ entre $f$ et $g$ telle que
\[
h'_i = h_i \quad \mbox{si} \quad 1 \leq i \leq n \quad \mbox{et} \quad h_{n+1} = h'_{n+1} + 
\delta_{Hoch}(\rho_n).
\]
\end{lemme}

\dem Nous raisonnons comme dans le lemme \ref{lemme2_unite}.  Posons $F = Bf$, $G = Bg$ et
$H = Bh : BA \ra BA'$ l'homotopie entre $F$ et $G$. Soit
$R$ la ($F$,$G$)-codÈrivation de degrÈ $-2$
qui est donnÈe (\ref{cogebres_tensorielles}) par la suite
\[
(0,\hdots ,0,s\rho_n\si \tp n,0,\hdots).
\]
Soit $H'$ dÈfini par l'ÈgalitÈ
\[
H' = H - b^{A'}R + Rb^A.
\]
C'est une ($F$,$G$)-codÈrivation qui est clairement une homotopie entre $F$ et $G$. Nous vÈrifions
qu'elle correspond ‡ une homotopie $h'$ entre $f$ et $g$ telle que 
\[
h'_i = h_i \quad \mbox{si} \quad 1 \leq i \leq n \quad \mbox{et} \quad h_{n+1} = h'_{n+1} + 
\delta_{Hoch}(\rho_n).
\]
\findem \\

{\em DÈmonstration du thÈorËme \ref{theoreme_unite_homotopies} : }\\

Nous cherchons une homotopie $h$ entre $f$ et $g$ telle que  
les morphismes $h_i$, $i\geq 1$, vÈrifient
\[
h_i (\Id \tp j \ts \unite \ts \Id \tp l) = 0, \quad i \geq 2 \mbox{ et } j+1+l = i.
\]
Construisons les $h_i$, $1 \leq i \leq n$,  par rÈcurrence
sur $n.$ Soit $n\geq 1$. Supposons qu'on a un morphisme
$h$, tel que les morphismes $h_i$, $2 \leq i\leq n$, vÈrifient la condition prÈcitÈe.
En utilisant les mÍmes arguments que pour le thÈorËme
\ref{theoreme_unite_ai-algebre} dans lesquels nous remplaÁons le complexe
$C(A,A)$ par le complexe $C(A,A')$ (voir \ref{section_obstruction_unite}) et le lemme d'obstruction
\ref{lemme_obstruction_structure_minimale} par le lemme
\ref{lemme_obstruction_structure_minimale_homotopie}, nous trouvons
qu'il existe deux ÈlÈments
\[
h'_{n+1} \in \Hom_{\sz \gr \sf C} (\b A \tp{n+1} ,A')\quad \mbox{et} \quad
\rho_{n} \in \Hom_{\sz \gr \sf C} (A \tp{n} ,A')
\]
tels que
\[
\delta_{Hoch} (\rho_n) + h_{n+1} = h'_{n} \quad \mbox{et} \quad
\delta_{Hoch}(h'_{n+1}) + r(h_2, \cdots , h_{n}) = 0.
\]
Par le lemme \ref{lemme3_unite}, il existe une homotopie $h'$ entre $f$ et $g$ telle
que nous avons les Èquations
\[
h'_i (\Id \tp j \ts \unite \ts \Id \tp l) = 0, \quad i \geq 2 \mbox{ et } j+1+l = i.
\]
\findem\\

Nous dÈduisons des thÈorËmes \ref{theoreme_unite_ai-algebre},
 \ref{theoreme_unite_ai-morphisme} et \ref{theoreme_unite_homotopies}
le corollaire suivant :

\begin{corollaire}\label{corollaire_unite_homotopies}
Soit $A$ et $A'$ des $\ai$-algËbres minimales strictement unitaires et
 $f : A \ra A'$ une Èquivalence d'homotopie strictement unitaire.
Il existe un inverse ‡ homotopie prËs $g$ de $f$ qui est strictement unitaire
et des homotopies $h$ et $h'$ strictement unitaires entre $\Id_{A'}$ et $f \circ g$,
et entre $\Id_{A}$ et $g \circ f$. \findem
\end{corollaire}

\subsection{ModËle minimal d'une $\ai$-algËbre strictement unitaire}
\label{section_modele_minimal_strictement_unitaire}

Le corollaire (\ref{corollaire_unite_ai-algebre}) montre
que toute $\ai$-algËbre homologiquement unitaire $A$ admet un modËle minimal
strictement unitaire $A'$ tel que l'$\ai$-quasi-isomorphisme
\[
f : A' \ra A
\]
vÈrifie $f \circ \unite = \unite.$
Le but de cette section est de montrer la proposition suivante :

\begin{proposition} \label{proposition_modele_minimal_strictement_unitaire} 
\index{modele minimal@{modËle minimal}}
Toute $\ai$-algËbre strictement unitaire $A$ admet un modËle minimal
strictement unitaire $A'$ tel que l'$\ai$-quasi-isomorphisme
\[
f : A' \ra A
\]
est strictement unitaire.
\end{proposition}

Notre dÈmonstration est basÈe sur le lemme de perturbation (voir \cite{Huebschmann91a},
\cite{Gugenheim86},  \cite{Gugenheim89}, \cite{Gugenheim91}, \cite{Merkulov99} et \cite{Kontsevich01}).\\

\dem
Posons $V = H^* A$.
Soit $i : (V,0) \ra (A,m_1)$ 
un morphisme de complexes qui induit l'identitÈ en homologie et tel que $i \circ \unite = \unite$.
Soit $p : A \ra K$ le conoyau de $i$. Le complexe $K$ est contractile. La suite de complexes
$(i,p)$ est donc scindable. Choisissons une rÈtraction $\rho$ et une section $\sigma$ telles que
\[
\rho \circ \sigma = 0 \quad \mbox{et} \quad i \circ \rho + \sigma \circ p = \Id_A.
\]
Soit $h$ une homotopie contractante de $K$ tel que
$h^2 = 0$.
Soit $A' = V^\delta$ l'$\ai$-algËbre (de complexe sous-jacent $V$) et
$f = f^\delta$
le morphisme d'$\ai$-algËbres construits ‡ partir de ces donnÈes dans (\ref{lemme_modele_minimal_I}).
Montrons que $A'$ est une $\ai$-algËbre strictement unitaire et que l'$\ai$-morphisme
$f$ est strictement unitaire. Nous utilisons les notations de la dÈmonstration
de (\ref{lemme_modele_minimal_I}). Nous avons clairement les ÈgalitÈs
\[
m'_1 \circ \unite = 0, \quad m'_2(\unite \ts \Id) = m'_2(\Id \ts \unite) = \Id \quad
\mbox{et} \quad f_1 \circ \unite = \unite.
\]
Il reste ‡ montrer que la composition de $f_i$, $i\geq 2$, et $m'_i$, $i\geq 3$, par
\[
\unite_\alpha = (\Id \tp \alpha \ts \unite \ts \Id \tp{i - 1 - \alpha}), \quad 0 \leq \alpha < i,
\]
est nulle. Il suffit de montrer que les compositions
\[
m_{i,T} \circ \unite_\alpha
\quad \mbox{et} \quad f_{i,T} \circ \unite_\alpha, \quad T \in \mathcal T,
\]
sont nulles.
Remarquons que ces compositions proviennent d'arbres $\b T$, coloriÈs comme pour
$m_{i,T}$ (resp.~$f_{i,T}$) sauf en une feuille qui est maintenant de couleur $\unite$.
Comme $A$ est strictement unitaire, nous avons
\[
m_j \circ \eta_\beta = 0, \quad j \geq 3,\quad  0 \leq \beta < j.
\]
Il suffit donc de vÈrifier la nullitÈ des compositions
provenant d'arbres coloriÈs dont un sous-arbre coloriÈ est de la forme
\begin{center}
\input{konts3.pstex_t} \hspace{2cm} 
 \input{konts4.pstex_t}.
\end{center}
Dans les deux premiers cas,
$m'_{i,T} \circ \unite_\alpha$ et $f_{i,T} \circ \unite_\alpha$, s'annulent
car $H^2 = 0$, dans les autres cas, car $i \circ H = 0.$
\findem 

\begin{remarque}{\em \label{remarque_proposition_modele_minimal_strictement_unitaire}
Nous vÈrifions de la mÍme maniËre que le morphisme $q^\delta$ et l'homotopie $H^\delta$ de la remarque
(\ref{remarque3_lemme_perturbation}) sont aussi strictement unitaires. Le lemme de perturbation
produit donc une contraction dans la catÈgorie des $\ai$-algËbres strictement unitaires.
}\end{remarque}

Soit $\big(\aia\big)_u$ (resp.~$\big(\aia\big)_{su}$) \indexnotation{aiau} \indexnotation{aiasu}
la catÈgorie des $\ai$-algËbres strictement unitaires dont les
espaces de morphismes sont formÈs des morphismes homologiquement unitaires (resp.~strictement unitaires).
Notons $\sim_{u}$ (resp.~$\sim_u$) la relation d'homotopie relativement aux
homotopies au sens de \ref{definition_homotopie_ai-morphismes} (resp.~aux homotopies
strictement unitaires).

\begin{proposition} \label{proposition_ai-algebres_strictement_unitaires}
L'inclusion
\[
\big(\aia\big)_{su} \hookrightarrow \big(\aia\big)_{u}
\]
induit une Èquivalence
\[
J : \big(\aia\big)_{su}/\!\sim_{su}\, \ra \big(\aia\big)_{u}/\!\sim_u\,.
\]
\end{proposition}

\dem
La remarque (\ref{remarque_proposition_modele_minimal_strictement_unitaire})
montre qu'il suffit de montrer que $J$ induit un isomorphisme dans les
espaces de morphismes dont le but et la source sont des $\ai$-algËbres minimales strictement unitaires.
Nous strictifions les $\ai$-morphismes, puis les homotopies entre
$\ai$-morphismes strictement unitaires gr‚ce aux thÈorËmes (\ref{theoreme_unite_ai-morphisme})
et (\ref{theoreme_unite_homotopies}).
\findem

\begin{corollaire}\label{corollaire_ai-algebres_strictement_unitaires}
La sous-catÈgorie $\big(\aia \big)_{hu} \subset \aia$ \indexnotation{aiahu}
des $\ai$-algËbres homologiquement unitaires dont les morphismes
sont les $\ai$-morphismes homologiquement unitaires et la catÈgorie $\big(\aia \big)_{su}$
deviennent Èquivalentes aprËs passage ‡
l'homotopie. \findem
\end{corollaire}

{\noindent \bf (Co)fibrations triviales strictement unitaires}\\

Nous finissons cette section par des rÈsultats qui nous seront utiles
dans la section (\ref{section_categorie_derivee_strictement_unitaire}).

\begin{lemme} \label{lemme_(co)fibrations_triviales_strictement_unitaires}
Soit $A$ et $A'$ des $\ai$-algËbres strictement unitaires.
\english \begin{itemize}
\item[{\it a.}] 
Soit $i : A \ra A'$ une cofibration triviale strictement unitaire.
Il existe un $\ai$-morphisme $p : A' \ra A$ strictement unitaire tel
que $p \circ i = \Id_A$.
\item[{\it b.}]
Soit $q : A' \ra A$ une fibration triviale strictement unitaire.
Il existe un $\ai$-morphisme $j : A \ra A'$ strictement unitaire tel
que $q \circ j = \Id_A$.
\end{itemize}\francais
\end{lemme}

\dem Les arguments de la dÈmonstration des deux points Ètant duaux nous
ne prouvons que le point {\it a}.
Supposons donnÈ un $\ai$-morphisme strictement unitaire $p'$ tel que
la composition $\alpha = p' \circ i$ est un automorphisme de $A$.
Comme $\alpha$ est la composÈe d'$\ai$-morphismes strictement unitaires,
il est strictement unitaire. Le lemme
(\ref{lemme_inverse_ai-isomorphisme_strictement_unitaire}) ci-dessous
montre l'$\ai$-morphisme $\alpha^{-1}$ est aussi strictement unitaire.
Posons $p = \alpha^{-1} \circ p'$ et nous avons le rÈsultat car $p \circ i = \Id_A$.

Nous devons donc trouver un $\ai$-morphisme $p'$ strictement unitaire tel que
$p' \circ i$ est un automorphisme de $A$.

{\em Premier cas : l'unitÈ $\unite$ est un bord de $A'$.}\\
Dans cette situation, l'unitÈ est nulle dans la cohomologie. Il en rÈsulte que
$A$ et $A'$ sont faiblement Èquivalentes ‡ $0.$
DÈfinissons $p'_1$ comme un scindage de $i_1$. Il vÈrifie l'ÈgalitÈ $p'_1 \circ \unite = \unite$.
Les morphismes $p'_i$, $i \geq 2$, sont dÈfinis par rÈcurrence sur $i$.
Soit $h$ une homotopie contractante de $A.$
Posons
\[
p'_i = -h \circ r(p'_1, \hdots ,p'_{i-1}), \quad i\geq 2,
\]
o˘ $r(p'_1, \hdots ,p'_{i-1})$ est le cycle du lemme (\ref{extension_nstructure_morphismes}).
Nous vÈrifions (par rÈcurrence) que $r(p'_1, \hdots ,p'_{i-1})$ composÈ avec
\[
\Id \tp \alpha \ts \unite \ts \Id \tp \beta, \quad \alpha + 1 + \beta = i+1,
\]
est nul. Les morphismes $p'_i$, $i \geq 1$, ainsi construits dÈfinissent bien un $\ai$-morphisme
gr‚ce au lemme (\ref{extension_nstructure_morphismes}). Il est strictement unitaire et,
comme nous avons l'ÈgalitÈ
\[
(p' \circ i)_1 = p'_1 \circ i_1 = \Id
\]
$p' \circ i$ est un automorphisme de $A$.

{\em DeuxiËme cas : l'unitÈ $\unite$ n'est pas un bord de $A'$.}\\
Comme $i$ est une cofibration triviale, l'axiome (CM4) de la catÈgorie $\aia$ (voir \ref{theoreme_cmf_aia})
nous donne un $\ai$-morphisme $q : A' \ra A$ tel que $q \circ i = \Id_A.$ 
L'$\ai$-morphisme $q$ est clairement homologiquement unitaire qui vÈrifie l'ÈgalitÈ
$q_1 \circ \unite = \unite$. Comme $A$ et $A'$ sont strictement
unitaires, il existe (\ref{proposition_ai-algebres_strictement_unitaires}) un
$\ai$-morphisme strictement unitaire $q' : A' \ra A$ homotope ‡ $q$.
Comme l'unitÈ $\unite$ n'est pas un bord de $A'$, il existe une rÈtraction de complexes de
de $\unite : e \ra A'$. Il induit un scindage $A' = e \oplus \b A'$.
Nous savons que le morphisme $q_1 - q'_1$ est homotope ‡ zÈro et qu'il s'annule
sur $e$. Il se factorise en $z \circ t$, o˘ $t$ est la projection $A' \ra \b A'$.
Comme cette projection est scindÈe dans la catÈgorie des complexes, $z$ est homotope ‡ zÈro.  
Il existe donc une homotopie $h_1$ entre $q_1$ et $q'_1$ telle que $h_1 \circ \unite = 0$ et
nous avons l'ÈgalitÈ $q'_1 \circ i_1 = \Id_A + \delta (h_1)\circ i_1$.

Construisons les morphismes $p'_i$, $i \geq 1$, 
‡ partir des morphismes $q'_j$, $j\geq 1$, par rÈcurrence sur $i$ :
Posons
\[
p'_1 = q'_1 - \delta(h_1)
\]
et, pour $i \geq 2$,
\[
p'_i = q'_i - \sum (-1)^s
m_{r+1+t} (p'_{i_1}\ts \hdots \ts p'_{i_r} \ts h_{1} \ts q'_{j_1}\ts \hdots \ts q'_{i_t})
+ \sum h_1 \circ m_i,
\]
o˘ $s$ est dÈfini en (\ref{definition_homotopie_ai-morphismes}).
Les morphismes $p'_i$ , $i \geq 1$, dÈfinissent ainsi un $\ai$-morphisme
strictement unitaire $A' \ra A$ tel que la suite
\[
(h_1, 0, \hdots)
\]
est une homotopie entre $q'$ et $p'.$
La composition $p' \circ i$ est un automorphisme car
\[
(p' \circ i)_1 = (q'_1 - \delta(h_1)) \circ i_1 = q'_1 \circ i_1 - \delta (h_1) \circ i_1 = \Id_A + \delta (h_1)\circ i_1 - \delta (h_1)\circ i_1 = \Id_A.
\]
\findem

\begin{lemme} \label{lemme_inverse_ai-isomorphisme_strictement_unitaire}
Soit $A$ et $A'$ deux $\ai$-algËbres strictement unitaires.
Soit $\alpha : A \ra A'$ un $\ai$-isomorphisme strictement unitaire.
L'$\ai$-morphisme $\beta = \alpha^{-1}$ est strictement unitaire.
\end{lemme}

\dem On note $\unite$ l'unitÈ des $\ai$-algËbres.
Comme $\alpha_1 \circ \unite = \unite$, nous avons l'ÈgalitÈ
$\beta_1 \circ \unite = \unite$.
Nous savons que le morphisme
\[
\alpha_2 \circ (\beta_1 \ts \beta_1) + \alpha_1 \circ \beta_2 : A' \tp 2 \ra A
\]
est nul. Si nous le composons avec $\unite \ts \Id$ (resp.~$\Id \ts \unite$),
nous trouvons que
\[
\alpha_1 \circ \beta_2 (\unite \ts \Id), \quad (\mbox{resp.}\quad \alpha_1 \circ \beta_2 (\Id \ts \unite))
\]
est nul. Comme $\alpha_1$ est un isomorphisme, ceci implique que
\[
\beta_2 (\unite \ts \Id) = 0 \quad \mbox{et} \quad \beta_2 (\Id \ts \unite) = 0.
\]
Nous continuons par rÈcurrence sur $n.$ Supposons que $\beta_1 \unite = \unite$ et
\[
\beta_i (\Id \tp j \ts \unite \ts \Id \tp k) = 0, \quad j + 1 + k = i, \quad 2 \leq i \leq n.
\]
Nous en dÈduisons l'ÈgalitÈ
\[
(\alpha \circ \beta)_{n+1} (\Id \tp j \ts \unite \ts \Id \tp k) = \alpha_1 \circ \beta_{n+1} (\Id \tp j \ts \unite \ts \Id \tp k), \quad j + 1 + k = n+1.
\]
Comme le terme dÈfinissant $(\alpha \circ \beta)_{n+1}$ est nul, nous en dÈduisons
que  
\[
\beta_{n+1} (\Id \tp j \ts \unite \ts \Id \tp k) = 0, \quad j + 1 + k = n+1.
\]
\findem

\section{Strictification unitaire des polydules} \label{section_strictification_unitaire_(bi)polydules}

Cette section traite des diffÈrents types de compatibilitÈ aux
unitÈs des $\ai$-(bi)polydules. Les dÈmonstrations sont omises car elles sont
similaires ‡ celles de la section \ref{section_strictification_unitaire_aia}.

\subsection{Polydules homologiquement unitaires}

\label{strictification_unite_modules}

\begin{definition}{\em 
Soit $A$ une $\ai$-algËbre homologiquement unitaire. Un $A$-polydule $M$ est
{\em homologiquement unitaire} 
si $H^*M$ est un $H^*A$-module unitaire.
Si $M$ et $M'$ sont deux $A$-polydules homologiquement unitaires, un $\ai$-morphisme
$f : M \ra M'$ est toujours {\em homologiquement unitaire}, i.~e.~$f_1$ induit
un morphisme de $H^*A$-modules unitaires
\[
H^*M \ra  H^*M'.
\]
}\end{definition}

Soit $A$ une $\ai$-algËbre strictement unitaire. Un $A$-polydule strictement unitaire
(\ref{definition_polydules_strictement_unitaires}) est clairement homologiquement unitaire.\\



{\noindent \bf Les rÈsultats}\\

Soit $A$ une $\ai$-algËbre minimale et strictement unitaire.

\begin{theoreme} \label{theoreme1_strictification_unite_polydules}
Tout $A$-polydule minimal homologiquement unitaire est isomorphe ‡ un $A$-polydule
strictement unitaire. \findem
\end{theoreme}

\begin{corollaire}\label{corollaire_strictification_unite_polydules}
Tout $A$-polydule homologiquement unitaire est homotopiquement Èquivalent ‡ un $A$-polydule
strictement unitaire. \findem
\end{corollaire}

\begin{theoreme}
Soit $M$ et $M'$ deux
$A$-polydules minimaux strictement unitaires.
Tout $\ai$-morphisme $f : M \ra M'$ est homotope ‡ un $\ai$-morphisme strictement
unitaire. \findem
\end{theoreme}

\begin{theoreme}
Soit $M$ et $M'$ deux $A$-polydules minimaux strictement unitaires.
Si $f$ et $g$ sont deux $\ai$-morphismes strictement unitaires homotopes $M \ra M'$
il existe une homotopie strictement unitaire entre $f$ et $g.$
\findem
\end{theoreme}

\begin{corollaire}
Soit $M$ et $M'$ des $A$-polydules strictement unitaires et
 $f : M \ra M'$ une Èquivalence d'homotopie strictement unitaire.
Il existe un inverse ‡ homotopie prËs $g$ de $f$ qui est strictement unitaire
et des homotopies $h$ et $h'$ strictement unitaires entre $\Id_{M'}$ et $f \circ g$,
et entre $\Id_{M}$ et $g \circ f$. \findem
\end{corollaire}

Soit $A$ une $\ai$-algËbre strictement unitaire.

\begin{proposition}
Tout $A$-polydule strictement unitaire $M$ admet un modËle minimal
strictement unitaire $M'$ tel que l'$\ai$-quasi-isomorphisme
\[
f : M' \ra M
\]
est strictement unitaire. \findem
\end{proposition}

Soit $\big(\aiMod A\big)_u$ \indexnotation{(aiMod)u} la {\em sous-catÈgorie} pleine de $\aiMod A$
formÈe des $A$-polydules strictement unitaires. 

\begin{proposition} \label{proposition_strictication_unitaire_polydules}
L'inclusion
\[
\aiModu A \hookrightarrow \big(\aiMod A\big)_{u}
\]
induit une Èquivalence
\[
\aiModu A/\!\sim\, \arr{\sim} \big(\aiMod A\big)_{u}/\!\sim\,,
\]
o˘ les symboles $\sim$ dÈsignent la relation d'homotopie 
(\ref{definition_polydules_strictement_unitaires}) et
(\ref{definition_homotopie_ai-morphismes_mod}).
\findem
\end{proposition}

\subsection{Bipolydules homologiquement unitaires} \label{strictification_unite_bimodules}

Soit $C$ et $C'$ deux cogËbres diffÈrentielles graduÈes et soit
$N$ et $N'$ deux $C$-$C'$-bicomodules diffÈrentiels graduÈs.
Notons $\Delta^R$ et $\Delta^L$ la comultiplication ‡ droite et ‡ gauche de
ces bicomodules.

\begin{definition}{\em 
Une {\em codÈrivation de bicomodules} \index{coderivation@{codÈrivation}}est un morphisme
\[
K : N \ra N'
\]
tel que
\[
\Delta^L \circ K = (\Id \ts K) \circ \Delta^L \quad \mbox{et} \quad
\Delta^{R} \circ K = (K \ts \Id) \circ \Delta^R.  
\]
}\end{definition}

Soit $A$ et $A'$ deux algËbres graduÈes associatives unitaires
et $M$ un $A$-$A'$-bimodules graduÈs.
ConsidÈrons les comme des $\ai$-algËbres et comme un $A$-$A'$-bipolydule.
L'espace de codÈrivations de $\Ba A$-$\Ba A'$-bicomodules
\[
\coder (BM,BM)
\]
joue dans cette section le rÙle de l'espace
\[
\coder ((BA)\+,(BA)\+)
\]
de la section \ref{section_obstruction_unite}.\\

Soit $A$ et $A'$ deux $\ai$-algËbres strictement unitaires.
Soit $\big(\aiMod (A,A')\big)_u$ \indexnotation{(aiMod)uAA'}
la {\em sous-catÈgorie} pleine de $\aiMod (A,A')$
formÈe des $A$-$A'$-bipolydules strictement unitaires. 
Nous montrons de la mÍme maniËre que prÈcÈdemment la proposition suivante :

\begin{proposition} \label{proposition_strictication_unitaire_bipolydules}
L'inclusion
\[
\aiModu (A,A') \hookrightarrow \big(\aiMod (A,A')\big)_{u}
\]
induit une Èquivalence
\[
\aiModu (A,A')/\!\sim\, \arr{\sim} \big(\aiMod (A,A')\big)_{u}/\!\sim\,,
\]
o˘ les symboles $\sim$ dÈsignent les relations d'homotopies.\findem
\end{proposition}

%% file: konts3.pstex_t
\begin{picture}(0,0)%
\includegraphics{konts3.pstex}%
\end{picture}%
\setlength{\unitlength}{3947sp}%
\begingroup\makeatletter\ifx\SetFigFont\undefined%
\gdef\SetFigFont#1#2#3#4#5{%
  \reset@font\fontsize{#1}{#2pt}%
  \fontfamily{#3}\fontseries{#4}\fontshape{#5}%
  \selectfont}%
\fi\endgroup%
\begin{picture}(1602,1781)(702,-904)
\put(1201,314){\makebox(0,0)[lb]{\smash{\SetFigFont{12}{14.4}{\rmdefault}{\mddefault}{\updefault}{\color[rgb]{0,0,0}$H$}%
}}}
\put(976,-136){\makebox(0,0)[lb]{\smash{\SetFigFont{12}{14.4}{\rmdefault}{\mddefault}{\updefault}{\color[rgb]{0,0,0}$m_2$}%
}}}
\put(976,-511){\makebox(0,0)[lb]{\smash{\SetFigFont{12}{14.4}{\rmdefault}{\mddefault}{\updefault}{\color[rgb]{0,0,0}$H$}%
}}}
\put(1801,314){\makebox(0,0)[lb]{\smash{\SetFigFont{12}{14.4}{\rmdefault}{\mddefault}{\updefault}{\color[rgb]{0,0,0}$H$}%
}}}
\put(2026,-511){\makebox(0,0)[lb]{\smash{\SetFigFont{12}{14.4}{\rmdefault}{\mddefault}{\updefault}{\color[rgb]{0,0,0}$H$}%
}}}
\put(1036,-904){\makebox(0,0)[lb]{\smash{\SetFigFont{12}{14.4}{\rmdefault}{\mddefault}{\updefault}{\color[rgb]{0,0,0}$\vdots$}%
}}}
\put(2075,-901){\makebox(0,0)[lb]{\smash{\SetFigFont{12}{14.4}{\rmdefault}{\mddefault}{\updefault}{\color[rgb]{0,0,0}$\vdots$}%
}}}
\put(1430,712){\makebox(0,0)[lb]{\smash{\SetFigFont{12}{14.4}{\rmdefault}{\mddefault}{\updefault}{\color[rgb]{0,0,0}$\vdots$}%
}}}
\put(1729,709){\makebox(0,0)[lb]{\smash{\SetFigFont{12}{14.4}{\rmdefault}{\mddefault}{\updefault}{\color[rgb]{0,0,0}$\vdots$}%
}}}
\put(702,355){\makebox(0,0)[lb]{\smash{\SetFigFont{12}{14.4}{\rmdefault}{\mddefault}{\updefault}{\color[rgb]{0,0,0}$\unite$}%
}}}
\put(2251,355){\makebox(0,0)[lb]{\smash{\SetFigFont{12}{14.4}{\rmdefault}{\mddefault}{\updefault}{\color[rgb]{0,0,0}$\unite$}%
}}}
\put(2026,-125){\makebox(0,0)[lb]{\smash{\SetFigFont{12}{14.4}{\rmdefault}{\mddefault}{\updefault}{\color[rgb]{0,0,0}$m_2$}%
}}}
\end{picture}

%% file: konts4.pstex_t
\begin{picture}(0,0)%
\includegraphics{konts4.pstex}%
\end{picture}%
\setlength{\unitlength}{3947sp}%
\begingroup\makeatletter\ifx\SetFigFont\undefined%
\gdef\SetFigFont#1#2#3#4#5{%
  \reset@font\fontsize{#1}{#2pt}%
  \fontfamily{#3}\fontseries{#4}\fontshape{#5}%
  \selectfont}%
\fi\endgroup%
\begin{picture}(1403,1540)(740,-646)
\put(928,-646){\makebox(0,0)[lb]{\smash{\SetFigFont{12}{14.4}{\rmdefault}{\mddefault}{\updefault}{\color[rgb]{0,0,0}$\vdots$}%
}}}
\put(1820,-646){\makebox(0,0)[lb]{\smash{\SetFigFont{12}{14.4}{\rmdefault}{\mddefault}{\updefault}{\color[rgb]{0,0,0}$\vdots$}%
}}}
\put(740,729){\makebox(0,0)[lb]{\smash{\SetFigFont{12}{14.4}{\rmdefault}{\mddefault}{\updefault}{\color[rgb]{0,0,0}$\unite$}%
}}}
\put(2108,718){\makebox(0,0)[lb]{\smash{\SetFigFont{12}{14.4}{\rmdefault}{\mddefault}{\updefault}{\color[rgb]{0,0,0}$\unite$}%
}}}
\put(1666,595){\makebox(0,0)[lb]{\smash{\SetFigFont{12}{14.4}{\rmdefault}{\mddefault}{\updefault}{\color[rgb]{0,0,0}$i$}%
}}}
\put(1201,588){\makebox(0,0)[lb]{\smash{\SetFigFont{12}{14.4}{\rmdefault}{\mddefault}{\updefault}{\color[rgb]{0,0,0}$i$}%
}}}
\put(893,202){\makebox(0,0)[lb]{\smash{\SetFigFont{12}{14.4}{\rmdefault}{\mddefault}{\updefault}{\color[rgb]{0,0,0}$m_2$}%
}}}
\put(867,-192){\makebox(0,0)[lb]{\smash{\SetFigFont{12}{14.4}{\rmdefault}{\mddefault}{\updefault}{\color[rgb]{0,0,0}$H$}%
}}}
\put(1760,-196){\makebox(0,0)[lb]{\smash{\SetFigFont{12}{14.4}{\rmdefault}{\mddefault}{\updefault}{\color[rgb]{0,0,0}$H$}%
}}}
\put(1775,198){\makebox(0,0)[lb]{\smash{\SetFigFont{12}{14.4}{\rmdefault}{\mddefault}{\updefault}{\color[rgb]{0,0,0}$m_2$}%
}}}
\end{picture}

%% file: Cat_der.tex
{\noindent \bf Introduction}\\

\noindent Soit $A$ une $\ai$-algËbre {\em augmentÈe}.
Dans le chapitre \ref{chapitre_Homot_aimod}, nous avons montrÈ que la catÈgorie
dÈrivÈe $\cd_\infty A$ admet trois descriptions :
\[
\big(\aiModu A\big) [Qis^{-1}], \quad
\ch_\infty A = \aiModu A/\!\sim\, \quad \mbox{et} \quad \big(\aiModust A\big) [Qis^{-1}]
\]
o˘ $\sim$ est la relation d'homotopie.
Dans ce chapitre, nous dÈfinissons la catÈgorie dÈrivÈe $\cd_\infty A$ d'une $\ai$-algËbre $A$ quelconque.
Nous montrons que les trois descriptions ci-dessus restent valables si $A$ est {\em strictement unitaire}.\\

{\bf \noindent Plan du chapitre}\\

\noindent 
Soit $B$, $B'$ deux $\corps$-algËbres
associatives et $X$ un $B$-$B'$-bimodule.
Les foncteurs standard associÈs ‡ $X$ sont les foncteurs adjoints
\[
\Hom_{B'}(X,\?)\quad \mbox{et} \quad ? \ts_B X.
\]
Maintenant, soit $A$ et $A'$ des $\ai$-algËbres et $X$ un  $A$-$A'$-bipolydule.
Dans la section \ref{section_foncteurs_standard}, nous dÈfinissons les foncteurs standard
\[
\Homi_{A'}(X,\?)\quad \mbox{et} \quad ? \tsi_A X
\]
et nous montrons qu'ils forment une paire de foncteurs adjoints.

Dans la section \ref{section_categorie_derivee_generale}, nous dÈfinissons
la catÈgorie $\cd_\infty A$ d'une $\ai$-algËbre quelconque
et nous la dÈcrivons dans le cas o˘ $A$ est $H$-unitaire (\ref{proposition_categorie_derivee_H-unitaire}).
Dans la section \ref{section_categorie_derivee_strictement_unitaire}, nous montrons
(\ref{theoreme_categorie_derivee_strictement_unitaire}) que
si $A$ est strictement unitaire, la catÈgorie $\cd_\infty A$
telle que dÈfinie dans la section prÈcÈdente est Èquivalente aux catÈgories
\[
\big(\aiModu A\big) [Qis^{-1}],\quad \ch_\infty A \quad \mbox{et} \quad \big(\aiModust A\big) [Qis^{-1}].
\]
En particulier, si $A$ est une $\ai$-algËbre augmentÈe, les dÈfinitions de la catÈgorie dÈrivÈe du
chapitre \ref{chapitre_Homot_aimod} et de celui-ci sont Èquivalentes.
Dans la section \ref{section_categorie_derivee_bipolydules}, nous
Ètudions la catÈgorie dÈrivÈe $\cd_\infty (A,A')$,
o˘ $A$ et $A'$ sont deux $\ai$-algËbres.

\section{La catÈgorie dÈrivÈe des polydules}
\label{section_categorie_derivee_polydules}

\subsection{Les foncteurs standard}  \label{section_foncteurs_standard}

{\noindent \bf Notations}\\

Soit $\sf C$ une bicatÈgorie (voir \cite[Chap.~XII, ß6]{MacLane98}). \indexnotation{C2}
Supposons que, pour tous $\mathbb O,\mathbb O' \in \Obj \sf C$,
la {\em catÈgorie}
\[
\sf C(\mathbb O,\mathbb O') = \Hom_{\sf C}(\mathbb O,\mathbb O')
\]
est une $k$-catÈgorie de Grothendieck semi-simple et que le foncteur
de composition (associatif ‡ un isomorphisme donnÈ prËs)
\[
\sf C(\mathbb O',\mathbb O'') \times \sf C(\mathbb O,\mathbb O') \ra \sf C(\mathbb O,\mathbb O''),
\quad (M,N) \mapsto M \circ N,
\]
o˘ $\mathbb O,\mathbb O', \mathbb O'' \in \Obj \sf C$,
est $k$-bilinÈaire dans les espaces de morphismes.
Nous appelons {\em produit tensoriel au-dessus de $\mathbb O'$} ce foncteur
et notons \indexnotation{tsbbO}
\[
M \ts_{\mathbb O'} N = M \circ N.
\]
Supposons en outre que, pour tout objet $X$  de $\sf C(\mathbb O',\mathbb O'')$, le foncteur
\[
? \ts_{\mathbb O'} X : \sf C(\mathbb O,\mathbb O') \ra \sf C(\mathbb O,\mathbb O'')
\]
admet un adjoint ‡ droite
\[
\Hom_{\mathbb O''}(X,?) : \sf C(\mathbb O,\mathbb O'') \ra \sf C(\mathbb O,\mathbb O').
\]
Remarquons que le produit tensoriel au-dessus de $\mathbb O$
\[
\sf C(\mathbb O,\mathbb O) \times \sf C(\mathbb O,\mathbb O) \ra \sf C(\mathbb O,\mathbb O),
\quad (M,N) \mapsto M \ts_{\mathbb O} N,
\]
o˘ $\mathbb O\in \Obj \sf C$, munit la catÈgorie $\sf C(\mathbb O,\mathbb O)$
\indexnotation{C(O,O')}
d'une structure de catÈgorie monoÔdale. Notons $e_{\mathbb O}$ \indexnotation{ebbO} l'ÈlÈment neutre
pour le produit tensoriel.
Soit $\mathbb O', \mathbb O''$ des objets de $\sf C$.
La catÈgorie $\sf C(\mathbb O,\mathbb O)$ agit ‡ droite sur la catÈgorie
$\sf C(\mathbb O',\mathbb O)$ et ‡ gauche sur la catÈgorie $\sf C(\mathbb O,\mathbb O')$
par le produit tensoriel  $\ts_{\mathbb O}$.\\

L'exemple suivant apparaÓt naturellement dans l'Ètude des $\ai$-catÈgories
(\ref{section_categories_de_base_aicat}).

\begin{exemple}{\em  
La bicatÈgorie $\sf C$ a pour objets les ensembles considÈrÈs comme des catÈgories discrËtes.
Soit $\mathbb A$ et $\mathbb B$ deux ensembles.
Nous dÈfinissons $\sf C( \mathbb A,\mathbb B)$ comme la catÈgorie des foncteurs
\[
\mathbb B\op \times \mathbb A \ra \vect \corps.
\]
La composition de $\sf C$ est donnÈe par les produits
tensoriels au-dessus des catÈgories.
Le foncteur adjoint au foncteur
\[
\lrus{\mathbb A}{?}{\mathbb B} \ts_{\mathbb B} (\lrus{\mathbb{B}}{X}{\mathbb C})
: \sf C(\mathbb A,\mathbb B) \ra \sf C(\mathbb A,\mathbb C)
\]
se rÈcrit plus naturellement
\[
\Hom_{\mathbb C}(\lrus{\mathbb{B}}{X}{\mathbb C},\lrus{\mathbb A}{?}{\mathbb C}) :
\sf C(\mathbb A,\mathbb C) \ra \sf C(\mathbb A,\mathbb B).
\]
}\end{exemple}

{\noindent \bf Plan de la section}\\

Soit $\mathbb P$, $\mathbb O$ et $\mathbb O'$ des objets de $\sf C$.
Soit $A$ et $A'$ deux $\ai$-algËbres
dans $\sf C(\mathbb O,\mathbb O)$ et $\sf C(\mathbb O',\mathbb O')$
et $X$ un $A$-$A'$-bipolydule dans $\sf C(\mathbb O,\mathbb O')$.
Nous allons construire un couple de foncteurs adjoints
\[
\big(? \tsi_{A} X , \Homi_{A'}(X,\?)\big) : \aiMod A \ra \aiMod A'.
\]
o˘ $\aiMod A$ est la catÈgorie des $A$-polydules dans $\sf C(\mathbb P,\mathbb O)$
et $\aiModu A'$ est la catÈgorie des $A'$-polydules dans $\sf C(\mathbb P,\mathbb O').$\\

{\noindent \bf Le foncteur $\Homi_{A'}(X,\?) : \aiMod A' \ra \aiMod A$}\index{foncteur!standard}
\indexnotation{Homi} \\

Soit $N'$ un $A'$-polydule. Remarquons que
$SX \ts \ct SA'$ est un objet de la catÈgorie $\sf C(\mathbb O,\mathbb O')$ et
que $SN' \ts \ct SA'$ est un objet de $\sf C(\mathbb P,\mathbb O')$.
Nous dÈfinissons l'objet graduÈ de $\sf C(\mathbb P,\mathbb O)$ sous-jacent
‡ $\Homi_{A'}(X,N')$ comme
\[
\Hom_{\coModcu \ct SA'}(SX \ts \ct SA',SN' \ts \ct SA'),
\]
o˘ $\Hom$ dÈsigne le foncteur adjoint $\Hom_{\mathbb O'}$.
Sa diffÈrentielle est le morphisme
\[
\delta : F \mapsto b^{BN'} \circ F - (-1)^{|F|} F \circ b^{BX_{A'}}
\]
o˘ $BX_{A'} = SX \ts \ct SA'$ est la construction bar de $X$ en tant que $A'\+$-polydule
et o˘ le morphisme $F$ est de degrÈ $|F|.$
C'est un module diffÈrentiel graduÈ sur
l'algËbre diffÈrentielle graduÈe 
\[
\mathsf{End}(BX_{A'}) = \Big(\Hom_{\gr \ct SA'}(SX \ts \ct SA',SX \ts \ct SA'),\delta \Big).
\]
La structure de $A$-polydule est donnÈe par la restriction du $\mathsf{End}(BX_{A'})$-module
diffÈrentiel graduÈ 
$\Homi_{A'}(X,N')$ le long de l'$\ai$-morphisme
\[
A  \ra \mathsf{End}(BX_{A'})
\]
dÈfini dans le lemme clef (\ref{lemme_clef}). Explicitons cette structure.
Le morphisme
\[
m^{\mathsf H}_i : \Homi_{A'}(X,N') \ts A \tp{i-1} \ra \Homi_{A'}(X,N'), \quad i \geq 1,
\]
est donnÈ par la diffÈrentielle de l'espace si $i = 1$ et, sinon, par le morphisme
\[
\xymatrix{
S\Hom_{\coModcu \ct SA'}(SX \ts \ct SA',SN' \ts \ct SA') \ts (SA)\tp{i-1} \ar[d]^{b^{\mathsf H}_i} \\
 S\Hom_{\coModcu \ct SA'}(SX \ts \ct SA',SN' \ts \ct SA')
}\]
qui envoie un ÈlÈment $s\Gamma  \ts \phi \in S\Hom_{\coModcu \ct SA'}(SX \ts \ct SA',SN' \ts \ct SA') \ts (SA)\tp{i-1}$
sur
\[
b^{\mathsf{Comc}}_2 (s \Gamma \ts s \Phi) \in S\Hom_{\coModcu \ct SA'}(SX \ts \ct SA',SN' \ts \ct SA'),
\]
o˘ le morphisme $b^{\mathsf{Comc}}_2$ correspond ‡ la composition de la catÈgorie $\coModcu \ct SA'$
et $\Phi$ est dÈfini dans le lemme clef (\ref{lemme_clef}).
Un morphisme $f : N' \ra N''$ dans $\aiMod A'$
induit un morphisme de $\mathsf{End}(BX_{A'})$-modules
diffÈrentiels graduÈs
\[
F_* : \Hom (SX \ts \ct SA',SN' \ts \ct SA') \ra \Hom(SX \ts \ct SA',SN'' \ts \ct SA'),
\]
o˘ $F_*$ est induit par la construction bar $F$ de $f$. Ainsi, le morphisme $F_*$ est strict en tant que
morphisme de $A$-polydules.
Ceci nous fournit un foncteur
\[
\Homi_{A'}(X,\?) : \aiMod A' \ra \aiModst A \hookrightarrow \aiMod A.
\]

\begin{remarque}{\em
Si $A$ sont strictement unitaire et si $X$ est un $A$-$A'$-bipolydule
strictement unitaire pour $A$, i.~e.~si la composition
\[
m_{i,j}(\Id \tpo \alpha \tso \unite \tso \Id \tpo \beta \ts \Id_M \ts \Id \tp j), \quad i,j \geq 0,
\]
est nulle si $(i,j) \neq (1,0),$ Ègale ‡ $\Id$ sinon,  le
$A$-polydule $\Homi_{A'}(X, N)$ est strictement unitaire.
Nous obtenons alors un foncteur
\[
\Homi_{A'}(X,\?) : \aiMod A' \ra \aiModust A \hookrightarrow \big(\aiMod A\big)_u,
\]
o˘  $\big(\aiMod A\big)_{u}$ est la sous-catÈgorie pleine
de $\aiMod A$ formÈe des objets strictement unitaires.
}\end{remarque}

{\noindent \bf Le foncteur $? \tsi_{A} X : \aiMod A \ra \aiMod A'$} \index{foncteur!standard}
\indexnotation{tsi}\\

Soit $N$ un $A$-polydule.
L'objet graduÈ de  $\sf C(\mathbb P,\mathbb O')$ sous-jacent ‡ $N \tsi_{A} X$ est
\[
N \ts \ct SA \ts X.
\]
La structure de $A'$-polydule sur $N \ts \ct SA \ts X$ est donnÈe par une diffÈrentielle $b$
sur
$S(N \ts \ct SA \ts X) \ts \ct SA'$. La suspension de ce $\ct SA'$-comodule diffÈrentiel
graduÈ s'identifie au {\em produit cotensoriel} \indexnotation{square}
\[
(SN \ts \ct SA) \tbox_{\ct SA} (\ct SA \ts SX \ts \ct SA'),
\]
c'est-‡-dire au noyau
\[
\ker \big(BN \ts BX \arr{\Delta \ts \Id - \Id \ts \Delta^{L}} BN \ts \ct SA \ts BX\big),
\]
o˘ $BX = \ct SA \ts SX \ts \ct SA'$ est la construction bar de $X$ en tant que $A\+$-$A'\+$-bipolydule
strictement unitaire.
Un morphisme de $A$-polydules $f :  N \ra N'$ induit un morphisme strict
\[
\Big(\si \circ F \circ s\Big) \ts \Id_X : N \ts \ct SA \ts X \ra N' \ts \ct SA \ts X.
\]
Nous obtenons ainsi un foncteur
\[
? \tsi_{A} X : \aiMod A \ra \aiModst A' \hookrightarrow \aiMod A'.
\]

\begin{remarque}{\em \label{remarque_tsi_A_strictement_unitaire}
Si $A'$ sont strictement unitaire et si $X$ est un $A$-$A'$-bipolydule
strictement unitaire pour $A'$, i.~e.~si la composition
\[
m_{i,j}(\Id \tp i \ts \Id_M \ts \Id \tpo \alpha \tso \unite \tso \Id \tpo \beta ), \quad i,j \geq 0,
\]
est nulle si $(i,j) \neq (0,1),$ Ègale ‡ $\Id$ sinon, le
$A'$-polydule $N \tsi_{A} X$ est strictement unitaire.
Nous obtenons alors un foncteur
\[
? \tsi_{A} X : \aiMod A \ra \aiModust A' \hookrightarrow \big(\aiMod A'\big)_u.
\]
}\end{remarque}

\begin{lemme} \label{lemme_adjonction_foncteurs_standard}
Le foncteur $? \tsi_{A} X$ est adjoint ‡ gauche au foncteur \hbox{$\Homi_{A'}(X, ?)$}
\end{lemme}

\dem Soit $L$ un objet de $\aiMod A$ et $R$ un objet de $\aiMod A'.$
Soit des morphismes graduÈs dans $\sf C(\mathbb O',\mathbb O')$ de degrÈ $2-i$
\[
f_{j} : L \ts \ct SA  \ts X \ts A' \tp j \ra R, \quad j \geq 0.
\]
Soit $F_j$, $j \geq 0$, les morphismes donnÈs par les bijections $F_j \lra f_j.$
Ils sont donnÈs par des morphismes de degrÈ $0$
\[
F_{i,j} : SL \ts (SA) \tp i \ts X \ts (SA') \tp j \ra SR, \quad i,j \geq 0.
\]
Soit $i \geq 0.$ Soit $g_i$ le morphisme graduÈ de $\sf C(\mathbb P,\mathbb O)$ de degrÈ $1-i$
\[
g_i : L \ts A \tp i \ra \Hom_{\ct SA'}(SX \ts \ct SA', SR \ts \ct SA'))
\]
dÈfini par l'Èquation
\[
G_i (\lambda \ts \phi) = s (\Gamma) \in S\Hom_{\ct SA'}(SX \ts \ct SA', SR \ts \ct SA'))
\]
o˘ $\lambda \ts \phi$ est un ÈlÈment de  $SL \ts (SA) \tp i$ de degrÈ $r = |\lambda \ts \phi|$,
o˘ $G_i$ est donnÈ par les bijections
$g_i \lra G_i$ et o˘ le morphisme $\Gamma$ est l'unique morphisme (voir
\ref{lemme_comodules_colibres}) tel que la composition $p_1 \circ \Gamma$ a pour
composantes les morphismes
\[
\xymatrix{
SX \ts (SA') \tp j \ar[rr]^(.4){(-1)^{|r|} \lambda \ts \phi \ts \Id}& &
SN \ts (SA) \tp i \ts SX \ts (SA') \tp j \ar[r]^(.8){F'_{i,j}}&  SR ;
}
\]
ici le morphisme $F'_{i,j}$ est le morphisme $F_{i,j} \si.$
Nous devons montrer l'Èquivalence entre les deux points suivants.
\english \begin{itemize}
\item[{\it a.}]
Les morphismes $g_{j}$ dÈfinissent un $\ai$-morphisme d'$A$-polydules
\[
L \ra \Homi_{A'}(X, R).
\]
\item[{\it b.}]
Les morphismes $f_{j}$ dÈfinissent un $\ai$-morphisme d'$A'$-polydules
\[
L \tsi_{A} X \ra R.
\]
\end{itemize} \francais

Supposons que l'ÈnoncÈ {\it a} est vrai : on a les ÈgalitÈs
\[
\sum_{k+l+m=n} G_{k+1+m} (\Id \tp k \ts b_l \ts \Id \tp m) = 
\sum_{k+m=n} b^{\mathsf H}_{1+m} (G_k \ts \Id \tp m), \quad n \geq 1,
\]
o˘ les symboles $b_l$ doivent Ítre interprÈtÈs convenablement.
Nous allons montrer que cela est Èquivalent
aux Èquations dans les espaces de morphismes
\[
\Hom_{\sf C(\mathbb P,\mathbb O')}\Big(S(L \tsi_{A} X) \ts (SA') \tp{n-1}, SR\Big) ,\quad n \geq 0,
\]
\[
\sum_{k+l+m=t} F_{k+1+m} (\Id \tp k \ts b_l \ts \Id \tp m) = 
\sum_{k+m=t} b_{1+m} (F_k \ts \Id \tp m), \quad t  \geq 1.
\]
Soit $\lambda \ts \phi \in SL \ts (SA) \tp{n-1}$ et $\kappa \ts \phi' \in SX \ts (SA')\tp{t-1}$.
Calculons
\[
G_{k+1+m} (\Id \tp k \ts b_l \ts \Id \tp m)(\lambda \ts \phi)(\kappa \ts \phi').
\]
Dans le cas o˘ $k = 0$, on a
\[
\begin{array}{rl}
  & G_{1+m} (b^R_l \ts \Id \tp m)(\lambda \ts \phi)(\kappa \ts \phi')\\
= & G_{1+m} (b^R_l (\lambda \ts \phi_{l-1}) \ts \phi_m)(\kappa \ts \phi'))\\
= &  s \Gamma (\kappa \ts \phi')\\
= & (-1)^{|\lambda \ts \phi_{l-1}|+1+|\phi_m|} F'_{1+m,t-1} 
     (b^R_l (\lambda \ts \phi_{l-1}) \ts \phi_m \ts \kappa \ts \phi'))\\
= & (-1)^{|\lambda \ts \phi|+1} F'_{1+m,t-1} (b^R_l \ts \Id\tp m \ts \Id \ts \Id \tp {t-1})
     (\lambda \ts \phi_{l-1} \ts \phi_m \ts \kappa \ts \phi')\\
= & (-1)^{|\lambda \ts \phi|+1} F'_{1+m,t-1} (b^R_l \ts \Id\tp m \ts \Id \ts \Id\tp {t-1})
(\lambda \ts \phi \ts \kappa \ts \phi'),
\end{array}
\]
o˘ $\phi_1 \ts \phi_2 = \phi$ et dans le cas o˘ $k \neq 0$, on a
\[
\begin{array}{rl}
  & G_{k+1+m} (\Id \tp k \ts b^A_l \ts \Id \tp m)(\lambda \ts \phi)(\kappa \ts \phi')\\
= & (-1)^{|\lambda| + |\phi_{k-1}|} G_{k+1+m}  (\lambda \ts \phi_{k-1} \ts b^A_l(\phi_l)
   \ts \phi_m)(\kappa \ts \phi'))\\
= & (-1)^{|\lambda| + |\phi_{k-1}|} s \Gamma(\kappa \ts \phi')\\
= & (-1)^{|\lambda| +|\phi_{k-1}|+|\lambda \ts \phi| +1} F'_{k+1+m,t-1}(\lambda \ts \phi_{k-1}
     \ts b^A_l(\phi_{l}) \ts \phi_{m} \ts \kappa \ts \phi')\\
= & (-1)^{|\lambda| +|\phi_{k-1}|+|\lambda \ts \phi|+1} F'_{k+1+m,t-1}(\Id\tp k \ts b^A_l
     \ts \Id\tp m \ts \Id \ts \Id \tp{t-1})\\
  &  \hspace*{9cm} (\lambda \ts \phi \ts \kappa \ts \phi')\\
= & (-1)^{|\lambda \ts \phi|+1}F'_{k+1+m,t-1}(\Id\tp k \ts b^A_l \ts \Id\tp m \ts \Id \ts \Id \tp{t-1})
(\lambda \ts \phi \ts \kappa \ts \phi'),\\
\end{array}
\]
o˘ $\phi_1 \ts \phi_2 \ts \phi_3 =  \phi.$
Le terme
\[
b^{\mathsf H}_{1} (G_n)(\lambda \ts \phi)(\kappa \ts \phi')
\]
vaut
\[
\begin{array}{rl}
  & b^{\mathsf H}_{1} (G_n (\lambda \ts \phi))(\kappa \ts \phi')\\
= & b^{\mathsf H}_{1} (s\Gamma) (\kappa \ts \phi')\\
= & - s m^{\mathsf H}_{1}(\Gamma) (\kappa \ts \phi')\\
= & - s [b \circ \Gamma - (-1)^{|\lambda + \phi|+1}\Gamma\circ b] (\kappa \ts \phi')\\
= & (b^L \circ s\Gamma ) (\kappa \ts \phi')  - (-1)^{|\lambda + \phi|+1} 
    (s\Gamma\circ b^{X_{A'}}) (\kappa \ts \phi')\\
= & (-1)^{|\lambda + \phi|+1} \sum  b^L_{\beta+1} (F'_{n-1,\alpha} 
    (\lambda \ts \phi \ts \kappa \ts \phi'_{\alpha})
        \ts \phi'_\beta)\\
  & + s\Gamma \sum (\id \tp{\gamma_1} \ts b_{\gamma_2} \ts \Id \tp {\gamma_3})(\kappa \ts \phi'_{\gamma_1} \ts
     \phi'_{\gamma_2} \ts \phi'_{\gamma_3})\\
= & (-1)^{|\lambda + \phi|+1} \sum  b^L_{\beta+1} (F'_{n-1,\alpha} \ts \Id\tp {\beta})
         (\lambda \ts \phi \ts \kappa \ts \phi')\\ 
  & \sum F'_{n-1,\gamma_1 + \gamma 3}(\lambda \ts \phi \ts \Id \tp{\gamma_1} 
       \ts b_{\gamma_2} \ts \Id \tp {\gamma_3})(\kappa \ts \phi')\\
= & (-1)^{|\lambda + \phi|+1} \sum  b^L_{\beta+1} (F'_{n-1,\alpha} \ts \Id\tp {\beta})
         (\lambda \ts \phi \ts \kappa \ts \phi')\\ 
  & + (-1)^{|\lambda + \phi|} \sum F'_{n-1,\gamma_1 + \gamma 3}
 (\Id \ts \Id\tp{n-1} \ts \Id \tp{\gamma_1} \ts b_{\gamma_2} \ts \Id \tp{\gamma_3})
   (\lambda \ts \phi \ts \kappa \ts \phi'),
\end{array}
\]
o˘ $\phi'_\alpha \ts \phi'_\beta = \phi$ et
les indices de la premiËre somme sont tels que $\alpha + \beta = t-1$, o˘ les indices de la
seconde somme sont tels que $\gamma_1 + \gamma_2 + \gamma_3 = t$ et les symboles $b_{\gamma_2}$ doivent Ítre
interprÈtÈs selon leur place par $b^X_{0,\gamma_2-1}$ ou par $b^{A'}_{\gamma_2}$.
Le terme
\[
b^{\mathsf H}_{1+m} (G_k \ts \Id \tp m)(\lambda \ts \phi)(\kappa \ts \phi')
\]
vaut
\[
\begin{array}{rl}
  & b^{\mathsf H}_{1+m} (G_k(\lambda \ts \phi_{k-1}) \ts \phi_m)(\kappa \ts \phi')\\
= &  b^{\mathsf H}_{1+m} (s\Gamma_{k-1} \ts \phi_m)(\kappa \ts \phi')\\
= & b^{\mathsf{Comc}}_2 (s\Gamma_{k-1} \ts s\Phi_m)(\kappa \ts \phi')\\
= & (-1)^{|\lambda + \phi_{k-1}| + 1}b^{\mathsf{Comc}}_2 (s \ts s)(\Gamma_{k-1} \ts \Phi_l)(\kappa \ts \phi')\\
= & (-1)^{|\lambda + \phi_{k-1}| + 1}sm^{\mathsf{Comc}}_2 (\Gamma_{k-1} \ts \Phi_l)(\kappa \ts \phi')\\
= & (-1)^{|\lambda + \phi_{k-1}| + 1}(s\Gamma_{k-1} \circ \Phi_l)(\kappa \ts \phi')\\
= & (-1)^{|\lambda + \phi_{k-1}| + 1 + |\lambda + \phi_{k-1}| + |\phi_m|}
\sum F'_{k,\beta}( \lambda \ts \phi_{k-1} \ts  b^X_{m,\alpha}(\phi_m \ts \kappa \ts \phi'_{m'}) \ts \phi'_{\beta})\\
= & (-1)^{1+ |\phi_m|} \sum F'_{k,\beta}( \lambda \ts \phi_{k-1} \ts  b^X_{m,\alpha}(\phi_m \ts \kappa \ts \phi'_{\alpha}) \ts \phi'_{\beta})\\
= & (-1)^{|\lambda| +|\phi| + 1} \sum F'_{k,\beta}( \Id  \ts \Id\tp{k-1} \ts  b^X_{m,\alpha}(\Id\tp m \ts \Id \ts \Id \tp \alpha)
 \ts \Id\tp \beta)\\
  & \hspace*{9cm}( \lambda \ts \phi \ts \kappa \ts \phi'),\\
\end{array}
\]
o˘ les indices de la somme sont tels que $\alpha + \beta = t-1.$
L'ÈgalitÈ $F'_{i,j} = F_{i,j} \si$ nous donne l'Èquivalence entre les points {\it a} et {\it b.}\findem

\begin{remarque}{\em
Si $A$ et $A'$ sont strictement unitaires et si $X$ est un $A$-$A'$-bipolydule strictement
unitaire, l'adjonction
\[
(? \tsi_{A} X,\Homi_{A'}(X,?) : \aiMod A \ra \aiMod A'
\]
ne se restreint pas aux sous-catÈgories $\aiModu A$ et $\aiModu A'$.
Cependant, la proposition (\ref{proposition_strictication_unitaire_polydules})
montre que les foncteurs restreints
\[
? \tsi_{A} X :  \aiModu A \ra \aiModu A' \quad \mbox{et} \quad
\Homi_{A'}(X,?) : \aiModu A' \ra \aiModu A 
\]
induisent des foncteurs adjoints dans les catÈgories dÈrivÈes $\cd_\infty A$ et $\cd_\infty A'$
(dÈfinies en \ref{section_categorie_derivee_strictement_unitaire}).
}\end{remarque}

Soit $A$ une $\ai$-algËbre strictement unitaire.
ConsidÈrons $A$ comme un $A$-$A$-bipolydule strictement unitaire.
Notons aussi
\[
? \tsi_{A} X \quad \mbox{et}
 \quad
\Homi_{A}(X,?)
\]
les foncteurs standard restreints ‡ la sous-catÈgorie $\big(\aiMod A\big)_{u}$ (voir la dÈfinition
en \ref{proposition_strictication_unitaire_polydules}).

\begin{lemme} \label{lemme2_adjonction_foncteurs_standard}
ConsidÈrons la catÈgorie des endofoncteurs de la catÈgorie
$(\aiMod A)_{u}$.\vspace{-.2cm}
\english \begin{itemize}
\item[{\it a.}]
Il existe un morphisme canonique de foncteurs
$? \tsi_{A}A \ra \Id$
qui est un quasi-isomorphisme.
\item[{\it b.}]
Il existe un morphisme canonique de foncteurs $\Id  \ra \Homi_{A}(A,?)$
qui est un quasi-isomorphisme.
\end{itemize}\francais
\end{lemme}

\dem Par l'adjonction
\[
(? \tsi_{A}A , \Homi_{A}(A,?)) : \big(\aiMod A\big)_u \ra \big(\aiMod A\big)_u,
\]
il suffit de montrer le point {\it a}.
Soit $M$ un $A$-polydule strictement unitaire. 
Nous avons un $\ai$-morphisme de $A$-polydules
\[
g : M \tsi_{A} A \ra M
\]
dont les $g_j,$ $j\geq 1$, sont dÈfinis par les morphismes
\[
m^M_{i+2+j-1}(\Id \ts \si\tp{i} \ts \Id \tp{j}) : M \ts (SA)\tp i \ts A \ts A \tp {j-1} \ra M,
\quad i\geq 0,j\geq 1.
\]
Montrons que le cÙne du morphisme
\[
g_1 : M \tsi_{A} A \ra M
\]
est acyclique. Nous vÈrifions que le morphisme de degrÈ $-1$
\[
r : M \ts \ct SA \ts SA \ra M \ts \ct SA \ts SA
\]
donnÈ par le morphisme
\[
\Id \ts s\unite : M \ts (SA)\tp{i} \ts SA \ra M \ts (SA)\tp{i+1} \ts SA, \quad i\geq 0,
\]
o˘ $\unite$ est l'unitÈ de $A$, est une homotopie contractante du cÙne de $g_1$.
\findem\\

\begin{remarque}{\em \label{remarque_lemme2_adjonction_foncteurs_standard}
Le morphisme de $A$-polydules $g$
est clairement strictement unitaire.
Le morphisme de $A$-polydules
\[
f : M \ra \Homi_{A}(A,M)
\]
correspondant par l'adjonction au morphisme $g$
est dÈfini de maniËre analogue au morphisme
\[
f : A \ra \mathsf{End} = \Homi_{A}(A,A)
\]
du lemme clef (\ref{lemme_clef}) du chapitre \ref{chapitre_ai_cat}.
Il est aussi strictement unitaire.
}\end{remarque}

\subsection{La catÈgorie dÈrivÈe d'une $\ai$-algËbre} \label{section_categorie_derivee_generale}

Soit $\mathbb O$ et $\mathbb P$ deux objets de $\sf C.$
Soit $A$ une $\ai$-algËbre dans $\sf C(\mathbb O,\mathbb O)$.
On rappelle que la catÈgorie $\aiMod A$ des $A$-polydules dans $\sf C(\mathbb P,\mathbb O)$
est isomorphe ‡ la catÈgorie $\aiModu A\+$ des $A\+$-polydules
strictement unitaires o˘ $A\+$ est l'augmentation de $A$. 
ConsidÈrons l'objet $e = e_{\mathbb O}$ comme une $\ai$-algËbre augmentÈe dans $\sf C(\mathbb O,\mathbb O)$.
ConsidÈrons l'objet $e$ comme un $A\+$-$e$-bipolydule strictement unitaire gr‚ce
‡ l'augmentation $A\+ \ra e.$ Par la section \ref{section_foncteurs_standard}, nous avons un
foncteur
\[
? \tsi_{A\+} e : \aiModu A\+ \ra \aiModu e.
\]
Il induit un foncteur dans les catÈgories dÈrivÈes que nous notons de la mÍme maniËre.

\begin{definition}{\em \label{definition_categorie_derivee_generale}
La {\em catÈgorie dÈrivÈe} d'une $\ai$-algËbre est le noyau du foncteur
\index{categorie@{catÈgorie}!derivee@{dÈrivÈe}} \indexnotation{cdinftyA}
\[
? \tsi_{A\+} e : \cd_\infty A\+ \ra \cd_\infty e.
\]
}\end{definition}

\begin{remarque}{\em
Nous montrerons en (\ref{remarque1_polydules_H-unitaires})
qu'un $A\+$-polydule strictement unitaire est dans le noyau si et seulement si sa construction
bar est acyclique.
}\end{remarque}

\begin{remarque}{\em
Le thÈorËme (\ref{theoreme_categorie_derivee_strictement_unitaire}) ci-dessous
montrera que cette dÈfinition Ètend la dÈfinition
de la catÈgorie dÈrivÈe d'une $\ai$-algËbre augmentÈe (voir \ref{definition_categorie_derivee_augmentee}).
En particulier, nous montrerons que
si $A$ est elle-mÍme augmentÈe, nous avons une suite exacte de catÈgories triangulÈes
\[
\cd_\infty A \ra \cd_\infty A\+ \ra \cd_\infty e.
\]
}\end{remarque}

\begin{theoreme} \label{theoreme_restriction_quasi-isomorphisme}
Soit $A$ et $A'$ deux $\ai$-algËbres et $f : A \ra A'$ un $\ai$-quasi-isomorphisme.
La restriction le long de $f$ induit une Èquivalence de catÈgories
\[
\cd_\infty A' \ra \cd_\infty A.
\]
\end{theoreme}

\dem Soit $f\+ : A\+ \ra A'\+$ le morphisme augmentÈ associÈ ‡~$f.$ C'est
un $\ai$-quasi-isomorphisme. Les foncteurs
\[
(\Res^{f\+}? )\tsi_{A\+} e  \quad \mbox{et} \quad ? \tsi_{A'\+} e : \aiModu A'\+ \ra \aiModu e
\]
sont donc quasi-isomorphes. Il suffit donc de montrer que la restriction le long de $f\+$
induit une Èquivalence
\[
\cd_\infty A'\+ \ra \cd_\infty A\+.
\]
Le lemme (\ref{lemme_algebre_enveloppante}) implique que le morphisme entre les algËbres enveloppantes
\[
\E (f\+) : \E (A\+) \ra \E (A'\+)
\]
est un quasi-isomorphisme.
Il s'ensuit \cite[6.1]{Keller94} que la restriction de long de $\E (f\+)$ est une Èquivalence de catÈgories
\[
\cd \E(A'\+) \ra \cd \E (A\+).
\]
Nous dÈduisons le rÈsultat du lemme (\ref{lemme2_categorie_derivee}).
\findem \\

{\noindent \bf Cas des $\ai$-algËbres $H$-unitaires}

\begin{definition}{\em
Une $\ai$-algËbre {\em $H$-unitaire} \index{H-unitaire@{$H$-unitaire}} est une
$\ai$-algËbre dont la construction bar non augmentÈe est quasi-isomorphe ‡ $0.$
}\end{definition}

La notion d'algËbre $H$-unitaire est due ‡ M.~Wodzicki \cite{Wodzicki88}.
Il montre qu'une algËbre est $H$-unitaire si et seulement si elle satisfait
la propriÈtÈ d'excision (voir \cite{Wodzicki88}, \cite{Wodzicki89}).

\begin{lemme}
Une $\ai$-algËbre minimale (i.~e.~$m_1 = 0$) strictement unitaire est $H$-unitaire.
\end{lemme}

\dem Soit $(A,\unite)$ une $\ai$-algËbre minimale strictement unitaire.
Le morphisme de degrÈ~$-1$
\[
h : BA \ra BA
\]
donnÈ par les morphismes
\[
\Id \ts (s\unite) : (SA)\tp i \ra (SA)\tp i \ts SA
\]
dÈfinit une homotopie contractante de $BA$. \findem

\begin{corollaire} \label{lemme1_H-unitaire}
Une $\ai$-algËbre homologiquement unitaire 
(voir la dÈfinition dans la section \ref{section_definition_unites})
est $H$-unitaire.
\end{corollaire}

\dem
Soit $A$ une $\ai$-algËbre homologiquement unitaire.
Le corollaire (\ref{corollaire_unite_ai-algebre}) montre que
$A$ admet un modËle minimal strictement unitaire $A'$.
Comme $BA'$ est faiblement Èquivalent ‡ $BA$ et comme
les Èquivalences faibles sont des quasi-isomorphismes, nous avons
le rÈsultat. \findem\\

{\noindent \bf La sous-catÈgorie $\Tria A$}\\

Soit $x : \mathbb P \ra \mathbb O$
un morphisme de $\sf C$.
Le morphisme $x$ induit un foncteur
\[
x^* : \sf C(\mathbb P,\mathbb O) \ra \sf C(\mathbb P,\mathbb P), \quad M \mapsto M(x).
\]
On suppose que ce foncteur admet un adjoint ‡ gauche
\[
x_! : \sf C(\mathbb P,\mathbb P) \ra \sf C(\mathbb P,\mathbb O).
\]

\begin{exemple}{\em  
Regardons l'exemple apparaissant dans l'Ètude des $\ai$-catÈgories
(\ref{section_categories_de_base_aicat}). Soit $\mathbb P$  et $\mathbb O$ deux
ensembles et soit
\[
x : \mathbb P \ra \mathbb O, \quad p \mapsto x(p)
\]
une application.
Le foncteur $x^*$ envoie $M \in C(\mathbb P,\mathbb O)$ sur
\[
(p,p') \mapsto M(x(p),p').
\]
Le foncteur $x_!$ envoie un objet $V$ de $\sf C(\mathbb P,\mathbb P)$ sur le
$\mathbb P$-$\mathbb O$-bimodule
\[
(o,p) \mapsto V(?,p) \ts_{\mathbb P} e_{\mathbb O}(o,x(?)).
\]
Supposons maintenant que $\mathbb P$ est un ensemble ‡ un ÈlÈment.
L'application $x$ est dÈterminÈ par l'ÈlÈment de $o = x(p)$ de $\mathbb O.$
Soit $V = e_\mathbb P$. L'adjonction nous donne alors un isomorphisme
\[
\Hom_{\sf C(\mathbb P,\mathbb O)}(e_{\mathbb O}(?,o),M) \arr{\sim} M(o).
\]
\findem
}\end{exemple}

Soit $x : \mathbb P \ra \mathbb O$
un morphisme de $\sf C$. Soit $V$ un objet de $\sf C(\mathbb P,\mathbb P)$.
Munissons l'objet $x_!(V) \ts_{\mathbb O} A$ de la structure de $A$-polydule
donnÈe par les multiplications de $A$. Comme nous avons un isomorphisme
\[
\Hom_{\sf C(\mathbb P,\mathbb P)}(x_!(V),M) \arr{\sim} \Hom_{\aiMod A}(x_!(V) \ts_{\mathbb O} A ,M),
\quad M \in \aiMod A,
\]
nous avons une adjonction
\[
(x_!(?) \ts_{\mathbb O} A,x^*) : \sf C(\mathbb P,\mathbb P) \ra \aiMod A.
\]

Notons $\Tria A$ la plus petite sous-catÈgorie triangulÈe aux sommes
infinies de $\cd_\infty A\+$ contenant les
\[
x \pw = x_!(e_\mathbb P) \ts_{\mathbb O} A , \quad x \in {\sf C}(\mathbb P,\mathbb O).
\]

\begin{remarque}{\em 
Cette notation est justifiÈe par le fait suivant. 
Dans l'exemple apparaissant dans l'Ètude des $\ai$-catÈgories
(\ref{section_categories_de_base_aicat}),
si $\mathbb P$ est un ensemble ‡ un ÈlÈment,
et $x$ est l'application donnÈ par un ÈlÈment $x$ de $\mathbb O,$
le $A$-polydule $x \pw$ est l'$\ai$-foncteur reprÈsentÈ par $x$
\[
x \pw = A(?,x).
\]
}\end{remarque}

\begin{proposition} \label{proposition_categorie_derivee_H-unitaire}
Soit $A$ une $\ai$-algËbre $H$-unitaire.
Nous avons une suite exacte de catÈgories triangulÈes
\[
\Tria A \hookrightarrow \cd_\infty A\+ \ra \cd_\infty e.
\]
En particulier, la catÈgorie dÈrivÈe $\cd_\infty A$ est Ègale ‡ $\Tria A.$
\end{proposition}

Dans le cas des algËbres diffÈrentielles
graduÈes cette proposition est dÈmontrÈe dans \cite{Keller94a}.
Dans la dÈmonstration ci-dessous, nous utilisons une filtration qui
est adaptÈe de celle de J.~A.~Guccione et
J.~J.~Guccione \cite{Guccione96}. Elle leur permet de montrer astucieusement la propriÈtÈ d'excision
des algËbres diffÈrentielles graduÈes $H$-unitaires.\\

\dem Montrons que la composition
\[
\Tria A \hookrightarrow \cd_\infty A\+ \ra \cd_\infty e
\]
est nulle. Comme $x \pw$ est le $A$-polydule $x_!(e) \ts_{\mathbb O} A$,
il suffit de montrer que $A \tsi_{A\+} e$ est
quasi-isomorphe ‡ $0$ dans la catÈgorie $\sf C(\mathbb O,\mathbb O)$.
Nous dÈfinissons une filtration de $A \tsi_{A\+} e = A \ts \ct (SA\+) \ts e$ par
\[
F_p  = \Big[ \bigoplus_{0 \leq i < p} A \ts (SA\+) \tp i \Big] \oplus
\Big[ \bigoplus_{0 \leq r} A \ts (SA) \tp{r} \ts (SA\+)\tp p\Big], \quad p \geq 0.
\]
Les $F_p$, $p\geq 0$, sont des sous-complexes de $A \tsi_{A\+} e$. Les objets graduÈs
\[
\ogr_p{A \tsi_{A\+} e = A \ts \ct (SA\+) \ts e} = \bigoplus_{0 \leq r}A \ts (SA)\tp r \ts (Se) \tp p, \quad p \geq 0,
\]
sont isomorphes en tant que complexes ‡
\[
S^{-1}BA \ts (Se) \tp p, \quad p \geq 0.
\]
Ils sont donc acycliques, ce qui montre que $A \tsi_{A\+} e$ est acyclique.

Pour dÈmontrer qu'on a une suite exacte de catÈgories triangulÈes, nous allons montrer
que l'inclusion de $\Tria A$ dans $\cd_\infty A\+$ a pour
adjoint ‡ droite le foncteur
\[
? \tsi_{A\+} A : \cd_\infty A\+ \ra \Tria A. 
\]
Ceci revient ‡ montrer que pour chaque $X \in \aiModu A\+$,
le triangle
\[
X \tsi_{A\+} A \ra X \ra X \tsi_{A\+} e \ra S(X \tsi_{A\+} A)
\]
est tel que l'objet $X \tsi_{A\+} e \in \Tria e$ est $(\Tria A)$-local, i.~e.~
\[
\Hom_{\cd_\infty A\+}(L,X \tsi_{A\+} e ) = 0, \quad L \in \Tria A.
\]
Comme $A \tsi_{A\+} e$ est quasi-isomorphe ‡ $0$, la seconde flËche du triangle
de $\mathbb O$-$\mathbb O$-bimodules
\[
A \tsi_{A\+} e \ra A\+ \tsi_{A\+} e \ra e \tsi_{A\+} e \ra S(A \tsi_{A\+} e)
\]
est un isomorphisme dans la catÈgorie dÈrivÈe des $A\+$-polydules
dans $\sf C(\mathbb O,\mathbb O).$  Par ailleurs, le morphisme
\[
 A\+ \tsi_{A\+} e \ra e
\]
est un quasi-isomorphisme car son cÙne, qui est la construction bar $BA\+ = \ct (SA\+)$,
est acyclique (\ref{lemme1_H-unitaire}). Ceci implique que
\[
e \ra e \tsi_{A\+} e
\]
est un isomorphisme de $A\+$-$A\+$-bipolydules dans $\sf C(\mathbb O,\mathbb O)$.
Soit $X \in \cd_\infty A\+$.
Montrons que l'objet $X \tsi_{A\+} e \in \Tria e$ est $(\Tria A)$-local.
Soit $L$ un objet de $\Tria A$ et un morphisme 
\[
f : L \ra X \tsi_{A\+} e.
\]
Nous avons un diagramme commutatif
\[
\xymatrix{
L \ar[r] \ar[d] & X \tsi_{A\+}e \ar[d]^{\sim} \\
L \tsi_{A\+} e \ar[r] & X \tsi e_{A\+} \tsi_{A\+} e,}
\]
o˘ la flËche verticale de droite reprÈsente un isomorphisme de $\cd_\infty A\+$ et
o˘ $L \tsi_{A\+} e$ est quasi-isomorphe ‡ $0.$ Le morphisme $f$ est donc nul. \findem

\subsection{La catÈgorie dÈrivÈe d'une $\ai$-algËbre strictement unitaire}
\label{section_categorie_derivee_strictement_unitaire}

Soit $A$ une $\ai$-algËbre strictement unitaire.
Dans cette section, nous donnons  plusieurs descriptions de la catÈgorie
dÈrivÈe $\cd_\infty A$ de (\ref{definition_categorie_derivee_generale}).
Plus prÈcisÈment, nous allons montrer le thÈorËme suivant :

\begin{theoreme} \label{theoreme_categorie_derivee_strictement_unitaire}
Les catÈgories suivantes sont Èquivalentes~:
\english \begin{itemize}
\item[{\rm D1.}] la catÈgorie dÈrivÈe $\cd_\infty A$ de (\ref{definition_categorie_derivee_generale}), c'est-‡-dire,
la sous-catÈgorie triangulÈe $\Tria A$ de $\cd_\infty A\+$ (\ref{proposition_categorie_derivee_H-unitaire}),
\item[{\rm D2.}] la catÈgorie (dont nous montrerons qu'elle est bien dÈfinie) \indexnotation{chinftyA}
\[
\ch_\infty A : \aiModu A/\!\sim\,
\]
o˘ $\sim$ est la relation d'homotopie (\ref{definition_polydules_strictement_unitaires}),
\item[{\rm D3.}] la catÈgorie localisÈe
\[
\big(\aiModu A \big)[Qis^{-1}]
\]
o˘ $Qis$ est la classe des $\ai$-quasi-isomorphismes de $\aiModu A$,
\item[{\rm D4.}] la catÈgorie homotopique
\[
\Ho \aiModust A
\]
de la catÈgorie de modËles $\aiModust A$ (dÈfinie plus bas).
\end{itemize} \francais
\end{theoreme}

Il en rÈsulte de ce thÈorËme que si $A$ est augmentÈe, la dÈfinition de $\cd_\infty A$ donnÈe dans
(\ref{definition_categorie_derivee_augmentee}) est Èquivalente ‡ celle de
(\ref{definition_categorie_derivee_generale}).

\begin{remarque}{\em \label{remarque_fin_categorie_derivee_II}
Les diffÈrentes descriptions de $\cd_\infty A$ montrent que les rÈsultats de la
proposition (\ref{corollaire2_categorie_derivee})
restent valides.
}\end{remarque}

{\noindent \bf L'Èquivalence entre les catÈgories de D1 et D2}\\

Comme $A$ est strictement unitaire, nous avons un
$\ai$-morphisme strictement unitaire d'$\ai$-algËbres
\[
r = \Big[ \begin{array}{c} i \\ \unite
\end{array} \Big]: A\+ = A \oplus e \ra A
\]
o˘ $\unite$ est l'unitÈ de $A$.
On a un foncteur restriction
\[
\Res : \aiModu A \ra \aiModu A\+
\]
qui est fidËle.
Nous savons que l'isomorphisme de catÈgories (\ref{section_augmentation_reduction})
\[
\aiMod A \arr{\sim} \aiModu A\+
\]
est compatible ‡ l'homotopie. La proposition (\ref{proposition_strictication_unitaire_polydules})
montre que le foncteur restriction induit un isomorphisme
\[
\Hom_{\aiModu A}(M,M')/\!\sim \, \arr{} \Hom_{\aiModu A\+}(\Res M,\Res M')/\!\sim, \quad M,M' \in \aiModu A,
\]
o˘ $\sim$ est la relation d'homotopie (\ref{definition_polydules_strictement_unitaires}).
Le corollaire (\ref{corollaire2_categorie_derivee})
dit que la relation d'homotopie (\ref{definition_polydules_strictement_unitaires}) dans $\aiModu A\+$
est une relation d'Èquivalence compatible ‡ la composition.
Ceci montre que la relation d'homotopie dans dans $\aiModu A$
est une relation d'Èquivalence compatible ‡ la composition. Nous avons donc une
catÈgorie bien dÈfinie
\[
\ch_\infty A = \aiModu A/\!\sim\,
\]
et un foncteur pleinement fidËle
\[
J : \ch_\infty A \hookrightarrow \ch_\infty A\+ \iso \cd_\infty A\+.
\]

\begin{proposition} \label{propostion_D1_D2}
Le foncteur restriction
\[
\Res : \aiModu A \ra \aiModu A\+
\]
induit une Èquivalence de catÈgorie
\[
\ch_\infty A \ra \Tria A.
\]
\end{proposition}

CommenÁons par introduire quelques notions.

\begin{definition}{\em 
Un $A\+$-polydule est {\em $H$-unitaire} \index{H-unitaire@{$H$-unitaire}}si son image par le foncteur
\[
B : \aiModu A\+ \ra \coModcu \Ba A\+
\]
est quasi-isomorphe ‡ $0.$
}\end{definition}

\begin{remarque}{\em \label{remarque1_polydules_H-unitaires}
Un $A\+$-polydule $M$ est $H$-unitaire si et seulement si l'objet $M \tsi_{A\+} e$ est quasi-isomorphe
‡ 0. La sous-catÈgories des $A\+$-polydules $H$-unitaires est donc Ègale ‡ la catÈgorie $\Tria A$ par
la proposition (\ref{proposition_categorie_derivee_H-unitaire}).
}\end{remarque}

\begin{remarque}{\em
Dans le cas o˘ $A$ est une algËbre associative unitaire et
$M$ un module unitaire, le complexe $BM$ est le cÙne de l'augmentation $\rp M \ra M$,
o˘ $\rp M$ est la rÈsolution bar de $M$ (voir par exemple \cite[IX.6]{Cartan56}).
En particulier, tout $A$-module unitaire
est un $A\+$-module $H$-unitaire. 
}\end{remarque}

\begin{lemme}\label{lemme_homologiqement_unitaire_H-unitaire}
Un $A\+$-polydule est $H$-unitaire si et seulement si il est homologiquement unitaire en tant
que $A$-polydule.
\end{lemme}

\dem 
Soit $M$ un $A$-polydule homologiquement unitaire. Il existe
une structure de $A$-polydule (nÈcessairement homologiquement unitaire)
sur $H^*M$ et un $\ai$-quasi-isomorphisme $H^*M \ra M$. Par le corollaire
\ref{theoreme1_strictification_unite_polydules}, nous pouvons choisir
$H^* M$ strictement unitaire. Nous avons alors une Èquivalence faible
\[
B(H^*M) \ra BM.
\]
Comme les Èquivalences faibles sont des quasi-isomorphismes, il nous suffit de
montrer que $B(H^*M)$ est quasi-isomorphe ‡ $0$.
Nous vÈrifions que le morphisme
\[
r : SH^*M \ts \Ba (A\+) \ra SH^*M \ts \Ba (A\+),
\]
dÈfini par les morphismes
\[
(\id \ts s\unite) : SH^*M \ts (SA) \tp i \ra SH^*M \ts (SA) \tp{i} \ts SA, \quad i\geq 0,
\]
o˘ $\unite : e \ra A$ est l'unitÈ stricte de $A$, est une homotopie contractante
de $B(H^*M)$. 

Pour dÈmontrer la rÈciproque nous introduisons quelques notions supplÈmentaires.\\

{\noindent \bf Les cochaÓnes tordantes gÈnÈralisÈes}\\

Soit $C$ un objet de $\cocoga$ et $A'$ un objet de $\aiaa.$
Une {\em cochaÓne tordante gÈnÈralisÈe} \index{cochaÓne tordante!gÈnÈralisÈe}
$\tau : C \ra A'$ est
un morphisme graduÈ de degrÈ $+1$ qui s'annule sur la 
co-augmentation $\coaugmentation^C$, qui se factorise par 
$\ker (A\+ \ra e)$ et qui vÈrifie
\[
\sum_{i\geq 1} m_i\circ  (\tau \tp i) \circ \Delta^{(i)} = 0.
\]
Remarquons que la somme infinie est bien dÈfinie car
$\tau$ s'annule sur la co-augmentation et $C$ est cocomplËte.

Soit $M$ un objet de $\aiModu A'.$ Nous munissons le produit tensoriel
$M \ts C$ du morphisme de degrÈ $+1$ qui est la somme (bien dÈfinie) de la diffÈrentielle
du produit tensoriel et des morphismes
\[
M \ts C \arr{\Id \ts \Delta^{(i)}} M \ts C \tp i \arr{\Id \ts \tau \tp {i-1} \ts \Id}
M \ts A'{}\tp{i-1} \ts C \arr{m_{i} \ts \Id} M \ts C, \quad i\geq 1.
\]
Nous vÈrifions que ce morphisme de degrÈ $+1$ est une diffÈrentielle de $M \ts C$.
Nous notons $M \tw C$ le produit tensoriel muni de cette diffÈrentielle.
Soit $N$ un objet de $\coModcu C.$ Nous munissons le produit tensoriel
$N \ts A$ de la diffÈrentielle qui est la somme (bien dÈfinie) de la diffÈrentielle
du produit tensoriel et des morphismes
\[
(\Id \ts m_{i} ) \circ (\Id \ts \tau \tp {i-1} \ts \Id) \circ (\Delta^{(i)} \ts \Id ) :
N \ts A' \ra N\ts A', \quad i \geq 1.
\]
Nous munissons $N \ts A'$ du morphisme $m_1$ donnÈe par la diffÈrentielle ci-dessus et
des morphismes $m_i$, $i\geq 2$, valant $\Id_N \ts m_i^{A'}.$ Ces morphismes dÈfinissent
une structure de $A'$-polydule sur $N \ts A'.$ Notons ce $A'$-polydule $N \tw A'.$
Ceci nous donne deux foncteurs
\[
\? \tw A' : \coModcu C \ra \aiModu A' \quad \mbox{et} \quad \? \tw C : \aiModu A' \ra \coModcu C
\]
appelÈs les {\em produits tensoriels tordus gÈnÈralisÈs}. \index{produit tensoriel tordu}\\

{\em Fin de la dÈmonstration du lemme (\ref{lemme_homologiqement_unitaire_H-unitaire})} :
Soit $M$ un $A\+$-polydule $H$-unitaire. Nous voulons montrer qu'il est homologiquement unitaire
en tant que $A$-polydule.
Nous vÈrifions que la composition
\[
\Ba A\+ = \ct SA \arr{p_1} SA \arr{\si} A \hookrightarrow A\+
\]
est un ÈlÈment tordant gÈnÈralisÈ.
Nous avons un morphisme de $A\+$-polydules
\[
\unite^{\Ba A\+} \ts \coaugmentation^{A\+} : \Ba A\+ \tw A\+ \ra e
\]
donnÈ par l'unitÈ de $\Ba A\+$ et la co-augmentation de $A\+.$
Le morphisme
\[
M \tw (\unite^{\Ba A\+} \ts \coaugmentation^{A\+}) : M \tw \Ba A\+ \tw A\+ \ra M = M\tw e
\]
est un quasi-isomorphisme (l'homotopie contractante de la dÈmonstration du lemme 
\ref{lemme_acyclicite_cochaines_universelles}
dÈfinit une homotopie contractante de son cÙne).
La co-augmentation $A\+ \ra e$ induit une suite exacte
\[
0 \ra M \tw \Ba A\+ \tw A \arr{i} M \tw \Ba A\+ \tw A\+ \ra M\tw \Ba A\+ \tw e \ra 0.
\]
Le $A\+$-polydule $M$ Ètant $H$-unitaire,
l'objet $M\tw \Ba A\+ \tw e$ est quasi-isomorphe ‡ $0$ car isomorphe ‡  $S^{-1}BM$.
Il en rÈsulte que le morphisme $i$ est un quasi-isomorphisme. Le $A\+$-polydule
$M$ qui est quasi-isomorphe ‡  $M \tw \Ba A\+ \tw A\+ $ est donc quasi-isomorphe
‡ $M \tw \Ba A\+ \tw A$. Comme ce dernier est strictement unitaire sur $A$,
$M$ est homologiquement
unitaire sur $A$.
\findem \\

{\em DÈmonstration de la proposition (\ref{propostion_D1_D2})}\\
Nous savons que le foncteur
\[
J : \ch_\infty A \hookrightarrow \ch_\infty A\+ \iso \cd_\infty A\+
\]
est pleinement fidËle. Il faut montrer que son image est formÈe des objets
de $\Tria A$.
Le lemme (\ref{lemme_homologiqement_unitaire_H-unitaire})
montre que tout objet de $\aiModu A$ est dans $\Tria A$. RÈciproquement,
si un $A\+$-polydule $M$ est dans $\Tria A$, il est homologiquement unitaire
sur $A$. Il est donc (\ref{corollaire_strictification_unite_polydules})
quasi-isomorphe ‡ un objet strictement unitaire. \findem\\

Nous munissons la catÈgorie $\ch_\infty A$ de la structure triangulÈe
induite par l'Èquivalence
\[
\ch_\infty A \ra \Tria A.
\]

{\noindent \bf Equivalence entre les catÈgories de D2 et D3}\\

Le foncteur
\[
J : \ch_\infty A \ra \ch_\infty A\+
\]
est pleinement fidËle et nous avons un isomorphisme de catÈgories
(\ref{corollaire3_categorie_derivee})
\[
\ch_\infty A\+ \arr{\sim} \cd_\infty A\+.
\]
Les $\ai$-quasi-isomorphismes sont donc les isomorphismes dans $\ch_\infty A$.
Comme
\[
\aiModu A \ra \ch_\infty A
\]
est un foncteur localisation (par rapport aux Èquivalences d'homotopie),
nous avons un isomorphisme
\[
\big(\aiModu A\big)[Qis^{-1}] \arr{\sim} \ch_\infty A.
\]

{\noindent \bf Equivalence entre les catÈgories de D3 et D4}\\

CommenÁons par montrer quelques rÈsultats sur la catÈgorie dÈrivÈe d'une
algËbre diffÈrentielle graduÈe.

\begin{lemme} \label{lemme_categorie_derivee_algebre_differentielle_graduee_unitaire}
Soit $A$ une algËbre diffÈrentielle graduÈe unitaire.
L'inclusion
\[
J : \Modu A \ra \aiModu A 
\]
induit une Èquivalence
\[
\cd A \ra \big(\aiModu A\big)[Qis^{-1}].
\]
L'inverse est donnÈ par le foncteur $? \tsi_{A}A$.
\end{lemme}

\dem
ConsidÈrons $A$ comme un $A$-$A$-bipolydule.
Nous lui associons (\ref{remarque_tsi_A_strictement_unitaire}) le foncteur
\[
? \tsi_{A} A : \aiModu A \ra \Modu A.
\]
Nous savons par le lemme (\ref{lemme2_adjonction_foncteurs_standard}) que
l'$\ai$-morphisme
\[
g_M : M \tsi_{A} A \ra M, \quad M \in \aiModu A,
\]
est un $\ai$-quasi-isomorphisme. Si $M$ est un module diffÈrentiel graduÈ sur $A$,
les multiplications $m^M_i$, $i \geq 3$, sont nulles et l'$\ai$-morphisme $g_M$
(construit dans la dÈmonstration du lemme (\ref{lemme2_adjonction_foncteurs_standard}))
est strict. L'$\ai$-morphisme $g_M$ est alors un morphisme de $A$-modules diffÈrentiels graduÈs.
Ceci montre que les foncteurs $J$ et $? \tsi_{A} A$ induisent des foncteurs
quasi-inverses l'un de l'autre entre les catÈgories
\[
\cd A \quad \mbox{et} \quad \big(\aiModu A\big)[Qis^{-1}].
\]
\findem

\begin{definition}{\em
Soit $A$ une algËbre diffÈrentielle graduÈe (non nÈcessairement unitaire).
La {\em catÈgorie dÈrivÈe} $\cd A$ \index{categorie@{catÈgorie}!derivee@{dÈrivÈe}} est le
noyau de
\[
? \stackrel{\bf L}{\ts} e : \cd A\+ \ra \cd e. 
\]
}\end{definition}

\begin{remarque}{\em
Dans le cas o˘ $A$ est unitaire, la catÈgorie dÈrivÈe dÈfinie ci-dessus est Èquivalente ‡ la catÈgorie
dÈrivÈe dÈfinie en (\ref{structure_triangulee_Ho_coModcu}).
}\end{remarque}

\begin{corollaire} \label{corollaire_categorie_derivee_algebre_differentielle_graduee}
Soit $A$ une algËbre diffÈrentielle graduÈe (non nÈcessairement unitaire).
Les catÈgories dÈrivÈes $\cd_\infty A$ et $\cd A$ sont Èquivalentes.
\end{corollaire}

\dem C'est une consÈquence du lemme (\ref{lemme_categorie_derivee_algebre_differentielle_graduee_unitaire})
et du fait que le foncteur $? \stackrel{\bf L}{\ts} e$ est exactement le
foncteur $? \tsi_{A} e$. \findem\\

{\noindent \bf La catÈgorie de modËles $\aiModust A$}\\

Nous utilisons ci-dessous les notations et la terminologie standard des opÈrades
diffÈrentielles graduÈes (voir par exemple \cite{Hinich97c}).

Une {\em opÈrade asymÈtrique} est une suite d'objets $\co (n)$, $n \geq 0$, de $\cc \sf C$
munie d'une composition $\mu$ vÈrifiant les mÍmes conditions d'associativitÈ
que la composition d'une opÈrade au sens habituel. Notons $\frak S_n$, $n \geq 1,$ le groupe
symÈtrique. La suite $\corps \frak S_n \ts_{\corps} \co(n)$, $n\geq 0$,
est un $S$-module dans $\cc \sf C$ et $\mu$ induit une structure d'opÈrade sur ce $S$-module.
L'opÈrade $\Ass$ des algËbres associative est Ègale ‡ $\corps \frak S_n \ts_{\corps} \Ass'(n)$, $n\geq 0$,
o˘ $\Ass'$ est une opÈrade asymÈtrique.

Soit $\co$ {\em l'opÈrade asymÈtrique des $\ai$-algËbres strictement unitaires}.
On note $\E(\co,A) = \E(A)$ l'algËbre enveloppante \index{algebre@{algËbre}!enveloppante}
de $A$ relativement ‡ l'opÈrade $\co$.
La catÈgorie $\aiModust A$ des $A$-polydules strictement unitaires
dont les morphismes sont les $\ai$-morphismes stricts est bien s˚r isomorphe ‡ la
catÈgorie des modules (‡ droite) sur la $\co$-algËbre $A$.
Nous avons donc un isomorphisme de catÈgories
\[
\Modu \E(A) \arr{\sim} \aiModust A.
\]
Nous dÈduisons du thÈorËme (\ref{theoreme_cmf_Modu}) le rÈsultat suivant.

\begin{proposition}
Les trois classes de morphismes ci-dessous dÈfinissent une structure de catÈgorie de modËles
sur $\aiModust A$ :
\english \begin{itemize}
\item[-] la classe $\weq$ formÈe des $\ai$-quasi-isomorphismes stricts,
\item[-] la classe $\fib$ formÈe des morphismes $f : M \ra M'$ tels que
$f^n$ est un Èpimorphis\-me pour tout $n \in \Z,$
\item[-] la classe $\cof$ formÈe des morphismes qui ont la propriÈtÈ de relËvement
‡ gauche par rapport aux morphismes appartenant ‡ $Qis \cap \fib$.
\end{itemize} \francais
\findem
\end{proposition}

Nous rappelons que la catÈgorie dÈrivÈe $\cd \E (A)$ est isomorphe ‡ la catÈgorie localisÈe
\[
\Ho \big(\aiModust A\big).
\]

\begin{remarque}{\em
Si $A$ est une $\ai$-algËbre augmentÈe, l'algËbre enveloppante $\E(A)$ est isomorphe ‡
$\Oma \Ba A$ (voir \ref{definition_algebre_enveloppante}).
}\end{remarque}

\begin{proposition} \label{proposition_D3_D4}
Soit $A$ une $\ai$-algËbre strictement unitaire.
L'inclusion
\[
J : \aiModust A \ra \aiModu A
\]
induit une Èquivalence
\[
\Ho \big(\aiModust A\big) \ra \big(\aiModu A\big)[Qis^{-1}].
\]
\end{proposition}

\dem 

{\em Premier cas : $A$ est une algËbre diffÈrentielle graduÈe unitaire.}
La suite d'inclusions
\[
\Modu A \hookrightarrow \aiModust A \hookrightarrow \aiModu A
\]
induit une suite de foncteurs fidËles
\[
\cd A \ra \Ho \big(\aiModust A\big) \ra \big(\aiModu A\big)[Qis^{-1}].
\]
Le lemme (\ref{lemme_categorie_derivee_algebre_differentielle_graduee_unitaire})
nous donne la pleine fidÈlitÈ de la composition. Le second foncteur est
donc plein et nous avons le rÈsultat.

{\em DeuxiËme cas : $A$ est une $\ai$-algËbre strictement unitaire quelconque.} \\
D'aprËs la proposition (\ref{proposition_modele_differentiel_gradue}), il
existe un modËle diffÈrentiel graduÈ unitaire $A'$ et une cofibration triviale
\[
i : A \ra A'
\]
strictement unitaire.
Le lemme
(\ref{lemme_(co)fibrations_triviales_strictement_unitaires}) montre qu'il
existe une fibration triviale $q : A' \ra A$ telle que $q \circ i = \Id_A$
et $i \circ q$ est homotope ‡ $\Id_{A'}$.
Les foncteurs restriction $\Res^i$ et $\Res^q$ induisent des foncteurs $i^*$ et $q^*$
entre les catÈgories homotopiques
\[
\Ho \big(\aiModust A\big) \quad \mbox{et} \quad \Ho \big(\aiModust A'\big).
\]
Nous avons clairement $i^* \circ q^* = \Id$. Montrons que $q^* \circ i^*$ est isomorphe au foncteur
identitÈ de $\Ho \big(\aiModust A'\big).$

Soit $A'\+$ l'augmentation de $A'$. Son algËbre enveloppante $\E(A'\+)$
est l'algËbre diffÈrentielle graduÈe $\Oma \Ba A'\+$ (voir \ref{A-polydules_EA-modules}).
Soit $j : A'\+ \ra \E(A'\+)$ l'$\ai$-morphisme universel construit en \ref{lemme_algebre_enveloppante}.
Comme il est un $\ai$-quasi-isomorphisme strictement unitaire augmentÈ, il induit
une Èquivalence
\[
\cd_\infty \E(A'\+)  \ra \cd_{\infty} A'\+
\]
compatible aux foncteurs
\[
\cd_\infty A'\+ \ra \cd_\infty e \quad \mbox{et} \quad \cd_\infty \E(A'\+) \ra \cd_\infty e.
\]
La sous-catÈgorie $\cd_\infty A' = \Tria A'$ est ainsi Èquivalente ‡ la sous-catÈgorie
$\cd_\infty \b \E(A'\+) = \Tria \b \E(A'\+)$ (l'algËbre $\b \E(A'\+) = \Omega \Ba A'\+$ est la rÈduction
de $\E(A'\+)$).
Notons $f$ l'$\ai$-morphisme composÈ $i \circ q.$ Soit $f\+ : A'\+ \ra A'\+$
le morphisme augmentÈ associÈ ‡ $f$.  Notons $g$ le morphisme 
\[
\Oma \Ba f\+ : \E(A'\+) \ra \E(A'\+).
\]
Le morphisme $g$ est un morphisme d'algËbres diffÈrentielles graduÈes unitaires.
Pour montrer que $\Res^f$ induit un endofoncteur de $\Ho \aiModust A'$
qui est isomorphe au foncteur identitÈ, il suffit de montrer que $\Res^g$
induit un endofoncteur de $\cd \E(A'\+)$ isomorphe au foncteur identitÈ.
Le morphisme $g$ est clairement homotope ‡ $\Id$ dans
la catÈgorie $\aia$. Les morphismes $g$ et $\Id$ deviennent donc
Ègaux dans $\alg [Qis^{-1}]$ (voir \ref{corollaire_cmf_aia}).
Comme $\Oma \Ba A'\+$ est une algËbre presque libre co-augmentÈe,
elle est un objet fibrant et cofibrant de la catÈgorie de modËles $\alg$ (voir \ref{section_cmf_alg}).
Il existe donc une homotopie ‡ droite entre $\Id$ et $g$.
Le lemme (\ref{lemme_homotopie_foncteur_restriction}) ci-dessous montre que l'endofoncteur $g^*$
de $\cd \b \E(A'\+)$ induit par $\Res^g$ est isomorphe
‡ l'identitÈ. \findem

\begin{lemme} \label{lemme_homotopie_foncteur_restriction}
Soit $A$ et $B$ deux algËbres diffÈrentielles graduÈes unitaires.
Soit $g$ et $g'$ deux morphismes $A \ra B$ unitaires homotopes ‡ droite.
Les foncteurs restriction le long de $g$ et $g'$ induisent des foncteurs
isomorphes
\[
\cd B \ra \cd A.
\]
\end{lemme}

\dem 
On rappelle qu'une algËbre de chemins $B^I$, c'est-‡-dire un objet de chemins
pour $B$ dans la catÈgorie de modËles $\alg$, est un objet de $\alg$
muni de morphismes
\[
B \arr{i} B^I \arr{p} B_0 \times B_1,
\]
o˘ $B_0$ et $B_1$ sont Ègales ‡ $B$, tels que $i$ est une Èquivalence
faible et $p \circ i$ est une factorisation de la diagonale $B \ra B_0 \times B_1$.
Notons $p_0$ et $p_1$ les morphismes composÈs
\[
B^I \arr{p} B_0 \times B_1 \ra B_0 \quad \mbox{et} \quad B^I \arr{p} B_0 \times B_1 \ra B_1.
\]
Nous avons les ÈgalitÈs $p_0 \circ i = p_1 \circ i = \Id$.

Les morphismes $g$ et $g'$ sont homotopes ‡ droite relativement ‡ l'algËbre de chemins
$B^I$, il existe donc un morphisme $H : A \ra B^I$
tel que $p_0 \circ H = g$ et $p_1 \circ H = g'$.
Ceci montre que
\[
\Res^{g} = \Res^{H} \circ \Res^{p_0} \quad \mbox{et} \quad \Res^{g'} = \Res^{H} \circ \Res^{p_1}.
\]
Pour montrer que $\Res^{g}$ et $\Res^{g'}$ induisent des foncteurs isomorphes dans les catÈgories
dÈrivÈes, il suffit de montrer que $\Res^{p_0}$ et $\Res^{p_1}$ induisent des foncteurs
isomorphes dans les catÈgories dÈrivÈes. Nous avons les ÈgalitÈs
\[
\Id = \Res^i \circ \Res^{p_0} = \Res^i \circ \Res^{p_1}.
\]
Comme $i$ est un quasi-isomorphisme, $\Res^i$ induit une Èquivalence dans les catÈgories
dÈrivÈes. Nous en dÈduisons que $\Res^{p_0}$ et $\Res^{p_1}$ induisent des foncteurs
isomorphes dans les catÈgories dÈrivÈes. \findem

\section{La catÈgorie dÈrivÈe des bipolydules}
\label{section_categorie_derivee_bipolydules}

Les dÈmonstrations de cette section sont omises car similaires
‡ celles de la section \ref{section_categorie_derivee_polydules}.\\

{\noindent \bf Le foncteur $M \tsi \,? \tsi M''$}\\

Soit $\mathbb O$, $\mathbb O'$, $\mathbb O''$ et $\mathbb O'''$ des objets de $\sf C.$
Soit $A$ (resp.~$A'$, $A''$, $A'''$) une $\ai$-algËbre dans $\sf C(\mathbb O, \mathbb O)$
(resp.~$\sf C(\mathbb O', \mathbb O')$, $\sf C(\mathbb O'', \mathbb O'')$, $\sf C(\mathbb O''', \mathbb O''')$).
Soit $M$ (resp.~$M''$) un $A$-$A'$-bipolydule
(resp.~$A''$-$A'''$-bipolydule) dans $\sf C(\mathbb O, \mathbb O')$
(resp.~$\sf C(\mathbb O'', \mathbb O''')$).
On dÈfinit le foncteur
\[
\aiMod (A',A'') \ra \aiMod (A,A'''), \quad M' \mapsto M \tsi M' \tsi M,
\]
par l'ÈgalitÈ de $\Ba A$-$\Ba A'''$-bicomodules diffÈrentiels graduÈs 
\[
B(M \tsi M' \tsi M) = BM \tbox_{\Ba A'} BM' \tbox_{\Ba A''} BM,
\]
o˘ $\tbox$ dÈsigne le produit cotensoriel (voir \ref{section_foncteurs_standard}).\\

{\noindent \bf La catÈgorie dÈrivÈe $\cd_\infty (A,A')$}\\

Soit $e_\mathbb O$ et $e_{\mathbb O'}$ les ÈlÈments neutres de $\sf C(\mathbb O,\mathbb O)$
et $\sf C(\mathbb O',\mathbb O')$ considÈrÈs comme des $\ai$-algËbres augmentÈes.
ConsidÈrons $e_\mathbb O$ et $e_{\mathbb O'}$ comme un $e_\mathbb O$-$A\+$-bipolydule
et un $A'\+$-$e_{\mathbb O'}$-bipolydule.

\begin{definition}{\em \label{definition_categorie_derivee_bipolydules_generale}
La {\em catÈgorie dÈrivÈe}  $\cd_\infty (A',A'')$ est le noyau du foncteur
\index{categorie@{catÈgorie}!derivee@{dÈrivÈe}}
\[
e_\mathbb O \tsi_{A'\+} \, ? \tsi_{A''\+} e_{\mathbb O'} : \cd_\infty (A'\+,A''\+) \ra \cd_\infty 
(e_\mathbb O,e_{\mathbb O'}).
\]
}\end{definition}

\newpage
{\noindent \bf La sous-catÈgorie $\Tria (A , A')$}\\

Supposons que la catÈgorie $\sf C$ admet un objet final $\mathbb P.$
Soit $x : \mathbb P \ra \mathbb O$ un morphisme de $\sf C$.
Le morphisme $x$ induit un foncteur
\[
x_* : \sf C(\mathbb O,\mathbb P) \ra \sf C(\mathbb P,\mathbb P), \quad M \mapsto M(x).
\]
Nous supposons que ce foncteur admet un adjoint ‡ gauche
\[
\lrus{!}{x}{} : \sf C(\mathbb P,\mathbb P) \ra \sf C(\mathbb P,\mathbb O).
\]
Nous avons un $A$-polydule ‡ gauche
\[
x^\vee = A \ts_{\mathbb O} \lrus{!}{x}{}(e_{\mathbb P}),
\]
dont la structure est donnÈe par les multiplications de $A$.

\begin{remarque}{\em 
Cette notation est justifiÈe par le fait suivant. 
Dans l'exemple apparaissant dans l'Ètude des $\ai$-catÈgories
(\ref{section_categories_de_base_aicat}), un objet final est
un ensemble ‡ un ÈlÈment. Soit  $\mathbb P$ un tel ensemble et
$\mathbb O$ un ensemble.
Soit $x$ une l'application $\mathbb P \ra \mathbb O$ donnÈe par
un ÈlÈment (notÈ aussi $x$) de $\mathbb O.$
le $A$-polydule $x^\vee$ est l'$\ai$-foncteur coreprÈsentÈ par $x$
\[
x^\vee = A(x,?).
\]
\findem
}\end{remarque}

Soit $x : \mathbb P \ra \mathbb O$ et $y : \mathbb P \ra \mathbb O'$ des morphismes de $\sf C$.
Le $\mathbb O$-$\mathbb O'$-bimodule
\[
x^\vee \ts_{\mathbb P} y \pw = A \ts_{\mathbb O'} \lrus{!}{x}{}(e_{\mathbb P})
\ts_{\mathbb P} y_!(e_{\mathbb P}) \ts_{\mathbb O'} A'
\]
est un $A$-$A'$-bipolydule.
La catÈgorie $\Tria (A , A')$ est la sous-catÈgorie triangulÈe de $\cd_\infty (A\+,A'\+)$ engendrÈe par
les
\[
 x^\vee \ts_{\mathbb P} y \pw ,\quad x \in \sf C(\mathbb P,\mathbb O), \quad y \in \sf C(\mathbb P,\mathbb O').
\]

\begin{proposition} \label{proposition_categorie_derivee_bipolydules_H-unitaire}
Soit $A$ et $A'$ des $\ai$-algËbres $H$-unitaires.
Nous avons une suite exacte de catÈgories triangulÈes
\[
\Tria (A , A') \hookrightarrow \cd_\infty (A\+,A'\+) \ra \cd_\infty (e_\mathbb O,e_\mathbb O').
\]
En particulier, la catÈgorie dÈrivÈe $\cd_\infty A$ est Ègale ‡ $\Tria (A , A').$
\findem
\end{proposition}

\begin{theoreme} \label{theoreme_categorie_derivee_strictement_unitaire_bipolydules}
Soit $A$ et $A'$ deux $\ai$-algËbres strictement unitaires.
Les catÈgories suivantes sont Èquivalentes~:
\english \begin{itemize}
\item[{\rm D1.}] la catÈgorie dÈrivÈe $\cd_\infty (A,A')$ de
(\ref{definition_categorie_derivee_bipolydules_generale}), c'est-‡-dire,
la sous-catÈgorie triangulÈe $\Tria (A , A')$ de $\cd_\infty (A\+,A'\+)$
(\ref{proposition_categorie_derivee_H-unitaire}),
\item[{\rm D2.}] la catÈgorie (bien dÈfinie) \indexnotation{chinftyAA'}
\[
\ch_\infty (A,A') = \aiModu (A,A')/\!\sim\,
\]
o˘ $\sim$ est la relation d'homotopie,
\item[{\rm D3.}] la catÈgorie localisÈe \indexnotation{cdinftyAA'}
\[
\big(\aiModu (A,A') \big)[Qis^{-1}]
\]
o˘ $Qis$ est la classe des $\ai$-quasi-isomorphismes de $\aiModu (A,A')$,
\item[{\rm D4.}] la catÈgorie localisÈe
\[
\big(\aiModust (A,A')\big)[Qis^{-1}]
\]
de la catÈgorie $\aiModust (A,A').$
\end{itemize} \francais
\end{theoreme}

\dem
Les Èquivalences entre les catÈgories de D1, D2 et D3 se montrent de la mÍme maniËre
que dans le thÈorËme (\ref{theoreme_categorie_derivee_strictement_unitaire}).
L'Èquivalence entre les catÈgories de D3 et D4 dans
le cas o˘ $A$ et $A'$ sont des algËbres diffÈrentielles graduÈes
unitaires se prouve comme dans la proposition (\ref{proposition_D3_D4}).
Si $A$ et $A'$ dont des $\ai$-algËbres strictement unitaires quelconques, nous procÈdons
de la maniËre suivante. On montre comme dans la proposition (\ref{proposition_D3_D4}) que
l'inclusion
\[
\Modu (\E(A),\E(A')) \hookrightarrow \aiModu (A,A')
\]
induit une Èquivalence
\[
\Ho \big(\Modu (\E(A),\E(A'))\big) \ra \big(\aiModu (A,A') \big)[Qis^{-1}].
\]
Comme cette Èquivalence est la composition des foncteurs fidËles
\[
\Ho \big(\Modu (\E(A),\E(A'))\big) \ra \big(\aiModust (A,A') \big)[Qis^{-1}] \arr{K} \big(\aiModu (A,A') \big)[Qis^{-1}]
\]
le foncteur $K$ est plein. Il est donc une Èquivalence.
\findem

%% file: Ai_cat.tex
{\bf \noindent Plan du chapitre}\\

\noindent Une $\ai$-catÈgorie est une $\ai$-algËbre avec plusieurs objets, et rÈciproquement, une
$\ai$-algËbre est une $\ai$-catÈgorie avec un objet. Les problËmes soulevÈs
par l'augmentation du nombre d'objets sont nombreux et la gÈnÈralisation des
rÈsultats des chapitres prÈcÈdents est parfois trËs technique.

Dans la  section \ref{section_preliminaire_ai-categories},
nous fixons des notations qui codent la variation des ensembles d'objets des petites $\ai$-catÈgories.
Nous introduisons pour cela une bicatÈgorie
$\sf C$ dont les objets sont les ensembles, puis nous dÈfinissons
une petite $\ai$-catÈgorie dont l'ensemble des objets est en bijection avec
un ensemble $\mathbb O$
comme une $\ai$-algËbre dans la catÈgorie (monoÔdale)
$\sf C(\mathbb O,\mathbb O).$ Nous dÈfinissons ensuite les $\ai$-foncteurs. 

 Dans la section \ref{section_categorie_dg_polydules}, nous dÈfinissons
les catÈgories diffÈrentielles graduÈes des (bi)polydules sur des $\ai$-catÈgories.

Dans la section \ref{section_lemme_clef},
nous Ètablissons un lemme (dit {\em lemme clef}) qui sera fondamental dans la construction de
l'$\ai$-foncteur de Yoneda (\ref{definition_ai-foncteur_de_Yoneda}) et celle de l'$\ai$-foncteur
de Yoneda gÈnÈralisÈ (\ref{section_ai-foncteur_de_Yoneda_generalise}).

\section{DÈfinitions}
\label{section_preliminaire_ai-categories}

\subsection{Les catÈgories de base $\sf C({\mathbb O},{\mathbb O}')$ et $\sf C({\mathbb O})$} 
\label{section_categories_de_base_aicat}

Nous fixons des notations que nous utiliserons tout au long de cette partie.
Nous construisons une bicatÈgorie $\sf C$ \indexnotation{C3} dont les objets sont les ensembles 
(voir \cite[Chap.~XII, ß6]{MacLane98} pour les bicatÈgories).\\

Soit $\corps$ un corps. Le produit tensoriel au-dessus de $\corps$ est
notÈ $\ts$.
Soit ${\mathbb O}$ un ensemble. \indexnotation{bbO}
ConsidÈrons le comme la {\em petite catÈgorie} dont
les objets sont en bijection avec ${\mathbb O}$ et dont l'espace des
morphismes $o \ra o'$ est vide si $o \neq o'$, et contient
uniquement le morphisme identitÈ $\id_o$ sinon. 

Soit ${\mathbb O}$, ${\mathbb O}'$ et ${\mathbb O}''$ trois ensembles.
Un {\em ${\mathbb O}'$-${\mathbb O}$-bimodule} \index{bimodule@{${\mathbb A}$-${\mathbb B}$-bimodule}}
(resp.~un {\em ${\mathbb O}$-module} \index{module@{${\mathbb A}$-module}}
‡ droite) est un foncteur
\[
M : {\mathbb O} \op \times {\mathbb O}' \ra \vect \corps,\quad \Big(\mbox{resp.}\quad M : {\mathbb O} \op  \ra \vect \corps\Big)
\]
o˘ $\vect \corps$ est la catÈgorie des $\corps$-espaces vectoriels.
Un {\em morphisme} de bimodules (resp.~de modules) est un morphisme de foncteurs.
Nous notons $\sf C({\mathbb O},{\mathbb O}')$ et $\sf C({\mathbb O})$ ces catÈgories.
Soit $M$ un objet de $\sf C({\mathbb O},{\mathbb O}')$ et $N$ un objet de $\sf C({\mathbb O}',{\mathbb O}'')$.
Le {\em produit tensoriel $M \tso_{{\mathbb O}'} N$ \indexnotation{tso}
au-dessus de ${\mathbb O}'$} est l'objet de $\sf C ({\mathbb O},{\mathbb O}'')$
dÈfini par
\[
\Big(M \tso_{{\mathbb O}'} N\Big)(o'',o) = \bigoplus_{o'\in {\mathbb O}'}M(o',o) \ts  N(o'',o').
\]
Nous noterons simplement $\odot$ le tenseur au-dessus de ${\mathbb O}'$ lorsque cela ne prÍtera pas
‡ confusion.
Le produit tensoriel au-dessus de ${\mathbb O}'$ nous donne un foncteur
\[
\sf C({\mathbb O}',{\mathbb O}) \times \sf C({\mathbb O}'',{\mathbb O}')   \ra \sf C ({\mathbb O}'',{\mathbb O}'), \quad
(M, N )\mapsto M \tso_{{\mathbb O}'} N,
\]
et si ${\mathbb O}'''$ est un ensemble
et $T$ un objet de $\sf C({\mathbb O}'',{\mathbb O}'''),$ on a des contraintes d'associativitÈ
\[
(M \tso_{{\mathbb O}'} N) \tso_{{\mathbb O}''} T \arr{\sim} M \tso_{{\mathbb O}'} (N \tso_{{\mathbb O}''} T).
\]
Soit $f: {\mathbb O} \ra {\mathbb O}'$ une application.
On a un foncteur
\[
\sf C({\mathbb O}'',{\mathbb O}')  \longrightarrow   \sf C({\mathbb O}'',{\mathbb O}),
\]
qui envoie le ${\mathbb O}'$-${\mathbb O}''$-bimodule $M$ sur le ${\mathbb O}$-${\mathbb O}''$-bimodule
\[
\begin{array}{rcl} M_{f} : {\mathbb O} \op \times {\mathbb O}'' &\ra &\vect \corps\\
o \times  o'' & \mapsto & M(f o, o'').
\end{array}
\]
De maniËre similaire, si $g : {\mathbb O} \ra {\mathbb O}''$ est une application, on a un
foncteur
\[
\sf C({\mathbb O}'',{\mathbb O}')  \longrightarrow   \sf C({\mathbb O},{\mathbb O}'), \quad M \mapsto \lrus{g}{M}{}.
\]
La catÈgorie $\sf C({\mathbb O},{\mathbb O}')$ est $\corps$-linÈaire abÈlienne, semi-simple,
cocomplËte, aux co\-limites filtrantes exactes (i.e.~c'est une $\corps$-catÈgorie de Grothen\-dieck
semi-simple). 
Par la section \ref{categorie_de_base}, nous avons les catÈgories $\gr \sf C({\mathbb O},{\mathbb O}')$
des {\em bimodules graduÈs} et $\cc \sf C({\mathbb O},{\mathbb O}')$ des {\em bimodules diffÈrentiels graduÈs}.
Remarquons que le produit tensoriel $\tso_{\mathbb O}$ et le bimodule
\[
e_{\mathbb O}(\?,\?) = \corps  \Hom_{{\mathbb O}}(\?,\?)
\]
dÈfinissent une structure de catÈgorie monoÔdale sur $(\sf C({\mathbb O},{\mathbb O}),\tso ,e_{\mathbb O})$.
Le foncteur
\[
\sf C({\mathbb O}',{\mathbb O}') \ra \sf C({\mathbb O},{\mathbb O}), \quad M \mapsto \lrus{f}{M}{f},
\]
est compatible ‡ la structure monoÔdale.
Comme la catÈgorie $\sf C({\mathbb O})$ est isomorphe ‡  $\sf C(\{*\},{\mathbb O})$,
o˘ $\{*\}$ est un ensemble ‡ un ÈlÈment,
nous obtenons une action ‡ droite de la catÈgorie monoÔdale $C({\mathbb O},{\mathbb O})$ sur $C({\mathbb O})$
\[
\sf C({\mathbb O})  \times  \sf C({\mathbb O},{\mathbb O}) \ra \sf C ({\mathbb O}), \quad
(M, N)\mapsto M \tso N.
\]

\begin{remarque}{\em
Soit $\ca$ une petite $\corps$-catÈgorie dont l'ensemble des objets est en bijection avec
un ensemble ${\mathbb A}$.
Le ${\mathbb A}$-${\mathbb A}$-bimodule
\[
\Hom_\ca : A \times A' \mapsto \Hom_\ca(A,A'),
\]
muni des morphismes
\[
\mu : \Hom_\ca \tso \Hom_\ca \ra \Hom_\ca \quad \mbox{et}
\quad \unite : e_{\mathbb A} \ra \Hom_\ca, \quad \id_A \mapsto \Id_A,
\]
donnÈs par la composition de $\ca$ et par les morphismes identitÈ $\Id_A$ de $\ca$,
est une algËbre unitaire dans la catÈgorie des ${\mathbb A}$-${\mathbb A}$-bimodules.
RÈciproquement, une algËbre unitaire dans la catÈgorie des ${\mathbb A}$-${\mathbb A}$-bimodules dÈfinit une
petite $\corps$-catÈgorie dont l'ensemble des objets est en bijection avec ${\mathbb A}$.

Soit $\ca$ et $\cb$ deux petites $\corps$-catÈgories dont les ensembles des objets sont
en bijection avec des ensembles ${\mathbb A}$ et ${\mathbb B}.$
Soit $f : \ca \ra \cb$ un foncteur. On note
\[
\dot f  : \Obj \ca \ra \Obj \cb
\]
l'application qui envoie $A$ sur son image par le foncteur $f.$
Le foncteur $f$ induit un morphisme d'algËbres unitaires
\[
\Hom_{\ca} \ra \lrus{\dot f}{\Hom_{\cb}}{\dot f}, \quad x \mapsto f(x),
\]
RÈciproquement, si $\Lambda$ et $\Lambda'$ sont deux algËbres unitaires dans les catÈgories des
${\mathbb A}$-${\mathbb A}$-bimodules
et des ${\mathbb B}$-${\mathbb B}$-bimodules, une application  $\dot f : {\mathbb A} \ra {\mathbb B}$
et un morphisme d'algËbres unitaires
$\Lambda \ra \lrus{\dot f}{\Lambda'}{\dot f}$ dans la catÈgorie $\sf C({\mathbb A},{\mathbb A})$
dÈfinissent un foncteur entre les $\corps$-catÈgories correspondantes.
 }\end{remarque}

\begin{definition}{\em Soit ${\mathbb A}$ un ensemble.
Une {\em (petite) catÈgorie diffÈrentielle graduÈe sur ${\mathbb A}$}
\index{categorie@{catÈgorie}!diffÈrentielle graduÈe} est
une algËbre diffÈrentielle graduÈe unitaire dans $\sf C({\mathbb A},{\mathbb A}).$
}\end{definition}

%
%
%
%
%
%
%
%
%
%
%
%
%
%

\subsection{DÈfinitions}
\label{section_definitions_ai-categories}

\begin{definition}{\em \label{definition_ai-categorie}
Soit ${\mathbb A}$ un ensemble.  \indexnotation{ca} \indexnotation{bbA}
Une {\em $\ai$-catÈgorie sur ${\mathbb A}$}
\index{A-infini categorie@{$\ai$-catÈgorie}} est une $\ai$-algËbre dans la catÈgorie 
\[
(\gr \sf C({\mathbb A},{\mathbb A}),\tso, e_{\mathbb A}).
\]
}\end{definition}

\begin{remarque}{\em
Soit $\ca$ une $\ai$-catÈgorie.
Elle est dÈterminÈe par
\begin{itemize}
\item[-] un ensemble d'{\em objets} $\Obj \ca = {\mathbb A}$,
\item[-] pour tout couple $(A,A')$ d'objets de $\ca$, un espace graduÈ de {\em morphismes}
\[
\Hom_{\ca}(A,A') = \ca(A,A'),
\]
\item[-] pour toute suite $(A_0,\hdots,A_n)$ d'objets de $\ca$,
des {\em compositions}
\[
m_n : \ca(A_{n-1},A_n) \ts \hdots \ts \ca (A_0,A_1) \ra \ca(A_0,A_n),
\]
 vÈrifiant les Èquations $(*_n)$, $n \geq 1$, de la dÈfinition \ref{definition_ai-algebre}, 
\end{itemize}
Si $\ca$ est homologiquement unitaire (en tant que $\ai$-algËbre dans $\gr \sf C({\mathbb A},{\mathbb A})$)
alors, pour tout objet $A \in \ca$, nous avons un morphisme {\em identitÈ} \index{identitÈ!d'un objet}
$\id_A \in \ca(A,A)$  tel que
sa classe $[\id_A]$ dans $H^*\!\ca(A,A)$ vÈrifie
\[
\mu (f,[\id_A]) = f, \quad f \in H^*\!\ca(A,A') \quad \mbox{et} \quad 
\mu ([\id_A],g) = g,  \quad g \in H^*\!\ca(A',A),
\]
o˘ $\mu$ est la composition de $H^*\!\ca$ induite par $m_2$.
}\end{remarque}

\begin{remarque}{\em
La composition $m_2$ induit une
composition associative
\[
\mu : H^0 \ca \tso H^0 \ca \ra H^0 \ca.
\]
Si $\ca$ est homologiquement unitaire alors $H^0 \ca$ est une catÈgorie au sens classique. Le morphisme
identitÈ d'un objet $A \in H^0 \ca$ est la classe $[\id_A].$
}\end{remarque}

\begin{lemme}\label{lemme_ai-categorie_fBf} \indexnotation{fBf}
Soit ${\mathbb B}$ un ensemble, $\cb$ une $\ai$-catÈgorie sur ${\mathbb B}$ homologiquement unitaire
et
\[
f : {\mathbb A} \ra {\mathbb B}
\]
une application.
Le ${\mathbb A}$-${\mathbb A}$-bimodule graduÈ $\lrus{f}{\cb}{f}$ est une $\ai$-catÈgorie
homologiquement unitaire
pour les compositions et les morphismes identitÈ induits par ceux de $\ca.$
\findem
\end{lemme}

%
%
%
%
%
%
%
%

\begin{definition}{\em \label{definition_ai-foncteurs} \indexnotation{f} \indexnotation{dotf}
Soit ${\mathbb A}$ et ${\mathbb B}$ deux ensembles et $\ca$ et $\cb$ deux $\ai$-catÈgories
sur ${\mathbb A}$ et ${\mathbb B}$. Un {\em $\ai$-foncteur}
\index{A-infini foncteur@{$\ai$-foncteur}}
\[
f : \ca \ra \cb
\]
est la donnÈe d'un couple $(\dot f,f_{\Hom})$ formÈ d'une application
\[
\dot f : {\mathbb A}\ra {\mathbb B}
\]
et d'un $\ai$-morphisme dans la catÈgorie $\gr \sf C({\mathbb A},{\mathbb A})$
\[
f_{\Hom} : \ca \ra \lrus{\dot f}{\cb}{\dot f}.
\]
Nous noterons souvent ce dernier $f$ au lieu de $f_\Hom$.
L'$\ai$-foncteur {\em identitÈ de $\ca$} \index{identitÈ!le foncteur} \indexnotation{Id}
est notÈ
\[
\Id_\ca : \ca \ra \ca.
\]
Attention ‡ ne pas confondre ce symbole avec
$\id_{A}$, le morphisme identitÈ d'un objet $A \in \ca$. \indexnotation{id}
}\end{definition}

\begin{remarque}{\em
Soit $\ca$ et $\cb$ deux petites $\ai$-catÈgories.
Un $\ai$-foncteur $f : \ca \ra \cb$ est dÈterminÈ par
\begin{itemize}
\item[-] une application
\[
\dot f : \Obj \ca \ra \Obj \cb,
\]
\item[-]  pour toute suite $(A_0,\hdots,A_n)$ d'objets de $\ca$,
des morphismes
\[
f_n : \ca(A_{n-1},A_n) \ts \hdots \ts \ca (A_0,A_1) \ra \cb(\dot f A_0,\dot f A_n),
\]
vÈrifiant les Èquations $(**_n)$, $n \geq 1$, de la dÈfinition \ref{definition_ai-morphisme}.
\end{itemize}
}\end{remarque}

\begin{remarque}{\em \label{remarque_ai-morphismes_ai-foncteurs}
Soit $\ca$ et $\cb$ deux petites $\ai$-catÈgories sur $\mathbb A$.
Un $\ai$-morphisme $f : \ca \ra \cb$ dans $\mathsf C(\mathbb A,\mathbb A)$
donne un $\ai$-foncteur
\[
(\Id_{\mathbb A}, f) : \ca \ra \cb, \quad x \mapsto f(x).
\]
RÈciproquement, un $\ai$-foncteur $(\dot f, f)$
dont l'application sous-jacente $\dot f$ est Ègale ‡ $\Id_{\mathbb A}$ donne un $\ai$-morphisme
$f : \ca \ra \cb$.
}\end{remarque}

{\noindent \bf Rappel sur la construction bar}\\

Soit $\ca$ et $\cb$ deux $\ai$-catÈgories et $f : \ca \ra \cb$ un $\ai$-foncteur.
Rappelons que les bijections de la section \ref{construction_bar_cobar},
\[
m_i \lra b_i \quad \Big(\mbox{resp.} \quad f_i  \lra F_i \Big), \quad i\geq 1,
\]
entre les espaces de morphismes
\[
\Hom_{\gr \sf C({\mathbb A},{\mathbb A})} \big(\ca \tpo i,\ca\big) \quad \mbox{et} \quad
\Hom_{\gr \sf C({\mathbb A},{\mathbb A})}\big((S\ca)\tpo i,S\ca\big)
\]
\[
\Big(\mbox{resp.} \quad
\Hom_{\gr \sf C({\mathbb A},{\mathbb A})}\big(\ca \tpo i,\lrus{\dot f}{\cb}{\dot f}\big) \quad \mbox{et} \quad
\Hom_{\gr \sf C({\mathbb A},{\mathbb A})}\big((S\ca)\tpo i, \lrus{\dot f}{\cb}{\dot f}\big) \Big)
\]
sont dÈfinies par les relations
\[
\si \circ b_i = -   m_i \circ \si \tpo i \quad \Big(\mbox{resp.} \quad
\si \circ F_{i} = (-1)^{|F_{i}|}f_{i} \circ \si \tpo{j} \Big),
\]
o˘ $F_i$ est un morphisme graduÈ de degrÈ $|F_{i}|$ et $\si = s^{-1}$.
Un foncteur $f :\ca \ra \cb$ est la donnÈe d'une application $\dot f : \Obj \ca \ra \Obj \cb$
et d'un morphisme diffÈrentiel graduÈ
\[
F  = Bf : B\ca \ra B\lrus{\dot f}{\cb}{\dot f}
\]
dans la catÈgorie $\cocog ({\mathbb A},{\mathbb A})$ des cogËbres cocomplËtes diffÈrentielles graduÈes de $\sf C({\mathbb A},{\mathbb A})$.\\

%
%
%
%
%
%
%
%

\begin{definition}{\em
Soit ${\mathbb A}$ un ensemble et $\ca$ une $\ai$-catÈgorie sur ${\mathbb A}$.
Un {\em $\ca$-polydule} \index{polydule}
 est un $\ca$-polydule  dans $\gr \sf C({\mathbb A})$ (voir \ref{definition_ai-module}).
Il est donnÈ par un ${\mathbb A}$-module ‡ droite
\[
M : {\mathbb A}^{op} \ra \gr \sf C
\]
muni de morphismes graduÈs de ${\mathbb A}$-modules ‡ droite
\[
m_i : M \tso \ca ^{\tso i-1} \ra M, \quad i \geq 1,
\]
de degrÈ $2-i$, tels qu'une Èquation $(*'_n)$, $n \geq 1$, de la mÍme
forme que l'Èquation $(*_n)$, $n \geq 1$, de la dÈfinition
\ref{definition_ai-algebre} est vÈrifiÈe.
}\end{definition}

\begin{remarque} \label{remarque1_polydules}{\em
Soit $V$ un objet de $\sf C({\mathbb A})$. Le ${\mathbb A}$-module $V \tso \ca$ muni des morphismes
\[
\Id_V \ts m_i : \big(V \tso \ca\big) \tso \ca \tpo{i-1} \ra V \tso \ca, \quad i\geq 1,
\]
est un $\ca$-polydule.
En particulier, si $A$ est un objet de $\ca$, le ${\mathbb A}$-module
\[
\ca(\?,A) = e(\?,A) \tso \ca
\]
est un $\ca$-polydules dans $\sf C({\mathbb A})$. On le notera $A \pw$.
}\end{remarque}

%
%
%
%
%
%
%
%

\begin{definition}{\em  
Soit ${\mathbb A}$ et ${\mathbb B}$ deux ensembles et
$\ca$ et $\cb$ deux $\ai$-catÈgories sur ${\mathbb A}$ et ${\mathbb B}$.
Un $\ca$-$\cb$-{\em bipolydule} \index{bipolydule} est un $\ca$-$\cb$-bipolydule
dans $\gr \sf C({\mathbb A},{\mathbb B})$ (voir \ref{definition_bipolydules}).
}\end{definition}

%
%
%
%
%
%
%
%
%
%
%
%
%
%
%

\section{CatÈgories diffÈrentielles graduÈes des poly\-dules}

 \label{section_categorie_dg_polydules}

{\noindent \bf La catÈgorie $\cc_\infty \Ba \ca$}\\

Soit ${\mathbb A}$ un ensemble et $\ca$ une $\ai$-catÈgorie sur ${\mathbb A}$.
La {\em catÈgorie} $\cc_\infty \Ba \ca$
 a pour objets ceux de $\coModcu \Ba \ca$. Si $N$ et $N'$
sont deux objets de $\coModcu \Ba \ca$,
l'espace de morphismes 
\[
\Hom_{\cc_\infty \Ba \ca} (N,N')
\]
est l'espace des morphismes unitaires graduÈs de comodules $N \ra N'$ muni
de la diffÈrentielle
\[
\delta :  F \mapsto b^{N'} \circ F - (-1)^{|F|} F \circ b^{N},
\]
o˘ $F$ est de degrÈ $|F|.$ C'est une catÈgorie diffÈrentielle graduÈe.
Remarquons que la catÈgorie $\coModcu \Ba \ca$ est isomorphe ‡
la catÈgorie $Z^0\cc_\infty \Ba \ca$, i.~e.~la catÈgorie
dont les objets sont ceux de $\cc_\infty \Ba \ca$ et dont les morphismes
sont les zÈro-cycles des complexes de morphismes de $\cc_\infty \Ba \ca$.\\

{\noindent \bf La catÈgorie $\cn_\infty \ca$}\\

La {\em catÈgorie}  \indexnotation{cninftyA}
$\cn_\infty \ca$ est la catÈgorie diffÈrentielle graduÈe 
dont les objets sont les $\ca$-polydules  et dont les espaces
de morphismes sont dÈfinis par
\[
\Hom_{\sz \cn_\infty \ca}(M,M') = \Hom_{\cc_\infty \Ba \ca}(BM ,BM'), \quad M,M' \in \cn_\infty \ca.
\]
Un morphisme $f : M \ra M'$ de degrÈ $n$ est donc donnÈ par une suite
de morphismes graduÈs de ${\mathbb A}$-module
\[
f_i : M \tso \ca \tpo{i-1} \ra M'
\]
de degrÈ $1-i+n.$ Les $\ai$-morphismes $f : M \ra M'$
sont les zÈro-cycles de $\Hom_{\cc_\infty \ca}(M,M').$
(La lettre $\cn$ se rapporte au ``$\cn$on'' dans ``$\ca$-polydules non unitaires''.)

\begin{remarque}{\em
Soit $\cb$ une $\ai$-catÈgorie et $X$ un $\cb$-$\ca$-bipolydules.
Nous avons un isomorphisme de complexes
\[
\Homi_{\ca}(X,M) = \Hom_{\sz \cn_\infty \ca}(X_\ca,M), \quad M \in \cn_\infty \ca,
\]
o˘ $\Homi_{\ca}(X,M)$ est dÈfini en (\ref{section_foncteurs_standard}).
}\end{remarque}

{\noindent \bf La catÈgorie $\cc_\infty \ca$}\\

Supposons dÈsormais que $\ca$ est strictement unitaire.
Si $M$ et $M'$ sont deux $\ca$-polydules strictement unitaires, un morphisme $f : M \ra M'$
de degrÈ $n$ est {\em strictement unitaire} s'il vÈrifie les Èquations
\[
f_i ( \Id \tp \alpha \ts \unite \ts \Id \tp \beta) = 0, \quad i\geq 2.
\]
Nous notons $\big(\cn_\infty \ca\big)_u$ \indexnotation{cninftyAu}
la {\em sous-catÈgorie} pleine de $\cn_\infty \ca$ formÈe des $\ca$-polydules strictement unitaires et
$\cc_\infty \ca$ \indexnotation{ccinftyA}
la {\em sous-catÈgorie} non pleine de $\cn_\infty \ca$ formÈe des $\ca$-polydules strictement unitaires dont
les morphismes sont les morphismes strictement unitaires.
Remarquons que si $\ca$ est augmentÈe, nous avons un isomorphisme de catÈgories
\[
\cc_\infty \ca \arr{\sim} \cn_\infty \b \ca.
\]

\begin{remarque}{\em \label{remarque_categorie_differentielle_graduee_polydules}
La catÈgorie $H^0\cc_\infty \ca$ est clairement isomorphe ‡ $\ch_\infty \ca$
(voir la dÈfinition \ref{definition_categorie_derivee_generale}). Elle
est Èquivalente ‡ la catÈgorie $\cd_\infty \ca$ par le corollaire \ref{corollaire3_categorie_derivee}.
}\end{remarque}

\begin{proposition}\label{proposition_strictication_unitaire_polydules_graduee}
L'inclusion
\[
\cc_\infty \ca \ra \cn_\infty \ca
\]
induit un quasi-isomorphisme dans les espaces de morphismes.
\end{proposition}

\dem La dÈmonstration est la mÍme que celle de la proposition
(\ref{proposition_strictication_unitaire_polydules}). Au lieu de considÈrer
uniquement les $\ai$-morphismes, i.~e.~les morphismes de $\cn_\infty \ca$ qui sont des cycles de degrÈ zÈro
et les homotopies entre $\ai$-morphismes,
nous considÈrons les morphismes de degrÈ quelconque qui sont des cycles. \findem\\

{\noindent \bf La catÈgorie $\cc_\infty (\ca,\cb)$} \label{section_categorie_dg_bipolydules}\\

Soit $\mathbb A$ et $\mathbb B$ deux ensembles et $\ca$ et $\cb$ des $\ai$-catÈgories sur $\mathbb A$
et $\mathbb B$.
La catÈgorie $\cn_\infty (\ca,\cb)$ est construite de maniËre strictement analogue ‡ $\cn_\infty \ca$.
Soit la {\em catÈgorie} diffÈrentielle graduÈe $\cc_\infty (\Ba \ca,\Ba \cb)$ dont les objets
sont ceux de $\coModcu (\Ba \ca,\Ba \cb)$. 
La {\em catÈgorie} $\cn_\infty (\ca,\cb)$ \indexnotation{cninftyAA'}
est la catÈgorie diffÈrentielle graduÈe 
qui a les mÍmes objets que que $\aiMod (\ca,\cb)$ et dont les espaces
de morphismes sont dÈfinis par le sous-espaces
\[
\Hom_{\sz \cn_\infty (\ca,\cb)}(M,M') = \Hom_{\cc_\infty (\Ba \ca,\Ba \cb)}(BM ,BM'),
\quad M,M' \in \cc_\infty(\ca,\cb).
\]

Si $\ca$ et $\cb$ sont strictement unitaires, on dÈfinit les catÈgories \indexnotation{cninftyAA'u}
$\big(\cn_\infty (\ca,\cb)\big)_u$ et $\cc_\infty(\ca,\cb)$ \indexnotation{ccinftyAA'}
de maniËre analogue aux catÈgories
$\big(\cn_\infty \ca\big)_u$ et $\cc_\infty \ca.$
La catÈgorie $\aiModu (\ca,\cb)$ est isomorphe ‡ $Z^0\cc_\infty (\ca,\cb)$.

\begin{proposition} \label{proposition_strictication_unitaire_bipolydules_graduee}
L'inclusion
\[
\cc_\infty (\ca,\cb) \ra \cn_\infty (\ca,\cb)
\]
induit un quasi-isomorphisme dans les espaces de morphismes.
\findem
\end{proposition}

%
%
%
%
%
%
%
%
%
%
%
%
%

\section{Lemme clef}
\label{section_lemme_clef}

Le lemme ci-dessous sera utile pour la construction
de l'$\ai$-foncteur de Yoneda (voir \ref{definition_ai-foncteur_de_Yoneda}).\\

Soit ${\mathbb A}$ et ${\mathbb B}$ deux ensembles, $M$ un objet graduÈ de $\sf C({\mathbb A},{\mathbb B})$ et
$\ca$ et $\cb$ deux $\ai$-catÈgories sur ${\mathbb A}$ et ${\mathbb B}$. 
Soit une famille de morphismes graduÈs de ${\mathbb A}$-${\mathbb B}$-bimodules
\[
m_{i,j} : \ca \tpo{_{{\mathbb A}} i}  \tso_{{\mathbb A}} M \tso_{{\mathbb B}} \cb \tpo{_{{\mathbb B}} j} \ra M, \quad
i,j\geq 0,
\]
de degrÈ $1-i-j.$  
Munissons les cogËbres tensorielles co-augmentÈes $\ct S\cb$ et $\ct S\ca$ des diffÈrentielles $b^\ca$
et $b^\cb$ des constructions bar co-augmentÈes. Les morphismes
\[
b_{0,j} : SM \tso_{{\mathbb B}} (S\cb) \tpo{_{{\mathbb B}} j} \ra SM,\quad j\geq 0,
\]
donnÈs par les bijections $m_{0,j} \lra b_{0,j}$ de la section \ref{section_construction_bar_bipolydules},
se relËvent (voir \ref{lemme_comodules_colibres}) en une unique codÈrivation de comodules
graduÈs dans $\sf C({\mathbb A},{\mathbb B})$
\[
D : SM \tso_{{\mathbb B}} \ct (S\cb)  \ra SM \tso_{{\mathbb B}} \ct (S\cb).
\]
Notons $\mathsf{End} = \mathsf{End}\Big(\big(SM \tso_{{\mathbb B}} \ct (S\cb)\Big)$
l'algËbre des endomorphismes graduÈs de $\Ba \cb$-comodules dans la catÈgorie $C({\mathbb A},{\mathbb B})$.
Remarquons que cet un objet de la catÈgorie $\sf C({\mathbb A},{\mathbb A})$ est aussi dÈfini par
\[
\mathsf{End}(A,A') = \Homi_\cb\big(M ,M(\?,A')\big)(A), \quad A,A' \in \mathbb A,
\]
o˘ $\Homi$ est le foncteur dÈfini en (\ref{section_foncteurs_standard}).
Nous munissons $\mathsf{End}$
des trois morphismes
\[
\begin{array}{rclrcl}
 m_0 : e_{\mathbb A} & \ra & \mathsf{End},& 1 & \mapsto & -D^2\\
 m_1 : \mathsf{End}& \ra & \mathsf{End}, & f &\mapsto & D \circ f - (-1)^{r}f \circ D\\
 m_2 :  \mathsf{End} \tso_{\mathbb A} \mathsf{End} & \ra & \mathsf{End}, &
f \tso g & \mapsto & f \circ g,
\end{array}
\]
o˘ $f$ est un morphisme de degrÈ $r$. Ils vÈrifient les Èquations
\[
m_{1}\cdot m_{0} =0 , \quad m_{2}(m_{0}\tso \Id + \Id \tso m_{0}) + m_{1}^2
= 0
\]
\[
m_{2}(m_{1}\tso \Id +
\Id \tso m_{1}) - m_{1} m_{2} = 0 \quad \mbox{et} \quad m_2(\Id \tso m_2 - m_2 \tso \Id) = 0.
\]
Une algËbre diffÈrentielle graduÈe $(A,d, \mu)$ vÈrifie clairement ces Èquations pour $m_0 = 0$,
$m_1 = d$ et $m_2 = \mu.$ RÈciproquement, si $M$ est un objet graduÈ muni de morphismes
$m_0$, $m_1$ et $m_2$ vÈrifiant ces Èquations,  $(M,m_1,m_2)$ est
une algËbre diffÈrentielle graduÈ si $m_0$ est nul.

Soit les morphismes graduÈs de ${\mathbb A}$-${\mathbb A}$-bimodules
\[
f_{i} : \ca \tpo i \ra \mathsf{End}, \quad i\geq 1,
\]
de degrÈ $2-i$, dÈfinis par l'Èquation
\[
F_{i}(\phi) = s(\Phi) \in 
S\mathsf{End},
\]
o˘ $F_i$ est donnÈ par la bijection $f_i \lra F_i$, o˘ $\phi$
est un ÈlÈment de $(S\ca)\tpo i$ de degrÈ $|\phi|$
et o˘
le morphisme $\Phi$ est l'unique morphisme (voir \ref{lemme_comodules_colibres})
tel que la composition $p_1 \circ \Phi$ a pour composantes les morphismes
\[
\xymatrix{
SM \tso (S\cb)\tpo j \ar[rr]^(.4){(-1)^{|\phi|} \phi \tso \Id} & &
(S\ca)\tpo i \tso SM \tso (S\cb)\tpo j \ar[r]^(.75){b_{i,j}} &
SM,} \quad j \geq 0.
\]

\begin{lemme} \label{lemme_clef} \index{lemme clef}
Les ÈnoncÈs suivants sont Èquivalents.
\english \begin{itemize}
\item[a.] Le triplet $(\mathsf{End},m_1,m_2)$
est une algËbre diffÈrentielle graduÈe et
les morphismes $f_i$, $i \geq 1$, dÈfinissent un $\ai$-morphisme
\[
f : \ca \ra \mathsf{End},
\]
o˘ $\mathsf{End}$ est munie de l'$\ai$-structure de la remarque \ref{remarque_ai-structure_alg}.

\item[b.] Les morphismes $m_{i,j}$, $i,j \geq 0$, dÈfinissent une structure de
$\ca$-$\cb$-bipoly\-dule sur $M.$

\end{itemize} \francais
\end{lemme}

\dem
Supposons que l'ÈnoncÈ {\it a} est vrai.
Nous allons montrer qu'il est Èquivalent aux Èquations
\[
\sum_{k+\bullet+m=n} b_\bullet (\Id \tpo k \tso b_\bullet \tso \Id \tpo m) = 0, \quad n \geq 0,
\]
o˘ les symboles $b_\bullet$ doivent Ítre interprÈtÈs convenablement. Ces Èquations sont Èquivalentes
aux Èquations $(*''_n)$, $n+1+n'\geq 1$, de la dÈfinition \ref{definition_bipolydules}.

Comme $(\mathsf{End},m_1,m_2)$ est une algËbre diffÈrentielle graduÈe,
le morphisme $m_0$ est nul. Ceci veut dire que $D$ est une diffÈrentielle
de comodules. L'Èquation $D^2 = 0$ est Èquivalente aux Èquations
\[
\sum_{1+j+k=n} b_{0,l}(b_{0,j}\tso \Id \tpo k)
+ \sum_{k+j+m=n} b_{0,l}(\Id \tpo k \tso b^{\cb}_{j}\tso \Id \tpo m) =
0, \quad n\geq 0.
\]

En vertu de la section \ref{construction_bar_cobar}, le fait que $f$ est un $\ai$-morphisme
se traduit par le fait que la suite des morphismes $F_i$, $i\geq 1$,
dÈfinit un morphisme de cogËbres diffÈrentielles
graduÈes
\[
F : \Ba \ca \ra \Ba \mathsf{End}.
\]
Ceci Èquivaut aux Èquations $(**_n)$, $n \geq 1$:
\[
\sum_{i+j+k=n} F_l(\Id \tpo i \tso b^{\ca}_j \tso \Id \tpo k) -
b^{\mathsf{End}}_1(F_n) - \sum_{i + j =n}b^{\mathsf{End}}_2(F_{i}\tso F_{j}) = 0.
\]

On rappelle que la dÈfinition des bijections $m^{\mathsf{End}}_i \lra b^{\mathsf{End}}_i$, $i\geq 2$, implique
que
\[
b^{\mathsf{End}}_1 \circ s = - s \circ m_1^{\mathsf{End}}\quad \mbox{et} \quad
b_2^{\mathsf{End}} \circ s \tpo 2 = s \circ m_2^{\mathsf{End}}.
\]

Soit $m \tso y$ un ÈlÈment de de $SM \tso (S\cb)\tpo n$.
Calculons l'image de $m \tso y$ par $b^{\mathsf{End}}_2(F_{i}\tso F_{j})(\phi)$ o˘
$\phi = \phi_j \tso \phi_i$
:\\[.3cm]
$
b^{\mathsf{End}}_2(F_{i}\tso F_{j})(\phi) (m \tso y)  = \ \ b^{\mathsf{End}}_2
(s\Phi_i \tso s\Phi_j)(m \tso y)\\[.08cm]
\hspace{2cm} = (-1)^{|\Phi_i|}b^{\mathsf{End}}_2(s\tso s)(\Phi_i \tso \Phi_j)(m \tso y) 
\\[.08cm]     =
(-1)^{|\Phi_i|}sm_2^{\mathsf{End}}(\Phi_i \tso \Phi_j)(m \tso y)\\[.08cm]   =
\sum_{k+l=n}(-1)^{|\Phi_i|+|\phi_i|+|\phi_j|}
sb_{i,l}(\phi_i \tso b_{j,k}(\phi_j \tso \Id_{SM} \tso \Id \tp{k})\tso \Id
\tpo{l})(m\tso y)\\[.08cm]
     =
\sum_{k+l=n}(-1)^{|\phi_j|+1}sb_{i,l}(\phi_i \tso b_{j,k}(\phi_j \tso
\Id_{SM} \tso \Id \tpo{k})\tso \Id
\tpo{l})(m \tso y)\\[.08cm]
=  \sum_{k+l=n} (-1)^{|\phi|+1}sb_{i,l}(\Id_{S\ca}\tpo i \tso
b_{j,j}(\Id_{S\ca}\tpo j \tso \Id_{SM} \tso
\Id \tpo{k})\tso \Id
\tp{l})(\phi \tso m \tso y),
$\\[.3cm]
puis $ b^{\mathsf{End}}_1(F_n)(\phi)(m \tso y)$ :\\[.3cm]
$
     b^{\mathsf{End}}_1(F_n)(\phi)(y) =   b^{\mathsf{End}}_1(s\Phi)(m \tso y)\\[.08cm]  =
- sm_1^{\mathsf{End}}(\Phi)(m \tso y)\\[.08cm]   =
- s(b \cdot\Phi -(-1)^{|\Phi|} \Phi \cdot b )(m \tso y)\\[.08cm]  =
- s\Big[ \sum_{k+l=n} (-1)^{|\phi|} b_{0,l}(b_{i+j,k}(\phi \tso \Id_{SM} \tso \Id \tpo{k}) \tso
\Id \tpo{l}) +\\
{\ \ \ }\sum_{u+v+l = n} -  (-1)^{|\phi|}
b_{i+j,u+1+l}(\phi \tso \Id_{SM} \tso \Id  \tpo{u} \tso b_{v}^\cb \tso \Id
\tpo{l}) +\\
{\ \ \ } -  (-1)^{|\phi|}
b_{i+j,u+1+l}(\phi \tso b_{0,n}(\Id_{SM} \tso \Id  \tpo{n})\Big] (m \tso y)
\\[.08cm]  =
(-1)^{|\phi|+1}
s\Big[\sum_{k+l=n} b_{0,l}(b_{i+j,k}(\Id_{S\ca}\tpo{i+j} \tso \Id_{SM} \tso
\Id \tpo{k}) \tso \Id \tpo{l}) +\\
{\ \ \ }\sum_{u+v+l = n}
b_{i+j,u+1+l}(\Id_{S\ca}\tpo{i+j} \tso \Id_{SM} \tso \Id  \tpo{u} \tso
b_{v}^\cb \tso \Id  \tpo{l}) +\\
{\ \ \ }b_{i+j,u+1+l}(\Id_{S\ca}\tpo{i+j} \tso b_{0,n}(\Id_{SM} \tso \Id
     \tpo{n})\Big] (\phi \tso m \tso y)
$\\[.3cm]
et enfin $F_l(\Id \tpo i \tso b^{\ca}_j \tso \Id \tpo k)(\phi)(m \tso y)$
:\\[.3cm]
$
F_l(\Id \tpo i \tso b^{\ca}_j \tso \Id \tpo k)(\phi) (m \tso y) = \\[.08cm]
\sum_{u+v+t = i+j}(-1)^{|\phi|+1}b_{u+1+t,n}
(\Id_{S\ca}\tpo u \tso b^\ca_{v}\tso \Id_{S\ca} \tpo t \tso \Id_{SM} \tso \Id^n)
(\phi \tso m \tso y).
$\\[.3cm]
Les Èquations $(**_n)$, $n\geq 1$, et le fait que la codÈrivation $D$ est une diffÈrentielle
sont donc Èquivalents aux Èquations
\[
\sum b_{\bullet}(\Id \tpo u \tso b_{\bullet} \tso \Id \tpo v) = 0
\]
o˘ les $b_{\bullet}$ et les $\Id$ doivent Ítre interprÈtÈs convenablement.
\findem \\

Soit $M$ un $\ca$-$\cb$-bipolydule.
Soit $A$ un objet de $\ca$.
On munit le ${\mathbb A}$-module $M(\?,A)$ de la structure de $\cb$-polydule donnÈe par les
morphismes $m_j$, $j \geq 1$, de ${\mathbb B}$-modules
\[
m_{0,j-1}(\?,A) : \Big(M \tso \cb \tp {j-1}\Big) (\?,A) \ra M(\?,A), \quad j \geq 1.
\]

\begin{corollaire} \label{corollaire_lemme_clef}
L'application 
\[
\dot{\theta}_M : {\mathbb A} \ra \Obj \cn_\infty \cb, \quad A \mapsto M(\?,A)
\]
se complËte de maniËre canonique en un $\ai$-foncteur
\[
\theta_M : \ca \ra \cn_\infty \cb.
\]
\end{corollaire}

\dem 
Le ${\mathbb A}$-${\mathbb A}$-bimodule 
\[
\Hom_{\cn_\infty \cb}(\dot{\theta_M}\?,\dot{\theta_M}\?)
\]
est par dÈfinition
l'algËbre des endomorphismes
\[
\mathsf{End}_{\cn_\infty \Ba \cb}\big((?,\?) \ts \ct S\cb\big),
\]
c'est-‡-dire, l'algËbre $\mathsf{End}$ du lemme clef.
Le foncteur canoniquement associÈ ‡ $M$ est donnÈ par les morphismes
\[
f_i : \ca \tpo i \ra \Hom_{\cn_\infty \cb }(\dot{\theta_M}\?,\dot{\theta_M}\?), \quad i\geq 1,
\]
du lemme \ref{lemme_clef}. Ils dÈfinissent un $\ai$-foncteur car
\[
f : \ca \ra \mathsf{End}
\] 
est un $\ai$-morphisme.
\findem

\francais
\begin{corollaire} \label{corollaire2_lemme_clef}
L'application $M \mapsto \theta_M$ de la classe des $\ca$-$\cb$-bipolydules sur
la classe des $\ai$-foncteurs $\ca \ra \cn_\infty \cb$ est une bijection.
Son application inverse associe ‡ un
$\ai$-foncteur
\[
(\dot g,g) : \ca \ra \cn_\infty \cb
\]
le ${\mathbb A}$-${\mathbb B}$-bimodule
\[
M(A,B) = (\dot g(A))(B)
\]
muni des multiplications $m_{i,j}$, $i,j \geq 0$, donnÈes par
\[
m_{i,{j-1}} = (g_i)_j.
\]
\findem
\end{corollaire}
\francais

{\noindent \bf Le cas strictement unitaire}\\

Supposons dÈsormais que $\ca$ et $\cb$ sont des $\ai$-catÈgories strictement unitaires.

\begin{remarque}{\em
Soit $M$ un $\ca$-$\cb$-bipolydule.
L'$\ai$-morphisme
\[
f : \ca \ra \mathsf{End}
\]
du lemme clef (\ref{lemme_clef})
est strictement unitaire si et seulement si les compositions
\[
m^M_{i,j}(\Id \tpo \alpha \tso \unite \tso \Id \tpo \beta \ts \Id_M \ts \Id \tp j), \quad i,j \geq 0,
\]
sont nulles pour $(i,j) \neq (1,0)$ et si elle est l'identitÈ si $(i,j) = (1,0).$
}\end{remarque}

\begin{remarque}{\em \label{remarque2_lemme_clef}
Si $M$ est un $\ca$-$\cb$-bipolydule strictement
unitaire, le $\cb$-polydules $M(\?,A)$, $A \in \mathbb A$, est strictement unitaire et l'$\ai$-foncteur
\[
\ca \ra \cn_\infty \cb
\]
du corollaire (\ref{corollaire_lemme_clef}) se factorise par
un foncteur
\[ 
\ca \ra \cc_\infty \cb.
\]
}\end{remarque}

\begin{remarque}{\em \label{remarque3_lemme_clef}
La bijection $M \mapsto \theta_M$ du corollaire (\ref{corollaire2_lemme_clef})
se restreint en une bijection
de la classe des $\ca$-$\cb$-bipolydules strictement unitaires sur
la classe des $\ai$-foncteurs strictement unitaires $\ca \ra \cc_\infty \cb$.
}\end{remarque}

%% file: Torsion.tex
\noindent Dans les chapitres \ref{chapitre_objets_tordus} et \ref{chapitre_ai_func},
nous allons construire des $\ai$-catÈgories dont les compositions sont construites par un processus de
torsion que nous dÈcrivons dans ce chapitre. 

En thÈorie des
dÈformations des algËbres de Lie diffÈrentielles graduÈes (ou algËbres associatives
diffÈrentielles graduÈes), la technique de la torsion est bien connue
(pour un panorama, voir par exemple \cite{Huebschmann99}).
La version $\ai$ (et $\mathrm{L}_\infty$) a ÈtÈ introduite dans
\cite[Chap.~4]{Fukaya01} (voir aussi \cite{Fukaya01b}). Notre dÈmonstration
du fait que les compositions tordues dÈfinissent bien une structure $\ai$ est diffÈrente.
La torsion d'une $\ai$-algËbre $A$ par une solution ‡ l'Èquation de Maurer-Cartan gÈnÈralisÈe modifie
non seulement la diffÈrentielle
$m_1$ mais aussi toutes les multiplications supÈrieures.

Ce chapitre est divisÈ en deux sections. Nous traitons d'abord le cas simple o˘ la torsion
est tensoriellement nilpotente, puis le cas o˘ les $\ai$-structures sont topologiques. 
Nous montrons que si $f : \ca \ra \cb$ est un $\ai$-foncteur qui induit des 
quasi-isomorphismes dans les espaces
de morphismes, sa torsion induit aussi des quasi-isomorphismes dans les espaces
de morphismes (\ref{proposition_torsion_ai-foncteurs_equiv_faible_I}).

\section{Le cas tensoriellement nilpotent}

\subsection{ElÈments tordants}

Soit ${\mathbb A}$ un ensemble et $\ca$ une $\ai$-catÈgorie sur ${\mathbb A}$.
Munissons l'ÈlÈment neutre $e = e_{\mathbb A}$ pour le produit tensoriel $\tso_{\mathbb A}$
de la structure de cogËbre donnÈe par la contrainte d'unitaritÈ de la catÈgorie monoÔdale de base $\sf C({\mathbb A},{\mathbb A})$
(voir \ref{section_categories_de_base_aicat} et \ref{categorie_de_base})
\[
e \arr{\sim} e \tso e.
\]
ConsidÈrons $e$ comme une cogËbre diffÈrentielle graduÈe concentrÈe en degrÈ~$0.$

\begin{definition}{\em 
Un {\em ÈlÈment tordant (tensoriellement nilpotent)} \index{element tordant@{ÈlÈment tordant}}
est un morphisme graduÈ $x : e \ra \ca$ de degrÈ $+1$ tel que
\english \begin{itemize}
\item[(1)] la composÈe $s \circ x$ se relËve en un morphisme de cogËbres
\[
X : e \ra B\ca,
\]
\item[(2)] et ce morphisme $X$ est compatible
aux diffÈrentielles.
\end{itemize} \francais
}\end{definition}

\begin{remarque}{\em \label{remarque_nilpotence_tensorielle}
Notons $p_1$ la projection $ B\ca \ra S\ca$.
La composition avec la projection donne une bijection
\[
\Hom_{\sz \cog} (e , B\ca) \arr{\sim} \Hom_{{nil}} (e,S\ca),
\]
o˘ $\Hom_{{nil}} (e,S\ca)$ est l'ensemble des morphismes
graduÈs $\phi : e \ra S\ca$ de degrÈ $0$
tels que, pour tout $A \in {\mathbb A}$, il existe un $N$ tel que
$\phi \tp n \Delta^{n-1} (\id_A)= 0$ pour $n\geq N.$
Nous en dÈduisons qu'un morphisme graduÈ $x : e \ra \ca$ de degrÈ $+1$ est
un ÈlÈment tordant si et seulement si
\english \begin{itemize}
\item[(1)] il est {\em tensoriellement nilpotent}\index{tensoriellement nilpotent} :
pour tout objet
$A\in {\mathbb A}$, l'ÈlÈment  $x(\id_A)  \in \ca(A,A)$ de degrÈ 1
est tel que
$x(\id_A)\tpo {n}$ est nul
pour un certain $n>0$,
\item[(2)] il vÈrifie l'Èquation de Maurer-Cartan \index{Maurer-Cartan (Èquation)}
\[
\sum_{i=1}^\infty m_i\big(x(\id_A)\tso \hdots \tso x(\id_A)\big) = 0, \quad A \in {\mathbb A}.
\]
(La somme est finie gr‚ce ‡ la propriÈtÈ de nilpotence tensorielle).
\end{itemize} \francais
}\end{remarque}

\subsection{Torsion des $\ai$-catÈgories} \label{section_torsion_ai-categories_I}

Soit ${\mathbb A}$ un ensemble et $\ca$ une $\ai$-catÈgorie sur ${\mathbb A}$. 
Soit $x$ un ÈlÈment tensoriellement nilpotent de $\ca.$
Soit 
\[
g : \ct S\ca = e \oplus \ctr S\ca \ra S\ca 
\]
le morphisme de composantes $[sx , p_1]$, o˘ $p_1$ est
la projection $\ctr S\ca \ra S\ca$.
Soit le morphisme de ${\mathbb A}$-${\mathbb A}$-bimodules
\[
\phi_x : \ct S\ca \ra \ct S\ca = \bigoplus_{i\geq 0}(S\ca)\tpo i
\]
dont la composÈe avec la projection sur $(S\ca)\tpo i$ est le morphisme
\[
(g \tpo i) \circ \Delta^{(i)} \quad \mbox{si} \quad i \geq 1, \quad \Id_e \quad \mbox{sinon.}
\]
Il est clairement
un morphisme de cogËbres co-unitaire et il est bien dÈfini car
sa restriction au sous-objet $(S\ca)\tpo i \in \sf C({\mathbb A},{\mathbb A})$
est Ègale ‡ la somme (bien dÈfinie par la propriÈtÈ de nilpotence tensorielle)
\[
\sum_{l \geq 0} \sum ((sx) \tpo{l_0} \tso \Id_{S\ca} \tso (sx) \tpo{l_1} \tso \hdots \tso \Id_{S\ca} \tso
(sx) \tpo{l_{i-1}} \tso \Id_{S\ca} \tso (sx) \tpo{l_{i}}),
\]
o˘ $l_0 + \hdots + l_{i} = l$. Remarquons que la composition
\[
\phi_x \circ \coaugmentation : e \ra \ct S\ca = e \oplus \ctr S\ca,
\]
o˘ $\coaugmentation$ est la co-augmentation de
$\ct S\ca$, a pour composantes le morphisme $\Id_e$ et le relËvement $X$
de $s \circ x$.
La matrice de
\[
\phi_x : \bigoplus_{j\geq 0}(S\ca)\tpo j \ra \bigoplus_{i\geq 0}(S\ca)\tpo i
\]
est triangulaire infÈrieure et sa diagonale est celle de l'identitÈ.
Le morphisme $\phi_x$ est donc un automorphisme co-unitaire (non co-augmentÈ)
de la cogËbre graduÈe co-augmentÈe $\ct S\ca$.
La diffÈrentielle de la construction bar $B\ca$ nous donne une diffÈrentielle
\[
b : \ct S\ca \ra \ct S\ca
\]
qui s'annule sur la co-augmentation.
Soit la composÈe
\[
D_x = \phi_x^{-1} \circ b \circ \phi_x : \ct S\ca \ra \ct S\ca.
\]
Supposons que $x$ vÈrifie l'Èquation de Maurer-Cartan.
Le relËvement $X : e \ra \ctr S\ca$ de $s\circ x$ est diffÈrentiel graduÈ.
La composition
\[
b \circ \phi_x \circ \coaugmentation = b \circ \Big[\begin{array}{c}\Id_e \\ X \end{array}\Big] : e \ra \ct S\ca = e \oplus \ctr S\ca
\]
est donc nulle et on a $D_x \circ \coaugmentation = 0.$
Posons $b_x$ le morphisme donnÈ par la flËche verticale de droite du diagramme
de suites exactes
\[
\xymatrix{0 \ar[r] & e \ar[r]^(.3){\coaugmentation } \ar[d]_0 & \ct S\ca \ar[r]
\ar[d]_{D_x} &\ctr S\ca \ar[r] \ar[d]_{{b_x}} & 0 \\
              0 \ar[r] & e \ar[r]^(.3){\coaugmentation } & \ct S\ca \ar[r]
                    &\ctr S\ca \ar[r] & 0.  \\}
\]
Comme $D_x$ est une $(\Id,\Id)$-codÈrivation de $\ct S\ca$,
le morphisme $b_x$ est une $(\Id,\Id)$-codÈrivation de $\ctr S\ca$.
Comme $D_x^2 = 0$, la codÈrivation $b_x$ est une diffÈrentielle de la cogËbre $\ctr S\ca.$
Elle est dÈterminÈe (\ref{cogebres_tensorielles}) par les composantes
\[
(b_x)_i :  (S\ca)\tpo i \ra S\ca
\]
de sa composition avec la projection sur $S\ca$.

\begin{lemme} \label{lemme_torsion_ai-categories_I}
Soit $i \geq 1.$ Le morphisme $(b_x)_i$ est la somme
\[
\sum_{l} \sum b_{l+m}((sx) \tpo{l_0} \tso \Id_{S\ca} \tso (sx) \tpo{l_1} \tso \hdots \tso \Id_{S\ca} \tso
(sx) \tpo{l_{i-1}} \tso \Id_{S\ca} \tso (sx) \tpo{l_{i}}),
\]
o˘ $l_0 + \hdots + l_{i} = l$.
\end{lemme}

\dem Remarquons que $D_x$ restreint au sous-objet $\ctr S\ca$ de $\ct S\ca$ est Ègal ‡ $b_x.$
Nous devons calculer
\[
(p_1 \circ D_x) |_{(S\ca)\tpo i} = (p_1 \circ \phi_{x}^{-1} \circ b \circ \phi_x)|_{(S\ca)\tpo i} \quad i\geq 1.
\]
Comme la matrice des coefficients de
\[
\phi_x : \bigoplus_{j\geq 0}(S\ca)\tpo j \ra \bigoplus_{i\geq 0}(S\ca)\tpo i
\]
est triangulaire infÈrieure et comme sa diagonale est celle
de l'identitÈ, la matrice de $\phi_x^{-1}$ est de la mÍme forme. Ainsi,
les morphisme $p_1 \circ \phi_{x}^{-1} \circ b \circ \phi_x$ et
$p_1 \circ b \circ \phi_x$ restreints ‡ $(S\ca)\tpo i$ sont Ègaux.
Ceci dÈmontre le lemme. \findem\\

\begin{definition}[K.~Fukaya \cite{Fukaya01} (voir aussi \cite{Fukaya01b})]
{\em \label{definition_torsion_ai-categories_I}
L'{\em $\ai$-catÈgorie tordue $\ca_x$} \index{A-infini categorie@{$\ai$-catÈgorie}!tordue} \index{torsion}
\indexnotation{cax1}
sur ${\mathbb A}$ est le ${\mathbb A}$-${\mathbb A}$-bimodule $\ca_x = \ca$
dont la construction bar $B\ca_x$ est la cogËbre tensorielle rÈduite diffÈrentielle graduÈe
\[
(\ctr S\ca, b_x).
\]
Ses compositions
\[
m^x_i : \ca_x \tpo i \ra \ca_x, \quad i \geq 1,
\]
sont donc dÈfinies par la somme
\[
\sum_{l} \sum (-1)^{s} m^{\ca}_{l+i}(x \tpo{l_0} \tso \Id_\ca \tso x \tpo{l_1} \tso \hdots \tso \Id_\ca \tso
x \tpo{l_{i-1}} \tso \Id_\ca \tso x \tpo{l_{i}}),
\]
o˘ l'exposant du signe est $s = \sum_{1 \leq t \leq i} t \times l_t$
(Cette somme infinie dÈfinit bien un morphisme gr‚ce ‡ la propriÈtÈ de nilpotence tensorielle
de $x$).
}\end{definition}

\subsection{Torsion des $\ai$-foncteurs}\label{section_torsion_ai-foncteurs_I}

Soit ${\mathbb A}$ et ${\mathbb B}$ deux ensembles, $\ca$ et $\cb$ deux $\ai$-catÈgories sur ${\mathbb A}$ et ${\mathbb B}.$ Soit
\[
(\dot f, f) : \ca \ra \cb
\]
un $\ai$-foncteur et $x$ et $x'$ des ÈlÈments tordants de $\ca$ et $\cb$ vÈrifiant une relation de
compatibilitÈ avec $f$ qui sera prÈcisÈe plus bas.
Cette relation dit approximativement que l'image de $x$ par $f$
est $x'.$ Le but de cette section est de construire un $\ai$-foncteur tordu
\[
 \ca_x \ra \cb_{x'}.
\]

Munissons le ${\mathbb A}$-${\mathbb A}$-bimodule $\lrus{\dot f}{\cb}{\dot f}$ de la
structure de $\ai$-catÈgorie sur ${\mathbb A}$ du lemme \ref{lemme_ai-categorie_fBf}. Nous notons $\cb'$ cette
$\ai$-catÈgorie sur ${\mathbb A}$. L'ÈlÈment tordant
\[
x' : e_{\mathbb B} \ra \cb
\]
donne un ÈlÈment tordant de $\cb'$
\[
e_{\mathbb A} \ra \cb', \quad \id_A \mapsto x'(\dot fA).
\]
Nous le notons aussi $x'.$ Soit
\[
F : \Ba\ca \ra \Ba\cb'
\]
la construction bar co-augmentÈe de l'$\ai$-morphisme  $f : \ca \ra \cb'.$
Nous allons construire l'$\ai$-foncteur tordu de faÁon ‡ ce que
le morphisme
\[
G = \phi_{x'}^{-1} \circ F \circ \phi_x : \ct S\ca \arr{} \ct S\cb'
\]
soit sa construction bar co-augmentÈe.
Remarquons que pour des ÈlÈments tordants $x$ et $x'$ quelconques,
le morphisme $G$ est bien un morphisme diffÈrentiel graduÈ
\[
G : \Ba \ca_x \ra \Ba \cb'_{x'}.
\]
Il n'y a cependant aucune raison pour qu'il soit co-augmentÈ car $\phi_x$ et $\phi_{x'}$
ne le sont pas. Demander qu'il le soit nous donne des relations de compatibilitÈs entre
$x$ et $x'$ :
Supposons que $G$ est augmentÈ. Nous avons l'ÈgalitÈ
\[
\phi_{x'}^{-1} \circ F \circ \phi_{x} \circ \coaugmentation = \coaugmentation,
\]
ou en d'autre terme, nous avons
\[
F \circ \phi_{x} \circ \coaugmentation = \phi_{x'} \circ \coaugmentation.
\]
Comme les compositions $\phi_x \circ \coaugmentation$ et $\phi_{x'} \circ \coaugmentation$ sont
Ègales aux morphismes
\[
\Id_e + X : e \ra e \oplus \ctr  S\ca \quad \mbox{et} \quad \Id_e + X' : e \ra e \oplus \ctr  S\cb'
\]
o˘ $X$ et $X'$ sont les relËvements de $x : e \ra \ca$ et $x' : e \ra \cb',$
la compatibilitÈ entre $x$ et $x'$
affirme que la somme (bien dÈfinie par la propriÈtÈ de
nilpotence tensorielle de $x$)
\[
\sum_{i \geq 1} f_i (x \tpo i) : e \ra \cb'
\]
est Ègale ‡ l'ÈlÈment tordant $x'.$

Comme le morphisme $G$ est co-augmentÈ, il est la co-augmentation d'un
morphisme des cogËbres diffÈrentielles graduÈes rÈduites
\[
F_x : B\ca_x \ra B\cb'_{x'}.
\]

\begin{lemme}\label{lemme_torsion_ai-foncteurs_I}
Soit $i \geq 1.$ Le morphisme $(F_x)_i : (S\ca)\tpo i \ra S\cb'$ est la somme
\[
\sum_{l} \sum F_{l+m}((sx) \tpo{l_0} \tso \Id_{S\ca} \tso (sx) \tpo{l_1} \tso \hdots \tso \Id_{S\ca} \tso
(sx) \tpo{l_{i-1}} \tso \Id_{S\ca} \tso (sx) \tpo{l_{i}}),
\]
o˘ $l_0 + \hdots + l_{i} = l$.
\end{lemme}

\dem Similaire ‡ celle du lemme \ref{lemme_torsion_ai-categories_I}.\\

Remarquons que l'$\ai$-catÈgorie $\cb'_{x'}$ est Ègale ‡ $\lrus{\dot f}{(\cb_{x'})}{\dot f}$.

\begin{definition}{\em \label{definition_torsion_ai-foncteurs_I}
L'{\em $\ai$-foncteur tordu} \index{A-infini foncteur@{$\ai$-foncteur}!tordu} \index{torsion}
\indexnotation{fx1}
\[
(\dot f,f^x) : \ca_x \ra \cb_{x'}
\]
est le foncteur dont la construction bar est $F_x.$
}\end{definition}

Il est donc dÈfini par des morphismes
\[
f^x_i : \ca_x \tpo i \ra \lrus{\dot f}{(\cb_{x'})}{\dot f}, \quad i\geq 1,
\]
dÈfinis par les sommes
\[
\sum_{l} \sum (-1)^{s} f^{\ca}_{l+i}(x \tpo{l_0} \tso \Id_\ca \tso x \tpo{l_1} \tso \hdots \tso \Id_\ca \tso
x \tpo{l_{i-1}} \tso \Id_\ca \tso x \tpo{l_{i}}),
\]
o˘ l'exposant du signe est $s = \sum_{1 \leq t \leq i} t \times l_t$.\\

{\noindent \bf Torsion et Èquivalences faibles}

\begin{lemme} 
Soit $\ca$ une $\ai$-catÈgorie et $x$ un ÈlÈment tordant tensoriellement nilpotent. 
Soit $\ca$ une $\ai$-catÈgorie faiblement Èquivalente ‡ zÈro,
i.~e.~le morphisme dans $\mathsf C(\mathbb A,\mathbb A)$
\[
\ca \ra 0
\]
est un $\ai$-quasi-isomorphisme.
La catÈgorie tordue $\ca_x$ est faiblement Èquivalente ‡ zÈro.
\end{lemme}

\dem La catÈgorie ambiante du raisonnement ci-dessous est $\mathsf C(\mathbb A,\mathbb A)$.
Nous rappellons (\ref{remarque_ai-morphismes_ai-foncteurs})
qu'un $\ai$-morphisme $f$ entre deux $\ai$-algËbres de $\mathsf C(\mathbb A,\mathbb A)$
est un $\ai$-foncteur dont l'application sous-jacente $\dot f$ est l'identitÈ de $\mathbb A$.
Soit $K$ le complexe contractile $(\ca,m_1)$.
ConsidÈrons le comme une $\ai$-algËbre (\ref{remarque_ai-structure_complexe}).
Le lemme (\ref{lemme1_cmf_aia}) montre qu'il existe un $\ai$-(iso)morphisme
\[
f : \ca \ra K
\]
tel que $f_1 = \Id_K.$
Munissons l'$\ai$-algËbre $K$ de l'ÈlÈment tordant
\[
x' = \sum_{i \geq 1} f_i (x \tpo i).
\]
Le diagramme commutatif
\[
\xymatrix{
\Ba \ca \ar[r]^{F} \ar[d]_{\phi_x} & \Ba K \ar[d]^{\phi_{x'}}\\
\Ba \ca_x \ar[r]_G & \Ba K_{x'}
}
\]
montre que $G$ est un isomorphisme. En particulier, $\ca_x$ est $\ai$-quasi-isomorphe
‡ $K_{x'}.$ Il suffit donc de montrer que $K_{x'}$ est faiblement Èquivalent ‡ zÈro.
Par construction, les multiplications $m^K_i$, $i\geq 2$, sont nulles. Nous en dÈduisons que
\[
m^{K_{x'}}_1 = m^K_1 \quad \mbox {et} \quad  m^{K_{x'}}_i = 0 \quad i \geq 2.
\]
Ainsi, l'$\ai$-catÈgorie tordue $K_{x'}$ est Ègale ‡ $K$ et elle est faiblement
Èquivalente ‡ zÈro.
\findem.

\begin{proposition}\label{proposition_torsion_ai-foncteurs_equiv_faible_I}
Soit $\ca$ et $\cb$ des $\ai$-catÈgories sur $\mathbb A$ et $\mathbb B.$
Soit
\[
(\dot f,f) : \ca \ra \cb
\]
un $\ai$-foncteur qui
induit un quasi-isomorphisme dans les espaces de morphismes, i.~e.~
les morphismes
\[
f_1 : \ca(A,A') \ra \cb(\dot fA,\dot f A'), \quad A,A' \in \mathbb A,
\]
sont des quasi-isomorphismes. Soit $x$ et $x'$ des ÈlÈments tordant nilpotents
de $\ca$ et $\cb$ compatibles ‡ $f.$
L'$\ai$-foncteur tordu
\[
(\dot f,f^x) : \ca_x \ra \cb_{x'}
\]
induit un quasi-isomorphisme dans les espaces de morphismes.
\end{proposition}

\dem Notons $\cb'$ l'$\ai$-catÈgorie $\lrus{\dot f}{\cb}{\dot f}$ sur $\mathbb A$
(voir \ref{lemme_ai-categorie_fBf}).
L'$\ai$-foncteur $f$ induit un quasi-isomorphisme dans les espaces de morphismes
si et seulement si l'$\ai$-morphisme dans la catÈgorie des $\ai$-algËbres dans $\sf C(\mathbb A,\mathbb A)$
\[
f' : \ca \ra \cb'
\]
induit par $f$ est une Èquivalence faible.
Supposons donc que $f$ est un $\ai$-quasi-isomorphisme dans $\sf C(\mathbb A,\mathbb A)$.
La dÈmonstration de l'axiome de factorisation (CM5) {\it a.}
de la catÈgorie $\aia$ (\ref{theoreme_cmf_aia}) nous
donne une factorisation de $f$ en
\[
\xymatrix{
\ca \ar@{ >->}[r]^(.4){i}
&  \ca \prod C \ar@{->>}[r]^{} & \cb,}
\]
o˘ $\ca \prod C$ est le produit dans $\aia$ de $\ca$ du cÙne $C$
de l'identitÈ du complexe $(\cb,m_1)$ (considÈrÈ comme $\ai$-algËbre),
et $i$ a pour composantes $\Id_\ca$ et $0.$ Il suffit de montrer le rÈsultat dans le cas o˘
$f$ est Ègal ‡ $i$ et dans le cas o˘ il est une fibration triviale.
CommenÁons par la cofibration triviale $i$. Munissons
$\ca \prod C$ de l'ÈlÈment tordant
\[
x'' = \sum_{j \geq 1} i_j (x \tpo j).
\]
Nous avons les ÈgalitÈs
\[
(\ca \prod C)_{x''} = \ca_x \prod C \quad \mbox{et} \quad 
i^{x} = \big[\begin{array}{c}\Id_{\ca_{x}} \\ 0 \end{array}\big] : \ca_x \ra \ca_x \prod C.
\]
Il en rÈsulte que $i^x$ est une Èquivalence faible.
Supposons maintenant que $f$ est une fibration triviale.
Un scindage de $f_1$ dans la catÈgorie des complexes
nous donne un isomorphisme de complexes
\[
j : \ca \ra \cb \oplus K,
\]
o˘ $K$ est un complexe contractile.
Soit $\cb \prod K$ le produit dans $\aia$ de l'$\ai$-algËbre $\cb$ et
du complexe $K$ considÈrÈ comme $\ai$-algËbre.
La projection canonique $p : \cb \prod K \ra \cb$ est une fibration triviale.
La remarque (\ref{remarque_theoreme_cmf_aia}) appliquÈe ‡ l'axiome de relËvement (CM4) {\it a} nous donne
un $\ai$-isomorphisme
\[
\tilde f : \ca \ra \cb \prod K
\]
tel que $\tilde f_1 = j$ et $p \circ \tilde f = f$.
Munissons $\cb \prod K$ de l'ÈlÈment tordant
\[
x'' = \sum_{j \geq 1} \tilde f_j (x \tpo j).
\]
Nous avons l'ÈgalitÈ
\[
(\cb \prod K)_{x''} = \cb_{x'} \prod K
\]
et l'$\ai$-morphisme tordu $p^{x''}$ s'identifie ‡ la projection canonique
\[
\cb_{x'} \prod K \ra \cb_{x'}.
\]
Comme $K$ est contractile, $p^{x''}$ est une Èquivalence faible. L'ÈgalitÈ
 $f^x = p^{x''} \circ \tilde f^x$
montre que $f^x$ est une Èquivalence faible.
\findem

\subsection{Torsion des $\ca$-$\cb$-bipolydules} \label{section_torsion_bipolydules_I}

Les dÈtails sont omis car ils sont similaires aux deux derniËres sections.\\

Soit ${\mathbb A}$ et ${\mathbb B}$ deux ensembles, $\ca$ et $\cb$ deux $\ai$-catÈgories
sur ${\mathbb A}$ et ${\mathbb B}$ et
$M$ un $\ca$-$\cb$-bipolydule. Soit $x$ et $x'$ des
ÈlÈments tordants de $\ca$ et $\cb.$

\begin{definition}{\em \label{definition_torsion_bipolydules_I}
Le {\em $\ca_x$-$\cb_{x'}$-bipolydule} $\lrus{x}{M}{x'}$ \index{bipolydule!tordu} \index{torsion}
\indexnotation{xMx'1}
a pour multiplications les morphismes
\[
m^{x,x'}_{i,j} : \ca_x \tpo i \tso \lrus{x}{M}{x'} \tso \cb_{x'}
\ra \lrus{x}{M}{x'}, \quad i,j\geq 0,
\]
dÈfinis par les sommes
\[
\sum_{l,k\geq 0} \sum (-1)^{s} m_{i+l,j+k}
(x \tpo{l_0} \tso \Id_\ca \hdots \Id_\ca \tso x \tpo{l_i} \tso \Id_M
\tso {x'} \tpo{k_o} \tso \Id_\cb 
\hdots  \Id_\cb \tso {x'} \tpo{k_{j}}),
\]
o˘ l'exposant du signe est
\[
s = \Big(\sum_{1 \leq t \leq i} t \times l_t \Big) + \Big(\sum_{1 \leq t \leq j} (j+ t) \times l_t \Big) 
\]
(Les sommes infinies dÈfinissent bien des morphismes gr‚ce ‡ la propriÈtÈ de nilpotence tensorielle de $x$
et $x'$).
}\end{definition}

\begin{remarque}{\em
La diffÈrentielle $b_{x,x'}$ de la construction bar
du $\ca_x$-$\cb_{x'}$-bipolydule $\lrus{x}{M}{x'}$ est la composÈe
\[
\big(\phi_x^{-1} \tso \Id \tso \phi_{x'}^{-1}\big)
\circ b \circ 
\big(\phi_x \tso \Id \tso \phi_{x'}\big)
\]
o˘
\[
b :  \ct S \ca \tso SM \tso \ct  S \cb \ra \ct S \ca \tso SM \tso \ct S \cb
\]
est la diffÈrentielle de la construction bar du $\ca$-$\cb$-bipolydule $M$.
}\end{remarque}

\begin{remarque}{\em
Soit $f : \ca \ra \cb$ un $\ai$-foncteur. Supposons que les ÈlÈments tordants $x$ et $x'$
sont compatibles ‡ $f$ (voir \ref{section_torsion_ai-foncteurs_I}).
Soit 
\[
y : \cb \ra \cc_\infty \cb
\]
l'$\ai$-foncteur de Yoneda qui sera dÈfini en \ref{definition_ai-foncteur_de_Yoneda}.
Par le corollaire \ref{corollaire2_lemme_clef}, les deux compositions de $\ai$-foncteurs
\[
\ca \arr{f} \cb \arr{y} \cc_\infty \quad \mbox{et} \quad \ca_x \arr{f_x} \cb_{x'} \arr{y} \cc_\infty \cb_{x'}
\]
sont donnÈes par un $\ca$-$\cb$-bipolydule $M$ et un $\ca_x$-$\cb_{x'}$-bipolydule $N.$
Nous vÈrifions qu'on a
\[
\lrus{x}{M}{x'} = N. 
\]
}\end{remarque}

\section{Le cas topologique} \label{section_torsion_II}

Soit ${\mathbb A}$ un ensemble et $\ca$ une $\ai$-catÈgorie sur ${\mathbb A}.$
Nous traitons ici de la torsion de $\ca$ par un morphisme
$x : e \ra \ca$ qui n'est pas tensoriellement nilpotent.
La somme de gauche dans l'Èquation de Maurer-Cartan (voir \ref{remarque_nilpotence_tensorielle})
\[
\sum_{i\geq 1} m_i (x \tpo i) = 0
\]
appliquÈe ‡ $\id_A$ n'est plus finie
mais l'ÈgalitÈ a encore un sens : si $\ca$ est munie d'une topologie,
nous interprÈtons l'Èquation ci-dessus comme la convergence
de la sÈrie vers $0.$
Nous montrons ‡ l'aide d'un artifice algÈbrique que les formules
donnant les structures tordues dans le cas o˘ $x$ est un ÈlÈment tordant
tensoriellement nilpotent donnent aussi
des structures tordues dans le cas o˘ $\ca$ est topologique et $x$
vÈrifie l'Èquation de Maurer-Cartan.

\subsection{DÈfinitions} \label{section_definitions_topologie}

{\noindent \bf La terminologie des objets topologiques}\\

Soit $(\mathsf M, \ts, e)$ une $\corps$-catÈgorie abÈlienne monoÔdale.
Une {\em topologie} \index{topologie} sur un objet $V \in \mathsf M$ est une filtration
dÈcroissante
\[
V_0 \supset V_1 \supset V_2 \supset \cdots \supset V_i \supset \cdots
\]
(voir \cite[Chap.~III~ß2~n∞5]{Bourbaki61}). 
La topologie est
{\em sÈparÈe} si $\cap_{i \in \N} V_i = 0$. On dira alors que les sous-objets $V_i$, $i\geq 1$, sont
un systËme de {\em voisinages} de~$0.$
Un objet {\em topologique} de $\sf M$ est un objet $M$ muni d'une topologie.
Sa {\em complÈtion} \index{complÈtion} \indexnotation{hatV} est la limite
\[
\hat V = \lim_{i\geq 0} V/V_i.
\]
Un objet $V$ est complet si $V = \hat V .$ 
Soit $V$ et $V'$ deux objets topologiques.
Un {\em morphisme} $f : V \ra V'$ est un morphisme continu.
Il est {\em contractant} \index{contractant} s'il vÈrifie
\[
f(V_i) \subset V'_i,  \quad i \geq 1.
\]
L'ÈlÈment neutre $e$ pour le produit tensoriel est muni de la topologie
discrËte. Le produit tensoriel $V \ts V'$ est topologique pour
le systËme de voisinages
\[
\big(V \ts V'\big)_i = \sum_{i_1 + i_2 \geq i} V_{i_1} \ts V_{i_2}, \quad i\geq 0.
\]
La catÈgorie des objets topologiques de $\sf M$, munie du produit tensoriel topolo\-gique
et de l'objet neutre $e$ est une catÈgorie
monoÔdale.
Le {\em produit tensoriel complet} $V \hat \ts V'$ \indexnotation{hatts} est la limite
\[
V \hat \ts V' = \lim_{i \geq 0} (V \ts V')/(V \ts V')_i.
\]
La catÈgorie des objets complets de $\sf M$, munie du produit tensoriel
complet et de l'objet neutre $e$ est une catÈgorie
monoÔdale.\\

{\noindent \bf $\ai$-structures topologiques}\\

Soit $\sf C$ une catÈgorie de base (voir \ref{categorie_de_base}).

\begin{definition}{\em \label{definition_ai-algebres_topologiques}
Une $\ai$-algËbre $A$ dans $\sf C$ est {\em topologique}
\index{A-infini algebre@{$\ai$-algËbre}!topologique}
si $A$ est muni
d'une topologie sÈparÈe et si les multiplications $m_i : A \tp i \ra A$, $i \geq 1$,
sont des morphismes continus contractants. 
Soit $A$ et $A'$ des $\ai$-algËbres topologiques.
Un $\ai$-morphisme {\em topologique} $f : A \ra A'$
est un $\ai$-morphisme tel que les morphismes
$f_i$, $i \geq 1,$ sont des  morphismes continus contractants.
Nous dÈfinissons de maniËre similaire les homotopies entre $\ai$-morphismes.
}\end{definition}

Soit $\sf C'$ une catÈgorie de Grothendieck munie d'une action ‡ droite
de la catÈgorie monoÔdale $\sf C.$ Cette action s'Ètend aux catÈgories des
objets topologiques de $\sf C'$ et $\sf C.$

\begin{definition}{\em
Un $\ca$-polydule {\em topologique} \index{polydule!topologique}
dans $\sf C'$ est un objet topologique sÈparÈ $M$ dans $\sf C'$ muni
d'une structure de $\ca$-polydule dont les multiplications
$m^M_i,$ $ i\geq 1,$ sont des morphismes continus contractants. On dÈfinit de maniËre similaire les
$\ai$-morphismes et les homotopies entre $\ai$-morphismes.
}\end{definition}

\subsection{ElÈments tordants}

\begin{definition}{\em
Soit $A$ une $\ai$-algËbre topologique.
Un morphisme graduÈ $x : e \ra A$ de degrÈ $+1$ est un {\em ÈlÈment tordant (topologique)}
\index{element tordant@{ÈlÈment tordant}}
si son image est dans le voisinage $A_1$ et si la somme
\[
\sum_{i\geq 1} m_i (x \tp i) 
\]
converge vers $0$.
}\end{definition}

\begin{remarque}{\em
Cette somme converge vers une limite bien dÈfinie
car la topologie de $A$ est sÈparÈe, l'image de $x$ est dans $A_1$ et les multiplications
$m_i$, $i \geq 1,$ sont contractantes.
}\end{remarque}

\subsection{AlgËbres locales}

Soit $\mathcal R$ \indexnotation{calR}
la catÈgorie des $\corps$-algËbres commutatives locales $R$ de corps rÈsiduel $\corps$
et dont l'idÈal maximal $\mathfrak m$ est nilpotent.
Soit $R$ un objet de $\mathcal R$. Nous notons $\ce$ la catÈgorie des modules sur $R.$
Soit ${\mathbb O}$, ${\mathbb O}'$ et ${\mathbb O}''$ trois ensembles. Nous notons
$\sf C^{R}({\mathbb O},{\mathbb O}')$ la catÈgorie des foncteurs
\[
{\mathbb O}' \op \times {\mathbb O} \ra \ce
\]
et $\sf C^{R}({\mathbb O}')$ la catÈgorie $\sf C^{R}(\{*\},{\mathbb O}').$ 
Si $M$ et $N$ sont des objets de $\sf C^R({\mathbb O},{\mathbb O}')$ et  $\sf C^R({\mathbb O}',{\mathbb O}'')$,
nous notons $\tso_R$ le produit
tensoriel 
\[
\Big(M \tso_{R} N\Big)(o'',o) = \bigoplus_{o'\in {\mathbb O}'}M(o',o) \ts_R  N(o'',o').
\]
\begin{definition}{\em
Soit $\mathbb A$ un ensemble.
Une {\em $R$-$\ai$-catÈgorie} \index{R-A-infini-categorie@{$R$-$\ai$-catÈgorie}}
est un objet $M$ de $\sf C^{R}({\mathbb A},{\mathbb A})$, muni de morphismes
\[
m_i : M^{\tso_R i } \ra M, \quad i\geq 1,
\]
vÈrifiant l'Èquation $(*_n)$, $n\geq 1$, de la dÈfinition \ref{definition_ai-algebre}
Les {\em $R$-$\ai$-foncteurs}
sont dÈfinis comme en \ref{definition_ai-foncteurs}.
}\end{definition}
Soit $M$ et $M'$ des objets de $\sf C({\mathbb A},{\mathbb A})$ et $i$ un entier $\geq 1$.
Soit
\[
\varphi : M^{\tso i } \ra M'
\]
un morphisme graduÈ. Soit
\[
\varphi^R : (M \ts_\corps \mathfrak m) ^{\tso_R i } \ra M' \ts_\corps  \mathfrak m
\]
le morphisme de $\sf C^{R}({\mathbb A},{\mathbb A})$
dÈfini par la composition
\[
\varphi \ts \mu^{(i)}: (M \ts_\corps  \mathfrak m)^{\tso_R i } \iso (M  \tpo i) \ts_\corps
( \mathfrak m) ^{\ts_R i } \ra M' \ts_\corps  \mathfrak m.
\]
Remarquons que, comme $\mf m$ est nilpotent, il existe un entier $N_0$ tel que
\hbox{$\mf m^{N_0} = 0$.}
Donc le morphisme $\varphi^R$ est nul dËs que $i \geq N_0.$

\begin{remarque}{\em 
Soit $\ca$ un objet de $\sf C({\mathbb A},{\mathbb A})$ et
\[
m_i : \ca \tpo i \ra \ca,\quad  i\geq 1,
\]
des morphismes graduÈs de degrÈ $2-i.$
Nous vÈrifions que les morphismes $m_i$, $i\geq 1$, dÈfinissent une structure de $\ai$-catÈgorie sur $\ca$
si et seulement si, pour tout $R \in \mathcal R$, les morphismes $m_i^R$, $i \geq 1$, dÈfinissent
une structure de $R$-$\ai$-catÈgorie sur $\ca \ts_\corps \mf m.$

Soit $\ca$ une $\ai$-catÈgorie et $R$ un objet de $\mathcal R$. Nous notons $\ca^R$ la $R$-$\ai$-catÈgorie
$\ca \ts_{\corps} \mf m$ sur ${\mathbb A}$ associÈe ‡ $\ca.$

Soit ${\mathbb A}$ et ${\mathbb B}$ deux ensembles et
$\ca$ et $\cb$ deux $\ai$-catÈgories sur ${\mathbb A}$ et ${\mathbb B}.$
Nous vÈrifions que des morphismes graduÈs $f_i$, $i\geq 1$, de degrÈ $1 -i$, dÈfinissent un $\ai$-foncteur
\[
f : \ca \ra \cb
\]
si et seulement si, pour tout $R \in \mathcal R$,  les morphismes $f^R_i$, $i\geq 1$,
dÈfinissent un $R$-$\ai$-foncteur
\[
f^R : \ca^R \ra \cb^R.
\]
Notons que les morphismes $m^R_i$ et $f_i^R$ sont nuls dËs que $i$ excËde le degrÈ de nilpotence de
l'idÈal maximal de $R$.
}\end{remarque}

{\noindent \bf Construction bar $B^R$}\\

Soit $R$ un objet de $\mathcal R.$
Le lemme \ref{cogebres_tensorielles} reste valable dans la catÈgorie $\sf C^{R}({\mathbb A},{\mathbb A}).$
En particulier la construction bar dÈfinit un foncteur pleinement fidËle
\[
B^R : \mathsf{Alg}_\infty^R \ra \cocog^R,
\]
o˘ $\mathsf{Alg}_\infty^R$ et $\cocog^R$ sont les catÈgories $\aia$ et $\cocog$ dans $\sf C^{R}({\mathbb A},{\mathbb A}).$\\

{\noindent \bf Rappel sur la complÈtion}\\

Soit $R$ un objet de $\mathcal R$.
Soit $V$ et $W$ des ${\mathbb A}$-${\mathbb A}$-$R$-bimodules. Nous munissons
la {\em $R$-cogËbre tensorielle rÈduite}
\[
\ctr V = \bigoplus_{i\geq 1} V ^{\tso_R i}
\]
de la {\em topologie canonique} dont la base de voisinages de $0$ est
\[
\bigoplus_{i \geq n} V^{\tso_R i}, \quad n \geq 1.
\]
La comultiplication est un morphisme continu pour cette topologie.
Rappelons que $\ct V$ dÈsigne la cogËbre co-augmentÈe $(\ctr V)\+.$
Nous la munissons des voisinages dÈfinis de la mÍme maniËre.

\begin{remarque}{\em
Un morphisme de $\sf C^R({\mathbb A},{\mathbb A})$
\[
\ct V \ra \ct W \quad \Big( \mbox{resp.} \ctr V \ra \ctr W \Big)
\] 
est continu si et seulement si sa matrice des composantes
\[
\bigoplus_{j\geq 0} V ^{\tso_R j} \ra \bigoplus_{i\geq 0} W ^{\tso_R i}
\quad \Big( \mbox{resp.} \bigoplus_{j\geq 1} V ^{\tso_R j} \ra \bigoplus_{i\geq 1} W ^{\tso_R i}\Big)
\]
a un nombre fini de coefficients non nuls sur chaque ligne.
En particulier, un morphisme de cogËbres
$f$ (resp.~une $(f',f'')$-codÈrivation $h$, o˘ $f'$ et $f''$ sont des
morphismes de cogËbres)
\[
\ctr V \ra \ctr W
\]
 est continu si et seulement si
les morphismes $f_i$,  $i\geq 0,$ (resp.~les morphismes $f'_i$, $f''_i$ et $h_i$, $i\geq 0,$)
sont presque tous nuls.
}\end{remarque}

La {\em $R$-cogËbre tensorielle complËte rÈduite} \indexnotation{hatctr}
$\hat \ct V$ est la complÈtion de $\ctr V.$
Elle a pour espace topologique sous-jacent
\[
\prod_{i\geq 1} V ^{\tso_R i}.
\]
Chaque morphisme continu $\varphi : \ctr V \ra \ctr W$ de $\sf C^R({\mathbb A},{\mathbb A})$
donne un morphisme
\[
\hat f : \hat \ct V \ra \hat \ct W.
\]
La {\em cogËbre tensorielle complËte co-augmentÈe} $\hat \ct  \+ V$ est la
co-augmentation de $\hat \ct V.$

\begin{lemme} \label{cogebres_tensorielles_completes}
Soit $V$ un objet de $\gr \sf C^R({\mathbb A},{\mathbb A})$
et $C$ une cogËbre graduÈe topologique topologique dans $\sf C^R({\mathbb A},{\mathbb A}).$
Soit $f'$ et $f''$ deux morphismes continus de cogËbres
\[
C \ra  \hat \ct \+ V.
\]
Un morphisme continu co-unitaire de cogËbres
complËtes (resp. une $(f',f'')$-codÈrivation) $C \ra  \hat \ct \+ V$ est dÈterminÈ par sa composition avec
la projection $\hat \ct \+ V \ra V.$
\findem
\end{lemme}

\subsection{Torsion des $\ai$-catÈgories}
\label{section_torsion_ai-categories_II}

{\noindent \bf Torsion de la diffÈrentielle de $B\ca^R$}\\

Soit ${\mathbb A}$ un ensemble.
Soit $\ca$ une $\ai$-catÈgorie topologique sur ${\mathbb A}$, i.~e.~une $\ai$-algËbre topologique
dans $\sf C({\mathbb A},{\mathbb A}).$ Soit $x : e \ra \ca$ un ÈlÈment tordant (topologique) de $\ca.$

Soit $R$ un objet de $\mathcal R.$ Notons $N_0$ l'indice de nilpotence
de son idÈal maximal $\mf m.$
Soit $\ca^R$ la $R$-$\ai$-catÈgorie sur ${\mathbb A}$ associÈe ‡ $\ca.$ 
Soit $\ctr S \ca^R$ la $R$-cogËbre tensorielle rÈduite  et
$\hat \ct \+ S\ca^R$ la $R$-cogËbre co-augmentÈe associÈe ‡ sa complÈtion.
La diffÈrentielle de la construction bar $B^R\ca^R$
\[
b^R : \ctr S\ca^R \ra \ctr S\ca^R
\]
est continue car les morphismes $m_i^R$ sont nuls pour $i \geq N_0.$
Notons $\hat b^R$ la diffÈrentielle de $\hat \ct \+ S\ca^R$ induite par
$b^R.$
Soit
\[
x^R : e^R \ra \ca^R
\]
le morphisme induit par $x$ et
\[
g : e \oplus \hat \ct  S\ca^R = \hat \ct \+S\ca^R\+ \ra S\ca^R.
\]
le morphisme de composantes le morphisme $x^R$ et la projection $p_1$ sur $S\ca^R.$
Soit le morphisme de ${\mathbb A}$-${\mathbb A}$-$R$-bimodules
\[
\phi_x^R :  \hat \ct \+ \ca^R \ra \hat \ct \+ \ca^R
\]
dont la composition avec la projection sur $(S\ca^R)\tpo n$
vaut
\[
g \tpo n \circ \Delta^{(n)} : \hat \ct \+S\ca^R \ra (S\ca^R)\tpo n
\]
si $n \geq 1$ et $\Id_e$ sinon.
Comme le morphisme $\phi_x$ de la section \ref{section_torsion_ai-categories_I},
le morphisme $\phi_x^R$ est un automorphisme continu co-unitaire (non co-augmentÈ) de cogËbres graduÈes et
la matrice de ses coefficients
\[
\prod_{j\geq 0} (S\ca^R)^{\tso_R j} \ra \prod_{i\geq 0} (S\ca^R)^{\tso_R i}
\]
est triangulaire infÈrieure; sa diagonale est celle de l'identitÈ. Soit la composÈe 
\[
D^R_x = (\phi^R_x)^{-1} \circ \hat{b^R}\+ \circ \phi^R_x.
\]
Comme $x$ est un ÈlÈment tordant, nous avons
\[
\sum_{1 \leq i \leq N_0} \hat b^R_i \big((x^R)^{\tso_R i }\big) = 0.
\]
Notons que le dÈfaut de nilpotence tensorielle est palliÈ par l'annulation des
morphismes $b^R_i$ pour $i \geq N_0.$
Comme dans la section \ref{section_torsion_ai-categories_I}, la composÈe
$D^R_x \circ \coaugmentation$ est nulle.
Soit $b^R_x$ le morphisme donnÈ par la flËche verticale de droite du diagramme
de suites exactes
\[
\xymatrix{0 \ar[r] & e \ar[r]^(.3){\coaugmentation } \ar[d]_0 &\hat \ct \+S\ca^R
\ar[r]
\ar[d]_{D^R_x}
                    &\hat \ct S\ca^R \ar[r] \ar[d]_{{b^R_x}} & 0 \\
              0 \ar[r] & e \ar[r]^(.3){\coaugmentation } & \hat \ct \+S\ca^R \ar[r]
                    &\hat \ct S\ca^R \ar[r] & 0.  \\}
\]
C'est une diffÈrentielle de la cogËbre $\hat \ct S\ca^R.$

\begin{lemme} \label{lemme_torsion_ai-categories_II}
La sous-cogËbre $\ct S\ca^R$ de $\hat \ct S\ca^R$ est stable par la diffÈrentielle~$b^R_x$.
La composÈe $p^R_1 \circ b^R_x$ restreinte ‡ $(S\ca^R)\tpo i$
est Ègale ‡ la somme
\[
\sum_{l} \sum b^R_{l+m}((sx) \tpo{l_0} \tso \Id_{S\ca^R} \tso (sx) \tpo{l_1} \tso \hdots \tso \Id_{S\ca^R}
 \tso (sx) \tpo{l_{i-1}} \tso \Id_{S\ca^R} \tso (sx) \tpo{l_{i}}),
\]
o˘ $l_0 + \hdots + l_{i} = l$.
\end{lemme} 

\dem Identique ‡ celle du lemme \ref{lemme_torsion_ai-categories_I} \findem \\

{\noindent \bf $\ai$-catÈgorie tordue par $x$}\\

Soit ${\mathbb A}$ un ensemble, $\ca$ une $\ai$-catÈgorie topologique
sur ${\mathbb A}$ et $x : e \ra \ca$ un ÈlÈment tordant.
Soit les morphismes
\[
m^x_i : \ca \tpo i \ra \ca, \quad i \geq 1,
\]
dÈfinis par la somme
\[
\sum_{l} \sum (-1)^{s} m^{\ca}_{l+i}(x \tpo{l_0} \tso \Id_\ca \tso x \tpo{l_1} \tso \hdots \tso \Id_\ca \tso
x \tpo{l_{i-1}} \tso \Id_\ca \tso x \tpo{l_{i}}),
\]
o˘ l'exposant du signe est $s = \sum_{1 \leq t \leq i} t \times l_t$.
Remarquons que ces sommes convergent vers des limites bien dÈfinies car $\ca$ est topologiquement sÈparÈ,
l'image de $x$ est dans le voisinage $\ca_1$ et
les compositions $m_i$, $i \geq 1,$ sont des morphismes continus contractants.

\begin{lemme} \label{lemme2_torsion_ai-categories_II}
Les morphismes $m_i^x$, $i \geq 1,$ dÈfinissent une structure d'$\ai$-catÈgorie sur le ${\mathbb A}$-${\mathbb A}$-bimodule
sous-jacent ‡ $\ca$.
\end{lemme}

\dem
Le lemme sera valide si, pour tout objet $R \in \mathcal R$, les morphismes
$(m^x_i)^R$, $i \geq 1$, dÈfinissent une structure de $R$-$\ai$-catÈgories
sur le $R$-${\mathbb A}$-${\mathbb A}$-bimodule sous-jacent ‡ $\ca^R.$ 

Soit $R$ un objet de $\mathcal R.$
Nous vÈrifions que le morphisme $b_x^R$ du lemme \ref{lemme_torsion_ai-categories_II} est
la codÈrivation
\[
\ct (S\ca^R) \ra \ct (S\ca^R)
\]
construite ‡ partir des $(m_i^x)^R$, $i \geq 1$. Comme elle est une diffÈrentielle nous avons le
rÈsultat.\findem

\begin{definition}{\em \label{definition_torsion_ai-categories_II}
L'{\em $\ai$-catÈgorie (topologique) tordue $\ca_x$} \indexnotation{cax}
\index{A-infini categorie@{$\ai$-catÈgorie}!tordue} \index{torsion}
est le ${\mathbb A}$-${\mathbb A}$-bimo\-dule $\ca_x = \ca$ munie des compositions
\[
m^x_i : \ca_x \tpo i \ra \ca_x, \quad i \geq 1,
\]
dÈfinies ci-dessus.
}\end{definition}

\subsection{Torsion des $\ai$-foncteurs}
\label{section_torsion_ai-foncteurs_II}

Soit ${\mathbb A}$ et ${\mathbb B}$ deux ensembles, $\ca$ et $\cb$ deux
$\ai$-catÈgories topologiques sur ${\mathbb A}$ et ${\mathbb B}$ et $x$ et $x'$ des ÈlÈments tordants
de $\ca$ et $\cb$ tel que pour tout $A \in \mathbb A$,
\[
\sum_{i \geq 1} f_i (x \tpo i)(\id_A) = \id_{\dot f A}.
\]
Remarquons que la somme de gauche converge vers une limite bien dÈfinie car $\cb'$ est topologiquement sÈparÈ,
l'image de $x$ est dans le voisinage $\ca_1$ et car les morphismes
$f_i$, $i \geq 1$, sont contractants. L'ÈgalitÈ ci-dessus exprime la compatibilitÈ
de $x$ et $x'$ ‡ $f$ (voir \ref{section_torsion_ai-foncteurs_I}).
Reprenons les notations $\cb'$, $\cb'_{x'}$ de la section \ref{section_torsion_ai-foncteurs_I}.
Soit les morphismes
\[
f^x_i : \ca \tpo i \ra \cb', \quad i\geq 1,
\]
dÈfinis par la somme (convergente)
\[
\sum_{l} \sum (-1)^{s} f^{\ca}_{l+i}(x \tpo{l_0} \tso \Id_\ca \tso x \tpo{l_1} \tso \hdots \tso \Id_\ca \tso
x \tpo{l_{i-1}} \tso \Id_\ca \tso x \tpo{l_{i}}),
\]
o˘ l'exposant du signe est $s = \sum_{1 \leq t \leq i} t \times l_t$.

\begin{lemme}\label{lemme_torsion_ai-foncteurs}
Les morphismes $f^x_i$ $i \geq 1$, dÈfinissent un $\ai$-foncteur
\[
(\dot f,f_x) : \ca_x \ra \cb_{x'}.
\]
\end{lemme}

\dem Nous allons montrer que, pour tout objet $R \in \mathcal R$, les morphismes $f_i^R$, $i \geq 1$,
dÈfinissent un $\ai$-foncteur
\[
f_x^R : \ca_x^R \ra \cb'{}_{x'}^R,
\]
ou, de faÁon Èquivalente, un morphisme diffÈrentiel graduÈ de cogËbres
\[
F^R_x : B^R\ca^R_x \ra B^R \cb'{}_{x'}^R.
\]

Soit $R \in \mathcal R$. Gr‚ce ‡ la compatibilitÈ de $x$ et $x'$ ‡ $f$,
le morphisme diffÈrentiel graduÈ de cogËbres complËtes co-unitaires
\[
G^R = (\phi^R_{x'})^{-1} \circ \hat F\+ \circ \phi^R_x : \hat \ct \+S \ca_x^R \ra \hat \ct \+ S{\cb'}_{x'}^R
\]
est co-augmentÈ. Il induit donc un morphisme diffÈrentiel graduÈ
\[
F_x : (\hat \ct S\ca_x^R,\hat b_x^R) \ra (\hat \ct S{\cb'}_{x'}^R ,\hat{b'}^R_{x'}).
\]
Soit $i \geq 1$. Nous montrons de maniËre similaire ‡ la dÈmonstration du
lemme \ref{lemme_torsion_ai-categories_II} que la restriction de $F_x$ au
sous-objet  $(S\ca^R_x)\tpo i$
est Ègale ‡ la somme 
\[
\sum_{l} \sum F^R_{l+m}((sx) \tpo{l_0} \tso sa_1 \tso (sx) \tpo{l_1} \tso \hdots \tso sa_{i - 1} \tso
(sx) \tpo{l_{i-1}} \tso sa_i \tso (sx) \tpo{l_{i}}),
\]
o˘ $l_0 + \hdots + l_{i} = l$.
Cette somme est finie car les morphismes $F^R_i$ sont nuls si $i$ excËde le degrÈ de nilpotence
de l'idÈal maximal de $R$.
Nous obtenons ainsi un morphisme de cogËbres
\[
F^R_x : (\ctr S\ca^R_x,b^R_x) \ra (\ctr S\cb'{}_{x'}^R,{b'}^R_{x'})
\]
qui est diffÈrentiel graduÈ. Nous avons donc le rÈsultat. \findem

\begin{definition}{\em \label{definition_torsion_ai-foncteurs_II}
L'{\em $\ai$-foncteur tordu}\index{A-infini foncteur@{$\ai$-foncteur}!tordu} \index{torsion}
\indexnotation{fx}
\[
(\dot f,f^x) : \ca_x \ra \cb_{x'}
\]
est donnÈ par les morphismes $f^x_i$, $i\geq 1,$ dÈfinis ci-dessus.
}\end{definition}

La proposition (\ref{proposition_torsion_ai-foncteurs_equiv_faible_I}) reste clairement valide
dans le cas topologique.

\subsection{Torsion des $\ca$-$\cb$-bipolydules}\label{section_torsion_bipolydules_II}

Les dÈtails sont omis car ils sont similaires aux deux derniËres sections.\\

Soit ${\mathbb A}$ et ${\mathbb B}$ deux ensembles, $\ca$ et $\cb$ deux $\ai$-catÈgories topologiques sur ${\mathbb A}$ et ${\mathbb B}$ et
$M$ un $\ca$-$\cb$-bipolydules topologiques. Soit $x$ et $x'$ des
ÈlÈments tordants de $\ca$ et $\cb.$

\begin{definition}{\em \label{definition_torsion_bipolydules_II}
Le {\em $\ca_x$-$\cb_{x'}$-bipolydule} $\lrus{x}{M}{x'}$ \index{bipolydule!tordu}\index{torsion}
\indexnotation{xMx'}
a pour multiplications
\[
m^{x,x'}_{i,j} : \ca_x \tpo i \tso \lrus{x}{M}{x'} \tso \cb_{x'}
\ra \lrus{x}{M}{x'}, \quad i,j\geq 0,
\]
dÈfinies par la somme (convergente)
\[
\sum_{l,k\geq 0} \sum (-1)^{s} m_{i+l,j+k}
(x \tpo{l_0} \tso \Id_\ca \hdots \Id_\ca \tso x \tpo{l_i} \tso \Id_M
\tso {x'} \tpo{k_o} \tso \Id_\cb 
\hdots  \Id_\cb \tso {x'} \tpo{k_{j}}),
\]
o˘ l'exposant du signe est
\[
s = \Big(\sum_{1 \leq t \leq i} t \times l_t \Big) + \Big(\sum_{1 \leq t \leq j} (j+ t) \times l_t \Big) 
\]
}\end{definition}

%% file: Objets_tordus.tex
{\bf \noindent Introduction}\\

\noindent Soit ${\mathbb A}$ un ensemble et $\ca$ une $\ai$-catÈgorie strictement unitaire sur ${\mathbb A}.$
Notons $\gr (H^*\!\ca)$ la catÈgorie des $H^*\!\ca$-modules graduÈs dont les morphismes
sont les morphismes graduÈs.
Dans cette section, nous relevons le foncteur de Yoneda
\[
H^* \ca \ra \gr (H^*\!\ca), \quad  A \mapsto (H^*\!\ca)(\?,A),
\]
en un $\ai$-foncteur 
\[
y : \ca \ra \cc_\infty \ca, \quad A \mapsto \ca(\?,A).
\]
Nous montrons ensuite le rÈsulat principal de ce chapitre (\ref{theoreme_factorisation_Yoneda}) :
{\em l'$\ai$-foncteur $y$ se factorise en
\[
\ca \arr{y'} \TW \ca \arr{y''} \cc_\infty \ca
\]
o˘ $\TW \ca$ est l'$\ai$-catÈgorie des objets tordus,
$y'$ est  un $\ai$-foncteur strict et pleinement fidËle 
 et $y''$ induit une Èquivalence}
\[
H^0\TW \ca \arr{\sim} \tria \ca \subset \cd_{\infty}\ca.
\]
La contruction des objets tordus dans le cas o˘ $\ca$ est diffÈrentielle graduÈe est due
‡ A.~I.~Bondal et M.~M.~Kapranov
\cite{Bondal91}, sa gÈnÈralisation aux $\ai$-catÈgories ‡ M.~Kontsevich \cite{Kontsevich94}.
RÈcemment, K.~Fukaya a construit indÈpendamment l'$\ai$-foncteur de Yoneda \cite{Fukaya01a}.\\

{\bf \noindent Plan du chapitre}\\

\noindent Dans la section \ref{section_ai-foncteur_de_Yoneda}, nous dÈfinissons l'$\ai$-foncteur de Yoneda et
nous ÈnonÁons le thÈorËme principal (\ref{theoreme_factorisation_Yoneda}). Le reste du chapitre (sauf
la section \ref{section_modele_differentiel_gradue}) est dÈdiÈ ‡ la dÈmonstration de ce thÈorËme.
Dans la section \ref{section_ai-categorie_des_objets_tordus}, nous construisons l'$\ai$-catÈgorie $\TW \ca$
des objets tordus. Les compositions de l'$\ai$-catÈgorie $\TW \ca$ sont obtenues
par torsion (voir chapitre \ref{chapitre_torsion}).
Nous montrons ensuite que l'$\ai$-catÈgorie $\TW \ca$ jouit d'une propriÈtÈ universelle. Nous en
dÈduirons l'existence de
la factorisation $y'' \circ y'$ de $y$.
Dans la section \ref{section_ai-foncteur_y''}, nous construisons explicitement
l'$\ai$-foncteur $y''$. Dans la section \ref{section_equivalence_tria_ca_H^0TWca}, nous montrons que 
l'$\ai$-foncteur de Yoneda $y$ induit des quasi-isomorphismes
entre les espaces de morphismes et nous en dÈduisons l'Èquivalence
\[
H^0\TW \ca \iso \tria \ca \subset \cd_{\infty}\ca.
\]

Dans la section \ref{section_modele_differentiel_gradue}, nous montrons que toute $\ai$-catÈgorie
homologiquement unitaire $\ca$ admet un {\em modËle diffÈrentiel graduÈ strictement 
unitaire}, c'est-‡-dire un $\ai$-quasi-isomorphisme homologiquement unitaire
$f : \ca \ra \ca'$ vers une catÈgorie diffÈrentielle
graduÈe strictement unitaire. 

Dans la section \ref{section_categories_stables}, nous montrons que toute catÈgorie
triangulÈe algÈbrique qui est engendrÈe par un ensemble d'objets
est $\ai$-prÈ-triangulÈe, i.~e.~elle est Èquivalente
‡ $H^0\TW \ca$, pour une certaine $\ai$-catÈgorie $\ca$.

\section{Le plongement de Yoneda}
\label{section_ai-foncteur_de_Yoneda}

Comme $\ca$ est une $\ai$-catÈgorie, le ${\mathbb A}$-${\mathbb A}$-bimodule $\ca$,
muni des morphismes $m_{i,j} = m^{\ca}_{i+1+j}$, $i,j\geq 0$,
est un $\ca$-$\ca$-bipolydule.
Par la remarque \ref{remarque2_lemme_clef}, nous avons
un $\ai$-foncteur
\[
y : \ca \ra \cc_\infty \ca,
\]
dont l'application sous-jacente
\[
\dot y : {\mathbb A} \ra \cc_\infty \ca
\]
envoie un objet  $A\in \ca$ sur le {\em $\ca$-polydule} \indexnotation{Apw}
\[
A \pw = \ca(\? , A).
\]
Pour tout $i\geq 1,$ le morphisme graduÈ
\[
y_i : \ca \tpo i  \ra  \lrus{\dot f}{\big(\cc_\infty \ca\big)}{\dot f},
\]
envoie un ÈlÈment $x \in (\ca \tpo i)(A,A')$
sur la suite de morphismes de ${\mathbb A}$-modules graduÈs
\[
\begin{array}{rclr}
 \ca(\?,A) \tso \ca \tpo {j-1} & \ra  & \ca(\?,A') ,&  \quad  j \geq 1.\\
x' \tso x'' & \mapsto & (-1)^{|x|+1}m_{i+1+j}(x' \tso x \tso x'')&
\end{array}
\]

\begin{definition}{\em \label{definition_ai-foncteur_de_Yoneda} \indexnotation{y}
L''{\em $\ai$-foncteur de Yoneda} est l'$\ai$-foncteur
$y : \ca \ra \cc_\infty \ca$. \index{Yoneda (A-infini foncteur de)@{Yoneda ($\ai$-foncteur de)}}
\index{A-infini foncteur@{$\ai$-foncteur}!de Yoneda}
}\end{definition}

\begin{definition}{\em
Un $\ai$-foncteur strict $f$ est {\em pleinement fidËle} si
\index{A-infini foncteur@{$\ai$-foncteur}!pleinement fidËle}
\[
f_1 : \ca \ra \lrus{\dot f}{\cb}{\dot f}
\]
est un isomorphisme de complexes.
}\end{definition}

\begin{definition}{\em
Soit $\calt$ une catÈgorie triangulÈe et $\mathbb T'$ un sous-ensemble de l'ensemble $\mathbb T$
des objets de $\calt$.
Notons $\tria \mathbb T'$ la {\em plus petite sous-catÈgorie triangulÈe} de $\calt$
qui contient les objets de $\mathbb T'.$ Elle est stable par sommes finies.
Soit $\ca$ une $\ai$-catÈgorie strictement unitaire et $\cd_\infty \ca$ sa catÈgorie
dÈrivÈe (voir \ref{section_categorie_derivee_generale}). 
Notons $\tria \ca$ la plus petite sous-catÈgorie triangulÈe de $\cd_\infty \ca$ 
qui contient tous les $\ca$-polydules $A \pw$, $A \in \Obj \ca$.
}\end{definition}

Dans ce chapitre, nous allons montrer l'ÈnoncÈ de M.~Kontsevich \cite{Kontsevich94}, \cite{Kontsevich98} suivant :

\begin{theoreme}[voir aussi K.~Fukaya \cite{Fukaya01a}]{\label{theoreme_factorisation_Yoneda}} \index{objet tordu}
Soit une $\ai$-catÈgorie $\ca$ avec des identitÈs strictes.
Il existe une $\ai$-catÈgorie $\TW \ca$
et une factorisation de l'$\ai$-foncteur de Yoneda
\[
 \ca \arr{y'} \TW \ca \arr{y''} \cc_{\infty} \ca
\]
telle que l'$\ai$-foncteur $y'$ est strict et pleinement
fidËle et l'$\ai$-foncteur $y''$ induit une Èquivalence
\[
H^0\TW \ca \iso \tria \ca \subset \cd_{\infty}\ca .
\]
\end{theoreme}
\dem Voir les trois sections suivantes.\\

\section{L'$\ai$-catÈgorie des objets tordus}
\label{section_ai-categorie_des_objets_tordus}

Soit $\Lambda$ une algËbre associative unitaire (non graduÈe). 
Nous notons $\cc^b(\free \Lambda)$ la {\em sous-catÈgorie} de $\cc \Lambda$ formÈe
des complexes bornÈs de $\Lambda$-modules libres de rang fini.
L'image $\cd^b(\free \Lambda)$ de la catÈgorie $\cc^b(\free \Lambda)$ par le foncteur
\[
\cc \Lambda \ra \cd \Lambda
\]
est Èquivalente ‡ la catÈgorie $\tria \Lambda$.
Les objets de $\cc^b(\free \Lambda)$ sont fibrants et cofibrants dans la catÈgorie des complexes
$\cc \Lambda$. Si $M$ et $M'$ sont des objets de $\cd^b(\free \Lambda)$,
les morphismes $M \ra M'$ dans  $\tria \Lambda$
sont donc en bijection avec les classes d'homotopie de morphismes $M\ra M'$
de $\Modu \Lambda$. Cette description des morphismes  permet de faire des calculs dans $\tria \Lambda$.
Le but de cette section est de gÈnÈraliser la construction
\[
\Lambda \leadsto \cc^b(\free \Lambda)
\]
aux $\ai$-catÈgories. Soit $\ca$ une $\ai$-catÈgorie.
Le rÙle de la catÈgorie $\cc^b(\free \Lambda)$ sera jouÈ par l'$\ai$-catÈgorie $\TW \ca$
des objets tordus. L'Èquivalence entre $\cd^b (\free \Lambda)$ et $\tria \Lambda$ sera remplacÈe
par une Èquivalence
\[
H^0 \TW \ca \arr{\sim} \tria \ca \subset \cd_\infty \ca.
\]
La construction $\ca \leadsto \TW \ca$ est la gÈnÈralisation  aux $\ai$-catÈgories \cite{Kontsevich94}
de la construction due ‡ A.~I.~Bondal et M.~M.~Kapranov
\cite{Bondal91} qui associe ‡ une catÈgorie diffÈrentielle graduÈe la catÈgorie de 
ses objets tordus (voir \ref{remarque_Kapranov}).\\

Pour rendre la construction qui suit plus intuitive, commenÁons par rÈinterprÈter
les objets de $\cc^b(\free \Lambda).$
Un complexe bornÈ $M$ de $\Lambda$-modules libres de rang fini est donnÈ
par ses composantes
\[
(M_r,M_{r+1},\hdots ,M_{l-1},M_{l}), \quad r\leq l, \quad r,l \in \Z,
\]
o˘ chaque $M_i$, $r \leq i \leq l$, est la suspension itÈrÈe d'un $\Lambda$-module libre de rang fini,
et par un morphisme de degrÈ $+1$
\[
\delta : \bigoplus_{r \leq j \leq l} M_j \ra \bigoplus_{r\leq i\leq l} M_i
\]
dont la matrice est strictement triangulaire infÈrieure et telle que $\delta \circ \delta = 0$.

Supposons maintenant que $\Lambda$ est une algËbre diffÈrentielle graduÈe. 
Les extensions itÈrÈes
dans la catÈgorie des complexes munie de la structure exacte
donnÈe par les suites de complexes qui se scindent en tant que suites de $\Lambda$-modules
graduÈs sont dÈcrites de la maniËre suivante.
Soit $M_i$, $r \leq i \leq l,$ des objets de $\Modu \Lambda$
qui sont des sommes finies de suspensions itÈrÈes de $\Lambda$.
Notons $d$ la diffÈrentielle de la somme des $M_i$, $r \leq i \leq l$.
Une extension itÈrÈe des objets $M_i$, $r \leq i \leq l,$ est donnÈe par une matrice de
la mÍme forme que ci-dessus qui vÈrifie l'Èquation de Maurer-Cartan
\[
d \circ \delta + \delta \circ d + \delta^2 = 0.
\]
La diffÈrentielle de l'extension itÈrÈe $M = \bigoplus_{r \leq j \leq l} M_j $ est la somme $d + \delta$. \\

{\noindent \bf Saturation par dÈcalages de $\ca$}\\

Soit $\Z \ca$ {\em l'$\ai$-catÈgorie} dont les objets sont des couples $(A,n)$, o˘
$A$ est un objet de $\ca$ et $n$ un entier. Les espaces de morphismes sont
dÈfinis par 
\[
\Z\ca ((A,n),(B,m)) = S^{m-n}\ca (A,B).
\]
Les compositions $m^{\Z\ca}_{i}$, $i \geq 1,$
\[
\xymatrix{
\Z\ca ((A_{i-1},n_{i-1}),(A_{i},n_{i})) \ts \hdots \ts
\Z\ca ((A_{0},n_{0}),(A_{1},n_{1})) \ar[d]^{m^{\Z\ca}_{i}} \\
\Z\ca((A_{0},n_{0}),(A_{i},n_{i}))}
\]
sont dÈfinies par
\[
(-1)^{i (n_i - n_0)} s^{n_i-n_o}\circ m_i \circ \big((s^{n_i-n_{i-1}})^{-1} \tso \hdots \tso (s^{n_1-n_{0}})^{-1}\big)
\]
(un calcul montre que ces compositions dÈfinissent bien une $\ai$-catÈgorie).\\

{\noindent \bf Saturation par extensions de $\Z \ca$}

\begin{definition}{\em \index{extension}
Une {\em extension itÈrÈe} $M$ d'objets de $\Z\ca$ est une suite
\[
(M_r,M_{r+1},\hdots ,M_{l-1},M_{l}), \quad r\leq l, \quad r,l \in \Z,
\]
munie d'une matrice ‡ coefficients dans $\Z \ca$ de degrÈ $+1$
\[
\delta^M : \bigoplus_{r \leq j \leq l} M_j \ra \bigoplus_{r\leq i \leq l} M_i
\]
qui est strictement triangulaire infÈrieure et vÈrifie l'Èquation de Maurer-Cartan
\[
\sum_{i\geq 1} m^{\Z \ca}_i \big((\delta^M) \tpo i\big) = 0.
\]
Ici, le produit tensoriel $\tso$ est l'extension du produit tensoriel de $\sf C(\mathbb A,\mathbb A)$
aux espaces de matrices ‡ coefficients dans $\Z \ca$.
L'entier $l - n + 1$ s'appelle la {\em hauteur} de l'extension. \index{hauteur}
Une extension itÈrÈe $M$ est {\em dÈgÈnÈrÈe} ou {\em scindÈe}
\index{scindee (extension)@{scindÈe (extension)}} \index{extension scindÈe} si $\delta^M = 0.$
Les extensions itÈrÈes dÈgÈnÈrÈes peuvent Ítre considÈrÈes comme les sommes formelles
d'objets de $\Z \ca$. Nous notons $\mathbb E$ l'ensemble des extensions itÈrÈes de $\Z \ca$.
}\end{definition}

\begin{definition}{\em \label{definition_ai-categorie_extensions}
Soit $M $ et $M'$ deux extensions itÈrÈes de $\Z \ca$.
Notons $\mathsf{Mat}^{\Z \ca}(M,M')$ l'espace graduÈ des matrices ‡ coefficients dans $\Z \ca$
\[
f : \bigoplus_{r \leq j \leq l} M_j \ra \bigoplus_{r'\leq i \leq l'} M'_i.
\]
Les compositions $m_i^{\Z \ca}$, $i\geq 1$, de $\Z \ca$ s'Ètendent clairement
en des compositions de matrices ‡ coefficients dans $\Z \ca$.
Notons $\ce_\ca$ l'{\em $\ai$-catÈgorie} dont les objets sont les extensions
itÈrÈes d'objets de $\Z\ca$ et dont les espaces de morphismes sont
\[
\Hom_{\ce_\ca}(M,M') = \mathsf{Mat}^{\Z \ca}(M,M').
\]
Nous avons clairement une suite d'inclusions de $\ai$-catÈgories 
\[
\ca \subset \Z \ca \subset \ce_\ca.
\]
}\end{definition}

{\noindent \bf L'ÈlÈment tordant nilpotent de l'$\ai$-catÈgorie $\ce_\ca$}\\

Nous rappelons (\ref{section_categories_de_base_aicat})
que $\id_M$ est le gÈnÈrateur de l'espace $e_{\mathbb E}(M,M).$
Soit
\[
x : e_{\mathbb E} \ra \ce_\ca
\]
le morphisme de $\mathbb E$-$\mathbb E$-bimodules qui envoie $\id_M$, $M \in \mathbb E$, sur
\[
\delta^M \in \mathsf{Mat}^{\Z \ca}(M,M).
\]
Le morphisme $x$ est de degrÈ $+1$. Il vÈrifie la condition de nilpotence
tensorielle (\ref{remarque_nilpotence_tensorielle}) car les matrices $\delta^M$
sont strictement triangulaires infÈrieures. Comme les morphismes $\delta^M$, $M \in \mathbb E,$
vÈrifient l'Èquation de Maurer-Cartan, le morphisme $x$ est un ÈlÈment tordant tensoriellement nilpotent.\\

{\noindent \bf La catÈgorie $\TW \ca$}

\begin{definition}{\em
Un {\em objet tordu} est une extension itÈrÈe d'objets de $\Z \ca.$ \index{objet tordu}
Notons $\mathbb{TW}\ca$ l'ensemble de objets tordus. Il est Ègal ‡ l'ensemble $\mathbb E.$
La {\em catÈgorie} $\TW \ca$ \indexnotation{TWca}
des objets tordus est la catÈgorie tordue $\big(\ce_\ca\big)_x$ (voir \ref{section_torsion_ai-categories_I}),
o˘ $x$ est l'ÈlÈment tordant ci-dessus.
}\end{definition}
Si $M$ et $M'$ sont des objets tordus, l'espace de morphismes $M \ra M'$ est donc
\[
\Hom_{\TW \ca} (M,M') = \mathsf{Mat}^{\Z\ca}(M,M').
\]
Remarquons que sur la sous-$\ai$-catÈgorie formÈe des extensions dÈgÈnÈrÈes, les compositions
tordues $m_i^{\ce_x} = m_i^{\TW \ca}$, $i \geq 1$, sont Ègales aux compositions $m_i^\ce$, $i \geq 1$.
Soit $\mathbb E_1$ l'ensemble des extensions (forcÈment dÈgÈnÈrÈe) de hauteur 1 et 
soit
\[
\dot y' : \mathbb A \ra \mathbb E_1,
\]
l'application qui envoie $A$ sur l'extension dÈgÈnÈrÈe de hauteur 1 dont la suite sous-jacente est
le 1-uplet $((A,0))$. C'est une bijection et nous avons un isomorphisme
\[
y'_1 : \ca \arr{\sim} \lrus{\dot y'}{\mathsf{Mat}^{\Z \ca}}{\dot y'} = \lrus{\dot y'}{\TW \ca}{\dot y'}
\]
qui donne clairement un $\ai$-foncteur strict et pleinement fidËle
\[
y' : \ca \ra \TW \ca.
\]

{\noindent \bf La propriÈtÈ universelle de $\TW \ca$}\\

Nous nous inspirons de l'article \cite{Bondal91}.\\

Soit $f : \ca \ra \cb$ un $\ai$-foncteur.
Il induit clairement un $\ai$-foncteur
\[
f : \ce_\ca \ra \ce_\cb
\]
tel que les ÈlÈments tordants $x_\ca$ et $x_\cb$ des $\ai$-catÈgories $\ce_\ca$ et $\ce_\cb$
sont compatibles ‡ $f$ (voir \ref{section_torsion_ai-foncteurs_I}).
Nous obtenons donc un $\ai$-foncteur tordu (voir \ref{section_torsion_ai-foncteurs_I})
\[
\TW f : \TW \ca \ra \TW \cb.
\]
La construction qui associe ‡ une $\ai$-catÈgorie $\ca$ la catÈgorie
des objets tordus $\TW \ca$ est un foncteur
\[
\TW : \aicat \ra \aicat,
\]
o˘ $\aicat$ est la catÈgorie des petites $\ai$-catÈgories.
Nous allons construire un morphisme de foncteurs
\[
\mathsf{Tot} : \TW \circ \TW \ra \TW.
\]
Soit $\ca$ une petite $\ai$-catÈgorie. L'$\ai$-foncteur strict $\mathsf{Tot}(\ca)$
associe ‡ un objet $N$ de $\TW \circ \TW \ca$, donnÈ par une
suite d'objets de $\TW \ca$
\[
(N_r , \hdots , N_l), \quad r \leq l, \quad r,l \in \Z,
\]
et une matrice $\delta^N$ ‡ coefficients dans $\Z \TW \ca$,
l'objet tordu de $\ca$ dont la suite sous-jacente est la concatÈnation
des suites dÈfinissant les $N_i$, $r \leq i \leq l$, et dont la matrice
\[
\delta^{\mathsf{Tot}} : \mathsf{Tot}(N) = \bigoplus{(N_j)_k} \ra \mathsf{Tot}(N) = \bigoplus{(N_i)_l} 
\]
est construite ‡ partir de la matrice $\delta^N$ en remplaÁant les coefficients $\delta^N_{i,j}$ par
les blocs donnÈs par les matrices $\delta^{N_i}$.
Nous vÈrifions que les morphismes de foncteurs de $\aicat$
\[
\unite = y': \Id_{\aicat} \ra \TW \quad \mbox{et} \quad \mathsf{Tot} : \TW \circ \TW \ra \TW
\]
dÈfinissent une monade de la catÈgorie des $\ai$-catÈgories au sens de Quillen et Mac~Lane \cite{May72}.
On rappelle qu'une $\TW$-algËbre $\cg$ est une $\ai$-catÈgorie munie d'un $\ai$-foncteur
\[
\TW \cg \ra \cg
\]
compatible ‡ la structure de monade.
La catÈgorie $\TW \ca$ est clairement la $\TW$-algËbre libre sur $\ca$. En particulier, l'$\ai$-foncteur
$y' : \ca \ra \TW \ca$ est universel parmi les $\ai$-foncteurs
\[
\ca \ra \cg
\]
o˘ $\cg$ est une algËbre sur la monade.

\begin{remarque}{\em \label{remarque_Kapranov}
Si $\cg$ est une catÈgorie diffÈrentielle graduÈe, $\TW \cg$ est une catÈgorie diffÈrentielle
graduÈe. La construction $\cg \leadsto \TW \cg$ correspond ‡ la construction de A.~I.~Bondal
et M.~M.~Kapranov qui associe ‡ $\cg$ la catÈgorie $\mbox{Pr-Tr}^+\cg$ des objets tordus
unilatÈraux \cite[ß4]{Bondal91}.
}\end{remarque}

{\noindent \bf Existence de l'$\ai$-foncteur $y''$}\\

Soit $\ca$ une petite $\ai$-catÈgorie. Soit 
\[
\TW \cc_\infty \ca \ra \cc_\infty \ca
\]
l'$\ai$-foncteur strict qui associe ‡ une extension itÈrÈe $M$ la somme des $M_i$, $r \leq i \leq l$,
munie de la diffÈrentielle $d + \delta_M$, o˘ $d$ est la diffÈrentielle de la somme des $M_i.$
Cet $\ai$-foncteur dÈfinit une structure de $\TW$-algËbre sur $\cc_\infty \ca$.
En particulier, l'$\ai$-foncteur
\[
y : \ca \ra \cc_\infty \ca
\]
se factorise en $y = y'' \circ y'$, o˘ $y''$ est l'$\ai$-foncteur
\[
\TW \ca \ra \cc_\infty \ca
\]
donnÈ par la propriÈtÈ universelle de $\TW \ca.$

\section{L'$\ai$-foncteur $y'' : \TW \ca \ra \cc_\infty \ca$}
\label{section_ai-foncteur_y''}

Dans cette section, nous construisons explicitement l'$\ai$-foncteur
\[
y'' : \TW \ca \ra \cc_\infty \ca.
\]
Par la remarque \ref{remarque3_lemme_clef}, les $\ai$-foncteurs
\[
\TW \ca \ra \cc_\infty \ca
\]
sont en bijection avec les $\TW \ca$-$\ca$-bipolydules strictement unitaires.
Le $\TW \ca$-$\ca$-bipolydule $N''$ associÈ ‡ $y''$ est construit
en tordant (voir la section \ref{section_torsion_bipolydules_I})
un $\ce$-$\ca$-bipolydule $N$. L'$\ai$-foncteur
\[
f : \ce \ra \cc_\infty \ca
\]
associÈ ‡ $N$ est l'extension de l'$\ai$-foncteur de Yoneda $y : \ca \ra \cc_\infty \ca$.
Nous donnons les formules explicites pour les $\ai$-foncteurs $f$ et $y''$.\\

{\noindent \bf Construction de $f : \ce \ra \cc_\infty \ca$}\\

Nous rappelons (\ref{definition_ai-categorie_extensions}) que nous avons
une suite d'inclusions de $\ai$-catÈgories 
\[
\ca \subset \Z \ca \subset \ce
\]
et que $y : \ca \ra \cc_\infty \ca$
dÈsigne l'$\ai$-foncteur de Yoneda (\ref{definition_ai-foncteur_de_Yoneda}).
Ce dernier s'Ètend en un $\ai$-foncteur
\[
\Z \ca \ra \cc_\infty \ca, \quad (A,n) \mapsto S^n(\dot yA) = S^nA \pw
\]
qui envoie un ÈlÈment
\[
x \in \Z\ca ((A_{i-1},n_{i-1}),(A_{i},n_{i})) \ts \hdots \ts
\Z\ca ((A_{0},n_{0}),(A_{1},n_{1}))
\]
sur le morphisme de $\ca$-polydules $S^{n_0}A_0\pw \ra S^{n_i}A_i\pw$
dÈfini par l'ÈlÈment de
\[
\Hom_{\cc_\infty \ca}(S^{n_0}A_0\pw , S^{n_i}A_i\pw)
\iso S^{n_i-n_0}\Hom_{\cc_\infty \ca}(A_0\pw , A_i\pw)
\]
donnÈ par
\[
s^{n_i-n_o}\circ y_i \circ \big((s^{n_i-n_{i-1}})^{-1} \tso \hdots \tso (s^{n_1-n_{0}})^{-1}\big)(x).
\]
Nous notons aussi $y$ cet $\ai$-foncteur. Nous l'Ètendons maintenant en un $\ai$-foncteur
\[
\ce \ra \cc_\infty \ca.
\]
Nous dÈfinissons une application
\[
\dot f : \mathbb E \ra \Obj \cc_\infty \ca
\]
qui envoie une extension itÈrÈe $M$, donnÈe par une suite $M_i$, $r \leq i \leq l$, et une matrice
$\delta^M$, sur le $\mathbb A$-module qui est la somme
\[
\sum_{r \leq i \leq l} \dot y M_i.
\]
Sa structure de  $\ca$-polydule
est induite par celle de la remarque \ref{remarque1_polydules}.
Notons que la matrice $\delta^M$ n'intervient pas dans la dÈfinition de l'image de $M.$
Les morphismes $y_i : (\Z \ca) \tpo i \ra \cc_\infty \ca$ s'Ètendent clairement
en des morphismes
\[
\Big(\mathsf{Mat}^{\Z \ca}\Big)\tpo i \ra \cc_\infty \ca.
\]
Ceci nous fournit un $\ai$-foncteur que nous notons $f : \ce \ra \cc_\infty \ca$ et
nous avons clairement la factorisation $y = f \circ y'$.
Par la remarque \ref{remarque3_lemme_clef}, l'$\ai$-foncteur $f$ est donnÈ par un $\ce$-$\ca$-bipolydule
$N$ qui, en tant que  $\mathbb E$-$\mathbb A$-bimodule, est
\[
(A,M) \mapsto \bigoplus_{r \leq i \leq l} S^{n_i}\ca (A,A_i),
\]
o˘ $M_i = (A_i,n_i)$, $r \leq i \leq l$.
Notons
\[
m^{N}_{i,j} : \ce \tpo i \tso N \tso \ca \tpo j \ra N, \quad i,j\geq 0.
\]
les multiplications du $\ce$-$\ca$-bipolydule $N$.
Elles sont clairement induites par l'extension ‡ $\Z \ca$, puis ‡ $\ce$, des compositions
\[
m^\ca_{i,j} = m_{i+1+j}^\ca : \ca \tpo i \tso \ca \tso \ca \tpo j \ra \ca, \quad i,j \geq 0.
\]

{\noindent \bf L'$\ai$-foncteur $y'' : \TW \ca \ra \cc_\infty \ca$}\\

Nous rappelons (\ref{definition_ai-categorie_extensions}) que $x$ dÈsigne
l'ÈlÈment tordant (nilpotent) de $\ce.$
Par la section \ref{section_torsion_bipolydules_I}, nous pouvons tordre $N$ en un
$\ce_x$-$\ca$-polydule $\lrus{x}{N}{} = N''.$
Comme l'$\ai$-catÈgorie $\TW \ca$ est par dÈfinition l'$\ai$-catÈgorie
tordue $\ce_x$, nous obtenons ainsi un $\TW \ca$-$\ca$-bipolydule $N''$ et par la remarque
(\ref{remarque3_lemme_clef}), un $\ai$-foncteur 
\[
y'' : \TW \ca \ra \cc_\infty \ca.
\]
Nous donnons ci-dessous les formules explicites le dÈfinissant.
Le $\mathbb{TW} \ca$-$\mathbb A$-bimodule $N''$ est donnÈ par
\[
(A,M) \mapsto \bigoplus_{r \leq i \leq l} S^{n_i}\ca (A,A_i).
\]
Comme $\mathbb{TW} \ca = \mathbb E$, il est isomorphe en tant que $\mathbb E$-$\mathbb A$-bimodule ‡ $N.$
En tant que $\ce_x$-$\ca$-bipolydule, ses multiplications $m_{i,j}^{N''}$, $i,j \geq 0$ sont donnÈes
(\ref{definition_torsion_bipolydules_I}) par la somme
\[
\sum_{l,k\geq 0} \sum (-1)^{s} m^N_{i+l,j+k}
(x \tpo{l_0} \tso \Id_\ce \tso x \tpo{l_1} \tso \hdots \tso  \Id_\ce \tso x \tpo{l_i}
\tso \Id_N \tso \Id_\ca \hdots  \tso \Id_\ca),
\]
o˘ $\Id_\ce$ dÈsigne l'identitÈ de l'espace des matrices ${\mathsf{Mat}^{\Z \ca}}$ et l'exposant du signe est
\[
s = \sum_{1 \leq t \leq i} t \times l_t .
\]
DÈtaillons maintenant l'application sous-jacente ‡ l'$\ai$-foncteur $y''$
\[
\dot y'' :  \TW \ca \ra \cc_\infty \ca.
\]
Elle envoie une extension itÈrÈe $M$, donnÈe par une suite $M_i$, $r \leq i \leq l$, et une matrice
$\delta^M$, sur le $\mathbb A$-module qui est la somme
\[
\dot y''M = \sum_{r \leq i \leq l} \dot y M_i.
\]
Les multiplications $m_{j}^{\dot y''M}$, $j \geq 1,$ dÈfinissant sa structure de $\ca$-polydule
sont les morphismes $m_{0,j-1}^{N''}$, $j\geq 1$, c'est-‡-dire la somme
\[
\sum_{l \geq 0} m_{l,j-1}^N(x \tpo l \tso \Id_{\dot y''M} \tso \Id_{\ca} \tpo{j-1}) =
\sum_{l \geq 0} m_{l,j-1}^N
\Big([y (\delta^M)] \tpo l \tso \Id_{\dot y''M} \tso \Id_{\ca} \tpo{j-1}\Big).
\]
Remarquons que mÍme si $\dot y''M$ et $\dot fM$ sont isomorphes en tant que
$\mathbb A$-modules, ils diffËrent en tant que $\ca$-polydules.
Le $\ca$-polydule $\dot y''M$ doit Ítre considÈrÈ comme la torsion de $fM$ par $y(\delta^M).$
Regardons maintenant les morphismes $y''_i$, $i\geq 1$, de l'$\ai$-foncteur $y''$.
Ils sont dÈfinis (\ref{corollaire2_lemme_clef}) par la relation
\[
(y''_i)_j = m^{N''}_{i,j-1}.
\]
En d'autre termes, le morphisme $y''_i$, $i \geq 1,$ envoie un ÈlÈment de
\[
\TW \ca (M_{i-1},M_i) \ts \hdots \ts \TW \ca (M_0,M_1)
\]
sur le morphisme de $\ca$-polydules $\varphi : (\dot y''M_0) \ra (\dot y'' M_i)$ donnÈ par la
suite des morphismes $\varphi_j : (\dot y''M_0) \tso \ca \tpo {j-1} \ra (\dot y''M_i)$ valant
\[
\sum_{l\geq  0} \sum (-1)^{s} m^N_{i+l,j-1}
\Big([y(\delta^{M_i})] \tpo{l_0} \tso \Id_{\TW \ca}  \hdots \tso  \Id_{\TW \ca} \tso
[y(\delta^{M_0})] \tpo{l_i}
\tso \Id_{\dot y''M_0} \tso \Id_\ca \tpo{j-1}\Big),
\]
o˘ $\Id_{\TW \ca}$ dÈsigne l'identitÈ de l'espace des matrices ${\mathsf{Mat}^{\Z \ca}}$ et 
l'exposant du signe est
\[
s = \sum_{1 \leq t \leq i} t \times l_t .
\]

Remarquons que l'unitaritÈ stricte de $\ca$ n'est pas intervenue dans 
la dÈmonstration de la factorisation du thÈorËme \ref{theoreme_factorisation_Yoneda}.
Elle joue un rÙle essentiel dans la prochaine section.

\section{L'Èquivalence entre les catÈgories $\tria \ca$ et $H^0\TW \ca$}
\label{section_equivalence_tria_ca_H^0TWca}

Nous rappelons (\ref{remarque_categorie_differentielle_graduee_polydules}) que
les catÈgories $H^0 \cc_\infty \ca$ et $\cd_\infty \ca$ sont Èquivalentes.
Nous montrons ci-dessous que l'$\ai$-foncteur $y'' : \TW \ca \ra \cc_\infty \ca$
induit un foncteur pleinement fidËle
\[
H^0 \TW \ca \ra \cd_\infty \ca.
\]
dont l'image est la catÈgorie $\tria \ca$.\\

Il s'agit de montrer que le foncteur $H^0y''$ est pleinement fidËle.
Nous devons donc montrer que, pour tous $M$, $M'$ objets de $\TW \ca$, nous avons
\[
H^0\Hom_{\TW \ca}(M,M') \arr{\sim} H^0\Hom_{\cc_\infty \ca}(\dot y''M,\dot y''M').
\]
Une extension $M$, donnÈe par une suite
\[
(M_r , \hdots ,M_i, \hdots ,M_l), \quad r \leq i \leq l,
\]
et une matrice $\delta^M$, est clairement filtrÈe dans la catÈgorie
des objets tordus $\TW \ca$ par
\[
F_k = (M_{r+k} , \hdots ,M_{l}), \quad 0 \leq k \leq l-r,
\]
(Le morphisme $\delta^M : M \ra M$ est compatible ‡ cette filtration). 
Les objets graduÈs de cette filtration sont des extensions tordues dÈgÈnÈrÈes,
i.~e.~des sommes formelles finies de $\Z \ca$ considÈrÈe comme objets de $\TW \ca$.
Il nous suffit donc de montrer qu'on a un isomorphisme
\[
H^0\Hom_{\TW \ca}(M,M') = H^0\Hom_{\cc_\infty \ca}(\dot y''M,\dot y''M').
\]
o˘ $M$ et $M'$ sont des objets de $\Z \ca$ considÈrÈs comme objets de $\TW \ca$.
Nous devons donc montrer le lemme suivant

\begin{lemme} \label{lemme_Yoneda_quasi-isomorphisme}
Pour toute paire d'objets $A$ et $A'$ dans $\ca$, 
l'$\ai$-foncteur de Yoneda $y : \ca \ra \cc_\infty \ca$ induit
un isomorphisme
\[
H^*\Hom_{\ca}(A,A') = H^*\Hom_{\cc_\infty \ca}(A\pw,A'\pw).
\]
\end{lemme}

\dem Le foncteur
pleinement fidËle (\ref{proposition_categorie_derivee_H-unitaire})
\[
\cd_\infty \ca \ra \cd_\infty \ca\+
\]
induit un isomorphisme
\[
H^*\Hom_{\cc_\infty \ca}(A\pw,A'\pw) \arr{\sim} H^*\Hom_{\cc_\infty \ca\+}(A\pw,A'\pw).
\]
Il suffit donc de montrer l'isomorphisme
\[
H^*\ca(A,A') \arr{\sim} H^*\Hom_{\cc_\infty \ca\+}(A\pw,A'\pw).
\]
Nous avons les ÈgalitÈs
\[
\ca(A,A') = A'\pw(A) \quad \mbox{et} \quad \Homi_{\ca}(\ca,A'\pw)(A) = \Hom_{\cc_\infty \ca\+}(A\pw,A'\pw).
\]
Nous dÈduisons alors le rÈsultat du lemme (\ref{lemme2_adjonction_foncteurs_standard}) et de la
remarque (\ref{remarque_lemme2_adjonction_foncteurs_standard}) qui
montrent que
\[
A'\pw \ra \Homi_{\ca}(\ca,A'\pw)
\]
est un quasi-isomorphisme.\findem

\section{ModËle diffÈrentiel graduÈ} \label{section_modele_differentiel_gradue}

Dans cette section, la catÈgorie de base $\sf C$ est Ègale ‡ $\sf C(\mathbb A,\mathbb A).$

\begin{definition}{\em
Soit $\ca$ une $\ai$-algËbre dans $\sf C.$
Un {\em modËle differentiel graduÈ} \index{modele differentiel gradue@{modËle diffÈrentiel graduÈ}}
de $\ca$ est une algËbre diffÈrentielle graduÈe $\ca'$ munie d'un $\ai$-quasi-isomorphisme
\[
\ca \ra \ca'.
\]
}\end{definition}

\begin{proposition}\label{proposition_modele_differentiel_gradue}
Toute $\ai$-algËbre strictement unitaire $\ca$ admet
un modËle diffÈrentiel graduÈ unitaire tel que l'$\ai$-morphisme
\[
\ca \ra \ca'
\]
est strictement unitaire.
\end{proposition}

Remarquons que dans le cas o˘ $\ca$ est une $\ai$-algËbre augmentÈe, son algËbre enveloppante $\E \ca$
(\ref{lemme_algebre_enveloppante}) est un modËle diffÈrentiel graduÈ unitaire de $\ca$ qui est augmentÈ.\\

\dem 
Nous dÈfinissons $\ca'$ comme le $\mathbb A$-$\mathbb A$-bimodule
\[
(A_0,A_1) \mapsto \Hom_{\cc_\infty \ca}(A_0\pw,A_1\pw).
\]
La structure diffÈrentielle graduÈe est celle induite par la composition et la diffÈrentielle
de la catÈgorie diffÈrentielle graduÈe $\cc_\infty \ca$.
Gr‚ce au thÈorËme (\ref{theoreme_factorisation_Yoneda}), l'$\ai$-foncteur de Yoneda
nous donne un $\ai$-quasi-isomor\-phisme d'$\ai$-algËbres dans $\sf C(\mathbb A,\mathbb A)$
\[
\ca \ra \ca'
\]
qui est strictement unitaire. \findem

\begin{corollaire}\label{corollaire_modele_differentiel_gradue}
Toute $\ai$-algËbre homologiquement unitaire $\ca$ admet
un modËle diffÈrentiel graduÈ unitaire tel que l'$\ai$-morphisme
\[
f : \ca \ra \ca'
\]
est unitaire, i.~e.~$f \circ \unite = \unite$.
\end{corollaire}

\dem Soit $\ca$ une $\ai$-algËbre homologiquement unitaire. 
Nous rappelons (\ref{corollaire_unite_ai-algebre}) qu'on peut munir $H^*\ca$
d'une structure d'$\ai$-algËbre strictement unitaire. Comme
l'$\ai$-morphisme
\[
\ca \ra H^*\ca
\]
est unitaire et il est $\ai$-quasi-isomorphisme, nous avons le rÈsultat. \findem

\section{CatÈgories stables}
\label{section_categories_stables}

Dans cette section, nous montrons que toute catÈgorie
triangulÈe algÈbrique qui est engendrÈe par un ensemble d'objets
est $\ai$-prÈ-triangulÈe, i.~e.~elle est Èquivalente
‡ $H^0\TW \ca$, pour une certaine $\ai$-catÈgorie $\ca$. 
\index{A-infini pretriangulee (categorie)@{$\ai$-prÈ-triangulÈe (catÈgorie)}}
\index{categorie@{catÈgorie}!A-infini pretriangulee@{$\ai$-prÈ-triangulÈe}}

\begin{definition}{\em \index{categorie@{catÈgorie}!triangulÈe algÈbrique}
Une $\corps$-catÈgorie triangulÈe est {\em algÈbrique} si elle est Èquivalente ‡ la catÈgorie stable d'une
$\corps$-catÈgorie de Frobenius (voir \ref{structure_triangulee_Ho_coModcu}). 
}\end{definition}

\begin{definition}{\em \index{compact (objet)}
Soit $\calt$ une catÈgorie triangulÈe aux sommes infinies. Un objet $X \in \calt$ est {\em compact}
si le foncteur $\Hom_{\calt}(X,\?)$ commute aux sommes infinies.
}\end{definition}

\begin{definition}{\em \index{generateur@{gÈnÈrateur}}
Soit $\calt$ une catÈgorie triangulÈe et $\mathbb A$ un sous-ensemble de
l'ensemble $\mathbb T$ des objets de $\calt$. \indexnotation{triabbA}
Nous notons $\tria \mathbb A$ la {\em plus petite sous-catÈgorie triangulÈe strictement
pleine}
de $\calt$ qui contient la sous-catÈgorie pleine formÈe des objets de $\mathbb A.$ 
Elle est stable par sommes finies. 
Les objets de $\mathbb A$ {\em engendrent} $\calt$ en tant que catÈgorie triangulÈe
si $\calt = \tria \mathbb A.$
Si $\calt$ admet des sommes infinies,
nous notons $\Tria \mathbb A$  {\em plus petite sous-catÈgorie triangulÈe stables par sommes infinies}
de $\calt$ qui contient la sous-catÈgorie pleine formÈe des objets de $\mathbb A.$
\indexnotation{TriabbA}
Les objets de $\mathbb A$ {\em engendrent} $\calt$ en tant que catÈgorie triangulÈe aux
sommes infinies si $\calt = \Tria \mathbb A.$
}\end{definition}

\begin{theoreme} \label{theoreme_caracterisation_categorie_derivee_sommes_infinies}
Soit $\calt$ une $\corps$-catÈgorie triangulÈe algÈbrique aux sommes infinies
qui est engendrÈe, en tant que catÈgorie triangulÈe aux sommes infinies, par un ensemble $\mathbb A$
d'objets compacts.
Il existe une $\ai$-catÈgorie $\ca$ strictement unitaire et minimale sur $\mathbb A$
et une Èquivalence triangulÈe
\[
\cd_\infty  \ca \ra \calt,  \quad A\pw \mapsto A.
\] 
\end{theoreme}

\dem Par dÈfinition des catÈgories triangulÈes algÈbriques,
$\calt$ est la catÈgorie stable $\underline \ce$ d'une catÈgorie de Frobenius $\ce$.
Nous rappelons \cite[4.3]{Keller94} qu'il existe une catÈgorie diffÈrentielle graduÈe unitaire $\ca'$
sur $\mathbb A$ et une Èquivalence triangulÈe
\[
\cd \ca' \ra \underline \ce,  \quad A\pw \mapsto A.
\]
Nous rappelons que $\cd \ca'$ est engendrÈe par les $A$-modules libres $\ca'(\?,A),$ $A \in \mathbb A$.
Choisissons un modËle minimal $\ca$ de $\ca'$ qui est strictement unitaire
(\ref{proposition_modele_minimal_strictement_unitaire}).
Nous dÈduisons du thÈorËme (\ref{theoreme_restriction_quasi-isomorphisme}) que
la restriction le long de $\ca' \ra \ca$ induit une Èquivalence de catÈgories
\[
\cd_\infty \ca \ra \cd \ca'.
\]
Comme, pour tout $A \in \mathbb A$, le $\ca'$-polydule restreint
$A \pw = \ca(\?,A)$ est $\ai$-quasi-isomorphe ‡ $\ca'(\?,A)$,
nous avons une Èquivalence
\[
\cd_\infty \ca \ra \underline \ce, \quad A\pw \mapsto A.
\]
\findem

\begin{remarque}{\em
Par construction de la catÈgorie $\ca'$ dans \cite[4.3]{Keller94}, le $\mathbb A$-$\mathbb A$-bimodule
sous-jacent ‡ $\ca$ est donnÈ par
\[
(A,A') \mapsto \ca(A,A') = \bigoplus_{n \in \Z} \Hom_\calt (A, S^nA'), \quad A,A' \in \mathbb A,
\]
et $m^\ca_2$ par la composition de $\calt.$
}\end{remarque}

\begin{theoreme} \label{theoreme_caracterisation_categorie_derivee}
Soit $\calt$ une $\corps$-catÈgorie triangulÈe algÈbrique qui est engendrÈe par un ensemble
d'objets $\mathbb A$.
Il existe une $\ai$-catÈgorie $\ca$ strictement unitaire et minimale sur $\mathbb A$
et une Èquivalence triangulÈe
\[
\tria \ca \ra \calt,  \quad A\pw \mapsto A,
\] 
o˘ $\tria \ca$ est la sous-catÈgorie de $\cd_\infty \ca$ engendrÈe par les
objets libres $A\pw,$ $A \in \mathbb A$.
\end{theoreme}

\dem Par dÈfinition des catÈgories triangulÈes algÈbriques, 
$\calt$ est la catÈgorie stable $\underline \ce$ d'une catÈgorie de Frobenius $\ce$.
La construction de \cite[4.3]{Keller94} nous donne une catÈgorie diffÈrentielle graduÈe unitaire $\ca'$
sur $\mathbb A$ telle que nous avons une Èquivalence triangulÈe
\[
\tria \ca' \ra \underline \ce,  \quad A\pw \mapsto A ,
\]
o˘ $\tria \ca'$ est la sous-catÈgorie de $\cd \ca'$ engendrÈe par les $A$-modules libres
$\ca'(\?,A),$ $A \in \mathbb A$.
Choisissons un modËle minimal $\ca$ de $\ca'$ qui est strictement unitaire
(\ref{proposition_modele_minimal_strictement_unitaire}).
L'Èquivalence de catÈgories
\[
\cd \ca' \ra \cd_\infty \ca.
\]
induit une Èquivalence
\[
\tria \ca' \ra \tria \ca
\] 
car le $\ca$-polydule
$A \pw = \ca(\?,A),$ $A \in \mathbb A$ est $\ai$-quasi-isomorphe ‡ la restriction de $\ca'(\?,A)$.
Nous en dÈduisons qu'on a une Èquivalence (triangulÈe)
\[
\tria \ca \ra \underline \ce, \quad A\pw \mapsto A.
\]
\findem

\begin{corollaire}
Soit $\calt$ une $\corps$-catÈgorie triangulÈe algÈbrique 
telle que dans le thÈorËme (\ref{theoreme_caracterisation_categorie_derivee}).
Il existe une $\ai$-catÈgorie $\ca$ strictement unitaire et minimale sur $\mathbb A$
et une Èquivalence triangulÈe
\[
H^0 (\TW \ca) \ra \calt,  \quad A \mapsto A.
\] 
\end{corollaire}

\dem ImmÈdiat par les thÈorËmes (\ref{theoreme_factorisation_Yoneda}) et 
(\ref{theoreme_caracterisation_categorie_derivee}). \findem

%% file: Ai_func.tex
{\bf \noindent Introduction}\\

\noindent 
Le but de ce chapitre est de construire l'analogue $\ai$ de la $2$-catÈgorie $\mathsf{cat}$
des petites catÈgories.
Nous construisons une $2$-catÈgorie $\aicat$ dont les
objets sont les $\ai$-catÈgories strictement unitaires.
La catÈgorie des espaces de morphismes
\[
\aicat (\ca,\cb), \quad \ca,\cb \in \Obj \aicat,
\]
sera dÈfinie comme l'homologie $H^0\Func (\ca,\cb)$ d'une $\ai$-catÈgorie dont les objets
sont les $\ai$-foncteurs strictement unitaires.

La catÈgorie $\Func (\ca,\cb)$ a ÈtÈ construite indÈpendamment par
K.~Fukaya \cite{Fukaya01a}, V.~Lyubashenko \cite{Lyubashenko02}  et
M.~Kontsevich et Y.~Soibelman \cite{Kontsevich02}, \cite{Kontsevich02a}.
Les compositions de $\Func(\ca,\cb)$ de V.~Lyubashenko,
bien qu'obtenues par une mÈthode diffÈrente, sont les mÍmes que
les nÙtres.\\

{\bf \noindent Plan du chapitre}\\

Soit $\ca$ et $\cb$ deux petites $\ai$-catÈgories (non nÈcessairement unitaires).
Dans la section \ref{section_Nunc}, nous construisons une $\ai$-catÈgorie
$\Nunc (\ca,\cb)$ dont les objets sont les $\ai$-foncteurs $\ca \ra \cb$.
Les compositions de $\Nunc(\ca,\cb)$ seront
construites par un processus de torsion (voir le chapitre \ref{chapitre_torsion}).
Dans la section \ref{section_fonctorialite_Nunc}, nous montrons que $\Nunc(\ca,\cb)$ est fonctoriel
en $\ca$ et $\cb$ et nous dÈfinissons la catÈgorie $\ainat$ dont les objets
sont les $\ai$-catÈgories.
Dans la section \ref{section_Func}, nous montrons que toutes les constructions des deux sections
prÈcÈdentes sont
compatibles aux $\ai$-structures strictement unitaires ($\ai$-catÈgories, $\ai$-foncteurs...)
et nous dÈfinissons la $2$-catÈgorie $\aicat$ comme une sous-catÈgorie non pleine de $\ainat$.

Dans la section (\ref{section_theorie_homotopie_ai-foncteurs}), nous construisons un $\ai$-foncteur
\[
z : \Func (\ca, \cb) \ra \cc_{\infty}(\ca,\cb), \quad \ca, \cb \in \aicat,
\]
o˘ $\cc_{\infty}(\ca,\cb)$ est la catÈgorie diffÈrentielle
graduÈe des $\ca$-$\cb$-bipolydules strictement unitaires 
(\ref{section_ai-foncteur_de_Yoneda_generalise}).
Ce foncteur gÈnÈralise l'$\ai$-foncteur de Yoneda construit en
(\ref{definition_ai-foncteur_de_Yoneda}). Nous montrerons que'il induit des quasi-isomorphismes
dans les espaces de morphismes. Dans la section \ref{section_equiv_faible_ai-foncteurs},
nous dÈfinissons les {\em Èquivalences faibles} d'$\ai$-foncteurs
strictement unitaires (elles sont
l'analogue $\ai$-catÈgorique des {\em homotopies} entre $\ai$-morphismes) et nous les caractÈriserons
‡ l'aide de leurs images par l'$\ai$-foncteur~$z$.

\section{L'$\ai$-catÈgorie des $\ai$-foncteurs}
\label{section_ai-categorie_des_ai-foncteurs}

\subsection{L'$\ai$-catÈgorie $\Nunc (\ca,\cb)$} \label{section_Nunc}

Soit $\mathbb A$ et $\mathbb B$ deux ensembles et $\ca$ et $\cb$
des $\ai$-catÈgories sur $\mathbb A$ et $\mathbb B$.
Nous construisons dans cette section, l'$\ai$-catÈgorie $\Nunc (\ca,\cb)$
des $\ai$-foncteurs non nÈcessairement strictement unitaires. La lettre
$\mathsf N$ remplace la lettre $\sf F$ dans $\Func$ et
se rapporte au $\mathsf N$ de ``$\mathsf N$on unitaires''.\\

Soit $f_1$ et $f_2$ deux $\ai$-foncteurs $\ca \ra \cb$.
Nous rappelons que $\lrus{\dot f_2}{\cb}{\dot f_1}$ est le $\mathbb A$-$\mathbb A$-bimodule
\[
(A',A) \mapsto \cb(\dot f_1A',\dot f_2A).
\]

\begin{definition}{\em \label{definition_morphismes_ai-foncteurs}
Nous posons
\[
\Hom_{\Nunc}(f_1,f_2) =
\Hom_{\gr \sf C(\mathbb A,\mathbb A)}\big(\ct S \ca, \lrus{\dot f_2}{\cb}{\dot f_1}\big).
\]
Nous obtenons ainsi un objet graduÈ dans la catÈgorie de base $\vect \corps$.
}\end{definition}

\begin{remarque}{\em \label{remarque_ai-categorie_Func}
Soit $H$ un ÈlÈment de degrÈ $r$ de $\Hom_{\Nunc}(f_1,f_2)$.
Pour tout entier $i\geq 0$, nous notons $\mathsf{incl}$ l'inclusion de $(S\ca)\tpo i$ dans $\ct S \ca$.
Soit $H_i$, $i\geq 0$, la composition
\[
(S\ca) \tpo i \arr{\mathsf{incl}} \ct S \ca \arr{h} 
\lrus{\dot f_2}{\cb}{\dot f_1}.
\]
Nous dÈfinissons les morphismes 
\[
h_i : \ca \tpo i \ra \lrus{\dot f_2}{\cb}{\dot f_1}, \quad i\geq 0,
\]
par les relations
\[
H_i \circ (\si \tpo i)^{-1} = (-1)^r  h_i, \quad i\geq 0.
\]
Les applications $H_i \mapsto h_i$, $i\geq 0$, sont clairement des bijections.
Le morphisme $H$ est donc dÈterminÈ par des morphismes graduÈs
\[
h_i : \ca(A_{i-1},A_i) \ts \hdots \ts \ca (A_0,A_1) \ra \cb(\dot f_1A_0,\dot f_2A_i), \quad i \geq 0.
\]
de degrÈ $r-i$, pour toute suite $(A_0,\hdots,A_i)$ d'objets de $\ca$.
En particulier, si $i = 0$, nous avons un morphisme
\[
h_0 : e_{\mathbb A} \ra \lrus{\dot f_2}{\cb}{\dot f_1}, \quad \id_A \mapsto h_0(\id_A).
\]
Nous noterons souvent $h_A \in \Hom_{\cb}(\dot f_1A,\dot f_2A)$ l'ÈlÈment $h_0(\id_A)$.
}\end{remarque}

\begin{remarque}{\em \label{remarque1_ai-categorie_Func}
Soit  $f : \ca \ra \cb$ un $\ai$-foncteur. Posons  $h_i = f_i$ si $i\geq 1$ et $h_0 = 0.$
Ceci nous fournit un ÈlÈment $H$ de degrÈ $+1$ de $\Hom_{\Nunc}(f,f)$. 
Nous avons alors un diagramme commutatif
\[
\xymatrix{(S\ca)\tpo i \ar@{^(->}[r] \ar[rd]_{H_i} & \Ba \ca \ar[r]^F 
\ar[d]^H & \Ba \lrus{\dot f}{\cb}{\dot f} \ar[d]^{p_1} \\
 &  \lrus{\dot f}{\cb}{\dot f} &  S\lrus{\dot f}{\cb}{\dot f} \ar[l]^{\si}
}
\]
dont nous dÈduisons les ÈgalitÈs
$H_i = \si \circ  F_i$, o˘ $F$ est la construction bar co-augmentÈe de $f$.
}\end{remarque}

{\noindent \bf Compositions naÔves de morphismes d'$\ai$-foncteurs}\\

Nous construisons dans ce paragraphe une $\ai$-catÈgorie $\cf(\ca,\cb) = \cf$ \indexnotation{cfcacb}
dont les objets sont les $\ai$-foncteurs $\ca \ra \cb$ et dont les espaces graduÈs de morphismes
sont les
\[
\Hom_{\cf}(f_1,f_2) =
\Hom_{\gr \sf C(\mathbb A,\mathbb A)}\big(\ct S\ca, \lrus{\dot f_2}{\cb}{\dot f_1}\big).
\]
Nous montrons que $\cf$ est munie d'une topologie pour laquelle elle est
une $\ai$-catÈgorie topologique. Nous construisons ensuite
un ÈlÈment tordant topologique de $\cf$ (voir \ref{section_torsion_II}).

Au lieu de construire les compositions $m_i^\cf$, $i\geq 1$, nous allons construire des
morphismes (voir les bijections $m_i \lra b_i$ dans la section \ref{construction_bar_cobar})
\[
b^\cf_i : S\cf \tpo i \ra S\cf, \quad i \geq 1,
\]
puis nous vÈrifions que cela dÈfinit bien une $\ai$-catÈgorie.
Remarquons que nous avons un isomorphisme
\[
S\cf (f_1,f_2) \arr{\sim} \Hom_{\gr \sf C(\mathbb A,\mathbb A)}\big(\Ba \ca,
S\lrus{\dot f_2}{\cb}{\dot f_1}\big).
\]
Le morphisme
\[
b^\cf_1 : \Hom_{\gr \sf C(\mathbb A,\mathbb A)}\big(\Ba \ca, S\lrus{\dot f_2}{\cb}{\dot f_1}\big)
\ra \Hom_{\gr \sf C(\mathbb A,\mathbb A)}\big(\Ba \ca, S\lrus{\dot f_2}{\cb}{\dot f_1}\big)
\]
est la diffÈrentielle de l'espaces de morphismes graduÈs entre complexes : elle est dÈfinie par
\[
 \varphi \mapsto b^\cb_1 \circ \varphi - (-1)^{|\varphi|} \varphi \circ b^{\Ba \ca},
\]
o˘ $\varphi$ est de degrÈ $|\varphi|$.
Soit $i\geq 2$ et $(f_0,\hdots, f_i)$ une suite d'$\ai$-foncteurs $\ca \ra \cb$. 
Le morphisme $b_i^\cf$ envoie un ÈlÈment
\[
g_i \tso \hdots \tso g_1 \in \Hom_{\gr \sf C(\mathbb A,\mathbb A)}
\big(\Ba \ca, S\lrus{\dot f_{i}}{\cb}{\dot f_{i-1}}\big)
\tso \hdots \tso \Hom_{\gr \sf C(\mathbb A,\mathbb A)}\big(\Ba \ca, S\lrus{\dot f_{1}}{\cb}{\dot f_0}\big)
\]
sur la composition
\[
\Ba \ca \arr{\Delta^{(i)}} (\Ba \ca)\tpo i \arr{g_i \tso \hdots \tso g_1} S\lrus{\dot f_{i}}{\cb}{\dot f_{i-1}}
\tso \hdots \tso S\lrus{\dot f_{1}}{\cb}{\dot f_{0}} \arr{b_i^\cb} S\lrus{\dot f_{i}}{\cb}{\dot f_{0}}.
\]

\begin{lemme}
Les morphismes $m^\cf_i$, $i \geq 1$, dÈfinissent une structure de $\ai$-catÈgorie sur $\cf$.
\end{lemme}

\dem
Nous avons clairement $b_1^\cf \circ b_1^\cf = 0.$ Soit $n \geq 2$ et soit
$g_i$, $1 \leq  i \leq n$, des ÈlÈments de $S\cf$ de degrÈ $|g_i|.$
Les termes de la somme
\[
\Big[\sum_{j+k+l = n}  b^{\cf}_i (\id \tpo j \tso  b^{\cf}_k \tso \id \tpo l)\Big]
(g_n \tso \hdots \tso g_1)
\]
sont des trois types : ceux o˘ $i=n$ et $k=1$, ceux o˘ $i=1$ et 
$k=n$ et ceux o˘ $i,j \neq 1.$ \\[.2cm]

\begin{itemize}
\item Lorsque $i=n$ et $k=1$ nous trouvons
\[
\begin{array}{cl}
& \big[  b^{\cf}_n (\id \tpo j \tso  b^{\cf}_1 \tso \id \tpo l)\big]
( g_n \tso \hdots \tso g_1) \\\\
= & (-1)^{\sum_{r < l+1}|g_r|}b^{\cf}_n 
( g_n \tso \hdots \tso b^{\cf}_1( g_{l+1}) \tso \hdots \tso g_1)\\\\
= & (-1)^{\sum_{r < l+1}|g_r|}b^{\cf}_n 
( g_n \tso \hdots \tso b^{\cb}_1 g_{l+1} \tso \hdots \tso g_1)\\
& - (-1)^{\sum_{r \leq l+1}|g_r|}b^{\cf}_n 
( g_n \tso \hdots \tso  g_{l+1} b^{\Ba \ca} \tso \hdots \tso g_1)\\\\
= & b^{\cb}_n (\id \tpo j \tso  b^{\cb}_1 \tso \id \tpo l)( g_n \tso \hdots \tso g_1)\Delta^{(n)}\\
& - (-1)^{\sum_{r}|g_r|}b^{\cb}_n 
( g_n \tso \hdots \tso  g_{l+1} \tso \hdots \tso g_1)(\id \tpo j \tso  b^{\Ba \ca} \tso \id \tpo l)\Delta^{(n)}\\\\
= & b^{\cb}_n (\id \tpo j \tso  b^{\cb}_1 \tso \id \tpo l)( g_n \tso \hdots \tso g_1)\Delta^{(n)}\\
& - (-1)^{\sum_{r}|g_r|}b^{\cb}_n 
( g_n \tso \hdots \tso  g_{l+1} \tso \hdots \tso g_1)\Delta^{(n)} b^{\Ba \ca}\\\\
= & b^{\cb}_n (\id \tpo j \tso  b^{\cb}_1 \tso \id \tpo l)( g_n \tso \hdots \tso g_1)\Delta^{(n)}\\
& - (-1)^{\sum_{r}|g_r|}b^{\cf}_n 
( g_n \tso \hdots \tso  g_{l+1} \tso \hdots \tso g_1) b^{\Ba \ca}
\end{array}
\]

\item  Lorsque $i=1$ et $k=n$ nous trouvons
\[
\begin{array}{cl}
&  b^{\cf}_1 \cdot b^{\cf}_n ( g_n \tso \hdots \tso g_1) \\\\
= & b^{\cb}_1 (b^{\cf}_n ( g_n \tso \hdots \tso g_1))\\
  & - (-1)^{1+\sum_{r}|g_r|}b^{\cf}_n ( g_n \tso \hdots g_{l+1} \tso \hdots \tso g_1) b^{\Ba \ca}\\\\
= & b^{\cb}_1 (b^{\cb}_n ( g_n \tso \hdots \tso g_1)\Delta^{(n)})\\
  & - (-1)^{1+\sum_{r} |g_r|}b^{\cf}_n  ( g_n \tso \hdots g_{l+1} \tso \hdots \tso g_1) b^{\Ba \ca}
\end{array}
\]

\item Lorsque $i\neq 1$ et $k \neq n$ nous trouvons
\[\hspace{-1cm}
\begin{array}{cl}
& \big[  b^{\cf}_i (\id \tpo j \tso  b^{\cf}_k \tso \id \tpo l)\big]
( g_n \tso \hdots \tso g_1) \\\\
= & (-1)^{\sum_{r < l+1}|g_r|}b^{\cf}_i 
( g_n \tso \hdots b^{\cf}_j( g_{l+j+1} \tso \hdots \tso  g_{l+1}) \tso \hdots \tso g_1)\\\\
= & (-1)^{\sum_{r < l+1}|g_r|}b^{\cf}_i 
( g_n \tso \hdots (b^{\cb}_j( g_{l+j+1} \tso \hdots \tso  g_{l+1})\Delta^{(k)}) \tso \hdots \tso g_1)\\\\
= & (-1)^{\sum_{r < l+1}|g_r|}b^{\cb}_i 
( g_n \tso \hdots (b^{\cb}_j( g_{l+j+1} \tso \hdots \tso  g_{l+1})\Delta^{(k)}) \tso \hdots \tso g_1)
\Delta^{(i)}\\\\
= & b^{\cb}_i  (\id \tpo j \tso b^{\cb}_j \tso \id \tpo l)
( g_n \tso \hdots  \tso g_1)
(\id \tpo j \tso \Delta^{(k)} \tso \id \tpo l)\Delta^{(i)}\\\\
= & b^{\cb}_i (\id \tpo j \tso b^{\cb}_j \tso \id \tpo l)
( g_n \tso  \hdots \tso g_1)
\Delta^{(n)}
\end{array}
\]
\end{itemize}

Les derniËres lignes des deux premiers cas se compensent gr‚ce aux signes et
la somme de ce qui reste est nulle car $\cb$ est une $\ai$-catÈgorie. \findem

\begin{remarque}{\em
L'$\ai$-catÈgorie $\cf(\ca,\cb)$ ainsi construite est clairement fonctorielle en $\ca$ et $\cb.$
Si $f : \ca \ra \ca'$ est un $\ai$-foncteur, l'$\ai$-foncteur induit $\cf(\ca',\cb) \ra \cf(\ca,\cb)$
est strict. Il envoit $H \in \Hom_{\Nunc}(f_1,f_2)$ sur sa composition avec $Bf$.
Si $f : \cb \ra \cb'$ est un $\ai$-foncteur,  l'$\ai$-foncteur induit $\cf(\ca',\cb) \ra \cf(\ca,\cb)$
n'est plus strict. Notons le $g$. 
Soit $G$ sa construction bar. Le morphisme $G_1$ envoie $H \in \Hom_{\Nunc}(f_1,f_2)$
sur sa composition avec $F_1.$
Les formules dÈfinissant les $G_i$, $i\geq 2$, sont obtenues ‡ partir des formules dÈfinissant les
$b_i^\cf$, $i \geq 2$, en remplacant les $b^\cb_i$ par des $F_i$. Les problËmes de fonctorialitÈ seront
ÈtudiÈs plus
en dÈtails dans la section \ref{section_fonctorialite_Nunc}.
}\end{remarque}

{\noindent \bf Description concrËte}\\

Regardons ce que sont les compositions de morphismes d'$\ai$-catÈgories du point
de vue de la remarque \ref{remarque_ai-categorie_Func}. \\

Soit $H$ un ÈlÈment de $\Hom_\Nunc (f_1,f_2)$ de degrÈ $|H|$.
Le morphisme $m^\cf_1(H)$ est dÈterminÈ par des morphismes
\[
h'_i : \ca \tpo i \ra \lrus{f_2}{\cb}{f_1}, \quad i \geq 0.
\]
Nous vÈrifions que $h_i'$ est Ègal ‡ la somme
\[
m^\cb_1 \circ h_i - (-1)^{|H|}\sum (-1)^{l+kj} h_{j+1+l}(\Id \tpo j \tso m^{\ca}_k \tso \Id \tpo l).
\]

Soit $i\geq 2$.
Soit $f_0, \hdots ,f_i$ des $\ai$-foncteurs $\ca \ra \cb$. Pour tout $1 \leq t \leq i$, soit
$H_t$ un ÈlÈment de $\Hom_{\Nunc}(f_{t-1},f_t)$ de degrÈ $|H_t|$.
Notons $|H|$ la somme des degrÈs $|H_t|.$
Soit $H'$ l'ÈlÈment de $\Hom_\Nunc (f_0,f_i)$ Ègal ‡ $m_i^\cf (H_n \tso \hdots \tso H_1)$.
Alors $H'$ est donnÈ par des morphismes graduÈs
\[
h'_n : \ca(A_{n-1},A_n) \ts \hdots \ts \ca (A_0,A_1) \ra \cb(\dot f_0A_0,\dot f_iA_n), \quad n \geq 0.
\]
de degrÈ $|H|-n$, pour toute suite $(A_0,\hdots,A_n)$ d'objets de $\ca$.
Soit $x_k \in \ca(A_{k-1},A_k)$, $1 \leq k \leq n.$
Nous notons $\mathsf{incl}$ l'inclusion de $(S\ca)\tpo i$ dans $\Ba \ca$.
L'ÈlÈment $h'_n (x_n \tso \hdots \tso x_1)$ est Ègal ‡
\[
\begin{array}{cl}
  - \si \circ b^\cb_i \circ \big[(\si\tpo i)^{-1}(H_i \tso \hdots \tso H_1)\big] \circ \Delta^{(i)}
\circ \mathsf{incl} \circ (\si\tpo n)^{-1}  (x_n \tso \hdots \tso x_1)
\end{array}
\]
Prenons un exemple simple.

\begin{exemple}{\em \label{exemple}
Supposons que $i = 3$ et $n = 2$.
La composition $\Delta^{(3)}
\circ \mathsf{incl} \circ (\si\tpo 2)^{-1} (x_2\tso x_1)$ est Ègale ‡ la somme dans $\Ba \ca \tpo 3$
\[
\begin{array}{l}
\big[ \id_{A_2} \tso \id_{A_2} \tso (\si \tpo 2)^{-1} -  \id_{A_2} \tso (\si)^{-1} \tso (\si)^{-1} -\\
\hspace*{1cm} (\si)^{-1} \tso \id_{A_1} \tso (\si)^{-1} + \id_{A_2} \tso (\si \tpo 2)^{-1} \tso \id_{A_0}\\
\hspace*{1.5 cm} - (\si)^{-1} \tso (\si)^{-1} \tso \id_{A_0} 
+ (\si \tpo 2)^{-1} \tso \id_{A_0} \tso \id_{A_0}\big](x_2 \tso x_1).
\end{array}
\]
Nous trouvons donc que $m^\cf_3(h_3 \tso h_2 \tso h_1)(x_2 \tso x_1)$ est Ègal ‡ la somme des
ÈlÈments
\[
\begin{array}{l}
m^\cb_3 \Big(\pm (h_3)_{A_2} \tso (h_2)_{A_2} \tso (h_1)_2 \pm (h_3)_{A_2} \tso (h_2)_1 \tso (h_1)_1 \\
\hspace*{1cm} \pm (h_3)_{1} \tso (h_2)_{A_1} \tso (h_1)_{1} \pm (h_3)_{A_2} \tso (h_2)_{2} \tso (h_1)_{A_0} \\
\hspace*{1.5cm} \pm (h_3)_{1} \tso (h_2)_{1} \tso (h_1)_{A_0} \pm
(h_3)_{2} \tso (h_2)_{A_0} \tso (h_1)_{A_0}\Big)(x_2 \tso x_1).
\end{array}
\]
Le morphisme
\[
h'_2(x_2 \tso x_1) : f_0 A_0 \ra f_3 A_2
\]
est donc la somme des compositions (au signes prËs) des suites de morphismes reprÈsentÈes par
un chemin de flËches menant de $f_0A_0$ ‡ $f_3A_2$ dans le diagramme ci-dessous 
\begin{center}
\hspace{1cm} \input{figure1.pstex_t}\hspace{1.5cm}.
\end{center}
Remarquons qu'il n'y a aucune flËche verticale (qui correspondrait ‡ un $(f_j)_1(x_i)$ ou un
$(f_j)_2(x_2 \ts x_1)$) dans ces chemins de flËches.
}\end{exemple}

De faÁon gÈnÈrale, nous trouvons que l'ÈlÈment $H'$ de $\Hom_{\Nunc}(f_0,f_n)$ est donnÈ par 
\[
h'_n = \sum_{j_1 + .. + j_l = n} (-1)^{s} m^{\cb}_l\big((h_i)_{j_1} \tso \hdots \tso (h_1)_{j_l} \big), \quad n\geq 0,
\]
o˘ les entiers $j_\alpha$ sont $\geq 0$, et o˘ le signe est donnÈ par l'ÈgalitÈ
\[
(-1)^{s} \big((H_i)_{j_1} \tso \hdots \tso (H_1)_{j_l} \big) \circ (\si \tpo n)
= \big((h_i)_{j_1} \tso \hdots \tso (h_1)_{j_l} \big).
\]

\begin{remarque}{\em
Soit $H$ l'ÈlÈment de $\Hom_\Nunc (f, f)$ construit ‡ la remarque (\ref{remarque1_ai-categorie_Func}).
Si $f_t = f$, $0 \leq t \leq i$, et $H_t = H$, $ 1 \leq t \leq i$,
le signe $(-1)^s$ ci-dessus est le mÍme que le signe
$(-1)^s$ de l'Èquation $(**_n)$, $n\geq 1$, dans la dÈfinition
des $\ai$-foncteurs (\ref{definition_ai-foncteurs}). 
}\end{remarque}

{\noindent \bf Topologie sur $\cf$}\\

Nous munissons l'espace
\[
\Hom_{\cf}(f_1,f_2) = \Hom_{\gr \sf C(\mathbb A, \mathbb A)}(\Ba \ca, \lrus{\dot f_2}{\cb}{\dot f_1})
\]
de la topologie dÈfinie par la filtration dÈcroissante $F_i$, $i \geq 0$, o˘
\[
F_i = \Hom_{\gr \sf C(\mathbb A, \mathbb A)}\Big(\bigoplus_{j\geq i}(S\ca)\tpo j,
\lrus{\dot f_2}{\cb}{\dot f_1}\Big).
\]
Cette topologie est sÈparÈe. La description ci-dessus montre que
les compositions de
$\cf$ sont des morphismes continus contractants (voir \ref{section_definitions_topologie}).
L'$\ai$-catÈgorie $\cf$ est donc topologique (\ref{definition_ai-algebres_topologiques}).\\

{\noindent \bf ElÈment tordant de $\cf$}\\

Notons $\mathbb F$ l'ensemble des $\ai$-foncteurs $\ca \ra \cb.$
L'ÈlÈment tordant
\[
x : e_{\mathbb F} \ra \cf
\]
envoie le gÈnÈrateur $\id_{f}$ de $e_{\mathbb F}(f,f)$ sur l'ÈlÈment
$H$ de degrÈ $+1$ de $\Hom_\Nunc (f,f)$ construit ‡ partir de $f$ (voir \ref{remarque1_ai-categorie_Func}).

VÈrifions maintenant que $x$ est un ÈlÈment tordant topologique.
Comme le morphisme $h_0$ est nul, l'image de $x$ est dans le voisinage $\cf_1$.
La restriction de la somme
\[
\sum_{i \geq 1} m_i^{\cf} (H\tpo i)(\id_{f}) : \Ba \ca \ra \lrus{\dot f}{\cb}{\dot f}
\]
‡ $(S\ca)\tpo n$ est la somme
\[
- \sum (-1)^{jk+l}h_i(\Id^{\ts j}\ts m_k\ts \Id^{\ts l}) +
\sum_{j_1 + .. + j_l = n} (-1)^{s} m^{\cb}_l\big(h_{j_1} \tso \hdots \tso h_{j_l})
\]
Rappelons que $h_i = f_i$, $i\geq 1$. L'Èquation de Maurer-Cartan appliquÈe ‡ $\id_{f}$
est donc Èquivalente
‡ l'ensemble de Èquations $(**_n)$, $n \geq 1$, de la dÈfinition d'un $\ai$-foncteur
(\ref{definition_ai-foncteurs}).\\

{\noindent \bf L'$\ai$-catÈgorie $\Nunc (\ca,\cb)$}

\begin{definition}[Voir aussi \cite{Fukaya01a}, \cite{Lyubashenko02}  et
\cite{Kontsevich02}, \cite{Kontsevich02a}]
{\em  \label{definition_ai-cat_Nunc} \indexnotation{Nunc}
$\hspace{1cm}$ \newline L'{\em $\ai$-catÈgorie} $\Nunc (\ca, \cb)$ est la catÈgorie tordue $\cf_x$
(voir \ref{definition_torsion_ai-categories_II} pour la torsion).
}\end{definition}

Remarquons que les compositions $m^{\Nunc}_i$, $i\geq 1$, de \cite{Fukaya01a}, \cite{Lyubashenko02} sont
les mÍmes mais obtenues de maniËre diffÈrentes.\\

{\noindent \bf Description concrËte}\\

Donnons maintenant une description du morphisme 
\[
m^{\Nunc}_1 : \Hom_\Nunc (f_1,f_2) \ra \Hom_\Nunc (f_1,f_2).
\]
Soit $H$ un ÈlÈment de degrÈ $|H|$ de $\Hom_\Nunc (f_1,f_2).$
Le morphisme $H' = m^{\Nunc}_1(H)$ est dÈterminÈ par des morphismes
\[
h'_i : \ca \tpo i \ra \lrus{f_2}{\cb}{f_1}, \quad i \geq 0.
\]
Nous vÈrifions que $h'_i$ est Ègal ‡ la somme
\[
\begin{array}{l}
\sum_{j_1 + .. + j_l = n} (-1)^{s} m^\cb_l \big((f_2)_{j_1} \tso \hdots \tso (f_2)_{j_t} \tso h_{j_{t+1}}\tso (f_1)_{j_{t+1}}
\hdots \tso (f_1)_{j_l} \big) \\
\hspace*{2cm} - (-1)^{|h|+l+kj} h_{j+1+l}(\Id \tpo j \tso m^{\ca}_k \tso \Id \tpo l),
\end{array}
\]
o˘ l'exposant du signe $s$ est la somme du signe apparaissant dans la torsion
(\ref{section_torsion_ai-categories_I}) et du signe donnÈ par l'ÈgalitÈ
\[
\begin{array}{l}
(-1)^{*}  \big((\si F_2)_{j_1} \tso \hdots \tso (\si F_2)_{j_t} \tso H_{j_{t+1}}\tso (\si F_1)_{j_{t+1}} 
\hdots \tso (\si F_1)_{j_l} \big) \circ (\si \tpo n)\\
\hspace*{2cm}= \big((f_2)_{j_1} \tso \hdots \tso (f_2)_{j_t} \tso h_{j_{t+1}}\tso (f_1)_{j_{t+1}} \hdots \tso (f_1)_{j_l} \big).
\end{array}
\]

La description des compositions supÈrieures $m^{\Nunc}_i$, $i\geq 2$, se fait de faÁon similaire.
Reprenons l'exemple \ref{exemple} et posons
\[
H'' = m_3^{\Nunc}(h_3 \tso h_2 \tso h_1) \in \Hom_{\Nunc}(f_0,f_3).
\]
Le morphisme
\[
h''_2(x_2 \tso x_1) : f_0 A_0 \ra f_3 A_2
\]
est la somme des compositions (au signes prËs) des suites de morphismes reprÈsentÈes par
un chemin de flËches menant de $f_0A_0$ ‡ $f_3A_2$ dans le diagramme ci-dessous 
\begin{center}
\hspace{1cm} \input{figure2.pstex_t}\hspace{1.5cm}.
\end{center}
Graphiquement, la torsion consiste donc ‡ autoriser les flËches verticales dans les
chemins.

\begin{remarque}{\em
Si $\cb$ est une catÈgorie diffÈrentielle graduÈe, la catÈgorie $\Nunc(\ca,\cb)$ est
aussi une catÈgorie diffÈrentielle graduÈe car les compositions $m^{\Nunc}_i$, $i\geq 3$
sont nulles.
}\end{remarque}

\subsection{FonctorialitÈ de $\Nunc(\ca,\cb)$} \label{section_fonctorialite_Nunc}

{\noindent \bf FonctorialitÈ en $\ca$}\\

Soit $\ca$, $\ca'$, $\cb$ des petites $\ai$-catÈgories.
Soit $g \in \ca' \ra \ca$, $f_1, f_2 : \ca \ra \cb$ des $\ai$-foncteurs. 
Soit $H$ un ÈlÈment de $\Hom_\Nunc (f_1,f_2)$.
Nous dÈfinissons l'ÈlÈment
\[
H \star g \in \Hom_\Nunc((f_1 \circ g),(f_2 \circ g))
\]
comme la composition
\[
\Ba \ca' \arr{G} \Ba \lrus{\dot g}{\ca}{\dot g} \ra  \lrus{\dot f_2 \dot g}{\cb}{\dot f_1 \dot g}
\]
o˘ la seconde flËche est induite par $H.$ Comme $G$ est un morphisme de cogËbres diffÈrentielles
graduÈes, le morphisme de $\mathbb F$-$\mathbb F$-bimodules
\[
? \star g : \Nunc (f_1,f_2) \ra \Nunc((f_1 \circ g),(f_2 \circ g))
\]
est un $\ai$-foncteur strict.\\

{\noindent \bf FonctorialitÈ en $\cb$}\\

Soit $\ca$, $\cb$ et $\cb'$ des petites $\ai$-catÈgories.
Soit $g \in \cb \ra \cb'$, $f_1, f_2 : \ca \ra \cb$ des $\ai$-foncteurs.
Soit $H$ un ÈlÈment de $\Hom_\Nunc (f_1,f_2)$. Nous allons construire un ÈlÈment
\[
g \star H \in \Hom_\Nunc((g \circ f_1),(g \circ f_2)).
\]
Cela nous fournira un $\ai$-foncteur strict
\[
g \, \star \, ? : \Nunc (f_1,f_2) \ra \Nunc((g \circ f_1),(g \circ f_2))
\]
CommenÁons par introduire quelques notions.

Soit $M$ un $\mathbb A$-$\mathbb A$-bimodule diffÈrentiel graduÈ.
Soit $C$, $C_1$ et $C_2$ des cogËbres cocomplËtes
dans la catÈgorie des cogËbres diffÈrentielles graduÈes de la catÈgorie de base
$\sf C(\mathbb A,\mathbb A).$ Nous munissons le $\mathbb A$-$\mathbb A$-bimodule
$C_2 \tso M \tso C_1$ de la structure de $C_2$-$C_1$-bicomodule (cocomplet) induite par les
comultiplications de $C_2$ et $C_1$. Soit 
\[
F_1 : C \ra C_1 \quad \mbox{et} \quad F_2 : C \ra C_2
\]
des morphismes de cogËbres.

\begin{definition}{\em  \index{coderivation@{codÈrivation}}
Une {\em ($F_1$,$F_2$)-codÈrivation} est un morphisme de $\mathbb A$-$\mathbb A$-bimo\-dules
\[
K : C \ra C_2 \tso M \tso C_1
\]
tel que
\[
(\Delta^{C_2} \tso \Id \tso \Id) \circ K = (F_2 \tso K) \circ \Delta^C \quad \mbox{et} \quad
(\Id \tso \Id \tso \Delta^{C_1}) \circ K = (K \tso F_1) \circ \Delta^C.  
\]
}\end{definition}

\begin{lemme}
Soit $p_1$ la projection $C_2 \tso M \tso C_1$ sur $M.$
L'application $K \circ p_1 \circ K$ est une bijection de l'ensemble des $(F_1,F_2)$-codÈrivations
sur les morphismes de  $\mathbb A$-$\mathbb A$-bimodules $C \ra M.$
\findem
\end{lemme}
Soit $C_1$, $C_2$ et $C_3$ des cogËbres cocomplËtes
dans la catÈgorie des cogËbres diffÈrentielles graduÈes de la catÈgorie de base
$\sf C(\mathbb A,\mathbb A).$ Le {\em produit cotensoriel} d'un $C_1$-$C_2$-bicomodule $M$
avec un $C_2$-$C_3$-bicomodule $N$ est le noyau \indexnotation{boxdot}
\[
M \boxdot N = \ker \Big(M \tso N \arr{\Delta \tso \Id - \Id \tso \Delta} M \tso C_2 \tso N\Big).
\]

Reprenons la construction de $H \star g.$
Nous rappelons que les $\mathbb A$-$\mathbb A$-bimodules
$\lrus{\dot f_1}{\cb}{\dot f_1}$ et $\lrus{\dot f_2}{\cb}{\dot f_2}$ sont des 
$\ai$-catÈgories sur $\mathbb A$.
Soit
\[
F_1 : \Ba \ca \ra \Ba \lrus{\dot f_1}{\cb}{\dot f_1} \quad \mbox{et} \quad
F_2 : \Ba \ca \ra \Ba \lrus{\dot f_2}{\cb}{\dot f_2}
\]
les constructions bar co-augmentÈes de $f_1$ et $f_2.$
Le morphisme
\[
H : \Ba \ca \ra \lrus{\dot f_2}{\cb}{\dot f_1}
\]
se relËve en une $(F_1,F_2)$-codÈrivation de comodules
\[
K : \Ba \ca \ra \Ba \lrus{\dot f_2}{\cb}{\dot f_2} \tso \lrus{\dot f_2}{\cb}{\dot f_1} \tso \Ba
\lrus{\dot f_1}{\cb}{\dot f_1}.
\]
L'$\ai$-foncteur $g : \cb \ra \cb'$ induit un morphisme $G$ de degrÈ $0$
\[
 \Ba \lrus{\dot f_2}{\cb}{\dot f_2} \tso \lrus{\dot f_2}{\cb}{\dot f_1} \tso \Ba
\lrus{\dot f_1}{\cb}{\dot f_1} \ra
\Ba \lrus{\dot g \dot f_2}{\cb}{\dot g \dot f_2} \tso \lrus{\dot g \dot f_2}{\cb}{\dot g \dot f_1} \tso \Ba
\lrus{\dot g \dot f_1}{\cb}{\dot g \dot f_1}.
\]
Nous vÈrifions que la composition $G \circ K$ dÈfinit une $(GF_1,GF_2)$-codÈrivation
\[
\Ba \ca \ra \Ba \lrus{\dot g \dot f_2}{\cb}{\dot g \dot f_2} \tso
\lrus{\dot g \dot f_2}{\cb}{\dot g \dot f_1} \tso \Ba
\lrus{\dot g \dot f_1}{\cb}{\dot g \dot f_1}
\]
et nous dÈfinissons l'ÈlÈment $g \star H$ par la composition
\[
p_1 \circ (G \circ K) : \Ba \ca \ra \lrus{\dot g \dot f_2}{\cb}{\dot g \dot f_1}.
\]

Munissons  $\Ba \lrus{\dot f_2}{\cb}{\dot f_2} \tso \lrus{\dot f_2}{\cb}{\dot f_1} \tso \Ba
\lrus{\dot f_1}{\cb}{\dot f_1}$
de la diffÈrentielle induite par les $b_i^\cb,$ $i\geq 1,$ et notons $D(f_2,f_1)$
le bicomodule diffÈrentiel graduÈ obtenu ainsi.
Nous pouvons considÈrÈ $D(f_1,f_2)$ comme la construction bar du
$\lrus{\dot f_2}{\cb}{\dot f_2}$-$\lrus{\dot f_1}{\cb}{\dot f_1}$-bipolydule 
$\lrus{\dot f_2}{\cb}{\dot f_1}.$

\begin{remarque}{\em
Soit $H$ un ÈlÈment de $\Hom_\Nunc (f_1,f_2)$ et $K$ la codÈrivation associÈe.
L'ÈlÈment $m^\Nunc_1 (H)$ correspond ‡ la codÈrivation $\delta (K)$
dans l'espace diffÈrentiel graduÈ des morphismes graduÈs
\[
\Big( \Hom_{\gr \sf C(\mathbb A,\mathbb A)}
\big(\Ba \ca, D(f_2,f_1)\big), \delta\Big).
\]
}\end{remarque}

Soit $i\geq 2$.
Soit $f_0, \hdots ,f_i$ des $\ai$-foncteurs $\ca \ra \cb$. Pour tout $1 \leq t \leq i$, soit
$H_t$ un ÈlÈment de $\Hom_{\Nunc}(f_{t-1},f_t)$ de degrÈ $|H_t|$.
Notons $C_t$ la cogËbre diffÈrentielle graduÈe $\Ba \lrus{\dot f_t}{\cb}{\dot f_t}$.
Le $C_i$-$C_0$-bicomodule
\[
D(f_i,f_{i-1}) \boxdot \hdots \boxdot D(f_1,f_{0})
\]
est isomorphe en tant qu'objet graduÈ ‡
\[
C_i \tso \lrus{\dot f_i}{\cb}{\dot f_{i-1}} \tso 
C_{i-1}\tso \lrus{\dot f_{i-1}}{\cb}{\dot f_{i-2}} \tso C_{i-2} \tso \hdots \tso 
C_1\tso \lrus{\dot f_1}{\cb}{\dot f_0} \tso 
C_1.
\]
Nous le munissons de la diffÈrentielle induite par les $b_i^\cb,$ $i\geq 1.$
L'ÈlÈment
\[
m_i (H_i \tso \hdots \tso H_1) : \Ba \ca \ra \lrus{\dot f_i}{\cb}{\dot f_{0}}
\]
correspond ‡ la $F_i$-$F_1$-codÈrivation
\[
K : \Ba \ca \ra D(f_i,f_0)
\]
qui est le relËvement
\[
\Ba \ca \arr{\Delta^{(i)}} (\Ba \ca)\tpo i \arr{K_i \boxdot \hdots \boxdot K_1}
D(f_i,f_{i-1}) \boxdot \hdots \boxdot D(f_1,f_{0})
\arr{q} \lrus{\dot f_i}{\cb}{\dot f_0},
\]
o˘ $q$ est induit par les $b_i^\cb,$ $i\geq 1.$

L'$\ai$-foncteur $g$ induit des morphismes
\[
D(f_i,f_{i-1}) \boxdot \hdots \boxdot D(f_1,f_{0}) \ra D(gf_i,gf_{i-1}) \boxdot \hdots \boxdot D(gf_1,gf_{0})
\]
et un relËvement vers $D(gf_i,gf_0)$ de
\[
D(f_i,f_{i-1}) \boxdot \hdots \boxdot D(f_1,f_{0})
\arr{} \lrus{\dot g\dot f_i}{\cb}{\dot g\dot f_0}
\]
qui sont compatibles aux diffÈrentielles.
Nous en dÈduisons que le morphisme de $\mathbb F$-$\mathbb F$-bimodules
\[
g \, \star \,? : \Nunc (f_1,f_2) \ra \Nunc((g \circ f_1),(g \circ f_2))
\]
dÈfinit un $\ai$-foncteur strict.\\

{\noindent \bf La catÈgorie $\ainat$}\\

Soit $\ainat$ \indexnotation{ainat}
la {\em catÈgorie} dont les objets sont les petites $\ai$-catÈgories
(non nÈcessairement strictement unitaires), dont les espaces
de morphismes sont les catÈgories (sans unitÈs en gÈnÈral)
\[
\ainat (\ca,\cb) = H^0\Nunc(\ca,\cb).
\]
Il rÈsulte de la fonctorialitÈ de $\Nunc(\ca,\cb)$ que $\ainat$ est une ``$2$-catÈgorie
sans unitÈs pour les $2$-morphismes''.
La lettre $\mathsf n$ remplace la lettre $\mathsf c$ de $\aicat$ et exprime
le fait que les objets de $\ainat$ sont les $\ai$-``$\mathsf{cat}$''Ègories ``${\bf \mathsf n}$''on 
(nÈcessairement) strictement unitaires.

\begin{remarque}{\em
Soit $f_1$ et $f_2$ in $\Nunc(\ca,\cb)$. Soit $H$ un morphisme de
$\Hom_\Nunc(f_1,f_2)$ qui est un zÈro cycle. Soit $x$ un ÈlÈment de
$\ca (A_0,A_1)$. Comme $H$ est un cycle, nous avons la relation
\[
m^\cb_1(h_1(x)) - m^\cb_2(h_{A_1} \tso f_1x) +
m^\cb_2( f_2x \tso h_{A_0}) = 0.
\]
Nous avons donc un diagramme commutatif dans $H^0\cb$
\[
\xymatrix{
\dot f_1 A_0 \ar[r]^{f_1x} \ar[d]_{h_{A_0}}  & \dot f_1 A_1 \ar[d]^{h_{A_1}}\\
f_2A_0 \ar[r]_{f_2x} & f_2A_1.
}
\]
}\end{remarque}

\subsection{L'$\ai$-catÈgorie $\Func (\ca,\cb)$} \label{section_Func}

Reprenons les notations de la section \ref{section_Nunc} mais supposons dÈsormais
que $\ca$ et $\cb$ sont strictement unitaires. L'$\ai$-catÈgorie
$\Func (\ca,\cb)$ dont les objets sont les $\ai$-foncteurs strictement unitaires
est dÈfinie de la maniËre suivante :

Soit $f_1$ et $f_2$ deux $\ai$-foncteurs $\ca \ra \cb.$
Un ÈlÈment $H$ de $\Hom_{\Nunc}(f_1,f_2)$ est {\em strictement unitaire}
\index{strictement unitaire!element de $\Hom_{\Nunc}(f_1,f_2)$@{ÈlÈment de $\Hom_{\Nunc}(f_1,f_2)$}}
s'il vÈrifie
\[
h_i (\Id \tpo \alpha \tso \unite \tso \Id \tpo \beta) = 0, \quad i\geq 1. 
\]
Les $\ai$-foncteurs strictement unitaires et les morphismes strictement unitaires d'$\ai$-foncteurs
strictement unitaires forment une sous-$\ai$-catÈgorie de $\Nunc(\ca,\cb).$
Nous la notons $\Func(\ca,\cb).$ \indexnotation{Func}
Nous vÈrifions que $\Func(\ca,\cb)$ est fonctoriel
par rapport aux $\ai$-foncteurs strictement unitaires.\\

{\noindent \bf La 2-catÈgorie $\aicat$}\\

\begin{definition}{\em \label{definition_aicat} \indexnotation{aicat}
Soit $\aicat$ la {\em catÈgorie} dont les objets sont les petites $\ai$-catÈgories
strictement unitaires, dont les espaces
de morphismes sont les catÈgories 
\[
\aicat (\ca,\cb) = H^0\Func(\ca,\cb).
\]
Il rÈsulte de la fonctorialitÈ de $\Func(\ca,\cb)$ que $\aicat$ est une $2$-catÈgorie.
}\end{definition}

\begin{remarque}{\em
Soit $f_1$ et $f_2 \in \Func(\ca,\cb)$. Soit $H$ un morphisme de
$\Hom_\Func(f_1,f_2)$ qui est un zÈro cycle. Soit $\id_A$ le morphisme identitÈ de $A \in \ca$.
Comme $H$ est un cycle, nous avons la relation
\[
m^\cb_1(h_1(\id_A)) - m^\cb_2(h_{A} \tso f_1\id_A) +
m^\cb_2( f_2\id_A \tso h_{A}) = 
\]
\[
 - m^\cb_2(h_{A} \tso \id_{f_1A}) +
m^\cb_2( \id_{f_2A} \tso h_{A}) = 0.
\]
Nous avons donc un diagramme commutatif dans $H^0\cb$
\[
\xymatrix{
\dot f_1 A \ar[r]^{\id} \ar[d]_{h_{A}}  & \dot f_1 A \ar[d]^{h_{A}}\\
f_2A \ar[r]_{\id} & f_2A.
}
\]
}\end{remarque}

\section{ThÈorie de l'homotopie des $\ai$-foncteurs}
\label{section_theorie_homotopie_ai-foncteurs}

\noindent Cette section est divisÈe en deux sous-sections. 
Soit $\ca$ et $\cb$ deux $\ai$-catÈgories strictement unitaires sur $\mathbb A$ et $\mathbb B$.
Dans la premiËre, nous construisons une gÈnÈralisation de l'$\ai$-foncteur de
Yoneda $y$ (\ref{definition_ai-foncteur_de_Yoneda}) : nous dÈfinissons un $\ai$-foncteur
\[
z : \Func(\ca,\cb) \ra \cc_\infty (\ca,\cb), \quad \ca,\cb \in \aicat,
\]
qui nous redonne l'$\ai$-foncteur de
Yoneda pour $\ca$ Ègal ‡ $e_{\mathbb B}$. Nous montrons ensuite que l'$\ai$-foncteur
de Yoneda gÈnÈralisÈ $z$ induit un quasi-isomorphisme dans les espaces de morphismes.
Dans la seconde partie, nous dÈfinissons les Èquivalences
faibles de l'$\ai$-catÈgorie $\Func(\ca,\cb)$ 
(elles sont l'analogue $\ai$-catÈgorique des {\em homotopies} entre $\ai$-morphismes)
et nous les caractÈrisons
‡ l'aide de leurs images par l'$\ai$-foncteur $z$.

\subsection{L'$\ai$-foncteur de Yoneda gÈnÈralisÈ}
\label{section_ai-foncteur_de_Yoneda_generalise}

L'$\ai$-foncteur de Yoneda gÈnÈralisÈ \indexnotation{z}
\[
z : \Func(\ca,\cb) \ra \cc_\infty (\ca,\cb)
\]
est dÈfini comme la composition 
\index{Yoneda generalise (A-infini foncteur de)@{Yoneda gÈnÈralisÈ ($\ai$-foncteur de)}}
\index{A-infini foncteur@{$\ai$-foncteur}!de Yoneda gÈnÈralisÈ}
\[
\Func(\ca,\cb) \ra \Func(\ca,\cc_\infty \cb) \arr{\theta^{-1}}\cc_\infty (\ca,\cb)
\]
o˘ la premiËre flËche est induite par le foncteur de Yoneda $y : \cb \ra \cc_\infty \cb$
du chapitre \ref{chapitre_objets_tordus} et o˘ $\theta$ est dÈfinie dans la proposition ci-dessous.

\begin{proposition} \label{proposition_Func_aiBimodu}
Soit $\ca$ et $\cb$ deux $\ai$-catÈgories sur $\mathbb A$ et $\mathbb B$.
Il existe un isomorphisme fonctoriel de catÈgories diffÈrentielles graduÈes \indexnotation{theta}
\[
\theta : \cn_\infty(\ca,\cb) \arr{\sim} \Nunc(\ca,\cn_\infty \cb).
\]
Il se restreint en un isomorphisme 
\[
\theta : \cc_\infty(\ca,\cb) \arr{\sim} \Func(\ca,\cc_\infty \cb)
\]
si $\ca$ et $\cb$ sont strictement unitaires.
\end{proposition}

{\noindent \em DÈmonstration de la proposition (\ref{proposition_Func_aiBimodu}) :}\\

{\noindent \bf Le foncteur $\theta$}\\

Nous rappelons (\ref{corollaire2_lemme_clef}) que l'application
\[
\begin{array}{rcl}
\Obj \cn_\infty(\ca,\cb) & \ra & \Obj \Nunc(\ca,\cn_\infty \cb),\\
M & \mapsto & \theta_M,
\end{array}
\]
est une bijection.
Nous allons Ètendre cette application
en un isomorphisme de catÈgories diffÈrentielles graduÈes
\[
\theta : \cn_\infty (\ca,\cb) \ra \Nunc(\ca,\cn_\infty \cb).
\]
Soit $X$ et $X'$ deux $\ca$-$\cb$-bipolydules et
\[
f : X \ra X'
\]
un morphisme de $\cn(\ca,\cb)$. Il est donnÈ par des
morphismes
\[
f_{i,j} : \ca \tpo i \tso X \tso \cb \tpo j \ra X', \quad i,j \geq 0.
\]
Le morphisme
\[
\theta(f) \in \Hom_{\Nunc}(\theta_{X},\theta_{X'})
\]
est donnÈ par un morphisme
\[
\Ba \ca \ra \lrus{\theta_{X'}}{\big(\cn_\infty \cb \big)}{\theta_X} = \Hom_{\ct S\cb}(SX \tso \ct S\cb,SX' \tso \ct S\cb)
\]
qui envoie un ÈlÈment $\phi$ de $(S\ca)\tpo i$ de degrÈ $|\phi|$ sur l'unique
morphisme (voir \ref{lemme_comodules_colibres}) $\Upsilon$ tel que la composition $p_1 \circ \Upsilon$
a pour composantes les morphismes
\[
\xymatrix{
SX \tso (S\cb)\tpo j \ar[rr]^(.4){(-1)^{|\phi|} \phi \tso \Id} & &
(S\ca)\tpo i \tso SX \tso (S\cb)\tpo j \ar[r]^(.75){F_{i,j}} &
SX',} \quad j \geq 0.
\]
Remarquons que si $i = 0$, le morphisme
\[
\Upsilon : SX \tso \ct S\cb \ra SX' \tso \ct S\cb
\]
est le morphisme donnÈ par les morphismes $F_{0,j}$, $j \geq 0.$
Nous avons ainsi dÈfini un isomorphisme d'objets graduÈs
\[
\Hom_{\cn_\infty (\ca,\cb)}\big(X,X'\big) \ra \Hom_{\Nunc(\ca,\cn_\infty \cb)}\big(\theta_X,\theta_{X'}\big).
\]
Montrons que cet isomorphisme dÈfinit un isomorphisme de catÈgories diffÈrentielles graduÈes.
Soit $f$ de degrÈ $p.$
La compatibilitÈ ‡ la composition $m_2$ est immÈdiate.
Montrons la compatibilitÈ ‡ $m_1$. 
Soit $\phi \in (S\ca)\tpo n$ de degrÈ $|\phi|$ et soit $\kappa \tso \psi \in SX \tso (S\cb)\tpo{n'}$.
Nous avons les ÈgalitÈs
(le calcul est le mÍme que pour la dÈmonstration du lemme clef \ref{lemme_clef})
\[
\begin{array}{rl}
 & m^\cb_1(\theta(f))(\phi)(\kappa \tso \psi) \\
 = &  (-1)^{|\phi|+1}\big[ \sum b^{X'}_{0,\beta'}(f_{n,\beta} \tso \Id \tpo {\beta'})
   (\phi \tso \kappa \tso \psi)\\
 & -(-1)^{p}\sum f_{i,j}(\Id \tpo{n}\tso b_{0,\alpha}^X \tso \Id \tpo \beta ) (\phi \tso \kappa \tso \psi)\\
 & -(-1)^{p}\sum f_{i,j}(\Id \tpo{n}\tso \Id_X \ts \Id \tpo \alpha \tso b^\cb \tso \Id \tpo \beta )
     (\phi \tso \kappa \tso \psi)\\
 & -(-1)^{p} \sum f_{i,j}(\Id \tpo \alpha \tso b^\ca \tso \Id \tpo \beta \tso \Id_X \ts \Id \tpo{n'})
     (\phi \tso \kappa \tso \psi)\big],\\
\end{array}
\]
\[
\begin{array}{rl}
 & -m^\cb_2(\theta_{X'},\theta(f))(\phi)(\kappa \tso \psi) \\
 = &  (-1)^{|\phi|+1}\sum_{\alpha'>0} b^{X'}_{\alpha',\beta'}(\Id\tpo{\alpha'} \tso f_{\alpha,\beta} \tso \Id \tpo {\beta'})
   (\phi \tso \kappa \tso \psi),
\end{array}
\]
\[
\begin{array}{rl}
 & m^\cb_2(\theta(f),\theta_X)(\phi)(\kappa \tso \psi) \\
 = &  -(-1)^{p +|\phi|+1}
\sum_{\alpha >0} f_{\alpha',\beta'}(\Id\tpo{\alpha'} \tso b^{X}_{\alpha,\beta} \tso \Id \tpo {\beta'})
   (\phi \tso \kappa \tso \psi).
\end{array}
\]
Nous en dÈduisons l'ÈgalitÈ
\[
d(\theta(f)) = m_1^{\cf}(\theta(f)) - m^\cf_2(\theta_X' \tso \theta(f)) + m_2^\cf(\theta(f) \tso \theta_X) = \theta(d(f))
\]
et nous avons le rÈsultat.\\

{\noindent \bf CompatibilitÈ de $\theta$ ‡ la fonctorialitÈ}\\

Si $f : \ca' \ra \ca$ et $g : \cb \ra \cb'$ sont des $\ai$-foncteurs,
ils induisent clairement des morphismes qui rendent commutatifs les carrÈs
\[
\begin{array}{rl}
\hspace*{-.3cm}\xymatrix{
\cn_\infty(\ca,\cb) \ar[r]^{f^*} \ar[d]_\theta & \cn_\infty (\ca',\cb) \ar[d]^\theta \\
\Nunc (\ca,\cn_\infty \cb) \ar[r]_{f^*} & \Nunc (\ca',\cn_\infty \cb)
} &
\xymatrix{
\cn_\infty(\ca,\cb) \ar[r]^{g_*} \ar[d]_\theta & \cn_\infty (\ca,\cb') \ar[d]^\theta \\
\Nunc (\ca,\cn_\infty \cb) \ar[r]_{g_*} & \Nunc (\ca,\cn_\infty \cb').
}
\end{array}
\]

{\noindent \bf Le cas strictement unitaire}\\

Supposons maintenant que $\ca$ et $\cb$ sont des $\ai$-catÈgories strictement unitaires.
Nous avons les sous-catÈgories (\ref{section_categorie_dg_polydules})
\[
\cc_\infty (\ca,\cb) \subset \cn_\infty (\ca,\cb) \quad \mbox{et} \quad \Nunc(\ca,\cc_\infty \cb)
\subset \Func(\ca,\cn_\infty \cb).
\]
Par la remarque (\ref{remarque3_lemme_clef}), la bijection
\[
\begin{array}{rcl}
\Obj \cn_\infty(\ca,\cb) & \ra & \Obj \Nunc(\ca,\cn_\infty \cb),\\
M & \mapsto & \theta_M,
\end{array}
\]
se restreint en une bijection
\[
\Obj \cc_\infty(\ca,\cb)  \ra  \Obj \Func(\ca,\cc_\infty \cb)
\]
et il est clair que, pour $X$ et $X'$ dans $\cc_\infty (\ca,\cb)$, l'application
$f \mapsto \theta(f)$ induit un isomorphisme
\[
\Hom_{\cc_\infty (\ca,\cb)}\big(X,X'\big) \arr{\sim} \Hom_{\Func(\ca,\cc_\infty \cb)}\big(\theta_X,\theta_{X'}\big).
\]
Nous avons donc un isomorphisme de catÈgories diffÈrentielles graduÈes
\[
\theta : \cc_\infty (\ca,\cb) \arr{\sim} \Func(\ca,\cc_\infty \cb).
\]
\findem

\begin{theoreme} \label{theoreme_foncteur_de_Yoneda_generalise}
L'$\ai$-foncteur de Yoneda gÈnÈralisÈ
\[
z : \Func(\ca,\cb) \ra \cc_\infty(\ca,\cb) 
\]
induit un quasi-isomorphisme dans les espaces de morphismes.
\end{theoreme}

 CommenÁons par quelques lemmes.\\

Soit $\big(\Nunc(\ca,\cb)\big)_u$ la {\em sous-catÈgorie} pleine de $\Nunc(\ca,\cb)$ formÈe
des $\ai$-foncteurs strictement unitaires. \indexnotation{Nuncu}

\begin{lemme}\label{lemme1_theoreme_foncteur_de_Yoneda_generalise}
Le foncteur fidËle
\[
\Func(\ca,\cb) \ra \Nunc(\ca,\cb)
\]
induit un isomorphisme
\[
H^*\Func(\ca,\cb) \ra H^*\big(\Nunc(\ca,\cb)\big)_u.
\]
\end{lemme}

\dem Dans cette dÈmonstration, nous utilisons une filtration qui
est adaptÈe de celle de J.~A.~Guccione et J.~J.~Guccione \cite{Guccione96}.

Soit $f_1$ et $f_2$ deux $\ai$-foncteurs strictement unitaires $\ca \ra \cb.$
Nous rappelons que l'espace des ÈlÈments strictement unitaires de
\[
\Hom_{\Nunc}(f_1,f_2) = \Hom_{\sf C(\mathbb A,\mathbb A)}(\ct S\ca ,\lrus{\dot f_2}{\cb}{\dot f_1})
\]
est formÈ des $H$ qui se factorisent par $\ct S\b \ca$, o˘ $\b \ca$ est le conoyau de l'unitÈ de $\ca$.
Nous avons donc l'ÈgalitÈ
\[
\Hom_{\Func}(f_1,f_2) =  \Hom_{\sf C(\mathbb A,\mathbb A)}
\Big(\bigoplus_{0 \leq p} (S\b \ca)\tpo p, \lrus{\dot f_2}{\cb}{\dot f_1}\Big).
\]
Pour tout $p \geq 0$, nous posons
\[
F_p = \Hom_{\sf C(\mathbb A,\mathbb A)}
\Big(\bigoplus_{0 \leq i < p} (S\b \ca)\tpo i, \lrus{\dot f_2}{\cb}{\dot f_1}\Big)
\oplus \Hom_{\sf C(\mathbb A,\mathbb A)}\Big(\bigoplus_{0 \leq j}
(S\b \ca)\tpo p \tso (S\ca)\tpo j, \lrus{\dot f_2}{\cb}{\dot f_1}\Big).
\]
Nous avons clairement l'inclusion $F_{i+1} \subset F_i$, $i \geq 0$. 
La limite inverse des $F_p$, $p \geq 0,$ est l'espace  $\Hom_{\Func}(f_1,f_2)$
et $F_0$ est l'espace  $\Hom_{\Nunc}(f_1,f_2).$
Nous avons une injection d'espaces graduÈs
\[
J_p : F_p \hookrightarrow \Hom_{\sf C(\mathbb A,\mathbb A)}(\ct S\ca, \lrus{\dot f_2}{\cb}{\dot f_1}),
\quad p \geq 0.
\]
Munissons $\Hom_{\sf C(\mathbb A,\mathbb A)}(\ct S\ca, \lrus{\dot f_2}{\cb}{\dot f_1})$
de la diffÈrentielle $m_1^\Nunc$ et montrons qu'elle induit une
diffÈrentielle sur $F_p$, $ p\geq 1$.

Soit $p \geq 1$. Notons $Q_p$ la projection sur le conoyau de $J_p$. 
Soit $H \in \Hom_{\Nunc}(f_1,f_2)$ tel que $Q_p (H) = 0$. Cette condition est Èquivalente
au fait que
les morphismes $h_i$, $i \geq 0$, (dÈfinis en \ref{remarque_ai-categorie_Func}) vÈrifient les
Èquations
\[
h_i ((\Id \tpo \alpha \tso \unite \tso \Id \tpo \beta) \tso \Id \tpo \gamma) = 0,
\quad \alpha + 1 + \beta + \gamma = i, \quad \alpha + 1 + \beta \leq p.
\]
Nous dÈduisons de la description concrËte (\ref{definition_ai-cat_Nunc}) de l'ÈlÈment
$m_1^\Nunc (H)$ que la composition de $m_1^\Nunc (H)$ avec
\[
((\Id \tpo \alpha \tso \unite \tso \Id \tpo \beta) \tso \Id \tpo \gamma),
\quad \alpha + 1 + \beta + \gamma = i, \quad \alpha + 1 + \beta \leq p,
\]
s'annule.
Ceci montre que $Q_p(m_1^\Nunc (H)) = 0.$ 
Nous en dÈduisons que la diffÈrentielle $m_1^\Nunc$ induit une diffÈrentielle sur
l'objet graduÈ $F_p$, $p\geq 1$.

Montrons que le quotient de l'inclusion $F_{p+1} \subset F_p$, $p\geq 0$, est contractile.
Soit $G_p$ le conoyau de cette inclusion. Il est isomorphe ‡
\[
\Hom_{\sf C(\mathbb A,\mathbb A)}\Big(\bigoplus_{0 \leq j}
(Se)\tp p \tso (S\ca)\tpo j, \lrus{\dot f_2}{\cb}{\dot f_1}\Big) =
\Hom_{\sf C(\mathbb A,\mathbb A)}\Big((Se)\tpo p \tso \ct S\ca, \lrus{\dot f_2}{\cb}{\dot f_1}\Big)
\]
Soit $H$ un ÈlÈment de $F_i$ de degrÈ $|H|$. Nous dÈduisons de la description concrËte (\ref{definition_ai-cat_Nunc})
de $m_1^\Nunc (H)(\phi)$, o˘ $\phi$ est un ÈlÈment de $(Se)\tpo p \tso \ct S\ca$,
l'ÈgalitÈ
\[
m_1^{G_p}(H) = m_1^{\cf(\ca,\cb)}(H),
\]
o˘ $\cf (\ca,\cb)$ est la catÈgorie munie des compositions naÔves (\ref{section_Nunc}).
Par dÈfinition, l'ÈlÈment $m_1^{\cf(\ca,\cb)}(H)$ est Ègal ‡
\[
b^{\Ba \ca} \circ H -(-1)^{|H|} H \circ m^\cb_1.
\]
Comme l'$\ai$-catÈgorie
$\ca$ est strictement unitaire, elle est $H$-unitaire (\ref{lemme_homologiqement_unitaire_H-unitaire}).
Sa construction bar est donc quasi-isomorphe ‡ $0$. Nous en dÈduisons que $G_p$ est contractile.

Montrons que l'inclusion
\[
J : \Hom_{\Func}(f_1,f_2) \hookrightarrow \Hom_{\Nunc}(f_1,f_2)
\]
est un quasi-isomorphisme.
Les complexes $G_p$, $p \geq 0$, sont tous contractiles.
Nous en dÈduisons que le conoyau de l'injection $J_p$, $p\geq 0$, est isomorphe ‡
\[
\bigoplus_{0 \leq i \leq p} G_i.
\]
C'est un espace contractile.
L'espace $\Hom_{\Nunc}(f_1,f_2)$ est donc isomorphe ‡
\[
F_p \oplus \bigoplus_{0 \leq i \leq p} G_i, \quad p \geq 0.
\]
Le conoyau de $J$ est donc
\[
\prod_{0 \leq i} G_i.
\]
Il est clairement contractile, d'o˘ le rÈsultat. \findem

\begin{lemme} \label{lemme2_theoreme_foncteur_de_Yoneda_generalise}
Soit $\ca'$ et $\cb'$ des $\ai$-catÈgories sur $\mathbb A$ et $\mathbb B$ et
\[
g : \ca \ra \ca' \quad \mbox{et} \quad g' : \cb \ra \cb'
\]
des $\ai$-quasi-isomorphismes dans $\sf C(\mathbb A,\mathbb A)$ et $\sf C(\mathbb B,\mathbb B)$.
ConsidËrons les comme des $\ai$-foncteurs (\ref{remarque_ai-morphismes_ai-foncteurs}).
Les $\ai$-foncteurs
\[
g^* : \Nunc(\ca',\cb) \ra \Nunc(\ca,\cb) \quad \mbox{et} \quad g'_* : \Nunc(\ca,\cb) \ra \Nunc(\ca,\cb')
\]
induisent des quasi-isomorphismes dans les espaces de morphismes.
\end{lemme}

Nous dÈduisons de ce lemme et du lemme (\ref{lemme1_theoreme_foncteur_de_Yoneda_generalise})
le corollaire suivant :

\begin{corollaire} \label{corollaire_lemme2_theoreme_foncteur_de_Yoneda_generalise}
Reprenons les hypothËses du lemme (\ref{lemme2_theoreme_foncteur_de_Yoneda_generalise}).
Si les $\ai$-catÈgories $\ca$, $\ca'$, $\cb$ et $\cb'$ sont strictement unitaires et que
les $\ai$-morphismes $g$ et $g'$ sont strictement unitaires,
les $\ai$-foncteurs restreints
\[
\Func(\ca',\cb) \ra \Func(\ca,\cb) \quad \mbox{et} \quad \Func(\ca,\cb) \ra \Func(\ca,\cb')
\]
induisent des quasi-isomorphismes dans les espaces de morphismes.
\findem
\end{corollaire}

{\em DÈmonstration du lemme \ref{lemme2_theoreme_foncteur_de_Yoneda_generalise} :}\\
Par la proposition (\ref{proposition_torsion_ai-foncteurs_equiv_faible_I}), il suffit de montrer
que les $\ai$-foncteurs induits par $g$ et $g'$
\[
\cf(\ca',\cb) \ra \cf(\ca,\cb) \quad \mbox{et} \quad \cf(\ca,\cb) \ra \cf(\ca,\cb'),
\]
o˘ $\cf (\ca',\cb),$ $\cf(\ca,\cb),$ $\cf(\ca,\cb)$ et $\cf(\ca,\cb')$
sont les catÈgories munies des compositions naÔves (voir \ref{section_Nunc}),
donnent des quasi-isomorphismes dans les espaces de morphismes.
Les espaces de morphismes
\[
\Hom_{\cf(\ca,\cb)}(f_1,f_2) = \Hom_{\sf C(\mathbb A,\mathbb A)}\big(\ct S \ca,\lrus{\dot f_2}{\cb}{\dot f_1}\big)
\]
sont munies de la diffÈrentielle
\[
\delta  : H \mapsto m^\cb_1 \circ H - (-1)^{|H|} H \circ b^{\Ba \ca}.
\]
Comme les morphismes $g'_1 : \cb \ra \cb'$ et $\Ba g : \Ba \ca \ra \Ba \ca'$
sont des quasi-isomorphismes, nous avons le rÈsultat.
\findem\\

{\em DÈmonstration du thÈorËme (\ref{theoreme_foncteur_de_Yoneda_generalise}) : }\\

Nous allons d'abord montrer que nous pouvons se ramener au cas o˘ les $\ai$-catÈgories
strictement unitaires sont diffÈrentielles graduÈes unitaires, puis nous prouverons le rÈsultat en utilisant
des arguments d'algËbre homologique classique.\\

La proposition (\ref{proposition_modele_differentiel_gradue})
nous donne des modËles diffÈrentiels graduÈs  unitaires $\ca'$ et $\cb'$ munis d'$\ai$-quasi-isomorphismes
strictement unitaires
\[
\ca \ra \ca' \quad \mbox{et} \quad \cb \ra \cb'.
\]
Le lemme \ref{lemme2_theoreme_foncteur_de_Yoneda_generalise} et son corollaire
\ref{corollaire_lemme2_theoreme_foncteur_de_Yoneda_generalise} nous donne un diagramme
\[
\xymatrix{
\Func(\ca,\cb) \ar[r]^z  \ar[d]_{}&  
\cc_\infty(\ca,\cb) \ar[d]^{} \\
\Func(\ca,\cb') \ar[r]^z &  \cc_\infty(\ca,\cb')\\
\Func(\ca',\cb') \ar[r]^{z} \ar[u]_{} & \cc_\infty(\ca',\cb')\ar[u]_{}\\
}
\]
dont toutes les flËches verticales induisent des quasi-isomorphismes dans les espaces
de morphismes. Ils nous suffit donc de montrer que
\[
z : \Func(\ca,\cb)\ra  \cc_\infty(\ca,\cb)
\]
est un quasi-isomorphismes dans le cas o˘ $\ca$ et $\cb$ sont diffÈrentielles
graduÈes unitaires. Le lemme (\ref{lemme2_theoreme_foncteur_de_Yoneda_generalise})
et la proposition (\ref{proposition_strictication_unitaire_bipolydules_graduee})
montrent qu'il est Èquivalent de montrer que
\[
z : \big(\Nunc(\ca,\cb)\big)_u \ra  \big(\cn_\infty(\ca,\cb)\big)_u
\]
est un quasi-isomorphisme.
Soit $f_1$ et $f_2$ des $\ai$-foncteurs strictement unitaires $\ca \ra \cb$.
Nous avons un isomorphisme
\[
\xymatrix{
\Hom_{\sf C(\mathbb A,\mathbb A)}(\Ba \ca, \lrus{\dot f_2}{\cb}{\dot f_1})
\ar[r] &
\Hom_{\ca\op \tso \ca}(\ca \tso \Ba \ca \tso \ca, \lrus{\dot f_2}{\cb}{\dot f_1}).
}
\]
Rappelons (\ref{lemme_Yoneda_quasi-isomorphisme}) que l'$\ai$-foncteur de Yoneda
\[
y : \cb \ra \cc_\infty \cb
\]
induit un quasi-isomorphisme dans les espaces de morphismes.
Nous avons donc un quasi-isomorphisme
\[
\xymatrix{
\Hom_{\ca\op \tso \ca}(\ca \tso \Ba \ca \tso \ca, \lrus{\dot f_2}{\cb}{\dot f_1})
\ar[d] \\
\Hom_{\ca\op \tso \ca}\big(\ca \tso \Ba \ca \tso \ca, \Hom_{\cc_\infty \cb}(y \circ f_1,y \circ 
f_2)\big).
}
\]
Notons $\mathbf R \Hom$ le foncteur dÈrivÈ ‡ droite qui calcule les groupes $\Ext^*.$
Le dernier terme ci-dessus se rÈcrit
\[
\mathbf R\Hom_{\ca\op \tso \ca}\big(\ca, \mathbf R\Hom_{\cb}(y\circ  f_1,y\circ  f_2)\big).
\]
Il est isomorphe ‡
\[
\mathbf R\Hom_{\ca\op \tso \cb}\big(y\circ  f_1,y\circ  f_2\big)
\]
qui est isomorphe ‡ 
\[
\Hom_{\ca\op \tso \cb}
\big(\ca \tso \ct S\ca \tso  S(y \circ  f_1) \tso \ct S \cb \tso \cb, S(y \circ  f_2)\big) \iso
\]
\[
\Hom_{\sf C(\mathbb A,\mathbb B)}
\big(\ct S\ca \tso  S(y \circ  f_1) \tso \ct S \cb, S(y \circ  f_2)\big) \iso
\]
\[
\Hom_{(\ct S\ca)\op \tso (\ct S\cb) }
\big(B(y \circ  f_1), B(y \circ  f_2)\big).
\]
Comme nous avons des ÈgalitÈs de $\ca$-$\cb$-bipolydules
\[
y \circ f = z (f), \quad f \in \Nunc(\ca,\cb),
\]
le lemme (\ref{proposition_Func_aiBimodu}) montre que le dernier espace de morphismes ci-dessus est 
\[
\Hom_{\cn_{\infty}(\ca,\cb)}\big(z(f_1),z (f_2)\big).
\]
La composition de tous les (quasi-)isomorphismes ci-dessus Ètant le morphisme
\[
z(f_1,f_2) : \Hom_{\Nunc(\ca,\cb)}(f_1,f_2) \ra \Hom_{\cn_\infty(\ca,\cb)}(f_1,f_2),
\]
nous avons le rÈsultat.
\findem

\begin{remarque}{\em
Par construction, l'image de l'$\ai$-foncteur $z$ est formÈe
des $\ca$-$\cb$-bipolydules qui sont de la
forme
\[
\cb (?,\dot f\?), \quad f \in \Func (\ca,\cb).
\]
Ils sont libres en tant que $\cb$-polydules.
}\end{remarque}

\subsection{Equivalences faibles d'$\ai$-foncteurs}
\label{section_equiv_faible_ai-foncteurs}

\noindent Les Èquivalences faibles entre $\ai$-foncteurs sont
l'analogue $\ai$-catÈgorique des {\em homotopies} entre $\ai$-morphismes.

\begin{definition}{\em
Soit $\ca$ et $\cb$ deux $\ai$-catÈgories sur $\mathbb A$ et $\mathbb B$. Soit
$f$ et $g$ deux $\ai$-foncteurs $\ca \ra \cb$. 
Un ÈlÈment $H \in Z^0\Hom_{\Nunc}(f,g)$ est une {\em Èquivalence faible}
\index{equivalence faible@{Èquivalence faible}!d'A-infini foncteurs@{d'$\ai$-foncteurs}}
s'il devient un isomorphisme dans $H^0 \Nunc(\ca,\cb).$ Nous dirons alors que $f$ et $g$ sont
{\em faiblement Èquivalents} et Ècrirons $f \sim g.$
}\end{definition}

\begin{remarque}{\em
Supposons que $\ca$ et $\cb$ sont strictement unitaires et $f$ et $g$ sont des $\ai$-foncteurs
strictement unitaires. D'aprËs le lemme (\ref{lemme1_theoreme_foncteur_de_Yoneda_generalise})
$f$ et $g$ sont faiblement Èquivalents si et seulement si il existe un morphisme strictement
unitaire $H \in Z^0\Hom_{\Func}(f_1,f_2)$ qui devient
un isomorphisme dans $H^0 \Func(\ca,\cb).$
}\end{remarque}

\begin{proposition} \label{proposition_equivalence_faible_d'ai-foncteurs}
Soit $\ca$ et $\cb$ deux $\ai$-catÈgories strictement unitaires sur $\mathbb A$ et $\mathbb B$.
Soit $f$ et $g$ deux $\ai$-foncteurs strictement unitaires $\ca \ra \cb$. 
Un ÈlÈment $H \in Z^0\Hom_{\Nunc}(f,g)$ est une Èquivalence faible
si et seulement si $h_0 : e_{\mathbb A} \ra \Hom_{\cb}(\dot f?,\dot g\?)$ induit un
isomorphisme de foncteurs $H^0 f \ra H^0 g$ de $H^0\ca$ dans $H^0\cb.$
\end{proposition}

\dem
D'aprËs le thÈorËme (\ref{theoreme_foncteur_de_Yoneda_generalise}), nous avons un isomorphisme
\[
H^0\Func(\ca,\cb) \arr{\sim} H^0\cc_\infty (\ca,\cb).
\]
L'ÈlÈment $H$ est donc une Èquivalence
faible si et seulement si le morphisme de $\ca$-$\cb$-bipolydules
\[
z(H) : z(f) \ra z(g)
\]
est une Èquivalence d'homotopie dans $\cc_\infty (\ca,\cb),$ c'est-‡-dire
(voir l'Èquivalence entre D2 et D3 dans \ref{theoreme_categorie_derivee_strictement_unitaire})
si et
seulement si $z(H)$ est un $\ai$-quasi-isomorphisme de $\ca$-$\cb$-bipolydules.
Par dÈfinition des $\ai$-quasi-isomorphismes, ceci est Èquivalent au fait que le
morphisme de $\sf C(\mathbb A,\mathbb B)$
\[
S^{-1}(z(H))_{0,0} : \cb(?,f \?) \ra \cb(? ,g\?)
\]
est un quasi-isomorphisme, c'est-‡-dire qu'il devient un isomorphisme en cohomologie.
Comme le foncteur de Yoneda au sens classique
envoie la classe dans $H^*\cb$ de
\[
h_{A} = h_0(\id_A) : \dot fA \ra \dot gA, \quad A \in \mathbb A,
\]
sur
\[
S^{-1}(z(H))_{0,0} : H^*\cb(?,f A) \ra H^*\cb(? ,g A),
\]
$S^{-1}(z(H))_{0,0}$ est un quasi-isomorphisme si et seulement si 
$h_A$ induit un isomorphisme dans $H^*\cb$, ou de maniËre Èquivalente dans $H^0\cb$.
\findem

%% file: figure1.pstex_t
\begin{picture}(0,0)%
\includegraphics{figure1.pstex}%
\end{picture}%
\setlength{\unitlength}{3947sp}%
\begingroup\makeatletter\ifx\SetFigFont\undefined%
\gdef\SetFigFont#1#2#3#4#5{%
  \reset@font\fontsize{#1}{#2pt}%
  \fontfamily{#3}\fontseries{#4}\fontshape{#5}%
  \selectfont}%
\fi\endgroup%
\begin{picture}(3525,2295)(451,-1711)
\put(1126,239){\makebox(0,0)[lb]{\smash{\SetFigFont{12}{14.4}{\rmdefault}{\mddefault}{\updefault}{\color[rgb]{0,0,0}$f_0$}%
}}}
\put(601,-1486){\makebox(0,0)[lb]{\smash{\SetFigFont{12}{14.4}{\rmdefault}{\mddefault}{\updefault}{\color[rgb]{0,0,0}$A_2$}%
}}}
\put(451,-511){\makebox(0,0)[lb]{\smash{\SetFigFont{12}{14.4}{\rmdefault}{\mddefault}{\updefault}{\color[rgb]{0,0,0}$x_1$}%
}}}
\put(451,-1111){\makebox(0,0)[lb]{\smash{\SetFigFont{12}{14.4}{\rmdefault}{\mddefault}{\updefault}{\color[rgb]{0,0,0}$x_2$}%
}}}
\put(601,-286){\makebox(0,0)[lb]{\smash{\SetFigFont{12}{14.4}{\rmdefault}{\mddefault}{\updefault}{\color[rgb]{0,0,0}$A_0$}%
}}}
\put(601,-886){\makebox(0,0)[lb]{\smash{\SetFigFont{12}{14.4}{\rmdefault}{\mddefault}{\updefault}{\color[rgb]{0,0,0}$A_1$}%
}}}
\put(3976,-1486){\makebox(0,0)[lb]{\smash{\SetFigFont{12}{14.4}{\rmdefault}{\mddefault}{\updefault}{\color[rgb]{0,0,0}$f_3A_2$}%
}}}
\put(3301,389){\makebox(0,0)[lb]{\smash{\SetFigFont{12}{14.4}{\rmdefault}{\mddefault}{\updefault}{\color[rgb]{0,0,0}$h_3$}%
}}}
\put(2401,389){\makebox(0,0)[lb]{\smash{\SetFigFont{12}{14.4}{\rmdefault}{\mddefault}{\updefault}{\color[rgb]{0,0,0}$h_2$}%
}}}
\put(1576,389){\makebox(0,0)[lb]{\smash{\SetFigFont{12}{14.4}{\rmdefault}{\mddefault}{\updefault}{\color[rgb]{0,0,0}$h_1$}%
}}}
\put(2776,-61){\makebox(0,0)[lb]{\smash{\SetFigFont{12}{14.4}{\rmdefault}{\mddefault}{\updefault}{\color[rgb]{0,0,0}$f_2A_0$}%
}}}
\put(2026,239){\makebox(0,0)[lb]{\smash{\SetFigFont{12}{14.4}{\rmdefault}{\mddefault}{\updefault}{\color[rgb]{0,0,0}$f_1$}%
}}}
\put(2926,239){\makebox(0,0)[lb]{\smash{\SetFigFont{12}{14.4}{\rmdefault}{\mddefault}{\updefault}{\color[rgb]{0,0,0}$f_2$}%
}}}
\put(3826,239){\makebox(0,0)[lb]{\smash{\SetFigFont{12}{14.4}{\rmdefault}{\mddefault}{\updefault}{\color[rgb]{0,0,0}$f_3$}%
}}}
\put(1126,-61){\makebox(0,0)[lb]{\smash{\SetFigFont{12}{14.4}{\rmdefault}{\mddefault}{\updefault}{\color[rgb]{0,0,0}$f_0A_0$}%
}}}
\put(3451,-586){\makebox(0,0)[lb]{\smash{\SetFigFont{12}{14.4}{\rmdefault}{\mddefault}{\updefault}{\color[rgb]{0,0,0}$(h_3)_2(x_2 \tso x_1)$}%
}}}
\put(1801,-1711){\makebox(0,0)[lb]{\smash{\SetFigFont{12}{14.4}{\rmdefault}{\mddefault}{\updefault}{\color[rgb]{0,0,0}$f_1A_2$}%
}}}
\put(2326,-1636){\makebox(0,0)[lb]{\smash{\SetFigFont{12}{14.4}{\rmdefault}{\mddefault}{\updefault}{\color[rgb]{0,0,0}$(h_2)_{A_2}$}%
}}}
\put(2926,-1711){\makebox(0,0)[lb]{\smash{\SetFigFont{12}{14.4}{\rmdefault}{\mddefault}{\updefault}{\color[rgb]{0,0,0}$f_2A_2$}%
}}}
\end{picture}

%% file: figure2.pstex_t
\begin{picture}(0,0)%
\includegraphics{figure2.pstex}%
\end{picture}%
\setlength{\unitlength}{3947sp}%
\begingroup\makeatletter\ifx\SetFigFont\undefined%
\gdef\SetFigFont#1#2#3#4#5{%
  \reset@font\fontsize{#1}{#2pt}%
  \fontfamily{#3}\fontseries{#4}\fontshape{#5}%
  \selectfont}%
\fi\endgroup%
\begin{picture}(3675,2445)(451,-1861)
\put(1126,239){\makebox(0,0)[lb]{\smash{\SetFigFont{12}{14.4}{\rmdefault}{\mddefault}{\updefault}{\color[rgb]{0,0,0}$f_0$}%
}}}
\put(601,-1486){\makebox(0,0)[lb]{\smash{\SetFigFont{12}{14.4}{\rmdefault}{\mddefault}{\updefault}{\color[rgb]{0,0,0}$A_2$}%
}}}
\put(451,-511){\makebox(0,0)[lb]{\smash{\SetFigFont{12}{14.4}{\rmdefault}{\mddefault}{\updefault}{\color[rgb]{0,0,0}$x_1$}%
}}}
\put(451,-1111){\makebox(0,0)[lb]{\smash{\SetFigFont{12}{14.4}{\rmdefault}{\mddefault}{\updefault}{\color[rgb]{0,0,0}$x_2$}%
}}}
\put(601,-286){\makebox(0,0)[lb]{\smash{\SetFigFont{12}{14.4}{\rmdefault}{\mddefault}{\updefault}{\color[rgb]{0,0,0}$A_0$}%
}}}
\put(601,-886){\makebox(0,0)[lb]{\smash{\SetFigFont{12}{14.4}{\rmdefault}{\mddefault}{\updefault}{\color[rgb]{0,0,0}$A_1$}%
}}}
\put(3301,389){\makebox(0,0)[lb]{\smash{\SetFigFont{12}{14.4}{\rmdefault}{\mddefault}{\updefault}{\color[rgb]{0,0,0}$h_3$}%
}}}
\put(2401,389){\makebox(0,0)[lb]{\smash{\SetFigFont{12}{14.4}{\rmdefault}{\mddefault}{\updefault}{\color[rgb]{0,0,0}$h_2$}%
}}}
\put(1576,389){\makebox(0,0)[lb]{\smash{\SetFigFont{12}{14.4}{\rmdefault}{\mddefault}{\updefault}{\color[rgb]{0,0,0}$h_1$}%
}}}
\put(2026,239){\makebox(0,0)[lb]{\smash{\SetFigFont{12}{14.4}{\rmdefault}{\mddefault}{\updefault}{\color[rgb]{0,0,0}$f_1$}%
}}}
\put(2926,239){\makebox(0,0)[lb]{\smash{\SetFigFont{12}{14.4}{\rmdefault}{\mddefault}{\updefault}{\color[rgb]{0,0,0}$f_2$}%
}}}
\put(3826,239){\makebox(0,0)[lb]{\smash{\SetFigFont{12}{14.4}{\rmdefault}{\mddefault}{\updefault}{\color[rgb]{0,0,0}$f_3$}%
}}}
\put(1126,-61){\makebox(0,0)[lb]{\smash{\SetFigFont{12}{14.4}{\rmdefault}{\mddefault}{\updefault}{\color[rgb]{0,0,0}$f_0A_0$}%
}}}
\put(3751,-1636){\makebox(0,0)[lb]{\smash{\SetFigFont{12}{14.4}{\rmdefault}{\mddefault}{\updefault}{\color[rgb]{0,0,0}$f_3A_2$}%
}}}
\put(4126,-886){\makebox(0,0)[lb]{\smash{\SetFigFont{12}{14.4}{\rmdefault}{\mddefault}{\updefault}{\color[rgb]{0,0,0}$(f_3)_2(x_2 \tso x_1)$}%
}}}
\put(2701,-1861){\makebox(0,0)[lb]{\smash{\SetFigFont{12}{14.4}{\rmdefault}{\mddefault}{\updefault}{\color[rgb]{0,0,0}$f_3(x_2)$}%
}}}
\put(3751,-61){\makebox(0,0)[lb]{\smash{\SetFigFont{12}{14.4}{\rmdefault}{\mddefault}{\updefault}{\color[rgb]{0,0,0}$f_3A_0$}%
}}}
\end{picture}

%% file: Ai_equiv.tex
Ce chapitre est divisÈ en deux parties.
Dans la section \ref{section_ai-isomorphie},
nous dÈfinissons l'$\ai$-isomorphie dans une $\ai$-catÈgorie $\ca$ et nous montrerons
que cette notion est un relËvement $\ai$-catÈgorique de l'isomorphie dans $H^0\ca$ au sens classique.
Dans la section \ref{section_ai-equivalences},
nous dÈfinissons les $\ai$-Èquivalences et nous montrons qu'un $\ai$-foncteur $f$ est
une $\ai$-Èquivalence si et seulement si $f_1$ est un quasi-isomorphisme et
$H^0f_1 : H^0 \ca \ra H^0 \cb$ est une Èquivalence de catÈgories au sens classique.
Cette caractÈrisation des $\ai$-Èquivalences a ÈtÈ ÈnoncÈe par M.~Kontsevich \cite{Kontsevich98}.
K.~Fukaya l'a dÈmontrÈ de maniËre indÈpendante \cite[thm.~8.6]{Fukaya01a}, ainsi que
V.~Lyubashenko \cite{Lyubashenko02}.

\section{L'$\ai$-isomorphie}
\label{section_ai-isomorphie}

Soit $\mathbb O$ un ensemble. ConsidËrons comme une catÈgorie de la maniËre
suivante : les objets sont en bijection avec $\mathbb O$ et, pour $i,j \in \mathbb \mathbb O$,
l'espace de morphismes $\Hom_{\mathbb O}(i,j)$ contient un unique ÈlÈment
notÈ $(i,j)$. La composition est alors nÈcessairement donnÈe par
\[
(j,k) \circ (i,j) = (i,k), \quad i,j,k \in \mathbb O.
\]
En particulier, l'identitÈ de $i \in \mathbb O$ est le morphisme $(i,i)$ et
tous les morphismes $(i,j)$ sont des isomorphismes.

\begin{definition}{\em \indexnotation{In}
Soit $n \geq 1.$  ConsidËrons l'ensemble ‡ $n$ ÈlÈments $\{1,\hdots ,n\}.$
Soit $\I_n$ la $\corps$-catÈgorie engendrÈe par la catÈgorie $\{1,\hdots ,n\}.$
}\end{definition}

\begin{remarque}{\em \index{isomorphie}
Soit $n \geq 2$.
Soit $\ca$ une $\corps$-catÈgorie et des objets $A_i \in \Obj \ca$, $1 \leq i \leq n.$
Ils sont isomorphes si et seulement si il existe un foncteur
\[
f : \I_n \ra \ca
\]
qui envoie $i$ sur $A_i$. Nous disons alors que $f$ est un {\em foncteur d'isomorphie} pour
les objets $A_i \in \Obj \ca$, $1 \leq i \leq n.$
}\end{remarque}

\begin{definition}{\em \label{definition_ai-isomorphie}
Soit $n \geq 2$. \index{A-infini isomorphie@{$\ai$-isomorphie}}
Soit $\ca$ une $\ai$-catÈgorie strictement unitaire sur $\mathbb A$
et des objets $A_i \in \mathbb A$, $1 \leq i \leq n.$
Les objets $A_i \in \mathbb A$, $1 \leq i \leq n,$ sont {\em $\ai$-isomorphes} s'il existe
un $\ai$-foncteur strictement unitaire 
\[
f : \I_n \ra \ca
\]
qui envoie $i$ sur $A_i$. Nous disons alors que $f$ est un {\em $\ai$-foncteur d'$\ai$-isomorphie} pour
les objets $A_i \in \mathbb A$, $1 \leq i \leq n.$
}\end{definition}

Nous prouvons maintenant un lemme ÈnonÁÈ dans \cite{Kontsevich98} :

\begin{lemme} \label{lemme_ai-isomorphie}
Soit $\ca$ une $\ai$-catÈgorie strictement unitaire. Soit $n \geq 1$.
Des objets $A_i \in \mathbb A$, $1 \leq i \leq n,$ sont $\ai$-isomorphes dans $\ca$ si
et seulement ils sont isomorphes dans $H^0\ca$.
\end{lemme}

\dem
Comme $\ca$ est strictement unitaire, il existe (\ref{proposition_modele_minimal_strictement_unitaire})
un modËle minimal strictement unitaire $H^*\ca$
pour $\ca$ et des $\ai$-foncteurs strictement unitaires
(\ref{remarque_proposition_modele_minimal_strictement_unitaire})
\[
i : H^* \ca \ra \ca \quad \mbox{et} \quad q : \ca \ra H^*\ca.
\]
Nous en dÈduisons que des objets $A_i \in \mathbb A$, $1 \leq i \leq n,$ sont $\ai$-isomorphes dans $\ca$ si
et seulement si ils sont $\ai$-isomorphes dans $H^*\ca.$
Nous pouvons donc supposer que l'$\ai$-catÈgorie $\ca$ est minimale.

Soit $\ca$ une $\ai$-catÈgorie minimale. Montrons que l'$\ai$-isomorphie
dans $\ca$ entraine l'isomorphie dans $H^0\ca$.
Soit $f : \I_n \ra \ca$ un $\ai$-foncteur d'$\ai$-isomorphie pour $A_i \in \mathbb A$, $1 \leq i \leq n.$
Comme les $\ai$-catÈgories $\I_n$ et $\ca$ sont minimales,
$f_0 : \I_n \ra \ca^0 = H^0\ca$ dÈfinit un foncteur d'isomorphie pour
les objets $A_i \in \mathbb A$, $1 \leq i \leq n.$

Montrons que l'isomorphie dans $H^0\ca$ implique l'$\ai$-isomorphie dans $\ca$.
Soit $g : \I_n \ra H^0 \ca$ un foncteur d'isomorphie pour
les objets $A_i \in \Obj H^0\ca$, $1 \leq i \leq n.$
Nous cherchons un $\ai$-foncteur strictement unitaire
\[
f : \I_n \ra \ca
\]
tel que $f_1 = i \circ g$, o˘ $i$ est l'inclusion $\ca^0 \hookrightarrow \ca$.
D'aprËs le thÈorËme (\ref{theoreme_unite_ai-morphisme}), il suffit de construire
un $\ai$-foncteur $f'$ (non nÈcessairement strictement unitaire) tel que $f'_1 = f_1$.
Nous allons construire les $f'_r$, $r\geq 2$, par rÈcurrence sur $r$.
Supposons donnÈ des morphismes graduÈs $f'_i$, $1 \leq i \leq r$, de degrÈ $1-i$,
dÈfinissant un $\rm{A}_r$-foncteur $\I_n \ra \ca$.
Soit $f'_{r+1}$ un morphisme de degrÈ $-r$. 
Le lemme (\ref{lemme_obstruction_structure_minimale_morphisme}) affirme que
la suite des $f'_i$, $1 \leq i \leq r+1,$
dÈfinit un $\rm{A}_{r+1}$-foncteur si nous avons l'ÈgalitÈ
\[
\delta_{Hoch}(f'_{r+1}) = - r(f'_2, \hdots ,f'_r)
\]
o˘ $r(f'_2, \hdots ,f'_r)$ est un certain cycle du complexe de Hochschild
$C^*(\I_n, \lrus{\dot f'}{\ca}{\dot f'})$.
Comme la catÈgorie $\I_n$ est Èquivalente ‡ la catÈgorie triviale $\I_1$,
le complexe de Hochschild $C^*(\I_n, \lrus{\dot f}{\ca}{\dot f})$ est acyclique.
Il existe donc un morphisme $f'_{r+1}$
tel que les morphismes graduÈs $f'_i$, $1 \leq i \leq r+1$,
dÈfinissent un $\rm{A}_{r+1}$-foncteur $\I_n \ra \ca$.\findem

\section{La caractÈrisation des $\ai$-Èquivalences}
\label{section_ai-equivalences}

\begin{definition}{\em \label{definition_ai-equivalence} \index{A-infini equivalence@{$\ai$-Èquivalence}}
Deux $\ai$-catÈgories strictement unitaires $\ca$ et $\cb$ sur $\mathbb A$ et $\mathbb B$
sont {\em $\ai$-Èquivalentes} s'il existe des $\ai$-foncteurs strictement unitaires
\[
f : \ca \ra \cb \quad \mbox{et} \quad g : \cb \ra \ca
\]
tels que $f \circ g$ et $\Id_{\cb}$ sont $\ai$-isomorphes dans $\Func(\cb,\cb)$ et $g \circ f$ et
$\Id_\ca$ sont $\ai$-isomorphes dans $\Func(\ca,\ca).$
Nous dirons alors que $f$ (ou $g$) est une {\em $\ai$-Èquivalence} entre $\ca$ et $\cb$. 
}\end{definition}

\begin{definition}{\em \label{definition_quasi-equivalence}
Soit $\ca$ et $\cb$ deux catÈgories diffÈrentielles graduÈe  unitaires sur $\mathbb A$ et $\mathbb B$.
Elles sont {\em Èquivalentes} (au sens classique) s'il existe des foncteurs
\[
f : \ca \ra \cb \quad \mbox{et} \quad g : \cb \ra \ca
\]
et des isomorphismes des foncteurs
\[
\mu : f \circ g \ra \Id_\cb \quad \mbox{et} \quad  \nu : g \circ f  \ra \Id_\ca.
\]
Nous dirons alors que $f$ (ou $g$) est une {\em Èquivalence} entre $\ca$ et $\cb$. 
}\end{definition}

\begin{remarque}{\em Soit $\ca$ et $\cb$ deux catÈgories diffÈrentielles graduÈe  unitaires sur
$\mathbb A$ et $\mathbb B$. Supposons qu'elles sont Èquivalentes.
Soit $f$ une Èquivalence entre $\ca$ et $\cb$. Soit $g$, $\mu$ et $\nu$ comme dans la dÈfinition
(\ref{definition_quasi-equivalence}).
L'ÈlÈment $H \in \Hom_{\Func}(f \circ g,\Id_\cb)$ (resp.~$H' \in \Hom_{\Func}(g \circ f,\Id_\ca)$) dÈfini par
\[
h_0 = \mu, \quad h_i = 0, \quad i \geq 1, \quad
\Big(\mbox{resp.}\quad h'_0 = \mu, \quad h'_i = 0, \quad  i \geq 1\Big)
\]
est un cycle dans $\Func(\cb,\cb)$ (resp.~dans $\Func(\ca,\ca)$).
Il induit un isomorphisme dans $H^0\Func(\cb,\cb)$ (resp.~$H^0\Func(\ca,\ca)$).
Ceci montre que $\ca$ et $\cb$ sont $\ai$-Èquivalentes en tant que $\ai$-catÈgories.
}\end{remarque}

L'ÈnoncÈ du thÈorËme suivant est du ‡ M.~Kontsevich \cite{Kontsevich98}.

\begin{theoreme}[Voir aussi K.~Fukaya \cite{Fukaya01a} et V.~Lyubashenko \cite{Lyubashenko02}]
Soit $\ca$ et $\cb$ deux $\ai$-catÈgories strictement unitaires sur $\mathbb A$ et $\mathbb B$
et $f : \ca \ra \cb$ un $\ai$-foncteur strictement unitaire.
Les ÈnoncÈs suivants sont Èquivalents :
\english
\begin{itemize}
\item[a.] $f$ est une $\ai$-Èquivalence.
\item[b.] $f_1$ induit une Èquivalence $H^* \ca \ra H^* \cb,$ o˘ $H^*\ca$ et $H^*\cb$ sont la
cohomologie de $\ca$ et $\cb$ considÈrÈes comme $\corps$-catÈgories graduÈes.
\item[c.] $f_1$ est un quasi-isomorphisme et induit une Èquivalence $H^0 \ca \ra H^0 \cb.$
\end{itemize}\francais
\end{theoreme}

\dem

$ a \Ra b$ :  Supposons que $f$ est une $\ai$-Èquivalence. Soit $g : \cb \ra \ca$ vÈrifiant
les conditions de la dÈfinition (\ref{definition_ai-equivalence}). D'aprËs le lemme
(\ref{lemme_ai-isomorphie}), l'$\ai$-isomorphie dans $\Func(\cb,\cb)$ (resp.~dans $\Func(\ca,\ca)$)
est Èquivalente ‡ l'isomorphie dans $H^0 \Func(\cb,\cb)$ (resp.~dans $H^0\Func(\ca,\ca)$).
Comme $f \circ g$ et $\Id_\cb$ sont isormorphes dans $H^0\Func(\ca,\ca),$
il existe donc un ÈlÈment
\[
H  \in Z^0\Hom_{\Func}(g \circ f,\Id_\cb)
\]
induisant un isomorphisme dans $H^0 \Func(\cb,\cb).$
D'aprËs la proposition (\ref{proposition_equivalence_faible_d'ai-foncteurs}),
le morphisme $h_0$ induit un isomorphisme de foncteurs
\[
H^0(h_0) : H^*(g_1 \circ f_1) \ra H^*\Id_\cb.
\]
L'isomorphisme de foncteurs entre $H^*(f_1 \circ g_1)$ et $\Id_{H^*\ca}$ est construit de la mÍme maniËre.

$ b \Ra c$ : C'est clair.

$c \Ra a$ : Nous allons montrer cette implication dans deux cas particuliers
puis nous montrerons que cela implique le cas gÈnÈral.

{\em Premier cas : l'application $\dot f : \mathbb A \ra \mathbb B$ est une bijection.}\\
Nous pouvons considÈrer que $\mathbb A$ est Ègal ‡ $\mathbb B$ et que $\dot f$ est l'identitÈ de $\mathbb A.$
L'$\ai$-foncteur $f$ est ainsi (\ref{remarque_ai-morphismes_ai-foncteurs}) un $\ai$-morphisme
dans la catÈgorie $\sf C(\mathbb A,\mathbb A).$ D'aprËs le point {\it b} du corollaire
(\ref{corollaire_cmf_aia}), il existe un $\ai$-morphisme $g : \cb \ra \ca$ et des
homotopies $h$ et $h'$ de $f \circ g$ vers $\Id_\cb$ et de $g \circ f$ vers $\Id_\ca$.
Gr‚ce ‡ la proposition (\ref{proposition_ai-algebres_strictement_unitaires}), nous pouvons supposer
que l'$\ai$-morphisme $g$ et les homotopies $h$ et $h'$ sont strictement unitaires.
Soit $H$ l'ÈlÈment de $\Hom_{\Func}(f\circ g,\Id_\cb)$ donnÈ (voir \ref{remarque_ai-categorie_Func})
par les morphismes $h_i,$ $i\geq 1,$ et $h_A = \Id_A$, $A \in \mathbb A$.
Posons $Z = m_1^\Func (H)$. Il est donnÈ
par des morphismes $z_i$, $i \geq 0$. Montrons que $H$ est un cycle dans
$\Hom_{\Func}(\ca,\cb)$. Le morphisme $z_0$ est clairement nul.
Pour $n \geq 1$, nous vÈrifions (en utilisant le fait
que $f\circ g$, $\Id_\cb$ et $h$ sont strictement unitaires) que 
\[
\begin{array}{rl}
z_n = & (f \circ g)_n - (\Id_\cb)_n \\ &  - \sum (-1)^s
m_{r+1+t} ((f\circ g)_{i_1}\ts \hdots \ts (f \circ g)_{i_r} \ts h_{k} \ts (\Id_\cb)_{j_1}\ts \hdots \ts
(\Id_\cb)_{i_t})\\
& - \sum (-1)^{jk+l}h_i(\Id^{\ts j}\ts m_k\ts \Id^{\ts l}).
\end{array}
\]
o˘ $s$ est le signe intervenant dans l'Èquation $(***_n)$ de (\ref{definition_homotopie_ai-morphismes}).
Comme $h$ est une homotopie entre $f \circ g$ et $\Id_\cb$, le terme de droite est
nul. Ceci montre que $H$ est un cycle de $\Hom_{\Func}(f\circ g,\Id_\cb).$
Le morphisme $h_A$, $A \in \mathbb A$ valant $\Id_A,$ la proposition
(\ref{proposition_equivalence_faible_d'ai-foncteurs}) implique que $H$ induit un isomorphisme dans
$H^0\Func(\cb,\cb).$ Nous en dÈduisons (\ref{lemme_ai-isomorphie}) que les $\ai$-foncteurs $\Id_\cb$ et $f \circ g$
sont $\ai$-isomorphes dans $\Func(\cb,\cb).$
L'$\ai$-isomorphie entre $g \circ f$ et $\Id_\ca$ se montre de la mÍme maniËre.

\begin{remarque}{\em
En particulier, ceci implique qu'une $\ai$-catÈgorie strictement unitaire $\ca$ est $\ai$-Èquivalente 
‡ son modËle minimal (\ref{proposition_modele_minimal_strictement_unitaire}) et ‡
tout ses modËles diffÈrentiels graduÈs (\ref{proposition_modele_differentiel_gradue}).
}\end{remarque}

{\em DeuxiËme cas : $f$ est une inclusion $\ca \hookrightarrow \cb$ o˘ $\ca$ est une
sous-$\ai$-catÈgorie pleine de $\cb$.} Gr‚ce ‡ la remarque prÈcÈdente,
Nous pouvons supposer que $\cb$ est diffÈrentielle
graduÈe. Comme $H^0 f : H^0 \ca \ra H^0 \cb$ est une Èquivalence,
il suffit de montrer le thÈorËme dans le cas suivant : Choisissons dans chaque
classe d'isomorphie $[B]$ de $\cb$ un reprÈsentant $B_0$. Soit $\mathbb A$ l'ensemble de
ces rÈprÈsentants. Nous posons $\ca$ Ègale ‡ la sous-catÈgorie pleine de $\cb$ formÈe des objets
$A \in \mathbb A.$ Soit
\[
r : \mathbb B \ra \mathbb A, \quad B  \mapsto r(B) = B_0.
\]
Soit $\ca'$ la catÈgorie diffÈrentielle graduÈe
$\lrus{r}{\ca}{r}$ sur $\mathbb B$ (voir \ref{lemme_ai-categorie_fBf}).
Nous avons alors les ÈgalitÈs
\[
\ca'(B,B') = \ca(B_0,B'_0) = \cb(B_0,B'_0)
\]
et les catÈgories diffÈrentielles graduÈes $\ca$ et $\ca'$ sont Èquivalentes au
sens classique. Ils nous suffit donc de montrer que $\ca'$ et $\cb$ sont
$\ai$-Èquivalentes. Soit $i : \ca \ra \ca'$ l'inclusion. Nous allons construire
une $\ai$-Èquivalence
\[
g : \ca' \ra \cb
\]
tel que $f = g \circ i$. \\

{\em Construction de $g$ :}
Posons $\dot g = \Id_{\mathbb B}$.
L'$\ai$-foncteur $g$ est donc donnÈ par un $\ai$-morphisme
$\ca' \ra  \cb$ dans $\sf C(\mathbb B,\mathbb B).$
Par hypothËse, tout ÈlÈment $B \in \mathbb B$ est $\ai$-isomorphe ‡ $r(B).$
Pour chaque $B \in \mathbb B$, choisissons un ÈlÈment $\alpha_B$ de $\cb(r(B),B)$ qui
devient un isomorphisme dans $H^0\cb(r(B),B).$
ConsidÈrons le diagramme de $\mathbb B$-$\mathbb B$-bimodules diffÈrentiels graduÈs
\[
\diagnum{I}{
\xymatrix{ & \cb(\?,?) \ar[d]^{\alpha^*} \\
\ca'(\?,?) = \cb(r(\?),r(?)) \ar[r]_(.65){\alpha_*} & \cb(r(\?),?).}}
\]
L'$\ai$-foncteur de Yoneda $y : \cb \ra \cc_\infty \cb$ (\ref{definition_ai-foncteur_de_Yoneda})
envoie le diagramme $(I)$ sur un diagramme quasi-isomorphe de $\mathbb B$-$\mathbb B$-bimodules
\[
\diagnum{I'}{\xymatrix{ & \cc_\infty \cb(y\?,y?) \ar[d]^{(y\alpha)^*} \\
\cc_\infty \cb(yr(\?),yr(?)) \ar[r]_{(y\alpha)_*} & \cc_\infty \cb(yr(\?),y?).}}
\]
Pour chaque $B \in \mathbb B$, le morphisme $\alpha_B$ devient un isomorphisme
dans $H^0\cb$. Comme le foncteur de Yoneda induit un quasi-isomorphisme 
dans les espaces de morphismes (\ref{lemme_Yoneda_quasi-isomorphisme}),
il induit un isomorphisme $H^0 \cb \ra H^0\cc_\infty \cb.$ Nous en dÈduisons que
le morphisme $y\alpha_B$ est une Èquivalence d'homotopie dans $\cc_\infty \cb.$
D'aprËs l'Èquivalence entre les catÈgories de D3 et D4 (\ref{theoreme_categorie_derivee_strictement_unitaire}),
il est un quasi-isomorphisme.
Ceci implique que les flËches du diagramme $(I')$ sont des quasi-isomorphismes.
La catÈgorie $\cb$ Ètant diffÈrentielle graduÈe, les $\cb$-bipolydules $y(B)$, $B \in \mathbb B$,
sont des $\cb$-modules diffÈrentiels graduÈs et le morphisme $y\alpha_B : y(r(B)) \ra y(B)$ est un morphisme
de $\cb$-modules diffÈrentiels graduÈs.
L'axiome (CM5) de la catÈgorie $\Modu \cb$ nous donne
une factorisation de $y\alpha_B$ en une cofibration triviale et une fibration triviale
\[
\xymatrix{
yr(B) \ar@{ >->}^{i_B}[r] & m(B) \ar@{->>}[r]^{p_B} & yB.
}
\]
Gr‚ce ‡ l'axiome (CM4) de la catÈgorie $\Modu \cb$, il existe un quasi-isomorphisme $\sigma_B$
tel que $p_B \circ \sigma_B = \Id_{yB}.$
Le morphisme
\[
\cc_\infty \cb(y\?,y?) \ra \cc_\infty \cb(m\?,m?), \quad x \mapsto \sigma \circ x \circ p,
\]
est un quasi-isomorphisme d'algËbres diffÈrentielles graduÈes.
Le diagramme
\[
\diagnum{I''}{
\xymatrix{ & \cc_\infty \cb(m\?,m?) \ar@{->>}[d]^{i^*} \\
\cc_\infty \cb(yr(\?),yr(?)) \ar[r]_{(y\alpha)_*} & \cc_\infty \cb(yr(\?),y?)}}
\]
est ainsi quasi-isomorphe ‡ $(I').$
Les cofibrations Ètant des monomorphismes, la flÍche verticale du diagramme $(I'')$
est une surjection. 
Nous en dÈduisons que les projections canoniques
\[
\cc_\infty \cb(yr(\?),yr(?)) \la P \ra \cc_\infty \cb(m\?,m?),
\]
o˘ $P$ est le produit fibrÈ au-dessus du diagramme $(I'')$ sont des quasi-isomorphismes.
Comme $\cc_\infty \cb(yr(\?),yr(?))$ et $\cc_\infty \cb(m\?,m?)$ sont des algËbres diffÈrentielles
graduÈes unitaires,  $P$ est une algËbre diffÈrentielle graduÈe unitaire et les projections
canoniques ci-dessus sont des morphismes d'algËbres diffÈrentielles graduÈes unitaires.
Nous avons ainsi construit une suite de quasi-isomorphisme d'algËbres diffÈrentielles graduÈes unitaires dans
$\sf C(\mathbb B,\mathbb B)$
\[
\ca' \ra \cc_\infty \cb(yr(\?),yr(?)) \la P \ra \cc_\infty \cb(m\?,m?) \la \cc_\infty \cb(y\?,y?) \la \cb.
\]
Les quasi-isomorphismes d'algËbres Ètant inversibles ‡ homotopie prËs dans la catÈgorie
$\aia$, nous obtenons un $\ai$-quasi-isomorphisme homologiquement unitaire 
\[
g' : \ca' \ra \cb. 
\]
D'aprËs la proposition \ref{proposition_ai-algebres_strictement_unitaires}, il 
existe un $\ai$-morphisme strictement unitaire $g$ homotope ‡ $g'.$ En particulier, $g$ est un
$\ai$-quasi-isomorphisme. C'est une $\ai$-Èquivalence (voir le premier cas.)

{\em Le cas gÈnÈral :} Soit $\ca$ et $\cb$ deux $\ai$-catÈgories strictement unitaires sur
$\mathbb A$ et $\mathbb B$ et $f$ un $\ai$-foncteur tel que
$f_1$ est un quasi-isomorphisme et induit une Èquivalence $H^0 \ca \ra H^0 \cb$.
Choisissons dans chaque classe d'$\ai$-isomorphie $[A]$ de $\ca$ un reprÈsentant $A_0$
et notons $B_0$ son image par $\dot f$. Comme $H^0f : H^0 \ca \ra H^0 \cb$ est une
Èquivalence, nous dÈduisons du lemme (\ref{lemme_ai-isomorphie}) que toute classe
d'$\ai$-isomorphie $[B]$ dans $\cb$ admet un unique reprÈsentant parmi les $B_0$.
Notons $\ca'$ (resp.~$\cb'$) la sous-catÈgorie pleine de $\ca$ (resp.~$\cb$)
formÈe des $A_0$ (resp.~des $B_0$). L'inclusion 
\[
\ca' \ra \ca \quad \Big(\mbox{resp.}\quad \cb' \ra \cb\Big)
\]
est une $\ai$-Èquivalence (voir le deuxiËme cas).
Pour montrer que $f$ est une $\ai$-Èquivalence, il suffit donc de montrer que le foncteur
\[
f' : \ca' \ra \cb'
\]
induit par l'$\ai$-foncteur $f$ est une $\ai$-Èquivalence.
Son application sous-jacente $\dot f'$ est une bijection et $f'_1$ est
un quasi-isomorphisme. Nous sommes donc dans le premier cas et $f'$ est une $\ai$-Èquivalence.
 \findem

%% file: Appendice_A.tex
%
%
%
%

Dans cet appendice, nous rappelons la dÈfinition, due ‡ D. Quillen \cite{Quillen67},
d'une catÈgorie de modËles (fermÈe), quelques
notions fondamentales (objets fibrants, cofibrants, homotopies, foncteurs de Quillen) et
quelques ÈnoncÈs-clÈs. Nous rappelons ensuite
les exemples dont nous avons besoin dans ce manuscrit. Nous renvoyons au livre
de M. Hovey  \cite{Hovey99} et ‡ l'article de W.~Dwyer et J.~Spalinski \cite{DwyerSpalinski95}
pour plus de dÈtails.\\

{\noindent \bf DÈfinitions et propositions}

\begin{definition}{\em Soit $\sf E$ une catÈgorie. 
Un {\em relËvement (de $g$ relatif ‡ $f$)} \index{relËvement}
dans le diagramme
\[
(1) \quad \begin{array}{c}\xymatrix{ A \ar@{->}[d]_i \ar@{->}[r]^f & B
 \ar@{->}[d]^p \\ C \ar@{->}[r]_g  &  D }\end{array}
\]
est un morphisme $\alpha : C \ra B$ tel que les deux triangles du diagramme
\[
 \begin{array}{c} \xymatrix{ A
 \ar@{->}[d]_i \ar@{->}[r]^f & B \ar@{->}[d]^p \\ C \ar@{->}[r]_g
 \ar@{.>}[ur]^\alpha &  D }\end{array}
\]
sont commutatifs.
Soit $i$ et $p$ deux morphismes dans $\sf E.$ Nous dirons que {\em $p$ a
la propriÈtÈ de relËvement ‡ droite par rapport ‡ $i$}  et que {\em $i$ a
la propriÈtÈ de relËvement ‡ gauche par rapport ‡ $p$} si tout diagramme
de la forme $(1)$ admet un relËvement $\alpha$.

Soit $f : X \ra X'$ et $g : Y \ra Y'$ deux morphismes. Le morphisme $f$ est un {\em rÈtract} de $g$ s'il
existe un diagramme commutatif
\[
\xymatrix{
X \ar[d]_f \ar[r] & Y \ar[d]_g \ar[r] & X \ar[d]^f \\
X' \ar[r] & Y' \ar[r] & X' }
\]
tel que les compositions horizontales sont l'identitÈ de $X$ et l'identitÈ de $X'$.
}
\end{definition}

\begin{definition} \label{definition_cmf} {\em 
Une {\em catÈgorie de modËles} \index{categorie@{catÈgorie}!de modËles} est un quadruplet 
\[(\sf E, \weq ,\fib , \cof),\] \indexnotation{weq}
o˘
\english
\begin{itemize}
\item[-] $\sf E$ est une catÈgorie,
\item[-] $\weq$ est une classe de morphismes appelÈs {\em Èquivalences faibles}, 
\index{equivalence faible@{Èquivalence faible}}
\item[-] $\fib$  est une classe de morphismes appelÈs {\em fibrations} \index{fibration}
(elles sont reprÈsentÈes par des flËches ‡ double tÍte $\twoheadrightarrow$),
\item[-] $\cof$ est une classe de morphismes appelÈs {\em cofibrations} \index{cofibration}
(elles sont reprÈ\-sentÈes par des flËches avec une queue $\rightarrowtail$),
\end{itemize}
\francais
tel que les axiomes (CM1) -- (CM5) ci-dessous
sont vÈrifiÈs. Un morphisme appartenant ‡ $\weq \cap \cof$
sera appelÈ une {\em cofibration triviale} \index{triviale (cofibration)} et un morphisme
de $\weq \cap \fib$ sera appelÈ une {\em fibration triviale} \index{triviale (fibration)}.
\begin{description}
\item[](CM1) La catÈgorie $\sf E$ admet toutes les limites finies et toutes les colimites finies.
\item[](CM2) La classe des Èquivalences faibles est {\em saturÈe} \index{saturÈe (classe)}, i.e.~si deux
morphismes parmi $f,g,fg$ sont des  Èquivalences faibles, le
troisiËme l'est aussi.
\item[](CM3) Les trois classes de morphismes sont stables par rÈtracts.
\item[](CM4) {\em relËvement :} \\
{\em  a.} Les cofibrations
ont la propriÈtÈ de relËvement ‡ gauche par rapport aux fibrations triviales,\\
{\em  b.} Les fibrations
ont la propriÈtÈ de relËvement ‡ droite par rapport aux cofibrations triviales.
\item[](CM5) {\em factorisation :}\\
a. Tout morphisme $f:A \ra B$ se factorise en $f =pi$ o˘
$i:A \rightarrowtail A'$ est une cofibration triviale
et $p : A' \twoheadrightarrow B $ est une fibration.\\ {\em
b.} Tout morphisme $f:A \ra B$ se factorise en $f =pi$ o˘
$i:A \rightarrowtail B'$ est une cofibration  et
$p : B' \twoheadrightarrow B $ est une fibration
triviale.
\end{description}
}
\end{definition}
\begin{remarque}{\em 
Nous nous conformons ‡ la terminologie de  \cite{DwyerSpalinski95} en appelant ``catÈgorie de modËles'' ce que
Quillen \cite{Quillen67}, \cite{Quillen69} appelle ``catÈgorie de modËles fermÈe''.
Notons que les axiomes sont auto-duaux.}
\end{remarque}

Soit $(\sf E,\weq ,\fib , \cof)$ une catÈgorie de modËles. On a les propriÈtÈs suivantes :
\english \begin{itemize}
\item[-] La catÈgorie $\sf E$ a un objet initial $\emptyset$ et un
objet final $*$.
\item[-] Les fibrations sont exactement les morphismes ayant
la propriÈtÈ de relËvement ‡ droite par rapport aux cofibrations
triviales.
\item[-]  Les fibrations triviales sont exactement les morphismes ayant
la propriÈtÈ de relËvement ‡ droite par rapport aux cofibrations.
\item[-] Les cofibrations ont les propriÈtÈs de relËvement duales.
\end{itemize} \francais

\begin{definition} \label{definition_homotopie_cmf}
{\em Soit $X$ un objet de $\sf E.$
Un {\em cylindre} \index{cylindre} pour
$X$ est un objet $X\wedge I$  muni de morphismes
$i : X \coprod X \ra X \wedge I$ et $p : X\wedge I \ra X$ tels que
\begin{enumerate}
\item le morphisme $p$ est une Èquivalence faible,
\item la composition $p\circ i : X \coprod X \ra X\wedge I \ra X$ est
le morphisme 
\[
 [\Id, \Id] : X \coprod X \ra X.
\]
\end{enumerate}
Soit $X\wedge I$ un cylindre pour $X$.
Deux morphismes $f,g :X \ra Y$ de $\sf E$ sont {\em $X \wedge I$-homotopes ‡
gauche} si le morphisme
$[f, g] :X\coprod X \ra Y$ se factorise en 
\[
X\coprod X  \arr{i} X\wedge I \arr{H} Y
\]
pour un morphisme $H.$
Un tel morphisme $H$ est appelÈ
une {\em $X \wedge I$-homotopie ‡ gauche}\index{homotopie!‡ gauche} de $f$ ‡ $g$.
Les morphismes $f$ et $g$ sont {\em homotopes ‡ gauche} s'ils sont
$X \wedge I$-homotopes pour un cylindre $X\wedge I$ pour $X$. Nous
Ècrirons alors
\[
f \sim_l g.
\]
La dÈfinition d'un {\em objet de chemins} \index{objet de chemins}
pour $X$ est duale de celle d'un
cylindre pour $X$. 
La notion d'{\em homotopie ‡ droite}  (notÈe $\sim_r$)
est duale de celle d'homotopie ‡ gauche.

}
\end{definition}

\begin{definition}{\em
Un objet $X$ de $\sf E$ est {\em cofibrant} \index{cofibrant} si le morphisme
$\emptyset \ra X$ est une cofibration. Il est  {\em fibrant} \index{fibrant}
si le morphisme $X \ra *$ est une fibration. La sous-catÈgorie pleine des objets
fibrants est notÈe $\sf E_{\sf{f}}$, celle des objets cofibrants $\sf E_{\sf{c}}$  et
celle des objets fibrants et cofibrants est notÈe $\sf E_{\sf{cf}}$. \indexnotation{Ecf}
}
\end{definition}

\begin{definition}{\em
Soit $X$ un objet de $\sf E$. Une {\em rÈsolution cofibrante de $X$}
\index{rÈsolution cofibrante} est  une fibration triviale $X_{\sf c} \twoheadrightarrow
X$, o˘ $X_{\sf c}$ est cofibrant.
Une {\em rÈsolution fibrante de $X$} 
\index{rÈsolution fibrante}
est une cofibration triviale $Y \rightarrowtail
X_{\sf f}$, o˘  $X_{\sf f}$  est fibrant.
}
\end{definition}

Il rÈsulte de l'axiome (CM5) que tout objet admet une rÈsolution cofibrante et une
rÈsolution fibrante.

\begin{lemme} \label{lemme1_cmf}  Soit $X$ un objet
cofibrant et $Y$ un objet fibrant.
\english \begin{itemize}
\item[a.] La relation de $X\wedge I$-homotopie ‡ gauche ne dÈpend pas du choix du cylindre $X \wedge I.$
De mÍme, la relation de $P Y $-homotopie ‡ droite ne dÈpend pas du choix de l'objet de chemins $P Y.$

\item[b.] Les relations d'homotopie ‡ gauche et d'homotopie ‡ droite coÔncident
sur $\sf E(X,Y)$.
On dÈfinit la {\em relation d'homotopie} $\sim$ comme Ègale ‡ ces deux relations.

\item[c.] La relation d'homotopie est une relation d'Èquivalence sur $\sf E(X,Y)$.
\item[d.] Soit $X'$ un objet cofibrant et $Y'$ un objet fibrant.
La relation $f\sim g$ implique $fh \sim gh$ et $h'f \sim h'g$ quels que soient les morphismes
\[
h : X'\ra X\quad \mbox{ et } \quad h' : Y \ra Y'.
\]
\end{itemize} \francais
\findem
\end{lemme}

Le quotient
$\sf E_{\sf{cf}}/\sim$ est donc une catÈgorie. 
On dÈfinit la {\em catÈgorie homotopique}\index{categorie@{catÈgorie}!homotopique}
\indexnotation{HoE}
$\Ho \sf E$ comme la localisation $\sf  E [\weq^{-1}]$ 
de $\sf E$ par rapport ‡ la classe des Èquivalences faibles (voir \cite[I.1]{GabrielZisman67}).

\begin{proposition} \label{proposition_cmf}
\english \begin{itemize}
\item[a.] L'inclusion $\sf E_{\sf{cf}} \ra \sf E$ induit une Èquivalence
\[
\sf E_{\sf{cf}}/\sim \ \arr{} \Ho \sf E.
\]
\item[b.] Soit $X$ et $Y$ deux objets de $\sf E$. Soit $X_{\sf c} \twoheadrightarrow X$
une rÈsolution cofibrante de $X$ et $Y \rightarrowtail Y_{\sf f}$ une rÈsolution
fibrante de $Y$. On a une bijection canonique
\[
\Ho \sf E\,(X,Y) \iso \sf E\,(X_{\sf c},Y_{\sf f})/\sim.
\]
\findem
\end{itemize}\francais
\end{proposition}

{\noindent \bf Equivalence de Quillen}

\begin{definition}\label{foncteurs_de_Quillen}{\em
Soit $\sf E$ et $\sf F$ deux catÈgories de modËles.
Un foncteur $G :\sf E \ra \sf F$ est un {\em foncteur de Quillen ‡ gauche}
\index{foncteur!de Quillen}
s'il admet un adjoint ‡ droite et s'il prÈserve les cofibrations et les
cofibrations triviales. Un foncteur $D : \sf F \ra \sf E$ est un {\em foncteur
de Quillen ‡ droite} s'il admet un adjoint ‡ gauche et s'il prÈserve les fibrations
et les fibrations triviales.
Soit une paire de foncteurs adjoints $(G,D,\phi)$, c'est-‡-dire que $G$ est adjoint ‡ gauche
‡ $D$ et que $\phi$ est une bijection fonctorielle
\[
\Hom_{\sf F} (GX,Y) \arr{\sim} \Hom_{\sf E} (X,DY).
\]
On dira qu'elle est une {\em adjonction de Quillen}
\index{adjonction de Quillen} si $G$ est un foncteur de Quillen ‡ gauche. 
(Ceci implique que $D$ est un foncteur de Quillen ‡ droite.)
Une adjonction de Quillen est une {\em Èquivalence de Quillen} 
\index{equivalence de Quillen@{Èquivalence de Quillen}}
si, pour tout objet cofibrant
$X$ de $\sf E$ et tout objet fibrant $Y$ de $\sf F,$ un morphisme $f : GX \ra Y$ est une
Èquivalence faible si et seulement si $\phi{f} : X \ra DY$ est une Èquivalence faible. Nous
renvoyons ‡ \cite[Sect.~9]{DwyerSpalinski95} pour les dÈtails de la dÈfinition suivante.
}
\end{definition}

\begin{definition}\label{foncteurs_derives}\index{foncteur!derive@{dÈrivÈ}} {\em
Soit $G$ un foncteur de Quillen ‡ gauche. Le {\em foncteur  dÈrivÈ  ‡ gauche de $G$}
est le foncteur
\[
{\bf L} G : \Ho \sf E  \arr{} \Ho \sf F 
\] 
qui envoie un objet $X$ de $\sf E$ sur $GX_{\sf c}$, o˘ $X_{\sf c} \twoheadrightarrow X$ est
une rÈsolution cofibrante de $X$.
Soit $D$ un foncteur de Quillen ‡ droite. Le {\em foncteur dÈrivÈ  ‡ droite de $D$}
est le foncteur
\[
{\bf R} D : \Ho \sf F  \arr{} \Ho \sf E 
\] 
qui envoie un objet $Y$ de $\sf F$ sur $GY_{\sf f}$, o˘ $Y \rightarrowtail Y_{\sf f}$ est
une rÈsolution fibrante de $Y$.}
\end{definition}

\begin{remarque}{\em 
Notons que si un foncteur $G$ (resp.~$D$) comme dans la dÈfinition prÈserve
les Èquivalences faibles, alors il induit un foncteur entre les catÈgories
homotopiques, et ${\bf L} G$ (resp.~${\bf R} D$) est canoniquement isomorphe ‡ ce
foncteur induit.
}
\end{remarque}

\begin{proposition} \label{proposition_equiv_Quillen}
Soit $(G,D,\phi)$ une adjonction de Quillen de $\sf E$ dans $\sf F.$
Les propositions suivantes sont Èquivalentes
\begin{enumerate}
\item[a.] $(G,D,\phi)$ est une Èquivalence de Quillen.
\item[b.] Les foncteurs ${\bf L} G$ et ${\bf R} D$ sont des Èquivalences
inverses l'une de l'autre entre $\Ho \sf E$ et $ \Ho \sf F$.\\
\end{enumerate}
\end{proposition}

{\noindent \bf Exemples de catÈgories de modËles}\\

\begin{exemple}[Complexes de $\sf C$]{\em \label{cmf_complexes} 
Soit $\sf C$ la catÈgorie de base (\ref{categorie_de_base}).
La catÈgorie $\cc \sf C$ de (\ref{complexes}) admet une structure
de catÈgorie de modËles telle que 
\english \begin{itemize}
\item[-] la classe des Èquivalences faibles est la classe $Qis$ des quasi-isomorphis\-mes
(notons que ce sont exactement les morphismes qui sont inversibles ‡ homotopie prËs),
\item[-] les fibrations sont les Èpimorphismes (c'est-‡-dire, les
morphismes dont les composantes sont des Èpimorphismes),
\item[-] les cofibrations sont  les monomorphismes
(c'est-‡-dire,  les morphismes dont les composantes sont des
monomorphismes).
\end{itemize} \francais
Tous les complexes sont fibrants et cofibrants pour cette structure.
La catÈgorie homotopique associÈe est $\ch \sf C.$
}
\end{exemple}

\begin{exemple}[Complexes de chaÓnes non bornÈs]{\em 
Soit $R$ un anneau. Soit $\cc R$ la catÈgorie des complexes
de chaÓnes
\[
\cdots \ra M^{p-1} \ra M^p \ra M^{p+1} \ra \cdots , \quad p \in \Z,
\]
de $R$-modules ‡ droite. Les trois classes de morphismes
suivantes dÈfinissent une structure de catÈgorie de modËles
sur $\cc R$ (voir \cite[Chap. 2]{Hovey99}). 
\english \begin{itemize}
\item[-] Les Èquivalences faibles  sont les quasi-isomorphismes.
\item[-] Les fibrations sont les morphismes $f : X \ra Y$ tels que
$f^n$ est surjectif pour tout $n\in \Z.$
\item[-] Les cofibrations sont les morphismes qui ont la propriÈtÈ
de relËvement ‡ gauche par rapport aux fibrations triviales.
\end{itemize} \francais
Tous les complexes sont fibrants pour cette structure. Si un
complexe $X$ est cofibrant,
alors toutes ses composantes $X^n$, $n\in \Z$, sont
projectives. La rÈciproque est fausse. Cependant, 
si on suppose que le complexe $X$ est bornÈ ‡ droite et
que ses composantes  sont toutes projectives, alors il est
cofibrant.
}
\end{exemple}

%% file: Appendice_B.tex
%
%
%
%
%
%
%
%
%
%
\section{ThÈorie de l'obstruction pour les $\ai$-algËbres}
\label{section_obstruction}

Nous Ètudions la thÈorie de l'obstruction des $\ai$-algËbres.
Soit $\sf C$ une catÈgorie de base telle que dans le chapitre \ref{chapitre_Homot_aialg}.
$(A,m_1,\hdots , m_n)$ une $\rm{A}_n$-algËbre. Il s'agit de mesurer
l'obstruction ‡ l'existence d'un morphisme \index{obstruction}
$m_{n+1}:A^{\ts n+1}\ra A$ tel que $(A,m_1,\hdots ,m_{n+1})$ soit une $\rm{A}_{n+1}$-algËbre
(\ref{extension_nstructure_algebres}).
Soit $A$ et $A'$ deux $\rm{A}_{n+1}$-algËbres. Soit une famille de morphismes graduÈs
\[
f_i : A \tp i \ra B, \hspace{1cm} 1 \leq i \leq n,
\]
dÈfinissant un $\rm{A}_n$-morphisme $A \ra A'.$
Nous mesurons ensuite l'obstruction ‡ l'existence d'un morphisme $f_{n+1} : A\tp{n+1} \ra A'$
tel que les $f_i$, $1 \leq i \leq n+1$, dÈfinissent un $\rm A_{n+1}$-morphisme $A \ra A'$
(\ref{extension_nstructure_morphismes}).
Nous montrerons que cette obstruction est
fonctorielle  par rapport aux $\rm{A}_{n+1}$-morphismes stricts (\ref{fonctorialite_obstruction}).

L'Ètude des obstructions est un outil classique, voir par exemple T.~Kadei\-shvili \cite{Kadeishvili80},
A.~ProutÈ \cite{Proute85}. Elle doit son existence au fait que l'opÈrade des $\ai$-algËbres
est un modËle cofibrant minimal au sens de M.~Markl \cite{Markl96}
pour l'opÈrade des algËbres associatives. Nous n'adoptons pas ici ce point de vue
lui prÈfÈrant une approche naÔve.\\
%
%
%
%
%
%
%
%
%
%

{\bf Les $\ai$-algËbres}\\

Soit $V$ un objet graduÈ.
Soit des morphismes graduÈs
\[
b_i : V\tp i \ra V,\hspace{1cm}1 \leq i \leq n+1,
\]
de degrÈ $+1$.
Notons $b$ la codÈrivation de $\pctr{n+1}V$ donnÈe par la suite
\[
(b_1,\hdots ,b_n,b_{n+1}).
\]
Posons
\[
c(b_2,\hdots ,b_n) = \sum_{\substack{2 \leq i \leq n}} b_i(\Id \tp j \ts b_k \ts \Id \tp l)
\]
o˘ les entiers $j,k,l$ vÈrifient $j+k+l = n+1$ et $j +1+ l = i$.
Rappelons que $i_1$ et $p_{n+1}$ dÈsignent les morphismes canoniques
\[
V \arr{} \pctr{n+1}V \hspace{1cm}\mbox{et} \hspace{1cm} \pctr{n+1} V \arr{} V\tp{n+1}.
\]

\begin{lemme}
Supposons que la codÈrivation de la cogËbre $\pctr{n}V$ donnÈe par la suite
\[
(b_1,\hdots ,b_n)
\]
est une diffÈrentielle.
\begin{enumerate}
\item[a.]La codÈrivation
\[
b^2 : \pctr{n+1}V \arr{} \pctr{n+1}V
\]
est Ègale ‡ $i_1 \circ \zeta \circ p_{n+1},$ o˘ $\zeta : V \tp{n+1} \ra V$ est donnÈ
par 
\[
\zeta = b_1 b_{n+1} + b_{n+1}b_1 + c(b_2,\hdots ,b_{n}) ;
\]
ici la derniËre occurrence de $b_1$ dÈsigne la diffÈrentielle de $(V,b_1)\tp{n+1}.$
\item[b.] Le morphisme graduÈ $c(b_2,\hdots,b_{n})$ est un cycle de
\[
(\Hom_{\gr \sf C} (V\tp{n+1},V),\delta),
\]
o˘ la diffÈrentielle $\delta$ est induite par celle du complexe $(V,b_1).$
\end{enumerate}
En particulier, la codÈrivation $b$ est une diffÈrentielle si et seulement
si le cycle $c(b_2,\hdots , b_{n})$ est Ègal au bord $-\del{b_{n+1}}.$
\end{lemme}
\dem 
 {\it a.} Notre hypothËse implique que le carrÈ $b^2$ se factorise par $p_{n+1}$.
L'image de la comultiplication
$\Delta$ est incluse dans \[
\pctr{n} V \ts \pctr{n} V \subset \pctr{n+1} V \ts \pctr{n+1} V.
\]
On a donc l'ÈgalitÈ
\[
\Delta {b}^2 = (\Id \ts
{b}^2 + {b}^2 \ts \Id )\Delta = 0.
\]
On en dÈduit que l'image de $b^2$ est incluse dans $ \ker \Delta = V.$
Ceci nous donne la factorisation par~$i_1.$ Un calcul direct nous donne
la formule pour $\zeta.$

{\it b.}  D'aprËs le premier point, on a
\[
b_1 \circ {b}^2 =  {b} \circ {b}^2 = {b}^2
\circ {b} = {b}^2 \circ b_1,
\]
o˘ la derniËre occurrence de $b_1$ dÈsigne la diffÈrentielle de $(V,b_1)\tp{n+1}.$
Ceci montre que $\zeta$ est un cycle dans
le complexe
\[
(\Hom_{\gr \sf C} (V^{\ts n+1},V),\delta).
\]
Comme nous avons
\[
\zeta = \del{b_{n+1}} + c(b_2,\hdots,b_n)
\]
il en est de mÍme pour $c(b_2,\hdots,b_n)$.\hfill $\Box$\\
\begin{corollaire} \label{extension_nstructure_algebres}
Soit $(A,m_1)$ un complexe. Soit des morphismes graduÈs 
\[
m_i :  A \tp i \ra A, \hspace{1cm}2 \leq i \leq n+1
\]
de degrÈ
$2-i$. Supposons que les morphismes $m_i$, $1 \leq i \leq n$, dÈfinissent
une structure de $\rm{A}_n$-algËbre sur $A.$ La sous-expression
\[
\sum_{i,k \neq 1} (-1)^{jk+l}m_i(\Id^{\ts j}\ts m_k\ts \Id^{\ts l})
\]
de l'Èquation $(*_{n+1})$ de (\ref{definition_ai-algebre})
dÈfinit un cycle de $(\Hom_{\gr \sf C}(A\tp{n+1},A),\delta).$
Nous le notons $r(m_2,\hdots, m_{n}).$
L'Èquation $(*_{n+1})$ se rÈcrit alors
\[
r(m_2,\hdots, m_{n}) + \del{m_{n+1}} = 0.
\]
\end{corollaire}

\dem Nous appliquons le lemme prÈcÈdent ‡ l'espace graduÈ $V = SA$ et aux
morphismes graduÈs $b_i$ dÈfinis ‡ l'aide des bijections $b_i \lra m_i.$ 
Ces mÍmes bijections envoient
le morphisme $r(m_2,\hdots, m_{n})$ sur le morphisme
 $c(b_2, \hdots ,b_n)$ et le morphisme $\del{m_{n+1}}$ sur
$\del{b_{n+1}}.$\findem \\
%
%
%
%
%
%
%
%
%
%

{\bf Les $\ai$-morphismes d'$\ai$-algËbres}\\

Les lemmes suivant se montrent de maniËre similaire. \\

Soit $V$ et $W$ deux objets graduÈs. Soit $b$ et $b'$ des diffÈrentielles
de cogËbres sur les cogËbres $\pctr{n+1} V$ et $\pctr{n+1} W.$ Soit une
famille de morphismes graduÈs
\[
F_i : V \tp i \ra W, \hspace{1cm} 1 \leq i \leq n+1,
\]
de degrÈ $0.$
Soit $F$ le morphisme de cogËbres
\[
\pctr{n+1} V \arr{} \pctr{n+1} W
\]
qui relËve les $F_i.$ Posons
\[
c(F_1,\hdots,F_n) = \sum_{k\geq 2} F_i(\Id^{\ts j}\ts b_k\ts \Id^{\ts l}) - \sum_{r\geq 2} 
b'_r (F_{i_1}\ts \hdots \ts F_{i_r}),
\]
o˘ les entiers $j,k,l$ de la somme de gauche vÈrifient $j+k+l = n+1$ et $j+1+l =i$, et
o˘ la somme des entiers $i_r$ de la somme de droite vaut $n+1.$

\begin{lemme} 
Supposons que le morphisme 
\[
F \prim{n} : \pctr{n} V \ra \pctr{n} W
\]
induit par $F$ dans les $n$-primitifs est compatible aux diffÈrentielles.
\begin{enumerate}
\item[a.] La $(F,F)$-codÈrivation
\[
b'F - Fb : \pctr{n+1}V \arr{} \pctr{n+1}W
\]
est Ègale ‡ $i_1 \circ \zeta \circ p_{n+1}$, o˘ $\zeta : V \tp{n+1} \ra W$ est donnÈ
par 
\[
\zeta = b_1 F_{n+1} + F_{n+1}b_1 + c(F_1,\hdots ,F_{n}) ;
\]
ici la derniËre occurrence de $b_1$ dÈsigne la diffÈrentielle de $(V,b_1)\tp{n+1}.$
\item[b.] Le morphisme graduÈ $c(F_1,\hdots,F_{n})$ est un cycle de
\[
(\Hom_{\gr \sf C} (V\tp{n+1},W),\delta),
\]
o˘ la diffÈrentielle $\delta$ est induite par celles des complexes $(V,b_1)$ et $(W,b'_1).$
\end{enumerate}
En particulier, le morphisme $F$ est compatible aux diffÈrentielles de cogËbres si et
seulement si on a
\[
\del{F_{n+1}} + c(F_1,\hdots ,F_n) = 0.
\]
\findem
\end{lemme}
%
%
%
%
%
%
%
%
%
%

Regardons maintenant le comportement de l'obstruction par rapport ‡ la
composition des $\rm{A}_{n+1}$-morphismes.

Soit $V'$ et $W'$ deux objets graduÈs.
Soit $d$ et $d'$ deux diffÈrentielles de cogËbres sur les cogËbres $\pctr{n+1}V'$ et
$\pctr{n+1}W'$. Soit deux morphismes de cogËbres diffÈrentielles graduÈes
\[
G : \pctr{n+1}V'  \arr{} \pctr{n+1}V \hspace{1cm} \mbox{et} \hspace{1cm} H : \pctr{n+1}W
\arr{} \pctr{n+1}W'.
\]
Des calculs directs nous donnent le lemme suivant.

\begin{lemme} 
\begin{enumerate}
\item[a.] On a l'ÈgalitÈ
\[
\hspace{-1cm} c({F_1},\hdots,{F_n}) \circ {G_1}\tp{n+1} + {F_{1}}\circ
{\delta(G_{n+1})} = c({(FG)_1},\hdots ,{(FG)_n})
\]
de morphismes de $(V')\tp{n+1}$ dans $W.$
\item[b.]  On a l'ÈgalitÈ
\[
\hspace{-1cm}{\delta(H_{n+1})}\circ {F_{1}}\tp{n+1} + {H_1} \circ
c({F_1},\hdots,{F_n}) = c({(HF)_1},\hdots ,{(HF)_n})
\]
de morphismes de $V \tp{n+1}$ dans $W'.$
\findem
\end{enumerate}
\end{lemme}

\begin{corollaire} \label{extension_nstructure_morphismes}
Soit $A$ et $B$ deux $\rm{A}_{n+1}$-algËbres. Soit des morphismes graduÈs
\[
f_i : A \tp i \ra B,\hspace{1cm} 1 \leq i \leq n+1,
\]
de degrÈ $1-i$. Supposons que les morphismes $f_i$, $1\leq i \leq n$,
dÈfinissent un  $A_n$-morphisme $A \ra B$. La sous-expression
\[
\sum_{k\neq 1} (-1)^{jk+l}f_i(\Id^{\ts j}\ts m_k\ts \Id^{\ts l}) - \sum_{r\neq 1} (-1)^s
m_r (f_{i_1}\ts \hdots \ts f_{i_r})
\]
de l'Èquation $(**_{n+1})$ de (\ref{definition_ai-morphisme})
dÈfinit un cycle dans $(\Hom_{\gr \sf C} (A\tp{n+1},B),\delta).$ Nous le
notons $r(f_1,\hdots ,f_n).$
L'Èquation $(**_{n+1})$ se rÈcrit
\[
r(f_1,\hdots ,f_n) + \del{f_{n+1}} = 0.
\]
\findem
\end{corollaire}

\begin{corollaire} \label{fonctorialite_obstruction}
Soit $A'$ et $B'$ deux $\rm{A}_{n+1}$-algËbres. Soit $g : A' \ra A$ et 
$h : B \ra B'$ deux $\rm{A}_{n+1}$-morphismes stricts.
On a les ÈgalitÈs de morphismes
\begin{enumerate}
\item $ r({f_1},\hdots,{f_n}) \circ {g_1}\tp{n+1} =
r({(fg)_1},\hdots ,{(fg)_n}), $
\item $  {h_1} \circ r({f_1},\hdots,{f_n}) =
r({(hf)_1},\hdots ,{(hf)_n}).$
\end{enumerate}
\end{corollaire}
L'obstruction est donc fonctorielle par rapport aux morphismes
stricts. \\

\dem C'est la traduction du lemme \ref{fonctorialite_obstruction} appliquÈ
aux constructions bar des algËbres $A,A',B$ et $B'$. Les morphismes $g$ et $h$
Ètant stricts, on a $H_{n+1} = 0$ et $G_{n+1} = 0.$ 
Les Èquations de (\ref{fonctorialite_obstruction}) se traduisent alors par
celles du corollaire.
\findem

\section{ThÈorie de l'obstruction pour les polydules}
\label{section_obstruction_polydules}

Les dÈmonstrations de cette section Ètant presque identiques ‡
celles de la section \ref{construction_bar_cobar}, nous nous contentons
d'Ènoncer les rÈsultats. Soit $\sf C$ et $\sf C'$ les catÈgories de base de
la section \ref{section_reduction_augmentation}.
\\

\begin{lemme} \label{obstruction_extension_polydules}
Soit $A$ une $\rm{A}_n$-algËbre. Soit $(M,m_1^M)$ un complexe. Soit des morphismes
graduÈs
\[
m_i^M : M \ts A\tp{i-1} \ra M, \quad 2 \leq i \leq n+1,
\]
de degrÈ $2-i$. Supposons que les morphismes $m_i$, $1 \leq i \leq n,$ dÈfinissent une
structure de $\rm{A}_n$-module sur $M.$ La sous-expression
\[
\sum_{i,k \neq 1} (-1)^{jk+l}m_i(\Id^{\ts j}\ts m_k\ts \Id^{\ts l})
\]
de l'Èquation $(*'_{n+1})$ de (\ref{definition_ai-module})
dÈfinit un cycle de $(\Hom_{\gr \sf C'}(M \ts A\tp{n},M),\delta),$ o˘ $\delta$ est induit
par $m^A_1$ et $m^M_1.$
Nous le notons $r(m_2,\hdots, m_{n}).$
L'Èquation $(*'_{n+1})$ se rÈcrit alors
\[
r(m_2,\hdots, m_{n}) + \del{m_{n+1}} = 0.
\]
\findem
\end{lemme}

\begin{lemme} \label{extension_nstructure_morphismes_mod}
Soit $A$ une $\rm{A}_n$-algËbre.
Soit $M$ et $N$ deux $\rm{A}_{n+1}$-modules sur $A$. Soit des morphismes graduÈs
\[
f_i : M \ts A \tp{i-1} \ra N,\hspace{1cm} 1 \leq i \leq n+1,
\]
de degrÈ $1-i$. Supposons que les morphismes $f_i$, $1\leq i \leq n$,
dÈfinissent un  $A_n$-morphisme $M \ra N$. La sous-expression
\[
\sum_{k \neq 1} (-1)^{jk+l}f_i(\Id^{\ts j}\ts m_k\ts \Id^{\ts l}) = \sum_{s \neq 0}
m_{s+1} (f_{r}\ts \Id\tp s)
\]
de l'Èquation $(**'_{n+1})$ de (\ref{definition_ai-morphisme_mod})
dÈfinit un cycle dans 
\[(\Hom_{\gr \sf C'} (M \ts A\tp{n},N),\delta).
\]
Nous le
notons $r(f_1,\hdots ,f_n).$
L'Èquation $(**'_{n+1})$ se rÈcrit alors
\[
r(f_1,\hdots ,f_n) + \del{f_{n+1}} = 0.
\]
\findem
\end{lemme}

\begin{lemme} \label{fonctorialite_obstruction_mod}
Soit $A$ une $\rm{A}_n$-algËbre.
Soit $M'$ et $N'$ deux $\rm{A}_{n+1}$-modules. Soit $g : M' \ra M$ et 
$h : N \ra N'$ deux $\rm{A}_{n+1}$-morphismes stricts.
On a les ÈgalitÈs de morphismes
\begin{enumerate}
\item $ r({f_1},\hdots,{f_n}) \circ {g_1} \ts \Id \tp{n} =
r({(fg)_1},\hdots ,{(fg)_n}), $
\item $  {h_1} \circ r({f_1},\hdots,{f_n}) =
r({(hf)_1},\hdots ,{(hf)_n}).$
\end{enumerate}
\findem
\end{lemme}

\section{ThÈorie de l'obstruction pour les bipolydules}
\label{section_obstruction_bipolydules}

Les dÈmonstrations de cette section
sont omises car elles sont similaires ‡ celles de la section \ref{section_obstruction}.
Soit $\sf C$, $\sf C'$ et $\sf C''$ les catÈgories de base de la section
\ref{section_categorie_derivee_augmentee_bipolydules}.\\

Soit $A$ et $A''$ deux $\ai$-algËbres dans $\sf C$ et $\sf C''.$
Dans ce qui suit, $r$ et $t$ dÈsignent deux entiers $\geq 0$ et
$\ce$ dÈsigne l'ensemble des couples d'entiers $(i,j)$
tels que $0 \leq i \leq r$ et $0 \leq j \leq t-1$, ou, $0 \leq i \leq r-1$ et $0 \leq j \leq t$
(voir le graphique ci-dessous). L'ensemble $\ce'$ est Ègal ‡ $\ce \setminus (0,0).$\\
\begin{tabular}{ccc}
& & \\
\input{grid1.pstex_t}
&  &
\input{grid2.pstex_t} \\
ensemble $\ce$ & & ensemble $\ce'$
\end{tabular}\vspace{.7cm}

Soit $M$ un objet diffÈrentiel graduÈ de $\sf C'.$ On note sa diffÈrentielle
$m_{0,0}$. 
Soit
\[
m_{i,j} : A \tp i \ts M \ts A''\tp j \ra M \ra M, \quad 0 \leq j \leq t, \quad 0 \leq i \leq r,\quad
(i,j) \neq (0,0),
\]
un morphisme graduÈ de degrÈ $1-i-j$ dans $\sf C'.$

\begin{lemme}
Supposons que les morphismes $m_{i,j}$, $(i,j) \in \ce',$ vÈrifient les Èquations
$(*''_{k,l})$, $(k,l) \in \ce,$ de la dÈfinition \ref{definition_A_n-A_m-bimodules}.
La sous-expression
\[
\sum_{\ast \notin \{1, (0,0), (r,t)\}} 
(-1)^{j + i(|m_\ast|)} m_{\bullet,\bullet}(\Id \tp i \ts m_{\ast} \ts \Id \tp j)
\]
de l'Èquation $(*''_{r,t})$ dÈfinit un cycle de
\[
\Hom_{\gr \sf C'}(A \tp r \ts M \ts A'' \tp t, M).
\]
On le note $c\big(m_{i,j},\  (i,j) \in \ce'\big).$
Les morphismes $m_{i,j}$, $0 \leq j \leq t,$ $0 \leq i \leq r,$ vÈrifient l'Èquation
$(*''_{r,t})$ si et seulement si on a l'ÈgalitÈ
\[
\delta (m_{r,t}) = c\big(m_{i,j},\ (i,j) \in \ce'\big).
\]
\findem
\end{lemme}

Soit $M$ et $N$ deux $A$-$A''$-bipolydules dans $\sf C'$.
Soit
\[
f_{i,j} : A \tp i \ts M \ts A''\tp j \ra M \ra M, \quad 0 \leq j \leq t, \quad 0 \leq i \leq r,
\]
un morphisme graduÈ de degrÈ $-i-j$ dans $\gr \sf C'.$

\begin{lemme}
Supposons que les morphismes $f_{i,j}$, $(i,j) \in \ce,$ vÈrifient les Èquations
$(**''_{k,l})$, $(k,l) \in \ce,$ de la dÈfinition \ref{definition_A_n-A_m-bimodules}.
La sous-expression
\[
\begin{array}{c}
\sum_{(\alpha,\beta)\neq (0,0)}
(-1)^{\alpha(-i-j)} m_{\alpha,\beta}(\Id \tp \alpha \ts f_{k,l} \ts \Id \tp \beta) = \hspace*{4cm} \\
\hspace*{1cm} \sum_{\ast \notin \{1, (0,0)\}} (-1)^{j + i(|m_\ast|)}
f_{\bullet,\bullet}(\Id \tp i \ts m_{\ast} \ts \Id \tp j)
\end{array}
\]
de l'Èquation $(**''_{r,t})$ est un cycle de
\[
\Hom_{\gr \sf C'}(A \tp r \ts M \ts A'' \tp t, N).
\]
On le note $c'\big(f_{i,j},\ (i,j) \in \ce\big).$
Les morphismes $f_{i,j}$, $0 \leq j \leq t,$ $0 \leq i \leq r,$ vÈrifient l'Èquation
$(**''_{r,t})$ si et seulement si on a l'ÈgalitÈ
\[
\delta (f_{r,t}) = c'\big(f_{i,j},\  (i,j) \in \ce\big).
\] \findem
\end{lemme}

On regarde maintenant la compatibiltÈ de l'obstruction aux morphismes stricts.

Soit $M'$ et $N'$ deux $A$-$A'$-bipolydules et
\[
g : M' \ra M \quad \mbox{et} \quad h : N \ra N'
\]
deux $\ai$-morphismes stricts de bipolydules donnÈs par des morphismes graduÈs de degrÈ $0$ dans $\gr \sf C'$
\[
g_{0,0} : M' \ra M \quad \mbox{et} \quad h_{0,0} : N \ra N'.
\]
On dÈfinit les morphismes
\[
(f \circ g)_{i,j} \quad \mbox{et} \quad (h \circ f)_{i,j}, \quad 0 \leq j \leq t, \quad 0 \leq i \leq r,
\]
par les mÍmes formules que celles donnant la compositions des morphismes de bipolydules.

\begin{lemme} On a les ÈgalitÈs suivantes :
\english \begin{enumerate}
\item $c'\big(f_{i,j},\ (i,j) \in \ce\big) \circ (\Id \tp r \ts {g_{0,0}} \ts \Id \tp t) =
c'\big((f \circ g)_{i,j},\ (i,j) \in \ce\big), $
\item ${h_{0,0}} \circ c'\big(f_{i,j},\ (i,j) \in \ce\big) =
c'\big((h \circ f)_{i,j},\ (i,j) \in \ce\big).$
\end{enumerate}\francais \findem
\end{lemme}

\section{Cohomologie de Hochschild et thÈorie de l'obstruction pour les $\ai$-structures mini\-males}
\label{section_obstruction_unite}

Dans cette section, on rappelle la cohomologie de Hochschild d'une algËbre graduÈe ‡ coefficients
dans un bimodule graduÈ. Nous Ètablissons ensuite une thÈorie de l'obstruction des $\ai$-algËbres
minimales (resp.~des $\ai$-morphismes entre $\ai$-algËbres minimales et des homotopies entre ces
$\ai$-morphismes).\\

{\noindent \bf Rappel sur la cohomologie de Hochschild}\\

Soit $\sf C$ une catÈgorie de base telle que dans le chapitre \ref{chapitre_Homot_aialg}.
Soit $A \in \gr \sf C$ une algËbre associative.
On considËre $A$ comme une $\ai$-algËbre dont $m_2 = \mu^A$ et $m_i = 0$ pour tout $i \neq 2.$
Rappellons que $(BA)\+$ est la construction bar co-augmentÈe de $A$. Soit
$\coder ((BA)\+)$ l'espace des codÈrivations $(BA)\+ \ra (BA)\+$. Il est graduÈ par le degrÈ des codÈrivations.
L'application 
\[
\delta : D \mapsto b^A \circ D - (-1)^{|D|} D \circ b^A,
\]
o˘ $b^A$ est la diffÈrentielle de $(BA)\+$ et $D$ est de degrÈ $|D|$,
dÈfinit une diffÈrentielle sur $\coder ((BA)\+).$
Nous montrons (comme dans le lemme \ref{cogebres_tensorielles}) que nous avons une bijection naturelle
\[
\begin{array}{rcl}
\coder ((BA)\+) & \arr{\sim} & \Hom_{\sz \gr \sf C}((BA)\+,SA)\\
D & \mapsto & p_1 \circ D.
\end{array}
\]
Ainsi, une codÈrivation $D$ est dÈterminÈe par les composantes de $p_1 \circ D$
\[
D_i : (SA) \tp i \ra SA, \quad i\geq 0.
\]
Les bijections $b_i \lra m_i$, $i \geq 1$, de la section \ref{construction_bar_cobar} (complÈtÈe
de la bijection qui associe au morphisme $b_0 : e \ra SA$ le morphisme $m_0 = - \si b_0 : e \ra A$)
nous donnent une bijection
\[
\Hom_{\sz \gr \sf C}((BA)\+,SA) \arr{\sim} \prod_{i\geq 0} \Hom_{\sz \gr \sf C}(A\tp i,A).
\] 
Le {\em complexe de Hochschild} \index{Hochschild!complexe} est dÈfini par ces bijections comme
\[
C(A,A) = S^{-1}\prod_{i\geq 0}\Hom_{\sz \gr \sf C}(A\tp i,A).
\]
Sa diffÈrentielle $\delta_{Hoch}$ \indexnotation{delHoch}
envoie un morphisme $f : A \tp n \ra A$ de degrÈ $r$
sur le morphisme 
\[
\delta_{Hoch}(f) : A \tp{n+1}  \arr{}  A
\]
donnÈ par la somme
\[
\sum (-1)^{r+n+k}f_{i}(\Id \tp j \ts \mu \ts \Id \tp k) + (-1)^{r + n + 1}\mu (\Id \ts f_i) + (-1)^{r}
\mu (f_i \ts \Id).
\] 
Si le degrÈ de $f$ est nul, nous retrouvons la dÈfinition habituelle
(voir par exemple \cite[Chap.~IX]{Cartan56}).
Soit $M\in \gr \sf C$ un $A$-$A$-bimodule. Le {\em complexe de Hochschild ‡ coefficients dans $M$} est l'espace
\indexnotation{C(A,M)}
\[
C(A,M) = \prod_{i \geq 0} \Hom_{\sz \gr \sf C}(A \tp i, M),
\]
sa graduation est induite par la graduation de l'espace
\[
\prod_{i\geq 0}\Hom_{\sz \gr \sf C}((SA)\tp i,SM)
\]
et sa diffÈrentielle $\delta_{Hoch}$ \indexnotation{delta_Hoch}
est dÈfinie par la mÍme formule que prÈcÈdemment.
La {\em cohomologie de Hochschild de $A$ ‡ coefficients dans $M$} est la cohomologie de $C(A,M)$.
Si $A$ est unitaire,
le complexe $C(A,M)$ est homotopiquement Èquivalent au {\em sous-complexe de Hochschild rÈduit}
\index{Hochschild!complexe rÈduit} \indexnotation{bC(A,M)} (voir \cite[Chap.~IX]{Cartan56})
\[
\b C(A,M) = \prod_{i \geq 0} \Hom_{\sz \gr \sf C}(\b A \tp i, M),
\]
o˘ $\b A$ est le conoyau de l'unitÈ de $A.$ La diffÈrentielle de $\b C(A,M)$
est induite par celle de $C(A,M).$\\

{\noindent \bf Obstruction ‡ l'extension d'une $\rm{A}_n$-algËbre minimale en une
$\rm{A}_{n+1}$-algËbre minimale}

\begin{lemme}\label{lemme_obstruction_structure_minimale}
Soit $A$ une algËbre graduÈe de $\gr \sf C.$ Soit des morphismes graduÈs
\[
m_i : A \tp i \ra A, \quad 3 \leq i \leq n,
\]
de degrÈ $2-i.$ Nous posons $m_2 = \mu^A$.
Supposons que les morphismes $m_i$, $2 \leq i \leq n - 1,$ dÈfinissent une
structure de $\rm{A}_{n}$-algËbre minimale sur $A$. La sous-expression
\[
\sum_{i,k \notin \{1,2\}} (-1)^{j+kl}m_i(\Id^{\ts j}\ts m_k\ts \Id^{\ts l})
\]
de l'Èquation $(*_{n+1})$ de (\ref{definition_ai-algebre})
dÈfinit un cycle de Hochschild. Nous le notons $r(m_3,\hdots ,m_{n-1}).$
L'Èquation $(*_{n+1})$ se rÈcrit alors
\[
\delta_{Hoch}(m_{n}) + r(m_3,\hdots ,m_{n-1}) = 0.
\]
\end{lemme}

\dem 
Soit la suite des morphismes $b_i$, $2 \leq i \leq n$, donnÈs par les bijections
$b_i \lra m_i$ (voir \ref{construction_bar_cobar}).
 On note $D$ la codÈrivation de $(BA)\+$
telle que les composantes de $p_1 \circ D$ sont donnÈes par la suite
\[
(0,0,b_2, \hdots ,b_{n-1},b_{n},0,\hdots).
\]
Comme les $m_i$, $2 \leq i \leq n-1,$ dÈfinissent une
structure de $\rm{A}_{n}$-algËbre minimale, le carrÈ de la codÈrivation $D$ restreint
‡ la sous-cogËbre $\pctr{n}SA$ est nul. On en dÈduit que la composition
\[
\Delta \circ {D^2} = (\Id \ts D^2 + D^2 \ts \Id) \circ \Delta
\]
s'annule sur le sous-espace $(SA)\tp{n+1}.$ Il s'ensuit que l'image par $D^2$ du sous-espace
$(SA)\tp{n+1}$ est contenue dans $\ker \Delta = SA$ et que celle du
sous-espace $(SA)\tp{n+2}$ est contenue dans $(SA)\tp 2 \oplus SA$.
Ainsi, sur le sous-espace $(SA)\tp{n+2}$, on
a l'ÈgalitÈ
\[
D^2 \circ b_2 = D^3 = b_2 \circ D^2.
\]
Ceci montre que l'ÈlÈment
\[
D^2|_{(SA)\tp{n+1}} \in \Hom ((SA)\tp{n+1},SA)
\]
est un cycle. La premiËre assertion du lemme est dÈduite du fait que
l'ÈlÈment
\[
\si (D^2|_{(SA)\tp{n+1}}) 
\]
correspond ‡ l'ÈlÈment $r(m_3,\hdots,m_{n-1})$ par l'isomorphisme de complexes
\[
S^{-1}\Hom_{\sz \gr \sf C}((BA)\+,SA) \arr{\sim} C(A,A).
\]
La derniËre assertion du lemme est immÈdiate.\findem\\

{\noindent \bf Obstruction ‡ l'extension d'un $\rm{A}_n$-morphisme entre $\ai$-algËbres minimales
en un $\rm{A}_{n+1}$-morphisme}

\begin{lemme}\label{lemme_obstruction_structure_minimale_morphisme}
Soit $A$ et $A'$ deux $\ai$-algËbres minimales. Soit
\[
g : (A,m_2) \ra (A',m_2')
\]
un morphisme d'algËbres graduÈes. Soit des morphismes graduÈs
\[
f_i : A \tp i \ra A', \quad 2 \leq i \leq n,
\]
de degrÈ $1-i.$ Nous posons $f_1 = g$.
Supposons que les morphismes $f_i$, $1 \leq i \leq n-1,$ dÈfinissent un
 $\rm{A}_{n}$-morphisme $A \ra A'$. La sous-expression
\[
\sum_{k\notin \{1,2\}} (-1)^{jk+l}f_i(\Id^{\ts j}\ts m_k\ts \Id^{\ts l}) - \sum_{r\notin \{1,2\}} (-1)^s
m'_r (f_{i_1}\ts \hdots \ts f_{i_r})
\]
de l'Èquation $(**_{n+1})$ de (\ref{definition_ai-morphisme})
dÈfinit un cycle de Hochschild dans $C(A,A')$; la structure de  $A$-bimodule sur
$A'$ est donnÈe par $g$. Nous notons ce cycle $r(f_2,\hdots ,f_{n-1}).$
L'Èquation $(**_{n+1})$ se rÈcrit alors
\[
\delta_{Hoch}(f_{n}) + r(f_2,\hdots ,f_{n-1}) = 0.
\]
\findem
\end{lemme}

{\noindent \bf Obstruction ‡ l'extension d'une $\rm{A}_n$-homotopie entre $\ai$-morphismes
d'$\ai$-algËbres minimales
en une $\rm{A}_{n+1}$-homotopie}

\begin{lemme}\label{lemme_obstruction_structure_minimale_homotopie}
Soit $A$ et $A'$ deux $\ai$-algËbres minimales.
Soit $f$ et $g$ deux $\ai$-morphismes $A \ra A'$. Soit
des morphismes graduÈs
\[
h_i : A \tp i \ra A', \quad 2 \leq i \leq n,
\]
de degrÈ $-i.$ Posons $h_1 = 0.$
Supposons que les morphismes $h_i$, \hbox{$1 \leq i \leq n-1,$} dÈfinissent une
homotopie entre $f$ et $g$ considÈrÈs comme $\rm{A}_{n}$-morphisme \hbox{$A \ra A'$}
(on a alors $f_1 = g_1$).
La sous-expression
\[
 \Bigg( -  \sum_{r+1+t \not \in \{1,2\}} (-1)^s
m_{r+1+t} (f_{i_1}\ts \hdots \ts f_{i_r} \ts h_{k} \ts g_{j_1}\ts \hdots \ts g_{i_t}) +
\]
\[
 - \sum_{k \not \in \{1,2\}} (-1)^{jk+l}h_i(\Id^{\ts j}\ts m_k\ts \Id^{\ts l}) + f_{n+1} - g_{n+1} \Bigg)
\]
de l'Èquation $(***_{n+1})$ de la dÈfinition \ref{definition_homotopie_ai-morphismes}
dÈfinit un cycle de Hochschild dans $C(A,A')$; la structure de $A$-bimodule sur
$A'$ est donnÈe par $f_1$ et $g_1$. Nous notons ce cycle $r(h_2,\hdots ,h_{n-1}).$
L'Èquation $(***_{n+1})$ se rÈcrit alors
\[
\delta_{Hoch}(h_{n}) + r(h_2,\hdots ,h_{n-1}) = 0.
\]
\findem
\end{lemme}

%% file: grid1.pstex_t
\begin{picture}(0,0)%
\includegraphics{grid1.pstex}%
\end{picture}%
\setlength{\unitlength}{3947sp}%
\begingroup\makeatletter\ifx\SetFigFont\undefined%
\gdef\SetFigFont#1#2#3#4#5{%
  \reset@font\fontsize{#1}{#2pt}%
  \fontfamily{#3}\fontseries{#4}\fontshape{#5}%
  \selectfont}%
\fi\endgroup%
\begin{picture}(2625,1806)(376,-1486)
\put(526,-1486){\makebox(0,0)[lb]{\smash{\SetFigFont{12}{14.4}{\rmdefault}{\mddefault}{\updefault}{\color[rgb]{0,0,0}$0$}%
}}}
\put(2326,-1486){\makebox(0,0)[lb]{\smash{\SetFigFont{12}{14.4}{\rmdefault}{\mddefault}{\updefault}{\color[rgb]{0,0,0}$r$}%
}}}
\put(3001,-1336){\makebox(0,0)[lb]{\smash{\SetFigFont{12}{14.4}{\rmdefault}{\mddefault}{\updefault}{\color[rgb]{0,0,0}$i$}%
}}}
\put(526,164){\makebox(0,0)[lb]{\smash{\SetFigFont{12}{14.4}{\rmdefault}{\mddefault}{\updefault}{\color[rgb]{0,0,0}$j$}%
}}}
\put(376,-1336){\makebox(0,0)[lb]{\smash{\SetFigFont{12}{14.4}{\rmdefault}{\mddefault}{\updefault}{\color[rgb]{0,0,0}$0$}%
}}}
\put(376,-436){\makebox(0,0)[lb]{\smash{\SetFigFont{12}{14.4}{\rmdefault}{\mddefault}{\updefault}{\color[rgb]{0,0,0}$t$}%
}}}
\end{picture}

%% file: grid2.pstex_t
\begin{picture}(0,0)%
\includegraphics{grid2.pstex}%
\end{picture}%
\setlength{\unitlength}{3947sp}%
\begingroup\makeatletter\ifx\SetFigFont\undefined%
\gdef\SetFigFont#1#2#3#4#5{%
  \reset@font\fontsize{#1}{#2pt}%
  \fontfamily{#3}\fontseries{#4}\fontshape{#5}%
  \selectfont}%
\fi\endgroup%
\begin{picture}(2625,1806)(376,-1486)
\put(2326,-1486){\makebox(0,0)[lb]{\smash{\SetFigFont{12}{14.4}{\rmdefault}{\mddefault}{\updefault}{\color[rgb]{0,0,0}$r$}%
}}}
\put(3001,-1336){\makebox(0,0)[lb]{\smash{\SetFigFont{12}{14.4}{\rmdefault}{\mddefault}{\updefault}{\color[rgb]{0,0,0}$i$}%
}}}
\put(526,164){\makebox(0,0)[lb]{\smash{\SetFigFont{12}{14.4}{\rmdefault}{\mddefault}{\updefault}{\color[rgb]{0,0,0}$j$}%
}}}
\put(376,-436){\makebox(0,0)[lb]{\smash{\SetFigFont{12}{14.4}{\rmdefault}{\mddefault}{\updefault}{\color[rgb]{0,0,0}$t$}%
}}}
\put(376,-1036){\makebox(0,0)[lb]{\smash{\SetFigFont{12}{14.4}{\rmdefault}{\mddefault}{\updefault}{\color[rgb]{0,0,0}$1$}%
}}}
\put(826,-1486){\makebox(0,0)[lb]{\smash{\SetFigFont{12}{14.4}{\rmdefault}{\mddefault}{\updefault}{\color[rgb]{0,0,0}$1$}%
}}}
\end{picture}

%% file: Notation_index.tex
\begin{center}

{\em Les notations de base}

\begin{longtable}{p{3cm}p{8cm}p{2cm}}
${\corps}$ & corps de base & \pageref{corps} \\
${\mathsf C}$, $\mathsf C',$ $\sf C(\mathbb O,\mathbb O)$  & catÈgories monoÔdales ambiantes &  \pageref{C0}, 
\pageref{C1}, \pageref{C(O,O')} \\
$\ts$, $\ts_{\mathbb O}$  & produit tensoriel & \pageref{ts}, \pageref{tsbbO} \\
$e$, $e_{\mathbb O}$ & ÈlÈment neutre pour le produit tensoriel & \pageref{e}, \pageref{ebbO} \\
$\gr \mathsf C$ & catÈgorie des objets graduÈs de ${\mathsf C}$ & \pageref{grC}\\
$\cc \mathsf C$ & catÈgorie des objets diffÈrentiels graduÈs de ${\mathsf C}$ & \pageref{ccC}\\
${\mathsf C}$  & bicatÈgorie  ambiante (‡ partir du chapitre \ref{chapitre_Cat_der}) 
&  \pageref{C2}, \pageref{C3} \\
$S$ &  suspension des objets de $\gr \sf C$ et $\cc \sf C$ & \pageref{S}\\
$s$ & morphisme de foncteurs $ \Id \ra S$ & \pageref{s}\\
$\si = s^{-1} $ & morphismes de foncteur $S \ra \Id$ & \pageref{si} \\
$C(f)$ & cÙne d'un morphisme $f$ & \pageref{Cf} \\
$\del f$ & bord d'un morphisme graduÈ entre deux complexes & \pageref{del}\\
\end{longtable}

{\em Les catÈgories de modËles et catÈgories triangulÈes}

\begin{longtable}{p{3cm}p{8cm}p{2cm}}
$\weq$ & classe des Èquivalences faibles & \pageref{weq} \\
$\cof$ &  classe des cofibrations & \pageref{weq} \\
$\fib$ &  classe des fibrations & \pageref{weq} \\
$\mathsf{E}_{\sf c}$, $\mathsf{E}_{\sf f}$, $\mathsf{E}_{\sf{cf}}$ &
sous-catÈgories des objets cofibrants, des objets fibrants et de objets cofibrants et fibrants de $\mathsf E$
& \pageref{Ecf} \\
$\Ho \mathsf E$ & catÈgorie homotopique de $\mathsf E$ & \pageref{HoE} \\
$\tria \mathbb A$ & sous-catÈgorie triangulÈe engendrÈe par les objets de $\mathbb A$ & \pageref{triabbA} \\
$\Tria \mathbb A$ & sous-catÈgorie triangulÈe (aux sommes infinies)
engendrÈe par les objets de $\mathbb A$ & \pageref{TriabbA} \\
\end{longtable}

 {\em Les $\ai$-algËbres et les algËbres}  

\begin{longtable}{p{3cm}p{8cm}p{2cm}}
${\unite}$ & unitÈ  & \pageref{unite} \\
${\augmentation}$ & morphisme d'augmentation & \pageref{augmentation_mor} \\
$A^+$  &  augmentation de $A$ & \pageref{augmentation} \\
$\b A$  &  rÈduction d'une algËbre augmentÈe & \pageref{reduction} \\
$\b T V$, $TV $  &  algËbre tensorielle rÈduite et augmentÈe  & \pageref{bTV}, \pageref{TV} \\
$BA$ & construction bar rÈduite d'une $\ai$-algËbre & \pageref{B} \\
$\Ba A$ & construction bar augmentÈe d'une $\ai$-algËbre augmentÈe &  \pageref{Ba}, \pageref{Ba1} \\
$b^A$, $b$ & diffÈrentielle de la construction bar & \pageref{b} \\
$(*_m),\hspace*{.2cm} (**_m)$, \hbox{$(***_m)$} & Èquations de type $\ai$ & \pageref{definition_ai-algebre},
\pageref{definition_ai-morphisme}, \pageref{definition_homotopie_ai-morphismes} \\ 
$r(m_\bullet,\hdots, m_{n})$ & cycle mesurant les obstructions & 
\pageref{extension_nstructure_algebres}, \pageref{lemme_obstruction_structure_minimale} \\
$r(f_\bullet,\hdots ,f_n)$ & cycle mesurant les obstructions & \pageref{extension_nstructure_morphismes},
 \pageref{lemme_obstruction_structure_minimale_morphisme} \\
$\E A$ & algËbre enveloppante d'une $\ai$-algËbre $A$ & \pageref{EA}\\
${\alg}$ & catÈgorie des algËbres diffÈrentielles graduÈes & \pageref{alg}\\
${\alga}$ & catÈgorie des algËbres dg augmentÈes dont les morphismes sont augmentÈs & \pageref{alga}\\
${\aia}$ & catÈgorie des $\ai$-algËbres & \pageref{aia}\\
${\aiaa}$ & catÈgorie formÈe des $\ai$-algËbres augmentÈes dont les morphismes sont
augmentÈs & \pageref{aiaa}\\
$\big({\aia}\big)_{hu}$ & catÈgorie des $\ai$-algËbres homologiquement unitaires dont les
morphismes sont homologiquement unitaires & \pageref{aiahu}\\
$\big({\aia}\big)_u$ & catÈgorie des $\ai$-algËbres strictement unitaires dont les
morphismes sont homologiquement unitaires & \pageref{aiau}\\
$\big({\aia}\big)_{su}$ & catÈgorie des $\ai$-algËbres strictement unitaires dont les
morphismes sont strictement unitaires & \pageref{aiasu}\\
$A \rightarrowtail A\langle M,\alpha\rangle $  &  cofibration standard de $\alg$  & \pageref{A<M,a>} \\
$C^*(A,M)$ & complexe de Hochschild de $A$ ‡ coefficients dans $M$ & \pageref{C(A,M)} \\
$\b C^*(A,M)$ & complexe de Hochschild rÈduit & \pageref{bC(A,M)} \\
$\delta_{Hoch}$ & bord de Hochschild & \pageref{delHoch} \\
$\tau$ & cochaÓne tordante & \pageref{tau} \\
$\tau_A$, $\tau_C$ & cochaÓnes tordantes universelles & \pageref{tauA} \\
\end{longtable}

{\em Les $\ai$-cogËbres et les cogËbres }

\begin{longtable}{p{3cm}p{8cm}p{2cm}}
$C \prim n$ & $n$-primitifs de $C$ & \pageref{prim} \\
${\counite}$ & co-unitÈ & \pageref{counite} \\
${\coaugmentation}$ & morphisme de co-augmentation & \pageref{coaugmentation_mor} \\
$C\+$ & co-augmentation d'une cogËbre $C$ & \pageref{coaugmentation} \\
$\b C$ & rÈduction d'une cogËbre co-augmentÈe $C$ & \pageref{reduction_cog} \\
$\ctr V$, $\ct V$ & cogËbre tensorielle rÈduite et co-augmentÈe & \pageref{ctr}, \pageref{ct} \\
$\pctr n V $ & $n$-primitifs de la cogËbre $\ctr C$ & \pageref{pctr} \\
$\Omega C,$ $\Oma C$ & construction cobar rÈduite et co-augmentÈe & \pageref{Omega}, \pageref{Oma},
\pageref{Oma1} \\
${\cog}$  &  catÈgorie des cogËbres diffÈrentielles graduÈes & \pageref{cog} \\
${\cocog}$ &  catÈgorie des cogËbres dg cocomplËtes & \pageref{cogc} \\
$Qis$ & classe des quasi-isomorphismes & \pageref{Qis} \\
$\mbox{\it Qisf\,}$ & classe des quasi-isomorphismes filtrÈs & \pageref{Qisf} \\
${\aic}$  &   catÈgorie des $\ai$-cogËbres  & \pageref{aic} \\
\end{longtable}

 {\em Les polydules, les bipolydules et les modules}  

\begin{longtable}{p{3cm}p{8cm}p{2cm}}
$BM$ & construction bar d'un $A$-polydule & \pageref{BM} \\
$R_\tau M$, $RM$ & produit tensoriel tordu $M \tw C$ & \pageref{RM} \\
$M \tw C$ & produit tensoriel tordu  & \pageref{RM} \\
$r(m_2,\hdots, m_{n})$ & cycle mesurant les obstructions & \pageref{obstruction_extension_polydules} \\
$r(f_1,\hdots ,f_n)$ & cycle mesurant les obstructions & \pageref{extension_nstructure_morphismes_mod} \\
$\Homi_{A'} (X,\?)$ & foncteur standard & \pageref{Homi} \\
$? \tsi_A X$ & foncteur standard & \pageref{tsi} \\
${\Modu A}$ & catÈgorie des $A$-modules diffÈrentiels graduÈs unitaires & \pageref{Modu}\\
${\cd A}$ & catÈgorie dÈrivÈe de $\Modu A$ & \pageref{cdA} \\
${\aiMod A}$  &  catÈgorie des $A$-polydules (non nÈcessairement strictement unitaires) & \pageref{aiMod}\\
${\big(\aiMod A\big)_u}$  &  sous-catÈgorie pleine de $\aiMod A$ formÈe
des $A$-polydules strictement unitaires & \pageref{(aiMod)u}\\
${\aiModst A}$  &  catÈgorie ayant les mÍmes objets que $\aiMod A$ et dont les morphismes sont
  les $\ai$-morphismes stricts & \pageref{aiModst}\\
${\aiModu A}$  &  catÈgorie des $A$-polydules strictement unitaires  & \pageref{aiModu}\\
${\aiModust A}$  &  $\aiModst A \cap \aiModu A$  & \pageref{aiModust}\\
${\ch_\infty A}$ & catÈgorie $\aiModu A$ quotientÈe par la relation d'homotopie & \pageref{chinftyA} \\
${\cd_\infty A}$  & catÈgorie dÈrivÈe de $\aiModu A$ & \pageref{cdinftyA} \\
$\cn_\infty A$ & catÈgorie diffÈrentielle graduÈe des $A$-polydules (non nÈcessairement strictement
unitaires) & \pageref{cninftyA} \\
$\big(\cn_\infty A\big)_u$ & sous-catÈgorie de $\cn_\infty A$ formÈe 
des $A$-polydules strictement unitaires & \pageref{cninftyAu} \\
$\cc_\infty A$ & catÈgorie diffÈrentielle graduÈe des $A$-polydules strictement unitaires & 
\pageref{ccinftyA} \\
${\Modu (A,A')}$ & catÈgorie des $A$-$A'$-bimodules dg unitaires & \pageref{ModuAA'}\\
${\aiMod (A,A')}$ & catÈgorie des $A$-$A'$-polydules (non nÈcessairement
strictement unitaires) & \pageref{aiModAA'}\\
${\big(\aiMod (A,A')\big)_u}$  &  sous-catÈgorie pleine de $\aiMod (A,A')$ formÈe
des $A$-$A'$-bipolydules strictement unitaires & \pageref{(aiMod)uAA'}\\
${\aiModu (A,A')}$ & catÈgorie des $A$-$A'$-polydules strictement unitaires & \pageref{aiModuAA'}\\
${\aiModust (A,A')}$  &  catÈgorie ayant les mÍmes objets que $\aiMod A$ et dont les morphismes sont
les $\ai$-morphismes stricts & \pageref{aiModustAA'}\\
${\ch_\infty (A,A')}$ & catÈgorie $\aiModu (A,A')$ quotientÈe par la relation d'homotopie
& \pageref{chinftyAA'} \\
${\cd_\infty (A,A')}$  & catÈgorie dÈrivÈe de $\aiModu (A,A')$ & \pageref{cdinftyAA'} \\
$\cn_\infty (A,A')$ & catÈgorie diffÈrentielle graduÈe des $A$-$A'$-bipolydules (non nÈcessairement strictement
unitaires) & \pageref{cninftyAA'} \\
$\big(\cn_\infty (A,A')\big)_u$ & sous-catÈgorie pleine de $\cn_\infty (A,A')$ formÈe 
des $A$-$A'$-bipolydules strictement unitaires & \pageref{cninftyAA'u} \\
$\cc_\infty (A,A')$ & catÈgorie diffÈrentielle graduÈe des $A$-$A'$-bipolydules strictement unitaires
& \pageref{ccinftyAA'} \\
\end{longtable}

{\em Les comodules }

\begin{longtable}{p{3cm}p{8cm}p{2cm}}
$N \prim n$ & $n$-primitifs de $N$ & \pageref{primN} \\
$L_\tau N$, $L N$ & produit tensoriel tordu ${N \tw A}$ & \pageref{LN} \\
${N \tw A}$  &  produit tensoriel tordu &\pageref{produit_tens_tordu} \\
$\square_{C}$, $\square$ & produit cotensoriel au-dessus de $C$ & \pageref{square} \\
$\coModu C$ & catÈgorie des comodules dg unitaires & \pageref{comC} \\
${\coModcu C}$ &   catÈgorie des comodules dg cocomplets sur $C$ & \pageref{comcC} \\
$\cd C$ & catÈgorie dÈrivÈe de $\coModcu C$ & \pageref{cdC} \\
\end{longtable}

{\em Les $\ai$-catÈgories et les $\ai$-foncteurs }

\begin{longtable}{p{3cm}p{8cm}p{2cm}}
$\sf C$ & bicatÈgorie des ensembles & \pageref{C3} \\
$\mathbb O$, $\mathbb A$, $\mathbb B$ & ensembles : objets de $\sf C$ & \pageref{bbO} \\
$\ca$, $\cb$ & $\ai$-catÈgories & \pageref{ca} \\
$\mathbb A$, $\mathbb B$ & ensembles des objets des $\ai$-catÈgories $\ca$ et $\cb$ & \pageref{bbA} \\
$\tso$, $\tso_{\mathbb A}$ & produit tensoriel de $\sf C (\mathbb A,\mathbb A)$ & \pageref{tso} \\
$f$, $g$ & $\ai$-foncteurs & \pageref{f} \\
$\dot f$, $\dot g$ & applications sous-jacentes des $\ai$-foncteurs & \pageref{dotf} \\
$\lrus{\dot f}{\cb}{\dot f}$ & & \pageref{fBf} \\
$\Id_\ca$, $\Id$ & $\ai$-foncteur identitÈ de $\ca$ & \pageref{Id} \\
$\id_A$ & morphisme identitÈ d'un objet $A \in \mathbb A$ & \pageref{id} \\
$\ca_x$ & $\ai$-catÈgorie tordue par $x$ & \pageref{cax1}, \pageref{cax} \\
$f^x$ & $\ai$-foncteur tordu par $x$ & \pageref{fx1}, \pageref{fx} \\
$\lrus{x}{M}{x'}$ & bipolydule tordu par $x$ et $x'$ & \pageref{xMx'1}, \pageref{xMx'} \\
$\hat V$ & complÈtion d'un objet topologique & \pageref{hatV} \\
$\hat \ts$ & produit tensoriel complet & \pageref{hatts} \\
$\mathcal{R}$ & catÈgorie des algËbres locales commutatives  & \pageref{calR} \\
$\hat \ct V$ & cogËbre tensorielle complËte rÈduite & \pageref{hatctr} \\
$\TW \ca$ & $\ai$-catÈgorie des objets tordus de $\ca$ & \pageref{TWca} \\
$A \pw$ & polydule reprÈsentÈ $\ca (\?,A)$ & \pageref{Apw} \\
$y$ & $\ai$-foncteur de Yoneda & \pageref{y} \\
$\Nunc (\ca,\cb)$ & $\ai$-catÈgorie des $\ai$-foncteurs $\ca \ra \cb$ 
(non nÈcessairement strictement unitaires) & \pageref{Nunc} \\
$\cf(\ca,\cb)$ & $\ai$-catÈgorie $\Nunc (\ca,\cb)$ munie des compositions naÔves & \pageref{cfcacb} \\
$\big(\Nunc (\ca,\cb)\big)_u$ & sous-catÈgorie pleine de $\Nunc(\ca,\cb)$ formÈe des $\ai$-foncteurs
strictement unitaires & \pageref{Nuncu} \\
$\Func(\ca,\cb)$ & $\ai$-catÈgorie des $\ai$-foncteurs $\ca \ra \cb$ 
strictement unitaires & \pageref{Func} \\
$\boxdot $ & produit cotensoriel & \pageref{boxdot} \\
$\ainat$ & $2$-catÈgorie (non $2$-unitaire) des petites $\ai$-catÈgories (non nÈcessairement strictement
unitaires) & \pageref{ainat} \\
$\aicat$ & $2$-catÈgorie des petites $\ai$-catÈgories strictement unitaires
& \pageref{aicat} \\
$z$ & $\ai$-foncteur de Yoneda gÈnÈralisÈ & \pageref{z} \\
$\theta$ &  & \pageref{theta} \\
$\mathbf{I}_n$ &  & \pageref{In} \\
\end{longtable}

\end{center}

%% file: Page_dos.tex
\pagestyle{empty}

\begin{center}
{\large \bf Sur les $\ai$-catÈgories}\\
Kenji LefËvre-Hasegawa
\end{center}

{\small
{\bf RÈsumÈ : } 
Nous Ètudions les $\ai$-algËbres $\Z$-graduÈes (non nÈcessairement connexes)
et leurs $\ai$-modules. En utilisant les constructions bar et cobar ainsi que les outils de l'algËbre
homotopique de Quillen, nous dÈcrivons la localisation de la catÈgorie des
$\ai$-algËbres par rapport aux $\ai$-quasi-isomorphismes.
Nous adaptons ensuite ces mÈthodes pour dÈcrire la catÈgorie dÈrivÈe
$\cd_\infty A$ d'une $\ai$-algËbre augmentÈe $A$.
Le cas o˘ $A$ n'est pas muni d'une augmentation est traitÈ diffÈremment.
NÈanmoins, lorsque $A$ est strictement unitaire, sa catÈgorie dÈrivÈe peut Ítre dÈcrite
de la mÍme maniËre que dans le cas augmentÈ.
Nous Ètudions ensuite deux variantes de la notion
d'unitaritÈ pour les $\ai$-algËbres~: l'unitaritÈ stricte et l'unitaritÈ homologique. 
Nous montrons que d'un point de vue homotopique, il n'y a pas de diffÈrence
entre ces deux notions. Nous donnons ensuite un formalisme qui permet de
dÈfinir les $\ai$-catÈgories comme des $\ai$-algËbres dans certaines catÈgories monoÔdales.
Nous gÈnÈralisons ‡ ce cadre les constructions fondamentales de la thÈorie des catÈgories :
le foncteur de Yoneda,
les catÈgories de foncteurs, les Èquivalences de catÈgories...
Nous montrons que toute catÈgorie triangulÈe algÈbrique engendrÈe par un ensemble d'objets
est $\ai$-prÈtriangulÈe, c'est-‡-dire qu'elle est Èquivalente ‡ $H^0\TW \ca$, o˘ $\TW \ca$ est
l'$\ai$-catÈgorie des objets tordus d'une certaine $\ai$-catÈgorie $\ca$.\\

{\noindent \bf Discipline :} mathÈmatiques \\
{\bf Mots-clÈs : } $\ai$-catÈgorie, algËbre ‡ homotopie prËs, catÈgorie dÈrivÈe, algËbre homologique, catÈgorie triangulÈe, construction bar
}\\

\begin{center}
{\large \bf On $\ai$-categories}\\
Kenji LefËvre-Hasegawa
\end{center}

{\small
{\bf Abstract : } We study (not necessarily connected) $\Z$-graded $\ai$-algebras and their $\ai$-modules.
Using the cobar and the bar construction and Quillen's homotopical algebra, we describe the localisation
of the category of $\ai$-algebras with respect to $\ai$-quasi-isomorphisms. We then adapt these methods
to describe the derived category $\cd_\infty A$ of an augmented $\ai$-algebra $A$.
The case where $A$ is not endowed with an augmentation is treated differently.
Nevertheless, when $A$ is strictly unital, its derived category can be described in the same way as
in the augmented case. Next, we compare two different notions of $\ai$-unitarity :
strict unitarity and homological unitarity.
We show that, up to homotopy, there is no difference between these two notions.
We then establish a formalism which allows us to view $\ai$-categories as $\ai$-algebras in suitable
monoidal categories. We generalize the fundamental constructions of category theory to this setting : 
Yoneda embeddings, categories of
functors, equivalences of categories... We show that any algebraic
triangulated category $\calt$ which admits a set of
generators is $\ai$-pretriangulated,
that is to say, $\calt$ is equivalent to $H^0 \TW \ca$, where $\TW \ca$ is the $\ai$-category of twisted objets of a certain
$\ai$-category $\ca$. \\

{\bf \noindent Keywords: } $\ai$-categorie, homotopy algebra, derived category, homological algebra, triangulated category, bar construction
}\\[2cm]

{\noindent \bf Adresse : }
Kenji LefËvre-Hasegawa,
ThÈorie des Groupes,
Case 7012,
2 place Jussieu,
F-75251 Paris Cedex 05, France

{\noindent \bf Adresse Èlectronique :} {\tt lefevre@math.jussieu.fr}

\noindent